\documentclass[letterpaper,twoside]{article}
\usepackage[OT1,T1]{fontenc}
\usepackage[dvips]{graphicx,color}
\usepackage{citesort}
\usepackage{pst-3d}
\usepackage{pst-char}
\usepackage{pst-coil}
\usepackage{pst-eps}
\usepackage{pst-fill}
\usepackage{pst-grad}
\usepackage{pst-node}
\usepackage{pst-plot}
\usepackage{pst-text}
\usepackage{pst-tree}
\usepackage{paralist}

\usepackage{ccfonts}
\usepackage{beton}
\usepackage{fancybox}
\usepackage{amssymb}
\usepackage{amsmath,bm}
\usepackage[mathscr]{eucal}
\usepackage{amsfonts}
\usepackage{latexsym}
\usepackage{amsxtra}
\usepackage{amsbsy}
\usepackage{amsthm}
\usepackage{amscd}
\usepackage{amsopn}
\usepackage{amstext}
\usepackage{oldgerm}
\usepackage[varumlaut]{yfonts}
\newcounter{z0}
\newtheorem{ay}{Lemma}[section]

\newtheorem{bbbb}{Proposition}[subsection]
\newtheorem{cccc}{Lemma}[subsection]
\newtheorem{dddd}{Theorem}[subsection]
\newtheorem{ffff}{Corollary}[subsection]
\newtheorem{aaaaa}{Definition}[section]
\newtheorem{bbbbb}{Proposition}[section]
\newtheorem{ccccc}{Lemma}[section]

\newtheorem{fffff}{Corollary}[section]
\newtheorem*{rhp1}{RHP1}
\newtheorem*{rhp2}{RHP2}

\theoremstyle{definition}

\theoremstyle{definition}
\newtheorem{eeee}{Remark}[subsection]

\theoremstyle{definition}
\newtheorem{eeeee}{Remark}[section]

\theoremstyle{definition}
\newtheorem{notrem}{Notational Remark}[section]
\newcommand{\me}{\mathrm{e}}
\newcommand{\mi}{\mathrm{i}}
\newcommand{\md}{\mathrm{d}}
\renewcommand{\Im}{\mathrm{Im}}
\renewcommand{\Re}{\mathrm{Re}}
\newcommand{\id}{\pmb{\mi \md}}
\newcommand{\pvi}{\ensuremath{\int \hspace{-2.75mm} \rule[2.5pt]{2mm}{0.25mm}}}
\newcommand{\vip}{\ensuremath{\int \hspace{-3.35mm} \rule[2.5pt]{2mm}{0.25mm}}}
\newcommand{\norm}[1]{\lVert#1\rVert}
\numberwithin{equation}{section}
\allowdisplaybreaks[4]

\usepackage{relsize,exscale}
\usepackage{stmaryrd}
\usepackage{pxfonts}
\begin{document}
\fontsize{10}{12}\selectfont
\fontencoding{T1}\selectfont
\baselineskip=12pt
\frenchspacing
\title{Asymptotics of Laurent Polynomials of Even Degree Orthogonal with 
Respect to Varying Exponential Weights}
\author{K.~T.-R.~McLaughlin\thanks{\texttt{E-mail: mcl@math.arizona.edu}} \\
Department of Mathematics \\
The University of Arizona \\
617 N.~Santa Rita Ave. \\
P.~O.~Box 210089 \\
Tucson, Arizona 85721--0089 \\
U.~S.~A. \and
A.~H.~Vartanian\thanks{\texttt{E-mail: arthurv@math.ucf.edu}} \\
Department of Mathematics \\
University of Central Florida \\
P.~O.~Box 161364 \\
Orlando, Florida 32816--1364 \\
U.~S.~A. \and
X.~Zhou\thanks{\texttt{E-mail: zhou@math.duke.edu}} \\
Department of Mathematics \\
Duke University \\
Box 90320 \\
Durham, North Carolina 27708--0320 \\
U.~S.~A.}
\date{11 January 2006}
\maketitle
\begin{abstract}
\noindent
Let $\Lambda^{\mathbb{R}}$ denote the linear space over $\mathbb{R}$ spanned 
by $z^{k}$, $k \! \in \! \mathbb{Z}$. Define the real inner product (with 
varying exponential weights) $\langle \boldsymbol{\cdot},\boldsymbol{\cdot} 
\rangle_{\mathscr{L}} \colon \Lambda^{\mathbb{R}} \times \Lambda^{\mathbb{R}} 
\! \to \! \mathbb{R}$, $(f,g) \! \mapsto \! \int_{\mathbb{R}}f(s)g(s) \exp 
(-\mathscr{N} \, V(s)) \, \md s$, $\mathscr{N} \! \in \! \mathbb{N}$, where 
the external field $V$ satisfies: (i) $V$ is real analytic on $\mathbb{R} 
\setminus \{0\}$; (ii) $\lim_{\vert x \vert \to \infty}(V(x)/\ln (x^{2} \! + 
\! 1)) \! = \! +\infty$; and (iii) $\lim_{\vert x \vert \to 0}(V(x)/\ln (x^{
-2} \! + \! 1)) \! = \! +\infty$. Orthogonalisation of the (ordered) base 
$\lbrace 1,z^{-1},z,z^{-2},z^{2},\dotsc,z^{-k},z^{k},\dotsc \rbrace$ with 
respect to $\langle \boldsymbol{\cdot},\boldsymbol{\cdot} \rangle_{\mathscr{
L}}$ yields the even degree and odd degree orthonormal Laurent polynomials 
$\lbrace \phi_{m}(z) \rbrace_{m=0}^{\infty}$: $\phi_{2n}(z) \! = \! \xi^{(2n)
}_{-n}z^{-n} \! + \! \dotsb \! + \! \xi^{(2n)}_{n}z^{n}$, $\xi^{(2n)}_{n} \! 
> \! 0$, and $\phi_{2n+1}(z) \! = \! \xi^{(2n+1)}_{-n-1}z^{-n-1} \! + \! 
\dotsb \! + \! \xi^{(2n+1)}_{n}z^{n}$, $\xi^{(2n+1)}_{-n-1} \! > \! 0$. 
Define the even degree and odd degree monic orthogonal Laurent polynomials: 
$\boldsymbol{\pi}_{2n}(z) \! := \! (\xi^{(2n)}_{n})^{-1} \phi_{2n}(z)$ and 
$\boldsymbol{\pi}_{2n+1}(z) \! := \! (\xi^{(2n+1)}_{-n-1})^{-1}\phi_{2n+1}
(z)$. Asymptotics in the double-scaling limit as $\mathscr{N},n \! \to \! 
\infty$ such that $\mathscr{N}/n \! = \! 1 \! + \! o(1)$ of $\boldsymbol{
\pi}_{2n}(z)$ (in the entire complex plane), $\xi^{(2n)}_{n}$, $\phi_{2n}(z)$ 
(in the entire complex plane), and Hankel determinant ratios associated with 
the real-valued, bi-infinite, strong moment sequence $\left\lbrace c_{k} \! = 
\! \int_{\mathbb{R}}s^{k} \exp (-\mathscr{N} \, V(s)) \, \md s \right\rbrace_{
k \in \mathbb{Z}}$ are obtained by formulating the even degree monic 
orthogonal Laurent po\-l\-y\-n\-o\-m\-i\-a\-l problem as a matrix 
Riemann-Hilbert problem on $\mathbb{R}$, and then extracting the large-$n$ 
behaviour by applying the non-linear steepest-descent method introduced in 
\cite{a1} and further developed in \cite{a2,a3}.

\vspace{0.65cm}

47B36: (Secondary) 30C15, 30C70, 30E05, 30E10, 31A99, 41A20, 41A21, 41A60

\vspace{0.50cm}
{\bf Abbreviated Title.} Asymptotics of Even Degree Orthogonal Laurent 
Polynomials

\vspace{0.50cm}

Laurent-Jacobi matrices, Pad\'{e} approximants, parametrices, Riemann-Hilbert
problems, sing-

ular integral equations, strong moment problems, variational problems
\end{abstract}
\clearpage
\section{Introduction and Background}
Consider the \emph{classical Stieltjes} (resp., \emph{classical Hamburger}) 
\emph{moment problem} (SMP) (resp., HMP): given a simply-infinite (moment) 
sequence of real numbers $\{c_{n}\}_{n=0}^{\infty}$:
\begin{enumerate}
\item[(i)] find necessary and sufficient conditions for the existence of a 
non-negative Borel measure $\mu^{\text{S}}_{\text{MP}}$ (resp., $\mu^{\text{
H}}_{\text{MP}})$ on $[0,+\infty)$ (resp., $(-\infty,+\infty))$, and with 
infinite support, such that $c_{n} \! = \! \int_{0}^{+\infty}t^{n} \, \md 
\mu^{\text{S}}_{\text{MP}}(t)$, $n \! \in \! \mathbb{Z}_{0}^{+} \! := \! 
\{0\} \cup \mathbb{N}$ (resp., $c_{n} \! = \! \int_{-\infty}^{+\infty}t^{n} 
\, \md \mu^{\text{H}}_{\text{MP}}(t)$, $n \! \in \! \mathbb{Z}_{0}^{+})$, 
where the (improper) integral is to be understood in the Riemann-Stieltjes 
sense;
\item[(ii)] when there is a solution of the existence problem, in which case 
the SMP (resp., HMP) is \emph{determinate}, find conditions for the uniqueness 
of the solution; and
\item[(iii)]  when there is more than one solution, in which case the SMP 
(resp., HMP) is \emph{indeterminate}, describe the family of all solutions.
\end{enumerate}
The SMP was---first---treated in $1894/95$ by Stieltjes in the pioneering 
works \cite{a4}, and the HMP was introduced and solved in $1920/21$ by 
Hamburger in the landmark works \cite{a5}. The subsequent development of the 
theory of moment problems brought forth the profound fact that, over and 
above the indispensable utility afforded by the analytic theory of continued 
fractions, in particular, $S$- and real $J$-fractions, the theory of 
orthogonal polynomials \cite{a6} played a seminal, intimate and central 
r\^{o}le (see, for example, \cite{a7}).

Questions regarding two simply-infinite (moment) sequences $\{c_{n}\}_{n \in 
\mathbb{Z}_{0}^{+}}$ and $\{c_{-n}\}_{n \in \mathbb{N}}$ of real numbers, or, 
equivalently, doubly- or bi-infinite (moment) sequences $\{c_{n}\}_{n \in 
\mathbb{Z}}$ of real numbers, manifest, in various settings, purely 
mathematical and/or otherwise, as natural extensions of the foregoing. This 
generalisation is colloquially refered to as the \emph{strong Stieltjes} 
(resp., \emph{strong Hamburger}) \emph{moment problem} (SSMP) (resp., SHMP), 
namely, given a doubly- or bi-infinite (moment) sequence $\{c_{n}\}_{n \in 
\mathbb{Z}}$ of real numbers:
\begin{enumerate}
\item[(1)] find necessary and sufficient conditions for the existence of a 
non-negative measure $\mu^{\text{SS}}_{\text{MP}}$ (resp., $\mu^{\text{SH}}_{
\text{MP}})$ on $[0,+\infty)$ (resp., $(-\infty,+\infty))$, and with infinite 
support, such that $c_{n} \! = \! \int_{0}^{+\infty}t^{n} \, \md \mu^{\text{
SS}}_{\text{MP}}(t)$, $n \! \in \! \mathbb{Z}$ (resp., $c_{n} \! = \! \int_{-
\infty}^{+\infty}t^{n} \, \md \mu^{\text{SH}}_{\text{MP}}(t)$, $n \! \in \! 
\mathbb{Z})$, where the (improper) integral is to be understood in the sense 
of Riemann-Stieltjes;
\item[(2)] when there is a solution, in which case the SSMP (resp., SHMP) is 
determinate, find conditions for the uniqueness of the solution; and
\item[(3)] when there is more than one solution, in which case the SSMP 
(resp., SHMP) is indeterminate, describe the family of all solutions.
\end{enumerate}
The SSMP (resp., SHMP) was introduced in $1980$ (resp., $1981)$ by Jones 
\emph{et al.} \cite{a8} (resp., Jones \emph{et al.} \cite{a9}), and studied 
further in \cite{a10,a11,a12,a13,a14} (see, also, the review article 
\cite{a15}). Unlike the moment theory for the SMP and the HMP, wherein the 
theory of orthogonal polynomials, and the analytic theory of continued 
fractions, enjoyed a prominent r\^{o}le, the extension of the moment theory 
to the SSMP and the SHMP introduced a `rational generalisation' of the 
orthogonal polynomials, namely, the \emph{orthogonal Laurent} (or $L$-) 
\emph{polynomials} (as well as the introduction of special kinds of continued 
fractions commonly referred to as positive-$T$ fractions), which are discussed 
below \cite{a10,a11,a12,a13,a14,a15,a16,a17,a18,a19,a20,a21}. (The SHMP can 
also be solved using the spectral theory of unbounded self-adjoint operators 
in Hilbert space \cite{a22}; see, also, \cite{a23}.)

For any pair $(p,q) \! \in \! \mathbb{Z} \times \mathbb{Z}$, with $p \! 
\leqslant \! q$, let $\Lambda^{\mathbb{C}}_{p,q} \! := \! \left\lbrace 
\mathstrut f \colon \mathbb{C}^{\ast} \! \to \! \mathbb{C}; \, f(z) \! = \! 
\sum_{k=p}^{q} \widehat{\lambda}_{k}z^{k}, \, \widehat{\lambda}_{k} \! \in \! 
\mathbb{C}, \, k \! = \! p,\dotsc,q \right\rbrace$, where $\mathbb{C}^{\ast} 
\! := \! \mathbb{C} \setminus \{0\}$. For any $m \! \in \! \mathbb{Z}_{0}^{
+}$, set $\Lambda^{\mathbb{C}}_{2m} \! := \! \Lambda^{\mathbb{C}}_{-m,m}$, 
$\Lambda^{\mathbb{C}}_{2m+1} \! := \! \Lambda^{\mathbb{C}}_{-m-1,m}$, and 
$\Lambda^{\mathbb{C}} \! := \! \cup_{m \in \mathbb{Z}_{0}^{+}}(\Lambda^{
\mathbb{C}}_{2m} \cup \Lambda^{\mathbb{C}}_{2m+1})$. A function (or element) 
$f \! \in \! \Lambda^{\mathbb{C}}$ is called a \emph{Laurent} (or $L$-) 
\emph{polynomial}. (Note: the sets $\Lambda^{\mathbb{C}}_{p,q}$ and 
$\Lambda^{\mathbb{C}}$ form linear spaces over the field $\mathbb{C}$ with 
respect to the operations of addition and multiplication by a scalar.) Bases 
for each of the spaces $\Lambda^{\mathbb{C}}_{2m}$, $\Lambda^{\mathbb{C}}_{2
m+1}$, and $\Lambda^{\mathbb{C}}$, respectively, are $\lbrace z^{-m},\dotsc,
z^{m} \rbrace$, $\lbrace z^{-m-1},\dotsc,z^{m} \rbrace$, and $\lbrace 
\text{const.},z^{-1},z,z^{-2},z^{2},\dotsc,z^{-k},z^{k},\dotsc \rbrace$ (the 
basis for $\Lambda^{\mathbb{C}}$ corresponds to the \emph{cyclically-repeated 
pole sequence} $\lbrace \text{no pole},0,\infty,0,\infty,\dotsc,0,\infty,
\dotsc \rbrace)$. Furthermore, note that, for each $0 \! \not\equiv \! f \! 
\in \! \Lambda^{\mathbb{C}}$, there exists a unique $l \! \in \! \mathbb{Z}_{
0}^{+}$ such that $f \! \in \! \Lambda^{\mathbb{C}}_{l}$. For $l \! \in \! 
\mathbb{Z}_{0}^{+}$ and $0 \! \not\equiv \! f \! \in \! \Lambda^{\mathbb{C}
}_{l}$, the $L$-\emph{degree} of $f$, symbolically $LD(f)$, is defined as
\begin{equation*}
LD(f) \! := \! l.
\end{equation*}
For $\Lambda^{\mathbb{C}} \! \ni \! f \! = \! \sum_{j \in \mathbb{Z}} \widehat{
\lambda}_{j}z^{j}$, set $C_{j}(f) \! := \! \widehat{\lambda}_{j}$, $j \! \in 
\! \mathbb{Z}$. For each $l \! \in \! \mathbb{Z}_{0}^{+}$ and $0 \! \not\equiv 
\! f \! \in \! \Lambda^{\mathbb{C}}_{l}$, define the \emph{leading 
coefficient} of $f$, symbolically $LC(f)$, and the \emph{trailing coefficient} 
of $f$, symbolically $TC(f)$, as follows:
\begin{equation*}
LC(f) \! := \!
\begin{cases}
\widehat{\lambda}_{m}, &\text{$l \! = \! 2m$,} \\
\widehat{\lambda}_{-m-1}, &\text{$l \! = \! 2m \! + \! 1$,}
\end{cases}
\end{equation*}
and
\begin{equation*}
TC(f) \! := \!
\begin{cases}
\widehat{\lambda}_{-m}, &\text{$l \! = \! 2m$,} \\
\widehat{\lambda}_{m}, &\text{$l \! = \! 2m \! + \! 1$.}
\end{cases}
\end{equation*}
Thus, for $l \! \in \! \mathbb{Z}_{0}^{+}$ and $0 \! \not\equiv \! f \! \in \!
\Lambda^{\mathbb{C}}_{l}$, one writes, for $f \! := \! f_{l}(z)$: (1) if $l \!
= \! 2m$,
\begin{equation*}
f_{2m}(z) \! = \! TC(f)z^{-m} \! + \! \cdots \! + \! LC(f)z^{m};
\end{equation*}
and (2) if $l \! = \! 2m \! + \! 1$,
\begin{equation*}
f_{2m+1}(z) \! = \! LC(f)z^{-m-1} \! + \! \cdots \! + \! TC(f)z^{m}.
\end{equation*}
For $l \! \in \! \mathbb{Z}_{0}^{+}$, $0 \! \not\equiv \! f \! \in \! 
\Lambda^{\mathbb{C}}_{l}$ is called \emph{monic} if $LC(f) \! = \! 1$.

Consider the positive measure on $\mathbb{R}$ (oriented throughout this work, 
unless stated otherwise, {}from $-\infty$ to $+\infty)$ given by
\begin{equation*}
\md \widetilde{\mu}(z) \! = \! \widetilde{w}(z) \md z,
\end{equation*}
with varying exponential weight function of the form
\begin{equation*}
\widetilde{w}(z) \! = \! \exp (-\mathscr{N} \, V(z)), \quad \mathscr{N} \!
\in \! \mathbb{N},
\end{equation*}
where the \emph{external field} $V \colon \mathbb{R} \setminus \{0\} \! \to \! 
\mathbb{R}$ satisfies the following conditions:
\begin{gather}
V \, \, \text{is real analytic on} \, \, \mathbb{R} \setminus \{0\}; \tag{V1}
\\
\lim_{\vert x \vert \to \infty} \! \left(V(x)/\ln (x^{2} \! + \! 1) \right)
\! = \! +\infty; \tag{V2} \\
\lim_{\vert x \vert \to 0} \! \left(V(x)/\ln (x^{-2} \! + \! 1) \right) \! =
\! +\infty. \tag{V3}
\end{gather}
(For example, a rational function of the form $V(z) \! = \! \sum_{k=-2m_{1}}^{
2m_{2}} \! \varrho_{k}z^{k}$, with $\varrho_{k} \! \in \! \mathbb{R}$, $k \! 
= \! -2m_{1},\dotsc,2m_{2}$, $m_{1,2} \! \in \! \mathbb{N}$, and $\varrho_{-
2m_{1}},\varrho_{2m_{2}} \! > \! 0$ would suffice.) Define (uniquely) the 
\emph{strong moment linear functional} $\mathscr{L}$ by its action on the 
basis elements of $\Lambda^{\mathbb{C}}$: $\mathscr{L} \colon \Lambda^{
\mathbb{C}} \! \to \! \Lambda^{\mathbb{C}}$, $f \! = \! \sum_{k \in \mathbb{
Z}} \widehat{\lambda}_{k}z^{k} \! \mapsto \! \mathscr{L}(f) \! := \! \sum_{k
\in \mathbb{Z}} \widehat{\lambda}_{k}c_{k}$, where $c_{k} \! = \! \mathscr{L}
(z^{k}) \! = \! \int_{\mathbb{R}}s^{k} \exp (-\mathscr{N} \, V(s)) \, \md s$,
$(k,\mathscr{N}) \! \in \! \mathbb{Z} \times \mathbb{N}$. (Note that, as per 
the discussion above, $\left\lbrace c_{k} \! = \! \int_{\mathbb{R}}s^{k} \exp 
(-\mathscr{N} \, V(s)) \, \md s, \, \mathscr{N} \! \in \! \mathbb{N} 
\right\rbrace_{k \in \mathbb{Z}}$ is a bi-infinite, real-valued, \emph{strong 
moment sequence}: $c_{k}$ is called the \emph{$k$th strong moment of 
$\mathscr{L}$}.) Associated with the above-defined bi-infinite, real-valued, 
strong moment sequence $\{c_{k}\}_{k \in \mathbb{Z}}$ are the \emph{Hankel 
determinants} $H^{(m)}_{k}$, $(m,k) \! \in \! \mathbb{Z} \times \mathbb{N}$ 
\cite{a10,a11,a15,a17}:
\begin{equation}
H^{(m)}_{0} \! := \! 1 \qquad \quad \text{and} \qquad \quad 
H^{(m)}_{k} \! := \!
\begin{vmatrix}
c_{m} & c_{m+1} & \cdots & c_{m+k-2} & c_{m+k-1} \\
c_{m+1} & c_{m+2} & \cdots & c_{m+k-1} & c_{m+k} \\
c_{m+2} & c_{m+3} & \cdots & c_{m+k} & c_{m+k+1} \\
\vdots & \vdots & \ddots & \vdots & \vdots \\
c_{m+k-1} & c_{m+k} & \cdots & c_{m+2k-3} & c_{m+2k-2}
\end{vmatrix}.
\end{equation}

For any pair $(p,q) \! \in \! \mathbb{Z} \times \mathbb{Z}$, with $p \! 
\leqslant \! q$, let $\Lambda^{\mathbb{R}}_{p,q} \! := \! \left\lbrace 
\mathstrut f \colon \mathbb{C}^{\ast} \! \to \! \mathbb{C}; \, f(z) \! = \! 
\sum_{k=p}^{q} \widetilde{\lambda}_{k}z^{k}, \, \widetilde{\lambda}_{k} \! \in 
\! \mathbb{R}, \, k \! = \! p,\dotsc,q \right\rbrace$, and define, analogously 
as above, for $m \!\in \! \mathbb{Z}_{0}^{+}$, $\Lambda^{\mathbb{R}}_{2m} \! 
:= \! \Lambda^{\mathbb{R}}_{-m,m}$, $\Lambda^{\mathbb{R}}_{2m+1} \! := \! 
\Lambda^{\mathbb{R}}_{-m-1,m}$, and $\Lambda^{\mathbb{R}} \! := \! \cup_{m \in 
\mathbb{Z}_{0}^{+}}(\Lambda^{\mathbb{R}}_{2m} \cup \Lambda^{\mathbb{R}}_{2m+1}
)$. (Note: the sets $\Lambda^{\mathbb{R}}_{p,q}$ and $\Lambda^{\mathbb{R}}$ 
form linear spaces over the field $\mathbb{R}$ with respect to the operations 
of addition and multiplication by a scalar; furthermore, $\Lambda^{\mathbb{R}
}$ $(\subset \Lambda^{\mathbb{C}})$ is the linear space over $\mathbb{R}$ 
spanned by $z^{j}$, $j \! \in \! \mathbb{Z}$.) Hereafter, we shall be 
concerned only with (real) $L$-polynomials in $\Lambda^{\mathbb{R}}$: 
the---ordered---base for $\Lambda^{\mathbb{R}}$ is $\{1,z^{-1},z,z^{-2},z^{2},
\dotsc,z^{-k},z^{k},\dotsc\}$, corresponding to the cyclically-repeated pole 
sequence $\{\text{no pole},0,\infty,0,\infty,\dotsc,0,\infty,\dotsc\}$. Define 
the real bilinear form $\langle \boldsymbol{\cdot},\boldsymbol{\cdot} 
\rangle_{\mathscr{L}}$ as follows: $\langle \boldsymbol{\cdot},\boldsymbol{
\cdot} \rangle_{\mathscr{L}} \colon \Lambda^{\mathbb{R}} \times \Lambda^{
\mathbb{R}} \! \to \! \mathbb{R}$, $(f,g) \! \mapsto \! \langle f,g \rangle_{
\mathscr{L}} \! := \! \mathscr{L}(f(z)g(z)) \! = \! \int_{\mathbb{R}}f(s)g(s) 
\me^{-\mathscr{N} \, V(s)} \, \md s$, $\mathscr{N} \! \in \! \mathbb{N}$. It 
is a fact \cite{a10,a11,a15,a17} that the bilinear form $\langle \boldsymbol{
\cdot},\boldsymbol{\cdot} \rangle_{\mathscr{L}}$ thus defined is an inner 
product if and only if $H^{(-2m)}_{2m} \! > \! 0$ and $H^{(-2m)}_{2m+1} \! > 
\! 0$ $\forall \, \, m \! \in \! \mathbb{Z}_{0}^{+}$ (see Equations~(1.8) 
below, and Subsection~2.2, the proof of Lemma~2.2.1); and this fact is used, 
with little or no further reference, throughout this work (see, also, 
\cite{a24}).
\begin{eeeee}
These latter two (Hankel determinant) inequalities also appear when the 
question of the solvability of the SHMP is posed (in this case, the $c_{k}$, 
$k \! \in \! \mathbb{Z}$, which appear in Equations~(1.1) should be replaced 
by $c_{k}^{\text{\tiny SHMP}}$, $k \! \in \! \mathbb{Z})$: indeed, if these 
two inequalities are true $\forall \, \, m \! \in \! \mathbb{Z}_{0}^{+}$, then 
there is a non-negative measure $\mu^{\text{\tiny SH}}_{\text{\tiny MP}}$ (on 
$\mathbb{R})$ with the given (real) moments. For the case of the SSMP, there 
are four (Hankel determinant) inequalities (in this latter case, the $c_{k}$, 
$k \! \in \! \mathbb{Z}$, which appear in Equations~(1.1) should be replaced 
by $c_{k}^{\text{\tiny SSMP}}$, $k \! \in \! \mathbb{Z})$ which guarantee the 
existence of a non-negative measure $\mu^{\text{\tiny SS}}_{\text{\tiny MP}}$ 
(on $[0,+\infty))$ with the given moments, namely \cite{a8} (see, also, 
\cite{a10,a11}): for each $m \! \in \! \mathbb{Z}_{0}^{+}$, $H^{(-2m)}_{2m} 
\! > \! 0$, $H^{(-2m)}_{2m+1} \! > \! 0$, $H^{(-2m+1)}_{2m} \! > \! 0$, and 
$H^{(-2m-1)}_{2m+1} \! < \! 0$. It is interesting to note that the former 
solvability conditions do not automatically imply that the positive (real) 
moments $\{c_{k}^{\text{\tiny SHMP}}\}_{k \in \mathbb{Z}_{0}^{+}}$ determine 
a measure via the HMP: a similar statement holds true for the SMP (see the 
latter four solvability conditions). \hfill $\blacksquare$
\end{eeeee}
If $f \! \in \! \Lambda^{\mathbb{R}}$, then
\begin{equation*}
\norm{f(\cdot)}_{\mathscr{L}} \! := \! (\langle f,f \rangle_{\mathscr{L}})^{
1/2}
\end{equation*}
is called the \emph{norm of $f$ with respect to $\mathscr{L}$}: note that 
$\norm{f(\cdot)}_{\mathscr{L}} \! \geqslant \! 0 \, \, \forall \, \, f \! \in 
\! \Lambda^{\mathbb{R}}$, and $\norm{f(\cdot)}_{\mathscr{L}} \! > \! 0$ if $0 
\! \not\equiv \! f \! \in \! \Lambda^{\mathbb{R}}$. $\{\phi_{n}^{\flat}(z)\}_{
n \in \mathbb{Z}_{0}^{+}}$ is called a (real) orthonormal Laurent (or $L$-) 
polynomial sequence (ONLPS) with respect to $\mathscr{L}$ if, $\forall \, \, 
m,n \! \in \! \mathbb{Z}_{0}^{+}$:
\begin{enumerate}
\item[(i)] $\phi_{n}^{\flat} \! \in \! \Lambda^{\mathbb{R}}_{n}$, that is,
$LD(\phi_{n}^{\flat}) \! := \! n$;
\item[(ii)] $\langle \phi_{m}^{\flat},\phi_{n^{\prime}}^{\flat} \rangle_{
\mathscr{L}} \! = \! 0 \, \, \forall \, \, m \! \not= \! n^{\prime}$, or,
alternatively, $\langle f,\phi_{n}^{\flat} \rangle_{\mathscr{L}} \! = \! 0 \,
\, \forall \, \, f \! \in \! \Lambda^{\mathbb{R}}_{n-1}$;
\item[(iii)] $\langle \phi_{m}^{\flat},\phi_{m}^{\flat} \rangle_{\mathscr{L}}
\! =: \! \norm{\phi_{m}^{\flat}(\cdot)}^{2}_{\mathscr{L}} \! = \! 1$.
\end{enumerate}
Orthonormalisation of $\lbrace 1,z^{-1},z,z^{-2},z^{2},\dotsc,z^{-n},z^{n},
\dotsc \rbrace$, corresponding to the cyclically-repeated pole sequence 
$\lbrace \text{no pole},0,\infty,0,\infty,\dotsc,0,\infty,\dotsc \rbrace$, 
with respect to $\langle \boldsymbol{\cdot},\boldsymbol{\cdot} \rangle_{
\mathscr{L}}$ via the Gram-Schmidt orthogonalisation method, leads to the 
ONLPS, or, simply, orthonormal Laurent (or $L$-) polynomials (OLPs), $\{
\phi_{m}(z)\}_{m \in \mathbb{Z}_{0}^{+}}$, which, by suitable normalisation, 
may be written as, for $m \! = \! 2n$,
\begin{equation}
\phi_{2n}(z) \! = \! \xi^{(2n)}_{-n}z^{-n} \! + \! \dotsb \! + \!
\xi^{(2n)}_{n}z^{n}, \qquad \xi^{(2n)}_{n} \! > \! 0,
\end{equation}
and, for $m \! = \! 2n \! + \! 1$,
\begin{equation}
\phi_{2n+1}(z) \! = \! \xi^{(2n+1)}_{-n-1}z^{-n-1} \! + \! \dotsb \! + \!
\xi^{(2n+1)}_{n}z^{n}, \qquad \xi^{(2n+1)}_{-n-1} \! > \! 0.
\end{equation}
The $\phi_{n}$'s are normalised so that they all have real coefficients; in 
particular, the leading coefficients, $LC(\phi_{2n}) \! := \! \xi^{(2n)}_{n}$ 
and $LC(\phi_{2n+1}) \! := \! \xi^{(2n+1)}_{-n-1}$, $n \! \in \! \mathbb{Z}_{
0}^{+}$, are both positive, $\xi^{(0)}_{0} \! = \! 1$, and $\phi_{0}(z) \! 
\equiv \! 1$. Even though the leading coefficients, $\xi^{(2n)}_{n}$ and 
$\xi^{(2n+1)}_{-n-1}$, $n \! \in \! \mathbb{Z}_{0}^{+}$, are non-zero (in 
particular, they are positive), no such restriction applies to the trailing 
coefficients, $TC(\phi_{2n}) \! := \! \xi^{(2n)}_{-n}$ and $TC(\phi_{2n+1}) 
\! := \! \xi^{(2n+1)}_{n}$, $n \! \in \! \mathbb{Z}_{0}^{+}$. Furthermore, 
note that, by construction:
\begin{enumerate}
\item[(1)] $\langle \phi_{2n},z^{j} \rangle_{\mathscr{L}} \! = \! 0$, $j \! = 
\! -n,\dotsc,n \! - \! 1$;
\item[(2)] $\langle \phi_{2n+1},z^{j} \rangle_{\mathscr{L}} \! = \! 0$, $j \! 
= \! -n,\dotsc,n$;
\item[(3)] $\langle \phi_{j},\phi_{k} \rangle_{\mathscr{L}} \! = \! \delta_{j
k}$, $j,k \! \in \! \mathbb{Z}_{0}^{+}$, where $\delta_{jk}$ is the Kronecker 
delta.
\end{enumerate}
Moreover, if, for each $m \! \in \! \mathbb{Z}_{0}^{+}$, the orthonormal 
$L$-polynomials $\phi_{2m}(z)$ and $\phi_{2m+1}(z)$, respectively, are such 
that $TC(\phi_{2m}) \! := \! \xi^{(2m)}_{-m} \! \not= \! 0$ and $TC(\phi_{2
m+1}) \! := \! \xi^{(2m+1)}_{m} \! \not= \! 0$, then there are special 
Christoffel-Darboux formulae for the OLPs (see, for example, \cite{a12,a17};
see, also, \cite{a25}):
\begin{gather*}
\phi_{2m}(\zeta)(z \phi_{2m-1}(z) \! - \! \zeta \phi_{2m-1}(\zeta)) \! - \!
\zeta \phi_{2m-1}(\zeta)(\phi_{2m}(z) \! - \! \phi_{2m}(\zeta)) \! = \! (z \!
- \! \zeta) \dfrac{\xi^{(2m)}_{-m}}{\xi^{(2m-1)}_{-m}} \, \sum_{j=0}^{2m-1}
\phi_{j}(z) \phi_{j}(\zeta), \\
\phi_{2m}(\zeta)(z \phi_{2m+1}(z) \! - \! \zeta \phi_{2m+1}(\zeta)) \! - \!
\zeta \phi_{2m+1}(\zeta)(\phi_{2m}(z) \! - \! \phi_{2m}(\zeta)) \! = \! (z \!
- \! \zeta) \dfrac{\xi^{(2m+1)}_{m}}{\xi^{(2m)}_{m}} \, \sum_{j=0}^{2m} \phi_{
j}(z) \phi_{j}(\zeta),
\end{gather*}
where $\phi_{-1}(z) \! \equiv \! 0$, and (dividing by $z \! - \! \zeta$
and letting $\zeta \! \to \! z)$
\begin{gather*}
\phi_{2m}(z) \dfrac{\md}{\md z}(z \phi_{2m-1}(z)) \! - \! z \phi_{2m-1}(z)
\dfrac{\md}{\md z} \phi_{2m}(z) \! = \! \dfrac{\xi^{(2m)}_{-m}}{\xi^{(2m-1)}_{
-m}} \, \sum_{j=0}^{2m-1}(\phi_{j}(z))^{2}, \\
\phi_{2m}(z) \dfrac{\md}{\md z}(z \phi_{2m+1}(z)) \! - \! z \phi_{2m+1}(z)
\dfrac{\md}{\md z} \phi_{2m}(z) \! = \! \dfrac{\xi^{(2m+1)}_{m}}{\xi^{(2m)}_{
m}} \, \sum_{j=0}^{2m}(\phi_{j}(z))^{2}.
\end{gather*}

It is convenient to introduce the monic orthogonal Laurent (or $L$-) 
polynomials, $\boldsymbol{\pi}_{j}(z)$, $j \! \in \! \mathbb{Z}_{0}^{+}$: (i) 
for $j \! = \! 2n$, $n \! \in \! \mathbb{Z}_{0}^{+}$, with $\boldsymbol{\pi}_{
0}(z) \! \equiv \! 1$,
\begin{equation}
\boldsymbol{\pi}_{2n}(z) \! := \! \phi_{2n}(z)(\xi^{(2n)}_{n})^{-1} \! = \!
\nu^{(2n)}_{-n}z^{-n} \! + \! \dotsb \! + \! z^{n}, \qquad \quad \nu^{(2n)}_{
-n} \! := \! \xi^{(2n)}_{-n}/\xi^{(2n)}_{n};
\end{equation}
and (ii) for $j \! = \! 2n \! + \! 1$, $n \! \in \! \mathbb{Z}_{0}^{+}$,
\begin{equation}
\boldsymbol{\pi}_{2n+1}(z) \! := \! \phi_{2n+1}(z)(\xi^{(2n+1)}_{-n-1})^{-1}
\! = \! z^{-n-1} \! + \! \dotsb \! + \! \nu^{(2n+1)}_{n}z^{n}, \qquad \quad
\nu^{(2n+1)}_{n} \! := \! \xi^{(2n+1)}_{n}/\xi^{(2n+1)}_{-n-1}.
\end{equation}
The monic orthogonal $L$-polynomials, $\boldsymbol{\pi}_{j}(z)$, $j \! \in \! 
\mathbb{Z}_{0}^{+}$, possess the following properties:
\begin{enumerate}
\item[(1)] $\langle \boldsymbol{\pi}_{2n},z^{j} \rangle_{\mathscr{L}} \! = \!
0$, $j \! = \! -n,\dotsc,n \! - \! 1$;
\item[(2)] $\langle \boldsymbol{\pi}_{2n+1},z^{j} \rangle_{\mathscr{L}} \! =
\! 0$, $j \! = \! -n,\dotsc,n$;
\item[(3)] $\langle \boldsymbol{\pi}_{2n},\boldsymbol{\pi}_{2n} \rangle_{
\mathscr{L}} \! =: \! \norm{\boldsymbol{\pi}_{2n}(\cdot)}^{2}_{\mathscr{L}}
\! = \! (\xi^{(2n)}_{n})^{-2}$, whence $\xi^{(2n)}_{n} \! = \! 1/\norm{
\boldsymbol{\pi}_{2n}(\cdot)}_{\mathscr{L}}$ $(> \! 0)$;
\item[(4)] $\langle \boldsymbol{\pi}_{2n+1},\boldsymbol{\pi}_{2n+1} \rangle_{
\mathscr{L}} \! =: \! \norm{\boldsymbol{\pi}_{2n+1}(\cdot)}^{2}_{\mathscr{L}}
\! = \! (\xi^{(2n+1)}_{-n-1})^{-2}$, whence $\xi^{(2n+1)}_{-n-1} \! = \!
1/\norm{\boldsymbol{\pi}_{2n+1}(\cdot)}_{\mathscr{L}}$ $(> \! 0)$.
\end{enumerate}
Furthermore, in terms of the Hankel determinants, $H^{(m)}_{k}$, $(m,k) \! 
\in \! \mathbb{Z} \times \mathbb{N}$, associated with the real-valued, 
bi-infinite, strong moment sequence $\left\lbrace c_{k} \! = \! \int_{
\mathbb{R}}s^{k} \me^{-\mathscr{N} \, V(s)} \, \md s, \, \mathscr{N} \! \in 
\! \mathbb{N} \right\rbrace_{k \in \mathbb{Z}}$, the monic orthogonal 
$L$-polynomials, $\boldsymbol{\pi}_{j}(z)$, $j \! \in \! \mathbb{Z}_{0}^{+}$, 
are represented via the following determinantal formulae 
\cite{a10,a11,a15,a17} (see, also, Subsection~2.2, Proposition~2.2.1): for $m 
\! \in \! \mathbb{Z}_{0}^{+}$,
\begin{gather}
\boldsymbol{\pi}_{2m}(z) \! = \! \dfrac{1}{H^{(-2m)}_{2m}}
\begin{vmatrix}
c_{-2m} & c_{-2m+1} & \cdots & c_{-1} & z^{-m} \\
c_{-2m+1} & c_{-2m+2} & \cdots & c_{0} & z^{-m+1} \\
\vdots & \vdots & \ddots & \vdots & \vdots \\
c_{-1} & c_{0} & \cdots & c_{2m-2} & z^{m-1} \\
c_{0} & c_{1} & \cdots & c_{2m-1} & z^{m}
\end{vmatrix}, \\
\intertext{and}
\boldsymbol{\pi}_{2m+1}(z) \! = \! -\dfrac{1}{H^{(-2m)}_{2m+1}}
\begin{vmatrix}
c_{-2m-1} & c_{-2m} & \cdots & c_{-1} & z^{-m-1} \\
c_{-2m} & c_{-2m+1} & \cdots & c_{0} & z^{-m} \\
\vdots & \vdots & \ddots & \vdots & \vdots \\
c_{-1} & c_{0} & \cdots & c_{2m-1} & z^{m-1} \\
c_{0} & c_{1} & \cdots & c_{2m} & z^{m}
\end{vmatrix};
\end{gather}
moreover, it can be shown that (see, for example, \cite{a15,a17}), for $n \! 
\in \! \mathbb{Z}_{0}^{+}$,
\begin{gather}
\xi^{(2n)}_{n} \! \left(= \! \dfrac{1}{\norm{\boldsymbol{\pi}_{2n}(\cdot)}_{
\mathscr{L}}} \right) \! = \! \sqrt{\dfrac{H^{(-2n)}_{2n}}{H^{(-2n)}_{2n+1}}},
\qquad \xi^{(2n+1)}_{-n-1} \! \left(= \! \dfrac{1}{\norm{\boldsymbol{\pi}_{2n
+1}(\cdot)}_{\mathscr{L}}} \right) \! = \! \sqrt{\dfrac{H^{(-2n)}_{2n+1}}{H^{
(-2n-2)}_{2n+2}}}, \\
\nu^{(2n)}_{-n} \! \left(:= \! \dfrac{\xi^{(2n)}_{-n}}{\xi^{(2n)}_{n}} \right)
\! = \! \dfrac{H^{(-2n+1)}_{2n}}{H^{(-2n)}_{2n}}, \qquad \quad \nu^{(2n+1)}_{
n} \! \left(:= \! \dfrac{\xi^{(2n+1)}_{n}}{\xi^{(2n+1)}_{-n-1}} \right) \! =
\! -\dfrac{H^{(-2n-1)}_{2n+1}}{H^{(-2n)}_{2n+1}}.
\end{gather}

For each $m \! \in \! \mathbb{Z}_{0}^{+}$, the monic orthogonal $L$-polynomial 
$\boldsymbol{\pi}_{m}(z)$ and the index $m$ are called \emph{non-singular} if 
$0 \! \not= \! TC(\boldsymbol{\pi}_{m}) \! := \!
\begin{cases}
\nu^{(2n)}_{-n}, &\text{$m \! = \! 2n$,} \\
\nu^{(2n+1)}_{n}, &\text{$m \! = \! 2n \! + \! 1$;}
\end{cases}$ otherwise, $\boldsymbol{\pi}_{m}(z)$ and $m$ are \emph{singular}. 
{}From Equations~(1.9), it can be seen that, for each $m \! \in \! \mathbb{
Z}_{0}^{+}$:
\begin{enumerate}
\item[(i)] $\boldsymbol{\pi}_{2m}(z)$ is non-singular (resp., singular) if
$H^{(-2m+1)}_{2m} \! \not= \! 0$ (resp., $H^{(-2m+1)}_{2m} \! = \! 0)$;
\item[(ii)] $\boldsymbol{\pi}_{2m+1}(z)$ is non-singular (resp., singular) if
$H^{(-2m-1)}_{2m+1} \! \not= \! 0$ (resp., $H^{(-2m-1)}_{2m+1} \! = \! 0)$.
\end{enumerate}
For each $m \! \in \! \mathbb{Z}_{0}^{+}$, let $\mu_{2m} \! := \! 
\operatorname{card} \lbrace \mathstrut z; \, \boldsymbol{\pi}_{2m}(z) \! = \! 
0 \rbrace$ and $\mu_{2m+1} \! := \! \operatorname{card} \lbrace \mathstrut z; 
\, \boldsymbol{\pi}_{2m+1}(z) \! = \! 0 \rbrace$. It is an established fact 
\cite{a10,a11,a17} that, for $m \! \in \! \mathbb{Z}_{0}^{+}$:
\begin{enumerate}
\item[(1)] the zeros of $\boldsymbol{\pi}_{2m}(z)$ are real, simple, and 
non-zero, and $\mu_{2m} \! = \! 2m$ (resp., $2m \! - \! 1)$ if $\boldsymbol{
\pi}_{2m}(z)$ is non-singular (resp., singular);
\item[(2)] the zeros of $\boldsymbol{\pi}_{2m+1}(z)$ are real, simple, and 
non-zero, and $\mu_{2m+1} \! = \! 2m \! + \! 1$ (resp., $2m)$ if $\boldsymbol{
\pi}_{2m+1}(z)$ is non-singular (resp., singular).
\end{enumerate}
For each $m \! \in \! \mathbb{Z}_{0}^{+}$, it can be shown that, via a 
straightforward factorisation argument and using Equations~(1.6) and~(1.7):
\begin{enumerate}
\item[(i)] if $\boldsymbol{\pi}_{2m}(z)$ is non-singular, upon setting 
$\left\{\mathstrut \alpha^{(2m)}_{k}, \, k \! = \! 1,\dotsc,2m \right\} \! 
:= \! \left\{\mathstrut z; \, \boldsymbol{\pi}_{2m}(z) \! = \! 0 \right\}$,
\begin{equation*}
\prod_{k=1}^{2m} \alpha^{(2m)}_{k} \! = \nu^{(2m)}_{-m};
\end{equation*}
\item[(ii)] if $\boldsymbol{\pi}_{2m+1}(z)$ is non-singular, upon setting 
$\left\{\mathstrut \alpha^{(2m+1)}_{k}, \, k \! = \! 1,\dotsc,2m \! + \! 1 
\right\} \! := \! \left\{\mathstrut z; \, \boldsymbol{\pi}_{2m+1}(z) \! = \! 
0 \right\}$,
\begin{equation*}
\prod_{k=1}^{2m+1} \alpha^{(2m+1)}_{k} \! = \! - \! \left(\nu^{(2m+1)}_{m} 
\right)^{-1}.
\end{equation*}
\end{enumerate}

Unlike orthogonal polynomials, which satisfy a system of three-term recurrence 
relations, monic orthogonal, and orthonormal, $L$-polynomials may satisfy 
recurrence relations consisting of a pair of four-term recurrence relations 
\cite{a15}, a pair of systems of three- or five-term recurrence relations 
(which is guaranteed in the case when the corresponding monic orthogonal, 
and orthonormal, $L$-polynomials are non-singular) \cite{a15,a16,a17}, or 
a system consisting of four five-term recurrence relations \cite{a23}.
\begin{eeeee}
The non-vanishing of the leading and trailing coefficients of the OLPs 
$\lbrace \phi_{m}(z) \rbrace_{m=0}^{\infty}$, that is,
\begin{equation*}
LC(\phi_{m}) \! := \!
\begin{cases}
\xi^{(2n)}_{n}, &\text{$m \! = \! 2n$,} \\
\xi^{(2n+1)}_{-n-1}, &\text{$m \! = \! 2n \! + \! 1$,}
\end{cases}
\end{equation*}
and
\begin{equation*}
TC(\phi_{m}) \! := \!
\begin{cases}
\xi^{(2n)}_{-n}, &\text{$m \! = \! 2n$,} \\
\xi^{(2n+1)}_{n}, &\text{$m \! = \! 2n \! + \! 1$,}
\end{cases}
\end{equation*}
respectively, is of paramount importance: if both these conditions are not 
satisfied, then the `length' of the recurrence relations may be greater than 
three \cite{a16} (see, also, \cite{a24}). \hfill $\blacksquare$
\end{eeeee}
It can be shown that (see, for example, \cite{a17}, and Chapter~11 of 
\cite{a26}), if $\{\boldsymbol{\pi}_{m}(z)\}_{m \in \mathbb{Z}_{0}^{+}}$, as 
defined above, is a non-singular, monic orthogonal $L$-polynomial sequence, 
that is, $H^{(-2n+1)}_{2n} \! \not= \! 0$ $(m \! = \! 2n)$ and $H^{(-2n-1)}_{
2n+1} \! \not= \! 0$ $(m \! = \! 2n \! + \! 1)$, then $\{\boldsymbol{\pi}_{m}
(z)\}_{m \in \mathbb{Z}_{0}^{+}}$ satisfy the pair of three-term recurrence 
relations
\begin{align*}
\boldsymbol{\pi}_{2m+1}(z) &= \left(\dfrac{z^{-1}}{\beta_{2m}^{\natural}} \!
+ \! \beta_{2m+1}^{\natural} \right) \! \boldsymbol{\pi}_{2m}(z) \! + \!
\lambda_{2m+1}^{\natural} \boldsymbol{\pi}_{2m-1}(z), \\
\boldsymbol{\pi}_{2m+2}(z) &= \left(\dfrac{z}{\beta_{2m+1}^{\natural}} \! + \!
\beta_{2m+2}^{\natural} \right) \! \boldsymbol{\pi}_{2m+1}(z) \! + \! \lambda_{
2m+2}^{\natural} \boldsymbol{\pi}_{2m}(z),
\end{align*}
where $\boldsymbol{\pi}_{-1}(z) \! \equiv \! 0$,
\begin{gather*}
\beta_{2m}^{\natural}= \nu^{(2m)}_{-m}, \qquad \qquad \beta_{2m+1}^{\natural}=
\nu^{(2m+1)}_{m}, \\
\lambda_{2m+1}^{\natural} \! = \! -\dfrac{H^{(-2m-1)}_{2m+1} H^{(-2m+2)}_{2m-
1}}{H^{(-2m)}_{2m} H^{(-2m+1)}_{2m}} \quad (\not= \! 0), \qquad \qquad 
\lambda_{2m+2}^{\natural} \! = \! -\dfrac{H^{(-2m-1)}_{2m+2} H^{(-2m)}_{2m}
}{H^{(-2m)}_{2m+1} H^{(-2m-1)}_{2m+1}} \quad (\not= \! 0),
\end{gather*}
and $\lambda_{j} \beta_{j-1}/\beta_{j} \! > \! 0  \, \, \forall \, \, j \! 
\in \! \mathbb{N}$, with $\lambda_{1} \! := \! -c_{-1}$, leading to a 
\emph{tri-diagonal-type Laurent-Jacobi matrix} $\mathcal{F}$ for the `mixed' 
mapping
\begin{equation*}
\mathscr{F} \colon \Lambda^{\mathbb{R}} \! \to \! \Lambda^{\mathbb{R}}, \, \, 
f(z) \! \mapsto \! (z^{-1}(\oplus_{n=0}^{\infty} \operatorname{diag}(1,0)) \! 
+ \! z(\oplus_{n=0}^{\infty} \operatorname{diag}(0,1)))f(z),
\end{equation*}
where $\oplus_{n=0}^{\infty} \operatorname{diag}(1,0) \! := \! \operatorname{
diag}(1,0,\dotsc,1,0,\dotsc)$, and $\oplus_{n=0}^{\infty} \operatorname{diag}
(0,1) \! := \! \operatorname{diag}(0,1,\dotsc,0,1,\dotsc)$,
\begin{align*}
\setcounter{MaxMatrixCols}{13}
\mathcal{F} =& \, \operatorname{diag} \! \left(\beta_{0}^{\natural},\beta_{
1}^{\natural},\beta_{2}^{\natural},\dotsc \right) \!
\left(
\begin{smallmatrix}
-\beta_{1}^{\natural} & 1 & & & & & & & & & & & \\
-\lambda_{2}^{\natural} & -\beta_{2}^{\natural} & 1 & & & & & & & & & & \\
 & -\lambda_{3}^{\natural} & -\beta_{3}^{\natural} & 1 & & & & & & & & & \\
 & & -\lambda_{4}^{\natural} & -\beta_{4}^{\natural} & 1 & & & & & & & & \\
 & & & -\lambda_{5}^{\natural} & -\beta_{5}^{\natural} & 1 & & & & & & & \\
 & & & & -\lambda_{6}^{\natural} & -\beta_{6}^{\natural} & 1 & & & & & & \\
 & & & & & \ddots & \ddots & \ddots & & & & & \\
 & & & & & & -\lambda_{2m+1}^{\natural} & -\beta_{2m+1}^{\natural} & 1 & & & &
& \\
 & & & & & & & -\lambda_{2m+2}^{\natural} & -\beta_{2m+2}^{\natural} & 1 & & &
\\
 & & & & & & & & \ddots & \ddots & \ddots & &
\end{smallmatrix}
\right),
\end{align*}
with zeros outside the indicated diagonals (in terms of $\{\phi_{m}(z)\}_{m 
\in \mathbb{Z}_{0}^{+}}$, the pair of three-term recurrence relations reads 
\cite{a16}:
\begin{gather*}
\phi_{2m+1}(z) \! = \! (z^{-1} \! + \! \mathfrak{g}_{2m+1}) \phi_{2m}(z) \! +
\! \mathfrak{f}_{2m+1} \phi_{2m-1}(z), \\
\phi_{2m+2}(z) \! = \! (1 \! + \! \mathfrak{g}_{2m+2}z) \phi_{2m+1}(z) \! +
\! \mathfrak{f}_{2m+2} \phi_{2m}(z),
\end{gather*}
where $\mathfrak{f}_{2m+1},\mathfrak{f}_{2m+2} \! \not= \! 0$, $m \! \in \! 
\mathbb{Z}_{0}^{+}$, $\phi_{-1}(z) \! \equiv \! 0$, and $\phi_{0}(z) \! \equiv 
\! 1)$; otherwise, $\{\boldsymbol{\pi}_{m}(z)\}_{m \in \mathbb{Z}_{0}^{+}}$ 
satisfy the following pair of five-term recurrence relations \cite{a17}, with 
$\boldsymbol{\pi}_{-j}(z) \! \equiv \! 0$, $j \! = \! 1,2$,
\begin{align*}
\boldsymbol{\pi}_{2m+2}(z) =& \, \gamma_{2m+2,2m-2}^{\flat} \boldsymbol{\pi}_{
2m-2}(z) \! + \! \gamma_{2m+2,2m-1}^{\flat} \boldsymbol{\pi}_{2m-1}(z) \! + \!
(z \! + \! \gamma_{2m+2,2m}^{\flat}) \boldsymbol{\pi}_{2m}(z) \\
+& \, \gamma_{2m+2,2m+1}^{\flat} \boldsymbol{\pi}_{2m+1}(z), \\
\boldsymbol{\pi}_{2m+3}(z) =& \, \gamma_{2m+3,2m-1}^{\flat} \boldsymbol{\pi}_{
2m-1}(z) \! + \! \gamma_{2m+3,2m}^{\flat} \boldsymbol{\pi}_{2m}(z) \! + \!
(z^{-1} \! + \! \gamma_{2m+3,2m+1}^{\flat}) \boldsymbol{\pi}_{2m+1}(z) \\
+& \, \gamma_{2m+3,2m+2}^{\flat} \boldsymbol{\pi}_{2m+2}(z),
\end{align*}
where $\gamma_{l,k} \! = \! 0$, $k \! < \! 0$, $l \! \geqslant \! 2$, leading 
to a \emph{penta-diagonal-type Laurent-Jacobi matrix} $\mathcal{G}$ for the 
`mixed' mapping
\begin{equation*}
\mathscr{G} \colon \Lambda^{\mathbb{R}} \! \to \! \Lambda^{\mathbb{R}}, \, \, 
g(z) \! \mapsto \! (z(\oplus_{n=0}^{\infty} \operatorname{diag}(1,0)) \! + \! 
z^{-1}(\oplus_{n=0}^{\infty} \operatorname{diag}(0,1)))g(z),
\end{equation*}
\begin{align*}
\setcounter{MaxMatrixCols}{14}
&\mathcal{G} =
\left(
\begin{smallmatrix}
-\gamma_{2,0}^{\flat} & -\gamma_{2,1}^{\flat} & 1 & & & & & & & & & & & \\
-\gamma_{3,0}^{\flat} & -\gamma_{3,1}^{\flat} & -\gamma_{3,2}^{\flat} & 1 & &
& & & & & & & & \\
-\gamma_{4,0}^{\flat} & -\gamma_{4,1}^{\flat} & -\gamma_{4,2}^{\flat} &
-\gamma_{4,3}^{\flat} & 1 & & & & & & & & \\
 & -\gamma_{5,1}^{\flat} & -\gamma_{5,2}^{\flat} & -\gamma_{5,3}^{\flat} &
-\gamma_{5,4}^{\flat} & 1 & & & & & & & & \\
 & & -\gamma_{6,2}^{\flat} & -\gamma_{6,3}^{\flat} & -\gamma_{6,4}^{\flat} &
-\gamma_{6,5}^{\flat} & 1 & & & & & & & \\
 & & & \ddots & \ddots & \ddots & \ddots & \ddots & & & & & \\
 & & & & -\gamma_{2m+2,2m-2}^{\flat} & -\gamma_{2m+2,2m-1}^{\flat} & -\gamma_{
2m+2,2m}^{\flat} & -\gamma_{2m+2,2m+1}^{\flat} & 1 & & & & & \\
 & & & & & -\gamma_{2m+3,2m-1}^{\flat} & -\gamma_{2m+3,2m}^{\flat} & -\gamma_{
2m+3,2m+1}^{\flat} & -\gamma_{2m+3,2m+2}^{\flat} & 1 & & & & \\
 & & & & & & \ddots & \ddots & \ddots & \ddots & \ddots & & &
\end{smallmatrix}
\right),
\end{align*}
with zeros outside the indicated diagonals. The general form of these (system 
of) recurrence relations is a pair of three- and five-term recurrence 
relations \cite{a23}: for $n \! \in \! \mathbb{Z}_{0}^{+}$,
\begin{gather*}
z \phi_{2n+1}(z) \! = \! b_{2n+1}^{\sharp} \phi_{2n}(z) \! + \! a_{2n+1}^{
\sharp} \phi_{2n+1}(z) \! + \! b_{2n+2}^{\sharp} \phi_{2n+2}(z), \\
z \phi_{2n}(z) \! = \! c_{2n}^{\sharp} \phi_{2n-2}(z) \! + \! b_{2n}^{\sharp}
\phi_{2n-1}(z) \! + \! a_{2n}^{\sharp} \phi_{2n}(z) \! + \! b_{2n+1}^{\sharp}
\phi_{2n+1}(z) \! + \! c_{2n+2}^{\sharp} \phi_{2n+2}(z),
\end{gather*}
where all the coefficients are real, $c_{0}^{\sharp} \! = \! b_{0}^{\sharp} 
\! = \! 0$, and $c_{2k}^{\sharp} \! > \! 0$, $k \! \in \! \mathbb{N}$, and
\begin{gather*}
z^{-1} \phi_{2n}(z) \! = \! \beta_{2n}^{\sharp} \phi_{2n-1}(z) \! + \!
\alpha_{2n}^{\sharp} \phi_{2n}(z) \! + \! \beta_{2n+1}^{\sharp} \phi_{2n+1}
(z), \\
z^{-1} \phi_{2n+1}(z) \! = \! \gamma_{2n+1}^{\sharp} \phi_{2n-1}(z) \! + \!
\beta_{2n+1}^{\sharp} \phi_{2n}(z) \! + \! \alpha_{2n+1}^{\sharp} \phi_{2n+1}
(z) \! + \! \beta_{2n+2}^{\sharp} \phi_{2n+2}(z) \! + \! \gamma_{2n+3}^{
\sharp} \phi_{2n+3}(z),
\end{gather*}
where all the coefficients are real, $\beta_{0}^{\sharp} \! = \! \gamma_{1}^{
\sharp} \! = \! 0$, $\beta_{1}^{\sharp} \! > \! 0$, and $\gamma_{2l+1}^{
\sharp} \! > \! 0$, $l \! \in \! \mathbb{N}$, leading, respectively, to the 
real-symmetric, \emph{tri-penta-diagonal-type Laurent-Jacobi matrices}, 
$\mathcal{J}$ and $\mathcal{K}$, for the mappings
\begin{equation*}
\mathscr{J} \colon \Lambda^{\mathbb{R}} \! \to \! \Lambda^{\mathbb{R}}, \, \, 
j(z) \! \mapsto \! zj(z) \qquad \text{and} \qquad \mathscr{K} \colon \Lambda^{
\mathbb{R}} \! \to \! \Lambda^{\mathbb{R}}, \, \, k(z) \! \mapsto \! z^{-1}
k(z),
\end{equation*}
\begin{align*}
\setcounter{MaxMatrixCols}{18}
\mathcal{J} & = \!
\begin{pmatrix}
a_{0}^{\sharp} & b_{1}^{\sharp} & c_{2}^{\sharp} & & & & & & & & & & & & & \\
b_{1}^{\sharp} & a_{1}^{\sharp} & b_{2}^{\sharp} & & & & & & & & & & & & & \\
c_{2}^{\sharp} & b_{2}^{\sharp} & a_{2}^{\sharp} & b_{3}^{\sharp} & c_{4}^{
\sharp} & & & & & & & & & & & \\
 & & b_{3}^{\sharp} & a_{3}^{\sharp} & b_{4}^{\sharp} & & & & & & & & & & & \\
 & & c_{4}^{\sharp} & b_{4}^{\sharp} & a_{4}^{\sharp} & b_{5}^{\sharp} & c_{
6}^{\sharp} & & & & & & & & & \\
 & & & & b_{5}^{\sharp} & a_{5}^{\sharp} & b_{6}^{\sharp} & & & & & & & & & \\
 & & & & c_{6}^{\sharp} & b_{6}^{\sharp} & a_{6}^{\sharp} & b_{7}^{\sharp} &
c_{8}^{\sharp} & & & & & & & \\
 & & & & & & b_{7}^{\sharp} & a_{7}^{\sharp} & b_{8}^{\sharp} & & & & & & & \\
 & & & & & & c_{8}^{\sharp} & b_{8}^{\sharp} & a_{8}^{\sharp} & b_{9}^{\sharp}
& c_{10}^{\sharp} & & & & & \\
 & & & & & & & & \ddots & \ddots & \ddots & & & & & \\
 & & & & & & & & b_{2k+1}^{\sharp} & a_{2k+1}^{\sharp} & b_{2k+2}^{\sharp} & &
& & & \\
 & & & & & & & & c_{2k+2}^{\sharp} & b_{2k+2}^{\sharp} & a_{2k+2}^{\sharp} &
b_{2k+3}^{\sharp} & c_{2k+4}^{\sharp} & & & \\
 & & & & & & & & & & \ddots & \ddots & \ddots & & &
\end{pmatrix},
\end{align*}
and
\begin{align*}
\setcounter{MaxMatrixCols}{18}
\mathcal{K} &= \!
\begin{pmatrix}
\alpha_{0}^{\sharp} & \beta_{1}^{\sharp} & & & & & & & & & & & & & & \\
\beta_{1}^{\sharp} & \alpha_{1}^{\sharp} & \beta_{2}^{\sharp} & \gamma_{3}^{
\sharp}  & & & & & & & & & & & & \\
 & \beta_{2}^{\sharp} & \alpha_{2}^{\sharp} & \beta_{3}^{\sharp} & & & & & & &
& & & & & \\
 & \gamma_{3}^{\sharp} & \beta_{3}^{\sharp} & \alpha_{3}^{\sharp} & \beta_{4}^{
\sharp}  & \gamma_{5}^{\sharp} & & & & & & & & & & \\
 & & & \beta_{4}^{\sharp} & \alpha_{4}^{\sharp} & \beta_{5}^{\sharp} & & & & &
& & & & & \\
 & & & \gamma_{5}^{\sharp} & \beta_{5}^{\sharp} & \alpha_{5}^{\sharp} & \beta_{
6}^{\sharp} & \gamma_{7}^{\sharp} & & & & & & & & \\
 & & & & & \beta_{6}^{\sharp} & \alpha_{6}^{\sharp} & \beta_{7}^{\sharp} & & &
& & & & & \\
 & & & & & \gamma_{7}^{\sharp} & \beta_{7}^{\sharp} & \alpha_{7}^{\sharp} &
\beta_{8}^{\sharp} & \gamma_{9}^{\sharp} & & & & & & \\
 & & & & & & \beta_{8}^{\sharp} & \alpha_{8}^{\sharp} & \beta_{9}^{\sharp} & &
& & & & & \\
 & & & & & & \ddots & \ddots & \ddots & & & & & & & \\
 & & & & & & & & \gamma_{2k+1}^{\sharp} & \beta_{2k+1}^{\sharp} & \alpha_{2k+
1}^{\sharp} & \beta_{2k+2}^{\sharp} & \gamma_{2k+3}^{\sharp} & & & \\
 & & & & & & & & & & \beta_{2k+2}^{\sharp} & \alpha_{2k+2}^{\sharp} & \beta_{2
k+3}^{\sharp} & & & \\
 & & & & & & & & & & \ddots & \ddots & \ddots & & &
\end{pmatrix},
\end{align*}
with zeros outside the indicated diagonals; moreover, as shown in \cite{a23}, 
$\mathcal{J}$ and $\mathcal{K}$ are formal inverses, that is, $\mathcal{J} 
\mathcal{K} \! = \! \mathcal{K} \mathcal{J} \! = \! \operatorname{diag}(1,
\dotsc,1,\dotsc)$ (see, also, \cite{a27,a28,a29,a30,a31}).

It is convenient at this point to discuss, if only succinctly, a few of the 
multitudinous applications of $L$-polynomials (complete details may be found 
in the indicated references):
\begin{enumerate}
\item[(1)] as stated at the beginning of the Introduction, $L$-polynomials 
are intimately related with the solution of the SSMP and the SHMP. It is 
important to note \cite{a14} that the classical and strong moment problems 
(SMP, HMP, SSMP, and SHMP) are special cases of a more general theory, where 
moments corresponding to an arbitrary, countable sequence of (fixed) points 
are involved (in the classical and strong moment cases, respectively, the 
points are $\infty$ repeated and $0,\infty$ cyclically repeated), and where 
\emph{orthogonal rational functions} \cite{a26,a32,a33} play the r\^{o}le of 
orthogonal polynomials and orthogonal Laurent (or $L$-) polynomials; 
furthermore, since $L$-polynomials are rational functions with (fixed) poles 
at the origin and at the point at infinity, the step towards a more general 
theory where poles are at arbitrary, but fixed, positions/locations in 
$\mathbb{C} \cup \{\infty\}$ is natural, with applications to, say, 
multi-point Pad\'{e}, and Pad\'{e}-type, approximants 
\cite{a24,a34,a35,a36,a37,a38};
\item[(2)] in numerical analysis, the computation of integrals of the form 
$\int_{a}^{b}f(s) \, \md \mu (s)$, where $\mu$ is a positive measure on $[a,
b]$, and $-\infty \! \leqslant \! a \! < \! b \! \leqslant \! +\infty$, is an 
important problem. The most familiar quadrature formulae are the so-called 
Gauss-Christoffel formulae, that is, approximating the integral $\int_{a}^{b}
f(s) \, \md \mu (s)$ via a weighted-sum-of-products of function values of the 
form $\sum_{j=1}^{n} \mathscr{A}_{j,n}f(x_{j,n})$, $n \! \in \! \mathbb{N}$, 
where one chooses for the nodes $\lbrace x_{j,n} \rbrace_{j=1}^{n}$ the 
zeros/roots of $\varphi_{n}(z)$, the polynomial of degree $n$ orthogonal with 
respect to the inner product $\langle f,g \rangle \! = \! \int_{a}^{b}f(s) 
\overline{g(s)} \, \md \mu (s)$, and for the (positive) weights $\lbrace 
\mathscr{A}_{j,n} \rbrace_{j=1}^{n}$ the so-called Christoffel numbers 
\cite{a35}. When considering the computation of integrals of the form 
$\int_{-\pi}^{\pi}g(\me^{\mi \theta}) \, \md \mu (\theta)$, where $g$ is a 
complex-valued function on the unit circle $\mathbb{D} \! := \! \lbrace 
\mathstrut z \! \in \! \mathbb{C}; \, \lvert z \lvert = \! 1 \rbrace$ and 
$\mu$ is, say, a positive measure on $[-\pi,\pi]$, in 
p\-a\-r\-t\-i\-c\-u\-l\-ar, when $g$ is continuous on $\mathbb{D}$, 
keeping in mind that a function continuous on $\mathbb{D}$ can be uniformly 
approximated by $L$-polynomials, it is natural to consider, instead of 
orthogonal polynomials, Laurent polynomials, which are also related to 
the associated trigonometric moment problem \cite{a35,a39} (see, also, 
\cite{a40});
\item[(3)] for $V \colon \mathbb{R} \setminus \{0\} \! \to \! \mathbb{R}$ as 
described by conditions~(V1)--(V3), consider the function $g(z) \! = \! 
\int_{\mathbb{R}}(1 \! + \! sz)^{-1} \, \md \widetilde{\mu}(s)$, where $\md 
\widetilde{\mu}(s) \! = \! \exp (-\mathscr{N} \, V(s)) \, \md s$, $\mathscr{
N} \! \in \! \mathbb{N}$, which is holomorphic for $z \! \in \! \mathbb{C} 
\setminus \mathbb{R}$, with associated asymptotic expansions
\begin{gather*}
g(z) \underset{\mathbb{C} \setminus \mathbb{R} \ni z \to 0}{=} \, \sum_{m=
0}^{\infty}(-1)^{m}c_{m}z^{m} \! =: \! L^{0}(z) \qquad \text{and} \qquad g(z) 
\underset{\mathbb{C} \setminus \mathbb{R} \ni z \to \infty}{=} -\sum_{m=1}^{
\infty}(-1)^{m}c_{-m}z^{-m} \! =: \! L^{\infty}(z),
\end{gather*}
where $c_{l} \! = \! \int_{\mathbb{R}}s^{l} \me^{-\mathscr{N} \, V(s)} \, 
\md s$, $l \! \in \! \mathbb{Z}$, with respect to the (unbounded) domain $\{
\mathstrut z \! \in \! \mathbb{C}; \, \varepsilon \! \leqslant \! \vert 
\operatorname{Arg}(z) \vert \! \leqslant \! \pi \! - \! \varepsilon\}$, where 
$\operatorname{Arg}(\ast)$ denotes the principal argument of $\ast$, and 
$\varepsilon \! > \! 0$ is sufficiently small. Given the pair of formal power 
series $(L^{0}(z),L^{\infty}(z))$, the rational function $P_{k,n}(z)/Q_{k,n}
(z)$, where $P_{k,n}(z)$ belongs to the space of all polynomials of degree at 
most $n \! - \! 1$, and $Q_{k,n}(z)$ is a polynomial of degree exactly $n$ 
with $Q_{k,n}(0) \! \not= \! 0$, is said to be a \emph{$[k/n](z)$ two-point 
Pad\'{e} approximant} to $(L^{0}(z),L^{\infty}(z))$, $k \! \in \! \lbrace 0,
1,\dotsc,2n \rbrace$, if the following conditions are satisfied:
\begin{gather*}
L^{0}(z) \! - \! P_{k,n}(z)(Q_{k,n}(z))^{-1} \underset{z \to 0}{=} \mathcal{O}
(z^{k}), \\
L^{\infty}(z) \! - \! P_{k,n}(z)(Q_{k,n}(z))^{-1} \underset{z \to \infty}{=} 
\mathcal{O} \! \left((z^{-1})^{2n-k+1} \right).
\end{gather*}
The `balanced' situation corresponds to the case when $k \! = \! n$, in which 
case, the two-point Pad\'{e} approximants are denoted, simply, as $[n/n](z)$. 
An important, related problem of complex approximation theory is to study the 
convergence of sequences of two-point Pad\'{e} approximants constructed {}from 
the---formal---pair (of power series) $(L^{0}(z),L^{\infty}(z))$ to the 
function $g(z)$ on $\mathbb{C} \setminus \mathbb{R}$; in particular, denoting 
by $\mathscr{E}_{n}(z)$ the `error term' for the $[n/n](z)$ approximant, that 
is, $\mathscr{E}_{n}(z) \! := \! g(z) \! - \! [n/n](z)$, it can be shown that, 
following \cite{a41},
\begin{equation}
\mathscr{E}_{n}(z) \! = \! \left(\phi_{n}(-1/z) \right)^{-1} \int_{\mathbb{R}}
\dfrac{\phi_{n}(s) \me^{-\mathscr{N} \, V(s)}}{1 \! + \! sz} \, \md s, \quad
z \! \in \! \mathbb{C} \setminus \mathbb{R}, \tag{TPA1}
\end{equation}
where $\{\phi_{m}(z)\}_{m \in \mathbb{Z}_{0}^{+}}$ are the orthonormal 
$L$-polynomials defined in Equations~(1.2) and~(1.3). The main question 
regarding the convergence of two-point Pad\'{e} approximants for this class 
of functions is with which rate it takes place, that is, the so-called 
\emph{quantitative result} \cite{a42}: this necessitates obtaining results 
for the asymptotic behaviour (as $n \! \to \! \infty)$ of the orthonormal 
$L$-polynomials $\phi_{n}(z)$ in the entire complex plane. The theory of 
orthogonal $L$-pol\-y\-n\-o\-m\-i\-a\-l\-s is a natural framework for 
developing the theory of two-point Pad\'{e} approximants, for both the 
scalar and matrix cases \cite{a24,a41,a42,a43,a44};
\item[(4)] it turns out that, unlike the (finite) non-relativistic Toda 
lattice, whose direct and inverse spectral transform was constructed by Moser 
\cite{a45}, and which is based on the theory of orthogonal polynomials and 
tri-diagonal Jacobi matrices, the direct and inverse scattering transform for 
the (finite) relativistic Toda lattice, introduced by Ruijsenaars \cite{a46}, 
is based on the theory of orthogonal $L$-polynomials and pairs of bi-diagonal 
matrices \cite{a47} (see, also, \cite{a48}); and
\item[(5)] for a finite, countable or uncountable index set $K$, let $\lbrace 
\varsigma_{p}, \, p \! \in \! K \rbrace \subset \mathbb{C}_{+} \! := \! 
\lbrace \mathstrut z \! \in \! \mathbb{C}; \, \Im (z) \! > \! 0 \rbrace$, with 
$\varsigma_{p} \! \neq \! \varsigma_{q} \, \, \forall \, \, p \! \neq \! q \! 
\in \! K$, and $\lbrace \varpi_{p}, \, p \! \in \! K \rbrace \subset \mathbb{
C}$ be given point sets. A function $\mathfrak{F}(z)$ which is analytic for $z 
\! \in \! \mathbb{C}_{+}$, with $\Im (\mathfrak{F}(z)) \! \geqslant \! 0$, is 
called a \emph{Nevanlinna function}. The \emph{Pick-Nevanlinna problem} is: 
find a Nevanlinna function $\mathfrak{F}(z)$ so that $\mathfrak{F}(\varsigma_{
p}) \! = \! \varpi_{p} \, \, \forall \, \, p \! \in \! K$. A variant of this 
problem arises when, for $K \! = \! \mathbb{N}$, the points $\varsigma_{p}$, 
$p \! \in \! \mathbb{N}$, coalesce into the two points $0$ and $\infty$ (the 
point at infinity) according to the rule $\varsigma_{2i} \! = \! 0$, $i \! 
\in \! \mathbb{N}$, $\varsigma_{2j+1} \! = \! \infty$, $j \! \in \! \mathbb{ 
Z}_{0}^{+}$; then, the corresponding modification of the Pick-Nevanlinna 
problem is: given the bi-infinite sequence of numbers $\lbrace \breve{c}_{p} 
\rbrace_{p \in \mathbb{Z}}$, find a Nevanlinna function $\mathfrak{F}(z)$ with 
the asymptotic expansions $\mathfrak{F}(z) \! \sim_{z \to \infty} \! \sum_{k=
0}^{\infty} \breve{c}_{k}z^{-k}$ and $\mathfrak{F}(z) \! \sim_{z \to 0} \! 
\sum_{k=1}^{\infty} \breve{c}_{-k}z^{k}$ in every angular region $\lbrace 
\mathstrut z \! \in \! \mathbb{C}; \, \check{\delta} \! < \! \operatorname{Arg}
(z) \! < \! \pi \! - \! \check{\delta} \rbrace$, with $\check{\delta} \! > \! 
0$. This modified problem is equivalent to the SHMP \cite{a49}.
\end{enumerate}

Now that the principal objects have been defined, namely, the monic OLPs, $\{
\boldsymbol{\pi}_{m}(z)\}_{m \in \mathbb{Z}_{0}^{+}}$, and OLPs, $\{\phi_{m}
(z)\}_{m \in \mathbb{Z}_{0}^{+}}$, it is time to state what is actually 
studied in this work; in fact, this work constitutes the first part of a 
three-fold series of works devoted to asymptotics in the double-scaling limit 
as $\mathscr{N},n \! \to \! \infty$ such that $z_{o} \! := \! \mathscr{N}/n \! 
= \! 1 \! + \! o(1)$ (the simplified `notation' $n \! \to \! \infty$ will be 
adopted) of $L$-polynomials and related quantities. {}From the discussion 
above, an understanding of the large-$n$ (asymptotic) behaviour of the 
$L$-polynomials, as well as of the coefficients of the respective three- and 
five-term recurrence relations, is seminal in using the $L$-polynomials in 
several, seemingly disparate, applications: the purpose of the present series 
of works is, precisely, to analyse the $n \! \to \! \infty$ behaviour of the 
$L$-polynomials $\boldsymbol{\pi}_{n}(z)$ and $\phi_{n}(z)$ in $\mathbb{C}$, 
orthogonal with respect to the varying exponential measure\footnote{Note that 
$LD(\boldsymbol{\pi}_{m}) \! = \! LD(\phi_{m}) \! = \!
\begin{cases}
2n, &\text{$m \! = \! \text{even}$,} \\
2n \! + \! 1, &\text{$m \! = \! \text{odd}$,}
\end{cases}$ coincides with the parameter in the measure of orthogonality: 
the large parameter, $n$, enters simultaneously into the $L$-degree of the 
$L$-polynomials and the (varying exponential) weight; thus, asymptotics of 
the $L$-polynomials are studied along a `diagonal strip' of a doubly-indexed 
sequence.} $\md \mu (z) \! = \! \exp (-n \widetilde{V}(z)) \, \md z$, where 
$\widetilde{V}(z) \! := \! z_{o}V(z)$, and the (`scaled') external 
field\footnote{For real non-analytic external fields, see the recent work
\cite{a50}.} $\widetilde{V} \colon \mathbb{R} \setminus \{0\} \! \to \! 
\mathbb{R}$ satisfies conditions~(2.3)--(2.5) (see Subsection~2.2), as well 
as of the associated norming constants and coefficients of the (system of) 
recurrence relations; more precisely, then:
\begin{enumerate}
\item[\textbf{(i)}] in this work (Part~I), asymptotics (as $n \! \to \! 
\infty)$ of $\boldsymbol{\pi}_{2n}(z)$ (in the entire complex plane) and 
$\xi^{(2n)}_{n}$, thus $\phi_{2n}(z)$ (cf. Equation~(1.4)), and the Hankel 
determinant ratio $H^{(-2n)}_{2n}/H^{(-2n)}_{2n+1}$ (cf. Equations~(1.8)) 
are obtained;
\item[\textbf{(ii)}] in Part~II \cite{a51}, asymptotics (as $n \! \to \! 
\infty)$ of $\boldsymbol{\pi}_{2n+1}(z)$ (in the entire complex plane) and 
$\xi^{(2n+1)}_{-n-1}$, thus $\phi_{2n+1}(z)$ (cf. Equation~(1.5)), and the 
Hankel determinant ratio $H^{(-2n)}_{2n+1}/H^{(-2n-2)}_{2n+2}$ (cf. 
Equations~(1.8)) are obtained;
\item[\textbf{(iii)}] in Part~III \cite{a52}, asymptotics (as $n \! \to \! 
\infty)$ of $\nu^{(2n)}_{-n}$ $(= \! H^{(-2n+1)}_{2n}/H^{(-2n)}_{2n})$ and 
$\xi^{(2n)}_{-n}$, $\nu^{(2n+1)}_{n}$ $(= \! -H^{(-2n-1)}_{2n+1}/H^{(-2n)}_{2 
n+1})$ and $\xi^{(2n+1)}_{n}$, $\prod_{k=1}^{2n} \alpha^{(2n)}_{k}$ $(=\nu^{(
2n)}_{-n})$, and $\prod_{k=1}^{2n+1} \alpha^{(2n+1)}_{k}$ $(= \! -(\nu^{(2n+
1)}_{n})^{-1})$, as well as of the (elements of the) Laurent-Jacobi matrices, 
$\mathcal{J}$ and $\mathcal{K}$, and other, related, quantities constructed 
{}from the coefficients of the three- and five-term recurrence relations, 
are obtained.
\end{enumerate}

The above-mentioned asymptotics (as $n \! \to \! \infty)$ are obtained by 
reformulating, \emph{\`{a} la} Fokas-Its-Kitaev \cite{a53,a54}, the 
corresponding even degree and odd degree monic $L$-polynomial problems as 
(matrix) Riemann-Hilbert problems (RHPs) on $\mathbb{R}$, and then studying 
the large-$n$ behaviour of the corresponding solutions. The paradigm for the 
asymptotic (as $n \! \to \! \infty)$ analysis of the respective (matrix) 
RHPs is a union of the Deift-Zhou (DZ) non-linear steepest-descent method 
\cite{a1,a2}, used for the asymptotic analysis of undulatory RHPs, and the 
extension of Deift-Venakides-Zhou \cite{a3}, incorporating into the DZ method 
a non-linear analogue of the WKB method, making the asymptotic analysis of 
fully non-linear problems tractable (it should be mentioned that, in this 
context, the equilibrium measure \cite{a55} plays an absolutely crucial 
r\^{o}le in the analysis \cite{a56}); see, also, the multitudinous extensions 
and applications of the DZ method \cite{a57,a58,a59,a60,a61,a62,a63,a64,a65,%
a66,a67,a68,a69,a70,a71,a72,a73,a74,a75,a76,a77,a78,a79}. It is worth 
mentioning that asymptotics for Laurent-type polynomials and their zeros have 
been obtained in \cite{a42,a80} (see, also, \cite{a81,a82,a83}).

This article is organised as follows. In Section~2, necessary facts {}from 
the theory of compact Riemann surfaces are given, the respective `even 
degree' and `odd degree' RHPs on $\mathbb{R}$ are stated and the corresponding 
variational problems for the associated equilibrium measures are discussed, 
and the main results of this work, namely, asymptotics (as $n \! \to \! 
\infty)$ of $\boldsymbol{\pi}_{2n}(z)$ (in $\mathbb{C})$, and $\xi^{(2n)}_{n}$ 
and $\phi_{2n}(z)$ (in $\mathbb{C})$ are stated in Theorems~2.3.1 and~2.3.2, 
respectively. In Section~3, the detailed analysis of the `even degree' 
variational problem and the associated equilibrium measure is undertaken, 
including the construction of the so-called $g$-function, and the RHP 
formulated in Section~2 is reformulated as an equivalent, auxiliary RHP, 
which, in Sections~4 and~5, is augmented, by means of a sequence of contour 
deformations and transformations \emph{\`{a} la} Deift-Venakides-Zhou, into 
simpler, `model' (matrix) RHPs which, as $n \! \to \! \infty$, and in 
conjunction with the Beals-Coifman construction \cite{a84} (see, also, the 
extension of Zhou \cite{a85}) for the integral representation of the solution 
of a matrix RHP on an oriented contour, are solved explicitly (in closed 
form) in terms of Riemann theta functions (associated with the underlying 
finite-genus hyperelliptic Riemann surface) and Airy functions, {}from 
which the final asymptotic (as $n \! \to \! \infty)$ results stated in 
Theorems~2.3.1 and~2.3.2 are proved. The paper concludes with an Appendix.
\begin{eeeee}
The even degree OLPs, $\phi_{2n}(z)$, $n \! \in \! \mathbb{Z}_{0}^{+}$, are 
related, in a way, to the polynomials orthogonal with respect to the varying 
weight $\widehat{w}(z) \! := \! z^{-2n} \exp (-\mathscr{N} \, V(z))$, 
$\mathscr{N} \! \in \! \mathbb{N}$: this follows directly {}from the 
orthogonality relation satisfied by $\phi_{2n}(z)$. This does not help with 
any of the algebraic relations, such as the system of three- and five-term 
recurrence relations; however, this does provide for an alternative approach 
to computing large-$n$ asymptotics for $\phi_{2n}(z)$. The connection is not 
so clear for the odd degree OLPs, $\phi_{2n+1}(z)$, $n \! \in \! \mathbb{Z}_{
0}^{+}$. Indeed, in this latter case, the associated (density of the) measure 
for the orthogonal polynomials would take the form $\md \widehat{\mu}(z) \! 
:= \! z^{-2n-1} \exp (-\mathscr{N} \, V(z)) \, \md z$, and this measure 
changes signs, which causes a number of difficulties in the large-$n$ 
asymptotic analysis. In this paper, these connections are not used, and 
a complete asymptotic analysis of the even degree OLPs is carried out, 
directly. \hfill $\blacksquare$
\end{eeeee}
\section{Hyperelliptic Riemann Surfaces, The Riemann-Hilbert 
P\-r\-o\-b\-l\-e\-m\-s, a\-n\-d S\-u\-m\-m\-a\-r\-y of Results}
In this section, necessary facts {}from the theory of hyperelliptic Riemann 
surfaces are given (see Subsection~2.1), the respective RHPs on $\mathbb{R}$ 
for the even degree and odd degree monic orthogonal $L$-polynomials are 
formulated and the corresponding variational problems for the associated 
equilibrium measures are discussed (see Subsection~2.2), and asymptotics (as 
$n \! \to \! \infty)$ for $\boldsymbol{\pi}_{2n}(z)$ (in the entire complex 
plane), and $\xi^{(2n)}_{n}$ and $\phi_{2n}(z)$ (in the entire complex plane) 
are given in Theorems~2.3.1 and~2.3.2, respectively (see Subsection~2.3).

Before proceeding, however, the notation/nomenclature used throughout this 
work is summarised.
\begin{center}
\Ovalbox{\textsc{Notational Conventions}}
\end{center}
\begin{compactenum}
\item[(1)] $\mathrm{I} \! = \!
\left(
\begin{smallmatrix}
1 & 0 \\
0 & 1
\end{smallmatrix}
\right)$ is the $2 \times 2$ identity matrix, $\sigma_{1} \! = \!
\left(
\begin{smallmatrix}
0 & 1 \\
1 & 0
\end{smallmatrix}
\right)$, $\sigma_{2} \! = \!
\left(
\begin{smallmatrix}
0 & -\mi \\
\mi & 0
\end{smallmatrix}
\right)$ and $\sigma_{3} \! = \!
\left(
\begin{smallmatrix}
1 & 0 \\
0 & -1
\end{smallmatrix}
\right)$ are the Pauli matrices, $\sigma_{+} \! = \!
\left(
\begin{smallmatrix}
0 & 1 \\
0 & 0
\end{smallmatrix}
\right)$ and $\sigma_{-} \! = \!
\left(
\begin{smallmatrix}
0 & 0 \\
1 & 0
\end{smallmatrix}
\right)$ are, respectively, the raising and lowering matrices, $\mathbf{0} \! 
= \!
\left(
\begin{smallmatrix}
0 & 0 \\
0 & 0
\end{smallmatrix}
\right)$, $\mathbb{R}_{\pm} \! := \! \{\mathstrut x \! \in \! \mathbb{R}; \, 
\pm x \! > \! 0\}$, $\mathbb{C}_{\pm} \! := \! \{\mathstrut z \! \in \! 
\mathbb{C}; \, \pm \Im (z) \! > \! 0\}$, and $\mathrm{sgn}(x) \! := \! 0$ 
if $x \! = \! 0$ and $x \vert x \vert^{-1}$ if $x \! \not= \! 0;$
\item[(2)] for a scalar $\omega$ and a $2 \! \times \! 2$ matrix $\Upsilon$, 
$\omega^{\mathrm{ad}(\sigma_{3})} \Upsilon \! := \! \omega^{\sigma_{3}} 
\Upsilon \omega^{-\sigma_{3}}$;
\item[(3)] a contour $\mathcal{D}$ which is the finite union of 
piecewise-smooth, simple curves (as closed sets) is said to be 
\emph{orientable} if its complement $\mathbb{C} \setminus \mathcal{D}$ can 
always be divided into two, possibly disconnected, disjoint open sets $\mho^{
+}$ and $\mho^{-}$, either of which has finitely many components, such that
$\mathcal{D}$ admits an orientation so that it can either be viewed as a 
positively oriented boundary $\mathcal{D}^{+}$ for $\mho^{+}$ or as a 
negatively oriented boundary $\mathcal{D}^{-}$ for $\mho^{-}$ \cite{a85}, that 
is, the (possibly disconnected) components of $\mathbb{C} \setminus \mathcal{
D}$ can be coloured by $+$ or $-$ in such a way that the $+$ regions do not 
share boundary with the $-$ regions, except, possibly, at finitely many points
\cite{a86};
\item[(4)] for each segment of an oriented contour $\mathcal{D}$, according to 
the given orientation, the ``+'' side is to the left and the ``-'' side is to 
the right as one traverses the contour in the direction of orientation, that 
is, for a matrix $\mathcal{A}_{ij}(z)$, $i,j \! = \! 1,2$, $(\mathcal{A}_{ij}
(z))_{\pm}$ denote the non-tangential limits $(\mathcal{A}_{ij}(z))_{\pm} \! 
:= \! \lim_{\genfrac{}{}{0pt}{2}{z^{\prime} \, \to \, z}{z^{\prime} \, \in \, 
\pm \, \mathrm{side} \, \mathrm{of} \, \mathcal{D}}} \mathcal{A}_{ij}(z^{
\prime})$;
\item[(5)] for $1 \! \leqslant \! p \! < \! \infty$ and $\mathcal{D}$ some 
point set,
\begin{equation*}
\mathcal{L}^{p}_{\mathrm{M}_{2}(\mathbb{C})}(\mathcal{D}) \! := \! \left\{
\mathstrut f \colon \mathcal{D} \! \to \! \mathrm{M}_{2}(\mathbb{C}); \, 
\vert \vert f(\cdot) \vert \vert_{\mathcal{L}^{p}_{\mathrm{M}_{2}(\mathbb{C})
}(\mathcal{D})} \! := \! \left(\int_{\mathcal{D}} \vert f(z) \vert^{p} \, 
\vert \md z \vert \right)^{1/p} \! < \! \infty \right\},
\end{equation*}
where, for $\mathcal{A}(z) \! \in \! \operatorname{M}_{2}(\mathbb{C})$, $\vert 
\mathcal{A}(z) \vert \! := \! (\sum_{i,j=1}^{2} \overline{\mathcal{A}_{ij}(z)} 
\, \mathcal{A}_{ij}(z))^{1/2}$ is the Hilbert-Schmidt norm, with $\overline{
\bullet}$ denoting complex conjugation of $\bullet$, for $p \! = \! \infty$,
\begin{equation*}
\mathcal{L}^{\infty}_{\mathrm{M}_{2}(\mathbb{C})}(\mathcal{D}) \! := \!
\left\{\mathstrut g \colon \mathcal{D} \! \to \! \mathrm{M}_{2}(\mathbb{C});
\, \vert \vert g(\cdot) \vert \vert_{\mathcal{L}^{\infty}_{\mathrm{M}_{2}
(\mathbb{C})}(\mathcal{D})} \! := \! \max_{i,j = 1,2} \, \, \sup_{z \in
\mathcal{D}} \vert g_{ij}(z) \vert \! < \! \infty \right\},
\end{equation*}
and, for $f \! \in \! \mathrm{I} \! + \! \mathcal{L}^{2}_{\mathrm{M}_{2}
(\mathbb{C})}(\mathcal{D}) \! := \! \left\{\mathstrut \mathrm{I} \! + \! h;
\, h \! \in \! \mathcal{L}^{2}_{\mathrm{M}_{2}(\mathbb{C})}(\mathcal{D})
\right\}$,
\begin{equation*}
\vert \vert f(\cdot) \vert \vert_{\mathrm{I}+\mathcal{L}^{2}_{\mathrm{M}_{2}
(\mathbb{C})}(\mathcal{D})} \! := \! \left(\vert \vert f(\infty) \vert \vert_{
\mathcal{L}^{\infty}_{\mathrm{M}_{2}(\mathbb{C})}(\mathcal{D})}^{2} \! + \!
\vert \vert f(\cdot) \! - \! f(\infty) \vert \vert_{\mathcal{L}^{2}_{\mathrm{
M}_{2}(\mathbb{C})}(\mathcal{D})}^{2} \right)^{1/2};
\end{equation*}
\item[(6)] for a matrix $\mathcal{A}_{ij}(z)$, $i,j \! = \! 1,2$, to have 
boundary values in the $\mathcal{L}^{2}_{\mathrm{M}_{2}(\mathbb{C})}(\mathcal{
D})$ sense on an oriented contour $\mathcal{D}$, it is meant that $\lim_{
\genfrac{}{}{0pt}{2}{z^{\prime} \, \to \, z}{z^{\prime} \, \in \, \pm \,
\mathrm{side} \, \mathrm{of} \, \mathcal{D}}} \int_{\mathcal{D}} \vert 
\mathcal{A}(z^{\prime}) \! - \! (\mathcal{A}(z))_{\pm} \vert^{2} \, \vert 
\md z \vert \! = \! 0$ (e.g., if $\mathcal{D} \! = \! \mathbb{R}$ is oriented
{}from $+\infty$ to $-\infty$, then $\mathcal{A}(z)$ has $\mathcal{L}^{2}_{
\mathrm{M}_{2}(\mathbb{C})}(\mathcal{D})$ boundary values on $\mathcal{D}$
means that $\lim_{\varepsilon \downarrow 0} \int_{\mathbb{R}} \vert \mathcal{A}
(x \! \mp \! \mi \varepsilon) \! - \! (\mathcal{A}(x))_{\pm} \vert^{2} \, \md
x \! = \! 0)$;
\item[(7)] for a $2 \times 2$ matrix-valued function $\mathfrak{Y}(z)$, the 
notation $\mathfrak{Y}(z) \! =_{z \to z_{0}} \! \mathcal{O}(\ast)$ means 
$\mathfrak{Y}_{ij}(z) \! =_{z \to z_{0}} \! \mathcal{O}(\ast_{ij})$, $i,j \! 
= \! 1,2$ (\emph{mutatis mutandis} for $o(1))$;
\item[(8)] $\vert \vert \mathscr{F}(\cdot) \vert \vert_{\cap_{p \in J} 
\mathcal{L}^{p}_{\mathrm{M}_{2}(\mathbb{C})}(\ast)} \! := \! \sum_{p \in J} 
\vert \vert \mathscr{F}(\cdot) \vert \vert_{\mathcal{L}^{p}_{\mathrm{M}_{2}
(\mathbb{C})}(\ast)}$, with $\mathrm{card}(J) \! < \! \infty$;
\item[(9)] $\mathcal{M}_{1}(\mathbb{R})$ denotes the set of all non-negative,
bounded, unit Borel measures on $\mathbb{R}$ for which all moments exist,
\begin{equation*}
\mathcal{M}_{1}(\mathbb{R}) \! := \! \left\{\mathstrut \mu; \, \int_{\mathbb{
R}} \md \mu (s) \! = \! 1, \, \int_{\mathbb{R}}s^{m} \, \md \mu (s) \! < \!
\infty, \, m \! \in \! \mathbb{Z} \setminus \{0\} \right\};
\end{equation*}
\item[(10)] for $(\mu,\nu) \! \in \! \mathbb{R} \times \mathbb{R}$, denote the 
function $(\bullet \! - \! \mu)^{\mi \nu} \colon \mathbb{C} \setminus (-\infty,
\mu) \! \to \! \mathbb{C}$, $\bullet \! \mapsto \! \exp (\mi \nu \ln (\bullet 
-\mu))$, where $\ln$ denotes the principal branch of the logarithm;
\item[(11)] for $\widetilde{\boldsymbol{\gamma}}$ a null-homologous path in a 
region $\mathscr{D} \subset \mathbb{C}$, $\mathrm{int}(\widetilde{\boldsymbol{
\gamma}}) \! := \! \left\{\mathstrut \zeta \! \in \! \mathscr{D} \setminus 
\widetilde{\boldsymbol{\gamma}}; \, \mathrm{ind}_{\widetilde{\boldsymbol{
\gamma}}}(\zeta) \! := \! \int_{\widetilde{\boldsymbol{\gamma}}} \tfrac{1}{z
-\zeta} \, \tfrac{\md z}{2 \pi \mi} \! \not=\! 0 \right\}$;
\item[(12)] for some point set $\mathcal{D} \subset \mathcal{X}$, with 
$\mathcal{X} \! = \! \mathbb{C}$ or $\mathbb{R}$, $\overline{\mathcal{D}} 
:= \! \mathcal{D} \cup \partial \mathcal{D}$, and $\mathcal{D}^{c} \! := \! 
\mathcal{X} \setminus \overline{\mathcal{D}}$.
\end{compactenum}
\subsection{Riemann Surfaces: Preliminaries}
In this subsection, the basic elements associated with the construction of 
hyperelliptic and finite genus (compact) Riemann surfaces are presented 
(for further details and proofs, see, for example, \cite{a87,a88}).

\begin{eeee}
The superscripts ${}^{\pm}$, and sometimes subscripts ${}_{\pm}$, in this 
subsection should not be confused with the subscripts ${}_{\pm}$ appearing 
in the various RHPs (this is a general comment which applies, unless stated 
otherwise, throughout the entire text). Although $\overline{\mathbb{C}}$ 
(or $\mathbb{C} \mathbb{P}^{1})$ $:= \! \mathbb{C} \cup \{\infty\}$ (resp., 
$\overline{\mathbb{R}} \! := \! \mathbb{R} \cup \{-\infty\} \cup \{+\infty\})$
is the standard definition for the (closed) Riemann sphere (resp., closed 
real line), the simplified, and somewhat abusive, notation $\mathbb{C}$ 
(resp., $\mathbb{R})$ is used to denote both the (closed) Riemann sphere, 
$\overline{\mathbb{C}}$ (resp., closed real line, $\overline{\mathbb{R}})$, 
and the (open) complex field, $\mathbb{C}$ (resp., open real line, $\mathbb{
R})$, and the context(s) should make clear which object(s) the notation 
$\mathbb{C}$ (resp., $\mathbb{R})$ represents. \hfill $\blacksquare$
\end{eeee}

Let $N \! \in \! \mathbb{N}$ (with $N \! < \! \infty$ assumed throughout) and 
$\varsigma_{k} \! \in \! \mathbb{R} \setminus \{0,\pm \infty\}$, $k \! = \! 1,
\dotsc,2N \! + \! 2$, be such that $\varsigma_{i} \! \not= \! \varsigma_{j} 
\, \, \forall \, \, i \! \not= \! j \! = \! 1,\dotsc,2N \! + \! 2$, and 
enumerated/ordered according to $\varsigma_{1} \! < \! \varsigma_{2} \! < \! 
\cdots \! < \! \varsigma_{2N+2}$. Let $R(z) \! := \! \prod_{j=1}^{N}(z \! - 
\! \varsigma_{2j-1})(z \! - \! \varsigma_{2j})$ $\in \! \mathbb{R}[z]$ (the 
algebra of polynomials in $z$ with coefficients in $\mathbb{R})$ be the 
(unital) polynomial of even degree $\mathrm{deg}(R) \! = \! 2N \! + \! 2$ 
$(\mathrm{deg}(R) \! = \! 0$ $(\mathrm{mod} 2))$ whose (simple) zeros/roots 
are $\{\varsigma_{j}\}_{j=1}^{2N+2}$. Denote by $\mathscr{R}$ the 
hyperelliptic Riemann surface of genus $N$ defined by the equation $y^{2} \! 
= \! R(z)$ and realised as a two-sheeted branched (ramified) covering of the 
Riemann sphere such that its two sheets are two identical copies of $\mathbb{
C}$ with branch cuts along the intervals $(\varsigma_{1},\varsigma_{2})$, 
$(\varsigma_{3},\varsigma_{4})$, $\dotsc$, $(\varsigma_{2N+1},\varsigma_{2N+
2})$, and glued/pasted to each other `crosswise' along the opposite banks of 
the corresponding cuts $(\varsigma_{2j-1},\varsigma_{2j})$, $j \! = \! 1,
\dotsc,N \! + \! 1$. Denote the two sheets of $\mathscr{R}$ by $\mathscr{R}^{
+}$ (the first/upper sheet) and $\mathscr{R}^{-}$ (the second/lower sheet): to 
indicate that $z$ lies on the first (resp., second) sheet, one writes $z^{+}$ 
(resp., $z^{-})$; of course, as points in the plane $\mathbb{C}$, $z^{+} \! 
= \! z^{-} \! = \! z$. For points $z$ on the first (resp., second) sheet 
$\mathscr{R}^{+}$ (resp., $\mathscr{R}^{-})$, one has $z^{+} \! = \! (z,
+(R(z))^{1/2})$ (resp., $z^{-} \! = \! (z,-(R(z))^{1/2}))$, where the 
single-valued branch of the square root is chosen such that $z^{-(N+1)}
(R(z))^{1/2} \! \sim_{\underset{z \in \mathscr{R}^{\pm}}{z \to \infty}} \! 
\pm 1$.

Let $\mathscr{E}_{j} \! := \! (\varsigma_{2j-1},\varsigma_{2j})$, $j \! = \! 
1,\dotsc,N \! + \! 1$, and set $\mathscr{E} \! = \! \cup_{j=1}^{N+1} \mathscr{
E}_{j}$ (note that $\mathscr{E}_{i} \cap \mathscr{E}_{j} \! = \! \varnothing$, 
$i \! \not= \! j \! = \! 1,\dotsc,N \! + \! 1)$. Denote by $\mathscr{E}_{j}^{
+}$ $(\subset \mathscr{R}^{+})$ (resp., $\mathscr{E}_{j}^{-}$ $(\subset 
\mathscr{R}^{-}))$ the upper (resp., lower) bank of the interval $\mathscr{
E}_{j}$, $j \! = \! 1,\dotsc,N \! + \! 1$, forming $\mathscr{E}$, and oriented 
in accordance with the orientation of $\mathscr{E}$ as the boundary of 
$\mathbb{C} \setminus \mathscr{E}$, namely, the domain $\mathbb{C} \setminus 
\mathscr{E}$ is on the left as one proceeds along the upper bank of the $j$th 
interval {}from $\varsigma_{2j-1}$ to the point $\varsigma_{2j}$ and back 
along the lower bank from $\varsigma_{2j}$ to $\varsigma_{2j-1}$; thus, 
$\mathscr{E}_{j}^{\pm} \! := \! (\varsigma_{2j-1},\varsigma_{2j})^{\pm}$, $j 
\! = \! 1,\dotsc,N \! + \! 1$, are two (identical) copies of $(\varsigma_{2j
-1},\varsigma_{2j}) \subset \mathbb{R}$ `lifted' to $\mathscr{R}^{\pm}$. Set 
$\Gamma \! := \! \cup_{j=1}^{N+1} \Gamma_{j}$ $(\subset \mathscr{R})$, where 
$\Gamma_{j} \! := \! \mathscr{E}_{j}^{+} \cup \mathscr{E}_{j}^{-}$, $j \! = 
\! 1,\dotsc,N \! + \! 1$ $(\Gamma \! = \! \mathscr{E}^{+} \cup \mathscr{E}^{
-})$: note that $\Gamma$, as a curve on $\mathscr{R}$ (defined by the equation 
$y^{2} \! = \! R(z))$, consists of a finitely denumerable number of disjoint 
analytic closed Jordan curves, $\Gamma_{j}$, $j \! = \! 1,\dotsc,N \! + \! 1$, 
which are \emph{cycles} on $\mathscr{R}$, and that correspond to the intervals 
$\mathscr{E}_{j}$. {}From the above construction, it is clear that $\mathscr{
R} \! = \! \mathscr{R}^{+} \cup \mathscr{R}^{-} \cup \Gamma$; furthermore, the 
canonical projection of $\Gamma$ onto $\mathbb{C}$ $(\boldsymbol{\pi} \colon 
\mathscr{R} \! \to \! \mathbb{C})$ is $\mathscr{E}$, that is, $\boldsymbol{
\pi}(\Gamma) \! = \! \mathscr{E}$ (also, $\boldsymbol{\pi}(\mathscr{R}^{+}) 
\! = \! \boldsymbol{\pi}(\mathscr{R}^{-}) \! = \! \mathbb{C} \setminus 
\mathscr{E}$, or, alternately, $\boldsymbol{\pi}(z^{+}) \! = \! \boldsymbol{
\pi}(z^{-}) \! = \! z)$. One moves in the `positive $(+)$' (resp., `negative 
$(-)$') direction along the (closed) contour $\Gamma \subset \mathscr{R}$ if 
the domain $\mathscr{R}^{+}$ is on the left (resp., right) and the domain 
$\mathscr{R}^{-}$ is on the right (resp., left): the corresponding notation 
is (see above) $\Gamma^{+}$ (resp., $\Gamma^{-})$. For a function $f$ defined 
on the two-sheeted hyperelliptic Riemann surface $\mathscr{R}$, one defines 
the non-tangential boundary values, provided they exist, of $f(z)$ as $z \! 
\in \! \mathscr{R}^{+}$ (resp., $z \! \in \! \mathscr{R}^{-})$ approaches 
$\lambda \! \in \! \Gamma$, denoted $\lambda_{+}$ (resp., $\lambda_{-})$, by 
$f(\lambda_{\pm}) \! := \! f_{\pm}(\lambda) \! := \! \lim_{\underset{z \in 
\Gamma^{\pm}}{z \to \lambda}}f(z)$.

One takes the first $N$ contours among the (closed) contours $\Gamma_{j}$ for 
basis $\boldsymbol{\alpha}$-cycles $\{\boldsymbol{\alpha}_{j}, \, j \! = \! 1,
\dotsc,N\}$ and then completes/supplements this in the standard way with 
$\boldsymbol{\beta}$-cycles $\{\boldsymbol{\beta}_{j}, \, j \! = \! 1,\dotsc,
N\}$ so that the \emph{intersection matrix} has the (canonical) form 
$\boldsymbol{\alpha}_{k} \circ \boldsymbol{\alpha}_{j} \! = \! 
\boldsymbol{\beta}_{k} \circ \boldsymbol{\beta}_{j} \! = \! 0 \, \, \forall \, 
\, k \! \not= \! j \! = \! 1,\dotsc,N$, and $\boldsymbol{\alpha}_{k} \circ 
\boldsymbol{\beta}_{j} \! = \! \delta_{kj}$: the cycles $\{\boldsymbol{
\alpha}_{j},\boldsymbol{\beta}_{j}\}$, $j \! = \! 1,\dotsc,N$, form the 
\emph{canonical 1-homology basis} on $\mathscr{R}$, namely, any cycle 
$\widehat{\boldsymbol{\gamma}} \subset \mathscr{R}$ is homologous to an 
integral linear combination of $\{\boldsymbol{\alpha}_{j},\boldsymbol{\beta}_{
j}\}$, that is, $\widehat{\boldsymbol{\gamma}} \! = \! \sum_{j=1}^{N}(n_{j} 
\boldsymbol{\alpha}_{j} \! + \! m_{j} \boldsymbol{\beta}_{j})$, where $(n_{j},
m_{j}) \! \in \! \mathbb{Z} \times \mathbb{Z}$, $j \! = \! 1,\dotsc,N$. The 
$\boldsymbol{\alpha}$-cycles $\{\boldsymbol{\alpha}_{j}, \, j \! = \! 1,
\dotsc,N\}$, in the present case, are the intervals $(\varsigma_{2j-1},
\varsigma_{2j})$, $j \! = \! 1,\dotsc,N$, `going twice', that is, along the 
upper ({}from $\varsigma_{2j-1}$ to $\varsigma_{2j})$ and lower ({}from 
$\varsigma_{2j}$ to $\varsigma_{2j-1})$ banks $(\alpha_{j} \! = \! \mathscr{
E}_{j}^{+} \cup \mathscr{E}_{j}^{-}$, $j \! = \! 1,\dotsc,N)$, and the 
$\boldsymbol{\beta}$-cycles $\{\boldsymbol{\beta}_{j}, \, j \! = \! 1,\dotsc,
N\}$ are as follows: the $j$th $\boldsymbol{\beta}$-cycle consists of the 
$\boldsymbol{\alpha}$-cycles $\boldsymbol{\alpha}_{k}$, $k \! = \! j \! + \! 
1,\dotsc,N$, and the cycles `linked' with them and consisting of (the gaps) 
$(\varsigma_{2k},\varsigma_{2k+1})$, $k \! = \! 1,\dotsc,N$, `going twice', 
that is, {}from $\varsigma_{2k}$ to $\varsigma_{2k+1}$ on the first sheet and 
in the reverse direction on the second sheet. For an arbitrary holomorphic 
Abelian differential (one-form) $\boldsymbol{\omega}$ on $\mathscr{R}$, the 
function $\int^{z} \boldsymbol{\omega}$ is defined uniquely modulo its 
$\boldsymbol{\alpha}$- and $\boldsymbol{\beta}$-periods, 
$\oint_{\boldsymbol{\alpha}_{j}} \boldsymbol{\omega}$ and 
$\oint_{\boldsymbol{\beta}_{j}} \boldsymbol{\omega}$, $j \! = \! 1,\dotsc,N$, 
respectively. It is well known that the canonical $1$-homology basis 
$\{\boldsymbol{\alpha}_{j},\boldsymbol{\beta}_{j}\}$, $j \! = \! 1,\dotsc,N$, 
constructed above `generates', on $\mathscr{R}$, the corresponding 
$\boldsymbol{\alpha}$-normalised basis of holomorphic Abelian differentials 
(one-forms) $\{\omega_{1},\omega_{2},\dotsc,\omega_{N}\}$, where $\omega_{j} 
\! := \! \sum_{k=1}^{N} \tfrac{c_{jk}z^{N-k}}{\sqrt{\smash[b]{R(z)}}} \, 
\md z$, $c_{jk} \! \in \! \mathbb{C}$, $j \! = \! 1,\dotsc,N$, and 
$\oint_{\boldsymbol{\alpha}_{k}} \omega_{j} \! = \! \delta_{kj}$, 
$k, \, j \! = \! 1,\dotsc,N$: the associated $N \times N$ matrix of 
$\boldsymbol{\beta}$-periods, $\tau \! = \! (\tau_{ij})_{i,j=1,\dotsc,N} \! 
:= \! \left(\oint_{\boldsymbol{\beta}_{j}} \omega_{i} \right)_{i,j=1,\dotsc,
N}$, is a \emph{Riemann matrix}, that is, it is symmetric $(\tau_{ij} \! = 
\! \tau_{ji})$, pure imaginary, and $-\mi \tau$ is positive definite $(\Im 
(\tau_{jj}) \! > \! 0)$; moreover, $\tau$ is non-degenerate $(\det (\tau) 
\! \not= \! 0)$. {}From the condition that the basis of the differentials 
$\omega_{l}$, $l \! = \! 1,\dotsc,N$, is canonical, with respect to the given 
basis cycles $\{\boldsymbol{\alpha}_{j},\boldsymbol{\beta}_{j}\}$, it is seen 
that this implies that each $\omega_{l}$ is real valued on $\mathscr{E} \! = 
\! \cup_{j=1}^{N+1}(\varsigma_{2j-1},\varsigma_{2j})$ and has exactly one 
(real) root/zero in any interval (band) $(\varsigma_{2j-1},\varsigma_{2j})$, 
$j \! = \! 1,\dotsc,N \! + \! 1$, $j \! \not= \! l$; moreover, in the `gaps' 
$(\varsigma_{2j},\varsigma_{2j+1})$, $j \! = \! 1,\dotsc,N$, these 
differentials take non-zero, pure imaginary values.

Fix the `standard basis' $\boldsymbol{e}_{1},\boldsymbol{e}_{2},\dotsc,
\boldsymbol{e}_{N}$ in $\mathbb{R}^{N}$, that is, $(\boldsymbol{e}_{j})_{k} 
\! = \! \delta_{jk}$, $j,k \! = \! 1,\dotsc,N$ (these standard basis vectors 
should be viewed as column vectors): the vectors $\boldsymbol{e}_{1},
\boldsymbol{e}_{2},\dotsc,\boldsymbol{e}_{N},\tau \boldsymbol{e}_{1},\tau 
\boldsymbol{e}_{2},\dotsc,\tau \boldsymbol{e}_{N}$ are linearly independent 
over $\mathbb{R}$, and form a `basis' in $\mathbb{C}^{N}$. The quotient space 
$\mathbb{C}^{N}/\{N \! + \! \tau M\}$, $(N,M) \! \in \! \mathbb{Z}^{N} \times 
\mathbb{Z}^{N}$, where $\mathbb{Z}^{N} \! := \! \{\mathstrut (m_{1},m_{2},
\dotsc,m_{N}); \, m_{j} \! \in \! \mathbb{Z}, \, j \! = \! 1,\dotsc,N\}$, is 
a $2N$-dimensional real torus $\mathbb{T}^{2N}$, and is referred to as the 
\emph{Jacobi variety}, symbolically $\mathrm{Jac}(\mathscr{R})$, of the 
two-sheeted (hyperelliptic) Riemann surface $\mathscr{R}$ of genus $N$. 
Let $z_{0}$ be a fixed point in $\mathscr{R}$. A vector-valued function 
$\boldsymbol{\mathscr{A}}(z) \! = \! (\mathscr{A}_{1}(z),\mathscr{A}_{2}
(z),\dotsc,\mathscr{A}_{N}(z)) \! \in \! \mathrm{Jac}(\mathscr{R})$ with 
co-ordinates $\mathscr{A}_{k}(z) \! \equiv \! \int_{z_{0}}^{z} \omega_{k}$, 
$k \! = \! 1,\dotsc,N$, where, hereafter, unless stated otherwise and/or 
where confusion may arise, $\equiv$ denotes `congruence modulo the period 
lattice', defines the \emph{Abel map} $\boldsymbol{\mathscr{A}} \colon 
\mathscr{R} \! \to \! \mathrm{Jac}(\mathscr{R})$. The unordered set of points 
$z_{1},z_{2},\dotsc,z_{N}$, with $z_{k} \! \in \! \mathscr{R}$, form the $N$th 
symmetric power of $\mathscr{R}$, symbolically $\mathscr{R}^{N}_{\mathrm{symm}
}$ (or $\mathscr{S}^{N} \mathscr{R})$. The vector function $\boldsymbol{
\mathfrak{U}} \! = \! (\mathfrak{U}_{1},\mathfrak{U}_{2},\dotsc,\mathfrak{
U}_{N})$ with co-ordinates $\mathfrak{U}_{j} \! = \! \mathfrak{U}_{j}(z_{1},
z_{2},\dotsc,z_{N}) \! \equiv \! \sum_{k=1}^{N} \mathscr{A}_{j}(z_{k}) \! 
\equiv \! \sum_{k=1}^{N} \int_{z_{0}}^{z_{k}} \omega_{j}$, $j \! = \! 1, 
\dotsc,N$, that is, $(z_{1},z_{2},\dotsc,z_{N}) \! \to \! (\sum_{k=1}^{N} 
\int_{z_{0}}^{z_{k}} \omega_{1},\sum_{k=1}^{N} \int_{z_{0}}^{z_{k}} \omega_{
2},\dotsc,\sum_{k=1}^{N} \int_{z_{0}}^{z_{k}} \omega_{N})$, is also referred 
to as the \emph{Abel map}, $\boldsymbol{\mathfrak{U}} \colon \mathscr{R}^{
N}_{\mathrm{symm}} \! \to \! \mathrm{Jac}(\mathscr{R})$ (or $\boldsymbol{
\mathfrak{U}} \colon \mathscr{S}^{N} \mathscr{R} \! \to \! \mathrm{Jac}
(\mathscr{R}))$. It is known (see, for example, \cite{a89}) that the Abel map 
$\boldsymbol{\mathfrak{U}} \colon \mathscr{R}^{N}_{\mathrm{symm}} \! \to \! 
\mathrm{Jac}(\mathscr{R})$ is surjective and locally biholomorphic, but not 
injective globally. The \emph{dissected} Riemann surface, symbolically 
$\widetilde{\mathscr{R}}$, is obtained {}from $\mathscr{R}$ by `cutting' 
(canonical dissection) along the cycles of the canonical $1$-homology basis 
$\boldsymbol{\alpha}_{k},\boldsymbol{\beta}_{k}$, $k \! = \! 1,\dotsc,N$, of 
the original surface, namely, $\widetilde{\mathscr{R}} \! = \! \mathscr{R} 
\setminus (\cup_{j=1}^{N}(\boldsymbol{\alpha}_{j} \cup \boldsymbol{\beta}_{
j}))$; the surface $\widetilde{\mathscr{R}}$ is not only connected, as one 
can `pass' {}from one sheet to the other `across' $\Gamma_{N+1}$, but also 
simply connected (a $4N$-sided polygon ($4N$-gon) of a canonical dissection 
of $\mathscr{R}$ associated with the given canonical $1$-homology basis for 
$\mathscr{R})$. For a given vector $\vec{\boldsymbol{v}} \! = \! (\upsilon_{
1},\upsilon_{2},\dotsc,\upsilon_{N}) \! \in \! \mathrm{Jac}(\mathscr{R})$, 
the problem of finding an unordered collection of points $z_{1},z_{2},\dotsc,
z_{N}$, $z_{j} \! \in \! \mathscr{R}$, $j \! = \! 1,\dotsc,N$, for which 
$\mathfrak{U}_{k}(z_{1},z_{2},\dotsc,z_{N}) \! \equiv \! \upsilon_{k}$, $k \! 
= \! 1,\dotsc,N$, is called the \emph{Jacobi inversion problem} for Abelian 
integrals: as is well known, the Jacobi inversion problem is always solvable; 
but not, in general, uniquely.

By a \emph{divisor} on the Riemann surface $\mathscr{R}$ is meant a formal 
`symbol' $\boldsymbol{d} \! = \! z_{1}^{n_{f}(z_{1})}z_{2}^{n_{f}(z_{2})} 
\cdots z_{m}^{n_{f}(z_{m})}$, where $z_{j} \! \in \! \mathscr{R}$ and $n_{f}
(z_{j}) \! \in \! \mathbb{Z}$, $j \! = \! 1,\dotsc,m$: the number $\vert 
\boldsymbol{d} \vert \! := \! \sum_{j=1}^{m}n_{f}(z_{j})$ is called the 
\emph{degree} of the divisor $\boldsymbol{d}$: if $z_{i} \! \not= \! z_{j} \, 
\, \forall \, \, i \! \not= \! j \! = \! 1,\dotsc,m$, and if $n_{f}(z_{j}) \! 
\geqslant \! 0$, $j \! = \! 1,\dotsc,m$, then the divisor $\boldsymbol{d}$ 
is said to be \emph{integral}. Let $g$ be a meromorphic function defined on 
$\mathscr{R}$: for an arbitrary point $a \! \in \! \mathscr{R}$, one denotes 
by $n_{g}(a)$ (resp., $p_{g}(a))$ the multiplicity of the zero (resp., pole) 
of the function $g$ at this point if $a$ is a zero (resp., pole), and sets 
$n_{g}(a) \! = \! 0$ (resp., $p_{g}(a) \! = \! 0)$ otherwise; thus, $n_{g}(a),
\, p_{g}(a) \! \geqslant \! 0$. To a meromorphic function $g$ on $\mathscr{
R}$, one assigns the divisor $(g)$ of zeros and poles of this function as 
$(g) \! = \! z_{1}^{n_{g}(z_{1})}z_{2}^{n_{g}(z_{2})} \cdots z_{l_{1}}^{n_{g}
(z_{l_{1}})} \lambda_{1}^{-p_{g}(\lambda_{1})} \lambda_{2}^{-p_{g}(\lambda_{
2})} \cdots \lambda_{l_{2}}^{-p_{g}(\lambda_{l_{2}})}$, where $z_{i}, \, 
\lambda_{j} \! \in \! \mathscr{R}$, $i \! = \! 1,\dotsc,l_{1}$, $j \! = \! 
1,\dotsc,l_{2}$, are the zeros and poles of $g$ on $\mathscr{R}$, and $n_{g}
(z_{i}), \, p_{g}(\lambda_{j}) \! \geqslant \! 0$ are their multiplicities 
(one can also write $\{\mathstrut (a,n_{g}(a),-p_{g}(a)); \, a \! \in \! 
\mathscr{R}, \, n_{g}(a), \, p_{g}(a) \! \geqslant \! 0\}$ for the divisor 
$(g)$ of $g)$: these divisors are said to be \emph{principal}.

Associated with the Riemann matrix of $\boldsymbol{\beta}$-periods, $\tau$, 
is the \emph{Riemann theta function}, defined by
\begin{equation*}
\boldsymbol{\theta}(z;\tau) \! =: \! \boldsymbol{\theta}(z) \! = \! \sum_{m
\in \mathbb{Z}^{N}} \me^{2 \pi \mi (m,z)+\pi \mi (m,\tau m)}, \quad z \! \in
\! \mathbb{C}^{N},
\end{equation*}
where $(\boldsymbol{\cdot},\boldsymbol{\cdot})$ denotes the---real---Euclidean 
inner/scalar product (for $\mathbf{A} \! = \! (A_{1},A_{2},\dotsc,A_{N}) \! 
\in \! \mathbb{E}^{N}$ and $\mathbf{B} \! = \! (B_{1},B_{2},\dotsc,B_{N}) \! 
\in \! \mathbb{E}^{N}$, $(A,B) \! := \! \sum_{k=1}^{N}A_{k}B_{k})$, with the 
following evenness and (quasi-) periodicity properties,
\begin{equation*}
\boldsymbol{\theta}(-z) \! = \! \boldsymbol{\theta}(z), \qquad \boldsymbol{
\theta}(z \! + \! e_{j}) \! = \! \boldsymbol{\theta}(z), \qquad \mathrm{and}
\qquad \boldsymbol{\theta}(z \! \pm \! \tau_{j}) \! = \! \me^{\mp 2 \pi \mi
z_{j}-\mi \pi \tau_{jj}} \boldsymbol{\theta}(z),
\end{equation*}
where $e_{j}$ is the standard (basis) column vector in $\mathbb{C}^{N}$ with 
$1$ in the $j$th entry and $0$ elsewhere (see above), and $\tau_{j} \! := \! 
\tau e_{j}$ $(\in \! \mathbb{C}^{N})$, $j \! = \! 1,\dotsc,N$.

It turns out that, for the analysis of this work, the following multi-valued 
functions are essential:
\begin{enumerate}
\item[$\boldsymbol{\bullet}$] $(R_{e}(z))^{1/2} \! := \! (\prod_{k=0}^{N}(z \! 
- \! b_{k}^{e})(z \! - \! a_{k+1}^{e}))^{1/2}$, where, with the identification 
$a_{N+1}^{e} \! \equiv \! a_{0}^{e}$ (as points on the complex sphere, 
$\overline{\mathbb{C}})$ and with the point at infinity lying on the (open) 
interval $(a_{0}^{e},b_{0}^{e})$, $-\infty \! < \! a_{0}^{e} \! < \! b_{0}^{e} 
\! < \! a_{1}^{e} \! < \! b_{1}^{e} \! < \! \cdots \! < \! a_{N}^{e} \! < \! 
b_{N}^{e} \! < \! +\infty$, $a_{0}^{e}$ $(\equiv \! a_{N+1}^{e})$ $\not= \! 
-\infty,0$, and $b_{N}^{e} \! \not= \! 0,+\infty$ (see Figure~1);
\begin{figure}[tbh]
\begin{center}
\vspace{0.35cm}
\begin{pspicture}(0,0)(12,3)
\psset{xunit=1cm,yunit=1cm}
\psline[linewidth=0.9pt,linestyle=solid,linecolor=black]{o-o}(0.5,2)(2,2)
\psline[linewidth=0.9pt,linestyle=solid,linecolor=black]{o-o}(3,2)(4.5,2)
\psline[linewidth=0.9pt,linestyle=solid,linecolor=black]{o-o}(6.5,2)(8,2)
\psline[linewidth=0.9pt,linestyle=solid,linecolor=black]{o-o}(10,2)(11.5,2)
\psline[linewidth=0.7pt,linestyle=dotted,linecolor=darkgray](4.65,2)(6.35,2)
\psline[linewidth=0.7pt,linestyle=dotted,linecolor=darkgray](8.15,2)(9.85,2)
\rput(0.5,1.7){\makebox(0,0){$a_{0}^{e}$}}
\rput(0.5,0.9){\makebox(0,0){$a_{N+1}^{e}$}}
\rput{90}(0.5,1.3){\makebox(0,0){$\equiv$}}
\rput(1.25,2){\makebox(0,0){$\pmb{\times}$}}
\rput(1.25,2.3){\makebox(0,0){$\infty$}}
\rput(2,1.7){\makebox(0,0){$b_{0}^{e}$}}
\rput(3,1.7){\makebox(0,0){$a_{1}^{e}$}}
\rput(4.5,1.7){\makebox(0,0){$b_{1}^{e}$}}
\rput(6.5,1.7){\makebox(0,0){$a_{j}^{e}$}}
\rput(8,1.7){\makebox(0,0){$b_{j}^{e}$}}
\rput(10,1.7){\makebox(0,0){$a_{N}^{e}$}}
\rput(11.5,1.7){\makebox(0,0){$b_{N}^{e}$}}
\end{pspicture}
\end{center}
\vspace{-0.95cm}
\caption{Union of (open) intervals in the complex $z$-plane}
\end{figure}
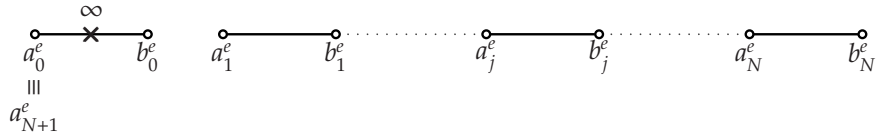
\item[$\boldsymbol{\bullet}$] $(R_{o}(z))^{1/2} \! := \! (\prod_{k=0}^{N}(z \! 
- \! b_{k}^{o})(z \! - \! a_{k+1}^{o}))^{1/2}$, where, with the identification 
$a_{N+1}^{o} \! \equiv \! a_{0}^{o}$ (as points on the complex sphere, 
$\overline{\mathbb{C}})$ and with the point at infinity lying on the (open) 
interval $(a_{0}^{o},b_{0}^{o})$, $-\infty \! < \! a_{0}^{o} \! < \! b_{0}^{o} 
\! < \! a_{1}^{o} \! < \! b_{1}^{o} \! < \! \cdots \! < \! a_{N}^{o} \! < \! 
b_{N}^{o} \! < \! +\infty$, $a_{0}^{o}$ $(\equiv \! a_{N+1}^{o})$ $\not= \! 
-\infty,0$, and $b_{N}^{o} \! \not= \! 0,+\infty$ (see Figure~2).
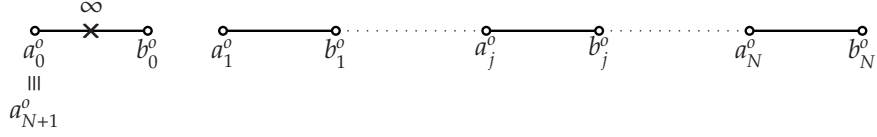
\begin{figure}[tbh]
\begin{center}
\vspace{0.35cm}
\begin{pspicture}(0,0)(12,3)
\psset{xunit=1cm,yunit=1cm}
\psline[linewidth=0.9pt,linestyle=solid,linecolor=black]{o-o}(0.5,2)(2,2)
\psline[linewidth=0.9pt,linestyle=solid,linecolor=black]{o-o}(3,2)(4.5,2)
\psline[linewidth=0.9pt,linestyle=solid,linecolor=black]{o-o}(6.5,2)(8,2)
\psline[linewidth=0.9pt,linestyle=solid,linecolor=black]{o-o}(10,2)(11.5,2)
\psline[linewidth=0.7pt,linestyle=dotted,linecolor=darkgray](4.65,2)(6.35,2)
\psline[linewidth=0.7pt,linestyle=dotted,linecolor=darkgray](8.15,2)(9.85,2)
\rput(0.5,1.7){\makebox(0,0){$a_{0}^{o}$}}
\rput(0.5,0.9){\makebox(0,0){$a_{N+1}^{o}$}}
\rput{90}(0.5,1.3){\makebox(0,0){$\equiv$}}
\rput(1.25,2){\makebox(0,0){$\pmb{\times}$}}
\rput(1.25,2.3){\makebox(0,0){$\infty$}}
\rput(2,1.7){\makebox(0,0){$b_{0}^{o}$}}
\rput(3,1.7){\makebox(0,0){$a_{1}^{o}$}}
\rput(4.5,1.7){\makebox(0,0){$b_{1}^{o}$}}
\rput(6.5,1.7){\makebox(0,0){$a_{j}^{o}$}}
\rput(8,1.7){\makebox(0,0){$b_{j}^{o}$}}
\rput(10,1.7){\makebox(0,0){$a_{N}^{o}$}}
\rput(11.5,1.7){\makebox(0,0){$b_{N}^{o}$}}
\end{pspicture}
\end{center}
\vspace{-0.95cm}
\caption{Union of (open) intervals in the complex $z$-plane}
\end{figure}
\end{enumerate}

The functions $R_{e}(z)$ and $R_{o}(z)$, respectively, are unital polynomials
$(\in \! \mathbb{R}[z])$ of even degree $(\deg (R_{e}(z)) \! = \! \deg (R_{o} 
(z)) \! = \! 2(N \! + \! 1))$ whose (simple) roots/zeros are $\lbrace 
b_{j-1}^{e},a_{j}^{e} \rbrace_{j=1}^{N+1}$ $(a_{N+1}^{e} \! \equiv \! a_{0}^{
e})$ and $\lbrace b_{j-1}^{o},a_{j}^{o} \rbrace_{j=1}^{N+1}$ $(a_{N+1}^{o} \! 
\equiv \! a_{0}^{o})$. The basic ingredients associated with the construction 
of the hyperelliptic Riemann surfaces of genus $N$ corresponding, 
respectively, to the multi-valued functions $y^{2} \! = \! R_{e}(z)$ and 
$y^{2} \! = \! R_{o}(z)$ was given above. One now uses the above construction; 
but particularised to the cases of the polynomials $R_{e}(z)$ and $R_{o}(z)$, 
to arrive at the following:
\begin{enumerate}
\item[\shadowbox{$\sqrt{\smash[b]{R_{e}(z)}}$}] Let $\mathcal{Y}_{e}$ denote 
the two-sheeted Riemann surface of genus $N$ associated with $y^{2} \! = \! 
R_{e}(z)$, with $R_{e}(z)$ as characterised above: the first/upper (resp., 
second/lower) sheet of $\mathcal{Y}_{e}$ is denoted by $\mathcal{Y}_{e}^{+}$ 
(resp., $\mathcal{Y}_{e}^{-})$, points on the first/upper (resp., 
second/lower) sheet are represented as $z^{+} \! := \! (z,+(R_{e}(z))^{1/2})$ 
(resp., $z^{-} \! := \! (z,-(R_{e}(z))^{1/2}))$, where, as points on the plane 
$\mathbb{C}$, $z^{+} \! = \! z^{-} \! = \! z$, and the single-valued branch 
for the square root of the (multi-valued) function $(R_{e}(z))^{1/2}$ is 
chosen such that $z^{-(N+1)}(R_{e}(z))^{1/2} \! \sim_{\underset{z \in 
\mathcal{Y}_{e}^{\pm}}{z \to \infty}} \! \pm 1$. $\mathcal{Y}_{e}$ is realised 
as a (two-sheeted) branched/rami\-f\-i\-e\-d covering of the Riemann sphere 
such that its two sheets are two identical copies of $\mathbb{C}$ with branch 
cuts (slits) along the intervals $(a_{0}^{e},b_{0}^{e}),(a_{1}^{e},b_{1}^{e}),
\dotsc,(a_{N}^{e},b_{N}^{e})$ and pasted/glued together along $\cup_{j=1}^{N+
1}(a_{j-1}^{e},b_{j-1}^{e})$ $(a_{0}^{e} \! \equiv \! a_{N+1}^{e})$ in such a 
way that the cycles $\boldsymbol{\alpha}_{0}^{e}$ and $\{\boldsymbol{\alpha}_{
j}^{e},\boldsymbol{\beta}_{j}^{e}\}$, $j \! = \! 1,\dotsc,N$, where the latter 
forms the canonical $\mathbf{1}$-homology basis for $\mathcal{Y}_{e}$, are 
characterised by the fact that (the closed contours) $\boldsymbol{\alpha}_{
j}^{e}$, $j \! = \! 0,\dotsc,N$, lie on $\mathcal{Y}_{e}^{+}$, and (the closed 
contours) $\boldsymbol{\beta}_{j}^{e}$, $j \! = \! 1,\dotsc,N$, pass {}from 
$\mathcal{Y}_{e}^{+}$ (starting {}from the slit $(a_{j}^{e},b_{j}^{e}))$, 
through the slit $(a_{0}^{e},b_{0}^{e})$ to $\mathcal{Y}_{e}^{-}$, and back 
again to $\mathcal{Y}_{e}^{+}$ through the slit $(a_{j}^{e},b_{j}^{e})$ (see 
Figure~3).
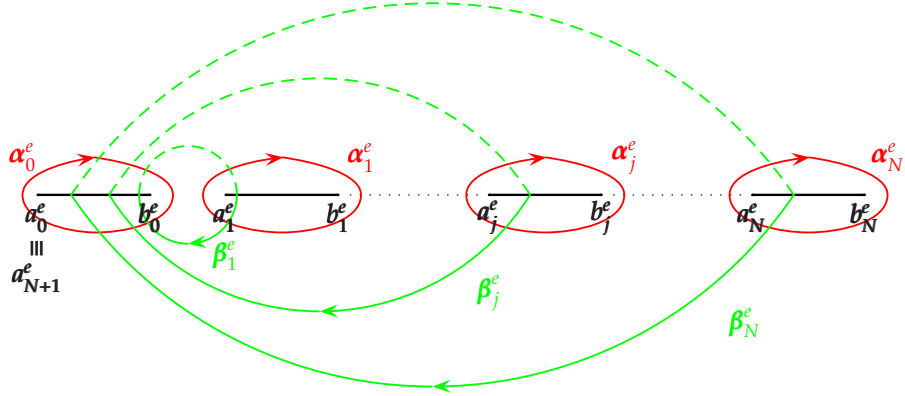
\begin{figure}[tbh]
\begin{center}
\vspace{0.45cm}
\begin{pspicture}(-2,0)(13,6)
\psset{xunit=1cm,yunit=1cm}
\psline[linewidth=0.9pt,linestyle=solid,linecolor=black]{c-c}(0.5,3)(2,3)
\psline[linewidth=0.9pt,linestyle=solid,linecolor=black]{c-c}(3,3)(4.5,3)
\psline[linewidth=0.9pt,linestyle=solid,linecolor=black]{c-c}(6.5,3)(8,3)
\psline[linewidth=0.9pt,linestyle=solid,linecolor=black]{c-c}(10,3)(11.5,3)
\psline[linewidth=0.7pt,linestyle=dotted,linecolor=darkgray](4.65,3)(6.35,3)
\psline[linewidth=0.7pt,linestyle=dotted,linecolor=darkgray](8.15,3)(9.85,3)
\pscurve[linewidth=0.7pt,linestyle=solid,linecolor=red,arrowsize=1.5pt 5]{->}%
(1.25,3.5)(2.3,3)(1.25,2.5)(0.3,3)(1.27,3.5)
\pscurve[linewidth=0.7pt,linestyle=solid,linecolor=red,arrowsize=1.5pt 5]{->}%
(3.75,3.5)(4.8,3)(3.75,2.5)(2.7,3)(3.77,3.5)
\pscurve[linewidth=0.7pt,linestyle=solid,linecolor=red,arrowsize=1.5pt 5]{->}%
(7.25,3.5)(8.3,3)(7.25,2.5)(6.2,3)(7.27,3.5)
\pscurve[linewidth=0.7pt,linestyle=solid,linecolor=red,arrowsize=1.5pt 5]{->}%
(10.75,3.5)(11.8,3)(10.75,2.5)(9.7,3)(10.77,3.5)
\rput(0.3,3.5){\makebox(0,0){{\color{red} $\boldsymbol{\alpha}_{0}^{e}$}}}
\rput(4.8,3.5){\makebox(0,0){{\color{red} $\boldsymbol{\alpha}_{1}^{e}$}}}
\rput(8.3,3.5){\makebox(0,0){{\color{red} $\boldsymbol{\alpha}_{j}^{e}$}}}
\rput(11.8,3.5){\makebox(0,0){{\color{red} $\boldsymbol{\alpha}_{N}^{e}$}}}
\psarc[linewidth=0.7pt,linestyle=solid,linecolor=green,arrowsize=1.5pt 5]{<-}%
(2.5,3){0.65}{270}{360}
\psarc[linewidth=0.7pt,linestyle=solid,linecolor=green](2.5,3){0.65}{180}{270}
\psarc[linewidth=0.7pt,linestyle=dashed,linecolor=green](2.5,3){0.65}{0}{180}
\psarc[linewidth=0.7pt,linestyle=solid,linecolor=green,arrowsize=1.5pt 5]%
{<-}(4.25,4.75){3.3}{270}{328}
\psarc[linewidth=0.7pt,linestyle=solid,linecolor=green](4.25,4.75){3.3}{212}%
{270}
\psarc[linewidth=0.7pt,linestyle=dashed,linecolor=green](4.25,1.25){3.3}{32}%
{148}
\psarc[linewidth=0.7pt,linestyle=solid,linecolor=green,arrowsize=1.5pt 5]%
{<-}(5.75,6.25){5.8}{270}{326}
\psarc[linewidth=0.7pt,linestyle=solid,linecolor=green](5.75,6.25){5.8}%
{214}{270}
\psarc[linewidth=0.7pt,linestyle=dashed,linecolor=green](5.75,-0.25){5.8}%
{34}{146}
\rput(3,2.2){\makebox(0,0){{\color{green} $\boldsymbol{\beta}_{1}^{e}$}}}
\rput(6.5,1.7){\makebox(0,0){{\color{green} $\boldsymbol{\beta}_{j}^{e}$}}}
\rput(9.9,1.3){\makebox(0,0){{\color{green} $\boldsymbol{\beta}_{N}^{e}$}}}
\rput(0.5,2.7){\makebox(0,0){$\pmb{a_{0}^{e}}$}}
\rput{90}(0.5,2.3){\makebox(0,0){$\pmb{\equiv}$}}
\rput(0.5,1.9){\makebox(0,0){$\pmb{a_{N+1}^{e}}$}}
\rput(2,2.7){\makebox(0,0){$\pmb{b_{0}^{e}}$}}
\rput(3,2.7){\makebox(0,0){$\pmb{a_{1}^{e}}$}}
\rput(4.5,2.7){\makebox(0,0){$\pmb{b_{1}^{e}}$}}
\rput(6.5,2.7){\makebox(0,0){$\pmb{a_{j}^{e}}$}}
\rput(8,2.7){\makebox(0,0){$\pmb{b_{j}^{e}}$}}
\rput(10,2.7){\makebox(0,0){$\pmb{a_{N}^{e}}$}}
\rput(11.5,2.7){\makebox(0,0){$\pmb{b_{N}^{e}}$}}
\end{pspicture}
\end{center}
\vspace{-0.55cm}
\caption{The Riemann surface $\mathcal{Y}_{e}$ of $y^{2} \! = \! \prod_{k=0
}^{N}(z \! - \! b_{k}^{e})(z \! - \! a_{k+1}^{e})$, $a_{N+1}^{e} \! \equiv 
\! a_{0}^{e}$. The solid (resp., dashed) lines are on the first/upper (resp., 
second/lower) sheet of $\mathcal{Y}_{e}$, denoted $\mathcal{Y}_{e}^{+}$ 
(resp., $\mathcal{Y}_{e}^{-})$.}
\end{figure}

\hspace*{0.50cm}
The canonical $\mathbf{1}$-homology basis $\{\boldsymbol{\alpha}_{j}^{e},
\boldsymbol{\beta}_{j}^{e}\}$, $j \! = \! 1,\dotsc,N$, generates, on $\mathcal{
Y}_{e}$, the (corresponding) $\boldsymbol{\alpha}^{e}$-normalised basis of 
holomorphic Abelian differentials (one-forms) $\{\omega_{1}^{e},\omega_{2}^{e},
\dotsc,\omega_{N}^{e}\}$, where $\omega_{j}^{e} \! := \! \sum_{k=1}^{N} \tfrac{
c_{jk}^{e}z^{N-k}}{\sqrt{\smash[b]{R_{e}(z)}}} \, \md z$, $c_{jk}^{e} \! \in 
\! \mathbb{C}$, $j \! = \! 1,\dotsc,N$, and $\oint_{\boldsymbol{\alpha}_{k}^{
e}} \omega_{j}^{e} \! = \! \delta_{kj}$, $k,j \! = \! 1,\dotsc,N$: $\omega_{
l}^{e}$, $l \! = \! 1,\dotsc,N$, is real valued on $\cup_{j=1}^{N+1}(a_{j-1}^{
e},b_{j-1}^{e})$, and has exactly one (real) root in any (open) interval $(a_{
j-1}^{e},b_{j-1}^{e})$, $j \! = \! 1,\dotsc,N \! + \! 1$; furthermore, in the 
intervals $(b_{j-1}^{e},a_{j}^{e})$, $j \! = \! 1,\dotsc,N$, $\omega_{l}^{e}$, 
$l \! = \! 1,\dotsc,N$, take non-zero, pure imaginary values. Let $\boldsymbol{
\omega}^{e} \! := \! (\omega_{1}^{e},\omega_{2}^{e},\dotsc,\omega_{N}^{e})$ 
denote the basis of holomorphic one-forms on $\mathcal{Y}_{e}$ as normalised 
above with the associated $N \! \times \! N$ Riemann matrix of $\boldsymbol{
\beta}^{e}$-periods, $\tau^{e} \! = \! (\tau^{e})_{i,j=1,\dotsc,N} \! := \! 
(\oint_{\boldsymbol{\beta}_{j}^{e}} \omega_{i}^{e})_{i,j=1,\dotsc,N}$: the 
Riemann matrix is symmetric $(\tau_{ij}^{e} \! = \! \tau_{ji}^{e})$ and pure 
imaginary, $-\mi \tau^{e}$ is positive definite $(\Im (\tau^{e}_{jj}) \! > \! 
0)$, and $\det (\tau^{e}) \! \not= \! 0$ (non-degenerate). For the holomorphic 
Abelian differential (one-form) $\boldsymbol{\omega}^{e}$ defined above, 
choose $a_{N+1}^{e}$ as the \emph{base point}, and set $\boldsymbol{u}^{e} 
\colon \mathcal{Y}_{e} \! \to \! \operatorname{Jac}(\mathcal{Y}_{e})$ $(:= 
\! \mathbb{C}^{N}/\{N \! + \! \tau^{e}M\}$, $(N,M) \! \in \! \mathbb{Z}^{N} 
\! \times \! \mathbb{Z}^{N})$, $z \! \mapsto \! \boldsymbol{u}^{e}(z) \! := \! 
\int_{a_{N+1}^{e}}^{z} \boldsymbol{\omega}^{e}$, where the integration {}from 
$a_{N+1}^{e}$ to $z$ $(\in \mathcal{Y}_{e})$ is taken along any path on 
$\mathcal{Y}_{e}^{+}$.
\begin{eeee}
{}From the representation $\omega_{j}^{e} \! = \! \sum_{k=1}^{N} \tfrac{c_{j
k}^{e}z^{N-k}}{\sqrt{\smash[b]{R_{e}(z)}}} \, \md z$, $j \! = \! 1,\dotsc,N$, 
and the normalisation condition $\oint_{\boldsymbol{\alpha}_{k}^{e}} \omega_{
j}^{e} \! = \! \delta_{kj}$, $k,j \! = \! 1,\dotsc,N$, one shows that $c_{j
k}^{e}$, $k,j \! = \! 1,\dotsc,N$, are obtained from
\begin{equation}
\begin{pmatrix}
c_{11}^{e} & c_{12}^{e} & \dotsb & c_{1N}^{e} \\
c_{21}^{e} & c_{22}^{e} & \dotsb & c_{2N}^{e} \\
\vdots     & \vdots     & \ddots & \vdots     \\
c_{N1}^{e} & c_{N2}^{e} & \dotsb & c_{NN}^{e}
\end{pmatrix} \! = \! \widetilde{\mathfrak{S}}_{e}^{-1} \tag{E1},
\end{equation}
where
\begin{equation}
\widetilde{\mathfrak{S}}_{e} \! := \!
\begin{pmatrix}
\oint_{\boldsymbol{\alpha}_{1}^{e}} \frac{\md s_{1}}{\sqrt{\smash[b]{R_{e}
(s_{1})}}} & \oint_{\boldsymbol{\alpha}_{2}^{e}} \frac{\md s_{2}}{\sqrt{
\smash[b]{R_{e}(s_{2})}}} & \dotsb & \oint_{\boldsymbol{\alpha}_{N}^{e}}
\frac{\md s_{N}}{\sqrt{\smash[b]{R_{e}(s_{N})}}} \\
\oint_{\boldsymbol{\alpha}_{1}^{e}} \frac{s_{1} \md s_{1}}{\sqrt{\smash[b]{
R_{e}(s_{1})}}} & \oint_{\boldsymbol{\alpha}_{2}^{e}} \frac{s_{2} \md s_{2}}{
\sqrt{\smash[b]{R_{e}(s_{2})}}} & \dotsb & \oint_{\boldsymbol{\alpha}_{N}^{e}}
\frac{s_{N} \md s_{N}}{\sqrt{\smash[b]{R_{e}(s_{N})}}} \\
\vdots & \vdots & \ddots & \vdots \\
\oint_{\boldsymbol{\alpha}_{1}^{e}} \frac{s_{1}^{N-1} \md s_{1}}{\sqrt{
\smash[b]{R_{e}(s_{1})}}} & \oint_{\boldsymbol{\alpha}_{2}^{e}} \frac{s_{
2}^{N-1} \md s_{2}}{\sqrt{\smash[b]{R_{e}(s_{2})}}} & \dotsb & \oint_{
\boldsymbol{\alpha}_{N}^{e}} \frac{s_{N}^{N-1} \md s_{N}}{\sqrt{\smash[b]{
R_{e}(s_{N})}}}
\end{pmatrix} \tag{E2}.
\end{equation}
For a (representation-independent) proof of the fact that $\det (\widetilde{
\mathfrak{S}}_{e}) \! \not= \! 0$, see, for example, Chapter~10, 
Section~10--2, of \cite{a87}. \hfill $\blacksquare$
\end{eeee}

\hspace*{0.50cm}
Set (see Section~4), for $z \! \in \! \mathbb{C}_{+}$, $\gamma^{e}(z) \! := \!
(\prod_{k=1}^{N+1}(z \! - \! b_{k-1}^{e})(z \! - \! a_{k}^{e})^{-1})^{1/4}$,
and, for $z \! \in \! \mathbb{C}_{-}$, $\gamma^{e}(z) \! := \! -\mi (\prod_{k
=1}^{N+1}(z \! - \! b_{k-1}^{e})(z \! - \! a_{k}^{e})^{-1})^{1/4}$. It is
shown in Section~4 that $\gamma^{e}(z) \! =_{\underset{z \in \mathcal{Y}_{
e}^{\pm}}{z \to \infty}} \! (-\mi)^{(1 \mp 1)/2} \linebreak[4]
\cdot (1 \! + \! \mathcal{O}(z^{-1}))$, and
\begin{equation*}
\left\{z_{j}^{e,\pm} \right\}_{j=1}^{N} \! = \! \left\{\mathstrut z^{\pm} \!
\in \! \mathcal{Y}_{e}^{\pm}; \, (\gamma^{e}(z) \! \mp \! (\gamma^{e}(z))^{-
1}) \vert_{z=z^{\pm}} \! = \! 0 \right\},
\end{equation*}
with $z_{j}^{e,\pm} \! \in \! (a_{j}^{e},b_{j}^{e})^{\pm}$ $(\subset \!
\mathcal{Y}_{e}^{\pm})$, $j \! = \! 1,\dotsc,N$, where, as points on the
plane, $z_{j}^{e,+} \! = \! z_{j}^{e,-} \! := \! z_{j}^{e}$, $j \! = \! 1,
\dotsc,N$ (of course, on the plane, $z_{j}^{e} \! \in \! (a_{j}^{e},b_{j}^{
e})$, $j \! = \! 1,\dotsc,N)$.

\hspace*{0.50cm}
Corresponding to $\mathcal{Y}_{e}$, define $\boldsymbol{d}_{e} \! := \!
-\boldsymbol{K}_{e} \! - \! \sum_{j=1}^{N} \int_{a_{N+1}^{e}}^{z_{j}^{e,-}}
\boldsymbol{\omega}^{e}$ $(\in \! \mathbb{C}^{N})$, where $\boldsymbol{K}_{e}$
is the associated (`even') vector of Riemann constants, and the integration
{}from $a_{N+1}^{e}$ to $z_{j}^{e,-}$, $j \! = \! 1,\dotsc,N$, is taken along
a fixed path in $\mathcal{Y}_{e}^{-}$. It is shown in Chapter~VII of
\cite{a88} that $\boldsymbol{K}_{e} \! = \! \sum_{j=1}^{N} \int_{a_{j}^{e}}^{
a_{N+1}^{e}} \boldsymbol{\omega}^{e}$; furthermore, $\boldsymbol{K}_{e}$
is a point of order $2$, that is, $2 \boldsymbol{K}_{e} \! = \! 0$ and $s
\boldsymbol{K}_{e} \! \not= \! 0$ for $0 \! < \! s \! < \! 2$. Recalling the
definition of $\boldsymbol{\omega}^{e}$ and that $z^{-(N+1)}(R_{e}(z))^{1/2}
\! \sim_{\underset{z \in \mathbb{C}_{\pm}}{z \to \infty}} \! \pm 1$, using the
fact that $\boldsymbol{K}_{e}$ is a point of order $2$, one arrives at
\begin{align*}
\boldsymbol{d}_{e} =& \, -\boldsymbol{K}_{e} \! - \! \sum_{j=1}^{N} \int_{a_{N
+1}^{e}}^{z_{j}^{e,-}} \boldsymbol{\omega}^{e} \! = \! \boldsymbol{K}_{e} \! -
\! \sum_{j=1}^{N} \int_{a_{N+1}^{e}}^{z_{j}^{e,-}} \boldsymbol{\omega}^{e} \!
= \! -\boldsymbol{K}_{e} \! + \! \sum_{j=1}^{N} \int_{a_{N+1}^{e}}^{z_{j}^{e,
+}} \boldsymbol{\omega}^{e} \! = \! \boldsymbol{K}_{e} \! + \! \sum_{j=1}^{N}
\int_{a_{N+1}^{e}}^{z_{j}^{e,+}} \boldsymbol{\omega}^{e} \\
=& \, -\sum_{j=1}^{N} \int_{a_{j}^{e}}^{z_{j}^{e,-}} \boldsymbol{\omega}^{e}
\! = \! \sum_{j=1}^{N} \int_{a_{j}^{e}}^{z_{j}^{e,+}} \boldsymbol{\omega}^{e}.
\end{align*}

\hspace*{0.50cm}
Associated with the Riemann matrix of $\boldsymbol{\beta}^{e}$-periods,
$\tau^{e}$, is the (`even') Riemann theta function:
\begin{equation}
\boldsymbol{\theta}(z;\tau^{e}) \! =: \! \boldsymbol{\theta}^{e}(z) \! = \!
\sum_{m \in \mathbb{Z}^{N}} \me^{2 \pi \mi (m,z)+\pi \mi (m,\tau^{e}m)}, \quad
z \! \in \! \mathbb{C}^{N};
\end{equation}
$\boldsymbol{\theta}^{e}(z)$ has the following evenness and (quasi-)
periodicity properties,
\begin{equation*}
\boldsymbol{\theta}^{e}(-z) \! = \! \boldsymbol{\theta}^{e}(z), \qquad
\boldsymbol{\theta}^{e}(z \! + \! e_{j}) \! = \! \boldsymbol{\theta}^{e}(z),
\qquad \mathrm{and} \qquad \boldsymbol{\theta}^{e}(z \! \pm \! \tau_{j}^{e})
\! = \! \me^{\mp 2 \pi \mi z_{j}-\mi \pi \tau_{jj}^{e}} \boldsymbol{\theta}^{e}
(z),
\end{equation*}
where $\tau_{j}^{e} \! := \! \tau^{e} e_{j}$ $(\in \! \mathbb{C}^{N})$, $j
\! = \! 1,\dotsc,N$. Extensive use of this entire apparatus will be made in
Section~4.

\item[\shadowbox{$\sqrt{\smash[b]{R_{o}(z)}}$}] Let $\mathcal{Y}_{o}$ denote
the two-sheeted Riemann surface of genus $N$ associated with $y^{2} \! = \!
R_{o}(z)$, with $R_{o}(z)$ as characterised above: the first/upper (resp.,
second/lower) sheet of $\mathcal{Y}_{o}$ is denoted by $\mathcal{Y}_{o}^{+}$
(resp., $\mathcal{Y}_{o}^{-})$, points on the first/upper (resp.,
second/lower) sheet are represented as $z^{+} \! := \! (z,+(R_{o}(z))^{1/2})$
(resp., $z^{-} \! := \! (z,-(R_{o}(z))^{1/2}))$, where, as points on the plane
$\mathbb{C}$, $z^{+} \! = \! z^{-} \! = \! z$, and the single-valued branch
for the square root of the (multi-valued) function $(R_{o}(z))^{1/2}$ is
chosen such that $z^{-(N+1)}(R_{o}(z))^{1/2} \! \sim_{\underset{z \in
\mathcal{Y}_{o}^{\pm}}{z \to \infty}} \! \pm 1$. $\mathcal{Y}_{o}$ is realised
as a (two-sheeted) branched/rami\-f\-i\-e\-d covering of the Riemann sphere
such that its two sheets are two identical copies of $\mathbb{C}$ with branch
cuts (slits) along the intervals $(a_{0}^{o},b_{0}^{o}),(a_{1}^{o},b_{1}^{o}),
\dotsc,(a_{N}^{o},b_{N}^{o})$ and pasted/glued together along $\cup_{j=1}^{N+1}
(a_{j-1}^{o},b_{j-1}^{o})$ $(a_{0}^{o} \! \equiv \! a_{N+1}^{o})$ in such a
way that the cycles $\boldsymbol{\alpha}_{0}^{o}$ and $\{\boldsymbol{\alpha}_{
j}^{o},\boldsymbol{\beta}_{j}^{o}\}$, $j \! = \! 1,\dotsc,N$, where the latter
forms the canonical $\mathbf{1}$-homology basis for $\mathcal{Y}_{o}$, are
characterised by the fact that (the closed contours) $\boldsymbol{\alpha}_{
j}^{o}$, $j \! = \! 0,\dotsc,N$, lie on $\mathcal{Y}_{o}^{+}$, and (the closed
contours) $\boldsymbol{\beta}_{j}^{o}$, $j \! = \! 1,\dotsc,N$, pass {}from
$\mathcal{Y}_{o}^{+}$ (starting {}from the slit $(a_{j}^{o},b_{j}^{o}))$,
through the slit $(a_{0}^{o},b_{0}^{o})$ to $\mathcal{Y}_{o}^{-}$, and back
again to $\mathcal{Y}_{o}^{+}$ through the slit $(a_{j}^{o},b_{j}^{o})$ (see
Figure~4).
\begin{figure}[htb]
\begin{center}
\vspace{0.45cm}
\begin{pspicture}(-2,0)(13,6)
\psset{xunit=1cm,yunit=1cm}
\psline[linewidth=0.9pt,linestyle=solid,linecolor=black]{c-c}(0.5,3)(2,3)
\psline[linewidth=0.9pt,linestyle=solid,linecolor=black]{c-c}(3,3)(4.5,3)
\psline[linewidth=0.9pt,linestyle=solid,linecolor=black]{c-c}(6.5,3)(8,3)
\psline[linewidth=0.9pt,linestyle=solid,linecolor=black]{c-c}(10,3)(11.5,3)
\psline[linewidth=0.7pt,linestyle=dotted,linecolor=darkgray](4.65,3)(6.35,3)
\psline[linewidth=0.7pt,linestyle=dotted,linecolor=darkgray](8.15,3)(9.85,3)
\pscurve[linewidth=0.7pt,linestyle=solid,linecolor=red,arrowsize=1.5pt 5]{->}%
(1.25,3.5)(2.3,3)(1.25,2.5)(0.3,3)(1.27,3.5)
\pscurve[linewidth=0.7pt,linestyle=solid,linecolor=red,arrowsize=1.5pt 5]{->}%
(3.75,3.5)(4.8,3)(3.75,2.5)(2.7,3)(3.77,3.5)
\pscurve[linewidth=0.7pt,linestyle=solid,linecolor=red,arrowsize=1.5pt 5]{->}%
(7.25,3.5)(8.3,3)(7.25,2.5)(6.2,3)(7.27,3.5)
\pscurve[linewidth=0.7pt,linestyle=solid,linecolor=red,arrowsize=1.5pt 5]{->}%
(10.75,3.5)(11.8,3)(10.75,2.5)(9.7,3)(10.77,3.5)
\rput(0.3,3.5){\makebox(0,0){{\color{red} $\boldsymbol{\alpha}_{0}^{o}$}}}
\rput(4.8,3.5){\makebox(0,0){{\color{red} $\boldsymbol{\alpha}_{1}^{o}$}}}
\rput(8.3,3.5){\makebox(0,0){{\color{red} $\boldsymbol{\alpha}_{j}^{o}$}}}
\rput(11.8,3.5){\makebox(0,0){{\color{red} $\boldsymbol{\alpha}_{N}^{o}$}}}
\psarc[linewidth=0.7pt,linestyle=solid,linecolor=green,arrowsize=1.5pt 5]{<-}%
(2.5,3){0.65}{270}{360}
\psarc[linewidth=0.7pt,linestyle=solid,linecolor=green](2.5,3){0.65}{180}{270}
\psarc[linewidth=0.7pt,linestyle=dashed,linecolor=green](2.5,3){0.65}{0}{180}
\psarc[linewidth=0.7pt,linestyle=solid,linecolor=green,arrowsize=1.5pt 5]{<-}%
(4.25,4.75){3.3}{270}{328}
\psarc[linewidth=0.7pt,linestyle=solid,linecolor=green](4.25,4.75){3.3}{212}%
{270}
\psarc[linewidth=0.7pt,linestyle=dashed,linecolor=green](4.25,1.25){3.3}{32}%
{148}
\psarc[linewidth=0.7pt,linestyle=solid,linecolor=green,arrowsize=1.5pt 5]{<-}%
(5.75,6.25){5.8}{270}{326}
\psarc[linewidth=0.7pt,linestyle=solid,linecolor=green](5.75,6.25){5.8}%
{214}{270}
\psarc[linewidth=0.7pt,linestyle=dashed,linecolor=green](5.75,-0.25){5.8}%
{34}{146}
\rput(3,2.2){\makebox(0,0){{\color{green} $\boldsymbol{\beta}_{1}^{o}$}}}
\rput(6.5,1.7){\makebox(0,0){{\color{green} $\boldsymbol{\beta}_{j}^{o}$}}}
\rput(9.9,1.3){\makebox(0,0){{\color{green} $\boldsymbol{\beta}_{N}^{o}$}}}
\rput(0.5,2.7){\makebox(0,0){$\pmb{a_{0}^{o}}$}}
\rput{90}(0.5,2.3){\makebox(0,0){$\pmb{\equiv}$}}
\rput(0.5,1.9){\makebox(0,0){$\pmb{a_{N+1}^{o}}$}}
\rput(2,2.7){\makebox(0,0){$\pmb{b_{0}^{o}}$}}
\rput(3,2.7){\makebox(0,0){$\pmb{a_{1}^{o}}$}}
\rput(4.5,2.7){\makebox(0,0){$\pmb{b_{1}^{o}}$}}
\rput(6.5,2.7){\makebox(0,0){$\pmb{a_{j}^{o}}$}}
\rput(8,2.7){\makebox(0,0){$\pmb{b_{j}^{o}}$}}
\rput(10,2.7){\makebox(0,0){$\pmb{a_{N}^{o}}$}}
\rput(11.5,2.7){\makebox(0,0){$\pmb{b_{N}^{o}}$}}
\end{pspicture}
\end{center}
\vspace{-0.55cm}
\caption{The Riemann surface $\mathcal{Y}_{o}$ of $y^{2} \! = \! \prod_{k=0
}^{N}(z \! - \! b_{k}^{o})(z \! - \! a_{k+1}^{o})$, $a_{N+1}^{o} \! \equiv \!
a_{0}^{o}$. The solid (resp., dashed) lines are on the first/upper (resp.,
second/lower) sheet of $\mathcal{Y}_{o}$, denoted $\mathcal{Y}_{o}^{+}$
(resp., $\mathcal{Y}_{o}^{-})$.}
\end{figure}
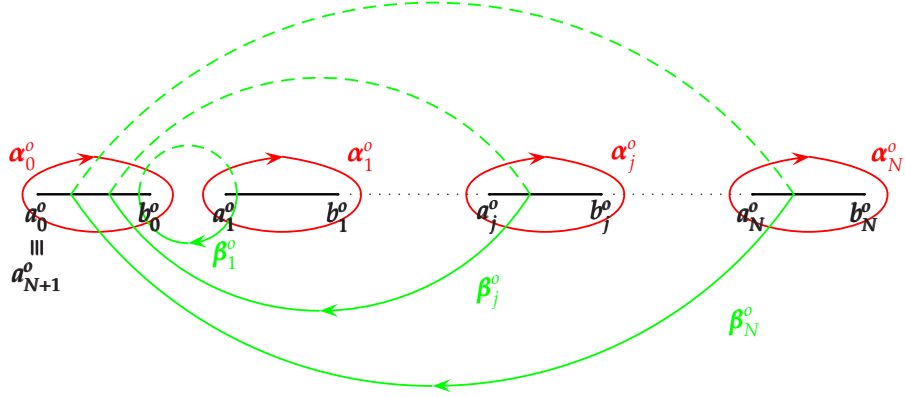

\hspace*{0.50cm}
The canonical $\mathbf{1}$-homology basis $\{\boldsymbol{\alpha}_{j}^{o},
\boldsymbol{\beta}_{j}^{o}\}$, $j \! = \! 1,\dotsc,N$, generates, on $\mathcal{
Y}_{o}$, the (corresponding) $\boldsymbol{\alpha}^{o}$-normalised basis of
holomorphic Abelian differentials (one-forms) $\{\omega_{1}^{o},\omega_{2}^{o},
\dotsc,\omega_{N}^{o}\}$, where $\omega_{j}^{o} \! := \! \sum_{k=1}^{N} \tfrac{
c_{jk}^{o}z^{N-k}}{\sqrt{\smash[b]{R_{o}(z)}}} \, \md z$, $c_{jk}^{o} \! \in
\! \mathbb{C}$, $j \! = \! 1,\dotsc,N$, and $\oint_{\boldsymbol{\alpha}_{k}^{
o}} \omega_{j}^{o} \! = \! \delta_{kj}$, $k,j \! = \! 1,\dotsc,N$: $\omega_{
l}^{o}$, $l \! = \! 1,\dotsc,N$, is real valued on $\cup_{j=1}^{N+1}(a_{j-1}^{
o},b_{j-1}^{o})$, and has exactly one (real) root in any (open) interval $(a_{j
-1}^{o},b_{j-1}^{o})$, $j \! = \! 1,\dotsc,N \! + \! 1$; furthermore, in the
intervals $(b_{j-1}^{o},a_{j}^{o})$, $j \! = \! 1,\dotsc,N$, $\omega_{l}^{o}$,
$l \! = \! 1,\dotsc,N$, take non-zero, pure imaginary values. Let $\boldsymbol{
\omega}^{o} \! := \! (\omega_{1}^{o},\omega_{2}^{o},\dotsc,\omega_{N}^{o})$
denote the basis of holomorphic one-forms on $\mathcal{Y}_{o}$ as normalised
above with the associated $N \! \times \! N$ Riemann matrix of $\boldsymbol{
\beta}^{o}$-periods, $\tau^{o} \! = \! (\tau^{o})_{i,j=1,\dotsc,N} \! := \!
(\oint_{\boldsymbol{\beta}_{j}^{o}} \omega_{i}^{o})_{i,j=1,\dotsc,N}$: the
Riemann matrix is symmetric $(\tau_{ij}^{o} \! = \! \tau_{ji}^{o})$ and
pure imaginary, $-\mi \tau^{o}$ is positive definite $(\Im (\tau^{o}_{jj})
\! > \! 0)$, and $\det (\tau^{o}) \! \not= \! 0$ (non-degenerate). For the
holomorphic Abelian differential (one-form) $\boldsymbol{\omega}^{o}$ defined
above, choose $a_{N+1}^{o}$ as the base point, and set $\boldsymbol{u}^{o}
\colon \mathcal{Y}_{o} \! \to \! \operatorname{Jac}(\mathcal{Y}_{o})$ $(:= \!
\mathbb{C}^{N}/\{N \! + \! \tau^{o}M\}$, $(N,M) \! \in \! \mathbb{Z}^{N} \!
\times \! \mathbb{Z}^{N})$, $z \! \mapsto \! \boldsymbol{u}^{o}(z) \! := \!
\int_{a_{N+1}^{o}}^{z} \boldsymbol{\omega}^{o}$, where the integration {}from
$a_{N+1}^{o}$ to $z$ $(\in \mathcal{Y}_{o})$ is taken along any path on
$\mathcal{Y}_{o}^{+}$.
\begin{eeee}
{}From the representation $\omega_{j}^{o} \! = \! \sum_{k=1}^{N} \tfrac{c_{j
k}^{o}z^{N-k}}{\sqrt{\smash[b]{R_{o}(z)}}} \, \md z$, $j \! = \! 1,\dotsc,N$,
and the normalisation condition $\oint_{\boldsymbol{\alpha}_{k}^{o}} \omega_{
j}^{o} \! = \! \delta_{kj}$, $k,j \! = \! 1,\dotsc,N$, one shows that $c_{j
k}^{o}$, $k,j \! = \! 1,\dotsc,N$, are obtained from
\begin{equation}
\begin{pmatrix}
c_{11}^{o} & c_{12}^{o} & \dotsb & c_{1N}^{o} \\
c_{21}^{o} & c_{22}^{o} & \dotsb & c_{2N}^{o} \\
\vdots     & \vdots     & \ddots & \vdots     \\
c_{N1}^{o} & c_{N2}^{o} & \dotsb & c_{NN}^{o}
\end{pmatrix} \! = \! \widetilde{\mathfrak{S}}_{o}^{-1} \tag{O1},
\end{equation}
where
\begin{equation}
\widetilde{\mathfrak{S}}_{o} \! := \!
\begin{pmatrix}
\oint_{\boldsymbol{\alpha}_{1}^{o}} \frac{\md s_{1}}{\sqrt{\smash[b]{R_{o}
(s_{1})}}} & \oint_{\boldsymbol{\alpha}_{2}^{o}} \frac{\md s_{2}}{\sqrt{
\smash[b]{R_{o}(s_{2})}}} & \dotsb & \oint_{\boldsymbol{\alpha}_{N}^{o}}
\frac{\md s_{N}}{\sqrt{\smash[b]{R_{o}(s_{N})}}} \\
\oint_{\boldsymbol{\alpha}_{1}^{o}} \frac{s_{1} \md s_{1}}{\sqrt{\smash[b]{
R_{o}(s_{1})}}} & \oint_{\boldsymbol{\alpha}_{2}^{o}} \frac{s_{2} \md s_{2}}{
\sqrt{\smash[b]{R_{o}(s_{2})}}} & \dotsb & \oint_{\boldsymbol{\alpha}_{N}^{o}}
\frac{s_{N} \md s_{N}}{\sqrt{\smash[b]{R_{o}(s_{N})}}} \\
\vdots & \vdots & \ddots & \vdots \\
\oint_{\boldsymbol{\alpha}_{1}^{o}} \frac{s_{1}^{N-1} \md s_{1}}{\sqrt{
\smash[b]{R_{o}(s_{1})}}} & \oint_{\boldsymbol{\alpha}_{2}^{o}} \frac{s_{
2}^{N-1} \md s_{2}}{\sqrt{\smash[b]{R_{o}(s_{2})}}} & \dotsb & \oint_{
\boldsymbol{\alpha}_{N}^{o}} \frac{s_{N}^{N-1} \md s_{N}}{\sqrt{\smash[b]{
R_{o}(s_{N})}}}
\end{pmatrix} \tag{O2}.
\end{equation}
For a (representation-independent) proof of the fact that $\det (\widetilde{
\mathfrak{S}}_{o}) \! \not= \! 0$, see, for example, Chapter~10,
Section~10--2, of \cite{a87}. \hfill $\blacksquare$
\end{eeee}

\hspace*{0.50cm}
Set (see \cite{a51}), for $z \! \in \! \mathbb{C}_{+}$, $\gamma^{o}(z) \! :=
\! (\prod_{k=1}^{N+1}(z \! - \! b_{k-1}^{o})(z \! - \! a_{k}^{o})^{-1})^{
1/4}$, and, for $z \! \in \! \mathbb{C}_{-}$, $\gamma^{o}(z) \! := \! -\mi
(\prod_{k=1}^{N+1}(z \! - \! b_{k-1}^{o})(z \! - \! a_{k}^{o})^{-1})^{1/4}$.
It is shown in \cite{a51} that $\gamma^{o}(z) \! =_{\underset{z \in \mathcal{
Y}_{o}^{\pm}}{z \to 0}} \! (-\mi)^{(1 \mp 1)/2} \gamma^{o}(0)(1 \! + \!
\mathcal{O}(z))$, where $\gamma^{o}(0) \! := \! (\prod_{k=1}^{N+1}b_{k-1}^{o}
/a_{k}^{o})^{1/4}$ $> \! 0$, and a set of $N$ upper-edge and lower-edge
finite-length-gap roots/zeros are
\begin{equation*}
\left\{z_{j}^{o,\pm} \right\}_{j=1}^{N} \! = \! \left\{\mathstrut z^{\pm} \!
\in \! \mathcal{Y}_{o}^{\pm}; \, ((\gamma^{o}(0))^{-1} \gamma^{o}(z) \! \mp
\! \gamma^{o}(0)(\gamma^{o}(z))^{-1}) \vert_{z=z^{\pm}} \! = \! 0 \right\},
\end{equation*}
with $z_{j}^{o,\pm} \! \in \! (a_{j}^{o},b_{j}^{o})^{\pm}$ $(\subset \!
\mathcal{Y}_{o}^{\pm})$, $j \! = \! 1,\dotsc,N$, where, as points on the
plane, $z_{j}^{o,+} \! = \! z_{j}^{o,-} \! := \! z_{j}^{o}$, $j \! = \!
1,\dotsc,N$ (of course, on the plane, $z_{j}^{o} \! \in \! (a_{j}^{o},
b_{j}^{o})$, $j \! = \! 1,\dotsc,N)$.

\hspace*{0.50cm}
Corresponding to $\mathcal{Y}_{o}$, define $\boldsymbol{d}_{o} \! := \! -
\boldsymbol{K}_{o} \! - \! \sum_{j=1}^{N} \int_{a_{N+1}^{o}}^{z_{j}^{o,-}}
\boldsymbol{\omega}^{o}$ $(\in \! \mathbb{C}^{N})$, where $\boldsymbol{K}_{o}$
is the associated (`odd') vector of Riemann constants, and the integration
{}from $a_{N+1}^{o}$ to $z_{j}^{o,-}$, $j \! = \! 1,\dotsc,N$, is taken along
a fixed path in $\mathcal{Y}_{o}^{-}$. It is shown in Chapter~VII of
\cite{a88} that $\boldsymbol{K}_{o} \! = \! \sum_{j=1}^{N} \int_{a_{j}^{o}}^{
a_{N+1}^{o}} \boldsymbol{\omega}^{o}$; furthermore, $\boldsymbol{K}_{o}$ is a
point of order $2$. Recalling the definition of $\boldsymbol{\omega}^{o}$ and
that $z^{-(N+1)}(R_{o}(z))^{1/2} \! \sim_{\underset{z \in \mathbb{C}_{\pm}}{z
\to \infty}} \! \pm 1$, using the fact that $\boldsymbol{K}_{o}$ is a point of
order $2$, one arrives at
\begin{align*}
\boldsymbol{d}_{o} =& \, -\boldsymbol{K}_{o} \! - \! \sum_{j=1}^{N} \int_{a_{N
+1}^{o}}^{z_{j}^{o,-}} \boldsymbol{\omega}^{o} \! = \! \boldsymbol{K}_{o} \! -
\! \sum_{j=1}^{N} \int_{a_{N+1}^{o}}^{z_{j}^{o,-}} \boldsymbol{\omega}^{o} \!
= \! -\boldsymbol{K}_{o} \! + \! \sum_{j=1}^{N} \int_{a_{N+1}^{o}}^{z_{j}^{o,
+}} \boldsymbol{\omega}^{o} \! = \! \boldsymbol{K}_{o} \! + \! \sum_{j=1}^{N}
\int_{a_{N+1}^{o}}^{z_{j}^{o,+}} \boldsymbol{\omega}^{o} \\
=& \, -\sum_{j=1}^{N} \int_{a_{j}^{o}}^{z_{j}^{o,-}} \boldsymbol{\omega}^{o}
\! = \! \sum_{j=1}^{N} \int_{a_{j}^{o}}^{z_{j}^{o,+}} \boldsymbol{\omega}^{o}.
\end{align*}

\hspace*{0.50cm}
Associated with the Riemann matrix of $\boldsymbol{\beta}^{o}$-periods,
$\tau^{o}$, is the (`odd') Riemann theta function:
\begin{equation*}
\boldsymbol{\theta}(z;\tau^{o}) \! =: \! \boldsymbol{\theta}^{o}(z) \! = \!
\sum_{m \in \mathbb{Z}^{N}} \me^{2 \pi \mi (m,z)+\pi \mi (m,\tau^{o}m)}, \quad
z \! \in \! \mathbb{C}^{N};
\end{equation*}
$\boldsymbol{\theta}^{o}(z)$ has the following evenness and (quasi-)
periodicity properties,
\begin{equation*}
\boldsymbol{\theta}^{o}(-z) \! = \! \boldsymbol{\theta}^{o}(z), \qquad
\boldsymbol{\theta}^{o}(z \! + \! e_{j}) \! = \! \boldsymbol{\theta}^{o}(z),
\qquad \mathrm{and} \qquad \boldsymbol{\theta}^{o}(z \! \pm \! \tau_{j}^{o})
\! = \! \me^{\mp 2 \pi \mi z_{j}-\mi \pi \tau_{jj}^{o}} \boldsymbol{\theta}^{
o}(z),
\end{equation*}
where $\tau_{j}^{o} \! := \! \tau^{o} e_{j}$ $(\in \! \mathbb{C}^{N})$, $j \!
= \! 1,\dotsc,N$. This entire latter apparatus is used extensively in 
\cite{a51}.
\end{enumerate}
\subsection{The Riemann-Hilbert Problems for the Monic OLPs}
In this subsection, the RHPs corresponding to the even degree and odd degree 
monic OLPs $\boldsymbol{\pi}_{2n}(z)$ and $\boldsymbol{\pi}_{2n+1}(z)$, 
defined, respectively, in Equations~(1.4) and~(1.5), are formulated 
\emph{\`{a} la} Fokas-Its-Kitaev \cite{a53,a54}. Furthermore, integral 
representations for the even degree and odd degree monic OLPs are also 
obtained.

Consider the varying exponential measure $\widetilde{\mu}$ $(\in \! \mathcal{
M}_{1}(\mathbb{R}))$ given by $\md \widetilde{\mu}(z) \! = \! \me^{-\mathscr{
N} \, V(z)} \, \md z$, $\mathscr{N} \! \in \! \mathbb{N}$, where (the external
field) $V \colon \mathbb{R} \setminus \{0\} \! \to \! \mathbb{R}$ satisfies
conditions~(V1)--(V3). The RHPs which characterise the even degree and odd
degree monic OLPs are now stated.
\begin{rhp1}
Let $V \colon \mathbb{R} \setminus \{0\} \! \to \! \mathbb{R}$ satisfy 
conditions~{\rm (V1)--(V3)}. Find $\overset{e}{\mathrm{Y}} \colon \mathbb{C} 
\setminus \mathbb{R} \! \to \! \mathrm{SL}_{2}(\mathbb{C})$ solving: {\rm (i)} 
$\overset{e}{\mathrm{Y}}(z)$ is holomorphic for $z \! \in \! \mathbb{C} 
\setminus \mathbb{R};$ {\rm (ii)} the boundary values $\overset{e}{\mathrm{
Y}}_{\pm}(z) \! := \! \lim_{\underset{\pm \Im (z^{\prime})>0}{z^{\prime} \to 
z}} \overset{e}{\mathrm{Y}}(z^{\prime})$ satisfy the jump condition
\begin{equation*}
\overset{e}{\mathrm{Y}}_{+}(z)= \overset{e}{\mathrm{Y}}_{-}(z) \! \left(
\mathrm{I} \! + \! \me^{-\mathscr{N} \, V(z)} \sigma_{+} \right), \quad z \!
\in \! \mathbb{R};
\end{equation*}
{\rm (iii)} $\overset{e}{\mathrm{Y}}(z)z^{-n \sigma_{3}} \! =_{\underset{z \in
\mathbb{C} \setminus \mathbb{R}}{z \to \infty}} \! \mathrm{I} \! + \! \mathcal{
O}(z^{-1});$ and {\rm (iv)} $\overset{e}{\mathrm{Y}}(z)z^{n \sigma_{3}} \! =_{
\underset{z \in \mathbb{C} \setminus \mathbb{R}}{z \to 0}} \! \mathcal{O}(1)$.
\end{rhp1}
\begin{rhp2}
Let $V \colon \mathbb{R} \setminus \{0\} \! \to \! \mathbb{R}$ satisfy 
conditions~{\rm (V1)--(V3)}. Find $\overset{o}{\mathrm{Y}} \colon \mathbb{C} 
\setminus \mathbb{R} \! \to \! \mathrm{SL}_{2}(\mathbb{C})$ solving: {\rm (i)} 
$\overset{o}{\mathrm{Y}}(z)$ is holomorphic for $z \! \in \! \mathbb{C} 
\setminus \mathbb{R};$ {\rm (ii)} the boundary values $\overset{o}{\mathrm{
Y}}_{\pm}(z) \! := \! \lim_{\underset{\pm \Im (z^{\prime})>0}{z^{\prime} \to 
z}} \overset{o}{\mathrm{Y}}(z^{\prime})$ satisfy the jump condition
\begin{equation*}
\overset{o}{\mathrm{Y}}_{+}(z)=\overset{o}{\mathrm{Y}}_{-}(z) \! \left(
\mathrm{I} \! + \! \me^{-\mathscr{N} \, V(z)} \sigma_{+} \right), \quad z \!
\in \! \mathbb{R};
\end{equation*}
{\rm (iii)} $\overset{o}{\mathrm{Y}}(z)z^{n \sigma_{3}} \! =_{\underset{z \in
\mathbb{C} \setminus \mathbb{R}}{z \to 0}} \! \mathrm{I} \! + \! \mathcal{O}
(z);$ and {\rm (iv)} $\overset{o}{\mathrm{Y}}(z)z^{-(n+1) \sigma_{3}} \! =_{
\underset{z \in \mathbb{C} \setminus \mathbb{R}}{z \to \infty}} \! \mathcal{O}
(1)$.
\end{rhp2}
\begin{cccc}
Let $\overset{e}{\mathrm{Y}} \colon \mathbb{C} \setminus \mathbb{R} \! \to
\! \mathrm{SL}_{2}(\mathbb{C})$ solve {\rm \pmb{RHP1}}. {\rm \pmb{RHP1}}
possesses a unique solution given by: {\rm (i)} for $n \! = \! 0$,
\begin{equation*}
\overset{e}{\mathrm{Y}}(z) \! = \!
\begin{pmatrix}
1 & \int_{\mathbb{R}} \frac{\exp (-\mathscr{N} \, V(s))}{s-z} \, \frac{\md
s}{2 \pi \mi} \\
0 & 1
\end{pmatrix}, \quad z \! \in \! \mathbb{C} \setminus \mathbb{R},
\end{equation*}
where $\bm{\pi}_{0}(z) \! := \! \overset{e}{\mathrm{Y}}_{11}(z) \! \equiv \!
1$, with $\overset{e}{\mathrm{Y}}_{11}(z)$ the $(1 \, 1)$-element of
$\overset{e}{\mathrm{Y}}(z);$ and {\rm (ii)} for $n \! \in \! \mathbb{N}$,
\begin{equation}
\overset{e}{\mathrm{Y}}(z) \! = \!
\begin{pmatrix}
\boldsymbol{\pi}_{2n}(z) & \int_{\mathbb{R}} \frac{\boldsymbol{\pi}_{2n}(s)
\exp (-\mathscr{N} \, V(s))}{s-z} \, \frac{\md s}{2 \pi \mi} \\
\overset{e}{\mathrm{Y}}_{21}(z) & \int_{\mathbb{R}} \frac{\overset{e}{\mathrm{
Y}}_{21}(s) \exp (-\mathscr{N} \, V(s))}{s-z} \, \frac{\md s}{2 \pi \mi}
\end{pmatrix}, \quad z \! \in \! \mathbb{C} \setminus \mathbb{R},
\end{equation}
where $\overset{e}{\mathrm{Y}}_{21} \colon \mathbb{C}^{\ast} \! \to \!
\mathbb{C}$ denotes the $(2 \, 1)$-element of $\overset{e}{\mathrm{Y}}(z)$,
and $\boldsymbol{\pi}_{2n}(z)$ is the even degree monic {\rm OLP} defined in
Equation~{\rm (1.4)}.
\end{cccc}

\emph{Proof.} Set $\widetilde{w}(z) \! := \! \exp (-\mathscr{N} \, V(z))$,
$\mathscr{N} \! \in \! \mathbb{N}$, where $V \colon \mathbb{R} \setminus
\lbrace 0 \rbrace \! \to \! \mathbb{R}$ satisfies conditions~(V1)--(V3).
Since $\int_{\mathbb{R}}s^{j} \widetilde{w}(s) \, \md s \! < \! \infty$, $j \!
\in \! \mathbb{Z}$, it follows via a straightforward application of the
Sokhotski-Plemelj formula that, for $n \! = \! 0$, \textbf{RHP1} has the
(unique) upper-triangular solution
\begin{equation*}
\overset{e}{\mathrm{Y}}(z) \! = \!
\begin{pmatrix}
1 & \int_{\mathbb{R}} \frac{\widetilde{w}(s)}{s-z} \, \frac{\md s}{2 \pi
\mi} \\
0 & 1
\end{pmatrix}, \quad z \! \in \! \mathbb{C} \setminus \mathbb{R},
\end{equation*}
where $\bm{\pi}_{0}(z) \! := \! \overset{e}{\mathrm{Y}}_{11}(z) \! \equiv \!
1$. Hereafter, $n \! \in \! \mathbb{N}$ is considered.

If $\overset{e}{\mathrm{Y}} \colon \mathbb{C} \setminus \mathbb{R} \! \to \! 
\operatorname{SL}_{2}(\mathbb{C})$ solves \textbf{RHP1}, then it follows 
{}from the jump condition~(ii) of \textbf{RHP1} that, for the elements of the 
first column of $\overset{e}{\mathrm{Y}}(z)$,
\begin{equation*}
\left(\overset{e}{\mathrm{Y}}_{j1}(z) \right)_{+} \! = \! \left(\overset{e}{
\mathrm{Y}}_{j1}(z) \right)_{-} \! := \! \overset{e}{\mathrm{Y}}_{j1}(z),
\quad j \! = \! 1,2,
\end{equation*}
and, for the elements of the second row,
\begin{equation*}
\left(\overset{e}{\mathrm{Y}}_{j2}(z) \right)_{+} \! - \! \left(\overset{e}{
\mathrm{Y}}_{j2}(z) \right)_{-} \! = \! \overset{e}{\mathrm{Y}}_{j1}(z)
\widetilde{w}(z), \quad j \! = \! 1,2.
\end{equation*}
{}From condition~(i), the normalisation condition~(iii), and the boundedness
condition~(iv) of \textbf{RHP1}, in particular, $\overset{e}{\mathrm{Y}}_{11}
(z)z^{-n} \! =_{\underset{z \in \mathbb{C} \setminus \mathbb{R}}{z \to \infty}
} \! 1 \! + \! \mathcal{O}(z^{-1})$, $\overset{e}{\mathrm{Y}}_{11}(z)z^{n} \!
=_{\underset{z \in \mathbb{C} \setminus \mathbb{R}}{z \to 0}} \! \mathcal{O}
(1)$, $\overset{e}{\mathrm{Y}}_{21}(z)z^{-n} \! =_{\underset{z \in \mathbb{C}
\setminus \mathbb{R}}{z \to \infty}} \! \mathcal{O}(z^{-1})$, and $\overset{
e}{\mathrm{Y}}_{21}(z)z^{n} \! =_{\underset{z \in \mathbb{C} \setminus
\mathbb{R}}{z \to 0}} \! \mathcal{O}(1)$, and the fact that $\overset{e}{
\mathrm{Y}}_{11}(z)$ and $\overset{e}{\mathrm{Y}}_{21}(z)$ have no jumps
throughout the $z$-plane, it follows that $\overset{e}{\mathrm{Y}}_{11}(z)$
is a monic rational function with a pole at the origin and at the point at
infinity, with representation $\overset{e}{\mathrm{Y}}_{11}(z) \! = \! \sum_{
l=-n}^{n} \nu_{l}z^{l}$, where $\nu_{n} \! = \! 1$, and $\overset{e}{\mathrm{
Y}}_{21}(z)$ is a rational function with a pole at the origin and at the point
at infinity, with representation $\overset{e}{\mathrm{Y}}_{21}(z) \! = \!
\sum_{l=-n}^{n-1} \nu_{l}^{\sharp}z^{l}$. Application of the Sokhotski-Plemelj
formula to the jump relations for $\overset{e}{\mathrm{Y}}_{j2}(z)$, $j \! =
\! 1,2$, gives rise to the following Cauchy-type integral representations:
\begin{equation}
\overset{e}{\mathrm{Y}}_{j2}(z) \! = \! \int_{\mathbb{R}} \dfrac{\overset{e}{
\mathrm{Y}}_{j1}(s) \widetilde{w}(s)}{s \! - \! z} \dfrac{\md s}{2 \pi \mi},
\quad j \! = 1,2, \quad z \! \in \! \mathbb{C} \setminus \mathbb{R}. \tag{CA1}
\end{equation}
One now studies $\overset{e}{\mathrm{Y}}_{j1}(z)$, $j \! = \! 1,2$, in more
detail. {}From the normalisation condition~(iii) of \textbf{RHP1}, in
particular, $\overset{e}{\mathrm{Y}}_{12}(z)z^{n} \! =_{\underset{z \in
\mathbb{C} \setminus \mathbb{R}}{z \to \infty}} \! \mathcal{O}(z^{-1})$ and
$\overset{e}{\mathrm{Y}}_{22}(z)z^{n} \! =_{\underset{z \in \mathbb{C}
\setminus \mathbb{R}}{z \to \infty}} \! 1 \! + \! \mathcal{O}(z^{-1})$,
the formulae~(CA1), the fact that $\int_{\mathbb{R}}s^{j} \widetilde{w}(s) \,
\md s \! < \! \infty$, $j \! \in \! \mathbb{Z}$, and the expansion (for $\vert
s/z \vert \! \ll \! 1)$ $\tfrac{1}{s-z} \! = \! -\sum_{k=0}^{l} \tfrac{s^{k}
}{z^{k+1}} \! + \! \tfrac{s^{l+1}}{z^{l+1}(s-z)}$, $l \! \in \! \mathbb{Z}_{
0}^{+}$, it follows that
\begin{equation*}
\int_{\mathbb{R}} \overset{e}{\mathrm{Y}}_{11}(s)s^{k} \widetilde{w}(s) \, \md
s \! = \! 0, \quad k \! = \! 0,1,\dotsc,n \! - \! 1, \qquad \text{and} \qquad
\int_{\mathbb{R}} \overset{e}{\mathrm{Y}}_{11}(s)s^{n} \widetilde{w}(s) \, \md
s \! = \! -2 \pi \mi \mathfrak{p}^{e},
\end{equation*}
for some (pure imaginary) $\mathfrak{p}^{e}$ of the form $\mathfrak{p}^{e} \!
= \! \mi \mathfrak{q}^{e}$, with $\mathfrak{q}^{e} \! > \! 0$ (see below), and
\begin{equation*}
\int_{\mathbb{R}} \overset{e}{\mathrm{Y}}_{21}(s)s^{j} \widetilde{w}(s) \, \md
s \! = \! 0, \quad j \! = \! 0,1,\dotsc,n \! - \! 2, \qquad \text{and} \qquad
\int_{\mathbb{R}} \overset{e}{\mathrm{Y}}_{21}(s)s^{n-1} \widetilde{w}(s) \,
\md s \! = \! -2 \pi \mi;
\end{equation*}
and, {}from the boundedness condition~(iv) of \textbf{RHP1}, in particular,
$\overset{e}{\mathrm{Y}}_{12}(z)z^{-n} \! =_{\underset{z \in \mathbb{C}
\setminus \mathbb{R}}{z \to 0}} \! \mathcal{O}(1)$ and $\overset{e}{\mathrm{
Y}}_{22}(z) \linebreak[4]
\cdot z^{-n} \! =_{\underset{z \in \mathbb{C} \setminus \mathbb{R}}{z \to
0}} \! \mathcal{O}(1)$, the formulae~(CA1), the fact that $\int_{\mathbb{R}}
s^{j} \widetilde{w}(s) \, \md s \! < \! \infty$, $j \! \in \! \mathbb{Z}$, and
the expansion (for $\vert z/s \vert \! \ll \! 1)$ $\tfrac{1}{z-s} \! = \!
-\sum_{k=0}^{l} \tfrac{z^{k}}{s^{k+1}} \! + \! \tfrac{z^{l+1}}{s^{l+1}(z-s)}$,
$l \! \in \! \mathbb{Z}_{0}^{+}$, it follows that
\begin{equation*}
\int_{\mathbb{R}} \overset{e}{\mathrm{Y}}_{11}(s)s^{-k} \widetilde{w}(s) \,
\md s \! = \! 0, \quad k \! = \! 1,2,\dotsc,n, \qquad \text{and} \qquad \int_{
\mathbb{R}} \overset{e}{\mathrm{Y}}_{21}(s)s^{-j} \widetilde{w}(s) \, \md s \!
= \! 0, \quad j \! = \! 1,2,\dotsc,n:
\end{equation*}
these give rise to $2n \! + \! 1$ conditions for $\overset{e}{\mathrm{Y}}_{11}
(z)$, and $2n$ conditions for $\overset{e}{\mathrm{Y}}_{21}(z)$. Consider,
first, the $2n$ conditions for $\overset{e}{\mathrm{Y}}_{21}(z)$. Recalling
that the strong moments are defined by $c_{j} \! := \! \int_{\mathbb{R}}s^{j}
\widetilde{w}(s) \, \md s$, $j \! \in \! \mathbb{Z}$, it follows {}from the
representation (established above) $\overset{e}{\mathrm{Y}}_{21}(z) \! = \!
\sum_{l=-n}^{n-1} \nu_{l}^{\sharp}z^{l}$ and the $2n$ conditions for
$\overset{e}{\mathrm{Y}}_{21}(z)$ that
\begin{equation*}
\sum_{l=-n}^{n-1} \nu_{l}^{\sharp}c_{l+k} \! = \! 0, \quad k \! = \! -n,-(n \!
- \! 1),\dotsc,n \! - \! 2, \qquad \quad \text{and} \qquad \quad \sum_{l=-n}^{
n-1} \nu_{l}^{\sharp}c_{l+n-1} \! = \! -2 \pi \mi,
\end{equation*}
that is,
\begin{equation*}
\begin{pmatrix}
c_{-2n} & c_{-2n+1} & \dotsb & c_{-2} & c_{-1} \\
c_{-2n+1} & c_{-2n+2} & \dotsb & c_{-1} & c_{0} \\
\vdots & \vdots & \ddots & \vdots & \vdots \\
c_{-2} & c_{-1} & \dotsb & c_{2n-4} & c_{2n-3} \\
c_{-1} & c_{0} & \dotsb & c_{2n-3} & c_{2n-2}
\end{pmatrix} \!
\begin{pmatrix}
\nu_{-n}^{\sharp} \\
\nu_{-n+1}^{\sharp} \\
\vdots \\
\nu_{n-2}^{\sharp} \\
\nu_{n-1}^{\sharp}
\end{pmatrix} \! = \!
\begin{pmatrix}
0 \\
0 \\
\vdots \\
0 \\
-2 \pi \mi
\end{pmatrix}.
\end{equation*}
This linear system of $2n$ equations for the $2n$ unknowns $\nu_{l}^{\sharp}$,
$l \! = \! -n,-(n \! - \! 1),\dotsc,n \! - \! 1$, admits a unique solution if,
and only if, the determinant of the coefficient matrix, in this case $H^{(-2n)
}_{2n}$ (cf.~Equations~(1.1)), is non-zero; in fact, it will be shown that
$H^{(-2n)}_{2n} \! > \! 0$. An integral representation for the Hankel
determinants $H^{(m)}_{k}$, $(m,k) \! \in \! \mathbb{Z} \times \mathbb{N}$,
is now obtained; then the substitutions $m \! = \! -2n$ and $k \! = \! 2n$ are
made. In the calculations that follow, $\mathfrak{S}_{k}$ denotes the $k!$
permutations $\boldsymbol{\sigma}$ of $\lbrace 1,2,\dotsc,k \rbrace$.
Recalling that $c_{j} \! := \! \int_{\mathbb{R}}s^{j} \, \md \widetilde{\mu}
(s)$, $j \! \in \! \mathbb{Z}$, where $\md \widetilde{\mu}(z) \! = \!
\widetilde{w}(z) \, \md z \! = \! \exp (-\mathscr{N} \, V(z)) \, \md z$, and
using the multi-linearity property of the determinant, via Equations~(1.1),
one proceeds thus (recall that $H^{(m)}_{0} \! := \! 1)$:
\begin{align*}
H^{(m)}_{k} :=& \,
\begin{vmatrix}
c_{m} & c_{m+1} & \dotsb & c_{m+k-1} \\
c_{m+1} & c_{m+2} & \dotsb & c_{m+k} \\
\vdots & \vdots & \ddots & \vdots \\
c_{m+k-2} & c_{m+k-1} & \dotsb & c_{m+2k-3} \\
c_{m+k-1} & c_{m+k} & \dotsb & c_{m+2k-2}
\end{vmatrix} \\
=& \,
\begin{vmatrix}
\int_{\mathbb{R}}s_{1}^{m} \, \md \widetilde{\mu}(s_{1}) & \int_{\mathbb{R}}
s_{2}^{m+1} \, \md \widetilde{\mu}(s_{2}) & \dotsb & \int_{\mathbb{R}}s_{k}^{
m+k-1} \, \md \widetilde{\mu}(s_{k}) \\
\int_{\mathbb{R}}s_{1}^{m+1} \, \md \widetilde{\mu}(s_{1}) & \int_{\mathbb{R}}
s_{2}^{m+2} \, \md \widetilde{\mu}(s_{2}) & \dotsb & \int_{\mathbb{R}}s_{k}^{
m+k} \, \md \widetilde{\mu}(s_{k}) \\
\vdots & \vdots & \ddots & \vdots \\
\int_{\mathbb{R}}s_{1}^{m+k-2} \, \md \widetilde{\mu}(s_{1}) & \int_{\mathbb{
R}}s_{2}^{m+k-1} \, \md \widetilde{\mu}(s_{2}) & \dotsb & \int_{\mathbb{R}}s_{
k}^{m+2k-3} \, \md \widetilde{\mu}(s_{k}) \\
\int_{\mathbb{R}}s_{1}^{m+k-1} \, \md \widetilde{\mu}(s_{1}) & \int_{\mathbb{
R}}s_{2}^{m+k} \, \md \widetilde{\mu}(s_{2}) & \dotsb & \int_{\mathbb{R}}s_{
k}^{m+2k-2} \, \md \widetilde{\mu}(s_{k})
\end{vmatrix} \\
=& \, \underbrace{\int_{\mathbb{R}} \int_{\mathbb{R}} \dotsi \int_{\mathbb{R}}
}_{k} \, \, \md \widetilde{\mu}(s_{1}) \, \md \widetilde{\mu}(s_{2}) \, \dotsb
\, \md \widetilde{\mu}(s_{k})
\begin{vmatrix}
s_{1}^{m} & s_{2}^{m+1} & \dotsb & s_{k}^{m+k-1} \\
s_{1}^{m+1} & s_{2}^{m+2} & \dotsb & s_{k}^{m+k} \\
\vdots & \vdots & \ddots & \vdots \\
s_{1}^{m+k-2} & s_{2}^{m+k-1} & \dotsb & s_{k}^{m+2k-3} \\
s_{1}^{m+k-1} & s_{2}^{m+k} & \dotsb & s_{k}^{m+2k-2}
\end{vmatrix} \\
=& \, \underbrace{\int_{\mathbb{R}} \int_{\mathbb{R}} \dotsi \int_{\mathbb{R}}
}_{k} \, \, \md \widetilde{\mu}(s_{1}) \, \md \widetilde{\mu}(s_{2}) \, \dotsb
\, \md \widetilde{\mu}(s_{k}) \, s_{1}^{m}s_{2}^{m+1} \dotsb s_{k}^{m+k-1}
\underbrace{\begin{vmatrix}
1 & 1 & \dotsb & 1 \\
s_{1} & s_{2} & \dotsb & s_{k} \\
\vdots & \vdots & \ddots & \vdots \\
s_{1}^{k-2} & s_{2}^{k-2} & \dotsb & s_{k}^{k-2} \\
s_{1}^{k-1} & s_{2}^{k-1} & \dotsb & s_{k}^{k-1}
\end{vmatrix}}_{=: \, \mathrm{V}(s_{1},s_{2},\dotsc,s_{k})} \\
=& \, \dfrac{1}{k!} \sum_{\boldsymbol{\sigma} \in \mathfrak{S}_{k}} \,
\underbrace{\int_{\mathbb{R}} \int_{\mathbb{R}} \dotsi \int_{\mathbb{R}}}_{k}
\, \, \md \widetilde{\mu}(s_{\sigma (1)}) \, \md \widetilde{\mu}(s_{\sigma
(2)}) \, \dotsb \, \md \widetilde{\mu}(s_{\sigma (k)}) \prod_{j=1}^{k}s_{
\sigma (j)}^{m}s_{\sigma (j)}^{j-1} \, \mathrm{V} \! \left(s_{\sigma (1)},
s_{\sigma (2)},\dotsc,s_{\sigma (k)} \right) \\
=& \, \dfrac{1}{k!} \sum_{\boldsymbol{\sigma} \in \mathfrak{S}_{k}} \,
\underbrace{\int_{\mathbb{R}} \int_{\mathbb{R}} \dotsi \int_{\mathbb{R}}}_{k}
\, \, \md \widetilde{\mu}(s_{1}) \, \md \widetilde{\mu}(s_{2}) \, \dotsb \,
\md \widetilde{\mu}(s_{k}) \, s_{1}^{m}s_{2}^{m} \dotsb s_{k}^{m} \! \left(
\operatorname{sgn}(\boldsymbol{\sigma}) \prod_{j=1}^{k}s_{\sigma (j)}^{j-1}
\right) \\
\times& \, \mathrm{V} \! \left(s_{1},s_{2},\dotsc,s_{k} \right) \\
=& \, \dfrac{1}{k!} \underbrace{\int_{\mathbb{R}} \int_{\mathbb{R}} \dotsi
\int_{\mathbb{R}}}_{k} \, \, \md \widetilde{\mu}(s_{1}) \, \md \widetilde{\mu}
(s_{2}) \, \dotsb \, \md \widetilde{\mu}(s_{k}) \, s_{1}^{m}s_{2}^{m} \dotsb
s_{k}^{m} \, \mathrm{V} \! \left(s_{1},s_{2},\dotsc,s_{k} \right) \\
\times& \, \underbrace{\sum_{\boldsymbol{\sigma} \in \mathfrak{S}_{k}}
\operatorname{sgn}(\boldsymbol{\sigma})s_{\sigma (1)}^{0}s_{\sigma (2)}^{1}
\dotsb s_{\sigma (k)}^{k-1}}_{= \, \mathrm{V}(s_{1},s_{2},\dotsc,s_{k})} 
\quad \Rightarrow \\
H^{(m)}_{k} =& \, \dfrac{1}{k!} \underbrace{\int_{\mathbb{R}} \int_{\mathbb{
R}} \dotsi \int_{\mathbb{R}}}_{k} \, \, \md \widetilde{\mu}(s_{1}) \, \md
\widetilde{\mu}(s_{2}) \, \dotsb \, \md \widetilde{\mu}(s_{k}) \, s_{1}^{m}
s_{2}^{m} \dotsb s_{k}^{m} \! \left(\mathrm{V}(s_{1},s_{2},\dotsc,s_{k})
\right)^{2};
\end{align*}
using the well-known determinantal formula $\mathrm{V}(s_{1},s_{2},\dotsc,
s_{k}) \! = \! \prod_{\underset{j<i}{i,j=1}}^{k}(s_{i} \! - \! s_{j})$, one
arrives at
\begin{equation}
H^{(m)}_{k} \! = \! \dfrac{1}{k!} \underbrace{\int_{\mathbb{R}} \int_{\mathbb{
R}} \dotsi \int_{\mathbb{R}}}_{k} \, \, s_{1}^{m}s_{2}^{m} \dotsb s_{k}^{m}
\prod_{\underset{l<i}{i,l=1}}^{k}(s_{i} \! - \! s_{l})^{2} \, \md \widetilde{
\mu}(s_{1}) \, \md \widetilde{\mu}(s_{2}) \, \dotsb \, \md \widetilde{\mu}(s_{
k}), \quad (m,k) \! \in \! \mathbb{Z} \times \mathbb{N} \tag{HA1}.
\end{equation}
Letting $m \! = \! -2n$ and $k \! = \! 2n$, it follows from the formula~(HA1)
that
\begin{equation*}
H^{(-2n)}_{2n} \! = \! \dfrac{1}{(2n)!} \underbrace{\int_{\mathbb{R}} \int_{
\mathbb{R}} \dotsi \int_{\mathbb{R}}}_{2n} \, \, s_{1}^{-2n}s_{2}^{-2n} \dotsb
s_{2n}^{-2n} \prod_{\underset{l<i}{i,l=1}}^{2n}(s_{i} \! - \! s_{l})^{2} \,
\md \widetilde{\mu}(s_{1}) \, \md \widetilde{\mu}(s_{2}) \, \dotsb \, \md
\widetilde{\mu}(s_{2n}) \! > \! 0,
\end{equation*}
whence the existence (and uniqueness) of $\overset{e}{\mathrm{Y}}_{21}(z)$.

Similarly, it follows, {}from the representation (established above)
$\overset{e}{\mathrm{Y}}_{11}(z) \! = \! \sum_{l=-n}^{n} \nu_{l}z^{l}$, with
$\nu_{n} \! = \! 1$, and the $2n \! + \! 1$ conditions for $\overset{e}{
\mathrm{Y}}_{11}(z)$, that
\begin{equation*}
\sum_{l=-n}^{n} \nu_{l}c_{l+k} \! = \! 0, \quad k \! = \! -n,-(n \! - \! 1),
\dotsc,n \! - \! 1, \qquad \quad \text{and} \qquad \quad \sum_{l=-n}^{n} 
\nu_{l}c_{l+n} \! = \! -2 \pi \mi \mathfrak{p}^{e},
\end{equation*}
that is,
\begin{equation*}
\begin{pmatrix}
c_{-2n} & \dotsb & c_{-n}  & \dotsb  & c_{-1}     & 0 \\
\vdots  & \ddots & \vdots  & \ddots  & \vdots     & \vdots \\
c_{-n}  & \dotsb & c_{0}   & \dotsb  & c_{n-1}    & 0 \\
\vdots  & \ddots & \vdots  & \ddots  & \vdots     & \vdots \\
c_{-1}  & \dotsb & c_{n-1} & \dotsb  & c_{2(n-1)} & 0 \\
c_{0}   & \dotsb & c_{n}   &\dotsb   & c_{2n-1}   & 2 \pi \mi
\end{pmatrix} \!
\begin{pmatrix}
\nu_{-n} \\
\vdots \\
\nu_{0} \\
\vdots \\
\nu_{n-1} \\
\mathfrak{p}^{e}
\end{pmatrix} \! = \!
\begin{pmatrix}
-c_{0} \\
\vdots \\
-c_{n} \\
\vdots \\
-c_{2n-1} \\
-c_{2n}
\end{pmatrix}.
\end{equation*}
This linear system of $2n \! + \! 1$ equations for the $2n \! + \! 1$ unknowns
$\nu_{l}$, $l \! = \! -n,-(n \! - \! 1),\dotsb,n \! - \! 1$, and $\mathfrak{
p}^{e}$ admits a unique solution if, and only if, the determinant of
the coefficient matrix, in this case $2 \pi \mi H^{(-2n)}_{2n}$, is non-zero;
but, it was shown above that $H^{(-2n)}_{2n} \! > \! 0$. Furthermore, via
Cramer's Rule:
\begin{equation*}
\mathfrak{p}^{e} \! = \! \dfrac{
\begin{vmatrix}
c_{-2n} & \dotsb & c_{-n}  & \dotsb  & c_{-1}     & -c_{0} \\
\vdots  & \ddots & \vdots  & \ddots  & \vdots     & \vdots \\
c_{-n}  & \dotsb & c_{0}   & \dotsb  & c_{n-1}    & -c_{n} \\
\vdots  & \ddots & \vdots  & \ddots  & \vdots     & \vdots \\
c_{-1}  & \dotsb & c_{n-1} & \dotsb  & c_{2(n-1)} & -c_{2n-1} \\
c_{0}   & \dotsb & c_{n}   &\dotsb   & c_{2n-1}   & -c_{2n}
\end{vmatrix}}{2 \pi \mi H^{(-2n)}_{2n}} \! = \! -\dfrac{1}{2 \pi \mi}
\dfrac{H^{(-2n)}_{2n+1}}{H^{(-2n)}_{2n}}.
\end{equation*}
Using the Hankel determinant formula~(HA1) with the substitutions $m \! = \!
-2n$ and $k \! = \! 2n \! + \! 1$, one arrives at
\begin{equation*}
H^{(-2n)}_{2n+1} \! = \! \dfrac{1}{(2n \! + \! 1)!} \underbrace{\int_{\mathbb{
R}} \int_{\mathbb{R}} \dotsi \int_{\mathbb{R}}}_{2n+1} \, \, s_{1}^{-2n}
s_{2}^{-2n} \dotsb s_{2n+1}^{-2n} \prod_{\underset{l<i}{i,l=1}}^{2n+1}(s_{i}
\! - \! s_{l})^{2} \, \md \widetilde{\mu}(s_{1}) \, \md \widetilde{\mu}(s_{2})
\, \dotsb \, \md \widetilde{\mu}(s_{2n+1}) \! > \! 0;
\end{equation*}
hence, $H^{(-2n)}_{2n+1}/H^{(-2n)}_{2n} \! > \! 0$. Using, now, the fact that
$\int_{\mathbb{R}} \overset{e}{\mathrm{Y}}_{11}(s)s^{k} \widetilde{w}(s) \,
\md s \! = \! 0$, $k \! = \! -n,-(n \! - \! 1),\dotsc,n \! - \! 1$, and the
relation $\int_{\mathbb{R}} \overset{e}{\mathrm{Y}}_{11}(s)s^{n} \widetilde{w}
(s) \, \md s \! = \! -2 \pi \mi \mathfrak{p}^{e}$, one notes, via the above
formula for $\mathfrak{p}^{e}$, that
\begin{align*}
\int_{\mathbb{R}} \overset{e}{\mathrm{Y}}_{11}(s)s^{n} \widetilde{w}(s) \,
\md s =& \, \int_{\mathbb{R}} \overset{e}{\mathrm{Y}}_{11}(s) \underbrace{
\left(s^{n} \! + \! \nu_{n-1}s^{n-1} \! + \! \dotsb \! + \! \nu_{-n}s^{-n}
\right)}_{= \, \overset{e}{\mathrm{Y}}_{11}(s)} \, \widetilde{w}(s) \, \md 
s \! = \! \int_{\mathbb{R}} \overset{e}{\mathrm{Y}}_{11}(s) \overset{e}{
\mathrm{Y}}_{11}(s) \widetilde{w}(s) \, \md s \\
=& \, -2 \pi \mi \mathfrak{p}^{e} \! = \! H^{(-2n)}_{2n+1}/H^{(-2n)}_{2n}
\quad (> \! 0);
\end{align*}
but the right-hand side of the latter expression (cf.~Equations~(1.8)) is
equal to $(\xi^{(2n)}_{n})^{-2} \! = \! \lvert \lvert \overset{e}{\mathrm{Y}
}_{11}(\pmb{\cdot}) \lvert \lvert_{\mathscr{L}}^{2}$ $(> \! 0)$: the existence
and uniqueness of $\overset{e}{\mathrm{Y}}_{11}(z) \! =: \! \boldsymbol{\pi}_{
2n}(z)$, the even degree monic OLP with respect to the inner product $\langle
\pmb{\cdot},\pmb{\cdot} \rangle_{\mathscr{L}}$, is thus established. \hfill
$\qed$
\begin{cccc}
Let $\overset{o}{\mathrm{Y}} \colon \mathbb{C} \setminus \mathbb{R} \! \to
\! \mathrm{SL}_{2}(\mathbb{C})$ solve {\rm \pmb{RHP2}}. {\rm \pmb{RHP2}}
possesses a unique solution given by: {\rm (i)} for $n \! = \! 0$,
\begin{equation*}
\overset{o}{\mathrm{Y}}(z) \! = \!
\begin{pmatrix}
z \bm{\pi}_{1}(z) & z \int_{\mathbb{R}} \frac{(s \bm{\pi}_{1}(s)) \exp (-
\mathscr{N} \, V(s))}{s(s-z)} \, \frac{\md s}{2 \pi \mi} \\
2 \pi \mi z & 1 \! + \! z \int_{\mathbb{R}} \frac{\exp (-\mathscr{N} \, V(s))
}{s-z} \, \md s
\end{pmatrix}, \quad z \! \in \! \mathbb{C} \setminus \mathbb{R},
\end{equation*}
where $\bm{\pi}_{1}(z) \! = \! \tfrac{1}{z} \! + \! \tfrac{\xi^{(1)}_{0}}{
\xi^{(1)}_{-1}}$, with $\tfrac{\xi^{(1)}_{0}}{\xi^{(1)}_{-1}} \! = \! -\int_{
\mathbb{R}}s^{-1} \exp (-\mathscr{N} \, V(s)) \, \md s$, $\mathscr{N} \! \in
\! \mathbb{N};$ and {\rm (ii)} for $n \! \in \! \mathbb{N}$,
\begin{equation*}
\overset{o}{\mathrm{Y}}(z) \! = \!
\begin{pmatrix}
z \boldsymbol{\pi}_{2n+1}(z) & z \int_{\mathbb{R}} \frac{(s \boldsymbol{\pi}_{
2n+1}(s)) \exp (-\mathscr{N} \, V(s))}{s(s-z)} \, \frac{\md s}{2 \pi \mi} \\
\overset{o}{\mathrm{Y}}_{21}(z) & z \int_{\mathbb{R}} \frac{\overset{o}{
\mathrm{Y}}_{21}(s) \exp (-\mathscr{N} \, V(s))}{s(s-z)} \, \frac{\md s}{2
\pi \mi}
\end{pmatrix}, \quad z \! \in \! \mathbb{C} \setminus \mathbb{R},
\end{equation*}
where $\overset{o}{\mathrm{Y}}_{21} \colon \mathbb{C}^{\ast} \! \to \!
\mathbb{C}$ denotes the $(2 \, 1)$-element of $\overset{o}{\mathrm{Y}}(z)$,
and $\boldsymbol{\pi}_{2n+1}(z)$ is the odd degree monic {\rm OLP} defined in
Equation~{\rm (1.5)}.
\end{cccc}

\emph{Proof.} See \cite{a51}, the proof of Lemma~2.2.2. \hfill $\qed$
\begin{ffff}
Let $V \colon \mathbb{R} \setminus \{0\} \! \to \! \mathbb{R}$ satisfy
conditions~{\rm (V1)--(V3)}. Let $\boldsymbol{\pi}_{2n}(z)$ and $\boldsymbol{
\pi}_{2n+1}(z)$ be the even degree and odd degree monic {\rm OLPs} with
respect to the inner product $\langle \pmb{\cdot},\pmb{\cdot} \rangle_{
\mathscr{L}}$ defined, respectively, in Equations~{\rm (1.4)} and~{\rm (1.5)},
and let $\xi^{(2n)}_{n}$ and $\xi^{(2n+1)}_{-n-1}$ be the corresponding `even'
and `odd' norming constants, respectively. Then, $\xi^{(2n)}_{n}$ and $\xi^{
(2n+1)}_{-n-1}$ have the following representations:
\begin{equation*}
\dfrac{\xi^{(2n)}_{n}}{\sqrt{\smash[b]{2n \! + \! 1}}} = \sqrt{\dfrac{
\underbrace{\int_{\mathbb{R}} \int_{\mathbb{R}} \dotsi \int_{\mathbb{R}}}_{2n}
\, \, s_{1}^{-2n}s_{2}^{-2n} \dotsb s_{2n}^{-2n} \mathlarger{\prod_{\underset{
l<i}{i,l=1}}^{2n}}(s_{i} \! - \! s_{l})^{2} \, \md \widetilde{\mu}(s_{1}) \,
\md \widetilde{\mu}(s_{2}) \, \dotsb \, \md \widetilde{\mu}(s_{2n})}{
\underbrace{\int_{\mathbb{R}} \int_{\mathbb{R}} \dotsi \int_{\mathbb{R}}}_{2
n+1} \, \, \lambda_{1}^{-2n} \lambda_{2}^{-2n} \dotsb \lambda_{2n+1}^{-2n}
\mathlarger{\prod_{\underset{l<i}{i,l=1}}^{2n+1}}(\lambda_{i} \! - \!
\lambda_{l})^{2} \, \md \widetilde{\mu}(\lambda_{1}) \, \md \widetilde{\mu}
(\lambda_{2}) \, \dotsb \, \md \widetilde{\mu}(\lambda_{2n+1})}},
\end{equation*}
\begin{equation*}
\dfrac{\xi^{(2n+1)}_{-n-1}}{\sqrt{\smash[b]{2(n \! + \! 1)}}} = \sqrt{
\dfrac{\underbrace{\int_{\mathbb{R}} \int_{\mathbb{R}} \dotsi \int_{\mathbb{
R}}}_{2n+1} \, \, \varpi_{1}^{-2n} \varpi_{2}^{-2n} \dotsb \varpi_{2n+1}^{-2n}
\mathlarger{\prod_{\underset{l<i}{i,l=1}}^{2n+1}}(\varpi_{i} \! - \! \varpi_{
l})^{2} \, \md \widetilde{\mu}(\varpi_{1}) \, \md \widetilde{\mu}(\varpi_{2})
\, \dotsb \, \md \widetilde{\mu}(\varpi_{2n+1})}{\underbrace{\int_{\mathbb{R}}
\int_{\mathbb{R}} \dotsi \int_{\mathbb{R}}}_{2n+2} \, \, \varsigma_{1}^{-2n-2}
\varsigma_{2}^{-2n-2} \dotsb \varsigma_{2n+2}^{-2n-2} \mathlarger{\prod_{
\underset{l<i}{i,l=1}}^{2n+2}}(\varsigma_{i} \! - \! \varsigma_{l})^{2} \, \md
\widetilde{\mu}(\varsigma_{1}) \, \md \widetilde{\mu}(\varsigma_{2}) \, \dotsb
\, \md \widetilde{\mu}(\varsigma_{2n+2})}},
\end{equation*}
where $\md \widetilde{\mu}(z) \! := \! \exp (-\mathscr{N} \, V(z)) \, \md z$,
$\mathscr{N} \! \in \! \mathbb{N}$.
\end{ffff}

\emph{Proof.} Consider, without loss of generality, the representation for 
$\xi^{(2n)}_{n}$. Recall that (cf. Equations (1.8)) $(\xi^{(2n)}_{n})^{2} 
\! = \! H^{(-2n)}_{2n}/H^{(-2n)}_{2n+1}$ $(> \! 0)$: using the integral 
representations for $H^{(-2n)}_{2n}$ and $H^{(-2n)}_{2n+1}$ derived in (the 
course of) the proof of Lemma~2.2.1, and taking positive square roots of 
both sides of the resulting equality, one arrives at the representation for 
$\xi^{(2n)}_{n}$. See \cite{a51}, Corollary~2.2.1, for the proof of the 
representation for $\xi^{(2n+1)}_{-n-1}$. \hfill $\qed$
\begin{bbbb}
Let $V \colon \mathbb{R} \setminus \{0\} \! \to \! \mathbb{R}$ satisfy
conditions~{\rm (V1)--(V3)}. Let $\boldsymbol{\pi}_{2n}(z)$ and $\boldsymbol{
\pi}_{2n+1}(z)$ be the even degree and odd degree monic {\rm OLPs} with
respect to the inner product $\langle \pmb{\cdot},\pmb{\cdot} \rangle_{
\mathscr{L}}$ defined, respectively, in Equations~{\rm (1.4)}
and~{\rm (1.5)}. Then, $\boldsymbol{\pi}_{2n}(z)$ and $\boldsymbol{\pi}_{2n+1}
(z)$ have, respectively, the following integral representations:
\begin{align*}
\boldsymbol{\pi}_{2n}(z) =& \, \dfrac{z^{-n}}{(2n)!H^{(-2n)}_{2n}} \,
\underbrace{\int_{\mathbb{R}} \int_{\mathbb{R}} \dotsi \int_{\mathbb{R}}}_{2n}
\, \, s_{0}^{-2n}s_{1}^{-2n} \dotsb s_{2n-1}^{-2n} \, \prod_{\underset{l<i}{i,
l=0}}^{2n-1}(s_{i} \! - \! s_{l})^{2} \, \prod_{j=0}^{2n-1}(z \! - \! s_{j}) \\
\times& \, \md \widetilde{\mu}(s_{0}) \, \md \widetilde{\mu}(s_{1}) \, \dotsb
\, \md \widetilde{\mu}(s_{2n-1}),
\end{align*}
\begin{align*}
\boldsymbol{\pi}_{2n+1}(z) =& \, -\dfrac{z^{-n-1}}{(2n \! + \! 1)!H^{(-2n)}_{2
n+1}} \, \underbrace{\int_{\mathbb{R}} \int_{\mathbb{R}} \dotsi \int_{\mathbb{
R}}}_{2n+1} \, \, s_{0}^{-2n-1}s_{1}^{-2n-1} \dotsb s_{2n}^{-2n-1} \, \prod_{
\underset{l<i}{i,l=0}}^{2n}(s_{i} \! - \! s_{l})^{2} \, \prod_{j=0}^{2n}(z \!
- \! s_{j}) \\
\times& \, \md \widetilde{\mu}(s_{0}) \, \md \widetilde{\mu}(s_{1}) \, \dotsb
\, \md \widetilde{\mu}(s_{2n}),
\end{align*}
where
\begin{gather*}
H^{(-2n)}_{2n} \! = \! \dfrac{1}{(2n)!} \underbrace{\int_{\mathbb{R}} \int_{
\mathbb{R}} \dotsi \int_{\mathbb{R}}}_{2n} \, \, \lambda_{1}^{-2n} \lambda_{
2}^{-2n} \dotsb \lambda_{2n}^{-2n} \prod_{\underset{l<i}{i,l=1}}^{2n}
(\lambda_{i} \! - \! \lambda_{l})^{2} \, \md \widetilde{\mu}(\lambda_{1}) \,
\md \widetilde{\mu}(\lambda_{2}) \, \dotsb \, \md \widetilde{\mu}(\lambda_{2
n}), \\
H^{(-2n)}_{2n+1} \! = \! \dfrac{1}{(2n \! + \! 1)!} \underbrace{\int_{\mathbb{
R}} \int_{\mathbb{R}} \dotsi \int_{\mathbb{R}}}_{2n+1} \, \, \lambda_{1}^{-2n}
\lambda_{2}^{-2n} \dotsb \lambda_{2n+1}^{-2n} \prod_{\underset{l<i}{i,l=1}}^{2
n+1}(\lambda_{i} \! - \! \lambda_{l})^{2} \, \md \widetilde{\mu}(\lambda_{1})
\, \md \widetilde{\mu}(\lambda_{2}) \, \dotsb \, \md \widetilde{\mu}(\lambda_{
2n+1}),
\end{gather*}
with $\md \widetilde{\mu}(z) \! := \! \exp (-\mathscr{N} \, V(z)) \, \md z$,
$\mathscr{N} \! \in \! \mathbb{N}$.
\end{bbbb}

\emph{Proof.} Consider, without loss of generality, the integral
representation for the even degree monic OLP $\boldsymbol{\pi}_{2n}(z)$. Let
$\mathfrak{S}_{k}$ denote the $k!$ permutations $\boldsymbol{\sigma}$ of
$\lbrace 0,1,\dotsc,k \! - \! 1 \rbrace$. Recalling that $c_{j} \! := \!
\int_{\mathbb{R}}s^{j} \, \md \widetilde{\mu}(s)$, $j \! \in \! \mathbb{Z}$,
where $\md \widetilde{\mu}(z) \! = \! \widetilde{w}(z) \, \md z \! = \! \exp
(-\mathscr{N} \, V(z)) \, \md z$, $\mathscr{N} \! \in \! \mathbb{N}$, with
$V \colon \mathbb{R} \setminus \{0\} \! \to \! \mathbb{R}$ satisfying
conditons~(V1)--(V3), and using the multi-linearity property of the
determinant, via the determinantal representation for $\boldsymbol{\pi}_{2n}
(z)$ given in Equation~(1.6), one proceeds thus:
\begin{align*}
\boldsymbol{\pi}_{2n}(z) =& \, \dfrac{1}{H^{(-2n)}_{2n}}
\begin{vmatrix}
c_{-2n} & c_{-2n+1} & \dotsb & c_{-1}  & z^{-n} \\
c_{-2n+1} & c_{-2n+2} & \dotsb & c_{0} & z^{-n+1} \\
\vdots & \vdots & \ddots & \vdots & \vdots \\
c_{-1} & c_{0} & \dotsb & c_{2n-2} & z^{n-1} \\
c_{0} & c_{1} & \dotsb & c_{2n-1} & z^{n}
\end{vmatrix} \\
=& \, \dfrac{z^{-n}}{H^{(-2n)}_{2n}}
\begin{vmatrix}
c_{-2n} & c_{-2n+1} & \dotsb & c_{-1}  & c_{0} \\
c_{-2n+1} & c_{-2n+2} & \dotsb & c_{0} & c_{1} \\
\vdots & \vdots & \ddots & \vdots & \vdots \\
c_{-1} & c_{0} & \dotsb & c_{2n-2} & c_{2n-1} \\
z^{0} & z^{1} & \dotsb & z^{2n-1} & z^{2n}
\end{vmatrix} \\
=& \, \dfrac{z^{-n}}{H^{(-2n)}_{2n}}
\begin{vmatrix}
\int_{\mathbb{R}}s_{0}^{-2n} \, \md \widetilde{\mu}(s_{0}) & \int_{\mathbb{R}}
s_{0}^{-2n+1} \, \md \widetilde{\mu}(s_{0}) & \dotsb & \int_{\mathbb{R}}s_{
0}^{0} \, \md \widetilde{\mu}(s_{0}) \\
\int_{\mathbb{R}}s_{1}^{-2n+1} \, \md \widetilde{\mu}(s_{1}) & \int_{\mathbb{
R}}s_{1}^{-2n+2} \, \md \widetilde{\mu}(s_{1}) & \dotsb & \int_{\mathbb{R}}
s_{1}^{1} \, \md \widetilde{\mu}(s_{1}) \\
\vdots & \vdots & \ddots & \vdots \\
\int_{\mathbb{R}}s_{2n-1}^{-1} \, \md \widetilde{\mu}(s_{2n-1}) & \int_{
\mathbb{R}}s_{2n-1}^{0} \, \md \widetilde{\mu}(s_{2n-1}) & \dotsb & \int_{
\mathbb{R}}s_{2n-1}^{2n-1} \, \md \widetilde{\mu}(s_{2n-1}) \\
z^{0} & z^{1} & \dotsb & z^{2n}
\end{vmatrix} \\
=& \, \dfrac{z^{-n}}{H^{(-2n)}_{2n}} \, \underbrace{\int_{\mathbb{R}} \int_{
\mathbb{R}} \dotsi \int_{\mathbb{R}}}_{2n} \, \, \md \widetilde{\mu}(s_{0})
\md \widetilde{\mu}(s_{1}) \, \dotsb \, \md \widetilde{\mu}(s_{2n-1})
\begin{vmatrix}
s_{0}^{-2n} & s_{0}^{-2n+1} & \dotsb & s_{0}^{-1} & s_{0}^{0} \\
s_{1}^{-2n+1} & s_{1}^{-2n+2} & \dotsb & s_{1}^{0} & s_{1}^{1} \\
\vdots & \vdots & \ddots & \vdots & \vdots \\
s_{2n-1}^{-1} & s_{2n-1}^{0} & \dotsb & s_{2n-1}^{2n-2} & s_{2n-1}^{2n-1} \\
z^{0} & z^{1} & \dotsb & z^{2n-1} & z^{2n}
\end{vmatrix} \\
=& \, \dfrac{z^{-n}}{H^{(-2n)}_{2n}} \, \underbrace{\int_{\mathbb{R}} \int_{
\mathbb{R}} \dotsi \int_{\mathbb{R}}}_{2n} \, \, \md \widetilde{\mu}(s_{0})
\md \widetilde{\mu}(s_{1}) \, \dotsb \, \md \widetilde{\mu}(s_{2n-1})s_{0}^{-
2n}s_{1}^{-2n+1} \dotsb s_{2n-1}^{-1} \\
\times& \,
\begin{vmatrix}
s_{0}^{0} & s_{0}^{1} & \dotsb & s_{0}^{2n-1} & s_{0}^{2n} \\
s_{1}^{0} & s_{1}^{1} & \dotsb & s_{1}^{2n-1} & s_{1}^{2n} \\
\vdots & \vdots & \ddots & \vdots & \vdots \\
s_{2n-1}^{0} & s_{2n-1}^{1} & \dotsb & s_{2n-1}^{2n-1} & s_{2n-1}^{2n} \\
z^{0} & z^{1} & \dotsb & z^{2n-1} & z^{2n}
\end{vmatrix} \\
=& \, \dfrac{z^{-n}}{H^{(-2n)}_{2n}(2n)!} \, \sum_{\boldsymbol{\sigma} \in
\mathfrak{S}_{2n}} \, \underbrace{\int_{\mathbb{R}} \int_{\mathbb{R}} \dotsi
\int_{\mathbb{R}}}_{2n} \, \, \md \widetilde{\mu}(s_{\sigma (0)}) \md
\widetilde{\mu}(s_{\sigma (1)}) \, \dotsb \, \md \widetilde{\mu}(s_{\sigma
(2n-1)})s_{\sigma (0)}^{-2n}s_{\sigma (1)}^{-2n} \dotsb s_{\sigma (2n-1)}^{-
2n} \\
\times& \, s_{\sigma (0)}^{0}s_{\sigma (1)}^{1} \dotsb s_{\sigma (2n-1)}^{2n
-1}
\begin{vmatrix}
s_{\sigma (0)}^{0} & s_{\sigma (0)}^{1} & \dotsb & s_{\sigma (0)}^{2n-1} &
s_{\sigma (0)}^{2n} \\
s_{\sigma (1)}^{0} & s_{\sigma (1)}^{1} & \dotsb & s_{\sigma (1)}^{2n-1} &
s_{\sigma (1)}^{2n} \\
\vdots & \vdots & \ddots & \vdots & \vdots \\
s_{\sigma (2n-1)}^{0} & s_{\sigma (2n-1)}^{1} & \dotsb & s_{\sigma (2n-1)}^{2
n-1} & s_{\sigma (2n-1)}^{2n} \\
z^{0} & z^{1} & \dotsb & z^{2n-1} & z^{2n}
\end{vmatrix} \\
=& \, \dfrac{z^{-n}}{H^{(-2n)}_{2n}(2n)!} \, \underbrace{\int_{\mathbb{R}}
\int_{\mathbb{R}} \dotsi \int_{\mathbb{R}}}_{2n} \, \, \md \widetilde{\mu}
(s_{0}) \md \widetilde{\mu}(s_{1}) \, \dotsb \, \md \widetilde{\mu}(s_{2n-1})
s_{0}^{-2n}s_{1}^{-2n} \dotsb s_{2n-1}^{-2n} \\
\times& \, \left(\sum_{\boldsymbol{\sigma} \in \mathfrak{S}_{2n}}
\operatorname{sgn}(\boldsymbol{\sigma})s_{\sigma (0)}^{0}s_{\sigma (1)}^{1}
\dotsb s_{\sigma (2n-1)}^{2n-1} \right) \!
\begin{vmatrix}
s_{0}^{0} & s_{0}^{1} & \dotsb & s_{0}^{2n-1} & s_{0}^{2n} \\
s_{1}^{0} & s_{1}^{1} & \dotsb & s_{1}^{2n-1} & s_{1}^{2n} \\
\vdots & \vdots & \ddots & \vdots & \vdots \\
s_{2n-1}^{0} & s_{2n-1}^{1} & \dotsb & s_{2n-1}^{2n-1} & s_{2n-1}^{2n} \\
z^{0} & z^{1} & \dotsb & z^{2n-1} & z^{2n}
\end{vmatrix} \\
=& \, \dfrac{z^{-n}}{H^{(-2n)}_{2n}(2n)!} \, \underbrace{\int_{\mathbb{R}}
\int_{\mathbb{R}} \dotsi \int_{\mathbb{R}}}_{2n} \, \, \md \widetilde{\mu}
(s_{0}) \md \widetilde{\mu}(s_{1}) \, \dotsb \, \md \widetilde{\mu}(s_{2n-1})
s_{0}^{-2n}s_{1}^{-2n} \dotsb s_{2n-1}^{-2n} \\
\times& \,
\begin{vmatrix}
s_{0}^{0} & s_{0}^{1} & \dotsb & s_{0}^{2n-1} \\
s_{1}^{0} & s_{1}^{1} & \dotsb & s_{1}^{2n-1} \\
\vdots & \vdots & \ddots & \vdots \\
s_{2n-1}^{0} & s_{2n-1}^{1} & \dotsb & s_{2n-1}^{2n-1}
\end{vmatrix}
\begin{vmatrix}
s_{0}^{0} & s_{0}^{1} & \dotsb & s_{0}^{2n-1} & s_{0}^{2n} \\
s_{1}^{0} & s_{1}^{1} & \dotsb & s_{1}^{2n-1} & s_{1}^{2n} \\
\vdots & \vdots & \ddots & \vdots & \vdots \\
s_{2n-1}^{0} & s_{2n-1}^{1} & \dotsb & s_{2n-1}^{2n-1} & s_{2n-1}^{2n} \\
z^{0} & z^{1} & \dotsb & z^{2n-1} & z^{2n}
\end{vmatrix};
\end{align*}
but a straightforward calculation shows that
\begin{equation*}
\begin{vmatrix}
s_{0}^{0} & s_{0}^{1} & \dotsb & s_{0}^{2n-1} & s_{0}^{2n} \\
s_{1}^{0} & s_{1}^{1} & \dotsb & s_{1}^{2n-1} & s_{1}^{2n} \\
\vdots & \vdots & \ddots & \vdots & \vdots \\
s_{2n-1}^{0} & s_{2n-1}^{1} & \dotsb & s_{2n-1}^{2n-1} & s_{2n-1}^{2n} \\
z^{0} & z^{1} & \dotsb & z^{2n-1} & z^{2n}
\end{vmatrix} \! = \! 
\begin{vmatrix}
s_{0}^{0} & s_{0}^{1} & \dotsb & s_{0}^{2n-1} \\
s_{1}^{0} & s_{1}^{1} & \dotsb & s_{1}^{2n-1} \\
\vdots & \vdots & \ddots & \vdots \\
s_{2n-1}^{0} & s_{2n-1}^{1} & \dotsb & s_{2n-1}^{2n-1}
\end{vmatrix} 
\prod_{j=0}^{2n-1}(z \! - \! s_{j}),
\end{equation*}
whence
\begin{align*}
\boldsymbol{\pi}_{2n}(z) =& \, \dfrac{z^{-n}}{H^{(-2n)}_{2n}(2n)!} \,
\underbrace{\int_{\mathbb{R}} \int_{\mathbb{R}} \dotsi \int_{\mathbb{R}}}_{2n}
\, \, \md \widetilde{\mu}(s_{0}) \md \widetilde{\mu}(s_{1}) \, \dotsb \, \md
\widetilde{\mu}(s_{2n-1})s_{0}^{-2n}s_{1}^{-2n} \dotsb s_{2n-1}^{-2n} \\
\times& \, \prod_{j=0}^{2n-1}(z \! - \! s_{j})
\begin{vmatrix}
s_{0}^{0} & s_{0}^{1} & \dotsb & s_{0}^{2n-1} \\
s_{1}^{0} & s_{1}^{1} & \dotsb & s_{1}^{2n-1} \\
\vdots & \vdots & \ddots & \vdots \\
s_{2n-1}^{0} & s_{2n-1}^{1} & \dotsb & s_{2n-1}^{2n-1}
\end{vmatrix}^{2} \\
=& \, \dfrac{z^{-n}}{H^{(-2n)}_{2n}(2n)!} \, \underbrace{\int_{\mathbb{R}}
\int_{\mathbb{R}} \dotsi \int_{\mathbb{R}}}_{2n} \, \, \md \widetilde{\mu}
(s_{0}) \md \widetilde{\mu}(s_{1}) \, \dotsb \, \md \widetilde{\mu}(s_{2n-1})
s_{0}^{-2n}s_{1}^{-2n} \dotsb s_{2n-1}^{-2n} \\
\times& \, \prod_{j=0}^{2n-1}(z \! - \! s_{j})
\underbrace{\begin{vmatrix}
1 & 1 & \dotsb & 1 \\
s_{0}^{1} & s_{1}^{1} & \dotsb & s_{2n-1}^{1} \\
\vdots & \vdots & \ddots & \vdots \\
s_{0}^{2n-1} & s_{1}^{2n-1} & \dotsb & s_{2n-1}^{2n-1}
\end{vmatrix}^{2}}_{= \, \prod_{\underset{l<i}{i,l=0}}^{2n-1}(s_{i}-s_{l})^{
2}};
\end{align*}
hence the integral representation for $\boldsymbol{\pi}_{2n}(z)$ stated in the
Proposition, with the integral representation for $H^{(-2n)}_{2n}$ derived in
the proof of Lemma~2.2.1. See \cite{a51}, Proposition~2.2.1, for the proof of
the integral representation for the odd degree monic OLP $\boldsymbol{\pi}_{2
n+1}(z)$. \hfill $\qed$
\begin{eeee}
For the purposes of the ensuing asymptotic analysis, it is convenient to
re-write $\md \widetilde{\mu}(z) \! = \! \exp (-\mathscr{N} \, V(z)) \, \md z
\! = \! \exp (-n \widetilde{V}(z)) \, \md z \! =: \! \md \mu (z)$, $n \! \in
\! \mathbb{N}$, where
\begin{equation*}
\widetilde{V}(z) \! = \! z_{o}V(z),
\end{equation*}
with
\begin{equation*}
z_{o} \colon \mathbb{N} \times \mathbb{N} \! \to \! \mathbb{R}_{+}, \,
(\mathscr{N},n) \! \mapsto \! z_{o} \! := \mathscr{N}/n,
\end{equation*}
and where the `scaled' external field $\widetilde{V} \colon \mathbb{R}
\setminus \{0\} \! \to \! \mathbb{R}$ satisfies the following conditions:
\begin{gather}
\widetilde{V} \, \, \text{is real analytic on} \, \, \mathbb{R} \setminus
\{0\}; \\
\lim_{\vert x \vert \to \infty} \! \left(\widetilde{V}(x)/\ln (x^{2} \! + \!
1) \right) \! = \! +\infty; \\
\lim_{\vert x \vert \to 0} \! \left(\widetilde{V}(x)/\ln (x^{-2} \! + \! 1)
\right) \! = \! +\infty.
\end{gather}
(For example, a rational function of the form $\widetilde{V}(z) \! = \! \sum_{
k=-2m_{1}}^{2m_{2}} \widetilde{\varrho}_{k}z^{k}$, with $\widetilde{\varrho}_{
k} \! \in \! \mathbb{R}$, $k \! = \! -2m_{1},\dotsc,2m_{2}$, $m_{1,2} \! \in
\! \mathbb{N}$, and $\widetilde{\varrho}_{-2m_{1}},\widetilde{\varrho}_{2
m_{2}} \! > \! 0$ would satisfy conditions~(2.3)--(2.5).)

Hereafter, the double-scaling limit as $\mathscr{N},n \! \to \! \infty$ such 
that $z_{o} \! = \! 1 \! + \! o(1)$ is studied (the simplified `notation' $n 
\! \to \! \infty$ will be adopted). \hfill $\blacksquare$
\end{eeee}

It is, by now, a well-known, if not established, mathematical fact that 
variational conditions for minimisation problems in logarithmic potential 
theory, via the \emph{equilibrium measure} \cite{a55,a56,a90,a91,a92}, play 
a crucial r\^{o}le in the asymptotic analysis of (matrix) RHPs associated 
with (continuous and discrete) orthogonal polynomials, their roots, and 
corresponding recurrence relation coefficients (see, for example, 
\cite{a58,a59,a61,a65,a75}). The situation with respect to the large-$n$ 
asymptotic analysis for the monic OLPs, $\boldsymbol{\pi}_{n}(z)$, is 
analogous; but, unlike the asymptotic analysis for the orthogonal polynomials 
case, the asymptotic analysis for $\boldsymbol{\pi}_{n}(z)$ requires the 
consideration of two different families of RHPs, one for even degree 
(\textbf{RHP1}) and one for odd degree (\textbf{RHP2}). Thus, one must 
consider two sets of variational conditions for two (suitably posed) 
minimisation problems.

The following discussion is decomposed into two parts: one part corresponding
to the RHP for $\overset{e}{\mathrm{Y}} \colon \mathbb{C} \setminus \mathbb{R}
\! \to \! \operatorname{SL}_{2}(\mathbb{C})$ formulated as \textbf{RHP1},
denoted by $\pmb{\mathrm{P}_{1}}$, and the other part corresponding to the RHP
for $\overset{o}{\mathrm{Y}} \colon \mathbb{C} \setminus \mathbb{R} \! \to \!
\operatorname{SL}_{2}(\mathbb{C})$ formulated as \textbf{RHP2}, denoted by
$\pmb{\mathrm{P}_{2}}$.
\begin{compactenum}
\item[\shadowbox{$\pmb{\mathrm{P}_{1}}$}] Let $\widetilde{V} \colon \mathbb{R}
\setminus \{0\} \! \to \! \mathbb{R}$ satisfy conditions~(2.3)--(2.5). Let
$\mathrm{I}_{V}^{e}[\mu^{e}] \colon \mathcal{M}_{1}(\mathbb{R}) \! \to \!
\mathbb{R}$ denote the functional
\begin{equation*}
\mathrm{I}_{V}^{e}[\mu^{e}] \! = \! \iint_{\mathbb{R}^{2}} \ln \! \left(
\dfrac{\lvert st \rvert}{\lvert s \! - \! t \rvert^{2}} \right) \md \mu^{e}(s)
\, \md \mu^{e}(t) \! + \! 2 \int_{\mathbb{R}} \widetilde{V}(s) \, \md \mu^{e}
(s),
\end{equation*}
and consider the associated minimisation problem,
\begin{equation*}
E_{V}^{e} \! = \! \inf \lbrace \mathstrut \mathrm{I}_{V}^{e}[\mu^{e}]; \,
\mu^{e} \! \in \! \mathcal{M}_{1}(\mathbb{R}) \rbrace.
\end{equation*}
The infimum is finite, and there exists a unique measure $\mu_{V}^{e}$,
referred to as the `even' equilibrium measure, achieving the infimum (that is,
$\mathcal{M}_{1}(\mathbb{R}) \! \ni \! \mu_{V}^{e} \! = \! \inf \lbrace
\mathstrut \mathrm{I}_{V}^{e}[\mu^{e}]; \, \mu^{e} \! \in \! \mathcal{M}_{1}
(\mathbb{R}) \rbrace)$. Furthermore, $\mu_{V}^{e}$ has the following
`regularity' properties (all of these results are proven in this work):
\begin{compactenum}
\item[\textbullet] the `even' equilibrium measure has compact support which
consists of the disjoint union of a finite number of bounded real intervals;
in fact, as shown in Section~3 (see Lemma~3.5), $\mathrm{supp}(\mu_{V}^{e})
\! =: \! J_{e}$\footnote{It would be more usual, {}from the outset, for the
bounded (and closed) set $\overline{J_{e}} \! := \! \cup_{j=1}^{N+1}[b_{j-1}^{
e},a_{j}^{e}]$ to denote the support of $\mu_{V}^{e}$; however, the open (and
bounded) set $J_{e}$ provides an effective description of (the interior of)
the support of $\mu_{V}^{e}$: for this reason, $J_{e}$ (and at other times
$\overline{J_{e}})$ is used to denote $\operatorname{supp}(\mu_{V}^{e})$;
\emph{mutatis mutandis} for $J_{o}$ and $\overline{J_{o}}$ (see
\shadowbox{$\mathrm{P}_{2}$} below). This should not cause confusion for the
reader.} $\! = \! \cup_{j=1}^{N+1}(b_{j-1}^{e},a_{j}^{e})$ $(\subset \!
\mathbb{R} \setminus \{0\})$, where $\{b_{j-1}^{e},a_{j}^{e}\}_{j=1}^{N+1}$,
with $b_{0}^{e} \! := \! \min \{\mathrm{supp} \linebreak[4]
(\mu_{V}^{e})\} \! \notin \! \lbrace -\infty,0 \rbrace$, $a_{N+1}^{e} \! := \!
\max \{\mathrm{supp}(\mu_{V}^{e})\} \! \notin \! \lbrace 0,+\infty \rbrace$,
and $-\infty \! < \! b_{0}^{e} \! < \! a_{1}^{e} \! < \! b_{1}^{e} \! < \!
a_{2}^{e} \! < \! \cdots \! < \! b_{N}^{e} \! < \! a_{N+1}^{e} \! < \!
+\infty$, constitute the end-points of the support of $\mu_{V}^{e}$;
\item[\textbullet] the end-points $\{b_{j-1}^{e},a_{j}^{e}\}_{j=1}^{N+1}$ are
not arbitrary; rather, they satisfy the $n$-dependent and (locally) solvable
system of $2(N \! + \! 1)$ \emph{moment conditions} (transcendental equations)
given in Lemma~3.5;
\item[\textbullet] the `even' equilibrium measure is absolutely continuous
with respect to Lebesgue measure. The \emph{density} is given by
\begin{equation*}
\md \mu_{V}^{e}(x) \! := \! \psi_{V}^{e}(x) \, \md x \! = \! \dfrac{1}{2 \pi
\mi}(R_{e}(x))^{1/2}_{+}h_{V}^{e}(x) \pmb{1}_{J_{e}}(x) \, \md x,
\end{equation*}
where
\begin{equation*}
(R_{e}(z))^{1/2} \! := \! \left(\prod_{j=1}^{N+1}(z \! - \! b_{j-1}^{e})(z \!
- \! a_{j}^{e}) \right)^{1/2},
\end{equation*}
with $(R_{e}(x))^{1/2}_{\pm} \! := \! \lim_{\varepsilon \downarrow 0}(R_{e}
(x \! \pm \! \mi \varepsilon))^{1/2}$ and the branch of the square root is
chosen, as per the discussion in Subsection~2.1, such that $z^{-(N+1)}(R_{e}
(z))^{1/2} \! \sim_{\underset{z \in \mathbb{C}_{\pm}}{z \to \infty}} \! \pm
1$, $h_{V}^{e}(z) \! := \! \tfrac{1}{2} \oint_{C_{\mathrm{R}}^{e}}(R_{e}(s)
)^{-1/2}(\tfrac{\mi}{\pi s} \! + \! \tfrac{\mi \widetilde{V}^{\prime}(s)}{2
\pi})(s \! - \! z)^{-1} \, \md s$ (real analytic for $z \! \in \! \mathbb{R}
\setminus \{0\})$, where ${}^{\prime}$ denotes differentiation with respect to
the argument, $C_{\mathrm{R}}^{e}$ $(\subset \mathbb{C}^{\ast})$ is the union
of two circular contours, one outer one of large radius $R^{\natural}$
traversed clockwise and one inner one of small radius $r^{\natural}$ traversed
counter-clockwise, with the numbers $0 \! < \! r^{\natural} \! < \! R^{
\natural} \! < \! +\infty$ chosen such that, for (any) non-real $z$ in the
domain of analyticity of $\widetilde{V}$ (that is, $\mathbb{C}^{\ast})$,
$\mathrm{int}(C_{\mathrm{R}}^{e}) \! \supset \! J_{e} \cup \{z\}$, and $\pmb{
1}_{J_{e}}(x)$ denotes the indicator (characteristic) function of the set $J_{
e}$. (Note that $\psi_{V}^{e}(x) \! \geqslant \! 0 \, \, \forall \, \, x \!
\in \! \overline{J_{e}} \! := \! \cup_{j=1}^{N+1}[b_{j-1}^{e},a_{j}^{e}]$: it
vanishes like a square root at the end-points of the support of the `even'
equilibrium measure, that is, $\psi_{V}^{e}(s) \! =_{s \downarrow b_{j-1}^{e}}
\! \mathcal{O}((s \! - \! b_{j-1}^{e})^{1/2})$ and $\psi_{V}^{e}(s) \! =_{s
\uparrow a_{j}^{e}} \! \mathcal{O}((a_{j}^{e} \! - \! s)^{1/2})$, $j \! = \!
1,\dotsc,N \! + \! 1$.);
\item[\textbullet] the `even' equilibrium measure and its (compact) support
are uniquely characterised by the following Euler-Lagrange variational
equations: there exists $\ell_{e} \! \in \! \mathbb{R}$, the `even' Lagrange
multiplier, and $\mu^{e} \! \in \! \mathcal{M}_{1}(\mathbb{R})$ such that
\begin{gather*}
4 \int_{J_{e}} \ln (\vert x \! - \! s \vert) \, \md \mu^{e} (s) \! - \! 2 \ln
\vert x \vert \! - \! \widetilde{V}(x) \! - \! \ell_{e} \! = \! 0, \quad x \!
\in \! \overline{J_{e}}, \tag{$\mathrm{P}_{1}^{(a)}$} \\
4 \int_{J_{e}} \ln (\vert x \! - \! s \vert) \, \md \mu^{e} (s) \! - \! 2 \ln
\vert x \vert \! - \! \widetilde{V}(x) \! - \! \ell_{e} \! \leqslant \! 0,
\quad x \! \in \! \mathbb{R} \setminus \overline{J_{e}}; \tag{$\mathrm{P}_{
1}^{(b)}$}
\end{gather*}
\item[\textbullet] the Euler-Lagrange variational equations can be 
conveniently recast in terms of the complex potential $g^{e}(z)$ of 
$\mu_{V}^{e}$:
\begin{equation*}
g^{e}(z) \! := \! \int_{J_{e}} \! \ln \! \left((z \! - \! s)^{2}(zs)^{-1}
\right) \md \mu_{V}^{e}(s), \quad z \! \in \! \mathbb{C} \setminus (-\infty,
\max \{0,a_{N+1}^{e}\}).
\end{equation*}
The function $g^{e} \colon \mathbb{C} \setminus (-\infty,\max \{0,a_{N+1}^{e}
\}) \! \to \! \mathbb{C}$ so defined satisfies:
\begin{compactenum}
\item[$(\mathrm{P}_{1}^{(1)})$] $g^{e}(z)$ is analytic for $z \! \in \!
\mathbb{C} \setminus (-\infty,\max \{0,a_{N+1}^{e}\})$;
\item[$(\mathrm{P}_{1}^{(2)})$] $g^{e}(z) \! =_{\underset{z \in \mathbb{C}
\setminus \mathbb{R}}{z \to \infty}} \! \ln (z) \! + \! \mathcal{O}(1)$;
\item[$(\mathrm{P}_{1}^{(3)})$] $g^{e}_{+}(z) \! + \! g^{e}_{-}(z) \! -
\! \widetilde{V}(z) \! - \! \ell_{e} \! + \! 2Q_{e} \! = \! 0$, $z \! \in \!
\overline{J_{e}}$, where $g^{e}_{\pm}(z) \! := \! \lim_{\varepsilon \downarrow
0}g^{e}(z \! \pm \! \mi \varepsilon)$, and $Q_{e} \! := \! \int_{J_{e}} \ln
(s) \, \md \mu_{V}^{e}(s) \! = \! \int_{J_{e}} \ln (\lvert s \rvert) \, \md
\mu_{V}^{e}(s) \! + \! \mi \pi \int_{J_{e} \cap \mathbb{R}_{-}} \md \mu_{V}^{
e}(s)$;
\item[$(\mathrm{P}_{1}^{(4)})$] $g^{e}_{+}(z) \! + \! g^{e}_{-}(z) \! - \!
\widetilde{V}(z) \! - \! \ell_{e} \! + \! 2Q_{e} \! \leqslant \! 0$, $z \! \in
\! \mathbb{R} \setminus \overline{J_{e}}$, where equality holds for at most a
finite number of points;
\item[$(\mathrm{P}_{1}^{(5)})$] $g^{e}_{+}(z) \! - \! g^{e}_{-}(z) \! = \!
\mi f_{g^{e}}^{\mathbb{R}}(z)$, $z \! \in \! \mathbb{R}$, where $f_{g^{e}}^{
\mathbb{R}} \colon \mathbb{R} \! \to \! \mathbb{R}$, and, in particular, $g^{
e}_{+}(z) \! - \! g^{e}_{-}(z) \! = \! \mi \operatorname{const.}$, $z \! \in
\! \mathbb{R} \setminus \overline{J_{e}}$, with $\operatorname{const.} \! \in
\! \mathbb{R}$;
\item[$(\mathrm{P}_{1}^{(6)})$] $\mi (g^{e}_{+}(z) \! - \! g^{e}_{-}(z))^{
\prime} \! \geqslant \! 0$, $z \! \in \! J_{e}$, where equality holds for at
most a finite number of points.
\end{compactenum}
\end{compactenum}
\item[\shadowbox{$\pmb{\mathrm{P}_{2}}$}] Let $\widetilde{V} \colon \mathbb{
R} \setminus \{0\} \! \to \! \mathbb{R}$ satisfy conditions~(2.3)--(2.5). Let
$\mathrm{I}_{V}^{o}[\mu^{o}] \colon \mathcal{M}_{1}(\mathbb{R}) \! \to \!
\mathbb{R}$ denote the functional
\begin{equation*}
\mathrm{I}_{V}^{o}[\mu^{o}] \! = \! \iint_{\mathbb{R}^{2}} \ln \! \left(\dfrac{
\lvert st \rvert}{\lvert s \! - \! t \rvert^{2+\frac{1}{n}}} \right) \md \mu^{
o}(s) \, \md \mu^{o}(t) \! + \! 2 \int_{\mathbb{R}} \widetilde{V}(s) \, \md
\mu^{o}(s), \quad n \! \in \! \mathbb{N},
\end{equation*}
and consider the associated minimisation problem,
\begin{equation*}
E_{V}^{o} \! = \! \inf \lbrace \mathstrut \mathrm{I}_{V}^{o}[\mu^{o}]; \,
\mu^{o} \! \in \! \mathcal{M}_{1}(\mathbb{R}) \rbrace.
\end{equation*}
The infimum is finite, and there exists a unique measure $\mu_{V}^{o}$,
referred to as the `odd' equilibrium measure, achieving the infimum (that
is, $\mathcal{M}_{1}(\mathbb{R}) \! \ni \! \mu_{V}^{o} \! = \! \inf \lbrace
\mathstrut \mathrm{I}_{V}^{o}[\mu^{o}]; \, \mu^{o} \! \in \! \mathcal{M}_{1}
(\mathbb{R}) \rbrace)$. Furthermore, $\mu_{V}^{o}$ has the following
`regularity' properties (see \cite{a51} for complete details and proofs):
\begin{compactenum}
\item[\textbullet] the `odd' equilibrium measure has compact support which
consists of the disjoint union of a finite number of bounded real intervals;
in fact, as shown in \cite{a51}, $\mathrm{supp}(\mu_{V}^{o}) \! =: \! J_{o}
\! = \! \cup_{j=1}^{N+1}(b_{j-1}^{o},a_{j}^{o})$ $(\subset \! \mathbb{R}
\setminus \{0\})$, where $\{b_{j-1}^{o},a_{j}^{o}\}_{j=1}^{N+1}$, with $b_{
0}^{o} \! := \! \min \{\mathrm{supp} (\mu_{V}^{o})\} \! \notin \! \lbrace
-\infty,0 \rbrace$, $a_{N+1}^{o} \! := \! \max \{\mathrm{supp}(\mu_{V}^{o})\}
\! \notin \! \lbrace 0,+\infty \rbrace$, and $-\infty \! < \! b_{0}^{o} \! <
\! a_{1}^{o} \! < \! b_{1}^{o} \! < \! a_{2}^{o} \! < \! \cdots \! < \! b_{
N}^{o} \! < \! a_{N+1}^{o} \! < \! +\infty$, constitute the end-points of the
support of $\mu_{V}^{o}$; (The number of intervals, $N \! + \! 1$, is the same
in the `odd' case as in the `even' case, which can be established by a lengthy
analysis similar to that contained in \cite{a92}.)
\item[\textbullet] the end-points $\{b_{j-1}^{o},a_{j}^{o}\}_{j=1}^{N+1}$ are
not arbitrary; rather, they satisfy an $n$-dependent and (locally) solvable
system of $2(N \! + \! 1)$ moment conditions (transcendental equations; see
\cite{a51}, Lemma~3.5);
\item[\textbullet] the `odd' equilibrium measure is absolutely continuous
with respect to Lebesgue measure. The density is given by
\begin{equation*}
\md \mu_{V}^{o}(x) \! := \! \psi_{V}^{o}(x) \, \md x \! = \! \dfrac{1}{2 \pi
\mi}(R_{o}(x))^{1/2}_{+}h_{V}^{o}(x) \pmb{1}_{J_{o}}(x) \, \md x,
\end{equation*}
where
\begin{equation*}
(R_{o}(z))^{1/2} \! := \! \left(\prod_{j=1}^{N+1}(z \! - \! b_{j-1}^{o})(z \!
- \! a_{j}^{o}) \right)^{1/2},
\end{equation*}
with $(R_{o}(x))^{1/2}_{\pm} \! := \! \lim_{\varepsilon \downarrow 0}(R_{o}
(x \! \pm \! \mi \varepsilon))^{1/2}$ and the branch of the square root is
chosen, as per the discussion in Subsection~2.1, such that $z^{-(N+1)}(R_{o}
(z))^{1/2} \! \sim_{\underset{z \in \mathbb{C}_{\pm}}{z \to \infty}} \! \pm
1$, $h_{V}^{o}(z) \! := \! (2 \! + \! \tfrac{1}{n})^{-1} \oint_{C_{\mathrm{
R}}^{o}}(R_{o}(s))^{-1/2}(\tfrac{\mi}{\pi s} \! + \! \frac{\mi \widetilde{
V}^{\prime}(s)}{2 \pi})(s \! - \! z)^{-1} \, \md s$ (real analytic for $z \!
\in \! \mathbb{R} \setminus \{0\})$, where $C_{\mathrm{R}}^{o}$ $(\subset
\mathbb{C}^{\ast})$ is the union of two circular contours, one outer one of
large radius $R^{\flat}$ traversed clockwise and one inner one of small radius
$r^{\flat}$ traversed counter-clockwise, with the numbers $0 \! < \! r^{\flat}
\! < \! R^{\flat} \! < \! +\infty$ chosen such that, for (any) non-real $z$ in
the domain of analyticity of $\widetilde{V}$ (that is, $\mathbb{C}^{\ast})$,
$\mathrm{int}(C_{\mathrm{R}}^{o}) \! \supset \! J_{o} \cup \{z\}$, and $\pmb{
1}_{J_{o}}(x)$ denotes the indicator (characteristic) function of the set $J_{
o}$. (Note that $\psi_{V}^{o}(x) \! \geqslant \! 0 \, \, \forall \, \, x \!
\in \! \overline{J_{o}} \! := \! \cup_{j=1}^{N+1}[b_{j-1}^{o},a_{j}^{o}]$:
it vanishes like a square root at the end-points of the support of the `odd'
equilibrium measure, that is, $\psi_{V}^{o}(s) \! =_{s \downarrow b_{j-1}^{o}
} \! \mathcal{O}((s \! - \! b_{j-1}^{o})^{1/2})$ and $\psi_{V}^{o}(s) \! =_{s
\uparrow a_{j}^{o}} \! \mathcal{O}((a_{j}^{o} \! - \! s)^{1/2})$, $j \! = \!
1,\dotsc,N \! + \! 1$.);
\item[\textbullet] the `odd' equilibrium measure and its (compact) support are
uniquely characterised by the following Euler-Lagrange variational equations:
there exists $\ell_{o} \! \in \! \mathbb{R}$, the `odd' Lagrange multiplier,
and $\mu^{o} \! \in \! \mathcal{M}_{1}(\mathbb{R})$ such that
\begin{gather*}
2 \! \left(2 \! + \! \dfrac{1}{n} \right) \! \int_{J_{o}} \ln (\vert x \! - \!
s \vert) \, \md \mu^{o} (s) \! - \! 2 \ln \vert x \vert \! - \! \widetilde{V}
(x) \! - \! \ell_{o} \! - \! 2 \left(2 \! + \! \dfrac{1}{n} \right) \!
\widetilde{Q}_{o} = \! 0, \quad x \! \in \! \overline{J_{o}},
\tag{$\mathrm{P}_{2}^{(a)}$} \\
2 \! \left(2 \! + \! \dfrac{1}{n} \right) \! \int_{J_{o}} \ln (\vert x \! - \!
s \vert) \, \md \mu^{o} (s) \! - \! 2 \ln \vert x \vert \! - \! \widetilde{V}
(x) \! - \! \ell_{o} \! - \! 2 \left(2 \! + \! \dfrac{1}{n} \right) \!
\widetilde{Q}_{o} \! \leqslant \! 0, \quad x \! \in \! \mathbb{R} \setminus
\overline{J_{o}}, \tag{$\mathrm{P}_{2}^{(b)}$}
\end{gather*}
where $\widetilde{Q}_{o} \! := \! \int_{J_{o}} \ln (\lvert s \rvert) \, \md
\mu^{o}(s)$;
\item[\textbullet] the Euler-Lagrange variational equations can be 
conveniently recast in terms of the complex potential $g^{o}(z)$ of 
$\mu_{V}^{o}$:
\begin{equation*}
g^{o}(z) \! := \! \int_{J_{o}} \ln \! \left((z \! - \! s)^{2+\frac{1}{n}}
(zs)^{-1} \right) \md \mu_{V}^{o}(s), \quad z \! \in \! \mathbb{C} \setminus
(-\infty,\max \{0,a_{N+1}^{o}\}).
\end{equation*}
The function $g^{o} \colon \mathbb{C} \setminus (-\infty,\max \{0,a_{N+1}^{o}
\}) \! \to \! \mathbb{C}$ so defined satisfies:
\begin{compactenum}
\item[$(\mathrm{P}_{2}^{(1)})$] $g^{o}(z)$ is analytic for $z \! \in \!
\mathbb{C} \setminus (-\infty,\max \{0,a_{N+1}^{o}\})$;
\item[$(\mathrm{P}_{2}^{(2)})$] $g^{o}(z) \! =_{\underset{z \in \mathbb{C}
\setminus \mathbb{R}}{z \to 0}} \! -\ln (z) \! + \! \mathcal{O}(1)$;
\item[$(\mathrm{P}_{2}^{(3)})$] $g^{o}_{+}(z) \! + \! g^{o}_{-}(z) \! - \!
\widetilde{V}(z) \! - \! \ell_{o} \! - \! \mathfrak{Q}_{\mathscr{A}}^{+} \! -
\! \mathfrak{Q}_{\mathscr{A}}^{-} \! = \! 0$, $z \! \in \! \overline{J_{o}}$,
where $g^{o}_{\pm}(z) \! := \! \lim_{\varepsilon \downarrow 0}g^{o}(z \! \pm
\! \mi \varepsilon)$, and $\mathfrak{Q}_{\mathscr{A}}^{\pm} \! := \! (1 \! +
\! \tfrac{1}{n}) \int_{J_{o}} \ln (\lvert s \rvert) \, \md \mu_{V}^{o}(s) \!
- \! \mi \pi \int_{J_{o} \cap \mathbb{R}_{-}} \md \mu_{V}^{o}(s) \! \pm \! \mi
\pi (2 \! + \! \tfrac{1}{n}) \int_{J_{o} \cap \mathbb{R}_{+}} \md \mu_{V}^{o}
(s)$;
\item[$(\mathrm{P}_{2}^{(4)})$] $g^{o}_{+}(z) \! + \! g^{o}_{-}(z) \! - \!
\widetilde{V}(z) \! - \! \ell_{o} \! - \! \mathfrak{Q}_{\mathscr{A}}^{+} \! -
\! \mathfrak{Q}_{\mathscr{A}}^{-} \! \leqslant \! 0$, $z \! \in \! \mathbb{R}
\setminus \overline{J_{o}}$, where equality holds for at most a finite number
of points;
\item[$(\mathrm{P}_{2}^{(5)})$] $g^{o}_{+}(z) \! - \! g^{o}_{-}(z) \! - \!
\mathfrak{Q}_{\mathscr{A}}^{+} \! + \! \mathfrak{Q}_{\mathscr{A}}^{-} \! = \!
\mi f_{g^{o}}^{\mathbb{R}}(z)$, $z \! \in \! \mathbb{R}$, where $f_{g^{o}}^{
\mathbb{R}} \colon \mathbb{R} \! \to \! \mathbb{R}$, and, in particular, $g^{
o}_{+}(z) \! - \! g^{o}_{-}(z) \! - \! \mathfrak{Q}_{\mathscr{A}}^{+} \! + \!
\mathfrak{Q}_{\mathscr{A}}^{-} \! = \! \mi \operatorname{const.}$, $z \! \in
\! \mathbb{R} \setminus \overline{J_{o}}$, with $\operatorname{const.} \! \in
\! \mathbb{R}$;
\item[$(\mathrm{P}_{2}^{(6)})$] $\mi (g^{o}_{+}(z) \! - \! g^{o}_{-}(z) \! -
\! \mathfrak{Q}_{\mathscr{A}}^{+} \! + \! \mathfrak{Q}_{\mathscr{A}}^{-} )^{
\prime} \! \geqslant \! 0$, $z \! \in \! J_{o}$, where equality holds for at
most a finite number of points.
\end{compactenum}
\end{compactenum}
\end{compactenum}

In this three-fold series of works on asymptotics of OLPs and related 
quantities, the so-called `regular case' is studied, namely:
\begin{compactenum}
\item[\textbullet] $\md \mu_{V}^{e}$, or $\widetilde{V} \colon \mathbb{R}
\setminus \{0\} \! \to \! \mathbb{R}$ satisfying conditions~(2.3)--(2.5), is
\emph{regular} if: (i) $h_{V}^{e}(x) \! \not\equiv \! 0$ on $\overline{J_{e}
}$; (ii) $4 \int_{J_{e}} \ln (\vert x \! - \! s \vert) \, \md \mu_{V}^{e}(s)
\! - \! 2 \ln \vert x \vert \! - \! \widetilde{V}(x) \! - \! \ell_{e} \!
< \! 0$, $x \! \in \! \mathbb{R} \setminus \overline{J_{e}}$; and (iii)
inequalities~$(\mathrm{P}_{1}^{(4)})$ and~$(\mathrm{P}_{1}^{(6)})$ in
$\pmb{\mathrm{P}_{1}}$ are strict, that is, $\leqslant$ (resp., $\geqslant)$
is replaced by $<$ (resp., $>)$;
\item[\textbullet] $\md \mu_{V}^{o}$, or $\widetilde{V} \colon \mathbb{R}
\setminus \{0\} \! \to \! \mathbb{R}$ satisfying conditions~(2.3)--(2.5), is
regular if: (i) $h_{V}^{o}(x) \! \not\equiv \! 0$ on $\overline{J_{o}}$; (ii)
$2(2 \! + \! \tfrac{1}{n}) \int_{J_{o}} \ln (\vert x \! - \! s \vert) \, \md
\mu_{V}^{o}(s) \! - \! 2 \ln \vert x \vert \! - \! \widetilde{V}(x) \! - \!
\ell_{o} \! - \! 2(2 \! + \! \tfrac{1}{n})Q_{o} \! < \! 0$, $x \! \in \!
\mathbb{R} \setminus \overline{J_{o}}$, where $Q_{o} \! := \! \int_{J_{o}}
\ln (\lvert s \rvert) \, \md \mu_{V}^{o}(s)$; and (iii)
inequalities~$(\mathrm{P}_{2}^{(4)})$ and~$(\mathrm{P}_{2}^{(6)})$ in
$\pmb{\mathrm{P}_{2}}$ are strict, that is, $\leqslant$ (resp., $\geqslant)$
is replaced by $<$ (resp., $>)$\footnote{There are three distinct situations
in which these conditions may fail: (i) for at least one $\widetilde{x}_{e} \!
\in \! \mathbb{R} \setminus \overline{J_{e}}$ (resp., $\widetilde{x}_{o} \!
\in \! \mathbb{R} \setminus \widetilde{J_{o}})$, $4 \int_{J_{e}} \ln (\lvert
\widetilde{x}_{e} \! - \! s \rvert) \, \md \mu_{V}^{e}(s) \! - \! 2 \ln \lvert
\widetilde{x}_{e} \rvert \! - \! \widetilde{V}(\widetilde{x}_{e}) \! - \!
\ell_{e} \! = \! 0$ (resp., $2(2 \! + \! \tfrac{1}{n}) \int_{J_{o}} \ln
(\lvert \widetilde{x}_{o} \! - \! s \rvert) \, \md \mu_{V}^{o}(s) \! - \! 2
\ln \lvert \widetilde{x}_{o} \rvert \! - \! \widetilde{V}(\widetilde{x}_{o})
\! - \! \ell_{o} \! - \! 2(2 \! + \! \tfrac{1}{n})Q_{o} \! = \! 0)$, that is,
for $n$ even (resp., $n$ odd) equality is attained for at least one point
$\widetilde{x}_{e}$ (resp., $\widetilde{x}_{o})$ in the complement of the
closure of the support of the `even' (resp., `odd') equilibrium measure $\mu_{
V}^{e}$ (resp., $\mu_{V}^{o})$, which corresponds to the situation in which a
`band' has just closed, or is about to open, about $\widetilde{x}_{e}$ (resp.,
$\widetilde{x}_{o})$; (ii) for at least one $\widehat{x}_{e}$ (resp.,
$\widehat{x}_{o})$, $h_{V}^{e}(\widehat{x}_{e}) \! = \! 0$ (resp., $h_{V}^{o}
(\widehat{x}_{o}) \! = \! 0)$, that is, for $n$ even (resp., $n$ odd) the
function $h_{V}^{e}$ (resp., $h_{V}^{o})$ vanishes for at least one point
$\widehat{x}_{e}$ (resp., $\widehat{x}_{o})$ within the support of the `even'
(resp., `odd') equilibrium measure $\mu_{V}^{e}$ (resp., $\mu_{V}^{o})$, which
corresponds to the situation in which a `gap' is about to open, or close,
about $\widehat{x}_{e}$ (resp., $\widehat{x}_{o})$; and (iii) there exists at
least one $j \! \in \! \lbrace 1,\dotsc,N \! + \! 1 \rbrace$, denoted $j_{e}$
(resp., $j_{o})$, such that $h_{V}^{e}(b_{j_{e}-1}^{e}) \! = \! 0$ and/or $h_{
V}^{e}(a_{j_{e}}^{e}) \! = \! 0$ (resp., $h_{V}^{o}(b_{j_{o}-1}^{o}) \! = \!
0$ and/or $h_{V}^{o}(a_{j_{o}}^{o}) \! = \! 0)$. Each of these three cases
can occur only a finite number of times due to the fact that $\widetilde{V}
\colon \mathbb{R} \setminus \lbrace 0 \rbrace \! \to \! \mathbb{R}$ satisfies
conditions~(2.3)--(2.5) \cite{a58,a92}.}.
\end{compactenum}
The (density of the) `even' and `odd' equilibrium measures $\md \mu_{V}^{e}$
and $\md \mu_{V}^{o}$, respectively, together with the corresponding
variational problems, emerge naturally in the asymptotic analyses of
\textbf{RHP1} and \textbf{RHP2}.
\begin{eeee}
The following correspondences should also be noted:
\begin{compactenum}
\item[\textbullet] $g^{e} \colon \mathbb{C} \setminus (-\infty,
\max \{0,a_{N+1}^{e}\}) \! \to \! \mathbb{C}$ solves the
\emph{phase conditions}~$(\mathrm{P}_{1}^{(1)})$--$(\mathrm{P}_{1}^{(6)}) \!
\Leftrightarrow \! \mathcal{M}_{1}(\mathbb{R}) \! \ni \! \mu_{V}^{e}$ solves
the variational conditions~$(\mathrm{P}_{1}^{(a)})$ and~$(\mathrm{P}_{1}^{
(b)})$;
\item[\textbullet] $g^{o} \colon \mathbb{C} \setminus (-\infty,
\max \{0,a_{N+1}^{o}\}) \! \to \! \mathbb{C}$ solves the phase
conditions~$(\mathrm{P}_{2}^{(1)})$--$(\mathrm{P}_{2}^{(6)}) \!
\Leftrightarrow \! \mathcal{M}_{1}(\mathbb{R}) \! \ni \! \mu_{V}^{o}$ solves
the variational conditions~$(\mathrm{P}_{2}^{(a)})$ and~$(\mathrm{P}_{2}^{
(b)})$. \hfill $\blacksquare$
\end{compactenum}
\end{eeee}

Since the main results of this paper are asymptotics (as $n \! \to \! \infty)$ 
for $\boldsymbol{\pi}_{2n}(z)$ $(z \! \in \! \mathbb{C})$, $\xi^{(2n)}_{n}$ 
and $\phi_{2n}(z)$ $(z \! \in \! \mathbb{C})$, which are, via Lemma~2.2.1, 
Equation~(2.2), and Equations~(1.2) and~(1.4), related to \textbf{RHP1} for 
$\overset{e}{\mathrm{Y}} \colon \mathbb{C} \setminus \mathbb{R} \! \to \! 
\operatorname{SL}_{2}(\mathbb{C})$, no further reference, henceforth, to 
\textbf{RHP2} (and Lemma~2.2.2) for $\overset{o}{\mathrm{Y}} \colon \mathbb{C} 
\setminus \mathbb{R} \! \to \! \operatorname{SL}_{2}(\mathbb{C})$ will be 
made (see \cite{a51} for the complete details of the asymptotic analysis of 
\textbf{RHP2}). In the ensuing analysis, the large-$n$ behaviour of the 
solution of \textbf{RHP1} (see Lemma~2.2.1, Equation~(2.2)), hence asymptotics 
for $\boldsymbol{\pi}_{2n}(z)$ (in the entire complex plane), $\xi^{(2n)}_{n}$ 
and $\phi_{2n}(z)$ (in the entire complex plane), are extracted. 
\subsection{Summary of Results}
In this subsection, the final results of this work are presented (see
Sections~3--5 for the detailed analyses and proofs). Before doing so, however,
some notational preamble is necessary. For $j \! = \! 1,\dotsc,N \! + \! 1$,
let
\begin{equation*}
\Phi_{a_{j}}^{e}(z) \! := \! \left(\dfrac{3n}{2} \int_{a_{j}^{e}}^{z}(R_{e}
(s))^{1/2}h_{V}^{e}(s) \, \md s \right)^{2/3} \quad \text{and} \quad \Phi_{
b_{j-1}}^{e}(z) \! := \! \left(-\dfrac{3n}{2} \int_{z}^{b_{j-1}^{e}}(R_{e}(s)
)^{1/2}h_{V}^{e}(s) \, \md s \right)^{2/3},
\end{equation*}
where $(R_{e}(z))^{1/2}$ and $h_{V}^{e}(z)$ are defined in Theorem~2.3.1,
Equations~(2.8) and~(2.9). Define the `small', mutually disjoint open discs
about the end-points of the support of the `even' equilibrium measure, $\{
b_{j-1}^{e},a_{j}^{e}\}_{j=1}^{N+1}$, as follows: for $j \! = \! 1,\dotsc,N \!
+ \! 1$,
\begin{equation*}
\mathbb{U}^{e}_{\delta_{a_{j}}} \! := \! \left\{\mathstrut z \! \in \!
\mathbb{C}; \, \vert z \! - \! a_{j}^{e} \vert \! < \! \delta_{a_{j}}^{e}
\right\} \quad \quad \text{and} \quad \quad \mathbb{U}^{e}_{\delta_{b_{j-1}}}
\! := \! \left\{\mathstrut z \! \in \! \mathbb{C}; \, \vert z \! - \! b_{j-
1}^{e} \vert \! < \! \delta^{e}_{b_{j-1}} \right\},
\end{equation*}
where $(0,1) \! \ni \! \delta^{e}_{a_{j}}$ (resp., $(0,1) \! \ni \! \delta^{
e}_{b_{j-1}})$ are chosen `sufficiently small' so that $\Phi^{e}_{a_{j}}(z)$
(resp., $\Phi^{e}_{b_{j-1}}(z))$, which are bi-holomorphic, conformal, and
orientation preserving (resp., bi-holomorphic, conformal, and non-orientation
preserving), map $\mathbb{U}^{e}_{\delta_{a_{j}}}$ (resp., $\mathbb{U}^{e}_{
\delta_{b_{j-1}}})$, as well as the oriented skeletons (see Figure~5) $\cup_{
l=1}^{4} \Sigma^{e,l}_{a_{j}}$ (resp., $\cup_{l=1}^{4} \Sigma^{e,l}_{b_{j-
1}}$ (see Figure~6)), injectively onto open (and convex), $n$-dependent
neighbourhoods of $0$ such that:
\begin{compactenum}
\item[\pmb{(i)}] $\Phi^{e}_{a_{j}}(a_{j}^{e}) \! = \! 0$ (resp., $\Phi^{e}_{
b_{j-1}}(b_{j-1}^{e}) \! = \! 0)$;
\item[\pmb{(ii)}] $\Phi^{e}_{a_{j}} \colon \mathbb{U}^{e}_{\delta_{a_{j}}} \!
\to \! \widehat{\mathbb{U}}^{e}_{\delta_{a_{j}}} \! := \! \Phi^{e}_{a_{j}}
(\mathbb{U}^{e}_{\delta_{a_{j}}})$ (resp., $\Phi^{e}_{b_{j-1}} \colon \mathbb{
U}^{e}_{\delta_{b_{j-1}}} \! \to \! \widehat{\mathbb{U}}^{e}_{\delta_{b_{j-1}}
} \! := \! \Phi^{e}_{b_{j-1}}(\mathbb{U}^{e}_{\delta_{b_{j-1}}}))$;
\item[\pmb{(iii)}] $\Phi^{e}_{a_{j}}(\mathbb{U}^{e}_{\delta_{a_{j}}} \cap
\Sigma^{e,l}_{a_{j}}) \! = \! \Phi^{e}_{a_{j}}(\mathbb{U}^{e}_{\delta_{a_{j}
}}) \cap \gamma^{e,l}_{a_{j}}$ (resp., $\Phi^{e}_{b_{j-1}}(\mathbb{U}^{e}_{
\delta_{b_{j-1}}} \cap \Sigma^{e,l}_{b_{j-1}}) \! = \! \Phi^{e}_{b_{j-1}}
(\mathbb{U}^{e}_{\delta_{b_{j-1}}}) \cap \gamma^{e,l}_{b_{j-1}})$;
\item[\pmb{(iv)}] $\Phi^{e}_{a_{j}}(\mathbb{U}^{e}_{\delta_{a_{j}}} \cap
\Omega^{e,l}_{a_{j}}) \! = \! \Phi^{e}_{a_{j}}(\mathbb{U}^{e}_{\delta_{a_{j}}
}) \cap \widehat{\Omega}^{e,l}_{a_{j}}$ (resp., $\Phi^{e}_{b_{j-1}}(\mathbb{
U}^{e}_{\delta_{b_{j-1}}} \cap \Omega^{e,l}_{b_{j-1}}) \! = \! \Phi^{e}_{b_{j
-1}}(\mathbb{U}^{e}_{\delta_{b_{j-1}}}) \cap \widehat{\Omega}^{e,l}_{b_{j-1}
})$, with $\widehat{\Omega}^{e,1}_{a_{j}}$ (and $\widehat{\Omega}^{e,1}_{b_{j
-1}})$ $= \! \lbrace \mathstrut \zeta \! \in \! \mathbb{C}; \, \arg (\zeta) \!
\in \! (0,2 \pi/3) \rbrace$, $\widehat{\Omega}^{e,2}_{a_{j}}$ (and $\widehat{
\Omega}^{e,2}_{b_{j-1}})$ $= \! \lbrace \mathstrut \zeta \! \in \! \mathbb{C};
\, \arg (\zeta) \! \in \! (2 \pi/3,\pi) \rbrace$, $\widehat{\Omega}^{e,3}_{a_{
j}}$ (and $\widehat{\Omega}^{e,3}_{b_{j-1}})$ $= \! \lbrace \mathstrut \zeta
\! \in \! \mathbb{C}; \, \arg (\zeta) \! \in \! (-\pi,-2 \pi/3) \rbrace$, and
$\widehat{\Omega}^{e,4}_{a_{j}}$ (and $\widehat{\Omega}^{e,4}_{b_{j-1}})$ $=
\! \lbrace \mathstrut \zeta \! \in \! \mathbb{C}; \, \arg (\zeta) \! \in \!
(-2 \pi/3,0) \rbrace$\footnote{The precise angles between the sectors are not
absolutely important; one could, for example, replace $2 \pi/3$ by any angle
strictly between $0$ and $\pi$ \cite{a2,a58,a59,a61,a90}.}.
\end{compactenum}
\begin{figure}[tbh]
\begin{center}
\vspace{0.55cm}
\begin{pspicture}(0,0)(14,8)
\psset{xunit=1cm,yunit=1cm}
\pscircle[linewidth=0.7pt,linestyle=solid,linecolor=red](3.5,4){2.5}
\pscircle[linewidth=0.7pt,linestyle=solid,linecolor=red](10.5,4){2.5}
\psarcn[linewidth=0.6pt,linestyle=solid,linecolor=blue,arrowsize=1.5pt 5]{->}%
(0.5,4){3}{65}{36}
\pstextpath[c]{\psarcn[linewidth=0.6pt,linestyle=solid,linecolor=blue](0.5,4)%
{3}{36}{0}}{\makebox(0,0){$\pmb{\Sigma^{e,1}_{a_{j}}}$}}
\psarc[linewidth=0.6pt,linestyle=solid,linecolor=green,arrowsize=1.5pt 5]{->}%
(0.5,4){3}{295}{324}
\pstextpath[c]{\psarc[linewidth=0.6pt,linestyle=solid,linecolor=green](0.5,4)%
{3}{324}{360}}{\makebox(0,0){$\pmb{\Sigma^{e,3}_{a_{j}}}$}}
\psline[linewidth=0.6pt,linestyle=solid,linecolor=cyan,arrowsize=1.5pt 5]{->}%
(0.5,4)(1.6,4)
\pstextpath[c]{\psline[linewidth=0.6pt,linestyle=solid,linecolor=cyan](1.6,4)%
(3.5,4)}{\makebox(0,0){$\pmb{\Sigma^{e,2}_{a_{j}}}$}}
\psline[linewidth=0.6pt,linestyle=solid,linecolor=magenta](5.5,4)(6.3,4)
\pstextpath[c]{\psline[linewidth=0.6pt,linestyle=solid,linecolor=magenta,%
arrowsize=1.5pt 5]{->}(3.5,4)(5.5,4)}{\makebox(0,0){%
$\pmb{\Sigma^{e,4}_{a_{j}}}$}}
\rput(3.5,7.5){\makebox(0,0){$z-\text{plane}$}}
\rput(3.5,0.5){\makebox(0,0){$\mathbb{U}^{e}_{\delta_{a_{j}}}$}}
\pszigzag[coilwidth=0.3cm,coilarm=0.25cm,coilaspect=45]{->}(3.7,0.7)(4.1,2.2)
\rput(4.65,5.1){\makebox(0,0){$\pmb{\Omega^{e,1}_{a_{j}}}$}}
\rput(2.35,5.1){\makebox(0,0){$\pmb{\Omega^{e,2}_{a_{j}}}$}}
\rput(2.35,2.9){\makebox(0,0){$\pmb{\Omega^{e,3}_{a_{j}}}$}}
\rput(4.65,2.9){\makebox(0,0){$\pmb{\Omega^{e,4}_{a_{j}}}$}}
\rput(7,7.2){\makebox(0,0){$\zeta \! = \! \Phi^{e}_{a_{j}}(z)$}}
\rput(7,0.8){\makebox(0,0){$z \! = \! (\Phi^{e}_{a_{j}})^{-1}(\zeta)$}}
\psarcn[linewidth=0.8pt,linestyle=solid,linecolor=black,arrowsize=1.5pt 5]%
{->}(7,5.3){1.5}{135}{45}
\psarc[linewidth=0.8pt,linestyle=solid,linecolor=black,arrowsize=1.5pt 5]%
{<-}(7,2.7){1.5}{225}{315}
\rput(10.5,7.5){\makebox(0,0){$\zeta -\text{plane}$}}
\rput(10.5,0.5){\makebox(0,0){$\widehat{\mathbb{U}}^{e}_{\delta_{a_{j}}} \!
:= \! \Phi^{e}_{a_{j}}(\mathbb{U}^{e}_{\delta_{a_{j}}})$}}
\pszigzag[coilwidth=0.3cm,coilarm=0.25cm,coilaspect=45]{->}(10.75,0.75)%
(10.15,2.2)
\psline[linewidth=0.6pt,linestyle=solid,linecolor=cyan,arrowsize=1.5pt 5]{->}%
(7.6,4)(8.6,4)
\pstextpath[c]{\psline[linewidth=0.6pt,linestyle=solid,linecolor=cyan](8.6,4)%
(10.5,4)}{\makebox(0,0){$\pmb{\gamma^{e,2}_{a_{j}}}$}}
\psline[linewidth=0.6pt,linestyle=solid,linecolor=magenta](12.5,4)(13.3,4)
\pstextpath[c]{\psline[linewidth=0.6pt,linestyle=solid,linecolor=magenta,%
arrowsize=1.5pt 5]{->}(10.5,4)(12.5,4)}{\makebox(0,0){$\pmb{\gamma^{e,4}_{%
a_{j}}}$}}
\psline[linewidth=0.6pt,linestyle=solid,linecolor=blue,arrowsize=1.5pt 5]{->}%
(8.4,6.1)(9.3,5.2)
\pstextpath[c]{\psline[linewidth=0.6pt,linestyle=solid,linecolor=blue]%
(9.3,5.2)(10.5,4)}{\makebox(0,0){$\pmb{\gamma^{e,1}_{a_{j}}}$}}
\psline[linewidth=0.6pt,linestyle=solid,linecolor=green,arrowsize=1.5pt 5]{->}%
(8.4,1.9)(9.3,2.8)
\pstextpath[c]{\psline[linewidth=0.6pt,linestyle=solid,linecolor=green]%
(9.3,2.8)(10.5,4)}{\makebox(0,0){$\pmb{\gamma^{e,3}_{a_{j}}}$}}
\rput(11.5,4.85){\makebox(0,0){$\pmb{\widehat{\Omega}^{e,1}_{a_{j}}}$}}
\rput(9,4.65){\makebox(0,0){$\pmb{\widehat{\Omega}^{e,2}_{a_{j}}}$}}
\rput(9,3.35){\makebox(0,0){$\pmb{\widehat{\Omega}^{e,3}_{a_{j}}}$}}
\rput(11.5,3.15){\makebox(0,0){$\pmb{\widehat{\Omega}^{e,4}_{a_{j}}}$}}
\psdots[dotstyle=*,dotscale=1.5](3.5,4)
\psdots[dotstyle=*,dotscale=1.5](10.5,4)
\rput(3.5,3.7){\makebox(0,0){$\pmb{a_{j}^{e}}$}}
\rput(10.5,3.7){\makebox(0,0){$\pmb{0}$}}
\end{pspicture}
\end{center}
\caption{The conformal mapping $\zeta \! = \! \Phi^{e}_{a_{j}}(z) \! := \!
(\tfrac{3n}{2} \int_{a_{j}^{e}}^{z}(R_{e}(s))^{1/2}h_{V}^{e}(s))^{2/3}$, $j
\! = \! 1,\dotsc,N \! + \! 1$, where $(\Phi^{e}_{a_{j}})^{-1}$ denotes the
inverse mapping}
\end{figure}
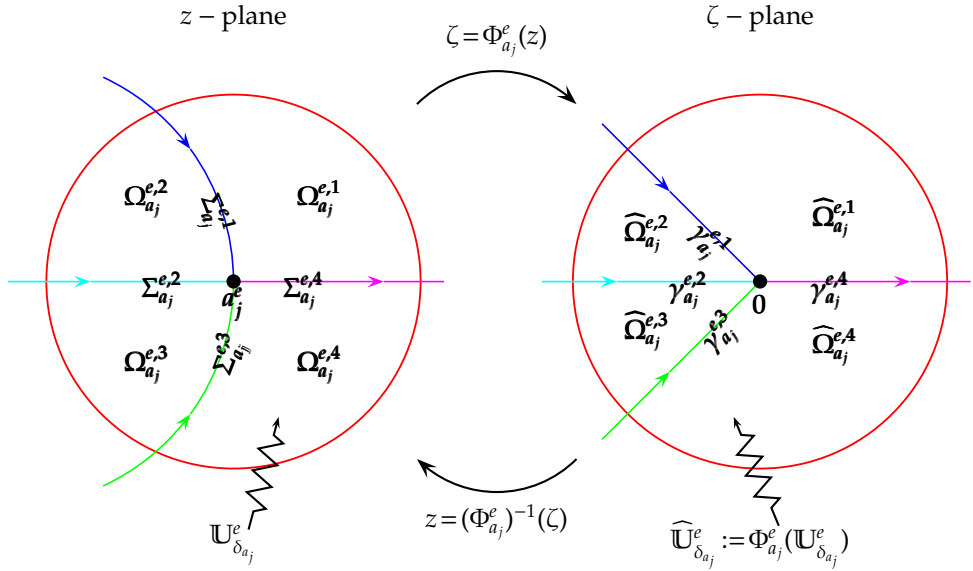
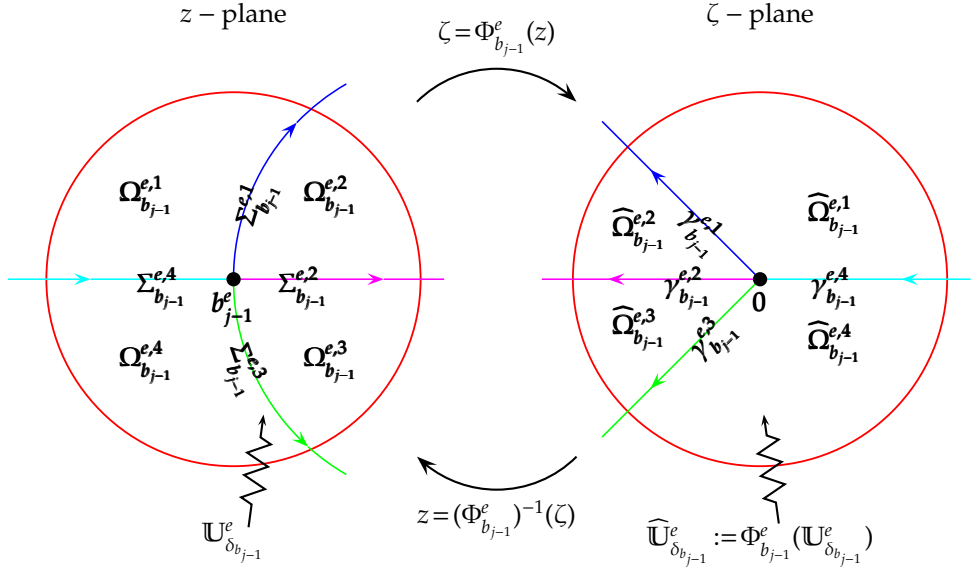
\begin{figure}[bht]
\begin{center}
\vspace{0.55cm}
\begin{pspicture}(0,0)(14,8)
\psset{xunit=1cm,yunit=1cm}
\pscircle[linewidth=0.7pt,linestyle=solid,linecolor=red](3.5,4){2.5}
\pscircle[linewidth=0.7pt,linestyle=solid,linecolor=red](10.5,4){2.5}
\psarcn[linewidth=0.6pt,linestyle=solid,linecolor=blue](6.5,4){3}{136}{120}
\pstextpath[c]{\psarcn[linewidth=0.6pt,linestyle=solid,linecolor=blue,%
arrowsize=1.5pt 5]{->}(6.5,4){3}{180}{136}}{\makebox(0,0){$\pmb{\Sigma^{e,%
1}_{b_{j-1}}}$}}
\psarc[linewidth=0.6pt,linestyle=solid,linecolor=green](6.5,4){3}{228}{240}
\pstextpath[c]{\psarc[linewidth=0.6pt,linestyle=solid,linecolor=green,%
arrowsize=1.5pt 5]{->}(6.5,4){3}{180}{228}}{\makebox(0,0){$\pmb{\Sigma^{e,%
3}_{b_{j-1}}}$}}
\psline[linewidth=0.6pt,linestyle=solid,linecolor=cyan,arrowsize=1.5pt 5]{->}%
(0.5,4)(1.6,4)
\pstextpath[c]{\psline[linewidth=0.6pt,linestyle=solid,linecolor=cyan](1.6,4)%
(3.5,4)}{\makebox(0,0){$\pmb{\Sigma^{e,4}_{b_{j-1}}}$}}
\psline[linewidth=0.6pt,linestyle=solid,linecolor=magenta](5.5,4)(6.3,4)
\pstextpath[c]{\psline[linewidth=0.6pt,linestyle=solid,linecolor=magenta,%
arrowsize=1.5pt 5]{->}(3.5,4)(5.5,4)}{\makebox(0,0){%
$\pmb{\Sigma^{e,2}_{b_{j-1}}}$}}
\rput(3.5,7.5){\makebox(0,0){$z-\text{plane}$}}
\rput(3.5,0.5){\makebox(0,0){$\mathbb{U}^{e}_{\delta_{b_{j-1}}}$}}
\pszigzag[coilwidth=0.3cm,coilarm=0.25cm,coilaspect=45]{->}(3.7,0.7)(3.9,2.2)
\rput(4.8,5.1){\makebox(0,0){$\pmb{\Omega^{e,2}_{b_{j-1}}}$}}
\rput(2.35,5.1){\makebox(0,0){$\pmb{\Omega^{e,1}_{b_{j-1}}}$}}
\rput(2.35,2.9){\makebox(0,0){$\pmb{\Omega^{e,4}_{b_{j-1}}}$}}
\rput(4.8,2.9){\makebox(0,0){$\pmb{\Omega^{e,3}_{b_{j-1}}}$}}
\rput(7,7.2){\makebox(0,0){$\zeta \! = \! \Phi^{e}_{b_{j-1}}(z)$}}
\rput(7,0.8){\makebox(0,0){$z \! = \! (\Phi^{e}_{b_{j-1}})^{-1}(\zeta)$}}
\psarcn[linewidth=0.8pt,linestyle=solid,linecolor=black,arrowsize=1.5pt 5]%
{->}(7,5.3){1.5}{135}{45}
\psarc[linewidth=0.8pt,linestyle=solid,linecolor=black,arrowsize=1.5pt 5]%
{<-}(7,2.7){1.5}{225}{315}
\rput(10.5,7.5){\makebox(0,0){$\zeta -\text{plane}$}}
\rput(10.5,0.5){\makebox(0,0){$\widehat{\mathbb{U}}^{e}_{\delta_{b_{j-1}}} \!
:= \! \Phi^{e}_{b_{j-1}}(\mathbb{U}^{e}_{\delta_{b_{j-1}}})$}}
\pszigzag[coilwidth=0.3cm,coilarm=0.25cm,coilaspect=45]{->}(10.75,0.75)%
(10.55,2.2)
\psline[linewidth=0.6pt,linestyle=solid,linecolor=magenta,arrowsize=1.5pt 5]%
{-<}(7.6,4)(8.7,4)
\pstextpath[c]{\psline[linewidth=0.6pt,linestyle=solid,linecolor=magenta]%
(8.6,4)(10.5,4)}{\makebox(0,0){$\pmb{\gamma^{e,2}_{b_{j-1}}}$}}
\psline[linewidth=0.6pt,linestyle=solid,linecolor=cyan](12.5,4)(13.3,4)
\pstextpath[c]{\psline[linewidth=0.6pt,linestyle=solid,linecolor=cyan,%
arrowsize=1.5pt 5]{-<}(10.5,4)(12.6,4)}{\makebox(0,0){$\pmb{\gamma^{e,4}_{%
b_{j-1}}}$}}
\psline[linewidth=0.6pt,linestyle=solid,linecolor=blue,arrowsize=1.5pt 5]{-<}%
(8.4,6.1)(9.2,5.3)
\pstextpath[c]{\psline[linewidth=0.6pt,linestyle=solid,linecolor=blue]%
(9.1,5.4)(10.5,4)}{\makebox(0,0){$\pmb{\gamma^{e,1}_{b_{j-1}}}$}}
\psline[linewidth=0.6pt,linestyle=solid,linecolor=green,arrowsize=1.5pt 5]{-<}%
(8.4,1.9)(9.2,2.7)
\pstextpath[c]{\psline[linewidth=0.6pt,linestyle=solid,linecolor=green]%
(9.1,2.6)(10.5,4)}{\makebox(0,0){$\pmb{\gamma^{e,3}_{b_{j-1}}}$}}
\rput(11.5,4.85){\makebox(0,0){$\pmb{\widehat{\Omega}^{e,1}_{b_{j-1}}}$}}
\rput(8.9,4.65){\makebox(0,0){$\pmb{\widehat{\Omega}^{e,2}_{b_{j-1}}}$}}
\rput(8.9,3.35){\makebox(0,0){$\pmb{\widehat{\Omega}^{e,3}_{b_{j-1}}}$}}
\rput(11.5,3.15){\makebox(0,0){$\pmb{\widehat{\Omega}^{e,4}_{b_{j-1}}}$}}
\psdots[dotstyle=*,dotscale=1.5](3.5,4)
\psdots[dotstyle=*,dotscale=1.5](10.5,4)
\rput(3.5,3.6){\makebox(0,0){$\pmb{b_{j-1}^{e}}$}}
\rput(10.5,3.7){\makebox(0,0){$\pmb{0}$}}
\end{pspicture}
\end{center}
\caption{The conformal mapping $\zeta \! = \! \Phi^{e}_{b_{j-1}}(z) \! := \!
(-\tfrac{3n}{2} \int_{z}^{b_{j-1}^{e}}(R_{e}(s))^{1/2}h_{V}^{e}(s))^{2/3}$,
$j \! = \! 1,\dotsc,N \! + \! 1$, where $(\Phi^{e}_{b_{j-1}})^{-1}$ denotes
the inverse mapping}
\end{figure}

Introduce, now, the Airy function, $\operatorname{Ai}(\pmb{\cdot})$, which 
appears in several of the final results of this work: $\operatorname{Ai}
(\pmb{\cdot})$ is determined (uniquely) as the solution of the second-order, 
non-constant coefficient, homogeneous ODE (see, for example, Chapter~10 of 
\cite{a93})
\begin{equation*}
\operatorname{Ai}^{\prime \prime}(z) \! - \! z \operatorname{Ai}(z) \! = \! 0,
\end{equation*}
with asymptotics (at infinity)
\begin{equation}
\begin{split}
\operatorname{Ai}(z) \underset{\underset{\vert \arg z \vert < \pi}{z \to
\infty}}{\sim}& \, \dfrac{1}{2 \sqrt{\smash[b]{\pi}}}z^{-1/4} \, \me^{-
\widehat{\zeta}(z)} \sum_{k=0}^{\infty}(-1)^{k}s_{k}(\widehat{\zeta}(z))^{-k},
\qquad \widehat{\zeta}(z) \! := \! \dfrac{2}{3}z^{3/2}, \\
\operatorname{Ai}^{\prime}(z) \underset{\underset{\vert \arg z \vert < \pi}{z
\to \infty}}{\sim}& \, -\dfrac{1}{2 \sqrt{\smash[b]{\pi}}}z^{1/4} \, \me^{-
\widehat{\zeta}(z)} \sum_{k=0}^{\infty}(-1)^{k}t_{k}(\widehat{\zeta}(z))^{-
k},
\end{split}
\end{equation}
where $s_{0} \! = \! t_{0} \! = \! 1$,
\begin{equation*}
s_{k} \! = \! \dfrac{\Gamma (3k \! + \! 1/2)}{54^{k}k! \Gamma (k \! + \! 1/2)}
\! = \! \dfrac{(2k \! + \! 1)(2k \! + \! 3) \cdots (6k \! - \! 1)}{216^{k}k!},
\qquad \quad t_{k} \! = \! -\left(\dfrac{6k \! + \! 1}{6k \! - \! 1} \right)
\! s_{k}, \quad k \! \in \! \mathbb{N},
\end{equation*}
and $\Gamma (\cdot)$ is the gamma (factorial) function.

In order to present the final asymptotic (as $n \! \to \! \infty)$ results, 
and for arbitrary $j \! = \! 1,\dotsc,N \! + \! 1$, consider the following 
decomposition (see Figure~7), into bounded and unbounded regions, of $\mathbb{
C}$ and the neighbourhoods of the end-points $b_{i-1}^{e},a_{i}^{e}$, $i \! 
= \! 1,\dotsc,N \! + \! 1$ (as per the discussion above, $\mathbb{U}^{e}_{
\delta_{b_{k-1}}} \cap \mathbb{U}^{e}_{\delta_{a_{k}}} \! = \! \varnothing$, 
$k \! = \! 1,\dotsc,N \! + \! 1)$.
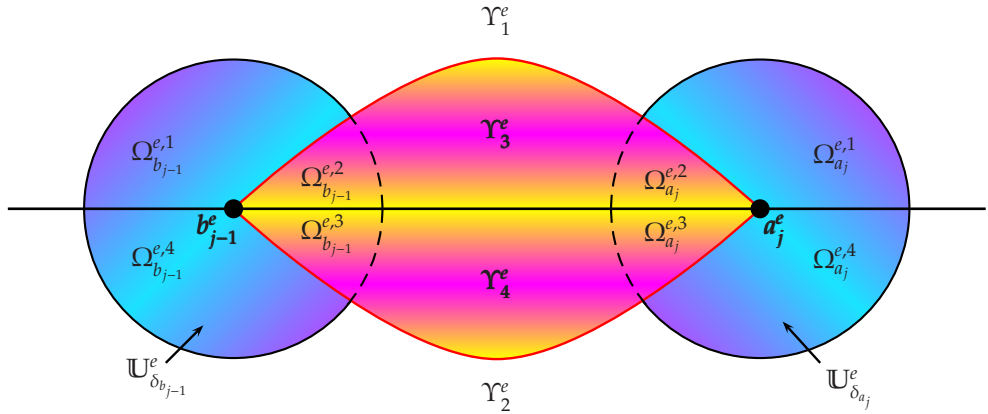
\begin{figure}[tbh]
\begin{center}
\vspace{0.55cm}
\begin{pspicture}(0,0)(14,6)
\psset{xunit=1cm,yunit=1cm}
\pscircle[fillstyle=gradient,gradangle=45,gradlines=600,gradbegin=magenta,%
gradend=cyan,gradmidpoint=0.5](3.5,3){2}
\pscircle[fillstyle=gradient,gradangle=135,gradlines=600,gradbegin=magenta,%
gradend=cyan,gradmidpoint=0.5](10.5,3){2}
\pscurve[linewidth=0.9pt,linestyle=solid,linecolor=red,fillstyle=gradient,%
gradangle=0,gradlines=600,gradbegin=yellow,gradend=magenta,gradmidpoint=0.5]%
(3.5,3)(7,5)(10.5,3)
\pscurve[linewidth=0.9pt,linestyle=solid,linecolor=red,fillstyle=gradient,%
gradangle=0,gradlines=600,gradbegin=yellow,gradend=magenta,gradmidpoint=0.5]%
(3.5,3)(7,1)(10.5,3)
\rput(7,5.5){\makebox(0,0){$\Upsilon^{e}_{1}$}}
\rput(7,0.5){\makebox(0,0){$\Upsilon^{e}_{2}$}}
\rput(7,4){\makebox(0,0){$\pmb{\Upsilon^{e}_{3}}$}}
\rput(7,2){\makebox(0,0){$\pmb{\Upsilon^{e}_{4}}$}}
\rput(2.5,0.75){\makebox(0,0){$\mathbb{U}^{e}_{\delta_{b_{j-1}}}$}}
\rput(11.70,0.60){\makebox(0,0){$\mathbb{U}^{e}_{\delta_{a_{j}}}$}}
\psline[linewidth=0.9pt,linestyle=solid,linecolor=black]{->}(2.6,0.95)%
(3.1,1.45)
\psline[linewidth=0.9pt,linestyle=solid,linecolor=black]{->}(11.3,0.95)%
(10.8,1.50)
\psline[linewidth=1.0pt,linestyle=solid,linecolor=black](0.5,3)(13.5,3)
\rput(2.5,3.7){\makebox(0,0){$\Omega^{e,1}_{b_{j-1}}$}}
\rput(2.5,2.3){\makebox(0,0){$\Omega^{e,4}_{b_{j-1}}$}}
\rput(4.75,3.35){\makebox(0,0){$\Omega^{e,2}_{b_{j-1}}$}}
\rput(4.75,2.65){\makebox(0,0){$\Omega^{e,3}_{b_{j-1}}$}}
\rput(9.25,3.35){\makebox(0,0){$\Omega^{e,2}_{a_{j}}$}}
\rput(9.25,2.65){\makebox(0,0){$\Omega^{e,3}_{a_{j}}$}}
\rput(11.5,3.7){\makebox(0,0){$\Omega^{e,1}_{a_{j}}$}}
\rput(11.5,2.3){\makebox(0,0){$\Omega^{e,4}_{a_{j}}$}}
\psarcn[linewidth=0.8pt,linestyle=dashed,linecolor=black](3.5,3){1.98}{40}{0}
\psarc[linewidth=0.8pt,linestyle=dashed,linecolor=black](3.5,3){1.98}{320}{360}
\psarc[linewidth=0.8pt,linestyle=dashed,linecolor=black](10.5,3){1.98}{140}%
{180}
\psarc[linewidth=0.8pt,linestyle=dashed,linecolor=black](10.5,3){1.98}{180}%
{220}
\psdots[dotstyle=*,dotscale=2](3.5,3)
\psdots[dotstyle=*,dotscale=2](10.5,3)
\rput(3.3,2.7){\makebox(0,0){$\pmb{b_{j-1}^{e}}$}}
\rput(10.7,2.7){\makebox(0,0){$\pmb{a_{j}^{e}}$}}
\end{pspicture}
\end{center}
\caption{Region-by-region decomposition of $\mathbb{C}$ and the neighbourhoods
surrounding the end-points of the support of the `even' equilibrium measure,
$\{b_{j-1}^{e},a_{j}^{e}\}_{j=1}^{N+1}$}
\end{figure}
Asymptotics (as $n \! \to \! \infty)$ for $\boldsymbol{\pi}_{2n}(z)$, with 
$z\! \in \! \cup_{j=1}^{4}(\Upsilon_{j}^{e} \cup (\cup_{k=1}^{N+1}(\Omega^{
e,j}_{b_{k-1}} \cup \Omega^{e,j}_{a_{k}})))$, are now presented. These 
asymptotic expansions are obtained via a union of the DZ non-linear 
steepest-descent method \cite{a1,a2} and the extension of Deift-Venakides-Zhou 
\cite{a3} (see, also, \cite{a57,a58,a59,a60,a61,a62,a63,a64,a65,a66,a67,a68,%
a69,a70,a71,a72,a73,a74,a75,a76,a79,a94}, and the detailed pedagogical 
exposition \cite{a90}).
\begin{eeee}
In order to eschew a flood of superfluous notation, the simplified `notation'
$\mathcal{O}(n^{-2})$ is maintained throughout Theorem~2.3.1 (see below), and
is to be understood in the following, \emph{normal} sense: for a compact
subset, $\mathfrak{D}$, say, of $\mathbb{C}$, and uniformly with respect to
$z \! \in \! \mathfrak{D}$, $\mathcal{O}(n^{-2}) \! := \! \mathcal{O}(c^{
\natural}(z,n)n^{-2})$, where $\norm{c^{\natural}(\boldsymbol{\cdot},n)}_{
\mathcal{L}^{p}(\mathfrak{D})} \! =_{n \to \infty} \! \mathcal{O}(1)$, $p \!
\in \! \lbrace 1,2,\infty \rbrace$, and $\exists \, \, K_{\mathfrak{D}} \! >
\! 0$ (and finite) such that, $\forall \, \, z \! \in \! \mathfrak{D}$, $\vert
c^{\natural}(z,n) \vert \! \leqslant_{n \to \infty} \! K_{\mathfrak{D}}$.
\hfill $\blacksquare$
\end{eeee}
\begin{dddd}
Let the external field $\widetilde{V} \colon \mathbb{R} \setminus \{0\} \! \to
\! \mathbb{R}$ satisfy conditions~{\rm (2.3)--(2.5)}. Set
\begin{equation}
\md \mu_{V}^{e}(x) \! := \! \psi_{V}^{e}(x) \, \md x \! = \! \dfrac{1}{2 \pi
\mi}(R_{e}(x))^{1/2}_{+}h_{V}^{e}(x) \pmb{1}_{J_{e}}(x) \, \md x,
\end{equation}
where
\begin{equation}
(R_{e}(z))^{1/2} \! := \! \left(\prod_{k=1}^{N+1}(z \! - \! b_{k-1}^{e})(z \!
- \! a_{k}^{e}) \right)^{1/2},
\end{equation}
with $(R_{e}(x))^{1/2}_{\pm} \! := \! \lim_{\varepsilon \downarrow 0}(R_{e}(x
\! \pm \! \mi \varepsilon))^{1/2}$, $x \! \in \! J_{e} \! := \! \operatorname{
supp}(\mu_{V}^{e}) \! = \! \cup_{j=1}^{N+1}(b_{j-1}^{e},a_{j}^{e})$ $(\subset
\mathbb{R} \setminus \{0\})$, $N \! \in \! \mathbb{N}$ (and finite), where
$b_{0}^{e} \! := \! \min \lbrace \operatorname{supp}(\mu_{V}^{e}) \rbrace \!
\notin \! \lbrace -\infty,0 \rbrace$, $a_{N+1}^{e} \! := \! \max \lbrace
\operatorname{supp}(\mu_{V}^{e}) \rbrace \! \notin \! \lbrace 0,+\infty
\rbrace$, and $-\infty \! < \! b_{0}^{e} \! < \! a_{1}^{e} \! < \! b_{1}^{e}
\! < \! a_{2}^{e} \! < \! \cdots \! < \! b_{N}^{e} \! < \! a_{N+1}^{e} \! < \!
+\infty$, the branch of the square root is chosen such that $z^{-(N+1)}(R_{e}
(z))^{1/2} \! \sim_{\underset{z \in \mathbb{C}_{\pm}}{z \to \infty}} \! \pm 1$,
\begin{equation}
h_{V}^{e}(z) \! := \! \dfrac{1}{2} \oint_{C_{\mathrm{R}}^{e}} \dfrac{\left(
\frac{\mi}{\pi s} \! + \! \frac{\mi \widetilde{V}^{\prime}(s)}{2 \pi} \right)
}{\sqrt{\smash[b]{R_{e}(s)}} \, (s \! - \! z)} \, \md s
\end{equation}
(real analytic for $z \! \in \! \mathbb{R} \setminus \{0\})$, $C_{\mathrm{R}
}^{e}$ $(\subset \mathbb{C}^{\ast})$ is the boundary of any open
doubly-connected annular region of the type $\lbrace \mathstrut z^{\prime} \!
\in \! \mathbb{C}; \, 0 \! < \! r^{\natural} \! < \! \vert z^{\prime} \vert \!
< \! R^{\natural} \! < \! +\infty \rbrace$, where the simple outer (resp.,
inner) boundary $\lbrace \mathstrut z^{\prime} \! = \! R^{\natural} \me^{\mi
\vartheta}, \, 0 \! \leqslant \! \vartheta \! \leqslant \! 2 \pi \rbrace$
(resp., $\lbrace \mathstrut z^{\prime} \! = \! r^{\natural} \me^{\mi
\vartheta}, \, 0 \! \leqslant \! \vartheta \! \leqslant \! 2 \pi \rbrace)$ is
traversed clockwise (resp., counter-clockwise), with the numbers $0 \! < \!
r^{\natural} \! < \! R^{\natural} \! < \! +\infty$ chosen so that, for (any)
non-real $z$ in the domain of analyticity of $\widetilde{V}$ (that is,
$\mathbb{C}^{\ast})$, $\mathrm{int}(C_{\mathrm{R}}^{e}) \! \supset \! J_{e}
\cup \{z\}$, $\pmb{1}_{J_{e}}(x)$ denotes the indicator (characteristic)
function of the set $J_{e}$, and $\{b_{j-1}^{e},a_{j}^{e}\}_{j=1}^{N+1}$
satisfy the following $n$-dependent and (locally) solvable system of $2(N \!
+ \! 1)$ moment conditions:
\begin{gather}
\begin{split}
\int_{J_{e}} \dfrac{(2s^{-1} \! + \! \widetilde{V}^{\prime}(s))s^{j}}{(R_{e}
(s))^{1/2}_{+}} \, \md s \! = \! 0, \quad j \! = \! 0,\dotsc,N, \qquad \qquad
\int_{J_{e}} \dfrac{(2s^{-1} \! + \! \widetilde{V}^{\prime}(s))s^{N+1}}{(R_{e}
(s))^{1/2}_{+}} \, \md s \! = \! -4 \pi \mi, \\
\int_{a_{j}^{e}}^{b_{j}^{e}} \! \left(\dfrac{\mi (R_{e}(s))^{1/2}}{2 \pi}
\int_{J_{e}} \dfrac{(2 \xi^{-1} \! + \! \widetilde{V}^{\prime}(\xi))}{(R_{e}
(\xi))^{1/2}_{+}(\xi \! - \! s)} \, \md \xi \right) \! \md s \! = \! \ln \!
\left\vert \dfrac{a_{j}^{e}}{b_{j}^{e}} \right\vert \! + \! \dfrac{1}{2} \!
\left(\widetilde{V}(a_{j}^{e}) \! - \! \widetilde{V}(b_{j}^{e}) \right), \quad
j \! = \! 1,\dotsc,N.
\end{split}
\end{gather}
Suppose, furthermore, that $\widetilde{V} \colon \mathbb{R} \setminus \{0\} \!
\to \! \mathbb{R}$ is regular, namely:
\begin{compactenum}
\item[{\rm (i)}] $h_{V}^{e}(x) \! \not\equiv \! 0$ on $\overline{J_{e}} := \!
J_{e} \cup \! \left(\cup_{k=1}^{N+1} \lbrace b_{k-1}^{e},a_{k}^{e} \rbrace
\right);$
\item[{\rm (ii)}]
\begin{equation}
4 \int_{J_{e}} \ln (\vert x \! - \! s \vert) \, \md \mu_{V}^{e}(s) \! - \! 2
\ln \vert x \vert \! - \! \widetilde{V}(x) \! - \! \ell_{e} \! = \! 0, \quad x
\! \in \! \overline{J_{e}},
\end{equation}
which defines the `even' variational constant $\ell_{e} \! \in \! \mathbb{R}$
(the same on each---compact---interval $[b_{j-1}^{e},a_{j}^{e}]$, $j \! = \! 
1,\dotsc,N \! + \! 1)$, and
\begin{equation*}
4 \int_{J_{e}} \ln (\vert x \! - \! s \vert) \, \md \mu_{V}^{e}(s) \! - \! 2
\ln \vert x \vert \! - \! \widetilde{V}(x) \! - \! \ell_{e} \! < \! 0, \quad x
\! \in \! \mathbb{R} \setminus \overline{J_{e}};
\end{equation*}
\item[{\rm (iii)}]
\begin{equation*}
g^{e}_{+}(z) \! + \! g^{e}_{-}(z) \! - \! \widetilde{V}(z) \! - \! \ell_{e} \!
+ \! 2Q_{e} \! < \! 0, \quad z \! \in \! \mathbb{R} \setminus \overline{J_{e}},
\end{equation*}
where
\begin{equation}
g^{e}(z) \! := \! \int_{J_{e}} \ln \! \left((z \! - \! s)^{2}(zs)^{-1} \right)
\md \mu_{V}^{e}(s), \quad z \! \in \! \mathbb{C} \setminus (-\infty,\max \{0,
a_{N+1}^{e}\}),
\end{equation}
with
\begin{align}
Q_{e} :=& \, \int_{J_{e}} \ln (s) \, \md \mu_{V}^{e}(s) \! = \! \int_{J_{e}}
\ln (\lvert s \rvert) \, \md \mu_{V}^{e}(s) \! + \! \mi \pi \int_{J_{e} \cap
\mathbb{R}_{-}} \md \mu_{V}^{e}(s) \nonumber \\
=& \, \int_{J_{e}} \ln (\lvert s \rvert) \, \md \mu_{V}^{e}(s) \! + \! \mi \pi
\begin{cases}
0, &\text{$J_{e} \! \subset \! \mathbb{R}_{+}$,} \\
1, &\text{$J_{e} \! \subset \! \mathbb{R}_{-}$,} \\
\int_{b_{0}^{e}}^{a_{j}^{e}} \md \mu_{V}^{e}(s), &\text{$(a_{j}^{e},b_{j}^{e})
\! \ni \! 0, \quad j \! = \! 1,\dotsc,N$;}
\end{cases}
\end{align}
\item[{\rm (iv)}]
\begin{equation*}
\mi (g^{e}_{+}(z) \! - \! g^{e}_{-}(z))^{\prime} \! > \! 0, \quad z \! \in \!
J_{e}.
\end{equation*}
\end{compactenum}

Set
\begin{equation}
\overset{e}{m}^{\raise-1.0ex\hbox{$\scriptstyle \infty$}}(z) \! = \!
\begin{cases}
\overset{e}{\mathfrak{M}}^{\raise-1.0ex\hbox{$\scriptstyle \infty$}}(z),
&\text{$z \! \in \! \mathbb{C}_{+}$,} \\
-\mi \, \overset{e}{\mathfrak{M}}^{\raise-1.0ex\hbox{$\scriptstyle \infty$}}
(z) \sigma_{2}, &\text{$z \in \! \mathbb{C}_{-}$,}
\end{cases}
\end{equation}
where $(\det (\overset{e}{m}^{\raise-1.0ex\hbox{$\scriptstyle \infty$}}(z))
\! = \! 1)$
\begin{gather}
\overset{e}{\mathfrak{M}}^{\raise-1.0ex\hbox{$\scriptstyle \infty$}}(z) \!
= \!
\begin{pmatrix}
\frac{(\gamma^{e}(z)+(\gamma^{e}(z))^{-1})}{2} \mathfrak{m}^{e}_{11}(z) &
-\frac{(\gamma^{e}(z)-(\gamma^{e}(z))^{-1})}{2 \mi} \mathfrak{m}^{e}_{12}(z)
\\
\frac{(\gamma^{e}(z)-(\gamma^{e}(z))^{-1})}{2 \mi} \mathfrak{m}^{e}_{21}(z)
& \frac{(\gamma^{e}(z)+(\gamma^{e}(z))^{-1})}{2} \mathfrak{m}^{e}_{22}(z)
\end{pmatrix}, \\
\gamma^{e}(z) \! := \! \left(\! \left(\dfrac{z \! - \! b_{0}^{e}}{z \! - \! a_{
N+1}^{e}} \right) \! \prod_{k=1}^{N} \! \left(\dfrac{z \! - \! b_{k}^{e}}{z \!
- \! a_{k}^{e}} \right) \right)^{1/4}, \\
\mathfrak{m}^{e}_{11}(z) \! := \! \dfrac{\boldsymbol{\theta}^{e}(\boldsymbol{
u}^{e}_{+}(\infty) \! + \! \boldsymbol{d}_{e}) \boldsymbol{\theta}^{e}
(\boldsymbol{u}^{e}(z) \! - \! \frac{n}{2 \pi} \boldsymbol{\Omega}^{e} \! + \!
\boldsymbol{d}_{e})}{\boldsymbol{\theta}^{e}(\boldsymbol{u}^{e}_{+}(\infty) \!
- \! \frac{n}{2 \pi} \boldsymbol{\Omega}^{e} \! + \! \boldsymbol{d}_{e})
\boldsymbol{\theta}^{e}(\boldsymbol{u}^{e}(z) \! + \! \boldsymbol{d}_{e})}, \\
\mathfrak{m}^{e}_{12}(z) \! := \! \dfrac{\boldsymbol{\theta}^{e}(\boldsymbol{
u}^{e}_{+}(\infty) \! + \! \boldsymbol{d}_{e}) \boldsymbol{\theta}^{e}(-
\boldsymbol{u}^{e}(z) \! - \! \frac{n}{2 \pi} \boldsymbol{\Omega}^{e} \! + \!
\boldsymbol{d}_{e})}{\boldsymbol{\theta}^{e}(\boldsymbol{u}^{e}_{+}(\infty) \!
- \! \frac{n}{2 \pi} \boldsymbol{\Omega}^{e} \! + \! \boldsymbol{d}_{e})
\boldsymbol{\theta}^{e}(-\boldsymbol{u}^{e}(z) \! + \! \boldsymbol{d}_{e})}, \\
\mathfrak{m}^{e}_{21}(z) \! := \! \dfrac{\boldsymbol{\theta}^{e}(\boldsymbol{
u}^{e}_{+}(\infty) \! + \! \boldsymbol{d}_{e}) \boldsymbol{\theta}^{e}
(\boldsymbol{u}^{e}(z) \! - \! \frac{n}{2 \pi} \boldsymbol{\Omega}^{e} \! - \!
\boldsymbol{d}_{e})}{\boldsymbol{\theta}^{e}(-\boldsymbol{u}^{e}_{+}(\infty)
\! - \! \frac{n}{2 \pi} \boldsymbol{\Omega}^{e} \! - \! \boldsymbol{d}_{e})
\boldsymbol{\theta}^{e}(\boldsymbol{u}^{e}(z) \! - \! \boldsymbol{d}_{e})}, \\
\mathfrak{m}^{e}_{22}(z) \! := \! \dfrac{\boldsymbol{\theta}^{e}(\boldsymbol{
u}^{e}_{+}(\infty) \! + \! \boldsymbol{d}_{e}) \boldsymbol{\theta}^{e}(-
\boldsymbol{u}^{e}(z) \! - \! \frac{n}{2 \pi} \boldsymbol{\Omega}^{e} \! - \!
\boldsymbol{d}_{e})}{\boldsymbol{\theta}^{e}(-\boldsymbol{u}^{e}_{+}(\infty)
\! - \! \frac{n}{2 \pi} \boldsymbol{\Omega}^{e} \! - \! \boldsymbol{d}_{e})
\boldsymbol{\theta}^{e}(\boldsymbol{u}^{e}(z) \! + \! \boldsymbol{d}_{e})},
\end{gather}
with
\begin{equation*}
\boldsymbol{u}^{e}(z) \! = \! \int_{a_{N+1}^{e}}^{z} \boldsymbol{\omega}^{e},
\qquad \qquad \boldsymbol{u}^{e}_{+}(\infty) \! = \! \int_{a_{N+1}^{e}}^{
\infty^{+}} \boldsymbol{\omega}^{e},
\end{equation*}
$\boldsymbol{\Omega}^{e} \! = \! (\Omega^{e}_{1},\Omega^{e}_{2},\dotsc,
\Omega^{e}_{N})^{\operatorname{T}}$ $(\in \! \mathbb{R}^{N})$, where
\begin{equation*}
\Omega^{e}_{j} \! := \! 4 \pi \int_{b_{j}^{e}}^{a_{N+1}^{e}} \psi_{V}^{e}(s)
\, \md s, \quad j \! = \! 1,\dotsc,N,
\end{equation*}
and
\begin{equation*}
\boldsymbol{d}_{e} \! \equiv \! \sum_{j=1}^{N} \int_{a_{j}^{e}}^{z_{j}^{e,+}}
\boldsymbol{\omega}^{e} \quad \left(\equiv \! -\sum_{j=1}^{N+1} \int_{a_{j}^{
e}}^{z_{j}^{e,-}} \boldsymbol{\omega}^{e} \right),
\end{equation*}
where
\begin{equation*}
\left\lbrace z_{j}^{e,\pm} \right\rbrace_{j=1}^{N} \! = \! \left\lbrace
\mathstrut z^{\pm} \! \in \! \mathbb{C}_{\pm}; \, (\gamma^{e}(z) \! \mp \!
(\gamma^{e}(z))^{-1}) \vert_{z=z^{\pm}} \! = \! 0 \right\rbrace,
\end{equation*}
with $z_{j}^{e,\pm} \! \in \! (a_{j}^{e},b_{j}^{e})^{\pm}$ $(\subset \!
\mathbb{C}_{\pm})$, $j \! = \! 1,\dotsc,N$.

Let $\overset{e}{\operatorname{Y}} \colon \mathbb{C} \setminus \mathbb{R}
\! \to \! \operatorname{SL}_{2}(\mathbb{C})$ be the unique solution of
{\rm \pmb{RHP1}} whose integral representations are given in
Lemma~{\rm 2.2.1;} in particular, $\boldsymbol{\pi}_{2n}(z) \! := \!
(\overset{e}{\mathrm{Y}}(z))_{11}$. Then:\\
{\rm \pmb{(1)}} for $z \! \in \! \Upsilon^{e}_{1}$ $(\subset \! \mathbb{C}_{
+})$,
\begin{align}
\boldsymbol{\pi}_{2n}(z) \underset{n \to \infty}{=}& \, \exp \! \left(n(g^{e}
(z) \! + \! Q_{e}) \right) \!
\left((\overset{e}{m}^{\raise-1.0ex\hbox{$\scriptstyle \infty$}}(z))_{11} \!
\left(1 \! + \! \dfrac{1}{n} \! \left(\mathscr{R}^{e}_{\infty}(z) \right)_{11}
\! + \! \mathcal{O} \! \left(\dfrac{1}{n^{2}} \right) \right) \right.
\nonumber \\
+&\left. \, (\overset{e}{m}^{\raise-1.0ex\hbox{$\scriptstyle \infty$}}(z))_{2
1} \! \left(\dfrac{1}{n} \! \left(\mathscr{R}^{e}_{\infty}(z) \right)_{12} \!
+ \! \mathcal{O} \! \left(\dfrac{1}{n^{2}} \right) \right) \right),
\end{align}
and
\begin{align}
\int_{\mathbb{R}} \dfrac{\boldsymbol{\pi}_{2n}(s) \me^{-n \widetilde{V}(s)}}{s
\! - \! z} \, \dfrac{\md s}{2 \pi \mi} \underset{n \to \infty}{=}& \, \exp \!
\left(-n(g^{e}(z) \! - \! \ell_{e} \! + \! Q_{e}) \right) \! \left(\!
(\overset{e}{m}^{\raise-1.0ex\hbox{$\scriptstyle \infty$}}(z))_{12} \! \left(
\! 1 \! + \! \dfrac{1}{n} \! \left(\mathscr{R}^{e}_{\infty}(z) \right)_{11}
\right. \right. \nonumber \\
+&\left. \left. \mathcal{O} \! \left(\dfrac{1}{n^{2}} \right) \right) \! +
\! (\overset{e}{m}^{\raise-1.0ex\hbox{$\scriptstyle \infty$}}(z))_{22} \!
\left(\dfrac{1}{n} \! \left(\mathscr{R}^{e}_{\infty}(z) \right)_{12} \! + \!
\mathcal{O} \! \left(\dfrac{1}{n^{2}} \right) \right) \right),
\end{align}
where
\begin{align}
\mathscr{R}^{e}_{\infty}(z):=& \, \sum_{j=1}^{N+1} \dfrac{1}{(z \! - \! b_{j-
1}^{e})} \! \left(\dfrac{\mathscr{A}^{e}(b_{j-1}^{e})}{\widehat{\alpha}_{0}^{e}
(b_{j-1}^{e})(z \! - \! b_{j-1}^{e})} \! + \! \dfrac{(\mathscr{B}^{e}(b_{j-1}^{
e}) \widehat{\alpha}_{0}^{e}(b_{j-1}^{e}) \! - \! \mathscr{A}^{e}(b_{j-1}^{e})
\widehat{\alpha}_{1}^{e}(b_{j-1}^{e}))}{(\widehat{\alpha}_{0}^{e}(b_{j-1}^{e}
))^{2}} \right) \nonumber \\
+& \, \sum_{j=1}^{N+1} \dfrac{1}{(z \! - \! a_{j}^{e})} \! \left(\dfrac{
\mathscr{A}^{e}(a_{j}^{e})}{\widehat{\alpha}_{0}^{e}(a_{j}^{e})(z \! - \! a_{
j}^{e})} \! + \! \dfrac{(\mathscr{B}^{e}(a_{j}^{e}) \widehat{\alpha}_{0}^{e}
(a_{j}^{e}) \! - \! \mathscr{A}^{e}(a_{j}^{e}) \widehat{\alpha}_{1}^{e}(a_{j}^{
e}))}{(\widehat{\alpha}_{0}^{e}(a_{j}^{e}))^{2}} \right),
\end{align}
with, for $j \! = \! 1,\dotsc,N \! + \! 1$,
\begin{gather}
\mathscr{A}^{e}(b_{j-1}^{e}) \! = \! -s_{1}(Q_{0}^{e}(b_{j-1}^{e}))^{-1} 
\me^{\mi n \mho_{j-1}^{e}} \! 
\begin{pmatrix}
\varkappa^{e}_{1}(b_{j-1}^{e}) \varkappa^{e}_{2}(b_{j-1}^{e}) & \mi 
(\varkappa^{e}_{1}(b_{j-1}^{e}))^{2} \\
\mi (\varkappa^{e}_{2}(b_{j-1}^{e}))^{2} & -\varkappa^{e}_{1}(b_{j-1}^{e}) 
\varkappa^{e}_{2}(b_{j-1}^{e})
\end{pmatrix}, \\
\mathscr{A}^{e}(a_{j}^{e}) \! = \! s_{1}Q_{0}^{e}(a_{j}^{e}) \me^{\mi n 
\mho_{j}^{e}} \!
\begin{pmatrix}
-\varkappa^{e}_{1}(a_{j}^{e}) \varkappa^{e}_{2}(a_{j}^{e}) & \mi 
(\varkappa^{e}_{1}(a_{j}^{e}))^{2} \\
\mi (\varkappa^{e}_{2}(a_{j}^{e}))^{2} & \varkappa^{e}_{1}(a_{j}^{e}) 
\varkappa^{e}_{2}(a_{j}^{e})
\end{pmatrix}, \\
\dfrac{\mathscr{B}^{e}(b_{j-1}^{e})}{\me^{\mi n \mho_{j-1}^{e}}} \! = \! 
\begin{pmatrix}
\boxed{\begin{matrix} \varkappa_{1}^{e}(b_{j-1}^{e}) \varkappa_{2}^{e}
(b_{j-1}^{e}) \! \left(-s_{1}(Q_{0}^{e}(b_{j-1}^{e}))^{-1} \right. \\
\left. \times \left\{\daleth^{1}_{1}(b_{j-1}^{e}) \! + \! \daleth^{1}_{-1}
(b_{j-1}^{e}) \! - \! Q_{1}^{e}(b_{j-1}^{e}) \right. \right. \\
\left. \left. \times \, (Q_{0}^{e}(b_{j-1}^{e}))^{-1} \right\} \! - \! t_{1} 
\! \left\{Q_{0}^{e}(b_{j-1}^{e}) \right. \right. \\
\left. \left. + \, (Q_{0}^{e}(b_{j-1}^{e}))^{-1} \aleph^{1}_{1}(b_{j-1}^{e}) 
\aleph^{1}_{-1} (b_{j-1}^{e}) \right\} \right. \\
\left. + \, \mi (s_{1} \! + \! t_{1}) \! \left\{\aleph^{1}_{-1}(b_{j-1}^{e}) 
\! - \! \aleph^{1}_{1}(b_{j-1}^{e}) \right\} \right)
\end{matrix}} & 
\boxed{\begin{matrix} (\varkappa_{1}^{e}(b_{j-1}^{e}))^{2} \! \left(-\mi 
s_{1}(Q_{0}^{e}(b_{j-1}^{e}))^{-1} \! \left\{2 \daleth^{1}_{1}(b_{j-1}^{e}) 
\right. \right. \\
\left. \left. -\, Q_{1}^{e}(b_{j-1}^{e})(Q_{0}^{e}(b_{j-1}^{e}))^{-1} \right\} 
\! + \! \mi t_{1} \! \left\{Q_{0}^{e}(b_{j-1}^{e}) \right. \right. \\
\left. \left. -\, (Q_{0}^{e}(b_{j-1}^{e}))^{-1}(\aleph^{1}_{1}(b_{j-1}^{e})
)^{2} \right\} \right. \\
\left. + \, 2(s_{1} \! - \! t_{1}) \aleph^{1}_{1}(b_{j-1}^{e}) \right)
\end{matrix}} \\
\boxed{\begin{matrix} (\varkappa_{2}^{e}(b_{j-1}^{e}))^{2} \! \left(-\mi 
s_{1}(Q_{0}^{e}(b_{j-1}^{e}))^{-1} \! \left\{2 \daleth^{1}_{-1}
(b_{j-1}^{e}) \right. \right. \\
\left. \left. -\, Q_{1}^{e}(b_{j-1}^{e})(Q_{0}^{e}(b_{j-1}^{e}))^{-1} \right\} 
\! + \! \mi t_{1} \! \left\{Q_{0}^{e}(b_{j-1}^{e}) \right. \right. \\
\left. \left. -\, (Q_{0}^{e}(b_{j-1}^{e}))^{-1}(\aleph^{1}_{-1}(b_{j-1}^{e}
))^{2} \right\} \right. \\
\left. - \, 2(s_{1} \! - \! t_{1}) \aleph^{1}_{-1}(b_{j-1}^{e}) \right)
\end{matrix}} & 
\boxed{\begin{matrix} \varkappa_{1}^{e}(b_{j-1}^{e}) \varkappa_{2}^{e}
(b_{j-1}^{e}) \! \left(s_{1}(Q_{0}^{e}(b_{j-1}^{e}))^{-1} \right. \\
\left. \times \left\{\daleth^{1}_{1}(b_{j-1}^{e}) \! + \! \daleth^{1}_{-1}
(b_{j-1}^{e}) \! - \! Q_{1}^{e}(b_{j-1}^{e}) \right. \right. \\
\left. \left. \times \, (Q_{0}^{e}(b_{j-1}^{e}))^{-1} \right\} \! + \! t_{1} 
\! \left\{Q_{0}^{e}(b_{j-1}^{e}) \right. \right. \\
\left. \left. + \, (Q_{0}^{e}(b_{j-1}^{e}))^{-1} \aleph^{1}_{1}(b_{j-1}^{e}) 
\aleph^{1}_{-1} (b_{j-1}^{e}) \right\} \right. \\
\left. + \, \mi (s_{1} \! + \! t_{1}) \! \left\{\aleph^{1}_{1}(b_{j-1}^{e}) 
\! - \! \aleph^{1}_{-1}(b_{j-1}^{e}) \right\} \right)
\end{matrix}}
\end{pmatrix}, \\
\dfrac{\mathscr{B}^{e}(a_{j}^{e})}{\me^{\mi n \mho_{j}^{e}}} \! = \! 
\begin{pmatrix}
\boxed{\begin{matrix} \varkappa_{1}^{e}(a_{j}^{e}) \varkappa_{2}^{e}
(a_{j}^{e}) \! \left(-s_{1} \! \left\{Q_{1}^{e}(a_{j}^{e}) \right. \right. \\
\left. \left. + \, Q_{0}^{e}(a_{j}^{e}) \left[\daleth^{1}_{1}(a_{j}^{e}) \! 
+ \! \daleth^{1}_{-1}(a_{j}^{e}) \right] \right\} \! - \! t_{1} \right. \\
\left. \times \left\{(Q_{0}^{e}(a_{j}^{e}))^{-1} \! + \! Q_{0}^{e}(a_{j}^{e}) 
\aleph^{1}_{1}(a_{j}^{e}) \aleph^{1}_{-1}(a_{j}^{e}) \right\} \right. \\
\left. + \, \mi (s_{1} \! + \! t_{1}) \! \left\{\aleph^{1}_{-1}(a_{j}^{e}) 
\! - \! \aleph^{1}_{1}(a_{j}^{e}) \right\} \right)
\end{matrix}} & 
\boxed{\begin{matrix} (\varkappa_{1}^{e}(a_{j}^{e}))^{2} \! \left(\mi s_{1} 
\! \left\{Q_{1}^{e}(a_{j}^{e}) \! + \! 2Q_{0}^{e}(a_{j}^{e}) \right. \right. 
\\
\left. \left. \times \, \daleth^{1}_{1}(a_{j}^{e}) \right\} \! + \! \mi t_{1} 
\! \left\{Q_{0}^{e}(a_{j}^{e})(\aleph^{1}_{1}(a_{j}^{e}))^{2} \right. \right. 
\\
\left. \left. - \, (Q_{0}^{e}(a_{j}^{e}))^{-1} \right\} \! - \! 2(s_{1} \! - 
\! t_{1}) \aleph^{1}_{1}(a_{j}^{e}) \right)
\end{matrix}} \\
\boxed{\begin{matrix} (\varkappa_{2}^{e}(a_{j}^{e}))^{2} \! \left(\mi s_{1} 
\! \left\{Q_{1}^{e}(a_{j}^{e}) \! + \! 2Q_{0}^{e}(a_{j}^{e}) \right. \right. 
\\
\left. \left. \times \, \daleth^{1}_{-1}(a_{j}^{e}) \right\} \! + \! \mi t_{1} 
\! \left\{Q_{0}^{e}(a_{j}^{e})(\aleph^{1}_{-1}(a_{j}^{e}))^{2} \right. 
\right. \\
\left. \left. - \, (Q_{0}^{e}(a_{j}^{e}))^{-1} \right\} \! + \! 2(s_{1} \! - 
\! t_{1}) \aleph^{1}_{-1}(a_{j}^{e}) \right)
\end{matrix}} & 
\boxed{\begin{matrix} \varkappa_{1}^{e}(a_{j}^{e}) \varkappa_{2}^{e}
(a_{j}^{e}) \! \left(s_{1} \! \left\{Q_{1}^{e}(a_{j}^{e}) \right. \right. \\
\left. \left. + \, Q_{0}^{e}(a_{j}^{e}) \left[\daleth^{1}_{1}(a_{j}^{e}) \! 
+ \! \daleth^{1}_{-1}(a_{j}^{e}) \right] \right\} \! + \! t_{1} \right. \\
\left. \times \left\{(Q_{0}^{e}(a_{j}^{e}))^{-1} \! + \! Q_{0}^{e}(a_{j}^{e}) 
\aleph^{1}_{1}(a_{j}^{e}) \aleph^{1}_{-1}(a_{j}^{e}) \right\} \right. \\
\left. + \, \mi (s_{1} \! + \! t_{1}) \! \left\{\aleph^{1}_{1}(a_{j}^{e}) \! 
- \! \aleph^{1}_{-1}(a_{j}^{e}) \right\} \right)
\end{matrix}}
\end{pmatrix}, \\
s_{1} \! = \! \dfrac{5}{72}, \qquad \qquad \quad t_{1} \! = \! -\dfrac{7}{72}, 
\qquad \qquad \quad \mho_{i}^{e} \! := \!
\begin{cases}
\Omega_{i}^{e}, &\text{$i \! = \! 1,\dotsc,N$,} \\
0, &\text{$i \! = \! 0,N \! + \! 1$,}
\end{cases} \\
Q_{0}^{e}(b_{0}^{e}) \! = \! -\mi \left(\! (a_{N+1}^{e} \! - \! b_{0}^{e})^{
-1} \prod_{k=1}^{N} \! \left(\dfrac{b_{k}^{e} \! - \! b_{0}^{e}}{a_{k}^{e} 
\! - \! b_{0}^{e}} \right) \right)^{1/2}, \\
Q_{1}^{e}(b_{0}^{e}) \! = \! \dfrac{1}{2} Q_{0}^{e}(b_{0}^{e}) \! \left(
\sum_{ k=1}^{N} \! \left(\dfrac{1}{b_{0}^{e} \! - \! b_{k}^{e}} \! - \! 
\dfrac{1}{b_{ 0}^{e} \! - \! a_{k}^{e}} \right) \! - \! \dfrac{1}{b_{0}^{e} 
\! - \! a_{N+1}^{e}} \right), \\
Q_{0}^{e}(a_{N+1}^{e}) \! = \! \left(\! (a_{N+1}^{e} \! - \! b_{0}^{e}) 
\prod_{k=1}^{N} \! \left(\dfrac{a_{N+1}^{e} \! - \! b_{k}^{e}}{a_{N+1}^{e} \! 
- \! a_{k}^{e}} \right) \right)^{1/2}, \\
Q_{1}^{e}(a_{N+1}^{e}) \! = \! \dfrac{1}{2}Q_{0}^{e}(a_{N+1}^{e}) \! \left(
\sum_{k=1}^{N} \! \left(\dfrac{1}{a_{N+1}^{e} \! - \! b_{k}^{e}} \! - \!
\dfrac{1}{a_{N+1}^{e} \! - \! a_{k}^{e}} \right) \! + \! \dfrac{1}{a_{N+1}^{e}
\! - \! b_{0}^{e}} \right), \\
Q_{0}^{e}(b_{j}^{e}) \! = \! -\mi \left(\! \dfrac{(b_{j}^{e} \! - \! b_{0}^{e}
)}{(a_{N+1}^{e} \! - \! b_{j}^{e})(b_{j}^{e} \! - \! a_{j}^{e})} \prod_{k=1}^{
j-1} \! \left(\dfrac{b_{j}^{e} \! - \! b_{k}^{e}}{b_{j}^{e} \! - \! a_{k}^{e}}
\right) \! \prod_{l=j+1}^{N} \! \left(\! \dfrac{b_{l}^{e} \! - \! b_{j}^{e}}{
a_{l}^{e} \! - \! b_{j}^{e}} \right) \right)^{1/2}, \\
Q_{1}^{e}(b_{j}^{e}) \! = \! \dfrac{1}{2} Q_{0}^{e}(b_{j}^{e}) \! \left(\sum_{
\substack{k=1\\k \not= j}}^{N} \! \left(\! \dfrac{1}{b_{j}^{e} \! - \! b_{k}^{
e}} \! - \! \dfrac{1}{b_{j}^{e} \! - \! a_{k}^{e}} \right) \! + \! \dfrac{1}{
b_{j}^{e} \! - \! b_{0}^{e}} \! - \! \dfrac{1}{b_{j}^{e} \! - \! a_{N+1}^{e}}
\! - \! \dfrac{1}{b_{j}^{e} \! - \! a_{j}^{e}} \right), \\
Q_{0}^{e}(a_{j}^{e}) \! = \! \left(\! \dfrac{(a_{j}^{e} \! - \! b_{0}^{e})(b_{
j}^{e} \! - \! a_{j}^{e})}{(a_{N+1}^{e} \! - \! a_{j}^{e})} \prod_{k=1}^{j-1}
\! \left(\dfrac{a_{j}^{e} \! - \! b_{k}^{e}}{a_{j}^{e} \! - \! a_{k}^{e}}
\right) \! \prod_{l=j+1}^{N} \! \left(\! \dfrac{b_{l}^{e} \! - \! a_{j}^{e}}{
a_{l}^{e} \! - \! a_{j}^{e}} \right) \right)^{1/2}, \\
Q_{1}^{e}(a_{j}^{e}) \! = \! \dfrac{1}{2}Q_{0}^{e}(a_{j}^{e}) \! \left(
\sum_{\substack{k=1\\k \not= j}}^{N} \! \left(\! \dfrac{1}{a_{j}^{e} \! - \!
b_{k}^{e}} \! - \! \dfrac{1}{a_{j}^{e} \! - \! a_{k}^{e}} \right) \! + \!
\dfrac{1}{a_{j}^{e} \! - \! b_{0}^{e}} \! - \! \dfrac{1}{a_{j}^{e} \! - \!
a_{N+1}^{e}} \! + \! \dfrac{1}{a_{j}^{e} \! - \! b_{j}^{e}} \right),
\end{gather}
where $\mi Q_{0}^{e}(b_{j-1}^{e}),Q_{0}^{e}(a_{j}^{e}) \! > \! 0$, $j \! = \! 
1,\dotsc,N \! + \! 1$,
\begin{gather}
\varkappa_{1}^{e}(\xi) \! = \! \dfrac{\bm{\theta}^{e}(\bm{u}^{e}_{+}(\infty) 
\! + \! \bm{d}_{e}) \bm{\theta}^{e}(\bm{u}^{e}_{+}(\xi) \! - \! \frac{n}{2 
\pi} \bm{\Omega}^{e} \! + \! \bm{d}_{e})}{\bm{\theta}^{e}(\bm{u}^{e}_{+}
(\infty) \! - \! \frac{n}{2 \pi} \bm{\Omega}^{e} \! + \! \bm{d}_{e}) \bm{
\theta}^{e}(\bm{u}^{e}_{+}(\xi) \! + \! \bm{d}_{e})}, \\
\varkappa_{2}^{e}(\xi) \! = \! \dfrac{\bm{\theta}^{e}(-\bm{u}^{e}_{+}(\infty) 
\! - \! \bm{d}_{e}) \bm{\theta}^{e}(\bm{u}^{e}_{+}(\xi) \! - \! \frac{n}{2 
\pi} \bm{\Omega}^{e} \! - \! \bm{d}_{e})}{\bm{\theta}^{e}(-\bm{u}^{e}_{+}
(\infty) \! - \! \frac{n}{2 \pi} \bm{\Omega}^{e} \! - \! \bm{d}_{e}) \bm{
\theta}^{e}(\bm{u}^{e}_{+}(\xi) \! - \! \bm{d}_{e})}, \\
\aleph^{\varepsilon_{1}}_{\varepsilon_{2}}(\xi) \! = \! -\dfrac{\mathfrak{
u}^{e}(\varepsilon_{1},\varepsilon_{2},\bm{0};\xi)}{\bm{\theta}^{e}
(\varepsilon_{1} \bm{u}^{e}_{+}(\xi) \! + \! \varepsilon_{2} \bm{d}_{e})} \! 
+ \! \dfrac{\mathfrak{u}^{e}(\varepsilon_{1},\varepsilon_{2},\bm{\Omega}^{e};
\xi)}{\bm{\theta}^{e}(\varepsilon_{1} \bm{u}^{e}_{+}(\xi) \! - \! \frac{
n}{2 \pi} \bm{\Omega}^{e} \! + \! \varepsilon_{2} \bm{d}_{e})}, \quad 
\varepsilon_{1},\varepsilon_{2} \! = \! \pm 1,
\end{gather}
\begin{align}
\daleth^{\varepsilon_{1}}_{\varepsilon_{2}}(\xi) =& \, -\dfrac{\mathfrak{
v}^{e}(\varepsilon_{1},\varepsilon_{2},\bm{0};\xi)}{\bm{\theta}^{e}
(\varepsilon_{1} \bm{u}^{e}_{+}(\xi) \! + \! \varepsilon_{2} \bm{d}_{e})} \! 
+ \! \dfrac{\mathfrak{v}^{e}(\varepsilon_{1},\varepsilon_{2},\bm{\Omega}^{e};
\xi)}{\bm{\theta}^{e}(\varepsilon_{1} \bm{u}^{e}_{+}(\xi) \! - \! \frac{n}{2 
\pi} \bm{\Omega}^{e} \! + \! \varepsilon_{2} \bm{d}_{e})} \! - \! \left(
\dfrac{\mathfrak{u}^{e}(\varepsilon_{1},\varepsilon_{2},\bm{0};\xi)}{\bm{
\theta}^{e}(\varepsilon_{1} \bm{u}^{e}_{+}(\xi) \! + \! \varepsilon_{2} 
\bm{d}_{e})} \right)^{2} \nonumber \\
+& \, \dfrac{\mathfrak{u}^{e}(\varepsilon_{1},\varepsilon_{2},\bm{0};\xi) 
\mathfrak{u}^{e}(\varepsilon_{1},\varepsilon_{2},\bm{\Omega}^{e};\xi)}{\bm{
\theta}^{e}(\varepsilon_{1} \bm{u}^{e}_{+}(\xi) \! + \! \varepsilon_{2} 
\bm{d}_{e}) \bm{\theta}^{e}(\varepsilon_{ 1} \bm{u}^{e}_{+}(\xi) \! - \! 
\frac{n}{2 \pi} \bm{\Omega}^{e} \! + \! \varepsilon_{2} \bm{d}_{e})},
\end{align}
\begin{gather}
\mathfrak{u}^{e}(\varepsilon_{1},\varepsilon_{2},\bm{\Omega}^{e},\xi) \! := 
\! 2 \pi \Lambda^{\raise-1.0ex\hbox{$\scriptstyle 1$}}_{e}(\varepsilon_{1},
\varepsilon_{2},\bm{\Omega}^{e},\xi), \qquad \mathfrak{v}^{e}(\varepsilon_{
1},\varepsilon_{2},\bm{\Omega}^{e},\xi) \! := \! -2 \pi^{2} 
\Lambda^{\raise-1.0ex\hbox{$\scriptstyle 2$}}_{e}(\varepsilon_{1},
\varepsilon_{2},\bm{\Omega}^{e},\xi), \\
\Lambda^{\raise-1.0ex\hbox{$\scriptstyle j^{\prime}$}}_{e}(\varepsilon_{1},
\varepsilon_{2},\bm{\Omega}^{e},\xi) \! = \! \sum_{m \in \mathbb{Z}^{N}}
(\mathfrak{r}_{e}(\xi))^{j^{\prime}} \me^{2 \pi \mi (m,\varepsilon_{1} 
\bm{u}^{e}_{+}(\xi)-\frac{n}{2 \pi} \bm{\Omega}^{e} \! + \! \varepsilon_{2} 
\bm{d}_{e})+ \pi \mi (m,\bm{\tau}^{e}m)}, \quad j^{\prime} \! = \! 1,2, \\
\mathfrak{r}_{e}(\xi) \! := \! \dfrac{2(m,\vec{\moo}_{e}(\xi))}{
\leftthreetimes^{\raise+0.3ex\hbox{$\scriptstyle e$}}(\xi)}, \qquad \qquad 
\vec{\moo}_{e}(\xi) \! = \! 
(\rightthreetimes^{\raise-0.9ex\hbox{$\scriptstyle e$}}_{1}(\xi),
\rightthreetimes^{\raise-0.9ex\hbox{$\scriptstyle e$}}_{2}(\xi),
\dotsc,\rightthreetimes^{\raise-0.9ex\hbox{$\scriptstyle e$}}_{N}(\xi)), \\
\rightthreetimes^{\raise-0.9ex\hbox{$\scriptstyle e$}}_{j^{\prime}}(\xi) \! 
:= \sum_{k=1}^{N} c_{j^{\prime}k}^{e} \xi^{N-k}, \quad j^{\prime} \! = \! 1,
\dotsc,N, \\
\leftthreetimes^{\raise+0.3ex\hbox{$\scriptstyle e$}}(b_{0}^{e}) \! = \! \mi 
(-1)^{N} \eta_{b_{0}^{e}}, \quad 
\leftthreetimes^{\raise+0.3ex\hbox{$\scriptstyle e$}}(a_{N+1}^{e}) \! = \! 
\eta_{a_{N+1}^{e}}, \quad 
\leftthreetimes^{\raise+0.3ex\hbox{$\scriptstyle e$}}(b_{j}^{e}) \! = \! \mi 
(-1)^{N-j} \eta_{b_{j}^{e}}, \quad 
\leftthreetimes^{\raise+0.3ex\hbox{$\scriptstyle e$}}(a_{j}^{e}) \! = \! 
(-1)^{N-j+1} \eta_{a_{j}^{e}}, \\
\eta_{b_{0}^{e}} \! := \! \left((a_{N+1}^{e} \! - \! b_{0}^{e}) \prod_{k=1}^{N}
(b_{k}^{e} \! - \! b_{0}^{e})(a_{k}^{e} \! - \! b_{0}^{e}) \right)^{1/2}, \\
\eta_{a_{N+1}^{e}} \! := \! \left((a_{N+1}^{e} \! - \! b_{0}^{e}) \prod_{k=
1}^{N}(a_{N+1}^{e} \! - \! b_{k}^{e})(a_{N+1}^{e} \! - \! a_{k}^{e}) \right)^{
1/2}, \\
\eta_{b_{j}^{e}} \! := \! \left(\! (b_{j}^{e} \! - \! a_{j}^{e})(a_{N+1}^{e} 
\! - \! b_{j}^{e})(b_{j}^{e} \! - \! b_{0}^{e}) \prod_{k=1}^{j-1}(b_{j}^{e} \! 
- \! b_{k}^{e})(b_{j}^{e} \! - \! a_{k}^{e}) \prod_{l=j+1}^{N}(b_{l}^{e} \! - 
\! b_{j}^{e})(a_{l}^{e} \! - \! b_{j}^{e}) \! \right)^{1/2}, \\
\eta_{a_{j}^{e}} \! := \! \left(\! (b_{j}^{e} \! - \! a_{j}^{e})(a_{N+1}^{e} 
\! - \! a_{j}^{e})(a_{j}^{e} \! - \! b_{0}^{e}) \prod_{k=1}^{j-1}(a_{j}^{e} \! 
- \! b_{k}^{e})(a_{j}^{e} \! - \! a_{k}^{e}) \prod_{l=j+1}^{N}(b_{l}^{e} \! - 
\! a_{j}^{e})(a_{l}^{e} \! - \! a_{j}^{e}) \! \right)^{1/2},
\end{gather}
where $c^{e}_{j^{\prime}k^{\prime}}$, $j^{\prime},k^{\prime} \! = \! 1,\dotsc,
N$, are obtained {}from Equations~{\rm (E1)} and~{\rm (E2)}, $\eta_{b_{j-1}^{
e}},\eta_{a_{j}^{e}} \! > \! 0$, $j \! = \! 1,\dotsc,N \! + \! 1$, and
\begin{align}
\widehat{\alpha}^{e}_{0}(b_{0}^{e}) \! &= \dfrac{4}{3} \mi (-1)^{N}h_{V}^{e}
(b_{0}^{e}) \eta_{b_{0}^{e}}, \\
\widehat{\alpha}^{e}_{1}(b_{0}^{e}) \! &= \mi (-1)^{N} \! \left(\dfrac{2}{5}
h_{V}^{e}(b_{0}^{e}) \eta_{b_{0}^{e}} \! \left(\sum_{l=1}^{N} \! \left(\dfrac{
1}{b_{0}^{e} \! - \! b_{l}^{e}} \! + \! \dfrac{1}{b_{0}^{e} \! - \! a_{l}^{e}} 
\right) \! + \! \dfrac{1}{b_{0}^{e} \! - \! a_{N+1}^{e}} \right) \! + \! 
\dfrac{4}{5}(h_{V}^{e}(b_{0}^{e}))^{\prime} \eta_{b_{0}^{e}} \right), \\
\widehat{\alpha}^{e}_{0}(a_{N+1}^{e}) \! &= \dfrac{4}{3}h_{V}^{e}(a_{N+1}^{e})
\eta_{a_{N+1}^{e}}, \\
\widehat{\alpha}^{e}_{1}(a_{N+1}^{e}) \! &= \dfrac{2}{5}h_{V}^{e}(a_{N+1}^{e}) 
\eta_{a_{N+1}^{e}} \! \left(\sum_{l=1}^{N} \! \left(\dfrac{1}{a_{N+1}^{e} \! - 
\! b_{l}^{e}} \! + \! \dfrac{1}{a_{N+1}^{e} \! - \! a_{l}^{e}} \right) \! + \! 
\dfrac{1}{a_{N+1}^{e} \! - \! b_{0}^{e}} \right) \nonumber \\
&+\dfrac{4}{5}(h_{V}^{e}(a_{N+1}^{e}))^{\prime} \eta_{a_{N+1}^{e}}, \\
\widehat{\alpha}^{e}_{0}(b_{j}^{e}) \! &= \dfrac{4}{3} \mi (-1)^{N-j}h_{V}^{e}
(b_{j}^{e}) \eta_{b_{j}^{e}}, \\
\widehat{\alpha}^{e}_{1}(b_{j}^{e}) \! &= \mi (-1)^{N-j} \! \left(\dfrac{2}{5}
h_{V}^{e}(b_{j}^{e}) \eta_{b_{j}^{e}} \! \left(\sum_{\substack{k=1\\k \not= 
j}}^{N} \! \left(\dfrac{1}{b_{j}^{e} \! - \! b_{k}^{e}} \! + \! \dfrac{1}{
b_{j}^{ e} \! - \! a_{k}^{e}} \right) \! + \! \dfrac{1}{b_{j}^{e} \! - \! 
a_{j}^{e}} \! + \! \dfrac{1}{b_{j}^{e} \! - \! b_{0}^{e}} \! + \! \dfrac{
1}{b_{j}^{e} \! - \! a_{N+1}^{e}} \right) \right. \nonumber \\
&\left. +\dfrac{4}{5}(h_{V}^{e}(b_{j}^{e}))^{\prime} \eta_{b_{j}^{e}} \right),
\\
\widehat{\alpha}^{e}_{0}(a_{j}^{e}) \! &= \dfrac{4}{3}(-1)^{N-j+1}h_{V}^{e}
(a_{j}^{e}) \eta_{a_{j}^{e}}, \\
\widehat{\alpha}^{e}_{1}(a_{j}^{e}) \! &= (-1)^{N-j+1} \! \left(\dfrac{2}{5}
h_{V}^{e}(a_{j}^{e}) \eta_{a_{j}^{e}} \! \left(\sum_{\substack{k=1\\k \not= 
j}}^{N} \! \left(\dfrac{1}{a_{j}^{e} \! - \! b_{k}^{e}} \! + \! \dfrac{1}{
a_{j}^{ e} \! - \! a_{k}^{e}} \right) \! + \! \dfrac{1}{a_{j}^{e} \! - \! 
b_{j}^{e}} \! + \! \dfrac{1}{a_{j}^{e} \! - \! a_{N+1}^{e}} \! + \! \dfrac{
1}{a_{j}^{e} \! - \! b_{0}^{e}} \right) \right. \nonumber \\
&\left. +\dfrac{4}{5}(h_{V}^{e}(a_{j}^{e}))^{\prime} \eta_{a_{j}^{e}}
\right);
\end{align}
{\rm \pmb{(2)}} for $z \! \in \! \Upsilon^{e}_{2}$ $(\subset \! \mathbb{C}_{
-})$,
\begin{align}
\boldsymbol{\pi}_{2n}(z) \underset{n \to \infty}{=}& \, \exp \! \left(n(g^{e}
(z) \! + \! Q_{e}) \right) \! \left(
(\overset{e}{m}^{\raise-1.0ex\hbox{$\scriptstyle \infty$}}(z))_{11} \! \left(
1 \! + \! \dfrac{1}{n} \! \left(\mathscr{R}^{e}_{\infty}(z) \right)_{11} \! +
\! \mathcal{O} \! \left(\dfrac{1}{n^{2}} \right) \right) \right. \nonumber \\
+&\left. \, (\overset{e}{m}^{\raise-1.0ex\hbox{$\scriptstyle \infty$}}(z))_{2
1} \! \left(\dfrac{1}{n} \! \left(\mathscr{R}^{e}_{\infty}(z) \right)_{12} \!
+ \! \mathcal{O} \! \left(\dfrac{1}{n^{2}} \right) \right) \right),
\end{align}
and
\begin{align}
\int_{\mathbb{R}} \dfrac{\boldsymbol{\pi}_{2n}(s) \me^{-n \widetilde{V}(s)}}{s
\! - \! z} \, \dfrac{\md s}{2 \pi \mi} \underset{n \to \infty}{=}& \, \exp
\! \left(-n(g^{e}(z) \! - \! \ell_{e} \! + \! Q_{e}) \right) \! \left(\!
(\overset{e}{m}^{\raise-1.0ex\hbox{$\scriptstyle \infty$}}(z))_{12} \! \left(
\! 1 \! + \! \dfrac{1}{n} \! \left(\mathscr{R}^{e}_{\infty}(z) \right)_{11}
\right. \right. \nonumber \\
+&\left. \left. \mathcal{O} \! \left(\dfrac{1}{n^{2}} \right) \right) \! +
\! (\overset{e}{m}^{\raise-1.0ex\hbox{$\scriptstyle \infty$}}(z))_{22} \!
\left(\dfrac{1}{n} \! \left(\mathscr{R}^{e}_{\infty}(z) \right)_{12} \! + \!
\mathcal{O} \! \left(\dfrac{1}{n^{2}} \right) \right) \right);
\end{align}
{\rm \pmb{(3)}} for $z \! \in \! \Upsilon^{e}_{3}$ $(\subset \! \cup_{j=1}^{N
+1} \left\lbrace \mathstrut z \! \in \! \mathbb{C}^{\ast}; \, \Re (z) \! \in
\! (b_{j-1}^{e},a_{j}^{e}), \, \inf_{q \in (b_{j-1}^{e},a_{j}^{e})} \vert z \!
- \! q \vert \! < \! 2^{-1/2} \min \{\delta_{b_{j-1}}^{e},\delta_{a_{j}}^{e}\}
\right\rbrace \! \subset \! \mathbb{C}_{+})$,
\begin{align}
\boldsymbol{\pi}_{2n}(z) \underset{n \to \infty}{=}& \, \exp \! \left(n \!
\left(g^{e}(z) \! + \! Q_{e} \right) \right) \! \left(\! \left(\!
(\overset{e}{m}^{\raise-1.0ex\hbox{$\scriptstyle \infty$}}(z))_{11} \! + \!
(\overset{e}{m}^{\raise-1.0ex\hbox{$\scriptstyle \infty$}}(z))_{12} \me^{-4n
\pi \mi \int_{z}^{a_{N+1}^{e}} \psi_{V}^{e}(s) \, \md s} \right) \! \! \left(
\! 1 \! + \! \dfrac{1}{n} \! \left(\mathscr{R}^{e}_{\infty}(z) \right)_{11}
\right. \right. \nonumber \\
+&\left. \left. \mathcal{O} \! \left(\dfrac{1}{n^{2}} \right) \! \right) \! +
\! \left(\! (\overset{e}{m}^{\raise-1.0ex\hbox{$\scriptstyle \infty$}}(z))_{2
1} \! + \! (\overset{e}{m}^{\raise-1.0ex\hbox{$\scriptstyle \infty$}}(z))_{22}
\me^{-4n \pi \mi \int_{z}^{a_{N+1}^{e}} \psi_{V}^{e}(s) \, \md s} \right) \!
\! \left(\dfrac{1}{n} \! \left(\mathscr{R}^{e}_{\infty}(z) \right)_{12} \! +
\! \mathcal{O} \! \left(\dfrac{1}{n^{2}} \right) \! \right) \! \right),
\end{align}
and
\begin{align}
\int_{\mathbb{R}} \dfrac{\boldsymbol{\pi}_{2n}(s) \me^{-n \widetilde{V}(s)}}{s
\! - \! z} \, \dfrac{\md s}{2 \pi \mi} \underset{n \to \infty}{=}& \, \exp
\! \left(-n \! \left(g^{e}(z) \! - \! \ell_{e} \! + \! Q_{e} \right) \right)
\! \left(\! (\overset{e}{m}^{\raise-1.0ex\hbox{$\scriptstyle \infty$}}(z))_{1
2} \! \left(\! 1 \! + \! \dfrac{1}{n} \! \left(\mathscr{R}^{e}_{\infty}(z)
\right)_{11} \! + \! \mathcal{O} \! \left(\dfrac{1}{n^{2}} \right) \right)
\right. \nonumber \\
+&\left. \, (\overset{e}{m}^{\raise-1.0ex\hbox{$\scriptstyle \infty$}}(z))_{2
2} \! \left(\dfrac{1}{n} \! \left(\mathscr{R}^{e}_{\infty}(z) \right)_{12} \!
+ \! \mathcal{O} \! \left(\dfrac{1}{n^{2}} \right) \right) \right);
\end{align}
{\rm \pmb{(4)}} for $z \! \in \! \Upsilon^{e}_{4}$ $(\subset \! \cup_{j=1}^{N
+1} \left\lbrace \mathstrut z \! \in \! \mathbb{C}^{\ast}; \, \Re (z) \! \in
\! (b_{j-1}^{e},a_{j}^{e}), \, \inf_{q \in (b_{j-1}^{e},a_{j}^{e})} \vert z \!
- \! q \vert \! < \! 2^{-1/2} \min \{\delta_{b_{j-1}}^{e},\delta_{a_{j}}^{e}\}
\right\rbrace \! \subset \! \mathbb{C}_{-})$,
\begin{align}
\boldsymbol{\pi}_{2n}(z) \underset{n \to \infty}{=}& \, \exp \! \left(n \!
\left(g^{e}(z) \! + \! Q_{e} \right) \right) \! \left(\! \left(\!
(\overset{e}{m}^{\raise-1.0ex\hbox{$\scriptstyle \infty$}}(z))_{11} \! - \!
(\overset{e}{m}^{\raise-1.0ex\hbox{$\scriptstyle \infty$}}(z))_{12} \me^{4n
\pi \mi \int_{z}^{a_{N+1}^{e}} \psi_{V}^{e}(s) \, \md s} \right) \! \left(1
\! + \! \dfrac{1}{n} \! \left(\mathscr{R}^{e}_{\infty}(z) \right)_{11} \right.
\right. \nonumber \\
+&\left. \left. \mathcal{O} \! \left(\dfrac{1}{n^{2}} \right) \right) \! + \!
\left((\overset{e}{m}^{\raise-1.0ex\hbox{$\scriptstyle \infty$}}(z))_{21} \!
- \! (\overset{e}{m}^{\raise-1.0ex\hbox{$\scriptstyle \infty$}}(z))_{22}
\me^{4n \pi \mi \int_{z}^{a_{N+1}^{e}} \psi_{V}^{e}(s) \, \md s} \right) \!
\left(\dfrac{1}{n} \! \left(\mathscr{R}^{e}_{\infty}(z) \right)_{12} \! + \!
\mathcal{O} \! \left(\dfrac{1}{n^{2}} \right) \right) \right),
\end{align}
and
\begin{align}
\int_{\mathbb{R}} \dfrac{\boldsymbol{\pi}_{2n}(s) \me^{-n \widetilde{V}(s)}}{s
\! - \! z} \, \dfrac{\md s}{2 \pi \mi} \underset{n \to \infty}{=}& \, \exp
\! \left(-n \! \left(g^{e}(z) \! - \! \ell_{e} \! + \! Q_{e} \right) \right)
\! \left(\! (\overset{e}{m}^{\raise-1.0ex\hbox{$\scriptstyle \infty$}}(z))_{1
2} \! \left(1 \! + \! \dfrac{1}{n} \! \left(\mathscr{R}^{e}_{\infty}(z)
\right)_{11} \! + \! \mathcal{O} \! \left(\dfrac{1}{n^{2}} \right) \right)
\right. \nonumber \\
+&\left. \, (\overset{e}{m}^{\raise-1.0ex\hbox{$\scriptstyle \infty$}}(z))_{2
2} \! \left(\dfrac{1}{n} \! \left(\mathscr{R}^{e}_{\infty}(z) \right)_{12} \!
+ \! \mathcal{O} \! \left(\dfrac{1}{n^{2}} \right) \right) \right);
\end{align}
{\rm \pmb{(5)}} for $z \! \in \! \Omega^{e,1}_{b_{j-1}}$ $(\subset \! \mathbb{
C}_{+} \cap \mathbb{U}^{e}_{\delta_{b_{j-1}}})$, $j \! = \! 1,\dotsc,N \! + \!
1$,
\begin{align}
\boldsymbol{\pi}_{2n}(z) \underset{n \to \infty}{=}& \, \exp \! \left(n \!
\left(g^{e}(z) \! + \! Q_{e} \right) \right) \! \left(\! (m^{b,1}_{p}(z))_{11}
\! \left(1 \! + \! \dfrac{1}{n} \! \left(\mathscr{R}^{e}_{\infty}(z) \! - \!
\widetilde{\mathscr{R}}^{e}_{\infty}(z) \right)_{11} \! + \! \mathcal{O} \!
\left(\dfrac{1}{n^{2}} \right) \right) \right. \nonumber \\
+&\left. \, (m^{b,1}_{p}(z))_{21} \! \left(\dfrac{1}{n} \! \left(\mathscr{R}^{
e}_{\infty}(z) \! - \! \widetilde{\mathscr{R}}^{e}_{\infty}(z) \right)_{12} \!
+ \! \mathcal{O} \! \left(\dfrac{1}{n^{2}} \right) \right) \right),
\end{align}
and
\begin{align}
\int_{\mathbb{R}} \dfrac{\boldsymbol{\pi}_{2n}(s) \me^{-n \widetilde{V}(s)}}{s
\! - \! z} \, \dfrac{\md s}{2 \pi \mi} \underset{n \to \infty}{=}& \, \exp \!
\left(-n \! \left(g^{e}(z) \! - \! \ell_{e} \! + \! Q_{e} \right) \right) \!
\left(\! (m^{b,1}_{p}(z))_{12} \! \left(1 \! + \! \dfrac{1}{n} \! \left(
\mathscr{R}^{e}_{\infty}(z) \! - \! \widetilde{\mathscr{R}}^{e}_{\infty}(z)
\right)_{11} \right. \right. \nonumber \\
+&\left. \left. \mathcal{O} \! \left(\dfrac{1}{n^{2}} \right) \right) \! + \!
(m^{b,1}_{p}(z))_{22} \! \left(\dfrac{1}{n} \! \left(\mathscr{R}^{e}_{\infty}
(z) \! - \! \widetilde{\mathscr{R}}^{e}_{\infty}(z) \right)_{12} \! + \!
\mathcal{O} \! \left(\dfrac{1}{n^{2}} \right) \right) \right),
\end{align}
where
\begin{align}
(m^{b,1}_{p}(z))_{11} :=& \, -\mi \sqrt{\smash[b]{\pi}} \, \me^{\frac{1}{2}n
\xi^{e}_{b_{j-1}}(z)} \! \left(\mi \! \left(\operatorname{Ai}(p_{b})(p_{b})^{
1/4} \! - \! \operatorname{Ai}^{\prime}(p_{b})(p_{b})^{-1/4} \right) \!
(\overset{e}{m}^{\raise-1.0ex\hbox{$\scriptstyle \infty$}}(z))_{11} \right.
\nonumber \\
-&\left. \left(\operatorname{Ai}(p_{b})(p_{b})^{1/4} \! + \!
\operatorname{Ai}^{\prime}(p_{b})(p_{b})^{-1/4} \right) \!
(\overset{e}{m}^{\raise-1.0ex\hbox{$\scriptstyle \infty$}}(z))_{12} \me^{-\mi
n \mho^{e}_{j-1}} \right), \\
(m^{b,1}_{p}(z))_{12} :=& \, \sqrt{\smash[b]{\pi}} \, \me^{-\frac{\mi \pi}{6}}
\me^{-\frac{1}{2}n \xi^{e}_{b_{j-1}}(z)} \! \left(\mi \! \left(-
\operatorname{Ai}(\omega^{2}p_{b})(p_{b})^{1/4} \! + \!
\omega^{2} \operatorname{Ai}^{\prime}(\omega^{2}p_{b})(p_{b})^{-1/4} \right)
\! (\overset{e}{m}^{\raise-1.0ex\hbox{$\scriptstyle \infty$}}(z))_{11} \right.
\nonumber \\
\times &\left. \me^{\mi n \mho^{e}_{j-1}} \! + \! \left(\operatorname{Ai}
(\omega^{2}p_{b})(p_{b})^{1/4} \! + \! \omega^{2} \operatorname{Ai}^{\prime}
(\omega^{2}p_{b})(p_{b})^{-1/4} \right) \!
(\overset{e}{m}^{\raise-1.0ex\hbox{$\scriptstyle \infty$}}(z))_{12}
\right), \\
(m^{b,1}_{p}(z))_{21} :=& \, -\mi \sqrt{\smash[b]{\pi}} \, \me^{\frac{1}{2}n
\xi^{e}_{b_{j-1}}(z)} \! \left(\mi \! \left(\operatorname{Ai}(p_{b})(p_{b})^{
1/4} \! - \! \operatorname{Ai}^{\prime}(p_{b})(p_{b})^{-1/4} \right) \!
(\overset{e}{m}^{\raise-1.0ex\hbox{$\scriptstyle \infty$}}(z))_{21} \right.
\nonumber \\
-&\left. \left(\operatorname{Ai}(p_{b})(p_{b})^{1/4} \! + \!
\operatorname{Ai}^{\prime}(p_{b})(p_{b})^{-1/4} \right) \!
(\overset{e}{m}^{\raise-1.0ex\hbox{$\scriptstyle \infty$}}(z))_{22} \me^{-
\mi n \mho^{e}_{j-1}} \right), \\
(m^{b,1}_{p}(z))_{22} :=& \, \sqrt{\smash[b]{\pi}} \, \me^{-\frac{\mi \pi}{6}}
\me^{-\frac{1}{2}n \xi^{e}_{b_{j-1}}(z)} \! \left(\mi \! \left(-
\operatorname{Ai}(\omega^{2}p_{b})(p_{b})^{1/4} \! + \! \omega^{2}
\operatorname{Ai}^{\prime}(\omega^{2}p_{b})(p_{b})^{-1/4} \right) \!
(\overset{e}{m}^{\raise-1.0ex\hbox{$\scriptstyle \infty$}}(z))_{21} \right.
\nonumber \\
\times &\left. \me^{\mi n \mho^{e}_{j-1}} \! + \! \left(\operatorname{Ai}
(\omega^{2}p_{b})(p_{b})^{1/4} \! + \! \omega^{2} \operatorname{Ai}^{\prime}
(\omega^{2}p_{b})(p_{b})^{-1/4} \right) \!
(\overset{e}{m}^{\raise-1.0ex\hbox{$\scriptstyle \infty$}}(z))_{22} \right),
\end{align}
with $\omega \! = \! \exp (2 \pi \mi/3)$,
\begin{gather}
\widetilde{\mathscr{R}}^{e}_{\infty}(z) \! := \! \sum_{j=1}^{N+1} \! \left(
\mathscr{R}^{\infty}_{b_{j-1}^{e}}(z) \pmb{1}_{\mathbb{U}^{e}_{\delta_{b_{j-
1}}}}(z) \! + \! \mathscr{R}^{\infty}_{a_{j}^{e}}(z) \pmb{1}_{\mathbb{U}^{e}_{
\delta_{a_{j}}}}(z) \right), \\
\xi^{e}_{b_{j-1}}(z) \! = \! -2 \int_{z}^{b^{e}_{j-1}}(R_{e}(s))^{1/2}h_{V}^{e}
(s) \, \md s, \qquad \quad p_{b} \! = \! \left(\dfrac{3}{4}n \xi^{e}_{b_{j-1}}
(z) \right)^{2/3}, \\
\mathscr{R}^{\infty}_{b_{j-1}^{e}}(z) \! = \! \dfrac{1}{\xi_{b_{j-1}}^{e}(z)}
\overset{e}{m}^{\raise-1.0ex\hbox{$\scriptstyle \infty$}}(z) \!
\begin{pmatrix}
-(s_{1}+t_{1}) & -\mi (s_{1}-t_{1}) \me^{\mi n \mho_{j-1}^{e}} \\
-\mi (s_{1}-t_{1}) \me^{-\mi n \mho_{j-1}^{e}} & (s_{1}+t_{1})
\end{pmatrix} \!
(\overset{e}{m}^{\raise-1.0ex\hbox{$\scriptstyle \infty$}}(z))^{-1}, \\
\mathscr{R}_{a_{j}^{e}}^{\infty}(z) \! = \! \dfrac{1}{\xi_{a_{j}}^{e}(z)}
\overset{e}{m}^{\raise-1.0ex\hbox{$\scriptstyle \infty$}}(z) \!
\begin{pmatrix}
-(s_{1}+t_{1}) & \mi (s_{1}-t_{1}) \me^{\mi n \mho_{j}^{e}} \\
\mi (s_{1}-t_{1}) \me^{-\mi n \mho_{j}^{e}} & (s_{1}+t_{1})
\end{pmatrix} \!
(\overset{e}{m}^{\raise-1.0ex\hbox{$\scriptstyle \infty$}}(z))^{-1}, \\
\xi_{a_{j}}^{e}(z) \! = \! 2 \int_{a_{j}^{e}}^{z}(R_{e}(s))^{1/2}h_{V}^{e}(s)
\, \md s,
\end{gather}
and $\pmb{1}_{\mathbb{U}^{e}_{\delta_{b_{j-1}}}}(z)$ (resp., $\pmb{1}_{\mathbb{
U}^{e}_{\delta_{a_{j}}}}(z))$ the indicator (characteristic) function of the
set $\mathbb{U}^{e}_{\delta_{b_{j-1}}}$ (resp., $\mathbb{U}^{e}_{\delta_{a_{
j}}});$\\
{\rm \pmb{(6)}} for $z \! \in \! \Omega^{e,2}_{b_{j-1}}$ $(\subset \! \mathbb{
C}_{+} \cap \mathbb{U}^{e}_{\delta_{b_{j-1}}})$, $j \! = \! 1,\dotsc,N \! + \!
1$,
\begin{align}
\boldsymbol{\pi}_{2n}(z) \underset{n \to \infty}{=}& \, \exp \! \left(n \!
\left(g^{e}(z) \! + \! Q_{e} \right) \right) \! \left(\! \left(\! (m^{b,2}_{p}
(z))_{11} \! + \! (m^{b,2}_{p}(z))_{12} \me^{-4n \pi \mi \int_{z}^{a_{N+1}^{e}
} \psi_{V}^{e}(s) \, \md s} \right) \! \left(1 \! + \! \dfrac{1}{n} \! \left(
\mathscr{R}^{e}_{\infty}(z) \right. \right. \right. \nonumber \\
-&\left. \left. \left. \widetilde{\mathscr{R}}^{e}_{\infty}(z) \right)_{11} \!
+ \! \mathcal{O} \! \left(\dfrac{1}{n^{2}} \right) \right) \! + \! \left((m^{
b,2}_{p}(z))_{21} \! + \! (m^{b,2}_{p}(z))_{22} \me^{-4n \pi \mi \int_{z}^{a_{
N+1}^{e}} \psi_{V}^{e}(s) \, \md s} \right) \! \left(\dfrac{1}{n} \! \left(
\mathscr{R}^{e}_{\infty}(z) \right. \right. \right. \nonumber \\
-&\left. \left. \left. \widetilde{\mathscr{R}}^{e}_{\infty}(z) \right)_{12} \!
+ \! \mathcal{O} \! \left(\dfrac{1}{n^{2}} \right) \right) \right),
\end{align}
and
\begin{align}
\int_{\mathbb{R}} \dfrac{\boldsymbol{\pi}_{2n}(s) \me^{-n \widetilde{V}(s)}}{s
\! - \! z} \, \dfrac{\md s}{2 \pi \mi} \underset{n \to \infty}{=}& \, \exp \!
\left(-n \! \left(g^{e}(z) \! - \! \ell_{e} \! + \! Q_{e} \right) \right) \!
\left(\! (m^{b,2}_{p}(z))_{12} \! \left(1 \! + \! \dfrac{1}{n} \! \left(
\mathscr{R}^{e}_{\infty}(z) \! - \! \widetilde{\mathscr{R}}^{e}_{\infty}(z)
\right)_{11} \right. \right. \nonumber \\
+&\left. \left. \mathcal{O} \! \left(\dfrac{1}{n^{2}} \right) \right) \! + \!
(m^{b,2}_{p}(z))_{22} \! \left(\dfrac{1}{n} \! \left(\mathscr{R}^{e}_{\infty}
(z) \! - \! \widetilde{\mathscr{R}}^{e}_{\infty}(z) \right)_{12} \! + \!
\mathcal{O} \! \left(\dfrac{1}{n^{2}} \right) \right) \right),
\end{align}
where
\begin{align}
(m^{b,2}_{p}(z))_{11} :=& \, -\mi \sqrt{\smash[b]{\pi}} \, \me^{\frac{1}{2}n
\xi^{e}_{b_{j-1}}(z)} \! \left(\mi \! \left(-\omega \operatorname{Ai}(\omega
p_{b})(p_{b})^{1/4} \! + \! \omega^{2} \operatorname{Ai}^{\prime}(\omega
p_{b})(p_{b})^{-1/4} \right) \!
(\overset{e}{m}^{\raise-1.0ex\hbox{$\scriptstyle \infty$}}(z))_{11} \right.
\nonumber \\
+&\left. \left(\omega \operatorname{Ai}(\omega p_{b})(p_{b})^{1/4} \! + \!
\omega^{2} \operatorname{Ai}^{\prime}(\omega p_{b})(p_{b})^{-1/4} \right) \!
(\overset{e}{m}^{\raise-1.0ex\hbox{$\scriptstyle \infty$}}(z))_{12} \me^{-\mi
n \mho^{e}_{j-1}} \right), \\
(m^{b,2}_{p}(z))_{12} :=& \, \sqrt{\smash[b]{\pi}} \, \me^{-\frac{\mi \pi}{6}}
\me^{-\frac{1}{2}n \xi^{e}_{b_{j-1}}(z)} \! \left(\mi \! \left(-\operatorname{
Ai}(\omega^{2}p_{b})(p_{b})^{1/4} \! + \! \omega^{2} \operatorname{Ai}^{\prime}
(\omega p_{b})(p_{b})^{-1/4} \right) \!
(\overset{e}{m}^{\raise-1.0ex\hbox{$\scriptstyle \infty$}}(z))_{11} \right.
\nonumber \\
\times &\left. \me^{\mi n \mho^{e}_{j-1}} \! + \! \left(\operatorname{Ai}
(\omega^{2} p_{b})(p_{b})^{1/4} \! + \! \omega^{2} \operatorname{Ai}^{\prime}
(\omega p_{b})(p_{b})^{-1/4} \right) \!
(\overset{e}{m}^{\raise-1.0ex\hbox{$\scriptstyle \infty$}}(z))_{12} \right), \\
(m^{b,2}_{p}(z))_{21} :=& \, -\mi \sqrt{\smash[b]{\pi}} \, \me^{\frac{1}{2}n
\xi^{e}_{b_{j-1}}(z)} \! \left(\mi \! \left(-\omega \operatorname{Ai}(\omega
p_{b})(p_{b})^{1/4} \! + \! \omega^{2} \operatorname{Ai}^{\prime}(\omega
p_{b})(p_{b})^{-1/4} \right) \!
(\overset{e}{m}^{\raise-1.0ex\hbox{$\scriptstyle \infty$}}(z))_{21} \right.
\nonumber \\
+&\left. \left(\omega \operatorname{Ai}(\omega p_{b})(p_{b})^{1/4} \! + \!
\omega^{2} \operatorname{Ai}^{\prime}(\omega p_{b})(p_{b})^{-1/4} \right) \!
(\overset{e}{m}^{\raise-1.0ex\hbox{$\scriptstyle \infty$}}(z))_{22} \me^{-\mi
n \mho^{e}_{j-1}} \right), \\
(m^{b,2}_{p}(z))_{22} :=& \, \sqrt{\smash[b]{\pi}} \, \me^{-\frac{\mi \pi}{6}}
\me^{-\frac{1}{2}n \xi^{e}_{b_{j-1}}(z)} \! \left(\mi \! \left(-\operatorname{
Ai}(\omega^{2}p_{b})(p_{b})^{1/4} \! + \! \omega^{2} \operatorname{Ai}^{\prime}
(\omega p_{b})(p_{b})^{-1/4} \right) \!
(\overset{e}{m}^{\raise-1.0ex\hbox{$\scriptstyle \infty$}}(z))_{21} \right.
\nonumber \\
\times &\left. \me^{\mi n \mho^{e}_{j-1}} \! + \! \left(\operatorname{Ai}
(\omega^{2}p_{b})(p_{b})^{1/4} \! + \! \omega^{2} \operatorname{Ai}^{\prime}
(\omega p_{b})(p_{b})^{-1/4} \right) \!
(\overset{e}{m}^{\raise-1.0ex\hbox{$\scriptstyle \infty$}}(z))_{22}
\right);
\end{align}
{\rm \pmb{(7)}} for $z \! \in \! \Omega^{e,3}_{b_{j-1}}$ $(\subset \! \mathbb{
C}_{-} \cap \mathbb{U}^{e}_{\delta_{b_{j-1}}})$, $j \! = \! 1,\dotsc,N \! + \!
1$,
\begin{align}
\boldsymbol{\pi}_{2n}(z) \underset{n \to \infty}{=}& \, \exp \! \left(n \!
\left(g^{e}(z) \! + \! Q_{e} \right) \right) \! \left(\! \left(\! (m^{b,3}_{p}
(z))_{11} \! - \! (m^{b,3}_{p}(z))_{12} \me^{4n \pi \mi \int_{z}^{a_{N+1}^{e}}
\psi_{V}^{e}(s) \, \md s} \right) \! \left(1 \! + \! \dfrac{1}{n} \! \left(
\mathscr{R}^{e}_{\infty}(z) \right. \right. \right. \nonumber \\
-&\left. \left. \left. \widetilde{\mathscr{R}}^{e}_{\infty}(z) \right)_{11} \!
+ \! \mathcal{O} \! \left(\dfrac{1}{n^{2}} \right) \right) \! + \! \left((m^{
b,3}_{p}(z))_{21} \! - \! (m^{b,3}_{p}(z))_{22} \me^{4n \pi \mi \int_{z}^{a_{
N+1}^{e}} \psi_{V}^{e}(s) \, \md s} \right) \! \left(\dfrac{1}{n} \! \left(
\mathscr{R}^{e}_{\infty}(z) \right. \right. \right. \nonumber \\
-&\left. \left. \left. \widetilde{\mathscr{R}}^{e}_{\infty}(z) \right)_{12}
\! + \! \mathcal{O} \! \left(\dfrac{1}{n^{2}} \right) \right) \right),
\end{align}
and
\begin{align}
\int_{\mathbb{R}} \dfrac{\boldsymbol{\pi}_{2n}(s) \me^{-n \widetilde{V}(s)}}{s
\! - \! z} \, \dfrac{\md s}{2 \pi \mi} \underset{n \to \infty}{=}& \, \exp \!
\left(-n \! \left(g^{e}(z) \! - \! \ell_{e} \! + \! Q_{e} \right) \right) \!
\left(\! (m^{b,3}_{p}(z))_{12} \! \left(1 \! + \! \dfrac{1}{n} \! \left(
\mathscr{R}^{e}_{\infty}(z) \! - \! \widetilde{\mathscr{R}}^{e}_{\infty}(z)
\right)_{11} \right. \right. \nonumber \\
+&\left. \left. \mathcal{O} \! \left(\dfrac{1}{n^{2}} \right) \right) \! + \!
(m^{b,3}_{p}(z))_{22} \! \left(\dfrac{1}{n} \! \left(\mathscr{R}^{e}_{\infty}
(z) \! - \! \widetilde{\mathscr{R}}^{e}_{\infty}(z) \right)_{12} \! + \!
\mathcal{O} \! \left(\dfrac{1}{n^{2}} \right) \right) \right),
\end{align}
where
\begin{align}
(m^{b,3}_{p}(z))_{11} :=& \, -\mi \sqrt{\smash[b]{\pi}} \, \me^{\frac{1}{2}n
\xi^{e}_{b_{j-1}}(z)} \! \left(\mi \! \left(-\omega^{2} \operatorname{Ai}
(\omega^{2} p_{b})(p_{b})^{1/4} \! + \! \omega \operatorname{Ai}^{\prime}
(\omega^{2} p_{b})(p_{b})^{-1/4} \right) \!
(\overset{e}{m}^{\raise-1.0ex\hbox{$\scriptstyle \infty$}}(z))_{11} \right.
\nonumber \\
+&\left. \left(\omega^{2} \operatorname{Ai}(\omega^{2}p_{b})(p_{b})^{1/4} \! +
\! \omega \operatorname{Ai}^{\prime}(\omega^{2} p_{b})(p_{b})^{-1/4} \right)
\! (\overset{e}{m}^{\raise-1.0ex\hbox{$\scriptstyle \infty$}}(z))_{12} \me^{
\mi n \mho^{e}_{j-1}} \right), \\
(m^{b,3}_{p}(z))_{12} :=& \, \sqrt{\smash[b]{\pi}} \, \me^{-\frac{\mi \pi}{6}}
\me^{-\frac{1}{2}n \xi^{e}_{b_{j-1}}(z)} \! \left(\mi \! \left(\omega^{2}
\operatorname{Ai}(\omega p_{b})(p_{b})^{1/4} \! - \! \operatorname{Ai}^{\prime}
(\omega p_{b})(p_{b})^{-1/4} \right) \!
(\overset{e}{m}^{\raise-1.0ex\hbox{$\scriptstyle \infty$}}(z))_{11} \right.
\nonumber \\
\times &\left. \me^{-\mi n \mho^{e}_{j-1}} \! - \! \left(\omega^{2}
\operatorname{Ai}(\omega p_{b})(p_{b})^{1/4} \! + \! \operatorname{Ai}^{\prime}
(\omega p_{b})(p_{b})^{-1/4} \right) \!
(\overset{e}{m}^{\raise-1.0ex\hbox{$\scriptstyle \infty$}}(z))_{12} \right), \\
(m^{b,3}_{p}(z))_{21} :=& \, -\mi \sqrt{\smash[b]{\pi}} \, \me^{\frac{1}{2}n
\xi^{e}_{b_{j-1}}(z)} \! \left(\mi \! \left(-\omega^{2} \operatorname{Ai}
(\omega^{2}p_{b})(p_{b})^{1/4} \! + \! \omega \operatorname{Ai}^{\prime}
(\omega^{2}p_{b})(p_{b})^{-1/4} \right) \!
(\overset{e}{m}^{\raise-1.0ex\hbox{$\scriptstyle \infty$}}(z))_{21} \right.
\nonumber \\
+&\left. \left(\omega^{2} \operatorname{Ai}(\omega^{2}p_{b})(p_{b})^{1/4} \! +
\! \omega \operatorname{Ai}^{\prime}(\omega^{2}p_{b})(p_{b})^{-1/4} \right) \!
(\overset{e}{m}^{\raise-1.0ex\hbox{$\scriptstyle \infty$}}(z))_{22} \me^{\mi n
\mho^{e}_{j-1}} \right), \\
(m^{b,3}_{p}(z))_{22} :=& \, \sqrt{\smash[b]{\pi}} \, \me^{-\frac{\mi \pi}{6}}
\me^{-\frac{1}{2}n \xi^{e}_{b_{j-1}}(z)} \! \left(\mi \! \left(\omega^{2}
\operatorname{Ai}(\omega p_{b})(p_{b})^{1/4} \! - \! \operatorname{Ai}^{\prime}
(\omega p_{b})(p_{b})^{-1/4} \right) \!
(\overset{e}{m}^{\raise-1.0ex\hbox{$\scriptstyle \infty$}}(z))_{21} \right.
\nonumber \\
\times &\left. \me^{-\mi n \mho^{e}_{j-1}} \! - \! \left(\omega^{2}
\operatorname{Ai}(\omega p_{b})(p_{b})^{1/4} \! + \! \operatorname{Ai}^{\prime}
(\omega p_{b})(p_{b})^{-1/4} \right) \!
(\overset{e}{m}^{\raise-1.0ex\hbox{$\scriptstyle \infty$}}(z))_{22}
\right);
\end{align}
{\rm \pmb{(8)}} for $z \! \in \! \Omega^{e,4}_{b_{j-1}}$ $(\subset \! \mathbb{
C}_{-} \cap \mathbb{U}^{e}_{\delta_{b_{j-1}}})$, $j \! = \! 1,\dotsc,N \! + \!
1$,
\begin{align}
\boldsymbol{\pi}_{2n}(z) \underset{n \to \infty}{=}& \, \exp \! \left(n \!
\left(g^{e}(z) \! + \! Q_{e} \right) \right) \! \left(\! (m^{b,4}_{p}(z))_{11}
\! \left(1 \! + \! \dfrac{1}{n} \! \left(\mathscr{R}^{e}_{\infty}(z) \! - \!
\widetilde{\mathscr{R}}^{e}_{\infty}(z) \right)_{11} \! + \! \mathcal{O} \!
\left(\dfrac{1}{n^{2}} \right) \right) \right. \nonumber \\
+&\left. \, (m^{b,4}_{p}(z))_{21} \! \left(\dfrac{1}{n} \! \left(\mathscr{R}^{
e}_{\infty}(z) \! - \! \widetilde{\mathscr{R}}^{e}_{\infty}(z) \right)_{12} \!
+ \! \mathcal{O} \! \left(\dfrac{1}{n^{2}} \right) \right) \right),
\end{align}
and
\begin{align}
\int_{\mathbb{R}} \dfrac{\boldsymbol{\pi}_{2n}(s) \me^{-n \widetilde{V}(s)}}{s
\! - \! z} \, \dfrac{\md s}{2 \pi \mi} \underset{n \to \infty}{=}& \, \exp \!
\left(-n \! \left(g^{e}(z) \! - \! \ell_{e} \! + \! Q_{e} \right) \right) \!
\left(\! (m^{b,4}_{p}(z))_{12} \! \left(1 \! + \! \dfrac{1}{n} \! \left(
\mathscr{R}^{e}_{\infty}(z) \! - \! \widetilde{\mathscr{R}}^{e}_{\infty}(z)
\right)_{11} \right. \right. \nonumber \\
+&\left. \left. \mathcal{O} \! \left(\dfrac{1}{n^{2}} \right) \right) \! + \!
(m^{b,4}_{p}(z))_{22} \! \left(\dfrac{1}{n} \! \left(\mathscr{R}^{e}_{\infty}
(z) \! - \! \widetilde{\mathscr{R}}^{e}_{\infty}(z) \right)_{12} \! + \!
\mathcal{O} \! \left(\dfrac{1}{n^{2}} \right) \right) \right),
\end{align}
where
\begin{align}
(m^{b,4}_{p}(z))_{11} :=& \, -\mi \sqrt{\smash[b]{\pi}} \, \me^{\frac{1}{2}n
\xi^{e}_{b_{j-1}}(z)} \! \left(\mi \! \left(\operatorname{Ai}(p_{b})(p_{b})^{
1/4} \! - \! \operatorname{Ai}^{\prime}(p_{b})(p_{b})^{-1/4} \right) \!
(\overset{e}{m}^{\raise-1.0ex\hbox{$\scriptstyle \infty$}}(z))_{11} \right.
\nonumber \\
-&\left. \left(\operatorname{Ai}(p_{b})(p_{b})^{1/4} \! + \!
\operatorname{Ai}^{\prime}(p_{b})(p_{b})^{-1/4} \right) \!
(\overset{e}{m}^{\raise-1.0ex\hbox{$\scriptstyle \infty$}}(z))_{12} \me^{\mi n
\mho^{e}_{j-1}} \right), \\
(m^{b,4}_{p}(z))_{12} :=& \, \sqrt{\smash[b]{\pi}} \, \me^{-\frac{\mi \pi}{6}}
\me^{-\frac{1}{2}n \xi^{e}_{b_{j-1}}(z)} \! \left(\mi \! \left(\omega^{2}
\operatorname{Ai}(\omega p_{b})(p_{b})^{1/4} \! - \! \operatorname{Ai}^{\prime}
(\omega p_{b})(p_{b})^{-1/4} \right) \!
(\overset{e}{m}^{\raise-1.0ex\hbox{$\scriptstyle \infty$}}(z))_{11} \right.
\nonumber \\
\times &\left. \me^{-\mi n \mho^{e}_{j-1}} \! - \! \left(\omega^{2}
\operatorname{Ai}(\omega p_{b})(p_{b})^{1/4} \! + \! \operatorname{Ai}^{\prime}
(\omega p_{b})(p_{b})^{-1/4} \right) \!
(\overset{e}{m}^{\raise-1.0ex\hbox{$\scriptstyle \infty$}}(z))_{12} \right), \\
(m^{b,4}_{p}(z))_{21} :=& \, -\mi \sqrt{\smash[b]{\pi}} \, \me^{\frac{1}{2}n
\xi^{e}_{b_{j-1}}(z)} \! \left(\mi \! \left(\operatorname{Ai}(p_{b})(p_{b})^{
1/4} \! - \! \operatorname{Ai}^{\prime}(p_{b})(p_{b})^{-1/4} \right) \!
(\overset{e}{m}^{\raise-1.0ex\hbox{$\scriptstyle \infty$}}(z))_{21} \right.
\nonumber \\
-&\left. \left(\operatorname{Ai}(p_{b})(p_{b})^{1/4} \! + \!
\operatorname{Ai}^{\prime}(p_{b})(p_{b})^{-1/4} \right) \!
(\overset{e}{m}^{\raise-1.0ex\hbox{$\scriptstyle \infty$}}(z))_{22} \me^{\mi n
\mho^{e}_{j-1}} \right), \\
(m^{b,4}_{p}(z))_{22} :=& \, \sqrt{\smash[b]{\pi}} \, \me^{-\frac{\mi \pi}{6}}
\me^{-\frac{1}{2}n \xi^{e}_{b_{j-1}}(z)} \! \left(\mi \! \left(\omega^{2}
\operatorname{Ai}(\omega p_{b})(p_{b})^{1/4} \! - \! \operatorname{Ai}^{\prime}
(\omega p_{b})(p_{b})^{-1/4} \right) \!
(\overset{e}{m}^{\raise-1.0ex\hbox{$\scriptstyle \infty$}}(z))_{21} \right.
\nonumber \\
\times &\left. \me^{-\mi n \mho^{e}_{j-1}} \! - \! \left(\omega^{2}
\operatorname{Ai}(\omega p_{b})(p_{b})^{1/4} \! + \! \operatorname{Ai}^{\prime}
(\omega p_{b})(p_{b})^{-1/4} \right) \!
(\overset{e}{m}^{\raise-1.0ex\hbox{$\scriptstyle \infty$}}(z))_{22}
\right);
\end{align}
{\rm \pmb{(9)}} for $z \! \in \! \Omega^{e,1}_{a_{j}}$ $(\subset \! \mathbb{
C}_{+} \cap \mathbb{U}^{e}_{\delta_{a_{j}}})$, $j \! = \! 1,\dotsc,N \! + \!
1$,
\begin{align}
\boldsymbol{\pi}_{2n}(z) \underset{n \to \infty}{=}& \, \exp \! \left(n \!
\left(g^{e}(z) \! + \! Q_{e} \right) \right) \! \left(\! (m^{a,1}_{p}(z))_{11}
\! \left(1 \! + \! \dfrac{1}{n} \! \left(\mathscr{R}^{e}_{\infty}(z) \! - \!
\widetilde{\mathscr{R}}^{e}_{\infty}(z) \right)_{11} \! + \! \mathcal{O} \!
\left(\dfrac{1}{n^{2}} \right) \right) \right. \nonumber \\
+&\left. \, (m^{a,1}_{p}(z))_{21} \! \left(\dfrac{1}{n} \! \left(\mathscr{R}^{
e}_{\infty}(z) \! - \! \widetilde{\mathscr{R}}^{e}_{\infty}(z) \right)_{12} \!
+ \! \mathcal{O} \! \left(\dfrac{1}{n^{2}} \right) \right) \right),
\end{align}
and
\begin{align}
\int_{\mathbb{R}} \dfrac{\boldsymbol{\pi}_{2n}(s) \me^{-n \widetilde{V}(s)}}{s
\! - \! z} \, \dfrac{\md s}{2 \pi \mi} \underset{n \to \infty}{=}& \, \exp \!
\left(-n \! \left(g^{e}(z) \! - \! \ell_{e} \! + \! Q_{e} \right) \right) \!
\left(\! (m^{a,1}_{p}(z))_{12} \! \left(1 \! + \! \dfrac{1}{n} \! \left(
\mathscr{R}^{e}_{\infty}(z) \! - \! \widetilde{\mathscr{R}}^{e}_{\infty}(z)
\right)_{11} \right. \right. \nonumber \\
+&\left. \left. \mathcal{O} \! \left(\dfrac{1}{n^{2}} \right) \right) \! + \!
(m^{a,1}_{p}(z))_{22} \! \left(\dfrac{1}{n} \! \left(\mathscr{R}^{e}_{\infty}
(z) \! - \! \widetilde{\mathscr{R}}^{e}_{\infty}(z) \right)_{12} \! + \!
\mathcal{O} \! \left(\dfrac{1}{n^{2}} \right) \right) \right),
\end{align}
where
\begin{align}
(m^{a,1}_{p}(z))_{11} :=& \, -\mi \sqrt{\smash[b]{\pi}} \, \me^{\frac{1}{2}
n \xi^{e}_{a_{j}}(z)} \! \left(\mi \! \left(\operatorname{Ai}(p_{a})(p_{a})^{
1/4} \! - \! \operatorname{Ai}^{\prime}(p_{a})(p_{a})^{-1/4} \right) \!
(\overset{e}{m}^{\raise-1.0ex\hbox{$\scriptstyle \infty$}}(z))_{11} \right.
\nonumber \\
+&\left. \left(\operatorname{Ai}(p_{a})(p_{a})^{1/4} \! + \!
\operatorname{Ai}^{\prime}(p_{a})(p_{a})^{-1/4} \right) \!
(\overset{e}{m}^{\raise-1.0ex\hbox{$\scriptstyle \infty$}}(z))_{12} \me^{-\mi
n \mho^{e}_{j}} \right), \\
(m^{a,1}_{p}(z))_{12} :=& \, \sqrt{\smash[b]{\pi}} \, \me^{-\frac{\mi \pi}{6}}
\me^{-\frac{1}{2}n \xi^{e}_{a_{j}}(z)} \! \left(\mi \! \left(\operatorname{Ai}
(\omega^{2}p_{a})(p_{a})^{1/4} \! - \! \omega^{2} \operatorname{Ai}^{\prime}
(\omega^{2}p_{a})(p_{a})^{-1/4} \right) \!
(\overset{e}{m}^{\raise-1.0ex\hbox{$\scriptstyle \infty$}}(z))_{11} \right.
\nonumber \\
\times &\left. \me^{\mi n \mho^{e}_{j}} \! + \! \left(\operatorname{Ai}
(\omega^{2}p_{a})(p_{a})^{1/4} \! + \! \omega^{2} \operatorname{Ai}^{\prime}
(\omega^{2}p_{a})(p_{a})^{-1/4} \right) \!
(\overset{e}{m}^{\raise-1.0ex\hbox{$\scriptstyle \infty$}}(z))_{12} \right), \\
(m^{a,1}_{p}(z))_{21} :=& \, -\mi \sqrt{\smash[b]{\pi}} \, \me^{\frac{1}{2}n
\xi^{e}_{a_{j}}(z)} \! \left(\mi \! \left(\operatorname{Ai}(p_{a})(p_{a})^{
1/4} \! - \! \operatorname{Ai}^{\prime}(p_{a})(p_{a})^{-1/4} \right) \!
(\overset{e}{m}^{\raise-1.0ex\hbox{$\scriptstyle \infty$}}(z))_{21} \right.
\nonumber \\
+&\left. \left(\operatorname{Ai}(p_{a})(p_{a})^{1/4} \! + \!
\operatorname{Ai}^{\prime}(p_{a})(p_{a})^{-1/4} \right) \!
(\overset{e}{m}^{\raise-1.0ex\hbox{$\scriptstyle \infty$}}(z))_{22} \me^{-\mi
n \mho^{e}_{j}} \right), \\
(m^{a,1}_{p}(z))_{22} :=& \, \sqrt{\smash[b]{\pi}} \, \me^{-\frac{\mi \pi}{6}}
\me^{-\frac{1}{2}n \xi^{e}_{a_{j}}(z)} \! \left(\mi \! \left(\operatorname{Ai}
(\omega^{2}p_{a})(p_{a})^{1/4} \! - \! \omega^{2} \operatorname{Ai}^{\prime}
(\omega^{2}p_{a})(p_{a})^{-1/4} \right) \!
(\overset{e}{m}^{\raise-1.0ex\hbox{$\scriptstyle \infty$}}(z))_{21} \right.
\nonumber \\
\times &\left. \me^{\mi n \mho^{e}_{j}} \! + \! \left(\operatorname{Ai}
(\omega^{2}p_{a})(p_{a})^{1/4} \! + \! \omega^{2} \operatorname{Ai}^{\prime}
(\omega^{2}p_{a})(p_{a})^{-1/4} \right) \!
(\overset{e}{m}^{\raise-1.0ex\hbox{$\scriptstyle \infty$}}(z))_{22}
\right),
\end{align}
with
\begin{equation}
p_{a} \! = \! \left(\dfrac{3}{4}n \xi^{e}_{a_{j}}(z) \right)^{2/3};
\end{equation}
{\rm \pmb{(10)}} for $z \! \in \! \Omega^{e,2}_{a_{j}}$ $(\subset \! \mathbb{
C}_{+} \cap \mathbb{U}^{e}_{\delta_{a_{j}}})$, $j \! = \! 1,\dotsc,N \! + \!
1$,
\begin{align}
\boldsymbol{\pi}_{2n}(z) \underset{n \to \infty}{=}& \, \exp \! \left(n \!
\left(g^{e}(z) \! + \! Q_{e} \right) \right) \! \left(\! \left(\! (m^{a,2}_{p}
(z))_{11} \! + \! (m^{a,2}_{p}(z))_{12} \me^{-4n \pi \mi \int_{z}^{a_{N+1}^{e}
} \psi_{V}^{e}(s) \, \md s} \right) \! \left(1 \! + \! \dfrac{1}{n} \! \left(
\mathscr{R}^{e}_{\infty}(z) \right. \right. \right. \nonumber \nonumber \\
-&\left. \left. \left. \widetilde{\mathscr{R}}^{e}_{\infty}(z) \right)_{11} \!
+ \! \mathcal{O} \! \left(\dfrac{1}{n^{2}} \right) \right) \! + \! \left((m^{
a,2}_{p}(z))_{21} \! + \! (m^{a,2}_{p}(z))_{22} \me^{-4n \pi \mi \int_{z}^{a_{
N+1}^{e}} \psi_{V}^{e}(s) \, \md s} \right) \! \left(\dfrac{1}{n} \! \left(
\mathscr{R}^{e}_{\infty}(z) \right. \right. \right. \nonumber \nonumber \\
-&\left. \left. \left. \widetilde{\mathscr{R}}^{e}_{\infty}(z) \right)_{12} \!
+ \! \mathcal{O} \! \left(\dfrac{1}{n^{2}} \right) \right) \right),
\end{align}
and
\begin{align}
\int_{\mathbb{R}} \dfrac{\boldsymbol{\pi}_{2n}(s) \me^{-n \widetilde{V}(s)}}{s
\! - \! z} \, \dfrac{\md s}{2 \pi \mi} \underset{n \to \infty}{=}& \, \exp
\! \left(-n \! \left(g^{e}(z) \! - \! \ell_{e} \! + \! Q_{e} \right) \right)
\! \left(\! (m^{a,2}_{p}(z))_{12} \! \left(1 \! + \! \dfrac{1}{n} \! \left(
\mathscr{R}^{e}_{\infty}(z) \! - \! \widetilde{\mathscr{R}}^{e}_{\infty}(z)
\right)_{11} \right. \right. \nonumber \nonumber \\
+&\left. \left. \mathcal{O} \! \left(\dfrac{1}{n^{2}} \right) \right) \! + \!
(m^{a,2}_{p}(z))_{22} \! \left(\dfrac{1}{n} \! \left(\mathscr{R}^{e}_{\infty}
(z) \! - \! \widetilde{\mathscr{R}}^{e}_{\infty}(z) \right)_{12} \! + \!
\mathcal{O} \! \left(\dfrac{1}{n^{2}} \right) \right) \right),
\end{align}
where
\begin{align}
(m^{a,2}_{p}(z))_{11} :=& \, -\mi \sqrt{\smash[b]{\pi}} \, \me^{\frac{1}{2}n
\xi^{e}_{a_{j}}(z)} \! \left(\mi \! \left(-\omega \operatorname{Ai}(\omega
p_{a})(p_{a})^{1/4} \! + \! \omega^{2} \operatorname{Ai}^{\prime}(\omega
p_{a})(p_{a})^{-1/4} \right) \!
(\overset{e}{m}^{\raise-1.0ex\hbox{$\scriptstyle \infty$}}(z))_{11} \right.
\nonumber \\
-&\left. \left(\omega \operatorname{Ai}(\omega p_{a})(p_{a})^{1/4} \! + \!
\omega^{2} \operatorname{Ai}^{\prime}(\omega p_{a})(p_{a})^{-1/4} \right) \!
(\overset{e}{m}^{\raise-1.0ex\hbox{$\scriptstyle \infty$}}(z))_{12} \me^{-\mi
n \mho^{e}_{j}} \right), \\
(m^{a,2}_{p}(z))_{12} :=& \, \sqrt{\smash[b]{\pi}} \, \me^{-\frac{\mi \pi}{6}}
\me^{-\frac{1}{2}n \xi^{e}_{a_{j}}(z)} \! \left(\mi \! \left(\operatorname{Ai}
(\omega^{2}p_{a})(p_{a})^{1/4} \! - \! \omega^{2} \operatorname{Ai}^{\prime}
(\omega^{2}p_{a})(p_{a})^{-1/4} \right) \!
(\overset{e}{m}^{\raise-1.0ex\hbox{$\scriptstyle \infty$}}(z))_{11} \right.
\nonumber \\
\times &\left. \me^{\mi n \mho^{e}_{j}} \! + \! \left(\operatorname{Ai}
(\omega^{2}p_{a})(p_{a})^{1/4} \! + \! \omega^{2} \operatorname{Ai}^{\prime}
(\omega^{2}p_{a})(p_{a})^{-1/4} \right) \!
(\overset{e}{m}^{\raise-1.0ex\hbox{$\scriptstyle \infty$}}(z))_{12} \right), \\
(m^{a,2}_{p}(z))_{21} :=& \, -\mi \sqrt{\smash[b]{\pi}} \, \me^{\frac{1}{2}n
\xi^{e}_{a_{j}}(z)} \! \left(\mi \! \left(-\omega \operatorname{Ai}(\omega
p_{a})(p_{a})^{1/4} \! + \! \omega^{2} \operatorname{Ai}^{\prime}(\omega
p_{a})(p_{a})^{-1/4} \right) \!
(\overset{e}{m}^{\raise-1.0ex\hbox{$\scriptstyle \infty$}}(z))_{21} \right.
\nonumber \\
-&\left. \left(\omega \operatorname{Ai}(\omega p_{a})(p_{a})^{1/4} \! + \!
\omega^{2} \operatorname{Ai}^{\prime}(\omega p_{a})(p_{a})^{-1/4} \right) \!
(\overset{e}{m}^{\raise-1.0ex\hbox{$\scriptstyle \infty$}}(z))_{22} \me^{-\mi
n \mho^{e}_{j}} \right), \\
(m^{a,2}_{p}(z))_{22} :=& \, \sqrt{\smash[b]{\pi}} \, \me^{-\frac{\mi \pi}{6}}
\me^{-\frac{1}{2}n \xi^{e}_{a_{j}}(z)} \! \left(\mi \! \left(\operatorname{Ai}
(\omega^{2}p_{a})(p_{a})^{1/4} \! - \! \omega^{2} \operatorname{Ai}^{\prime}
(\omega^{2}p_{a})(p_{a})^{-1/4} \right) \!
(\overset{e}{m}^{\raise-1.0ex\hbox{$\scriptstyle \infty$}}(z))_{21} \right.
\nonumber \\
\times &\left. \me^{\mi n \mho^{e}_{j}} \! + \! \left(\operatorname{Ai}
(\omega^{2}p_{a})(p_{a})^{1/4} \! + \! \omega^{2} \operatorname{Ai}^{\prime}
(\omega^{2}p_{a})(p_{a})^{-1/4} \right) \!
(\overset{e}{m}^{\raise-1.0ex\hbox{$\scriptstyle \infty$}}(z))_{22} \right);
\end{align}
{\rm \pmb{(11)}} for $z \! \in \! \Omega^{e,3}_{a_{j}}$ $(\subset \! \mathbb{
C}_{-} \cap \mathbb{U}^{e}_{\delta_{a_{j}}})$, $j \! = \! 1,\dotsc,N \! + \!
1$,
\begin{align}
\boldsymbol{\pi}_{2n}(z) \underset{n \to \infty}{=}& \, \exp \! \left(n \!
\left(g^{e}(z) \! + \! Q_{e} \right) \right) \! \left(\! \left(\! (m^{a,3}_{p}
(z))_{11} \! - \! (m^{a,3}_{p}(z))_{12} \me^{4n \pi \mi \int_{z}^{a_{N+1}^{e}}
\psi_{V}^{e}(s) \, \md s} \right) \! \left(1 \! + \! \dfrac{1}{n} \! \left(
\mathscr{R}^{e}_{\infty}(z) \right. \right. \right. \nonumber \\
-&\left. \left. \left. \widetilde{\mathscr{R}}^{e}_{\infty}(z) \right)_{11} \!
+ \! \mathcal{O} \! \left(\dfrac{1}{n^{2}} \right) \right) \! + \! \left((m^{
a,3}_{p}(z))_{21} \! - \! (m^{a,3}_{p}(z))_{22} \me^{4n \pi \mi \int_{z}^{a_{
N+1}^{e}} \psi_{V}^{e}(s) \, \md s} \right) \! \left(\dfrac{1}{n} \! \left(
\mathscr{R}^{e}_{\infty}(z) \right. \right. \right. \nonumber \\
-&\left. \left. \left. \widetilde{\mathscr{R}}^{e}_{\infty}(z) \right)_{12}
\! + \! \mathcal{O} \! \left(\dfrac{1}{n^{2}} \right) \right) \right),
\end{align}
and
\begin{align}
\int_{\mathbb{R}} \dfrac{\boldsymbol{\pi}_{2n}(s) \me^{-n \widetilde{V}(s)}}{s
\! - \! z} \, \dfrac{\md s}{2 \pi \mi} \underset{n \to \infty}{=}& \, \exp \!
\left(-n \! \left(g^{e}(z) \! - \! \ell_{e} \! + \! Q_{e} \right) \right) \!
\left(\! (m^{a,3}_{p}(z))_{12} \! \left(1 \! + \! \dfrac{1}{n} \! \left(
\mathscr{R}^{e}_{\infty}(z) \! - \! \widetilde{\mathscr{R}}^{e}_{\infty}(z)
\right)_{11} \right. \right. \nonumber \\
+&\left. \left. \mathcal{O} \! \left(\dfrac{1}{n^{2}} \right) \right) \! + \!
(m^{a,3}_{p}(z))_{22} \! \left(\dfrac{1}{n} \! \left(\mathscr{R}^{e}_{\infty}
(z) \! - \! \widetilde{\mathscr{R}}^{e}_{\infty}(z) \right)_{12} \! + \!
\mathcal{O} \! \left(\dfrac{1}{n^{2}} \right) \right) \right),
\end{align}
where
\begin{align}
(m^{a,3}_{p}(z))_{11} :=& \, -\mi \sqrt{\smash[b]{\pi}} \, \me^{\frac{1}{2}
n \xi^{e}_{a_{j}}(z)} \! \left(\mi \! \left(-\omega^{2} \operatorname{Ai}
(\omega^{2}p_{a})(p_{a})^{1/4} \! + \! \omega \operatorname{Ai}^{\prime}
(\omega^{2}p_{a})(p_{a})^{-1/4} \right) \!
(\overset{e}{m}^{\raise-1.0ex\hbox{$\scriptstyle \infty$}}(z))_{11} \right.
\nonumber \\
-&\left. \left(\omega^{2} \operatorname{Ai}(\omega^{2}p_{a})(p_{a})^{1/4} \! +
\! \omega \operatorname{Ai}^{\prime}(\omega^{2}p_{a})(p_{a})^{-1/4} \right) \!
(\overset{e}{m}^{\raise-1.0ex\hbox{$\scriptstyle \infty$}}(z))_{12} \me^{\mi
n \mho^{e}_{j}} \right), \\
(m^{a,3}_{p}(z))_{12} :=& \, \sqrt{\smash[b]{\pi}} \, \me^{-\frac{\mi \pi}{6}}
\me^{-\frac{1}{2}n \xi^{e}_{a_{j}}(z)} \! \left(\mi \! \left(-\omega^{2}
\operatorname{Ai}(\omega p_{a})(p_{a})^{1/4} \! + \! \operatorname{Ai}^{\prime}
(\omega p_{a})(p_{a})^{-1/4} \right) \!
(\overset{e}{m}^{\raise-1.0ex\hbox{$\scriptstyle \infty$}}(z))_{11} \right.
\nonumber \\
\times &\left. \me^{-\mi n \mho^{e}_{j}} \! - \! \left(\omega^{2}
\operatorname{Ai}(\omega p_{a})(p_{a})^{1/4} \! + \! \operatorname{Ai}^{\prime}
(\omega p_{a})(p_{a})^{-1/4} \right) \!
(\overset{e}{m}^{\raise-1.0ex\hbox{$\scriptstyle \infty$}}(z))_{12} \right), \\
(m^{a,3}_{p}(z))_{21} :=& \, -\mi \sqrt{\smash[b]{\pi}} \, \me^{\frac{1}{2}
n \xi^{e}_{a_{j}}(z)} \! \left(\mi \! \left(-\omega^{2} \operatorname{Ai}
(\omega^{2}p_{a})(p_{a})^{1/4} \! + \! \omega \operatorname{Ai}^{\prime}
(\omega^{2}p_{a})(p_{a})^{-1/4} \right) \!
(\overset{e}{m}^{\raise-1.0ex\hbox{$\scriptstyle \infty$}}(z))_{21} \right.
\nonumber \\
-&\left. \left(\omega^{2} \operatorname{Ai}(\omega^{2}p_{a})(p_{a})^{1/4} \! +
\! \omega \operatorname{Ai}^{\prime}(\omega^{2}p_{a})(p_{a})^{-1/4} \right) \!
(\overset{e}{m}^{\raise-1.0ex\hbox{$\scriptstyle \infty$}}(z))_{22} \me^{\mi n
\mho^{e}_{j}} \right), \\
(m^{a,3}_{p}(z))_{22} :=& \, \sqrt{\smash[b]{\pi}} \, \me^{-\frac{\mi \pi}{6}}
\me^{-\frac{1}{2}n \xi^{e}_{a_{j}}(z)} \! \left(\mi \! \left(-\omega^{2}
\operatorname{Ai}(\omega p_{a})(p_{a})^{1/4} \! + \! \operatorname{Ai}^{\prime}
(\omega p_{a})(p_{a})^{-1/4} \right) \!
(\overset{e}{m}^{\raise-1.0ex\hbox{$\scriptstyle \infty$}}(z))_{21} \right.
\nonumber \\
\times &\left. \me^{-\mi n \mho^{e}_{j}} \! - \! \left(\omega^{2}
\operatorname{Ai}(\omega p_{a})(p_{a})^{1/4} \! + \! \operatorname{Ai}^{\prime}
(\omega p_{a})(p_{a})^{-1/4} \right) \!
(\overset{e}{m}^{\raise-1.0ex\hbox{$\scriptstyle \infty$}}(z))_{22} \right);
\end{align}
and {\rm \pmb{(12)}} for $z \! \in \! \Omega^{e,4}_{a_{j}}$ $(\subset \!
\mathbb{C}_{-} \cap \mathbb{U}^{e}_{\delta_{a_{j}}})$, $j \! = \! 1,\dotsc,N
\! + \! 1$,
\begin{align}
\boldsymbol{\pi}_{2n}(z) \underset{n \to \infty}{=}& \, \exp \! \left(n \!
\left(g^{e}(z) \! + \! Q_{e} \right) \right) \! \left(\! (m^{a,4}_{p}(z))_{11}
\! \left(1 \! + \! \dfrac{1}{n} \! \left(\mathscr{R}^{e}_{\infty}(z) \! - \!
\widetilde{\mathscr{R}}^{e}_{\infty}(z) \right)_{11} \! + \! \mathcal{O} \!
\left(\dfrac{1}{n^{2}} \right) \right) \right. \nonumber \\
+&\left. \, (m^{a,4}_{p}(z))_{21} \! \left(\dfrac{1}{n} \! \left(\mathscr{R}^{
e}_{\infty}(z) \! - \! \widetilde{\mathscr{R}}^{e}_{\infty}(z) \right)_{12} \!
+ \! \mathcal{O} \! \left(\dfrac{1}{n^{2}} \right) \right) \right),
\end{align}
and
\begin{align}
\int_{\mathbb{R}} \dfrac{\boldsymbol{\pi}_{2n}(s) \me^{-n \widetilde{V}(s)}}{s
\! - \! z} \, \dfrac{\md s}{2 \pi \mi} \underset{n \to \infty}{=}& \, \exp \!
\left(-n \! \left(g^{e}(z) \! - \! \ell_{e} \! + \! Q_{e} \right) \right) \!
\left(\! (m^{a,4}_{p}(z))_{12} \! \left(1 \! + \! \dfrac{1}{n} \! \left(
\mathscr{R}^{e}_{\infty}(z) \! - \! \widetilde{\mathscr{R}}^{e}_{\infty}(z)
\right)_{11} \right. \right. \nonumber \\
+&\left. \left. \mathcal{O} \! \left(\dfrac{1}{n^{2}} \right) \right) \! + \!
(m^{a,4}_{p}(z))_{22} \! \left(\dfrac{1}{n} \! \left(\mathscr{R}^{e}_{\infty}
(z) \! - \! \widetilde{\mathscr{R}}^{e}_{\infty}(z) \right)_{12} \! + \!
\mathcal{O} \! \left(\dfrac{1}{n^{2}} \right) \right) \right),
\end{align}
where
\begin{align}
(m^{a,4}_{p}(z))_{11} :=& \, -\mi \sqrt{\smash[b]{\pi}} \, \me^{\frac{1}{2}n
\xi^{e}_{a_{j}}(z)} \! \left(\mi \! \left(\operatorname{Ai}(p_{a})(p_{a})^{
1/4} \! - \! \operatorname{Ai}^{\prime}(p_{a})(p_{a})^{-1/4} \right) \!
(\overset{e}{m}^{\raise-1.0ex\hbox{$\scriptstyle \infty$}}(z))_{11} \right.
\nonumber \\
+&\left. \left(\operatorname{Ai}(p_{a})(p_{a})^{1/4} \! + \!
\operatorname{Ai}^{\prime}(p_{a})(p_{a})^{-1/4} \right) \!
(\overset{e}{m}^{\raise-1.0ex\hbox{$\scriptstyle \infty$}}(z))_{12} \me^{\mi n
\mho^{e}_{j}} \right), \\
(m^{a,4}_{p}(z))_{12} :=& \, \sqrt{\smash[b]{\pi}} \, \me^{-\frac{\mi \pi}{6}}
\me^{-\frac{1}{2}n \xi^{e}_{a_{j}}(z)} \! \left(\mi \! \left(-\omega^{2}
\operatorname{Ai}(\omega p_{a})(p_{a})^{1/4} \! + \! \operatorname{Ai}^{\prime}
(\omega p_{a})(p_{a})^{-1/4} \right) \!
(\overset{e}{m}^{\raise-1.0ex\hbox{$\scriptstyle \infty$}}(z))_{11} \right.
\nonumber \\
\times &\left. \me^{-\mi n \mho^{e}_{j}} \! - \! \left(\omega^{2}
\operatorname{Ai}(\omega p_{a})(p_{a})^{1/4} \! + \! \operatorname{Ai}^{\prime}
(\omega p_{a})(p_{a})^{-1/4} \right) \!
(\overset{e}{m}^{\raise-1.0ex\hbox{$\scriptstyle \infty$}}(z))_{12} \right), \\
(m^{a,4}_{p}(z))_{21} :=& \, -\mi \sqrt{\smash[b]{\pi}} \, \me^{\frac{1}{2}n
\xi^{e}_{a_{j}}(z)} \! \left(\mi \! \left(\operatorname{Ai}(p_{a})(p_{a})^{
1/4} \! - \! \operatorname{Ai}^{\prime}(p_{a})(p_{a})^{-1/4} \right) \!
(\overset{e}{m}^{\raise-1.0ex\hbox{$\scriptstyle \infty$}}(z))_{21} \right.
\nonumber \\
+&\left. \left(\operatorname{Ai}(p_{a})(p_{a})^{1/4} \! + \!
\operatorname{Ai}^{\prime}(p_{a})(p_{a})^{-1/4} \right) \!
(\overset{e}{m}^{\raise-1.0ex\hbox{$\scriptstyle \infty$}}(z))_{22} \me^{\mi n
\mho^{e}_{j}} \right), \\
(m^{a,4}_{p}(z))_{22} :=& \, \sqrt{\smash[b]{\pi}} \, \me^{-\frac{\mi \pi}{6}}
\me^{-\frac{1}{2}n \xi^{e}_{a_{j}}(z)} \! \left(\mi \! \left(-\omega^{2}
\operatorname{Ai}(\omega p_{a})(p_{b})^{1/4} \! + \! \operatorname{Ai}^{\prime}
(\omega p_{a})(p_{a})^{-1/4} \right) \!
(\overset{e}{m}^{\raise-1.0ex\hbox{$\scriptstyle \infty$}}(z))_{21} \right.
\nonumber \\
\times &\left. \me^{-\mi n \mho^{e}_{j}} \! - \! \left(\omega^{2}
\operatorname{Ai}(\omega p_{a})(p_{a})^{1/4} \! + \! \operatorname{Ai}^{\prime}
(\omega p_{a})(p_{a})^{-1/4} \right) \!
(\overset{e}{m}^{\raise-1.0ex\hbox{$\scriptstyle \infty$}}(z))_{22} \right).
\end{align}
\end{dddd}
\begin{eeee}
Using limiting values, if necessary, all of the above (asymptotic) formulae 
for $\boldsymbol{\pi}_{2n}(z)$ and $\int_{\mathbb{R}} \boldsymbol{\pi}_{2n}(s) 
\me^{-n \widetilde{V}(s)}(s \! - \! z)^{-1} \, \tfrac{\md s}{2 \pi \mi}$ have 
a natural interpretation on the real and imaginary axes. \hfill $\blacksquare$
\end{eeee}
\begin{dddd}
Let all the conditions stated in Theorem~{\rm 2.3.1} be valid, and let 
$\overset{e}{\operatorname{Y}} \colon \mathbb{C} \setminus \mathbb{R} \! 
\to \! \operatorname{SL}_{2}(\mathbb{C})$ be the unique solution of 
{\rm \pmb{RHP1}}. Let $H^{(m)}_{k}$, $(m,k) \! \in \! \mathbb{Z} \times 
\mathbb{N}$, be the Hankel determinants associated with the bi-infinite, 
real-valued, strong moment sequence $\left\lbrace c_{k} \! = \! \int_{
\mathbb{R}}s^{k} \me^{-n \widetilde{V}(s)} \, \md s, \, n \! \in \! \mathbb{N} 
\right\rbrace_{k \in \mathbb{Z}}$ defined in Equations~{\rm (1.1)}, and let 
$\boldsymbol{\pi}_{2n}(z)$ be the even degree monic orthogonal $L$-polynomial 
defined in Lemma~{\rm 2.2.1}, that is, $\boldsymbol{\pi}_{2n}(z) \! := \! 
(\overset{e}{\mathrm{Y}}(z))_{11}$, with $n \! \to \! \infty$ asymptotics 
(in the entire complex plane) given by Theorem~{\rm 2.3.1}.
Then,
\begin{equation}
(\xi^{(2n)}_{n})^{2} \! = \! \dfrac{1}{\norm{\boldsymbol{\pi}_{2n}(\pmb{\cdot}
)}_{\mathscr{L}}^{2}} \! = \! \dfrac{H^{(-2n)}_{2n}}{H^{(-2n)}_{2n+1}}
\underset{n \to \infty}{=} \dfrac{\me^{-n \ell_{e}}}{\pi} \Xi^{\flat} \!
\left(1 \! + \! \dfrac{1}{n} \Xi^{\flat}(\mathfrak{Q}^{\flat})_{12} \! + \!
\mathcal{O} \! \left(\dfrac{c^{\flat}(n)}{n^{2}} \right) \right),
\end{equation}
where
\begin{gather}
\Xi^{\flat} \! := \! 2 \! \left(\sum_{k=1}^{N+1} \! \left(a_{k}^{e} \! - \!
b_{k-1}^{e} \right) \right)^{-1} \dfrac{\boldsymbol{\theta}^{e}(\boldsymbol{
u}^{e}_{+}(\infty) \! - \! \frac{n}{2 \pi} \boldsymbol{\Omega}^{e} \! + \!
\boldsymbol{d}_{e}) \boldsymbol{\theta}^{e}(-\boldsymbol{u}^{e}_{+}(\infty)
\! + \! \boldsymbol{d}_{e})}{\boldsymbol{\theta}^{e}(-\boldsymbol{u}^{e}_{+}
(\infty) \! - \! \frac{n}{2 \pi} \boldsymbol{\Omega}^{e} \! + \! \boldsymbol{
d}_{e}) \boldsymbol{\theta}^{e}(\boldsymbol{u}^{e}_{+}(\infty) \! + \!
\boldsymbol{d}_{e})}, \\
\mathfrak{Q}^{\flat} \! := \! 2 \mi \sum_{j=1}^{N+1} \! \left(\dfrac{(\mathscr{
B}^{e}(a_{j}^{e}) \widehat{\alpha}^{e}_{0}(a_{j}^{e}) \! - \! \mathscr{A}^{e}
(a_{j}^{e}) \widehat{\alpha}^{e}_{1}(a_{j}^{e}))}{(\widehat{\alpha}^{e}_{0}
(a_{j}^{e}))^{2}} \! + \! \dfrac{(\mathscr{B}^{e}(b_{j-1}^{e}) \widehat{
\alpha}^{e}_{0}(b_{j-1}^{e}) \! - \! \mathscr{A}^{e}(b_{j-1}^{e}) \widehat{
\alpha}^{e}_{1}(b_{j-1}^{e}))}{(\widehat{\alpha}^{e}_{0}(b_{j-1}^{e}))^{2}}
\right),
\end{gather}
$(\mathfrak{Q}^{\flat})_{12}$ denotes the $(1 \, 2)$-element of $\mathfrak{
Q}^{\flat}$, $c^{\flat}(n) \! =_{n \to \infty} \! \mathcal{O}(1)$, and all the 
relevant parameters are defined in Theorem~{\rm 2.3.1:} asymptotics for $\xi^{
(2n)}_{n}$ are obtained by taking the positive square root of both sides of 
Equation~{\rm (2.118)}. Furthermore, the $n \! \to \! \infty$ asymptotic 
expansion (in the entire complex plane) for the even degree orthonormal 
$L$-polynomial,
\begin{equation}
\phi_{2n}(z) \! = \! \xi^{(2n)}_{n} \boldsymbol{\pi}_{2n}(z),
\end{equation}
to $\mathcal{O}(n^{-2})$, is given by the (scalar) multiplication of the $n 
\! \to \! \infty$ asymptotics of $\boldsymbol{\pi}_{2n}(z)$ and $\xi^{(2n)}_{
n}$ stated, respectively, in Theorem~{\rm 2.3.1} and 
Equations~{\rm (2.118)--(2.120)}.
\end{dddd}
\begin{eeee}
Since, {}from general theory (cf. Section~1), and, by construction (cf. 
Equations~(1.2) and~(1.8)), $\xi^{(2n)}_{n} \! > \! 0$, it follows, 
incidentally, {}from Theorem~2.3.2, Equations~(2.118)--(2.120), that: (i) 
$\Xi^{\flat} \! > \! 0$; and (ii) $\Im ((\mathfrak{Q}^{\flat})_{12}) \! = \! 
0$. \hfill $\blacksquare$
\end{eeee}
\section{The Equilibrium Measure, the Variational Problem, and the 
Tr\-a\-n\-s\-f\-o\-r\-m\-e\-d RHP}
In this section, the detailed analysis of the `even degree' variational 
problem, and the associated `even' equilibrium measure, is undertaken (see 
Lemmas~3.1--3.3 and Lemma~3.5), including the discussion of the corresponding 
$g$-function, denoted, herein, as $g^{e}$, and \textbf{RHP1}, that is, 
$(\overset{e}{\mathrm{Y}}(z),\mathrm{I} \! + \! \exp (-n \widetilde{V}(z)) 
\sigma_{+},\mathbb{R})$, is reformulated as an equivalent\footnote{If 
there are two RHPs, $(\mathcal{Y}_{1}(z),\upsilon_{1}(z),\varGamma_{1})$ and 
$(\mathcal{Y}_{2}(z),\upsilon_{2}(z),\varGamma_{2})$, say, with $\varGamma_{2} 
\subset \varGamma_{1}$ and $\upsilon_{1}(z) \! \! \upharpoonright_{\varGamma_{
1} \setminus \varGamma_{2}} \! =_{n \to \infty} \! \mathrm{I} \! + \! o(1)$, 
then, within the BC framework \cite{a84}, and modulo $o(1)$ estimates, their 
solutions, $\mathcal{Y}_{1}$ and $\mathcal{Y}_{2}$, respectively, are 
(asymptotically) equal.}, auxiliary RHP (see Lemma~3.4). The proofs of 
Lemmas~3.1--3.3 are modelled on the calculations of Saff-Totik (\cite{a55}, 
Chapter~1), Deift (\cite{a90}, Chapter~6), and Johansson \cite{a91}.

One begins by establishing the existence of the `even' equilibrium measure, 
$\mu_{V}^{e}$ $(\in \! \mathcal{M}_{1}(\mathbb{R}))$.
\begin{ccccc}
Let the external field $\widetilde{V} \colon \mathbb{R} \setminus \{0\} \! \to 
\! \mathbb{R}$ satisfy conditions~{\rm (2.3)--(2.5)}, and set $w^{e}(z) \! := 
\! \me^{-\widetilde{V}(z)}$. For $\mu^{e} \! \in \! \mathcal{M}_{1}(\mathbb{R}
)$, define the weighted energy functional $\mathrm{I}_{V}^{e}[\mu^{e}] \colon 
\mathcal{M}_{1}(\mathbb{R}) \! \to \! \mathbb{R}$,
\begin{equation*}
\mathrm{I}_{V}^{e}[\mu^{e}] \! := \! \iint_{\mathbb{R}^{2}} \ln \! \left(
\vert s \! - \! t \vert^{2} \vert st \vert^{-1}w^{e}(s)w^{e}(t) \right)^{-1}
\md \mu^{e}(s) \, \md \mu^{e}(t),
\end{equation*}
and consider the minimisation problem
\begin{equation*}
E_{V}^{e} \! = \! \inf \left\lbrace \mathstrut \mathrm{I}_{V}^{e}[\mu^{e}]; 
\, \mu^{e} \! \in \! \mathcal{M}_{1}(\mathbb{R}) \right\rbrace.
\end{equation*}
Then: {\rm (1)} $E_{V}^{e}$ is finite; {\rm (2)} $\exists \, \, \mu_{V}^{e} 
\in \! \mathcal{M}_{1}(\mathbb{R})$ such that $\mathrm{I}_{V}^{e}[\mu_{V}^{e}] 
\! = \! E_{V}^{e}$ (the infimum is attained), and $\mu_{V}^{e}$ has finite 
weighted logarithmic energy $(-\infty \! < \! \mathrm{I}_{V}^{e}[\mu_{V}^{e}] 
\! < \! +\infty);$ and {\rm (3)} $J_{e} \! := \! \operatorname{supp}(\mu_{V}^{
e})$ is compact, $J_{e} \subset \{\mathstrut z; \, w^{e}(z) \! > \! 0\}$, and 
$J_{e}$ has positive logarithmic capacity, that is, $\operatorname{cap}(J_{e}) 
\! := \! \exp \! \left(-\inf \{\mathrm{I}_{V}^{e}[\mu^{e}]; \, \mu^{e} \! \in 
\! \mathcal{M}_{1}(J_{e})\} \right) \! > \! 0$.
\end{ccccc}

\emph{Proof.} Let $\mu^{e} \! \in \! \mathcal{M}_{1}(\mathbb{R})$, and set 
$w^{e}(z) \! := \! \exp (-\widetilde{V}(z))$, where $\widetilde{V} \colon 
\mathbb{R} \setminus \{0\} \! \to \! \mathbb{R}$ satisfies 
conditions~(2.3)--(2.5). {}From the definition of $\mathrm{I}_{V}^{e}[\mu^{e}
]$ given in the Lemma, one shows that
\begin{align*}
\mathrm{I}_{V}^{e}[\mu^{e}] &= \iint_{\mathbb{R}^{2}} \! \left(\ln (\vert s \!
- \! t \vert^{-1}) \! + \! \ln(\vert s^{-1} \! - \! t^{-1} \vert^{-1}) \right)
\md \mu^{e}(s) \, \md \mu^{e}(t) \! + \! 2 \int_{\mathbb{R}} \widetilde{V}(s)
\, \md \mu^{e}(s) \\
&=: \iint_{\mathbb{R}^{2}}K_{V}^{e}(s,t) \, \md \mu^{e}(s) \, \md \mu^{e}(t),
\end{align*}
where (the symmetric kernel)
\begin{equation*}
K_{V}^{e}(s,t) \! = \! K_{V}^{e}(t,s) \! := \! \ln (\vert s \! - \! t \vert^{-
1}) \! + \! \ln (\vert s^{-1} \! - \! t^{-1} \vert^{-1}) \! + \! \widetilde{V}
(s) \! + \! \widetilde{V}(t)
\end{equation*}
(of course, the definition of $\mathrm{I}_{V}^{e}[\mu^{e}]$ only makes sense 
provided both integrals exist and are finite). Recall the following 
inequalities (see, for example, Chapter~6 of \cite{a90}): $\vert s \! - \! t 
\vert \! \leqslant \! (1 \! + \! s^{2})^{1/2}(1 \! + \! t^{2})^{1/2}$ and 
$\vert s^{-1} \! - \! t^{-1} \vert \! \leqslant \! (1 \! + \! s^{-2})^{1/2}
(1 \! + \! t^{-2})^{1/2}$, $s,t \! \in \! \mathbb{R}$, whence $\ln (\vert s 
\! - \! t \vert^{-1}) \! \geqslant \! -\tfrac{1}{2} \ln (1 \! + \! s^{2}) \! 
- \! \tfrac{1}{2} \ln (1 \! + \! t^{2})$ and $\ln (\vert s^{-1} \! - \! t^{-
1} \vert^{-1}) \! \geqslant \! -\tfrac{1}{2} \ln (1 \! + \! s^{-2}) \! - \! 
\tfrac{1}{2} \ln (1 \! + \! t^{-2})$; thus,
\begin{equation*}
K_{V}^{e}(s,t) \! \geqslant \! \dfrac{1}{2} \! \left(2 \widetilde{V}(s) \!
- \! \ln (s^{2} \! + \! 1) \! - \! \ln (s^{-2} \! + \! 1) \right) \! + \!
\dfrac{1}{2} \! \left(2 \widetilde{V}(t) \! - \! \ln (t^{2} \! + \! 1) \! -
\! \ln (t^{-2} \! + \! 1) \right).
\end{equation*}
Recalling conditions~(2.3)--(2.5) for the external field $\widetilde{V} \colon
\mathbb{R} \setminus \{0\} \! \to \! \mathbb{R}$, in particular, $\exists \,
\, \delta_{1} \! > \! 0$ (resp., $\exists \, \, \delta_{2} \! > \! 0)$ such
that $\widetilde{V}(x) \! \geqslant \! (1 \! + \! \delta_{1}) \ln (x^{2} \! +
\! 1)$ (resp., $\widetilde{V}(x) \! \geqslant \! (1 \! + \! \delta_{2}) \ln
(x^{-2} \! + \! 1))$ for sufficiently large $\vert x \vert$ (resp., small
$\vert x \vert)$, it follows that $2 \widetilde{V}(x) \! - \! \ln (x^{2} \! +
\! 1) \! - \! \ln (x^{-2} \! + \! 1) \! \geqslant \! C_{V}^{e} \! > \!
-\infty$, whence $K_{V}^{e}(s,t) \! \geqslant \! C_{V}^{e}$ $(> \! -\infty)$,
which shows that $K_{V}^{e}(s,t)$ is bounded {}from below (on $\mathbb{R}^{2}
)$; hence,
\begin{equation*}
\mathrm{I}_{V}^{e}[\mu^{e}] \! \geqslant \! \iint_{\mathbb{R}^{2}}C_{V}^{e} \,
\md \mu^{e}(s) \, \md \mu^{e}(t) \! = \! C_{V}^{e} \underbrace{\int_{\mathbb{
R}} \md \mu^{e}(s)}_{= \, 1} \, \, \underbrace{\int_{\mathbb{R}} \md \mu^{e}
(t)}_{= \, 1} \! \geqslant \! C_{V}^{e} \quad (> \! -\infty).
\end{equation*}
It follows {}from the above inequality and the definition of $E_{V}^{e}$ 
stated in the Lemma that, $\forall \, \, \mu^{e} \! \in \! \mathcal{M}_{1}
(\mathbb{R})$, $E_{V}^{e} \! \geqslant \! C_{V}^{e} \! > \! -\infty$, which 
shows that $E_{V}^{e}$ is bounded {}from below. Let $\varepsilon$ be an 
arbitrarily fixed, sufficiently small positive real number, and set $\Sigma_{e,
\varepsilon} \! := \! \{\mathstrut z; \, w^{e}(z) \! \geqslant \! \varepsilon 
\}$; then, $\Sigma_{e,\varepsilon}$ is compact, and $\Sigma_{e,0} \! := \! 
\cup_{l=1}^{\infty} \Sigma_{e,1/l} \! = \! \cup_{l=1}^{\infty} \{\mathstrut 
z; \, w^{e}(z) \! > \! l^{-1}\} \! = \! \{\mathstrut z; \, w^{e}(z) \! > \! 
0\}$. Since, for $\widetilde{V} \colon \mathbb{R} \setminus \{0\} \! \to \! 
\mathbb{R}$ satisfying conditions (2.3)--(2.5), $w^{e}$ is an \emph{admissible 
weight} \cite{a55}, in which case $\Sigma_{e,0}$ has positive logarithmic 
capacity, that is, $\operatorname{cap}(\Sigma_{e,0}) \! = \! \exp (-\inf \{
\mathstrut \mathrm{I}_{V}^{e}[\mu^{e}]; \, \mu^{e} \! \in \! \mathcal{M}_{1}
(\Sigma_{e,0})\}) \! > \! 0$, it follows that $\exists \, \, l^{\ast} \! \in 
\! \mathbb{N}$ such that $\operatorname{cap}(\Sigma_{e,1/l^{\ast}}) \! = \! 
\exp (-\inf \{\mathstrut \mathrm{I}_{V}^{e}[\mu^{e}]; \, \mu^{e} \! \in \! 
\mathcal{M}_{1}(\Sigma_{e,1/l^{\ast}})\}) \! > \! 0$, which, in turn, means 
that there exists a probability measure, $\mu^{e}_{l^{\ast}}$, say, with 
$\operatorname{supp}(\mu^{e}_{l^{\ast}}) \subseteq \Sigma_{e,1/l^{\ast}}$, 
such that $\iint_{\Sigma_{e,1/l^{\ast}}^{2}} \ln (\vert s \! - \! t \vert^{-
2} \vert st \vert) \, \md \mu^{e}_{l^{\ast}}(s) \, \md \mu^{e}_{l^{\ast}}(t) 
\! < \! +\infty$, where $\Sigma_{e,1/l^{\ast}}^{2} \! = \! \Sigma_{e,1/l^{
\ast}} \times \Sigma_{e,1/l^{\ast}}$ $(\subseteq \mathbb{R}^{2})$. For $z \! 
\in \! \operatorname{supp}(\mu^{e}_{l^{\ast}}) \! \subseteq \! \Sigma_{e,
1/l^{\ast}}$, it follows that $w^{e}(z) \! \geqslant \! 1/l^{\ast}$, whence 
$\iint_{\Sigma_{e,1/l^{\ast}}^{2}} \ln (w^{e}(s)w^{e}(t))^{-1} \, \md \mu^{
e}_{l^{\ast}}(s) \, \md \mu^{e}_{l^{\ast}}(t) \! \leqslant \! 2 \ln (l^{\ast}) 
\! < \! +\infty$ $\Rightarrow$
\begin{equation*}
\mathrm{I}_{V}^{e}[\mu_{l^{\ast}}^{e}] \! = \! \iint_{\Sigma_{e,1/l^{\ast}}^{
2}} \ln \! \left(\vert s \! - \! t \vert^{2} \vert st \vert^{-1}w^{e}(s)w^{e}
(t) \right)^{-1} \md \mu^{e}_{l^{\ast}}(s) \, \md \mu^{e}_{l^{\ast}}(t) \! <
\! +\infty;
\end{equation*}
thus, it follows that $E_{V}^{e} \! := \! \inf \lbrace \mathstrut \mathrm{I}_{
V}^{e}[\mu^{e}]; \, \mu^{e} \! \in \! \mathcal{M}_{1}(\mathbb{R}) \rbrace$ is 
finite (see, also, below).

Choose a sequence of probability measures $\{\mu^{e}_{n}\}_{n=1}^{\infty}$ in 
$\mathcal{M}_{1}(\mathbb{R})$ such that $\mathrm{I}_{V}^{e}[\mu^{e}_{n}] \! 
\leqslant \! E_{V}^{e} \! + \! \tfrac{1}{n}$. {}From the analysis above, it 
follows that
\begin{align*}
\mathrm{I}_{V}^{e}[\mu_{n}^{e}] =& \iint_{\mathbb{R}^{2}}K_{V}^{e}(s,t) \, \md
\mu^{e}_{n}(s) \, \md \mu^{e}_{n}(t) \! \geqslant \! \iint_{\mathbb{R}^{2}} \!
\left(\tfrac{1}{2}(2 \widetilde{V}(s) \! - \! \ln (s^{2} \! + \! 1) \! - \!
\ln (s^{-2} \! + \! 1)) \right. \\
+& \left. \tfrac{1}{2}(2 \widetilde{V}(t) \! - \! \ln (t^{2} \! + \! 1) \! -
\! \ln (t^{-2} \! + \! 1) \right) \! \md \mu_{n}^{e}(s) \, \md \mu_{n}^{e}(t).
\end{align*}
Set
\begin{equation*}
\widehat{\psi}_{V}^{e}(z) \! := \! 2 \widetilde{V}(z) \! - \! \ln (z^{2} \! +
\! 1) \! - \! \ln (z^{-2} \! + \! 1).
\end{equation*}
Then $\mathrm{I}_{V}^{e}[\mu_{n}^{e}] \! \geqslant \! \int_{\mathbb{R}} 
\widehat{\psi}_{V}^{e}(s) \, \md \mu_{n}^{e}(s)$ $\Rightarrow$ $E_{V}^{e} \! 
+ \! \tfrac{1}{n} \! \geqslant \! \int_{\mathbb{R}} \widehat{\psi}_{V}^{e}(s) 
\, \md \mu_{n}^{e}(s)$. Recalling that $\exists \, \, \delta_{1} \! > \! 0$ 
(resp., $\exists \, \, \delta_{2} \! > \! 0)$ such that $\widetilde{V}(x) \! 
\geqslant \! (1 \! + \! \delta_{1}) \ln (x^{2} \! + \! 1)$ (resp., 
$\widetilde{V}(x) \! \geqslant \! (1 \! + \! \delta_{2}) \ln (x^{-2} \! + \! 
1))$ for sufficiently large $\vert x \vert$ (resp., small $\vert x \vert)$, 
it follows that, for any $b_{e} \! > \! 0$, $\exists \, \, M_{e} \! > \! 1$ 
such that $\widehat{\psi}_{V}^{e}(z) \! > \! b_{e} \, \, \forall \, \, z \! 
\in \! \{\vert z \vert \! \geqslant \! M_{e}\} \cup \{\vert z \vert \! 
\leqslant \! M_{e}^{-1}\} \! =: \! \mathfrak{D}_{e}$, which implies that
\begin{align*}
E_{V}^{e} \! + \! \tfrac{1}{n} \geqslant& \int_{\mathbb{R}} \widehat{\psi}_{
V}^{e}(s) \, \md \mu_{n}^{e}(s) \! = \! \int_{\mathfrak{D}_{e}} \underbrace{
\widehat{\psi}_{V}^{e}(s)}_{> \, b_{e}} \, \md \mu_{n}^{e}(s) \! + \! \int_{
\mathbb{R} \setminus \mathfrak{D}_{e}} \underbrace{\widehat{\psi}_{V}^{e}(s)
}_{\geqslant \, -\vert C_{V}^{e} \vert} \, \md \mu_{n}^{e}(s) \\
\geqslant& \, b_{e} \int_{\mathfrak{D}_{e}} \md \mu^{e}_{n}(s) \! - \! \vert
C_{V}^{e} \vert \underbrace{\int_{\mathbb{R} \setminus \mathfrak{D}_{e}} \md
\mu_{n}^{e}(s)}_{\in \, [0,1]} \! \geqslant \! b_{e} \int_{\mathfrak{D}_{e}}
\md \mu_{n}^{e}(s) \! - \! \vert C_{V}^{e} \vert;
\end{align*}
thus,
\begin{equation*}
\int_{\mathfrak{D}_{e}} \md \mu_{n}^{e}(s) \! \leqslant \! b_{e}^{-1} \!
\left(E_{V}^{e} \! + \! \vert C_{V}^{e} \vert \! + \! \dfrac{1}{n} \right),
\end{equation*}
whence
\begin{equation*}
\limsup_{n \to \infty} \int_{\mathfrak{D}_{e}} \md \mu_{n}^{e}(s) \! \leqslant
\! \limsup_{n \to \infty} \! \left(b_{e}^{-1} \! \left(E_{V}^{e} \! + \! \vert
C_{V}^{e} \vert \! + \! \dfrac{1}{n} \right) \right).
\end{equation*}
By the Archimedean property, it follows that, $\forall \, \, \epsilon_{o} \! 
> \! 0$, $\exists \, \, N \! \in \! \mathbb{N}$ such that, $\forall \, \, n 
\! > \! N \Rightarrow n^{-1} \! < \! \epsilon_{o}$; thus, choosing $b_{e} \! 
= \! \epsilon^{-1}(E_{V}^{e} \! + \! \vert C_{V}^{e} \vert \! + \! \epsilon_{
o})$, where $\epsilon$ is some arbitrarily fixed, sufficiently small positive 
real number, it follows that $\limsup_{n \to \infty} \int_{\mathfrak{D}_{e}} 
\md \mu_{n}^{e}(s) \! \leqslant \! \epsilon$ $\Rightarrow$ the sequence of 
probability measures $\{\mu_{n}^{e}\}_{n=1}^{\infty}$ in $\mathcal{M}_{1}
(\mathbb{R})$ is \emph{tight} \cite{a91} (that is, given $\epsilon \! > \! 
0$, $\exists \, \, M \! > \! 1$ such that $\limsup_{n \to \infty} \mu_{n}^{e}
(\{\vert s \vert \! \geqslant \! M\} \cup \{\vert s \vert \! \leqslant \! 
M^{-1}\}) \! := \! \limsup_{n \to \infty} \int_{\{\vert s \vert \geqslant M\} 
\cup \{\vert s \vert \leqslant M^{-1}\}} \md \mu_{n}^{e}(s) \! \leqslant \! 
\epsilon)$. Since the sequence of probabilty measures $\{\mu_{n}^{e}\}_{n=
1}^{\infty}$ in $\mathcal{M}_{1}(\mathbb{R})$ is tight, by a Helly Selection 
Theorem, there exists a $(\mathrm{weak}^{\ast}$ convergent) subsequence of 
probability measures $\{\mu^{e}_{n_{k}}\}_{k=1}^{\infty}$ in $\mathcal{M}_{1}
(\mathbb{R})$ converging (weakly) to a probability measure $\mu^{e} \! \in 
\! \mathcal{M}_{1}(\mathbb{R})$, symbolically $\mu^{e}_{n_{k}} \! \overset{
\ast}{\to} \! \mu^{e}$ as $k \! \to \! \infty$\footnote{A sequence of 
probability measures $\{\mu_{n}\}_{n=1}^{\infty}$ in $\mathcal{M}_{1}(D)$ 
is said to \emph{converge weakly} as $n \! \to \! \infty$ to $\mu \! \in \! 
\mathcal{M}_{1}(D)$, symbolically $\mu_{n} \! \overset{\ast}{\to} \! \mu$, 
if $\mu_{n}(f) \! := \! \int_{D}f(s) \, \md \mu_{n}(s) \! \to \! \int_{D}f(s) 
\, \md \mu (s) \! =: \mu (f)$ as $n \! \to \! \infty$ $\forall \, \, f \! \in 
 \! \mathbf{\mathrm{C}_{b}^{0}}(D)$, where $\mathbf{\mathrm{C}_{b}^{0}}(D)$ 
denotes the set of all bounded, continuous functions on $D$ with compact 
support.}. One now shows that, if $\mu_{n}^{e} \! \overset{\ast}{\to} \! 
\mu^{e}$, $\mu_{n}^{e},\mu^{e} \! \in \! \mathcal{M}_{1}(\mathbb{R})$, then 
$\liminf_{n \to \infty} \mathrm{I}_{V}^{e}[\mu_{n}^{e}] \! \geqslant \! 
\mathrm{I}_{V}^{e}[\mu^{e}]$. Since $w^{e}$ is continuous, thus upper 
semi-continuous \cite{a55}, there exists a sequence $\{w^{e}_{m}\}_{m=1}^{
\infty}$ (resp., $\{\widetilde{V}_{m}\}_{m=1}^{\infty})$ of continuous 
functions on $\mathbb{R}$ such that $w^{e}_{m+1} \! \leqslant \! w^{e}_{m}$ 
(resp., $\widetilde{V}_{m+1} \! \geqslant \! 
\widetilde{V}_{m})$\footnote{Adding a suitable constant, if necessary, which 
does not change $\mu_{m}^{e}$, or the regularity of $\widetilde{V} \colon 
\mathbb{R} \setminus \lbrace 0 \rbrace \! \to \! \mathbb{R}$, one may assume 
that $\widetilde{V} \! \geqslant \! 0$ and $\widetilde{V}_{m} \! \geqslant \! 
0$, $m \! \in \! \mathbb{N}$.}, $m \! \in \! \mathbb{N}$, and $w^{e}_{m}(z) 
\! \searrow \! w^{e}(z)$ (resp., $\widetilde{V}_{m}(z) \! \nearrow \! 
\widetilde{V}(z))$ as $m \! \to \! \infty$ for every $z \! \in \! \mathbb{
R}$; in particular,
\begin{equation*}
\mathrm{I}_{V}^{e}[\mu_{n}^{e}] \! = \! \iint_{\mathbb{R}^{2}}K_{V}^{e}(s,t)
\, \md \mu_{n}^{e}(s) \, \md \mu_{n}^{e}(t) \! \geqslant \! \iint_{\mathbb{
R}^{2}}K_{V_{m}}^{e}(s,t) \, \md \mu_{n}^{e}(s) \, \md \mu_{n}^{e}(t).
\end{equation*}
For arbitrary $q \! \in \! \mathbb{R}$, $\mathrm{I}_{V}^{e}[\mu_{n}^{e}] 
\! \geqslant \! \iint_{\mathbb{R}^{2}}p^{e}(s,t) \, \md \mu_{n}^{e}(s) \, 
\md \mu_{n}^{e}(t)$, where $p^{e}(s,t) \! = \! p^{e}(t,s) \! := \! \min 
\left\lbrace q,K_{V_{m}}^{e}(s,t) \right\rbrace$ (bounded and continuous on 
$\mathbb{R}^{2})$. Recall that $\{\mu_{n}^{e}\}_{n=1}^{\infty}$ is tight in 
$\mathcal{M}_{1}(\mathbb{R})$. For $M_{e} \! > \! 1$, let $h_{M}^{e}(x) \! 
\in \! \mathbf{\mathrm{C}_{b}^{0}}(\mathbb{R})$ be such that:
\begin{compactenum}
\item[(i)] $h_{M}^{e}(x) \! = \! 1$, $x \! \in \! [-M_{e},-M_{e}^{-1}] \cup
[M_{e}^{-1},M_{e}] \! =: \! \mathfrak{D}_{M_{e}}$;
\item[(ii)] $h_{M}^{e}(x) \! = \! 0$, $x \! \in \! \mathbb{R} \setminus
\mathfrak{D}_{M_{e}+1}$; and
\item[(iii)] $0 \! \leqslant \! h_{M}^{e}(x) \! \leqslant \! 1$, $x \! \in \!
\mathbb{R}$.
\end{compactenum}
Note the decomposition $\iint_{\mathbb{R}^{2}}p^{e}(t,s) \, \md \mu_{n}^{e}(t)
\, \md \mu_{n}^{e}(s) \! = \! I_{a} \! + \! I_{b} \! + \! I_{c}$, where
\begin{align*}
I_{a} :=& \iint_{\mathbb{R}^{2}}p^{e}(t,s)(1 \! - \! h_{M}^{e}(s)) \, \md
\mu^{e}_{n}(t) \, \md \mu_{n}^{e}(s), \\
I_{b} :=& \iint_{\mathbb{R}^{2}}p^{e}(t,s)h_{M}^{e}(s)(1 \! - \! h_{M}^{e}(t))
\, \md \mu^{e}_{n}(t) \, \md \mu_{n}^{e}(s), \\
I_{c} :=& \iint_{\mathbb{R}^{2}}p^{e}(t,s)h_{M}^{e}(t)h_{M}^{e}(s) \, \md
\mu^{e}_{n}(t) \, \md \mu_{n}^{e}(s).
\end{align*}
One shows that
\begin{align*}
\vert I_{a} \vert \leqslant& \iint_{\mathbb{R}^{2}} \vert p^{e}(t,s) \vert
(1 \! - \! h_{M}^{e}(s)) \, \md \mu_{n}^{e}(t) \, \md \mu_{n}^{e}(s) \\
\leqslant& \, \sup_{(t,s) \in \mathbb{R}^{2}} \vert p^{e}(t,s) \vert
\underbrace{\int_{\mathbb{R}} \md \mu_{n}^{e}(t)}_{= \, 1} \! \left(\int_{
\mathfrak{D}_{M_{e}}} \underbrace{(1 \! - \! h_{M}^{e}(s))}_{= \, 0} \, \md
\mu_{n}^{e}(s) \! + \! \int_{\mathbb{R} \setminus \mathfrak{D}_{M_{e}+1}}
(1 \! - \! \underbrace{h_{M}^{e}(s)}_{= \, 0}) \, \md \mu_{n}^{e}(s) \right),
\end{align*}
whence
\begin{equation*}
\limsup_{n \to \infty} \vert I_{a} \vert \! \leqslant \! \sup_{(t,s) \in
\mathbb{R}^{2}} \vert p^{e}(t,s) \vert \underbrace{\limsup_{n \to \infty}
\int_{\mathbb{R} \setminus \mathfrak{D}_{M_{e}+1}} \, \md \mu_{n}^{e}(s)}_{
\leqslant \, \epsilon} \! \leqslant \! \epsilon \sup_{(t,s) \in \mathbb{R}^{
2}} \vert p^{e}(t,s) \vert;
\end{equation*}
similarly,
\begin{equation*}
\limsup_{n \to \infty} \vert I_{b} \vert \! \leqslant \! \epsilon \sup_{(t,s)
\in \mathbb{R}^{2}} \vert p^{e}(t,s) \vert.
\end{equation*}
Since $p^{e}(t,s)$ is continuous and bounded on $\mathbb{R}^{2}$, there 
exists, by a generalisation of the Stone-Weierst\-r\-a\-s\-s Theorem (for the 
single-variable case), a polynomial in two variables, $p(t,s)$, say, with $p 
(t,s) \! = \! \sum_{i \geqslant i_{o}} \sum_{j \geqslant j_{o}} \gamma_{ij} 
t^{i}s^{j}$, such that $\vert p^{e}(t,s) \! - \! p(t,s) \vert \! \leqslant 
\! \epsilon$; thus,
\begin{equation*}
\vert h_{M}^{e}(t)h_{M}^{e}(s)p^{e}(t,s) \! - \! h_{M}^{e}(t)h_{M}^{e}(s)
p(t,s) \vert \! \leqslant \! \epsilon, \quad t,s \! \in \! \mathbb{R}.
\end{equation*}
Rewrite $I_{c}$ as
\begin{align*}
I_{c} &= \iint_{\mathbb{R}^{2}}h_{M}^{e}(s)h_{M}^{e}(t)p(t,s) \, \md
\mu_{n}^{e}(t) \, \md \mu_{n}^{e}(s) \! + \! \iint_{\mathbb{R}^{2}}h_{M}^{e}
(s)h_{M}^{e}(t)(p^{e}(t,s) \! - \! p(t,s)) \, \md \mu_{n}^{e}(t) \, \md \mu_{
n}^{e}(s) \\
&=: I_{c}^{\alpha} \! + \! I_{c}^{\beta}.
\end{align*}
One now shows that
\begin{align*}
\vert I_{c}^{\beta} \vert \leqslant& \iint_{\mathbb{R}^{2}}h_{M}^{e}(s)h_{M}^{
e}(t) \underbrace{\vert p^{e}(t,s) \! - \! p(t,s) \vert}_{\leqslant \,
\epsilon} \, \md \mu_{n}^{e}(t) \, \md \mu_{n}^{e}(s) \! \leqslant \! \epsilon
\int_{\mathbb{R}}h_{M}^{e}(s) \, \md \mu_{n}^{e}(s) \int_{\mathbb{R}}h_{M}^{e}
(t) \, \md \mu_{n}^{e}(t) \\
\leqslant& \, \epsilon \! \left(\int_{\mathfrak{D}_{M_{e}}} \underbrace{h_{
M}^{e}(s)}_{= \, 1} \, \md \mu_{n}^{e}(s) \! + \! \int_{\mathbb{R} \setminus
\mathfrak{D}_{M_{e}+1}} \underbrace{h_{M}^{e}(s)}_{= \, 0} \, \md \mu_{n}^{e}
(s) \right)^{2} \! \leqslant \! \epsilon \! \left(\int_{\mathfrak{D}_{M_{e}}}
\md \mu_{n}^{e}(s) \right)^{2} \\
\leqslant& \, \epsilon \left(\int_{\mathbb{R}} \md \mu_{n}^{e}(s) \right)^{2}
\! \leqslant \! \epsilon,
\end{align*}
and
\begin{align*}
I_{c}^{\alpha} &= \iint_{\mathbb{R}^{2}}h_{M}^{e}(s)h_{M}^{e}(t) \sum_{i
\geqslant i_{o}} \sum_{j \geqslant j_{o}} \gamma_{ij}t^{i}s^{j} \, \md \mu_{
n}^{e}(t) \, \md \mu_{n}^{e}(s) \\
&= \, \sum_{i \geqslant i_{o}} \sum_{j \geqslant j_{o}} \gamma_{ij} \! \left(
\int_{\mathbb{R}}h_{M}^{e}(t)t^{i} \, \md \mu_{n}^{e}(t) \right) \! \left(
\int_{\mathbb{R}}h_{M}^{e}(s)s^{j} \, \md \mu_{n}^{e}(s) \right) \\
&\to \, \sum_{i \geqslant i_{o}} \sum_{j \geqslant j_{o}} \gamma_{ij} \! \left(
\int_{\mathbb{R}}h_{M}^{e}(t)t^{i} \, \md \mu^{e}(t) \right) \! \left(\int_{
\mathbb{R}}h_{M}^{e}(s)s^{j} \, \md \mu^{e}(s) \right) \quad \quad
(\text{since} \, \, \mu_{n}^{e} \! \overset{\ast}{\to} \! \mu^{e} \, \,
\text{as} \, \, n \! \to \! \infty) \\
&= \, \iint_{\mathbb{R}^{2}} \! \left(\sum_{i \geqslant i_{o}} \sum_{j
\geqslant j_{o}} \gamma_{ij}t^{i}s^{j} \right) \! h_{M}^{e}(t)h_{M}^{e}(s) \,
\md \mu^{e}(t) \, \md \mu^{e}(s),
\end{align*}
whence, recalling that $p(t,s) \! = \! \sum_{i \geqslant i_{o}} \sum_{j
\geqslant j_{o}} \gamma_{ij}t^{i}s^{j}$, it follows that
\begin{equation*}
I_{c}^{\alpha} \! = \! \iint_{\mathbb{R}^{2}}p(t,s)h_{M}^{e}(t)h_{M}^{e}(s) \,
\md \mu^{e}(t) \, \md \mu^{e}(s).
\end{equation*}
Furthermore,
\begin{align*}
I_{c}^{\alpha} \leqslant& \iint_{\mathbb{R}^{2}}p^{e}(t,s)h_{M}^{e}(t)h_{M}^{e}
(s) \, \md \mu^{e}(t) \, \md \mu^{e}(s) \! + \! \epsilon \underbrace{\int_{
\mathbb{R}}h_{M}^{e}(t) \, \md \mu^{e}(t)}_{\leqslant \, 1} \, \, \underbrace{
\int_{\mathbb{R}}h_{M}^{e}(s) \, \md \mu^{e}(s)}_{\leqslant \, 1} \Rightarrow
\\
I_{c}^{\alpha} \leqslant& \, \iint_{\mathbb{R}^{2}}p^{e}(t,s) \vert 1 \! + \!
(h_{M}^{e}(t) \! - \! 1) \vert \vert 1 \! + \! (h_{M}^{e}(s) \! - \! 1) \vert
\, \md \mu^{e}(t) \, \md \mu^{e}(s) \! + \! \epsilon \\
\leqslant& \, \iint_{\mathbb{R}^{2}}p^{e}(t,s) \, \md \mu^{e}(t) \, \md \mu^{e}
(s) \! + \! \iint_{\mathbb{R}^{2}}p^{e}(t,s) \vert h_{M}^{e}(s) \! - \! 1
\vert \, \md \mu^{e}(t) \, \md \mu^{e}(s) \! + \! \epsilon \\
+& \, \iint_{\mathbb{R}^{2}}p^{e}(t,s) \vert h_{M}^{e}(t) \! - \! 1 \vert \,
\md \mu^{e}(t) \, \md \mu^{e}(s) \! + \! \iint_{\mathbb{R}^{2}}p^{e}(t,s)
\vert h_{M}^{e}(t) \! - \! 1 \vert \vert h_{M}^{e}(s) \! - \! 1 \vert \, \md
\mu^{e}(t) \, \md \mu^{e}(s) \\
\leqslant& \, \iint_{\mathbb{R}^{2}}p^{e}(t,s) \, \md \mu^{e}(t) \, \md \mu^{e}
(s) \! + \! 2 \sup_{(t,s) \in \mathbb{R}^{2}} \vert p^{e}(t,s) \vert
\underbrace{\int_{(\mathbb{R} \setminus \mathfrak{D}_{M_{e}}) \cup \mathfrak{
D}_{M_{e}}} \vert h_{M}^{e}(s) \! - \! 1 \vert \, \md \mu^{e}(s)}_{\leqslant
\, \epsilon} \, \, \underbrace{\int_{\mathbb{R}} \md \mu^{e}(t)}_{= \, 1} \\
+& \, \sup_{(t,s) \in \mathbb{R}^{2}} \vert p^{e}(t,s) \vert \left(\underbrace{
\int_{(\mathbb{R} \setminus \mathfrak{D}_{M_{e}}) \cup \mathfrak{D}_{M_{e}}}
\vert h_{M}^{e}(t) \! - \! 1 \vert \, \md \mu^{e}(t)}_{\leqslant \, \epsilon}
\right)^{2} \! + \! \epsilon \\
\leqslant& \, \iint_{\mathbb{R}^{2}}p^{e}(t,s) \, \md \mu^{e}(t) \, \md \mu^{
e}(s) \! + \! \epsilon \! \left(1 \! + \! 2 \sup_{(t,s) \in \mathbb{R}^{2}}
\vert p^{e}(t,s) \vert \right) \! + \! \mathcal{O}(\epsilon^{2}),
\end{align*}
whereupon, neglecting the $\mathcal{O}(\epsilon^{2})$ term, and setting
$\varkappa^{\flat} \! := \! 1 \! + \! 2 \sup_{(t,s) \in \mathbb{R}^{2}} \vert
p^{e}(t,s) \vert$, one obtains
\begin{equation*}
I_{c}^{\alpha} \! \leqslant \! \iint_{\mathbb{R}^{2}}p^{e}(t,s) \, \md \mu^{e}
(t) \, \md \mu^{e}(s) \! + \! \varkappa^{\flat} \epsilon.
\end{equation*}
Hence, assembling the above-derived bounds for $I_{a}$, $I_{b}$, $I_{c}^{
\beta}$, and $I_{c}^{\alpha}$, one arrives at, upon setting $\epsilon^{\flat}
\! := \! 2 \varkappa^{\flat} \epsilon$,
\begin{equation*}
\iint_{\mathbb{R}^{2}}p^{e}(t,s) \, \md \mu_{n}^{e}(t) \, \md \mu^{e}_{n}(s)
\! - \! \iint_{\mathbb{R}^{2}}p^{e}(t,s) \, \md \mu^{e}(t) \, \md \mu^{e}(s)
\! \leqslant \! \epsilon^{\flat};
\end{equation*}
thus,
\begin{equation*}
\iint_{\mathbb{R}^{2}}p^{e}(t,s) \, \md \mu_{n}^{e}(t) \, \md \mu^{e}_{n}(s)
\! \to \! \iint_{\mathbb{R}^{2}}p^{e}(t,s) \, \md \mu^{e}(t) \, \md \mu^{e}
(s) \quad \text{as} \, \, \, n \! \to \! \infty.
\end{equation*}
Recalling that $p^{e}(t,s) \! := \! \min \left\lbrace q,K_{V_{m}}^{e}(t,s)
\right\rbrace$, $(q,m) \! \in \! \mathbb{R} \times \mathbb{N}$, it follows
{}from the above analysis that
\begin{equation*}
\liminf_{n \to \infty} \mathrm{I}_{V}^{e}[\mu_{n}^{e}] \! \geqslant \! \iint_{
\mathbb{R}^{2}} \min \left\lbrace q,K_{V_{m}}^{e}(t,s) \right\rbrace \, \md
\mu^{e}(t) \, \md \mu^{e}(s):
\end{equation*}
letting $q \! \uparrow \! \infty$ and $m \! \to \! \infty$, and using the
Monotone Convergence Theorem, one arrives at, upon noting that $\min
\left\lbrace q,K_{V_{m}}^{e}(t,s) \right\rbrace \! \to \! K_{V}^{e}(t,s)$,
\begin{equation*}
\liminf_{n \to \infty} \mathrm{I}_{V}^{e}[\mu_{n}^{e}] \! \geqslant \! \iint_{
\mathbb{R}^{2}}K_{V}^{e}(t,s) \, \md \mu^{e}(t) \, \md \mu^{e}(s) \! = \!
\mathrm{I}_{V}^{e}[\mu^{e}], \quad \mu_{n}^{e},\mu^{e} \! \in \! \mathcal{M}_{
1}(\mathbb{R}).
\end{equation*}
Since, {}from the analysis above, it was shown that there exists a weakly 
$(\mathrm{weak}^{\ast})$ convergent subsequence (of probability measures) 
$\{\mu_{n_{k}}^{e}\}_{k=1}^{\infty}$ $(\subset \mathcal{M}_{1}(\mathbb{R}))$ 
of $\{\mu_{n}^{e}\}_{n=1}^{\infty}$ $(\subset \mathcal{M}_{1}(\mathbb{R}))$ 
with a weak limit $\mu^{e} \! \in \! \mathcal{M}_{1}(\mathbb{R})$, namely, 
$\mu_{n_{k}}^{e} \! \to \! \mu^{e}$ as $k \! \to \! \infty$, upon recalling 
that $\mathrm{I}_{V}^{e}[\mu_{n}^{e}] \! \leqslant \! E_{V}^{e} \! + \! 
\tfrac{1}{n}$, $n \! \in \! \mathbb{N}$, it follows that, in the limit as $n 
\! \to \! \infty$, $\mathrm{I}_{V}^{e}[\mu^{e}] \! \leqslant \! E_{V}^{e} \! 
:= \! \inf \lbrace \mathstrut \mathrm{I}_{V}^{e}[\mu^{e}]; \, \mu^{e} \! \in 
\! \mathcal{M}_{1}(\mathbb{R}) \rbrace$; {}from the latter two inequalities, 
it follows, thus, that $\exists \, \, \mu^{e} \! := \! \mu_{V}^{e} \! \in \! 
\mathcal{M}_{1}(\mathbb{R})$, the `even' equilibrium measure, such that 
$\mathrm{I}_{V}^{e}[\mu_{V}^{e}] \! = \! \inf \lbrace \mathstrut \mathrm{I}_{
V}^{e}[\mu^{e}]; \, \mu^{e} \! \in \! \mathcal{M}_{1}(\mathbb{R}) \rbrace$, 
that is, the infimum is attained (the uniqueness of $\mu_{V}^{e} \! \in \! 
\mathcal{M}_{1}(\mathbb{R})$ is proven in Lemma~3.3 below).

The compactness of $\operatorname{supp}(\mu_{V}^{e}) \! =: \! J_{e}$ is now 
established: actually, the following proof is true for any $\mu \! \in \! 
\mathcal{M}_{1}(\mathbb{R})$ achieving the above minimum; in particular, for 
$\mu \! = \! \mu_{V}^{e}$. Without loss of generality, therefore, let $\mu_{w} 
\! \in \! \mathcal{M}_{1}(\mathbb{R})$ be such that $\mathrm{I}_{V}^{e}[\mu_{
w}] \! = \! E_{V}^{e}$, and let $D$ be any proper subset of $\mathbb{R}$ for 
which $\mu_{w}(D) \! := \! \int_{D} \md \mu_{w}(s) \! > \! 0$. As in 
\cite{a91}, set
\begin{equation*}
\mu_{w}^{\varepsilon}(z) \! := \! \left(1 \! + \! \varepsilon \mu_{w}(D) 
\right)^{-1} \! \left(\mu_{w}(z) \! + \! \varepsilon (\mu_{w} \! \! 
\upharpoonright_{D})(z) \right), \quad \varepsilon \! \in \! (-1,1),
\end{equation*}
where $\mu_{w} \! \! \upharpoonright_{D}$ denotes the restriction of $\mu_{w}$ 
to $D$ (note, also, that $\mu^{\varepsilon}_{w} \! > \! 0$ and bounded, and 
$\int_{\mathbb{R}} \md \mu_{w}^{\epsilon}(s) \! = \! 1)$. Using the fact that 
$K_{V}^{e}(s,t) \! = \! K_{V}^{e}(t,s)$, one shows that
\begin{align*}
\mathrm{I}_{V}^{e}[\mu_{w}^{\varepsilon}] =& \iint_{\mathbb{R}^{2}}K_{V}^{e}(s,
t) \, \md \mu^{\varepsilon}_{w}(s) \, \md \mu^{\varepsilon}_{w}(t) \\
=& \, (1 \! + \! \varepsilon \mu_{w}(D))^{-2} \iint_{\mathbb{R}^{2}}K_{V}^{e}
(s,t)(\md \mu_{w}(s) \! + \! \varepsilon \md (\mu_{w} \! \! \upharpoonright_{
D})(s))(\md \mu_{w}(t) \! + \! \varepsilon \md (\mu_{w} \! \! \upharpoonright_{
D})(t)) \\
=& \, (1 \! + \! \varepsilon \mu_{w}(D))^{-2} \! \left(\mathrm{I}_{V}^{e}
[\mu_{w}] \! + \! 2 \varepsilon \iint_{\mathbb{R}^{2}}K_{V}^{e}(s,t) \, \md
\mu_{w}(s) \, \md (\mu_{w} \! \! \upharpoonright_{D})(t) \right. \\
+& \left. \varepsilon^{2} \iint_{\mathbb{R}^{2}}K_{V}^{e}(s,t) \, \md (\mu_{w}
\! \! \upharpoonright_{D})(t) \, \md (\mu_{w} \! \! \upharpoonright_{D})(s)
\right).
\end{align*}
(Note that all of the above integrals are finite due to the argument at 
the beginning of the proof.) By the minimal property of $\mu_{w} \! \in \!
\mathcal{M}_{1}(\mathbb{R})$, it follows that $\partial_{\varepsilon} 
\mathrm{I}_{V}^{e}[\mu_{w}^{\varepsilon}] \! = \! 0$, that is,
\begin{equation*}
\iint_{\mathbb{R}^{2}}(K_{V}^{e}(s,t) \! - \! \mathrm{I}_{V}^{e}[\mu_{w}]) \,
\md \mu_{w}(s) \, \md (\mu_{w} \! \! \upharpoonright_{D})(t) \! = \! 0;
\end{equation*}
but, recalling that, with $\widehat{\psi}_{V}^{e}(z) \! := \! 2 \widetilde{V}
(z) \! - \! \ln (z^{2} \! + \! 1) \! - \! \ln (z^{-2} \! + \! 1)$, $K_{V}^{e}
(t,s) \! \geqslant \! \tfrac{1}{2} \widehat{\psi}_{V}^{e}(s) \! + \! \tfrac{
1}{2} \widehat{\psi}_{V}^{e}(t)$, it follows {}from the above that
\begin{gather*}
\iint_{\mathbb{R}^{2}} \mathrm{I}_{V}^{e}[\mu_{w}] \, \md \mu_{w}(s) \, \md
(\mu_{w} \! \! \upharpoonright_{D})(t) \! \geqslant \! \iint_{\mathbb{R}^{2}}
\! \left(\tfrac{1}{2} \widehat{\psi}_{V}^{e}(s) \! + \! \tfrac{1}{2} \widehat{
\psi}_{V}^{e}(t) \right) \! \md \mu_{w}(s) \, \md (\mu_{w} \! \!
\upharpoonright_{D})(t) \Rightarrow \\
0 \! \geqslant \! \iint_{\mathbb{R}^{2}} \! \left(\tfrac{1}{2} \widehat{
\psi}_{V}^{e}(s) \! + \! \tfrac{1}{2} \widehat{\psi}_{V}^{e}(t) \! - \!
\mathrm{I}_{V}^{e}[\mu_{w}] \right) \! \md \mu_{w}(s) \, \md (\mu_{w} \! \!
\upharpoonright_{D})(t),
\end{gather*}
whence
\begin{equation*}
\int_{\mathbb{R}} \! \left(\widehat{\psi}_{V}^{e}(t) \! + \! \left(\int_{
\mathbb{R}} \widehat{\psi}_{V}^{e}(s) \, \md \mu_{w}(s) \right) \! - \! 2
\mathrm{I}_{V}^{e}[\mu_{w}] \right) \! \md (\mu_{w} \! \! \upharpoonright_{
D})(t) \! \leqslant \! 0.
\end{equation*}
Recalling that
\begin{equation*}
\widehat{\psi}_{V}^{e}(x) \! := \! 2 \widetilde{V}(x) \! - \! \ln (x^{2} \! 
+ \! 1) \! - \! \ln (x^{-2} \! + \! 1) \! = \!
\begin{cases}
+\infty, &\text{$\vert x \vert \! \to \! \infty$,} \\
+\infty, &\text{$\vert x \vert \! \to \! 0$,}
\end{cases}
\end{equation*}
it follows that, $\exists \, \, T_{m} \! > \! 1$ such that
\begin{equation*}
\widehat{\psi}_{V}^{e}(t) \! + \! \int_{\mathbb{R}} \widehat{\psi}_{V}^{e}(s)
\, \md \mu_{w}(s) \! - \! 2 \mathrm{I}_{V}^{e}[\mu_{w}] \! \geqslant \! 1
\quad \text{for} \quad t \! \in \! \left((-T_{m},-T_{m}^{-1}) \cup (T_{m}^{-
1},T_{m}) \right)^{c}
\end{equation*}
(note, also, that $+\infty \! > \! \mathrm{I}_{V}^{e}[\mu_{w}] \! = \! \iint_{
\mathbb{R}^{2}}K_{V}^{e}(t,s) \, \md \mu_{w}(t) \, \md \mu_{w}(s) \! = \! 
\int_{\mathbb{R}} \widehat{\psi}_{V}^{e}(\xi) \, \md \mu_{w}(\xi) \! =$ a 
finite real number). Hence, if $D$ $(\subset \mathbb{R})$ is such that $D 
\subset (\{\vert x \vert \! \geqslant \! T_{m}\} \cup \{\vert x \vert \! 
\leqslant \! T_{m}^{-1}\})$, $T_{m} \! > \! 1$, it follows {}from the above 
calculations that
\begin{equation*}
0 \! \geqslant \! \int_{\mathbb{R}} \! \left(\widehat{\psi}_{V}^{e}(t) \! + \!
\left(\int_{\mathbb{R}} \widehat{\psi}_{V}^{e}(s) \, \md \mu_{w}(s) \right)
\! - \! 2 \mathrm{I}_{V}^{e}[\mu_{w}] \right) \! \md (\mu_{w} \! \!
\upharpoonright_{D})(t) \!
\geqslant \! 1,
\end{equation*}
which is a contradiction; hence, $\operatorname{supp}(\mu_{w}) \subseteq [-T_{
m},-T_{m}^{-1}] \cup [T_{m}^{-1},T_{m}]$, $T_{m} \! > \! 1$; in particular, 
$J_{e} \! := \! \operatorname{supp}(\mu_{V}^{e}) \subseteq [-T_{m},-T_{m}^{-
1}] \cup [T_{m}^{-1},T_{m}]$, $T_{m} \! > \! 1$, which establishes the 
compactness of the support of the `even' equilibrium measure $\mu_{V}^{e} \! 
\in \! \mathcal{M}_{1}(\mathbb{R})$. Furthermore, it is worth noting that, 
since $J_{e} \! := \! \operatorname{supp}(\mu_{V}^{e}) \! =$ compact 
$(\subseteq \overline{\mathbb{R}} \setminus \{0,\pm \infty\})$, and 
$\widetilde{V} \colon \mathbb{R} \setminus \{0\} \! \to \! \mathbb{R}$ is 
real analytic on $J_{e}$,
\begin{align*}
+\infty \! > \! E_{V}^{e} \, (=& \, \mathrm{I}_{V}^{e}[\mu_{V}^{e}]) \,
\geqslant \, \iint_{\mathbb{R}^{2}} \ln \! \left(\vert s \! - \! t \vert^{2}
\vert st \vert^{-1}w^{e}(s)w^{e}(t) \right)^{-1} \md \mu_{V}^{e}(s) \, \md
\mu_{V}^{e}(t) \\
=& \, \iint_{J_{e}^{2}} \ln \! \left(\vert s \! - \! t \vert^{2} \vert st
\vert^{-1}w^{e}(s)w^{e}(t) \right)^{-1} \md \mu_{V}^{e}(s) \, \md \mu_{V}^{e}
(t) \! > \! -\infty;
\end{align*}
moreover, a straightforward consequence of the fact just established is that 
$J_{e}$ has positive logarithmic capacity, that is, $\operatorname{cap}(J_{
e}) \! = \! \exp (-E_{V}^{e}) \! > \! 0$. \hfill $\qed$
\begin{eeeee}
It is important to note {}from the latter part of the proof of Lemma~3.1 that 
$J_{e} \! \not\supseteq \! \{0,\pm \infty\}$. This can also be seen as 
follows. For $\varepsilon$ some arbitrarily fixed, sufficiently small positive 
real number and $\Sigma_{\varepsilon} \! := \! \{\mathstrut z; \, w^{e}(z) 
\! \geqslant \! \varepsilon\}$, if $(s,t) \! \notin \! \Sigma_{\varepsilon} 
\times \Sigma_{\varepsilon}$, then $\ln (\vert s \! - \! t \vert^{2} \vert 
st \vert^{-1}w^{e}(s)w^{e}(t))^{-1} \! =: \! K_{V}^{e}(s,t)$ $(= \! K_{V}^{e}
(t,s))$ $> \! E_{V}^{e} \! + \! 1$, which is a contradiction, since it was 
established above that the minimum is attained $\Leftrightarrow \! (s,t) \! 
\in \! \Sigma_{\varepsilon} \times \Sigma_{\varepsilon}$. Towards this end, 
it is enough to show that (see, for example, \cite{a55}), if $\{(s_{n},t_{n})
\}_{n=1}^{\infty}$ is a sequence with $\lim \min_{n \to \infty}\{w^{e}(s_{n}),
w^{e}(t_{n})\} \! = \! 0$, then $\lim_{n \to \infty} \ln (\vert s_{n} \! - \! 
t_{n} \vert^{2} \vert s_{n}t_{n} \vert^{-1}w^{e}(s_{n})w^{e}(t_{n}))^{-1} \! 
= \! \lim_{n \to \infty}K_{V}^{e}(s_{n},t_{n}) \! = \! +\infty$. Without loss 
of generality, one can assume that $s_{n} \! \to \! s$ and $t_{n} \! \to \! 
t$ as $n \! \to \! \infty$, where $s$, $t$, or both may be infinite; thus, 
there are several cases to consider:
\begin{compactenum}
\item[(i)] if $s$ and $t$ are finite, then, {}from $\lim \min_{n \to \infty}\{
w^{e}(s_{n}),w^{e}(t_{n})\} \! = \! \min \{w^{e}(s),w^{e}(t)\} \! = \! 0$, it
is clear that $\lim_{n \to \infty}K_{V}^{e}(s_{n},t_{n}) \! = \! +\infty$;
\item[(ii)] if $\vert s \vert \! = \! \infty$ (resp., $\vert t \vert \! = \!
\infty)$ but $t \! = \! \text{finite}$ (resp., $s \! = \! \text{finite})$,
then, due to the fact that $\widetilde{V} \colon \mathbb{R} \setminus \{0\} \!
\to \! \mathbb{R}$ satisfies the conditions
\begin{equation*}
2 \widetilde{V}(x) \! - \! \ln (x^{2} \! + \! 1) \! - \! \ln (x^{-2} \! + \!
1) \! = \!
\begin{cases}
+\infty, &\text{$\vert x \vert \! \to \! \infty$,} \\
+\infty, &\text{$\vert x \vert \! \to \! 0$,}
\end{cases}
\end{equation*}
it follows that $\lim_{n \to \infty}K_{V}^{e}(s_{n},t_{n}) \! = \! +\infty$;
\item[(iii)] if $\vert s \vert \! = \! 0$ (resp., $\vert t \vert \! = \! 0)$
but $t \! = \! \text{finite}$ (resp., $s \! = \! \text{finite})$, then, as a
result of the above conditions for $\widetilde{V}$, it follows that $\lim_{n
\to \infty}K_{V}^{e}(s_{n},t_{n}) \! = \! +\infty$;
\item[(iv)] if $\vert s \vert \! = \! \infty$ and $\vert t \vert \! = \!
\infty$, then, again due to the above conditions for $\widetilde{V}$, it
follows that $\lim_{n \to \infty}K_{V}^{e} \linebreak[4]
(s_{n},t_{n}) \! = \! +\infty$; and
\item[(v)] if $\vert s \vert \! = \! 0$ and $\vert t \vert \! = \! 0$, then,
again, as above, it follows that $\lim_{n \to \infty}K_{V}^{e}(s_{n},t_{n}) 
\! = \! +\infty$.
\end{compactenum}
Hence, $K_{V}^{e}(s,t) \! > \! E_{V}^{e} \! + \! 1$ if $(s,t) \! \notin \!
\Sigma_{\varepsilon} \times \Sigma_{\varepsilon}$, that is, if $s$, $t$, or
both $\in \! \{0,\pm \infty\}$ (which can not be the case, as the infimum $E_{
V}^{e}$ is attained $\Leftrightarrow \! (s,t) \! \in \! \Sigma_{\varepsilon}
\times \Sigma_{\varepsilon}$, whence $\operatorname{supp}(\mu_{V}^{e}) \! =:
\! J_{e} \! \not\supseteq \! \{0,\pm \infty\})$. \hfill $\blacksquare$
\end{eeeee}
In order to establish the uniqueness of the `even' equilibrium measure, 
$\mu_{V}^{e}$ $(\in \! \mathcal{M}_{1}(\mathbb{R}))$, the following lemma 
is requisite.
\begin{ccccc}
Let $\mu \! := \! \mu_{1} \! - \! \mu_{2}$, where $\mu_{1},\mu_{2}$ are 
non-negative, finite-moment $(\int_{\operatorname{supp}(\mu_{j})}s^{m} \, 
\md \mu_{j}(s) \! < \! \infty$, $m \! \in \! \mathbb{Z}$, $j \! = \! 1,2)$ 
measures on $\mathbb{R}$ supported on distinct sets $(\operatorname{supp}
(\mu_{1}) \cap \operatorname{supp}(\mu_{2}) \! = \! \varnothing)$, be the 
(unique) Jordan decomposition of the finite-moment signed measure on $\mathbb{
R}$ with mean zero, that is, $\int_{\operatorname{supp}(\mu)} \md \mu (s) \! 
= \! 0$, and with $\operatorname{supp}(\mu) \! = \! \text{compact}$. Suppose 
that $-\infty \! < \! \iint_{\mathbb{R}^{2}} \ln (\vert s \! - \! t \vert^{-
2} \vert st \vert) \, \md \mu_{j}(s) \, \md \mu_{j}(t) \! < \! +\infty$, $j 
\! = \! 1,2$. Then,
\begin{equation*}
\iint_{\mathbb{R}^{2}} \ln \! \left(\dfrac{\vert st \vert}{\vert s \! - \! t
\vert^{2}} \right) \! \md \mu (s) \, \md \mu (t) \! = \! \iint_{\mathbb{R}^{2}}
\ln \! \left(\dfrac{\vert s \! - \! t \vert^{2}}{\vert st \vert}w^{e}(s)w^{e}
(t) \right)^{-1} \md \mu (s) \, \md \mu (t) \! \geqslant \! 0,
\end{equation*}
where equality holds if, and only if, $\mu \! = \! 0$.
\end{ccccc}

\emph{Proof.} Recall the following identity \cite{a90} (see pg.~147, 
Equation~(6.44)): for $\xi \! \in \! \mathbb{R}$ and any $\varepsilon \! > 
\! 0$,
\begin{equation*}
\ln (\xi^{2} \! + \! \varepsilon^{2}) \! = \! \ln (\varepsilon^{2}) \! + \! 
2 \, \Im \! \left(\int_{0}^{+\infty} \! \left(\dfrac{\me^{\mi \xi v} \! - \! 
1}{\mi v} \right) \! \me^{-\varepsilon v} \, \md v \right);
\end{equation*}
thus, it follows that
\begin{align*}
\iint_{\mathbb{R}^{2}} \ln ((s \! - \! t)^{2} \! + \! \varepsilon^{2}) \, \md
\mu (s) \, \md \mu (t) =& \, \iint_{\mathbb{R}^{2}} \ln (\varepsilon^{2}) \,
\md \mu (s) \, \md \mu (t) \\
+& \, \iint_{\mathbb{R}^{2}} \! \left(2 \, \Im \! \left(\int_{0}^{+\infty} \!
\left(\dfrac{\me^{\mi (s-t)v} \! - \! 1}{\mi v} \right) \! \me^{-\varepsilon
v} \, \md v \right) \! \right) \! \md \mu (s) \, \md \mu (t), \\
\iint_{\mathbb{R}^{2}} \ln (s^{2} \! + \! \varepsilon^{2}) \, \md \mu (s) \,
\md \mu (t) =& \, \iint_{\mathbb{R}^{2}} \ln (\varepsilon^{2}) \, \md \mu (s)
\, \md \mu (t) \\
+& \, \iint_{\mathbb{R}^{2}} \! \left(2 \, \Im \! \left(\int_{0}^{+\infty} \!
\left(\dfrac{\me^{\mi sv} \! - \! 1}{\mi v} \right) \! \me^{-\varepsilon v} \,
\md v \right) \! \right) \! \md \mu (s) \, \md \mu (t), \\
\iint_{\mathbb{R}^{2}} \ln (t^{2} \! + \! \varepsilon^{2}) \, \md \mu (s) \,
\md \mu (t) =& \, \iint_{\mathbb{R}^{2}} \ln (\varepsilon^{2}) \, \md \mu (s)
\, \md \mu (t) \\
+& \, \iint_{\mathbb{R}^{2}} \! \left(2 \, \Im \! \left(\int_{0}^{+\infty} \!
\left(\dfrac{\me^{\mi tv} \! - \! 1}{\mi v} \right) \! \me^{-\varepsilon v} \,
\md v \right) \! \right) \! \md \mu (s) \, \md \mu (t);
\end{align*}
but, since $\iint_{\mathbb{R}^{2}} \md \mu (s) \, \md \mu (t) \! = \! \left(
\int_{\mathbb{R}} \md \mu (s) \right)^{2} \! = \! 0$, one obtains, after some
rearrangement,
\begin{align*}
\iint_{\mathbb{R}^{2}} \ln ((s \! - \! t)^{2} \! + \! \varepsilon^{2}) \,
\md \mu (s) \, \md \mu (t) =& 2 \, \Im \! \left(\int_{0}^{+\infty} \me^{-
\varepsilon v} \! \left(\iint_{\mathbb{R}^{2}} \! \left(\dfrac{\me^{\mi
(s-t)v} \! - \! 1}{\mi v} \right) \! \md \mu (s) \, \md \mu (t) \right) \!
\md v \right), \\
\iint_{\mathbb{R}^{2}} \ln (s^{2} \! + \! \varepsilon^{2}) \, \md \mu (s) \,
\md \mu (t) =& 2 \, \Im \left(\int_{0}^{+\infty} \me^{-\varepsilon v} \! 
\left(\iint_{\mathbb{R}^{2}} \! \left(\dfrac{\me^{\mi sv} \! - \! 1}{\mi v} 
\right) \! \md \mu (s) \, \md \mu (t) \right) \! \md v \right), \\
\iint_{\mathbb{R}^{2}} \ln (t^{2} \! + \! \varepsilon^{2}) \, \md \mu (s) \,
\md \mu (t) =& 2 \, \Im \left(\int_{0}^{+\infty} \me^{-\varepsilon v} \! 
\left(\iint_{\mathbb{R}^{2}} \! \left(\dfrac{\me^{\mi tv} \! - \! 1}{\mi v} 
\right) \! \md \mu (s) \, \md \mu (t) \right) \! \md v \right).
\end{align*}
Noting that
\begin{align*}
\iint_{\mathbb{R}^{2}} \! \left(\dfrac{\me^{\mi (s-t)v} \! - \! 1}{\mi v}
\right) \! \md \mu (s) \, \md \mu (t) =& \, \dfrac{1}{\mi v} \iint_{\mathbb{
R}^{2}} \me^{\mi (s-t)v} \, \md \mu (s) \, \md \mu (t) \! - \! \dfrac{1}{\mi
v} \underbrace{\iint_{\mathbb{R}^{2}} \md \mu (s) \, \md \mu (t)}_{= \, 0} \\
=& \, \dfrac{1}{\mi v} \int_{\mathbb{R}} \me^{\mi sv} \, \md \mu (s) \int_{
\mathbb{R}} \me^{-\mi tv} \, \md \mu (t),
\end{align*}
and setting $\widehat{\mu}(z) \! := \! \int_{\mathbb{R}} \me^{\mi \xi z} \,
\md \mu (\xi)$, one gets that
\begin{equation*}
\iint_{\mathbb{R}^{2}} \! \left(\dfrac{\me^{\mi (s-t)v} \! - \! 1}{\mi v}
\right) \! \md \mu (s) \, \md \mu (t) \! = \! \dfrac{1}{\mi v} \vert \widehat{
\mu}(v) \vert^{2}:
\end{equation*}
also,
\begin{equation*}
\iint_{\mathbb{R}^{2}} \! \left(\dfrac{\me^{\mi sv} \! - \! 1}{\mi v} \right)
\! \md \mu (s) \, \md \mu (t) \! = \! \dfrac{1}{\mi v} \int_{\mathbb{R}} \me^{
\mi sv} \, \md \mu (s) \underbrace{\int_{\mathbb{R}} \md \mu (t)}_{= \, 0}-
\dfrac{1}{\mi v} \underbrace{\int_{\mathbb{R}} \md \mu (s)}_{= \, 0} \, \,
\underbrace{\int_{\mathbb{R}} \md \mu (t)}_{= \, 0} \! = \! 0;
\end{equation*}
similarly,
\begin{equation*}
\iint_{\mathbb{R}^{2}} \! \left(\dfrac{\me^{\mi tv} \! - \! 1}{\mi v} \right)
\! \md \mu (s) \, \md \mu (t) \! = \! 0.
\end{equation*}
Hence,
\begin{gather*}
\iint_{\mathbb{R}^{2}} \ln ((s \! - \! t)^{2} \! + \! \varepsilon^{2}) \, \md
\mu (s) \, \md \mu (t) \! = \! 2 \, \Im \! \left(\int_{0}^{+\infty} \dfrac{
\vert \widehat{\mu}(v) \vert^{2}}{\mi v} \me^{-\varepsilon v} \, \md v
\right), \\
\iint_{\mathbb{R}^{2}} \ln (s^{2} \! + \! \varepsilon^{2}) \, \md \mu (s) \,
\md \mu (t) \! = \! \iint_{\mathbb{R}^{2}} \ln (t^{2} \! + \! \varepsilon^{2})
\, \md \mu (s) \, \md \mu (t) \! = \! 0.
\end{gather*}
Noting that $\widehat{\mu}(0) \! = \! \int_{\mathbb{R}} \md \mu (\xi) \! =
\! 0$, a Taylor expansion about $v \! = \! 0$ shows that $\widehat{\mu}(v) \!
=_{v \to 0} \! \widehat{\mu}^{\prime}(0)v \! + \! \mathcal{O}(v^{2})$, where
$\widehat{\mu}^{\prime}(0) \! := \! \partial_{v} \widehat{\mu}(v) \vert_{v=
0}$; thus, $v^{-1} \vert \widehat{\mu}(v) \vert^{2} \! =_{v \to 0} \! \vert
\widehat{\mu}^{\prime}(0) \vert^{2}v \! + \! \mathcal{O}(v^{2})$, which means
that there is no singularity in the integrand as $v \! \to \! 0$ (in fact,
$v^{-1} \vert \widehat{\mu}(v) \vert^{2}$ is real analytic in a neighbourhood
of the origin), whence
\begin{equation*}
\iint_{\mathbb{R}^{2}} \ln ((s \! - \! t)^{2} \! + \! \varepsilon^{2}) \, \md
\mu (s) \, \md \mu (t) \! = \! -2 \int_{0}^{+\infty}v^{-1} \vert \widehat{\mu}
(v) \vert^{2} \me^{-\varepsilon v} \, \md v.
\end{equation*}
Recalling that $\iint_{\mathbb{R}^{2}} \ln (\ast^{2} \! + \! \varepsilon^{2})
\, \md \mu (s) \, \md \mu (t) \! = \! 0$, $\ast \! \in \! \{s,t\}$, and
adding, it follows that
\begin{equation*}
\iint_{\mathbb{R}^{2}} \ln \! \left(\dfrac{(s^{2} \! + \! \varepsilon^{2})^{
1/2}(t^{2} \! + \! \varepsilon^{2})^{1/2}}{((s \! - \! t)^{2} \! + \!
\varepsilon^{2})} \right) \! \md \mu (s) \, \md \mu (t) \! = \! 2 \int_{0}^{+
\infty}v^{-1} \vert \widehat{\mu}(v) \vert^{2} \me^{-\varepsilon v} \, \md v.
\end{equation*}
Now, using the fact that $\ln ((s \! - \! t)^{2} \! + \! \varepsilon^{2})^{-
1}$ (resp., $\ln (s^{2} \! + \! \varepsilon^{2})^{1/2}$ and $\ln (t^{2} \!
+ \! \varepsilon^{2})^{1/2})$ is (resp., are) bounded below (resp., above)
uniformly with respect to $\varepsilon$ and that the measures have compact
support, letting $\varepsilon \! \downarrow \! 0$ and using the Monotone
Convergence Theorem, one arrives at
\begin{align*}
\iint_{\mathbb{R}^{2}} \ln \! \left(\dfrac{(s^{2} \! + \! \varepsilon^{2})^{
1/2}(t^{2} \! + \! \varepsilon^{2})^{1/2}}{((s \! - \! t)^{2} \! + \!
\varepsilon^{2})} \right) \! \md \mu (s) \, \md \mu (t) \underset{\varepsilon
\downarrow 0}{=}& \, \iint_{\mathbb{R}^{2}} \ln \! \left(\dfrac{\vert st
\vert}{\vert s \! - \! t \vert^{2}} \right) \! \md \mu (s) \, \md \mu (t) \\
=& \, 2 \int_{0}^{+\infty}v^{-1} \vert \widehat{\mu}(v) \vert^{2} \, \md v \!
\geqslant \! 0,
\end{align*}
where, trivially, equality holds if, and only if, $\mu \! = \! 0$.
Furthermore, noting that, since $\int_{\mathbb{R}} \md \mu (\xi) \! = \! 0$,
$\iint_{\mathbb{R}^{2}} \ln (w^{e}(\ast))^{-1} \, \md \mu (s) \, \md \mu (t)
\! = \! 0$, $\ast \! \in \! \{s,t\}$, letting $\varepsilon \! \downarrow \!
0$ and using monotone convergence, one also arrives at
\begin{align*}
\iint_{\mathbb{R}^{2}} \ln \! \left(\dfrac{(s^{2} \! + \! \varepsilon^{2})^{
1/2}(t^{2} \! + \! \varepsilon^{2})^{1/2}}{((s \! - \! t)^{2} \! + \!
\varepsilon^{2})w^{e}(s)w^{e}(t)} \right) \! \md \mu (s) \, \md \mu (t)
\underset{\varepsilon \downarrow 0}{=}& \, \iint_{\mathbb{R}^{2}} \ln \! \left(
\dfrac{\vert s \! - \! t \vert^{2}}{\vert st \vert}w^{e}(s)w^{e}(t) \right)^{-
1} \md \mu (s) \, \md \mu (t) \\
=& \, 2 \int_{0}^{+\infty}v^{-1} \vert \widehat{\mu}(v) \vert^{2} \, \md v \!
\geqslant \! 0,
\end{align*}
where, again, and trivially, equality holds if, and only if, $\mu \! = \! 0$.
\hfill $\qed$

The uniqueness of $\mu_{V}^{e}$ $(\in \! \mathcal{M}_{1}(\mathbb{R}))$ will 
now be established.
\begin{ccccc}
Let the external field $\widetilde{V} \colon \mathbb{R} \setminus \{0\} \! \to
\! \mathbb{R}$ satisfy conditions~{\rm (2.3)--(2.5)}. Set $w^{e}(z) \! := \!
\exp (-\widetilde{V}(z))$, and define
\begin{equation*}
\mathrm{I}_{V}^{e}[\mu^{e}] \colon \mathcal{M}_{1}(\mathbb{R}) \! \to \!
\mathbb{R},\, \, \mu^{e} \! \mapsto \! \iint_{\mathbb{R}^{2}} \ln \! \left(
\vert s \! - \! t \vert^{2} \vert st \vert^{-1}w^{e}(s)w^{e}(t) \right)^{-1}
\md \mu^{e}(s) \, \md \mu^{e}(t),
\end{equation*}
and consider the minimisation problem $E_{V}^{e} \! = \! \inf \lbrace
\mathstrut \mathrm{I}_{V}^{e}[\mu^{e}]; \, \mu^{e} \! \in \! \mathcal{M}_{1}
(\mathbb{R}) \rbrace$. Then, $\exists ! \, \, \mu_{V}^{e} \in \! \mathcal{M}_{
1}(\mathbb{R})$ such that $\mathrm{I}_{V}^{e}[\mu_{V}^{e}] \! = \! E_{V}^{e}$.
\end{ccccc}

\emph{Proof.} It was shown in Lemma~3.1 that $\exists \, \, \mu_{V}^{e} \! \in
\! \mathcal{M}_{1}(\mathbb{R})$, the `even' equilibrium measure, such that
$\mathrm{I}_{V}^{e}[\mu^{e}] \! = \! E_{V}^{e}$; therefore, it remains to
establish the uniqueness of the `even' equilibrium measure. Let $\widetilde{
\mu}_{V}^{e} \! \in \! \mathcal{M}_{1}(\mathbb{R})$ be a second probability
measure for which $\mathrm{I}_{V}^{e}[\widetilde{\mu}_{V}^{e}] \! = \! E_{V}^{
e} \! = \! \mathrm{I}_{V}^{e}[\mu_{V}^{e}]$: the argument in Lemma~3.1 shows
that $\widetilde{J}_{e} \! := \! \operatorname{supp}(\widetilde{\mu}_{V}^{e})
\! = \! \text{compact} \subseteq \! \overline{\mathbb{R}} \setminus \{0,\pm
\infty\}$, and that $\mathrm{I}_{V}^{e}[\widetilde{\mu}_{V}^{e}] \! < \!
+\infty$. Define the fi\-ni\-te-mo\-me\-nt signed measure $\mu^{\sharp} \! :=
\! \widetilde{\mu}_{V}^{e} \! - \! \mu_{V}^{e}$, where $\widetilde{\mu}_{V}^{
e},\mu_{V}^{e} \! \in \! \mathcal{M}_{1}(\mathbb{R})$, and $\widetilde{J}_{e}
\cap J_{e} \! = \! \varnothing$, with (cf. Lemma~3.1), $J_{e} \! = \! 
\operatorname{supp}(\mu_{V}^{e}) \! = \! \text{compact} \subseteq \overline{
\mathbb{R}} \setminus \{0,\pm \infty\}$; thus, {}from Lemma~3.2 (with $\mu 
\! \to \! \mu^{\sharp})$, namely,
\begin{equation*}
\iint_{\mathbb{R}^{2}} \ln \! \left(\vert s \! - \! t \vert^{-2} \vert st
\vert \right) \md \mu^{\sharp}(s) \, \md \mu^{\sharp}(t) \! = \! \iint_{
\mathbb{R}^{2}} \ln \! \left(\vert s \! - \! t \vert^{2} \vert st \vert^{-1}
w^{e}(s)w^{e}(t) \right)^{-1} \md \mu^{\sharp}(s) \, \md \mu^{\sharp}(t) \!
\geqslant \! 0,
\end{equation*}
it follows that
\begin{align*}
\iint_{\mathbb{R}^{2}} \ln \! \left(\vert st \vert \vert s \! - \! t \vert^{-
2} \right) \! \left(\md \widetilde{\mu}_{V}^{e}(s) \, \md \widetilde{\mu}_{
V}^{e}(t) \! + \! \md \mu_{V}^{e}(s) \, \md \mu_{V}^{e}(t) \right) \geqslant&
\, \iint_{\mathbb{R}^{2}} \ln \! \left(\vert st \vert \vert s \! - \! t \vert^{
-2} \right) \! \left(\md \widetilde{\mu}_{V}^{e}(s) \, \md \mu_{V}^{e}(t)
\right. \\
+& \left. \, \md \mu_{V}^{e}(s) \md \widetilde{\mu}_{V}^{e}(t) \right),
\end{align*}
or, via a straightforward symmetry argument,
\begin{align*}
\iint_{\mathbb{R}^{2}} \ln \! \left(\vert st \vert \vert s \! - \! t \vert^{-
2} \right) \! \left(\md \widetilde{\mu}_{V}^{e}(s) \, \md \widetilde{\mu}_{
V}^{e}(t) \! + \! \md \mu_{V}^{e}(s) \, \md \mu_{V}^{e}(t) \right) \geqslant&
\, 2 \iint_{\mathbb{R}^{2}} \ln \! \left(\vert st \vert \vert s \! - \! t
\vert^{-2} \right) \! \md \widetilde{\mu}_{V}^{e}(s) \, \md \mu_{V}^{e}(t) \\
=& \, 2 \iint_{\mathbb{R}^{2}} \ln \! \left(\vert st \vert \vert s \! - \! t
\vert^{-2} \right) \! \md \mu_{V}^{e}(s) \, \md \widetilde{\mu}_{V}^{e}(t).
\end{align*}
The above shows that (since both $\mathrm{I}_{V}^{e}[\mu_{V}^{e}]$ and 
$\mathrm{I}_{V}^{e}[\widetilde{\mu}_{V}^{e}] \! < \! +\infty)$ $\ln (\vert 
st \vert \vert s \! - \! t \vert^{-2})$ is integrable with respect to both 
$\md \widetilde{\mu}_{V}^{e}(s) \, \md \mu_{V}^{e}(t)$ and $\md \mu_{V}^{e}
(s) \, \md \widetilde{\mu}_{V}^{e}(t)$. {}From an argument on pg.~149 of 
\cite{a90}, it follows that $\ln (\vert st \vert \vert s \! - \! t \vert^{-
2})$ is integrable with respect to (the measure) $\md \mu_{t}^{e}(s) \, \md 
\mu_{t}^{e}(t^{\prime})$, where $\mu_{t}^{e}(z) \! := \! \mu_{V}^{e}(z) \! + 
\! t(\widetilde{\mu}_{V}^{e}(z) \! - \! \mu_{V}^{e}(z))$, $(z,t) \! \in \! 
\mathbb{R} \times [0,1]$. Set
\begin{equation*}
\mathscr{F}_{\mu}(t) \! := \! \iint_{\mathbb{R}^{2}} \ln \! \left(\vert st^{
\prime} \vert \vert s \! - \! t^{\prime} \vert^{-2}(w^{e}(s)w^{e}(t^{\prime})
)^{-1} \right) \md \mu_{t}^{e}(s) \, \md \mu_{t}^{e}(t^{\prime})
\end{equation*}
$(= \! \mathrm{I}_{V}^{e}[\mu_{t}^{e}])$. Noting that
\begin{align*}
\md \mu_{t}^{e}(s) \, \md \mu_{t}^{e}(t^{\prime})=& \, \md \mu_{V}^{e}(s) \,
\md \mu_{V}^{e}(t^{\prime}) \! + \! t \md \mu_{V}^{e}(s)(\md \widetilde{\mu}_{
V}^{e}(t^{\prime}) \! - \! \md \mu_{V}^{e}(t^{\prime})) \! + \! t \md \mu_{
V}^{e}(t^{\prime})(\md \widetilde{\mu}_{V}^{e}(s) \! - \! \md \mu_{V}^{e}(s))
\\
+& \, t^{2}(\md \widetilde{\mu}_{V}^{e}(s) \! - \! \md \mu_{V}^{e}(s))(\md
\widetilde{\mu}_{V}^{e}(t^{\prime}) \! - \! \md \mu_{V}^{e}(t^{\prime})),
\end{align*}
it follows that
\begin{align*}
\mathscr{F}_{\mu}(t) \! =& \, \mathrm{I}_{V}^{e}[\mu_{V}^{e}] \! + \! 2t
\iint_{\mathbb{R}^{2}} \ln \! \left(\dfrac{\vert st^{\prime} \vert}{\vert s
\! - \! t^{\prime} \vert^{2}}(w^{e}(s)w^{e}(t^{\prime}))^{-1} \right) \! \md
\mu_{V}^{e}(s)(\md \widetilde{\mu}_{V}^{e}(t^{\prime}) \! - \! \md \mu_{V}^{e}
(t^{\prime})) \\
+& \, t^{2} \iint_{\mathbb{R}^{2}} \ln \! \left(\dfrac{\vert st^{\prime}
\vert}{\vert s \! - \! t^{\prime} \vert^{2}}(w^{e}(s)w^{e}(t^{\prime}))^{-1}
\right) \! (\md \widetilde{\mu}_{V}^{e}(s) \! - \! \md \mu_{V}^{e}(s))(\md
\widetilde{\mu}_{V}^{e}(t^{\prime}) \! - \! \md \mu_{V}^{e}(t^{\prime})).
\end{align*}
Since $\mu^{\sharp} \! \in \! \mathcal{M}_{1}(\mathbb{R})$ is a finite-moment
signed measure with mean zero, that is, $\int_{\mathbb{R}} \md \mu^{\sharp}
(\xi) \! = \! \int_{\mathbb{R}} \md (\widetilde{\mu}_{V}^{e} \! - \! \mu_{V}^{
e})(\xi) \! = \! 0$, and compact support, it follows {}from the analysis
above and the result of Lemma~3.2 that $\mathscr{F}_{\mu}(t)$ is
convex\footnote{If $f$ is twice differentiable on $(a,b)$, then $f^{\prime
\prime}(x) \! \geqslant \! 0$ on $(a,b)$ is both a necessary and sufficient
condition that $f$ be convex on $(a,b)$.}; thus, for $t \! \in \! [0,1]$,
\begin{align*}
\mathrm{I}_{V}^{e}[\mu_{V}^{e}] \leqslant& \, \mathscr{F}_{\mu}(t) \! = \!
\mathrm{I}_{V}^{e}[\mu_{t}^{e}] \! = \! \mathscr{F}_{\mu}(t \! + \! (1 \! - \!
t)0) \! \leqslant \! t \mathscr{F}_{\mu}(1) \! + \! (1 \! - \! t) \mathscr{
F}_{\mu}(0) \\
=& \, t \mathrm{I}_{V}^{e}[\widetilde{\mu}_{V}^{e}] \! + \! (1 \! - \! t)
\mathrm{I}_{V}^{e}[\mu_{V}^{e}] \! = \! t \mathrm{I}_{V}^{e}[\mu_{V}^{e}] \! +
\! (1 \! - \! t) \mathrm{I}_{V}^{e}[\mu_{V}^{e}] \Rightarrow \\
\mathrm{I}_{V}^{e}[\mu_{V}^{e}] \leqslant& \, \mathrm{I}_{V}^{e}[\mu_{t}^{e}]
\! \leqslant \! \mathrm{I}_{V}^{e}[\mu_{V}^{e}],
\end{align*}
whence $\mathrm{I}_{V}^{e}[\mu_{t}^{e}] \! = \! \mathrm{I}_{V}^{e}[\mu_{V}^{e}
] \! := \! E_{V}^{e}$ $(= \! \text{const.})$. Since $\mathrm{I}_{V}^{e}
[\mu_{t}^{e}] \! = \! \mathscr{F}_{\mu}(t) \! = \! E_{V}^{e}$, it follows, in
particular, that $\mathscr{F}_{\mu}^{\prime \prime}(0) \! = \! 0 \Rightarrow$
\begin{align*}
0 =& \, \iint_{\mathbb{R}^{2}} \ln \! \left(\dfrac{\vert st^{\prime} \vert}{
\vert s \! - \! t^{\prime} \vert^{2}}(w^{e}(s)w^{e}(t^{\prime}))^{-1} \right)
\! (\md \widetilde{\mu}_{V}^{e}(s) \! - \! \md \mu_{V}^{e}(s))(\md \widetilde{
\mu}_{V}^{e}(t^{\prime}) \! - \! \md \mu_{V}^{e}(t^{\prime})) \\
=& \, \iint_{\mathbb{R}^{2}} \ln \! \left(\dfrac{\vert st^{\prime} \vert}{
\vert s \! - \! t^{\prime} \vert^{2}} \right) \! (\md \widetilde{\mu}_{V}^{e}
(s) \! - \! \md \mu_{V}^{e}(s))(\md \widetilde{\mu}_{V}^{e}(t^{\prime}) \! -
\! \md \mu_{V}^{e}(t^{\prime})) \\
+& \, 2 \int_{\mathbb{R}} \widetilde{V}(t^{\prime}) \, \md (\widetilde{\mu}_{
V}^{e} \! - \! \mu_{V}^{e})(t^{\prime}) \, \, \underbrace{\int_{\mathbb{R}}
\md (\widetilde{\mu}_{V}^{e} \! - \! \mu_{V}^{e})(s)}_{= \, 0} \Rightarrow \\
0 =& \, \iint_{\mathbb{R}^{2}} \ln \! \left(\dfrac{\vert st^{\prime} \vert}{
\vert s \! - \! t^{\prime} \vert^{2}} \right) \! \md (\widetilde{\mu}_{V}^{e}
\! - \! \mu_{V}^{e})(s) \, \md (\widetilde{\mu}_{V}^{e} \! - \! \mu_{V}^{e})
(t^{\prime});
\end{align*}
but, in Lemma~3.2, it was shown that
\begin{equation*}
\iint_{\mathbb{R}^{2}} \ln \! \left(\dfrac{\vert st^{\prime} \vert}{\vert s
\! - \! t^{\prime} \vert^{2}} \right) \! \md (\widetilde{\mu}_{V}^{e} \! - \!
\mu_{V}^{e})(s) \, \md (\widetilde{\mu}_{V}^{e} \! - \! \mu_{V}^{e})(t^{
\prime}) \! = \! 2 \int_{0}^{+\infty} \xi^{-1} \vert (\widehat{\widetilde{
\mu}_{V}^{e}} \! - \! \widehat{\mu_{V}^{e}})(\xi) \vert^{2} \, \md \xi \!
\geqslant \! 0,
\end{equation*}
whence $\int_{0}^{+\infty} \xi^{-1} \vert (\widehat{\widetilde{\mu}_{V}^{e}}
\! - \! \widehat{\mu_{V}^{e}})(\xi) \vert^{2} \, \md \xi \! = \! 0$
$\Rightarrow$ $\widehat{\widetilde{\mu}_{V}^{e}}(\xi) \! = \! \widehat{\mu_{
V}^{e}}(\xi)$, $\xi \! \geqslant \! 0$. Noting that
\begin{equation*}
\widehat{\widetilde{\mu}_{V}^{e}}(-\xi) \! = \! \int_{\mathbb{R}} \me^{\mi s
(-\xi)} \, \md \widetilde{\mu}_{V}^{e}(s) \! = \! \overline{\widehat{
\widetilde{\mu}_{V}^{e}}(\xi)} \qquad \text{and} \qquad \widehat{\mu_{V}^{e}}
(-\xi) \! = \! \int_{\mathbb{R}} \me^{\mi s(-\xi)} \, \md \mu_{V}^{e}(s) \! =
\! \overline{\widehat{\mu_{V}^{e}}(\xi)},
\end{equation*}
it follows {}from $\widehat{\widetilde{\mu}_{V}^{e}}(\xi) \! = \! \widehat{
\mu_{V}^{e}}(\xi)$, $\xi \! \geqslant \! 0$, via a complex-conjugation
argument, that $\widehat{\widetilde{\mu}_{V}^{e}}(-\xi) \! = \! \widehat{\mu_{
V}^{e}}(-\xi)$, $\xi \! \geqslant \! 0$; hence, $\widehat{\widetilde{\mu}_{
V}^{e}}(\xi) \! = \! \widehat{\mu_{V}^{e}}(\xi)$, $\xi \! \in \! \mathbb{R}$.
The latter relation shows that $\int_{\mathbb{R}} \me^{\mi s \xi} \, \md
(\widetilde{\mu}_{V}^{e} \! - \! \mu_{V}^{e})(s) \! = \! 0$ $\Rightarrow$
$\widetilde{\mu}_{V}^{e} \! = \! \mu_{V}^{e}$; thus the uniqueness of the
`even' equilibrium measure. \hfill $\qed$

Before proceeding to Lemma~3.4, the following observations, which are 
interesting, non-trivial and important results in their own right, should 
be noted. Let $\widetilde{V} \colon \mathbb{R} \setminus \{0\} \! \to \! 
\mathbb{R}$ satisfy conditions~(2.3)--(2.5). For each $n \! \in \! \mathbb{
N}$ and any $2n$-tuple $(x_{1},x_{2},\dotsc,x_{2n})$ of distinct, finite and 
non-zero real numbers, let
\begin{equation*}
\mathfrak{d}_{e,n}^{\widetilde{V}} \! := \! \dfrac{1}{2n(2n \! - \! 1)} \,
\inf_{\{x_{1},x_{2},\dotsc,x_{2n}\} \subset \mathbb{R} \setminus \{0\}} \left(
\sum_{\substack{j,k=1\\j \neq k}}^{2n} \ln \! \left(\left\vert x_{j} \! - \!
x_{k} \right\vert \! \left\vert x_{k}^{-1} \! - \! x_{j}^{-1} \right\vert
\right)^{-1} \! + \! 2(2n \! - \! 1) \sum_{i=1}^{2n} \widetilde{V}(x_{i})
\right).
\end{equation*}
For each $n \! \in \! \mathbb{N}$, a set $\left\lbrace x_{1}^{\sharp},x_{
2}^{\sharp},\dotsc,x_{2n}^{\sharp} \right\rbrace$ which realizes the above
infimum, that is,
\begin{equation*}
\mathfrak{d}_{e,n}^{\widetilde{V}} \! = \! \dfrac{1}{2n(2n \! - \! 1)} \!
\left(\sum_{\substack{j,k=1\\j \neq k}}^{2n} \ln \! \left(\left\vert x_{j}^{
\sharp} \! - \! x_{k}^{\sharp} \right\vert \! \left\vert (x_{k}^{\sharp})^{-1}
\! - \! (x_{j}^{\sharp})^{-1} \right\vert \right)^{-1} \! + \! 2(2n \! - \! 1)
\sum_{i=1}^{2n} \widetilde{V}(x_{i}^{\sharp}) \right),
\end{equation*}
will be called (with slight abuse of nomenclature) a \emph{generalised
weighted $2n$-Fekete set}, and the points $x_{1}^{\sharp},x_{2}^{\sharp},
\dotsc,x_{2n}^{\sharp}$ will be called \emph{generalised weighted Fekete
points}. For $\left\lbrace x_{1}^{\sharp},x_{2}^{\sharp},\dotsc,x_{2n}^{
\sharp} \right\rbrace$ a generalised weighted $2n$-Fekete set, denote by
\begin{equation*}
\mu_{\mathbf{x}^{\sharp}}^{e} \! := \! \dfrac{1}{2n} \sum_{j=1}^{2n} \delta_{
x_{j}^{\sharp}},
\end{equation*}
where $\delta_{x_{j}^{\sharp}}$, $j \! = \! 1,\dotsc,2n$, is the Dirac delta 
measure (atomic mass) concentrated at $x_{j}^{\sharp}$, the \emph{normalised 
counting measure}, that is, $\int_{\mathbb{R}} \md \mu_{\mathbf{x}^{\sharp}
}^{e}(s) \! = \! 1$. Then, mimicking the calculations in Chapter~6 of 
\cite{a90} and the techniques used to prove Theorem~1.34 in \cite{a56} (see, 
in particular, Section~2 of \cite{a56}), one proves that (the details are 
left to the interested reader):
\begin{compactenum}
\item[\textbullet] $\lim_{n \to \infty} \mathfrak{d}_{e,n}^{\widetilde{V}}$
exists, more precisely,
\begin{equation*}
\lim_{n \to \infty} \mathfrak{d}_{e,n}^{\widetilde{V}} \! = \! E_{V}^{e} \! =
\! \inf \lbrace \mathrm{I}_{V}^{e}[\mu^{e}]; \, \mu^{e} \! \in \! \mathcal{
M}_{1}(\mathbb{R}) \rbrace,
\end{equation*}
where (the functional) $\mathrm{I}_{V}^{e}[\mu^{e}] \colon \mathcal{M}_{1}
(\mathbb{R}) \! \to \! \mathbb{R}$ is defined in Lemma~3.1, and $\lim_{n \to
\infty} \exp (-\mathfrak{d}_{e,n}^{\widetilde{V}}) \! = \! \exp (-E_{V}^{e})$
is positive and finite;
\item[\textbullet] $\mu_{\mathbf{x}^{\sharp}}^{e}$ converges weakly (in the
weak-$\ast$ topology of measures) to the `even' equilibrium measure $\mu_{
V}^{e}$, that is, $\mu_{\mathbf{x}^{\sharp}}^{e} \overset{\ast}{\to} \mu_{V}^{
e}$ as $n \! \to \! \infty$.
\end{compactenum}

\textbf{RHP1}, that is, $(\overset{e}{\mathrm{Y}}(z),\mathrm{I} \! + \! \exp
(-n \widetilde{V}(z)) \sigma_{+},\mathbb{R})$, is now reformulated as an
equiv\-a\-l\-e\-n\-t, auxiliary RHP normalised at infinity.
\begin{notrem}
For completeness, the integrand appearing in the definition of $g^{e}(z)$ (see
Lemma~3.4 below) is defined as follows: $\ln ((z \! - \! s)^{2}(zs)^{-1}) \!
:= \! 2 \ln (z \! - \! s) \! - \! \ln z \! - \! \ln s$, where, for $s \!
< \! 0$, $\ln s \! := \! \ln \lvert s \rvert \! + \! \mi \pi$. \hfill
$\blacksquare$
\end{notrem}
\begin{ccccc}
Let the external field $\widetilde{V} \colon \mathbb{R} \setminus \{0\} \!
\to \! \mathbb{R}$ satisfy conditions~{\rm (2.3)--(2.5)}. For the associated
`even' equilibrium measure, $\mu_{V}^{e} \! \in \! \mathcal{M}_{1}(\mathbb{R}
)$, set $J_{e} \! := \! \operatorname{supp}(\mu_{V}^{e})$, where $J_{e}$
$(=$ compact) $\subset \! \overline{\mathbb{R}} \setminus \{0,\pm \infty\}$,
and let $\overset{e}{\mathrm{Y}} \colon \mathbb{C} \setminus \mathbb{R} \! \to
\! \mathrm{SL}_{2}(\mathbb{C})$ be the (unique) solution of {\rm \pmb{RHP1}}. 
Let
\begin{equation*}
\overset{e}{\mathscr{M}}(z) \! := \! \me^{-\frac{n \ell_{e}}{2} \mathrm{ad}
(\sigma_{3})} \overset{e}{\mathrm{Y}}(z) \me^{-n(g^{e}(z)+Q_{e}) \sigma_{3}},
\end{equation*}
where $g^{e}(z)$, the `even' $g$-function, is defined by
\begin{equation*}
g^{e}(z) \! := \! \int_{J_{e}} \ln \! \left((z \! - \! s)^{2}(zs)^{-1} \right)
\md \mu_{V}^{e}(s), \quad z \! \in \! \mathbb{C} \setminus (-\infty,\max \{0,
\max \{\operatorname{supp}(\mu_{V}^{e})\}\}),
\end{equation*}
$\ell_{e}$ $(\in \! \mathbb{R})$, the `even' variational constant, is given
in {\rm Lemma~3.6} below, and
\begin{equation*}
Q_{e} \! := \! \int_{J_{e}} \ln (s) \md \mu_{V}^{e}(s).
\end{equation*}
Then $\overset{e}{\mathscr{M}} \colon \mathbb{C} \setminus \mathbb{R} \! \to
\! \mathrm{SL}_{2}(\mathbb{C})$ solves the following (normalised at infinity)
{\rm RHP:} {\rm (i)} $\overset{e}{\mathscr{M}}(z)$ is holomorphic for $z \!
\in \! \mathbb{C} \setminus \mathbb{R};$ {\rm (ii)} the boundary values
$\overset{e}{\mathscr{M}}_{\pm}(z) := \! \lim_{\underset{\pm \Im (z^{\prime})
>0}{z^{\prime} \to z}} \overset{e}{\mathscr{M}}(z^{\prime})$ satisfy the jump 
condition
\begin{equation*}
\overset{e}{\mathscr{M}}_{+}(z) \! = \! \overset{e}{\mathscr{M}}_{-}(z) \!
\begin{pmatrix}
\me^{-n(g^{e}_{+}(z)-g^{e}_{-}(z))} & \me^{n(g^{e}_{+}(z)+g^{e}_{-}(z)-
\widetilde{V}(z)-\ell_{e}+2Q_{e})} \\
0 & \me^{n(g^{e}_{+}(z)-g^{e}_{-}(z))}
\end{pmatrix}, \quad z \! \in \! \mathbb{R},
\end{equation*}
where $g^{e}_{\pm}(z) \! := \! \lim_{\varepsilon \downarrow 0}g^{e}(z \! \pm
\! \mi \varepsilon);$ {\rm (iii)} $\overset{e}{\mathscr{M}}(z) \! =_{
\underset{z \in \mathbb{C} \setminus \mathbb{R}}{z \to \infty}} \! \mathrm{I}
\! + \! \mathcal{O}(z^{-1});$ and {\rm (iv)} $\overset{e}{\mathscr{M}}(z) \!
=_{\underset{z \in \mathbb{C} \setminus \mathbb{R}}{z \to 0}} \! \mathcal{O}
(1)$.
\end{ccccc}

\emph{Proof.} For (arbitrary) $z_{1},z_{2} \! \in \! \mathbb{C}_{\pm}$, note
that, {}from the definition of $g^{e}(z)$ stated in the Lemma, $g^{e}(z_{2})
\! - \! g^{e}(z_{1}) \! = \! \mi \pi \int_{z_{1}}^{z_{2}} \mathscr{F}^{e}(s)
\, \md s$, where
\begin{equation*}
\mathscr{F}^{e} \colon \mathbb{C} \setminus (\operatorname{supp}(\mu_{V}^{e})
\cup \{0\}) \! \to \! \mathbb{C},\, \, z \! \mapsto \! -\dfrac{1}{\pi \mi} \!
\left(\dfrac{1}{z} \! + \! 2 \int_{J_{e}} \dfrac{\md \mu_{V}^{e}(s)}{s \! - \!
z} \right),
\end{equation*}
with $\mathscr{F}^{e}(z) \! =_{z \to 0} \! -\tfrac{1}{\pi \mi z} \! + \!
\mathcal{O}(1)$ (since $\mu_{V}^{e} \! \in \! \mathcal{M}_{1}(\mathbb{R})$; in
particular, $\int_{\mathbb{R}}s^{m} \, \md \mu_{V}^{e}(s) \! < \! \infty$, $m
\! \in \! \mathbb{Z})$; thus, $\vert g^{e}(z_{2}) \! - \! g^{e}(z_{1}) \vert
\! \leqslant \! \pi \sup_{z \in \mathbb{C}_{\pm}} \vert \mathscr{F}^{e}(z)
\vert \vert z_{2} \! - \! z_{1} \vert$, that is, $g^{e}(z)$ is uniformly
Lipschitz continuous in $\mathbb{C}_{\pm}$. Thus, {}from the definition of
$g^{e}(z)$ stated in the Lemma:
\begin{compactenum}
\item[(1)] for $s \! \in \! J_{e}$, $z \! \in \! \mathbb{C} \setminus
(-\infty,\max \{0,\max \{\operatorname{supp}(\mu_{V}^{e})\}\})$, with $\vert
s/z \vert \! \ll \! 1$, and $\mu_{V}^{e} \! \in \! \mathcal{M}_{1}(\mathbb{R}
)$, in particular, $\int_{\mathbb{R}} \md \mu_{V}^{e}(s)$ $(= \! \int_{J_{e}}
\md \mu_{V}^{e}(s))$ $= \! 1$ and $\int_{\mathbb{R}} s^{m} \, \md \mu_{V}^{e}
(s)$ $(= \! \int_{J_{e}}s^{m} \, \md \mu_{V}^{e}(s))$ $< \! \infty$, $m \! \in
\! \mathbb{N}$, it follows {}from the expansions $\tfrac{1}{s-z} \! = \! -
\sum_{k=0}^{l} \tfrac{s^{k}}{z^{k+1}} \! + \! \tfrac{s^{l+1}}{z^{l+1}(s-z)}$,
$l \! \in \! \mathbb{Z}_{0}^{+}$, and $\ln (z \! - \! s) \! =_{\vert z \vert
\to \infty} \! \ln (z) \! - \! \sum_{k=1}^{\infty} \tfrac{1}{k}(\tfrac{s}{z}
)^{k}$, that
\begin{equation*}
g^{e}(z) \! \underset{\mathbb{C} \setminus \mathbb{R} \ni z \to \infty}{=} 
\ln (z) \! - \! Q_{e} \! + \! \mathcal{O}(z^{-1}),
\end{equation*}
where $Q_{e}$ is defined in the Lemma;
\item[(2)] for $s \! \in \! J_{e}$, $z \! \in \! \mathbb{C} \setminus (-
\infty,\max \{0,\max \{\operatorname{supp}(\mu_{V}^{e})\}\})$, with $\vert z/s
\vert \! \ll \! 1$, and $\mu_{V}^{e} \! \in \! \mathcal{M}_{1}(\mathbb{R})$,
in particular, $\int_{\mathbb{R}}s^{-m} \, \md \mu_{V}^{e}(s)$ $(= \! \int_{
J_{e}} s^{-m} \, \md \mu_{V}^{e}(s))$ $< \! \infty$, $m \! \in \! \mathbb{N}$,
it follows {}from the expansions $\tfrac{1}{z-s} \! = \! -\sum_{k=0}^{l}
\tfrac{z^{k}}{s^{k+1}} \! + \! \tfrac{z^{l+1}}{s^{l+1}(z-s)}$, $l \! \in \!
\mathbb{Z}_{0}^{+}$, and $\ln (s \! - \! z) \! =_{\vert z \vert \to 0} \! \ln
(s) \! - \! \sum_{k=1}^{\infty} \tfrac{1}{k}(\tfrac{z}{s})^{k}$, that
\begin{equation*}
g^{e}(z) \! \underset{\mathbb{C}_{\pm} \ni z \to 0}{=} \! -\ln(z) \! - \! 
Q_{e} \! + \! 2 \int_{J_{e}} \ln (\lvert s \rvert) \, \md \mu_{V}^{e}(s) \! 
\pm \! 2 \pi \mi \int_{J_{e} \cap \mathbb{R}_{+}} \md \mu_{V}^{e}(s) \! + \! 
\mathcal{O}(z),
\end{equation*}
where (see Lemma~3.5, item~(1), below)
\begin{equation*}
\int_{J_{e} \cap \mathbb{R}_{+}} \md \mu_{V}^{e}(s) \! = \!
\begin{cases}
0, &\text{$J_{e} \! \subset \! \mathbb{R}_{-}$,} \\
1, &\text{$J_{e} \! \subset \! \mathbb{R}_{+}$,} \\
\int_{b_{j}^{e}}^{a_{N+1}^{e}} \md \mu_{V}^{e}(s), &\text{$(a_{j}^{e},b_{j}^{
e}) \! \ni \! 0, \quad j \! = \! 1,\dotsc,N$.}
\end{cases}
\end{equation*}
\end{compactenum}
Items~(i)--(iv) now follow {}from the definitions of $\overset{e}{\mathscr{M}
}(z)$ (in terms of $\overset{e}{\mathrm{Y}}(z))$ and $g^{e}(z)$ stated in the
Lemma, and the above two asymptotic expansions. \hfill $\qed$

\begin{ccccc}
Let the external field $\widetilde{V} \colon \mathbb{R} \setminus \{0\} \! \to
\! \mathbb{R}$ satisfy conditions~{\rm (2.3)--(2.5)}. For $\mu_{V}^{e} \! \in
\! \mathcal{M}_{1}(\mathbb{R})$, the associated `even' equilibrium measure,
set $J_{e} \! := \! \operatorname{supp}(\mu_{V}^{e})$, where $J_{e}$ $(=$
compact) $\subset \! \overline{\mathbb{R}} \setminus \{0,\pm \infty\}$. Then: 
{\rm (1)} $J_{e} \! = \! \cup_{j=1}^{N+1}(b_{j-1}^{e},a_{j}^{e})$, with $N \! 
\in \! \mathbb{N}$ and finite, $b_{0}^{e} \! := \! \min \{\operatorname{supp}
(\mu_{V}^{e})\} \! \notin \! \lbrace -\infty,0 \rbrace$, $a_{N+1}^{e} \! := 
\! \max \{\operatorname{supp}(\mu_{V}^{e})\} \! \notin \! \lbrace 0,+\infty 
\rbrace$, and $-\infty \! \! < \! b_{0}^{e} \! < \! a_{1}^{e} \! < \! b_{1}^{
e} \! < \! a_{2}^{e} \! < \! \cdots \! < \! b_{N}^{e} \! < \! a_{N+1}^{e} \! 
< \! +\infty$, and $\lbrace b_{j-1}^{e},a_{j}^{e} \rbrace_{j=1}^{N+1}$ satisfy 
the $n$-dependent and (locally) solvable system of $2(N \! + \! 1)$ moment 
conditions
\begin{gather*}
\int_{J_{e}} \dfrac{(\frac{2 \mi}{\pi s} \! + \! \frac{\mi \widetilde{V}^{
\prime}(s)}{\pi})s^{j}}{(R_{e}(s))^{1/2}_{+}} \, \dfrac{\md s}{2 \pi \mi} \! =
\! 0, \quad j \! = \! 0,\dotsc,N, \qquad \qquad \int_{J_{e}} \dfrac{(\frac{2
\mi}{\pi s} \! + \! \frac{\mi \widetilde{V}^{\prime}(s)}{\pi})s^{N+1}}{(R_{e}
(s))^{1/2}_{+}} \, \dfrac{\md s}{2 \pi \mi} \! = \! \dfrac{2}{\pi \mi}, \\
\int_{a_{j}^{e}}^{b_{j}^{e}} \! \left(\mi (R_{e}(s))^{1/2} \int_{J_{e}}
\dfrac{(\frac{\mi}{\pi \xi} \! + \! \frac{\mi \widetilde{V}^{\prime}(\xi)}{2
\pi})}{(R_{e}(\xi))^{1/2}_{+}(\xi \! - \! s)} \, \dfrac{\md \xi}{2 \pi \mi}
\right) \! \md s \! = \! \dfrac{1}{2 \pi} \ln \! \left\vert \dfrac{a_{j}^{e}
}{b_{j}^{e}} \right\vert \! + \! \dfrac{1}{4 \pi} \! \left(\widetilde{V}(a_{
j}^{e}) \! - \! \widetilde{V}(b_{j}^{e}) \right), \quad j \! = \! 1,\dotsc,N,
\end{gather*}
where $(R_{e}(z))^{1/2}$ is defined in Theorem~{\rm 2.3.1},
Equation~{\rm (2.8)}, with $(R_{e}(z))^{1/2}_{\pm} \! := \! \lim_{\varepsilon
\downarrow 0}(R_{e}(z \! \pm \! \mi \varepsilon))^{1/2}$, and the branch of
the square root chosen so that $z^{-(N+1)}(R_{e}(z))^{1/2} \! \sim_{\underset{
z \in \mathbb{C}_{\pm}}{z \to \infty}} \! \pm 1;$ and {\rm (2)} the density of
the `even' equilibrium measure, which is absolutely continuous with respect to
Lebesgue measure, is given by
\begin{equation*}
\md \mu_{V}^{e}(x) \! := \! \psi_{V}^{e}(x) \, \md x \! = \! \dfrac{1}{2 \pi
\mi}(R_{e}(x))^{1/2}_{+}h_{V}^{e}(x) \pmb{1}_{J_{e}}(x) \, \md x,
\end{equation*}
where
\begin{equation*}
h_{V}^{e}(z) \! := \! \dfrac{1}{2} \oint_{C_{\mathrm{R}}^{e}} \dfrac{(\frac{
\mi}{\pi s} \! + \! \frac{\mi \widetilde{V}^{\prime}(s)}{2 \pi})}{(R_{e}(s))^{
1/2}(s \! - \! z)} \, \md s
\end{equation*}
(real analytic for $z \! \in \! \mathbb{R} \setminus \{0\})$, with $C_{
\mathrm{R}}^{e}$ $(\subset \mathbb{C}^{\ast})$ the boundary of any open
doubly-connected annular region of the type $\{\mathstrut z^{\prime} \! \in
\! \mathbb{C}; \, 0 \! < \! r \! < \! \vert z^{\prime} \vert \! < \! R \! < \!
+\infty\}$, where the simple outer (resp., inner) boundary $\{\mathstrut z^{
\prime} \! = \! R \me^{\mi \vartheta}, \, 0 \! \leqslant \! \vartheta \!
\leqslant \! 2 \pi\}$ (resp., $\{\mathstrut z^{\prime} \! = \! r \me^{\mi
\vartheta}, \, 0 \! \leqslant \! \vartheta \! \leqslant \! 2 \pi\})$ is
traversed clockwise (resp., counter-clockwise), with the numbers $0 \! < \! r
\! < \! R \! < \! +\infty$ chosen such that, for (any) non-real $z$ in the
domain of analyticity of $\widetilde{V}$ (that is, $\mathbb{C}^{\ast})$,
$\operatorname{int}(C_{\mathrm{R}}^{e}) \supset J_{e} \cup \{z\}$, $\pmb{1}_{
J_{e}}(x)$ is the indicator (characteristic) function of the set $J_{e}$, and
$\psi_{V}^{e}(x) \! \geqslant \! 0$ (resp., $\psi_{V}^{e}(x) \! > \! 0)$
$\forall \, \, x \! \in \! \overline{J_{e}}:= \! \cup_{j=1}^{N+1}[b_{j-1}^{e},
a_{j}^{e}]$ (resp., $\forall \, \, x \! \in \! J_{e})$.
\end{ccccc}

\emph{Proof.} One begins by showing that the support of the `even' equilibrium
measure, $\operatorname{supp}(\mu_{V}^{e}) \! =: \! J_{e}$, consists of the
union of a finite number of disjoint and bounded (real) intervals. Recall
{}from Lemma~3.1 that $J_{e} \! = \text{compact} \subset \overline{\mathbb{
R}} \setminus \lbrace 0,\pm \infty \rbrace$, and that $\widetilde{V}$ is real
analytic on $\mathbb{R} \setminus \{0\}$, thus real analytic on $J_{e}$, with
an analytic continuation to the following (open) neighbourhood of $J_{e}$,
$\mathbb{U} \! := \! \lbrace \mathstrut z \! \in \! \mathbb{C}; \, \inf_{q 
\in J_{e}} \vert z \! - \! q \vert \! < \! r \! \in \! (0,1) \rbrace \setminus 
\{0\}$. In analogy with Equation~(2.1) of \cite{a56}, for each $n \! \in \!
\mathbb{N}$ and any $2n$-tuple $(x_{1},x_{2},\dotsc,x_{2n})$ of distinct,
finite and non-zero real numbers, let
\begin{align*}
\mathrm{d}_{\widetilde{V},n}^{e} :=& \, \left(\sup_{\{x_{1},x_{2},\dotsc,
x_{2n}\} \subset \mathbb{R} \setminus \{0\}} \, \prod_{\substack{j,k=1\\j<
k}}^{2n} \left\vert x_{j} \! - \! x_{k} \right\vert^{2} \! \left\vert x_{
k}^{-1} \! - \! x_{j}^{-1} \right\vert^{2} \, \me^{-2 \widetilde{V}(x_{j})}
\me^{-2 \widetilde{V}(x_{k})} \right)^{\frac{1}{2n(2n-1)}} \\
=& \, \left(\sup_{\{x_{1},x_{2},\dotsc,x_{2n}\} \subset \mathbb{R} \setminus
\{0\}} \, \prod_{\substack{j,k=1\\j<k}}^{2n} \left\vert x_{j} \! - \! x_{k}
\right\vert^{2} \! \left\vert x_{k}^{-1} \! - \! x_{j}^{-1} \right\vert^{2}
\, \me^{-2(2n-1) \sum_{i=1}^{2n} \widetilde{V}(x_{i})} \right)^{\frac{1}{2n
(2n-1)}},
\end{align*}
where $\prod_{\substack{j,k=1\\j<k}}^{2n}(\star) \! = \! \prod_{j=1}^{2n-1}
\prod_{k=j+1}^{2n}(\star)$. Denote by $\left\lbrace x^{\ast}_{1},x^{\ast}_{2},
\dotsc,x^{\ast}_{2n} \right\rbrace$, with $ x^{\ast}_{i} \! < \! x^{\ast}_{j}
\, \, \forall \, \, i \! < \! j \! \in \! \lbrace 1,\dotsc,2n \rbrace$, the
associated generalised weighted $2n$-Fekete set (see the discussion preceding
Lemma~3.4), that is,
\begin{equation*}
\mathrm{d}_{\widetilde{V},n}^{e} \! = \! \left(\prod_{\substack{j,k=1\\j<k}}^{
2n} \left\vert x_{j}^{\ast} \! - \! x_{k}^{\ast} \right\vert^{2} \! \left\vert
(x_{k}^{\ast})^{-1} \! - \! (x_{j}^{\ast})^{-1} \right\vert^{2} \, \me^{-2(2n-
1) \sum_{i=1}^{2n} \widetilde{V}(x_{i}^{\ast})} \right)^{\frac{1}{2n(2n-1)}}.
\end{equation*}
Proceeding, now, as in the proof of Theorem~1.34, Equation~(1.35), of 
\cite{a56}, in particular, mimicking the calculations on pp.~408--413 of 
\cite{a56} (for the proofs of Lemmas~2.3 and~2.15 therein), namely, using 
those techniques to show that, in the present case, the nearest-neighbour 
distances $\lbrace x_{j+1}^{\ast} \! - \! x_{j}^{\ast} \rbrace_{j=1}^{2n-1}$ 
are not `too small' as $n \! \to \! \infty$, and the calculations on 
pp.~413--415 of \cite{a56} (for the proof of Lemma~2.26 therein), one shows 
that, for the regular case considered herein (cf. Subsection~2.2), the `even' 
equilibrium measure, $\mu_{V}^{e}$ $(\in \! \mathcal{M}_{1}(\mathbb{R}))$, is 
absolutely continuous with respect to Lebesgue measure, that is, the density 
of the `even' equilibrium measure has the representation $\md \mu_{V}^{e}(x) 
\! := \! \psi_{V}^{e}(x) \, \md x$, $x \! \in \! \operatorname{supp}(\mu_{
V}^{e})$, where $\psi_{V}^{e}(x) \! \geqslant \! 0$ on $\overline{J_{e}}$, 
with $\psi_{V}^{e}(\pmb{\cdot})$ determined (explicitly) below\footnote{The 
analysis of \cite{a56} is, in some sense, more complicated than the one of 
the present paper, because, unlike the `real-line' case considered herein, 
that is, $\operatorname{supp}(\mu_{V}^{e}) \! =: \! J_{e} \subset \overline{
\mathbb{R}} \setminus \lbrace 0,\pm \infty \rbrace$, the end-point effects 
at $\pm 1$ in \cite{a56} require special consideration (see, also, Section~4 
of \cite{a56}).}.

Set
\begin{equation}
\mathscr{H}^{e}(z) \! := \! (\mathscr{F}^{e}(z))^{2} \! - \! \int_{J_{e}}
\dfrac{(16 \mi \psi_{V}^{e}(\xi)(\mathcal{H} \psi_{V}^{e})(\xi) \! - \! 8 \mi
\psi_{V}^{e}(\xi)/\pi \xi)}{(\xi \! - \! z)} \, \dfrac{\md \xi}{2 \pi \mi},
\quad z \! \in \! \mathbb{C} \setminus (J_{e} \cup \{0\}),
\end{equation}
where, {}from the proof of Lemma~3.4,
\begin{equation}
\mathscr{F}^{e}(z) \! = \! -\dfrac{1}{\pi \mi} \! \left(\dfrac{1}{z} \! + \!
2 \int_{J_{e}} \dfrac{\md \mu_{V}^{e}(s)}{s-z} \right),
\end{equation}
with $\int_{J_{e}} \tfrac{\md \mu_{V}^{e}(s)}{s-z}$ the Stieltjes transform of
the `even' equilibrium measure, and
\begin{equation*}
\mathcal{H} \colon \mathcal{L}^{2}_{\mathrm{M}_{2}(\mathbb{C})} \! \to \!
\mathcal{L}^{2}_{\mathrm{M}_{2}(\mathbb{C})}, \, \, f \! \mapsto \! (\mathcal{
H}f)(z) \! := \! \vip_{\raise-1.95ex\hbox{$\scriptstyle{} \mathbb{R}$}}
\dfrac{f(s)}{z \! - \! s} \, \dfrac{\md s}{\pi}
\end{equation*}
denotes the Hilbert transform, with $\pvi_{}$ denoting the principle value
integral. Via the distributional identities $\tfrac{1}{x-(x_{0} \pm \mi 0)} \!
= \! \tfrac{1}{x-x_{0}} \! \pm \! \mi \pi \delta (x \! - \! x_{0})$, with
$\delta (\cdot)$ the Dirac delta function, and $\int_{\xi_{1}}^{\xi_{2}}f(\xi)
\delta (\xi \! - \! x) \md \xi \! = \!
\begin{cases}
f(x), &\text{$x \! \in \! (\xi_{1},\xi_{2})$,} \\
0, &\text{$x \! \in \! \mathbb{R} \setminus (\xi_{1},\xi_{2}),$}
\end{cases}$ it follows that
\begin{align*}
&\mathscr{H}^{e}_{\pm}(z) \! = \\
&\begin{cases}
(\mathscr{F}^{e}_{\pm}(z))^{2} \! - \! \pvi_{\raise-0.95ex\hbox{$\scriptstyle{}
J_{e}$}} \dfrac{(16 \mi \psi_{V}^{e}(\xi)(\mathcal{H} \psi_{V}^{e})(\xi) \! -
\! \frac{8 \mi \psi_{V}^{e}(\xi)}{\pi \xi})}{(\xi \! - \! z)} \, \dfrac{\md
\xi}{2 \pi \mi} \! \mp \! \dfrac{1}{2} \! \left(16 \mi \psi_{V}^{e}(z)
(\mathcal{H} \psi_{V}^{e})(z) \! - \! \dfrac{8 \mi \psi_{V}^{e}(z)}{\pi z}
\right), &\text{$z \! \in \! J_{e}$,} \\
(\mathscr{F}^{e}_{\pm}(z))^{2} \! - \! \int_{J_{e}} \dfrac{(16 \mi \psi_{V}^{e}
(\xi)(\mathcal{H} \psi_{V}^{e})(\xi) \! - \! \frac{8 \mi \psi_{V}^{e}(\xi)}{
\pi \xi})}{(\xi \! - \! z)} \, \dfrac{\md \xi}{2 \pi \mi}, &\text{$z \!
\notin \! J_{e}$,}
\end{cases}
\end{align*}
where $\star^{e}_{\pm}(z) \! := \! \lim_{\varepsilon \downarrow 0} \star^{e}
(z \! \pm \! \mi 0)$, $\star \! \in \! \{\mathscr{H},\mathscr{F}\}$. Recall
the definition of $g^{e}(z)$ given in Lemma~3.4:
\begin{equation*}
g^{e}(z) \! := \! \int_{J_{e}} \ln \! \left(\dfrac{(z \! - \! s)^{2}}{zs}
\right) \! \md \mu_{V}^{e}(s) \! = \! \int_{J_{e}} \ln \! \left(\dfrac{(z \! -
\! s)^{2}}{zs} \right) \! \psi_{V}^{e}(s) \, \md s, \quad z \! \in \! \mathbb{
C} \setminus (-\infty,\max \{0,\max \{J_{e}\}\});
\end{equation*}
using the above distributional identities and the fact that $\int_{J_{e}}
\psi_{V}^{e}(s) \, \md s \! = \! 1$, one shows that
\begin{equation*}
(g^{e}_{\pm}(z))^{\prime} \! := \! \lim_{\varepsilon \downarrow 0}(g^{e})^{
\prime}(z \! \pm \! \mi \varepsilon) \! = \!
\begin{cases}
-\tfrac{1}{z} \! - \! 2 \pvi_{\raise-1.05ex\hbox{$\scriptstyle{}J_{e}$}}
\tfrac{\psi_{V}^{e}(s)}{s-z} \, \md s \! \mp \! 2 \pi \mi \psi_{V}^{e}(z),
&\text{$z \! \in \! J_{e}$,} \\
-\tfrac{1}{z} \! - \! 2 \int_{J_{e}} \tfrac{\psi_{V}^{e}(s)}{s-z} \, \md s,
&\text{$z \! \notin \! J_{e}$,}
\end{cases}
\end{equation*}
whence one concludes that
\begin{align*}
(g^{e}_{+} \! + \! g^{e}_{-})^{\prime}(z) =& \, -\dfrac{2}{z} \! - \! 4 \vip_{
\raise-1.95ex\hbox{$\scriptstyle{}J_{e}$}} \dfrac{\psi_{V}^{e}(s)}{s \! - \!
z} \, \md s \! = \! -\dfrac{2}{z} \! + \! 4 \pi (\mathcal{H} \psi_{V}^{e})(z),
\quad z \! \in \! J_{e}, \\
(g^{e}_{+} \! - \! g^{e}_{-})^{\prime}(z) =& \,
\begin{cases}
-4 \pi \mi \psi_{V}^{e}(z), &\text{$z \! \in \! J_{e}$,} \\
0, &\text{$z \! \notin \! J_{e}$.}
\end{cases}
\end{align*}
Demanding that (see Lemma~3.6 below) $(g^{e}_{+} \! + \! g^{e}_{-})^{\prime}
(z) \! = \! \widetilde{V}^{\prime}(z)$, $z \! \in \! J_{e}$, one shows {}from
the above that, for $J_{e} \! \ni \! z$, $((g^{e}(z))^{\prime} \! + \! \tfrac{
1}{z})_{+} \! + \! ((g^{e}(z))^{\prime} \! + \! \tfrac{1}{z})_{-} \! = \! 4
\pi (\mathcal{H} \psi_{V}^{e})(z) \! = \! \tfrac{2}{z} \! + \! \widetilde{V}^{
\prime}(z)$ $\Rightarrow$
\begin{equation}
(\mathcal{H} \psi_{V}^{e})(z) \! = \! \dfrac{1}{2 \pi z} \! + \! \dfrac{
\widetilde{V}^{\prime}(z)}{4 \pi}, \quad z \! \in \! J_{e}.
\end{equation}
{}From Equation~(3.2) and the above distributional identities, one shows that
\begin{equation}
\mathscr{F}^{e}_{\pm}(z) \! := \! \lim_{\varepsilon \downarrow 0} \mathscr{
F}^{e}(z \! \pm \! \mi \varepsilon) \! = \!
\begin{cases}
-\tfrac{1}{\pi \mi z} \! - \! 2 \mi (\mathcal{H} \psi_{V}^{e})(z) \! \mp \! 2
\psi_{V}^{e}(z), &\text{$z \! \in \! J_{e}$,} \\
-\tfrac{1}{\pi \mi} \! \left(\tfrac{1}{z} \! + \! 2 \int_{J_{e}} \tfrac{\psi_{
V}^{e}(s)}{s-z} \, \md s \right), &\text{$z \! \notin \! J_{e}$;}
\end{cases}
\end{equation}
thus, for $z \! \in \! \mathbb{R} \setminus (J_{e} \cup \{0\})$, $\mathscr{F}^{
e}_{+}(z) \! = \! \mathscr{F}^{e}_{-}(z) \! = \! -\tfrac{1}{\pi \mi}(\tfrac{1}{
z} \! + \! 2 \int_{J_{e}} \tfrac{\psi_{V}^{e}(s)}{s-z} \, \md s)$. Hence, for
$z \! \notin \! J_{e} \cup \{0\}$, one deduces that $\mathscr{H}^{e}_{+}(z)
\! = \! \mathscr{H}^{e}_{-}(z)$. For $z \! \in \! J_{e}$, one notes that
\begin{equation*}
\mathscr{H}^{e}_{+}(z) \! - \! \mathscr{H}^{e}_{-}(z) \! = \! (\mathscr{F}^{
e}_{+}(z))^{2} \! - \! (\mathscr{F}^{e}_{-}(z))^{2} \! - \! 16 \mi \psi_{V}^{e}
(z)(\mathcal{H} \psi_{V}^{e})(z) \! + \! \dfrac{8 \mi \psi_{V}^{e}(z)}{\pi z},
\end{equation*}
and
\begin{equation*}
(\mathscr{F}^{e}_{\pm}(z))^{2} \! = \! -\dfrac{1}{\pi^{2}z^{2}} \! + \!
\dfrac{4(\mathcal{H} \psi_{V}^{e})(z)}{\pi z} \! \mp \! \dfrac{4 \mi \psi_{
V}^{e}(z)}{\pi z} \! - \! 4((\mathcal{H} \psi_{V}^{e})(z))^{2} \! \pm \! 8
\mi \psi_{V}^{e}(z)(\mathcal{H} \psi_{V}^{e})(z) \! + \! (2 \psi_{V}^{e}
(z))^{2},
\end{equation*}
whence $(\mathscr{F}^{e}_{+}(z))^{2} \! - \! (\mathscr{F}^{e}_{-}(z))^{2} \!
= \! -\tfrac{8 \mi \psi_{V}^{e}(z)}{\pi z} \! + \! 16 \mi \psi_{V}^{e}(z)
(\mathcal{H} \psi_{V}^{e})(z)$ $\Rightarrow$ $\mathscr{H}^{e}_{+}(z) \! - \!
\mathscr{H}^{e}_{-}(z) \! = \! 0$; thus, for $z \! \in \! J_{e}$, $\mathscr{
H}^{e}_{+}(z) \! = \! \mathscr{H}^{e}_{-}(z)$. The above argument shows,
therefore, that $\mathscr{H}^{e}(z)$ is analytic across $\mathbb{R} \setminus
\{0\}$; in fact, $\mathscr{H}^{e}(z)$ is entire for $z \! \in \! \mathbb{C}^{
\ast}$. Recalling that $\mu_{V}^{e} \! \in \! \mathcal{M}_{1}(\mathbb{R})$,
in particular, $\int_{J_{e}}s^{-m} \, \md \mu_{V}^{e}(s) \! = \! \int_{J_{e}}
s^{-m} \psi_{V}^{e}(s) \, \md s \! < \! \infty$, $m \! \in \! \mathbb{N}$, one
shows that, for $\vert z/s \vert \! \ll \! 1$, with $s \! \in \! J_{e}$ and $z
\! \notin \! J_{e}$, via the expansion $\tfrac{1}{z-s} \! = \! -\sum_{k=0}^{l}
\tfrac{z^{k}}{s^{k+1}} \! + \! \tfrac{z^{l+1}}{s^{l+1}(z-s)}$, $l \! \in \!
\mathbb{Z}_{0}^{+}$,
\begin{equation*}
(\mathscr{F}^{e}(z))^{2} \underset{z \to 0}{=} \! -\dfrac{1}{\pi^{2}z^{2}}
\! - \! \dfrac{1}{z} \! \left(\dfrac{4}{\pi^{2}} \int_{J_{e}}s^{-1} \, \md
\mu_{V}^{e}(s) \right) \! + \! \mathcal{O}(1),
\end{equation*}
whence, upon recalling the definition of $\mathscr{H}^{e}(z)$, in particular,
for $\vert z/\xi \vert \! \ll \! 1$, with $\xi \! \in \! J_{e}$ and $z \!
\notin \! J_{e}$, via the expansion $\tfrac{1}{z-\xi} \! = \! -\sum_{k=0}^{l}
\tfrac{z^{k}}{\xi^{k+1}} \! + \! \tfrac{z^{l+1}}{\xi^{l+1}(z-\xi)}$, $l \! \in
\! \mathbb{Z}_{0}^{+}$,
\begin{equation*}
\int_{J_{e}} \dfrac{(16 \mi \psi_{V}^{e}(\xi)(\mathcal{H} \psi_{V}^{e})(\xi)
\! - \! \frac{8 \mi \psi_{V}^{e}(\xi)}{\pi \xi})}{(\xi \! - \! z)} \, \dfrac{
\md \xi}{2 \pi \mi} \underset{z \to 0}{=} \! \mathcal{O}(1),
\end{equation*}
it follows that
\begin{equation*}
\mathscr{H}^{e}(z) \! \underset{z \to 0}{=} \! -\dfrac{1}{\pi^{2}z^{2}} \! -
\! \dfrac{1}{z} \! \left(\dfrac{4}{\pi^{2}} \int_{J_{e}}s^{-1} \, \md \mu_{
V}^{e}(s) \right) \! + \! \mathcal{O}(1),
\end{equation*}
which shows that $\mathscr{H}^{e}(z)$ has a pole of order $2$ at $z \! = \!
0$, with $\operatorname{Res}(\mathscr{H}^{e}(z);0) \! = \! -4 \pi^{-2} \int_{
J_{e}}s^{-1} \, \md \mu_{V}^{e}(s)$. One learns {}from the above analysis that
$z^{2} \mathscr{H}^{e}(z)$ is entire: look, in particular, at the behaviour
of $z^{2} \mathscr{H}^{e}(z)$ as $\vert z \vert \! \to \! \infty$. Recalling
Equations~(3.1) and~(3.2), one shows that, for $\mu_{V}^{e} \! \in \! \mathcal{
M}_{1}(\mathbb{R})$, in particular, $\int_{J_{e}} \md \mu_{V}^{e}(s) \! = \!
1$ and $\int_{J_{e}}s^{m} \, \md \mu_{V}^{e}(s) \! < \! \infty$, $m \! \in \!
\mathbb{N}$, for $\vert s/z \vert \! \ll \! 1$, with $s \! \in \! J_{e}$ and
$z \! \notin \! J_{e}$, via the expansion $\tfrac{1}{s-z} \! = \! -\sum_{k=
0}^{l} \tfrac{s^{k}}{z^{k+1}} \! + \! \tfrac{s^{l+1}}{z^{l+1}(s-z)}$, $l \!
\in \! \mathbb{Z}_{0}^{+}$,
\begin{align*}
z^{2} \mathscr{H}^{e}(z) \, +& \, \dfrac{1}{\pi^{2}} \! - \! \int_{J_{e}}s \!
\left(16 \mi \psi_{V}^{e}(s)(\mathcal{H} \psi_{V}^{e})(s) \! - \! \dfrac{8
\mi \psi_{V}^{e}(s)}{\pi s} \right) \dfrac{\md s}{2 \pi \mi} \\
-& \, z \int_{J_{e}} \! \left(16 \mi \psi_{V}^{e}(s)(\mathcal{H} \psi_{V}^{e})
(s) \! - \! \dfrac{8 \mi \psi_{V}^{e}(s)}{\pi s} \right) \dfrac{\md s}{2 \pi
\mi} \! \underset{\vert z \vert \to \infty}{=} \! \mathcal{O}(z^{-1});
\end{align*}
thus, due to the entirety of $\mathscr{H}^{e}(z)$, it follows, by a 
generalisation of Liouville's Theorem, that
\begin{align*}
z^{2} \mathscr{H}^{e}(z) \, +& \, \dfrac{1}{\pi^{2}} \! - \! \int_{J_{e}}
\! s \! \left(16 \mi \psi_{V}^{e}(s)(\mathcal{H} \psi_{V}^{e})(s) \! - \!
\dfrac{8 \mi \psi_{V}^{e}(s)}{\pi s} \right) \dfrac{\md s}{2 \pi \mi} \\
-& \, z \int_{J_{e}} \! \left(16 \mi \psi_{V}^{e}(s)(\mathcal{H} \psi_{V}^{e})
(s) \! - \! \dfrac{8 \mi \psi_{V}^{e}(s)}{\pi s} \right) \dfrac{\md s}{2 \pi
\mi} \! = \! 0.
\end{align*}
Substituting Equation~(3.1) into the above formula, one notes that
\begin{align*}
(\mathscr{F}^{e}(z))^{2} -& \, \int_{J_{e}} \dfrac{(16 \mi \psi_{V}^{e}(\xi)
(\mathcal{H} \psi_{V}^{e})(\xi) \! - \! \frac{8 \mi \psi_{V}^{e}(\xi)}{\pi
\xi})}{(\xi \! - \! z)} \, \dfrac{\md \xi}{2 \pi \mi} \! - \! \dfrac{1}{z}
\int_{J_{e}} \! \left(16 \mi \psi_{V}^{e}(\xi)(\mathcal{H} \psi_{V}^{e})(\xi)
\! - \! \dfrac{8 \mi \psi_{V}^{e}(\xi)}{\pi \xi} \right) \! \dfrac{\md \xi}{2
\pi \mi} \\
+& \, \dfrac{1}{\pi^{2}z^{2}} \! - \! \dfrac{1}{z^{2}} \int_{J_{e}} \xi \!
\left(16 \mi \psi_{V}^{e}(\xi)(\mathcal{H} \psi_{V}^{e})(\xi) \! - \! \dfrac{
8 \mi \psi_{V}^{e}(\xi)}{\pi \xi} \right) \! \dfrac{\md \xi}{2 \pi \mi} \! =
\! 0.
\end{align*}
Via Equation~(3.3), it follows that $16 \mi \psi_{V}^{e}(s)(\mathcal{H} \psi_{
V}^{e})(s) \! - \! \tfrac{8 \mi \psi_{V}^{e}(s)}{\pi s} \! = \! \tfrac{4 \mi
\psi_{V}^{e}(s) \widetilde{V}^{\prime}(s)}{\pi}$; substituting the latter
expression into the above equation, and re-arranging, one obtains,
\begin{equation}
(\mathscr{F}^{e}(z))^{2} \! - \! \dfrac{2}{\pi^{2}} \int_{J_{e}} \dfrac{
\widetilde{V}^{\prime}(\xi) \psi_{V}^{e}(\xi)}{\xi \! - \! z} \, \md \xi \! +
\! \dfrac{1}{\pi^{2}z^{2}} \! - \! \dfrac{2}{\pi^{2}z^{2}} \int_{J_{e}} \xi
\widetilde{V}^{\prime}(\xi) \psi_{V}^{e}(\xi) \, \md \xi \! - \! \dfrac{2}{
\pi^{2}z} \int_{J_{e}} \widetilde{V}^{\prime}(\xi) \psi_{V}^{e}(\xi) \, \md
\xi \! = \! 0.
\end{equation}
But
\begin{align*}
\dfrac{2}{\pi^{2}} \int_{J_{e}} \dfrac{\widetilde{V}^{\prime}(\xi) \psi_{V}^{e}
(\xi)}{\xi \! - \! z} \, \md \xi =& \, \dfrac{2}{\pi^{2}} \int_{J_{e}} \dfrac{
(\widetilde{V}^{\prime}(\xi) \! - \! \widetilde{V}^{\prime}(z)) \psi_{V}^{e}
(\xi)}{\xi \! - \! z} \, \md \xi \! + \! \dfrac{2}{\pi^{2}} \int_{J_{e}}
\dfrac{\widetilde{V}^{\prime}(z) \psi_{V}^{e}(\xi)}{\xi \! - \! z} \, \md \xi
\\
=& \, \dfrac{2}{\pi^{2}} \int_{J_{e}} \dfrac{(\widetilde{V}^{\prime}(\xi) \! -
\! \widetilde{V}^{\prime}(z)) \psi_{V}^{e}(\xi)}{\xi \! - \! z} \, \md \xi \!
+ \! \dfrac{\mi \widetilde{V}^{\prime}(z)}{\pi} \! \underbrace{\left(\dfrac{
2}{\pi \mi} \int_{J_{e}} \dfrac{\psi_{V}^{e}(\xi)}{\xi \! - \! z} \, \md \xi
\right)}_{= \, -\mathscr{F}^{e}(z)-(\mi \pi z)^{-1}} \\
=& \, \dfrac{2}{\pi^{2}} \int_{J_{e}} \dfrac{(\widetilde{V}^{\prime}(\xi) \!
- \! \widetilde{V}^{\prime}(z)) \psi_{V}^{e}(\xi)}{\xi \! - \! z} \, \md \xi
\! - \! \dfrac{\mi \widetilde{V}^{\prime}(z) \mathscr{F}^{e}(z)}{\pi} \! - \!
\dfrac{\widetilde{V}^{\prime}(z)}{\pi^{2}z}:
\end{align*}
substituting the above into Equation~(3.5), one arrives at, upon completing
the square and re-arr\-a\-n\-g\-i\-n\-g terms,
\begin{equation}
\left(\mathscr{F}^{e}(z) \! + \! \dfrac{\mi \widetilde{V}^{\prime}(z)}{2 \pi}
\right)^{2} \! + \! \dfrac{\mathfrak{q}_{V}^{e}(z)}{\pi^{2}} \! = \! 0,
\end{equation}
where
\begin{equation*}
\mathfrak{q}_{V}^{e}(z) \! := \! \left(\dfrac{\widetilde{V}^{\prime}(z)}{2}
\right)^{2} \! + \! \dfrac{\widetilde{V}^{\prime}(z)}{z} \! - \! 2 \int_{J_{
e}} \dfrac{(\widetilde{V}^{\prime}(\xi) \! - \! \widetilde{V}^{\prime}(z))
\psi_{V}^{e}(\xi)}{\xi \! - \! z} \, \md \xi \! + \! \dfrac{1}{z^{2}} \!
\left(1 \! - \! 2 \int_{J_{e}}(\xi \! + \! z) \widetilde{V}^{\prime}(\xi)
\psi_{V}^{e}(\xi) \, \md \xi \right).
\end{equation*}
(Equation~(3.6) above generalizes Equation~(3.5) for $q^{(0)}(x)$ in 
\cite{a58} for the case when $\widetilde{V} \colon \mathbb{R} \setminus 
\{0\} \! \to \! \mathbb{R}$ is real analytic; moreover, it is analogous to 
Equation~(1.37) of \cite{a56}.) Note that, since $\widetilde{V} \colon 
\mathbb{R} \setminus \{0\} \! \to \! \mathbb{R}$ satisfies 
conditions~(2.3)--(2.5), it follows {}from $\alpha^{l} \! - \! \beta^{l} \! 
= \! (\alpha \! - \! \beta)(\alpha^{l-1} \! + \! \alpha^{l-2} \beta \! + \! 
\cdots \! + \! \alpha \beta^{l-2} \! + \! \beta^{l-1})$, $l \! \in \! 
\mathbb{N}$, that $\mathfrak{q}_{V}^{e}(z)$ is real analytic on $J_{e}$ (and 
real analytic on $\mathbb{R} \setminus \{0\})$. For $x \! \in \! J_{e}$, set 
$z \! := \! x \! + \! \mi \varepsilon$, and consider the $\varepsilon \! 
\downarrow \! 0$ limit of Equation~(3.6): $\lim_{\varepsilon \downarrow 0}
(\mathscr{F}^{e}(x \! + \! \mi \varepsilon) \! + \! \tfrac{\mi \widetilde{
V}^{\prime}(x+\mi \varepsilon)}{2 \pi})^{2} \! = \! (\mathscr{F}^{e}_{+}(x) 
\! + \! \tfrac{\mi \widetilde{V}^{\prime}(x)}{2 \pi})^{2}$ (as $\widetilde{
V}$ is real analytic on $J_{e})$; recalling that $\mathscr{F}^{e}_{+}(x) \! 
= \! -\tfrac{1}{\pi \mi x} \! - \! 2 \mi (\mathcal{H} \psi_{V}^{e})(x) \! - 
\! 2 \psi_{V}^{e}(x)$, via Equation~(3.3), it follows that $\mathscr{F}^{e}_{
+}(x) \! = \! -\tfrac{\mi \widetilde{V}^{\prime}(x)}{2 \pi} \! - \! 2 \psi_{
V}^{e}(x)$ $\Rightarrow$ $(\mathscr{F}^{e}_{+}(x) \! + \! \tfrac{\mi 
\widetilde{V}^{\prime}(x)}{2 \pi})^{2} \! = \! (2 \psi_{V}^{e}(x))^{2}$, 
whence $(\psi_{V}^{e}(x))^{2} \! = \! -\mathfrak{q}_{V}^{e}(x)/(2 \pi)^{2}$, 
$x \! \in \! J_{e}$, whereupon, using the fact that (see above) $\psi_{V}^{e}
(x) \! \geqslant \! 0 \, \, \forall \, \, x \! \in \! \overline{J_{e}}$, it 
follows that $\mathfrak{q}_{V}^{e}(x) \! \leqslant \! 0$, $x \! \in \! J_{e}$; 
moreover, as a by-product, decomposing $\mathfrak{q}_{V}^{e}(x)$, for $x \! 
\in \! J_{e}$, into positive and negative parts, that is, $\mathfrak{q}_{V}^{
e}(x) \! = \! (\mathfrak{q}_{V}^{e}(x))^{+} \! - \! (\mathfrak{q}_{V}^{e}
(x))^{-}$, $x \! \in \! J_{e}$, where $(\mathfrak{q}_{V}^{e}(x))^{\pm} \! := 
\! \max \left\lbrace \pm \mathfrak{q}_{V}^{e}(x),0 \right\rbrace$ $(\geqslant 
\! 0)$, one learns {}from the above analysis that, for $x \! \in \! J_{e}$, 
$(\mathfrak{q}_{V}^{e}(x))^{+} \! \equiv \! 0$ and $\psi_{V}^{e}(x) \! = \! 
\tfrac{1}{2 \pi}((\mathfrak{q}_{V}^{e}(x))^{-})^{1/2}$; and, since $\int_{
J_{e}} \! \psi_{V}^{e}(s) \, \md s \! = \! 1$, it follows that $\tfrac{1}{2 
\pi} \int_{J_{e}}((\mathfrak{q}_{V}^{e}(s))^{-})^{1/2} \, \md s \! = \! 1$, 
which gives rise to the interesting fact that the function $(\mathfrak{q}_{
V}^{e}(x))^{-} \! \not\equiv \! 0$ on $J_{e}$. (Even though $(\mathfrak{q}_{
V}^{e}(x))^{-}$ depends on $\md \mu_{V}^{e}(x) \! = \! \psi_{V}^{e}(x) \, \md 
x$, and thus $\psi_{V}^{e}(x) \! = \! \tfrac{1}{2 \pi}((\mathfrak{q}_{V}^{e}
(x))^{-})^{1/2}$ is an implicit representation for $\psi_{V}^{e}$, it is 
still a useful relation which can be used to obtain additional, valuable 
information about $\psi_{V}^{e}$.) For $x \! \notin \! J_{e}$, set $z \! 
:= \! x \! + \! \mi \varepsilon$, and (again) study the $\varepsilon \! 
\downarrow \! 0$ limit of Equation~(3.6): in this case, $\lim_{\varepsilon 
\downarrow 0}(\mathscr{F}^{e}(x \! + \! \mi \varepsilon) \! + \! \tfrac{\mi 
\widetilde{V}^{\prime}(x+\mi \varepsilon)}{2 \pi})^{2} \! = \! (\mathscr{F}^{
e}_{+}(x) \! + \! \tfrac{\mi \widetilde{V}^{\prime}(x)}{2 \pi})^{2} \! = \! 
(\mathscr{F}^{e}(x) \! + \! \tfrac{\mi \widetilde{V}^{\prime}(x)}{2 \pi})^{
2}$; recalling that, for $x \! \notin \! J_{e}$, $\mathscr{F}^{e}(x) \! = \! 
-\tfrac{1}{\pi \mi}(\tfrac{1}{x} \! + \! 2 \int_{J_{e}} \tfrac{\psi_{V}^{e}
(s)}{s-x} \, \md s) \! = \! \tfrac{\mi}{\pi x} \! - \! 2 \mi (\mathcal{H} 
\psi_{V}^{e})(x)$, substituting the latter expression into Equation~(3.6), 
one arrives at $(\tfrac{1}{\pi x} \! - \! 2(\mathcal{H} \psi_{V}^{e})(x) \! 
+ \! \tfrac{\widetilde{V}^{\prime}(x)}{2 \pi})^{2} \! = \! \mathfrak{q}_{V}^{
e}(x)/\pi^{2}$, $x \! \notin \! J_{e}$ (since $\widetilde{V}^{\prime}$ is 
real analytic on $(\mathbb{R} \setminus \{0\}) \setminus J_{e}$, it follows 
that $\mathfrak{q}_{V}^{e}(x)$, too, is real analytic on $(\mathbb{R} 
\setminus \{0\}) \setminus J_{e}$, in which case, this latter relation merely 
states that, for $x \! = \! 0$, $+\infty \! = \! +\infty)$, whence 
$\mathfrak{q}_{V}^{e}(x) \! \geqslant \! 0 \, \, \forall \, \, x \! \notin 
\! J_{e}$.

Now, recalling that, on a compact subset of $\mathbb{R}$, an analytic function 
changes sign an at most countable number of times, it follows {}from the above 
argument, the fact that $\widetilde{V} \colon \mathbb{R} \setminus \{0\} 
\! \to \! \mathbb{R}$ satisfying conditions~(2.3)--(2.5) is regular (cf. 
Subsection~2.2), in particular, $\widetilde{V}$ is real analytic in the (open) 
neighbourhood $\mathbb{U} \! := \! \lbrace \mathstrut z \! \in \! \mathbb{C}; 
\, \inf_{q \in J_{e}} \vert z \! - \! q \vert \! < \! r \! \in \! (0,1) 
\rbrace \setminus \{0\}$, $\mu_{V}^{e}$ has compact support, and mimicking a 
part of the calculations subsumed in the proof of Theorem~1.38 in \cite{a56}, 
that $J_{e} \! := \! \operatorname{supp}(\mu_{V}^{e}) \! = \! \lbrace 
\mathstrut x \! \in \! \mathbb{R}; \, \mathfrak{q}_{V}^{e}(x) \! \leqslant \! 
0 \rbrace$ consists of the disjoint union of a finite number of bounded (real) 
intervals, with representation $J_{e} \! := \! \cup_{j=1}^{N+1}J_{j}^{e}$, 
where $J_{j}^{e} \! := \! [b_{j-1}^{e},a_{j}^{e}]$, with $N \! \in \! \mathbb{
N}$ and finite, $b_{0}^{e} \! := \! \min \lbrace J_{e} \rbrace \! \notin \! 
\lbrace -\infty,0 \rbrace$, $a_{N+1}^{e} \! := \! \max \lbrace J_{e} \rbrace 
\! \notin \! \lbrace 0,+\infty \rbrace$, and $-\infty \! < \! b_{0}^{e} \! < 
\! a_{1}^{e} \! < \! b_{1}^{e} \! < \! a_{2}^{e} \! < \! \cdots \! < \! b_{
N}^{e} \! < \! a_{N+1}^{e} \! < \! +\infty$. (One notes that $\widetilde{V}$ 
is real analytic in, say, the open neighbourhood $\widetilde{\mathbb{U}} \! 
:= \! \cup_{j=1}^{N+1} \widetilde{\mathbb{U}}_{j}$, where $\widetilde{\mathbb{
U}}_{j} \! := \! \lbrace \mathstrut z \! \in \! \mathbb{C}^{\ast}; \, \inf_{q 
\in J_{j}^{e}} \vert z \! - \! q \vert \! < \! r_{j} \! \in \! (0,1) \rbrace$, 
with $\widetilde{\mathbb{U}}_{i} \cap \widetilde{\mathbb{U}}_{j} \! = \! 
\varnothing$, $i \! \not= \! j \! = \! 1,\dotsc,N \! + \! 1$.) Furthermore, 
as a by-product of the above representation for $J_{e}$, it follows that, 
since $J_{i}^{e} \cap J_{j}^{e} \! = \! \varnothing$, $i \! \not= \! j \! = 
\! 1,\dotsc,N \! + \! 1$, $\text{meas}(J_{e}) \! = \! \sum_{j=1}^{N+1} \vert 
b_{j-1}^{e} \! - \! a_{j}^{e} \vert \! < \! +\infty$.

It remains, still, to determine the $2(N \! + \! 1)$ conditions satisfied by
the end-points of the support of the `even' equilibrium measure, $\lbrace b_{
j-1}^{e},a_{j}^{e} \rbrace_{j=1}^{N+1}$. Towards this end, one proceeds as
follows. {}From the formula for $\mathscr{F}^{e}(z)$ given in Equation~(3.2):
\begin{compactenum}
\item[(i)] for $\mu_{V}^{e} \! \in \! \mathcal{M}_{1}(\mathbb{R})$, in
particular, $\int_{\mathbb{R}} \md \mu_{V}^{e}(s) \! = \! 1$ and $\int_{
\mathbb{R}}s^{m} \, \md \mu_{V}^{e}(s) \! < \! \infty$, $m \! \in \! \mathbb{
N}$, $s \! \in \! J_{e}$ and $z \! \notin \! J_{e}$, with $\vert s/z \vert \!
\ll \! 1$ (e.g., $\vert z \vert \! \gg \! \max_{j=1,\dotsc,N+1} \lbrace \vert
b_{j-1}^{e} \! - \! a_{j}^{e} \vert \rbrace)$, via the expansion $\tfrac{1}{s
-z} \! = \! -\sum_{k=0}^{l} \tfrac{s^{k}}{z^{k+1}} \! + \! \tfrac{s^{l+1}}{z^{
l+1}(s-z)}$, $l \! \in \! \mathbb{Z}_{0}^{+}$, one gets that $\mathscr{F}^{e}
(z) \! =_{z \to \infty} \! \tfrac{1}{\pi \mi z} \! + \! \mathcal{O}(z^{-2})$;
\item[(ii)] for $\mu_{V}^{e} \! \in \! \mathcal{M}_{1}(\mathbb{R})$, in
particular, $\int_{\mathbb{R}}s^{-m} \, \md \mu_{V}^{e}(s) \! < \! \infty$,
$m \! \in \! \mathbb{N}$, $s \! \in \! J_{e}$ and $z \! \notin \! J_{e}$, with
$\vert z/s \vert \! \ll \! 1$ (e.g., $\vert z \vert \! \ll \! \min_{j=1,\dotsc,
N+1} \lbrace \vert b_{j-1}^{e} \! - \! a_{j}^{e} \vert \rbrace)$, via the
expansion $\tfrac{1}{z-s} \! = \! -\sum_{k=0}^{l} \tfrac{z^{k}}{s^{k+1}} \! +
\! \tfrac{z^{l+1}}{s^{l+1}(z-s)}$, $l \! \in \! \mathbb{Z}_{0}^{+}$, one gets
that $\mathscr{F}^{e}(z) \! =_{z \to 0} \! -\tfrac{1}{\pi \mi z} \! + \!
\mathcal{O}(1)$.
\end{compactenum}
Recalling, also, the formulae for $\mathscr{F}^{e}_{\pm}(z)$ given in
Equation~(3.4), one deduces that $\mathscr{F}^{e}_{+}(z) \! + \! \mathscr{F}^{
e}_{-}(z) \! = \! -\mi \widetilde{V}^{\prime}(z)/\pi$, $z \! \in \! J_{e}$,
and $\mathscr{F}^{e}_{+}(z) \! - \! \mathscr{F}^{e}_{-}(z) \! = \! 0$, $z \!
\notin \! J_{e}$; thus, one learns that $\mathscr{F}^{e} \colon \mathbb{C}
\setminus (J_{e} \cup \{0\}) \! \to \! \mathbb{C}$ solves the following
(scalar and homogeneous) RHP:
\begin{compactenum}
\item[(1)] $\mathscr{F}^{e}(z)$ is holomorphic (resp., meromorphic) for $z \!
\in \! \mathbb{C} \setminus (J_{e} \cup \{0\})$ (resp., $z \! \in \! \mathbb{
C} \setminus J_{e})$;
\item[(2)] $\mathscr{F}^{e}_{\pm}(z) \! := \! \lim_{\varepsilon \downarrow 0}
\mathscr{F}^{e}(z \! \pm \! \mi \varepsilon)$ satisfy the boundary condition
$\mathscr{F}^{e}_{+}(z) \! + \! \mathscr{F}^{e}_{-}(z) \! = \! -\mi
\widetilde{V}^{\prime}(z)/\pi$, $z \! \in \! J_{e}$, with $\mathscr{F}^{e}_{+}
(z) \! = \! \mathscr{F}^{e}_{-}(z) \! := \! \mathscr{F}^{e}(z)$ for $z \!
\notin \! J_{e}$;
\item[(3)] $\mathscr{F}^{e}(z) \! =_{\underset{z \in \mathbb{C} \setminus
\mathbb{R}}{z \to \infty}} \! \tfrac{1}{\pi \mi z} \! + \! \mathcal{O}(z^{-
2})$; and
\item[(4)] $\operatorname{Res}(\mathscr{F}^{e}(z);0) \! = \! -1/\pi \mi$.
\end{compactenum}
The solution of this RHP is (see, for example, \cite{a95})
\begin{equation*}
\mathscr{F}^{e}(z) \! = \! -\dfrac{1}{\pi \mi z} \! + \! (R_{e}(z))^{1/2}
\int_{J_{e}} \dfrac{(\frac{2}{\mi \pi s} \! + \! \frac{\widetilde{V}^{\prime}
(s)}{\mi \pi})}{(R_{e}(s))^{1/2}_{+}(s \! - \! z)} \, \dfrac{\md s}{2 \pi
\mi}, \quad z \! \in \! \mathbb{C} \setminus (J_{e} \cup \{0\}),
\end{equation*}
where $(R_{e}(z))^{1/2}$ is defined in the Lemma, with $(R_{e}(z))^{1/2}_{\pm}
\! := \! \lim_{\varepsilon \downarrow 0}(R_{e}(z \! \pm \! \mi \varepsilon))^{
1/2}$, and the branch of the square root is chosen so that $z^{-(N+1)}(R_{e}
(z))^{1/2} \! \sim_{\underset{z \in \mathbb{C}_{\pm}}{z \to \infty}} \! \pm
1$. (Note that $(R_{e}(z))^{1/2}$ is pure imaginary on $J_{e}$.) It follows
{}from the above integral representation for $\mathscr{F}^{e}(z)$ that, for
$s \! \in \! J_{e}$ and $z \! \notin \! J_{e}$, with $\vert s/z \vert \! \ll
\! 1$ (e.g., $\vert z \vert \! \gg \! \max_{j=1,\dotsc,N+1} \lbrace \vert b_{
j-1}^{e} \! - \! a_{j}^{e} \vert \rbrace)$, via the expansion $\tfrac{1}{s-z}
\! = \! -\sum_{k=0}^{l} \tfrac{s^{k}}{z^{k+1}} \! + \! \tfrac{s^{l+1}}{z^{l+1}
(s-z)}$, $l \! \in \! \mathbb{Z}_{0}^{+}$,
\begin{equation*}
\mathscr{F}^{e}(z) \! \underset{z \to \infty}{=} \! -\dfrac{1}{\mi \pi z} \! +
\! \dfrac{(z^{N+1} \! + \! \dotsb)}{z} \int_{J_{e}} \! \dfrac{(\frac{2 \mi}{
\pi s} \! + \! \frac{\mi \widetilde{V}^{\prime}(s)}{\pi})}{(R_{e}(s))^{1/2}_{
+}} \! \left(1 \! + \! \dfrac{s}{z} \! + \! \cdots \! + \! \dfrac{s^{N}}{z^{
N}} \! + \! \dfrac{s^{N+1}}{z^{N+1}} \! + \! \dotsb \right) \! \dfrac{\md
s}{2 \pi \mi}:
\end{equation*}
now, recalling {}from above that $\mathscr{F}^{e}(z) \! =_{z \to \infty} \!
\tfrac{1}{\pi \mi z} \! + \! \mathcal{O}(z^{-2})$, it follows that, upon
removing the secular (growing) terms,
\begin{equation*}
\int_{J_{e}} \! \left(\dfrac{2 \mi}{\pi s} \! + \! \dfrac{\mi \widetilde{V}^{
\prime}(s)}{\pi} \right) \! \dfrac{s^{j}}{(R_{e}(s))^{1/2}_{+}} \, \dfrac{\md
s}{2 \pi \mi} \! = \! 0, \quad j \! = \! 0,\dotsc,N
\end{equation*}
(which gives $N \! + \! 1$ (real) moment conditions), and, upon equating $z^{-
1}$ terms,
\begin{equation*}
\int_{J_{e}} \! \left(\dfrac{2 \mi}{\pi s} \! + \! \dfrac{\mi \widetilde{V}^{
\prime}(s)}{\pi} \right) \! \dfrac{s^{N+1}}{(R_{e}(s))^{1/2}_{+}} \, \dfrac{
\md s}{2 \pi \mi} \! = \! \dfrac{2}{\pi \mi};
\end{equation*}
it remains, therefore, to determine an additional $2(N \! + \! 1) \! - \! (N
\! + \! 1) \! - \! 1 \! = \! N$ (real) moment conditions. {}From the integral
representation for $\mathscr{F}^{e}(z)$, a residue calculus calculation shows
that
\begin{equation}
\mathscr{F}^{e}(z) \! = \! -\dfrac{\mi \widetilde{V}^{\prime}(z)}{2 \pi} \!
- \! \dfrac{(R_{e}(z))^{1/2}}{2} \oint_{C_{\mathrm{R}}^{e}} \! \dfrac{(\frac{
2 \mi}{\pi s} \! + \! \frac{\mi \widetilde{V}^{\prime}(s)}{\pi})}{(R_{e}(s))^{
1/2}(s \! - \! z)} \, \dfrac{\md s}{2 \pi \mi},
\end{equation}
where $C_{\mathrm{R}}^{e}$ $(\subset \mathbb{C}^{\ast})$ denotes the boundary
of any open doubly-connected annular region of the type $\lbrace \mathstrut
z^{\prime} \! \in \! \mathbb{C}; \, 0 \! < \! r \! < \! \vert z^{\prime} \vert
\! < \! R \! < \! +\infty \rbrace$, where the simple outer (resp., inner)
boundary $\lbrace \mathstrut z^{\prime} \! = \! R \me^{\mi \vartheta}, \, 0
\! \leqslant \! \vartheta \! \leqslant \! 2 \pi \rbrace$ (resp., $\lbrace
\mathstrut z^{\prime} \! = \! r \me^{\mi \vartheta}, \, 0 \! \leqslant \!
\vartheta \! \leqslant \! 2 \pi \rbrace)$ is traversed clockwise (resp.,
counter-clockwise), with the numbers $0 \! < \! r \! < \! R \! < \! +\infty$
chosen such that, for (any) non-real $z$ in the domain of analyticity of
$\widetilde{V}$ (that is, $\mathbb{C}^{\ast})$, $\mathrm{int}(C_{\mathrm{R}
}^{e}) \! \supset \! J_{e} \cup \{z\}$. Recall {}from Equation~(3.4) that, for
$z \! \in \! \mathbb{R} \setminus \overline{J_{e}}$ $(\supset \cup_{j=1}^{N}
(a_{j}^{e},b_{j}^{e}))$, $\mathscr{F}^{e}_{+}(z) \! = \! \mathscr{F}^{e}_{-}
(z) \! = \! -\tfrac{1}{\pi \mi}(\tfrac{1}{z} \! + \! 2 \int_{J_{e}} \tfrac{
\psi_{V}^{e}(s)}{s-z} \, \md s)$, whence $\mathscr{F}^{e}(z) \! + \! \tfrac{
1}{\pi \mi z} \! = \! -2 \mi (\mathcal{H} \psi_{V}^{e})(z)$; thus, using
Equation~(3.7), one arrives at
\begin{equation*}
(\mathcal{H} \psi_{V}^{e})(z) \! = \! \dfrac{\widetilde{V}^{\prime}(z)}{4 \pi}
\! + \! \dfrac{1}{2 \pi z} \! + \! \dfrac{\mi (R_{e}(z))^{1/2}}{2} \oint_{C_{
\mathrm{R}}^{e}} \! \dfrac{(\frac{1}{\pi \mi \xi} \! + \! \frac{\widetilde{
V}^{\prime}(\xi)}{2 \pi \mi})}{(R_{e}(\xi))^{1/2}(\xi \! - \! z)} \, \dfrac{
\md \xi}{2 \pi \mi}, \quad z \! \in \! \cup_{j=1}^{N}(a_{j}^{e},b_{j}^{e}).
\end{equation*}
A contour integration argument shows that
\begin{equation*}
\int_{a_{j}^{e}}^{b_{j}^{e}} \left((\mathcal{H} \psi_{V}^{e})(s) \! - \!
\dfrac{1}{2 \pi s} \! - \! \dfrac{\widetilde{V}^{\prime}(s)}{4 \pi} \right)
\md s \! = \! 0, \quad j \! = \! 1,\dotsc,N,
\end{equation*}
whence, using the above expression for $(\mathcal{H} \psi_{V}^{e})(z)$, $z \!
\in \! \cup_{j=1}^{N}(a_{j}^{e},b_{j}^{e})$, it follows that
\begin{equation}
\int_{a_{j}^{e}}^{b_{j}^{e}} \! \left(\dfrac{\mi (R_{e}(s))^{1/2}}{2} \oint_{
C_{\mathrm{R}}^{e}} \! \dfrac{(\frac{1}{\pi \mi \xi} \! + \! \frac{\widetilde{
V}^{\prime}(\xi)}{2 \pi \mi})}{(R_{e}(\xi))^{1/2}(\xi \! - \! s)} \, \dfrac{
\md \xi}{2 \pi \mi} \right) \! \md s \! = \! 0, \quad j \! = \! 1,\dotsc,N:
\end{equation}
now, `collapsing' the contour $C_{\mathrm{R}}^{e}$ down to $\mathbb{R}
\setminus \{0\}$ and using the Residue Theorem, one shows that
\begin{equation*}
\dfrac{\mi (R_{e}(z))^{1/2}}{2} \oint_{C_{\mathrm{R}}^{e}} \! \dfrac{(\frac{
1}{\pi \mi \xi} \! + \! \frac{\widetilde{V}^{\prime}(\xi)}{2 \pi \mi})}{(R_{e}
(\xi))^{1/2}(\xi \! - \! z)} \, \dfrac{\md \xi}{2 \pi \mi} \! = \! -\dfrac{1}{
2 \pi z} \! - \! \dfrac{\widetilde{V}^{\prime}(z)}{4 \pi} \! + \! \mi (R_{e}
(z))^{1/2} \int_{J_{e}} \! \dfrac{(\frac{1}{\pi \mi \xi} \! + \! \frac{
\widetilde{V}^{\prime}(\xi)}{2 \pi \mi})}{(R_{e}(\xi))^{1/2}_{+}(\xi \! - \!
z)} \, \dfrac{\md \xi}{2 \pi \mi};
\end{equation*}
substituting the latter relation into Equation~(3.8), one arrives at, after
straightforward integration and using the Fundamental Theorem of Calculus, for
$j \! = \! 1,\dotsc,N$,
\begin{equation*}
\int_{a_{j}^{e}}^{b_{j}^{e}} \left(\mi (R_{e}(s))^{1/2} \int_{J_{e}} \left(
\dfrac{\mi}{\pi \xi} \! + \! \dfrac{\mi \widetilde{V}^{\prime}(\xi)}{2 \pi}
\right) \! \dfrac{1}{(R_{e}(\xi))^{1/2}_{+}(\xi \! - \! s)} \, \dfrac{\md
\xi}{2 \pi \mi} \right) \md s \! = \! \dfrac{1}{2 \pi} \ln \! \left\vert
\dfrac{a_{j}^{e}}{b_{j}^{e}} \right\vert \! + \! \dfrac{1}{4 \pi} \! \left(
\widetilde{V}(a_{j}^{e}) \! - \! \widetilde{V}(b_{j}^{e}) \right),
\end{equation*}
which give the remaining $N$ moment conditions determining the end-points of
the support of the `even' equilibrium measure, $\lbrace b_{j-1}^{e},a_{j}^{e}
\rbrace_{j=1}^{N+1}$. Since $J_{e} \not\supseteq \{0,\pm \infty\}$ and
$\widetilde{V}$ is real analytic on $J_{e}$,
\begin{equation*}
(R_{e}(s))^{1/2} \underset{s \downarrow b_{j-1}^{e}}{=} \! \mathcal{O} \!
\left((s \! - \! b_{j-1}^{e})^{1/2} \right) \qquad \text{and} \qquad (R_{e}
(s))^{1/2} \underset{s \uparrow a_{j}^{e}}{=} \! \mathcal{O} \! \left(
(a_{j}^{e} \! - \! s)^{1/2} \right), \quad j \! = \! 1,\dotsc,N \! + \! 1,
\end{equation*}
which shows that all the integrals above constituting the $n$-dependent system
of $2(N \! + \! 1)$ moment conditions for the end-points of the support of
$\mu_{V}^{e}$ have removable singularities at $b_{j-1}^{e},a_{j}^{e}$, $j \!
= \! 1,\dotsc,N \! + \! 1$.

Recall {}from Equation~(3.4) that, for $z \! \in \! J_{e}$, $\mathscr{F}^{e}_{
\pm}(z) \! = \! -\tfrac{1}{\pi \mi z} \! - \! 2 \mi (\mathcal{H} \psi_{V}^{e})
(z) \! \mp \! 2 \psi_{V}^{e}(z)$: using the fact that, {}from Equation~(3.3),
for $z \! \in \! J_{e}$, $(\mathcal{H} \psi_{V}^{e})(z) \! = \! \tfrac{1}{2
\pi z} \! + \! \tfrac{\widetilde{V}^{\prime}(z)}{4 \pi}$, it follows that
\begin{equation*}
\mathscr{F}^{e}_{\pm}(z) \! = \! \dfrac{\widetilde{V}^{\prime}(z)}{2 \pi \mi}
\! \mp \! 2 \psi_{V}^{e}(z), \quad z \! \in \! J_{e}.
\end{equation*}
{}From Equation~(3.7), it follows that
\begin{equation*}
\mathscr{F}^{e}_{\pm}(z) \! = \! \dfrac{\widetilde{V}^{\prime}(z)}{2 \pi
\mi} \! + \! \dfrac{(R_{e}(z))^{1/2}_{\pm}}{2} \oint_{C_{\mathrm{R}}^{e}} \!
\dfrac{(\frac{2}{\pi \mi s} \! + \! \frac{\widetilde{V}^{\prime}(s)}{\mi \pi}
)}{(R_{e}(s))^{1/2}(s \! - \! z)} \, \dfrac{\md s}{2 \pi \mi};
\end{equation*}
thus, equating the above two expressions for $\mathscr{F}^{e}_{\pm}(z)$, one
arrives at $\psi_{V}^{e}(x) \! = \! \tfrac{1}{2 \pi \mi}(R_{e}(x))^{1/2}_{+}
h_{V}^{e}(x) \pmb{1}_{J_{e}}(x)$, where $h_{V}^{e}(z)$ is defined in the
Lemma, and $\pmb{1}_{J_{e}}(x)$ is the characteristic function of the set
$J_{e}$, which gives rise to the formula for the density of the `even'
equilibrium measure, $\md \mu_{V}^{e}(x) \! = \! \psi_{V}^{e}(x) \, \md x$
(the integral representation for $h_{V}^{e}(z)$ shows that it is analytic in
some open subset of $\mathbb{C}^{\ast}$ containing $J_{e})$. Now, recalling
that $\widetilde{V} \colon \mathbb{R} \setminus \{0\} \! \to \! \mathbb{R}$
satisfying conditions~(2.3)--(2.5) is regular, and that, for $s \! \in \!
J_{e}$ (resp., $s \! \in \! \overline{J_{e}})$, $\psi_{V}^{e}(s) \! > \! 0$
(resp., $\psi_{V}^{e}(s) \! \geqslant \! 0)$ and $(R_{e}(s))^{1/2}_{+} \! = \!
\mi (\vert R_{e}(s) \vert)^{1/2} \! \in \! \mi \mathbb{R}_{\pm}$ (resp.,
$(R_{e}(s))^{1/2}_{+} \! = \! \mi (\vert R_{e}(s) \vert)^{1/2} \! \in \! \mi
\mathbb{R})$, it follows {}from the formula $\psi_{V}^{e}(s) \! = \! \tfrac{
1}{2 \pi \mi}(R_{e}(s))^{1/2}_{+}h_{V}^{e}(s) \pmb{1}_{J_{e}}(s)$ and the
regularity assumption, namely, $h_{V}^{e}(z) \! \not\equiv \! 0$ for $z \!
\in \! \overline{J_{e}}$, that $(\vert R_{e}(s) \vert)^{1/2}h_{V}^{e}(s) \! >
\! 0$, $s \! \in \! J_{e}$ (resp., $(\vert R_{e}(s) \vert)^{1/2}h_{V}^{e}(s)
\! \geqslant \! 0$, $s \! \in \! \overline{J_{e}})$.

Finally, it will be shown that, if $\overline{J_{e}} \! := \! \cup_{j=1}^{N
+1}[b_{j-1}^{e},a_{j}^{e}]$, the end-points of the support of the `even' 
equilibrium measure, which satisfy the $n$-dependent system of $2(N \! + \! 
1)$ moment conditions stated in the Lemma, are (real) analytic functions of 
$z_{o}$, thus proving the (local) solvability of the $n$-dependent $2(N \! 
+ \! 1)$ moment conditions. Towards this end, one follows closely the idea 
of the proof of Theorem~1.3~(iii) in \cite{a92} (see, also, Section~8 of 
\cite{a56}, and \cite{a96}). Recall {}from Subsection~2.2 that 
$\widetilde{V}(z) \! := \! z_{o}V(z)$, where $z_{o} \colon \mathbb{N} \times 
\mathbb{N} \! \to \! \mathbb{R}_{+}$, $(n,\mathscr{N}) \! \mapsto \! z_{o} \! 
:= \! \mathscr{N}/n$, and, in the double-scaling limit as $\mathscr{N},n \! 
\to \! \infty$, $z_{o} \! = \! 1 \! + \! o(1)$. Furthermore, in the analysis 
above, it was shown that the end-points of the support of the `even' 
equilibrium measure were the simple zeros/roots of the function $\mathfrak{
q}_{V}^{e}(z)$, that is (up to re-arrangement), $\lbrace b_{0}^{e},b_{1}^{e},
\dotsc,b_{N}^{e},a_{1}^{e},a_{2}^{e},\dotsc,a_{N+1}^{e} \rbrace \! = \! 
\lbrace \mathstrut x \! \in \! \mathbb{R}; \, \mathfrak{q}_{V}^{e}(x) \! = 
\! 0 \rbrace$ (these are the only roots for the regular case studied in this 
work). The function $\mathfrak{q}_{V}^{e}(x) \! \in \! \mathbb{R}(x)$ (the 
algebra of rational functions in $x$ with coefficients in $\mathbb{R})$ is 
real rational on $\mathbb{R}$ and real analytic on $\mathbb{R} \setminus 
\{0\}$, it has analytic continuation to $\lbrace \mathstrut z \! \in \! 
\mathbb{C}; \, \inf_{p \in \mathbb{R}} \vert z \! - \! p \vert \! < \! r \! 
\in \! (0,1) \rbrace \setminus \{0\}$ (independent of $z_{o})$, and depends 
continuously on $z_{o}$; thus, its simple zeros/roots, that is, $b_{k-1}^{e} 
\! = \! b_{k-1}^{e}(z_{o})$ and $a_{k}^{e} \! = \! a_{k}^{e}(z_{o})$, $k \! 
= \! 1,\dotsc,N \! + \! 1$, are continuous functions of $z_{o}$.

Write the large-$z$ (e.g., $\vert z \vert \! \gg \! \max_{j=1,\dotsc,N+1} 
\lbrace \vert b_{j-1}^{e} \! - \! a_{j}^{e} \vert \rbrace)$ asymptotic 
expansion for $\mathscr{F}^{e}(z)$ given above as follows:
\begin{equation*}
\mathscr{F}^{e}(z) \underset{z \to \infty}{=} -\dfrac{1}{\mi \pi z} \! - \!
\dfrac{(R_{e}(z))^{1/2}}{2 \pi \mi z} \sum_{j=0}^{\infty} \mathcal{T}_{j}^{e}
z^{-j},
\end{equation*}
where
\begin{equation*}
\mathcal{T}_{j}^{e} \! := \! \int_{J_{e}} \! \left(\dfrac{2}{\mi \pi s} \! +
\! \dfrac{\widetilde{V}^{\prime}(s)}{\mi \pi} \right) \! \dfrac{s^{j}}{(R_{e}
(s))^{1/2}_{+}} \, \md s, \quad j \! \in \! \mathbb{Z}_{0}^{+}.
\end{equation*}
Set
\begin{equation*}
\mathcal{N}_{j}^{e} \! := \! \int_{a_{j}^{e}}^{b_{j}^{e}} \! \left(
(\mathcal{H} \psi_{V}^{e})(s) \! - \! \dfrac{1}{2 \pi s} \! - \! \dfrac{
\widetilde{V}^{\prime}(s)}{4 \pi} \right) \md s, \quad j \! = \! 1,\dotsc,N.
\end{equation*}
The $(n$-dependent) $2(N \! + \! 1)$ moment conditions are, thus, equivalent 
to $\mathcal{T}_{j}^{e} \! = \! 0$, $j \! = \! 0,\dotsc,N$, $\mathcal{T}_{N+
1}^{e} \! = \! -4$, and $\mathcal{N}_{j}^{e} \! = \! 0$, $j \! = \! 1,\dotsc,
N$. It will first be shown that, for regular $\widetilde{V} \colon \mathbb{R} 
\setminus \{0\} \! \to \! \mathbb{R}$ satisfying conditions~(2.3)--(2.5), the 
Jacobian of the transformation $\lbrace b_{0}^{e}(z_{o}),\dotsc,b_{N}^{e}(z_{
o}),a_{1}^{e}(z_{o}),\dotsc,a_{N+1}^{e}(z_{o}) \rbrace \! \mapsto \! \lbrace 
\mathcal{T}_{0}^{e},\dotsc,\mathcal{T}_{N+1}^{e},\mathcal{N}_{1}^{e},\dotsc,
\linebreak[4]
\mathcal{N}_{N}^{e} \rbrace$, that is, $\operatorname{Jac}(\mathcal{T}_{0}^{
e},\dotsc,\mathcal{T}_{N+1}^{e},\mathcal{N}_{1}^{e},\dotsc,\mathcal{N}_{N}^{
e}) \! := \! \tfrac{\partial (\mathcal{T}_{0}^{e},\dotsc,\mathcal{T}_{N+1}^{
e},\mathcal{N}_{1}^{e},\dotsc,\mathcal{N}_{N}^{e})}{\partial (b_{0}^{e},
\dotsc,b_{N}^{e},a_{1}^{e},\dotsc,a_{N+1}^{e})}$, is non-zero whenever $b_{j
-1}^{e} \! = \! b_{j-1}^{e}(z_{o})$ and $a_{j}^{e} \! = \! a_{j}^{e}(z_{o})$, 
$j \! = \! 1,\dotsc,N \! + \! 1$, are chosen so that $\overline{J_{e}} \! = 
\! \cup_{j=1}^{N+1}[b_{j-1}^{e},a_{j}^{e}]$. Using the equation $(\mathcal{H} 
\psi_{V}^{e})(z) \! = \! \tfrac{\mi}{2} \mathscr{F}^{e}(z) \! + \! \tfrac{1}{
2 \pi z}$ (cf. Equation~(3.2)), one follows the analysis on pp.~778--779 of
\cite{a92} to show that, for $k \! = \! 1,\dotsc,N \! + \! 1$:
\begin{gather}
\dfrac{\partial \mathcal{T}_{j}^{e}}{\partial b_{k-1}^{e}} \! = \! b_{k-1}^{e}
\dfrac{\partial \mathcal{T}_{j-1}^{e}}{\partial b_{k-1}^{e}} \! + \! \dfrac{
1}{2} \mathcal{T}_{j-1}^{e}, \quad j \! \in \! \mathbb{N}, \tag{T1} \\
\dfrac{\partial \mathcal{T}_{j}^{e}}{\partial a_{k}^{e}} \! = \! a_{k}^{e}
\dfrac{\partial \mathcal{T}_{j-1}^{e}}{\partial a_{k}^{e}} \! + \! \dfrac{1}{
2} \mathcal{T}_{j-1}^{e}, \quad j \! \in \! \mathbb{N}, \tag{T2} \\
\dfrac{\partial \mathscr{F}^{e}(z)}{\partial b_{k-1}^{e}} \! = \! -\dfrac{1}{2
\pi \mi} \! \left(\dfrac{\partial \mathcal{T}_{0}^{e}}{\partial b_{k-1}^{e}}
\right) \! \dfrac{(R_{e}(z))^{1/2}}{z \! - \! b_{k-1}^{e}}, \quad z \! \in \!
\mathbb{C} \setminus (\overline{J_{e}} \cup \{0\}), \tag{F1} \\
\dfrac{\partial \mathscr{F}^{e}(z)}{\partial a_{k}^{e}} \! = \! -\dfrac{1}{2
\pi \mi} \! \left(\dfrac{\partial \mathcal{T}_{0}^{e}}{\partial a_{k}^{e}}
\right) \! \dfrac{(R_{e}(z))^{1/2}}{z \! - \! a_{k}^{e}}, \quad z \! \in \!
\mathbb{C} \setminus (\overline{J_{e}} \cup \{0\}), \tag{F2} \\
\dfrac{\partial \mathcal{N}_{j}^{e}}{\partial b_{k-1}^{e}} \! = \! -\dfrac{1}{
4 \pi} \! \left(\dfrac{\partial \mathcal{T}_{0}^{e}}{\partial b_{k-1}^{e}}
\right) \! \int_{a_{j}^{e}}^{b_{j}^{e}} \dfrac{(R_{e}(s))^{1/2}}{s \! - \!
b_{k-1}^{e}} \, \md s, \quad j \! = \! 1,\dotsc,N, \tag{N1} \\
\dfrac{\partial \mathcal{N}_{j}^{e}}{\partial a_{k}^{e}} \! = \! -\dfrac{1}{4
\pi} \! \left(\dfrac{\partial \mathcal{T}_{0}^{e}}{\partial a_{k}^{e}} \right)
\! \int_{a_{j}^{e}}^{b_{j}^{e}} \dfrac{(R_{e}(s))^{1/2}}{s \! - \! a_{k}^{e}}
\, \md s, \quad j \! = \! 1,\dotsc,N; \tag{N2}
\end{gather}
furthermore, if one evaluates Equations~(T1) and~(T2) on the solution of the
$n$-dependent system of $2(N \! + \! 1)$ moment conditions, that is, $\mathcal{
T}_{j}^{e} \! = \! 0$, $j \! = \! 0,\dotsc,N$, $\mathcal{T}_{N+1}^{e} \! = \!
-4$, and $\mathcal{N}_{i}^{e} \! = \! 0$, $i \! = \! 1,\dotsc,N$, one arrives
at
\begin{equation}
\dfrac{\partial \mathcal{T}_{j}^{e}}{\partial b_{k-1}^{e}} \! = \! (b_{k-1}^{
e})^{j} \dfrac{\partial \mathcal{T}_{0}^{e}}{\partial b_{k-1}^{e}}, \qquad
\quad \dfrac{\partial \mathcal{T}_{j}^{e}}{\partial a_{k}^{e}} \! = \! (a_{
k}^{e})^{j} \dfrac{\partial \mathcal{T}_{0}^{e}}{\partial a_{k}^{e}}, \quad j
\! = \! 0,\dotsc,N \! + \! 1. \tag{S1}
\end{equation}
Via Equations~(N1), (N2), and~(S1), one now computes the Jacobian of the
transformation $\lbrace b_{0}^{e}(z_{o}),\dotsc,\linebreak[4]
b_{N}^{e}(z_{o}),a_{1}^{e}(z_{o}),\dotsc,a_{N+1}^{e}(z_{o}) \rbrace \! \mapsto
\! \lbrace \mathcal{T}_{0}^{e},\dotsc,\mathcal{T}_{N+1}^{e},\mathcal{N}_{1}^{
e},\dotsc,\mathcal{N}_{N}^{e} \rbrace$ on the solution of the $n$-dependent
system of $2(N \! + \! 1)$ moment conditions:
\begin{align*}
&\operatorname{Jac}(\mathcal{T}_{0}^{e},\dotsc,\mathcal{T}_{N+1}^{e},
\mathcal{N}_{1}^{e},\dotsc,\mathcal{N}_{N}^{e}) \! := \! \dfrac{\partial
(\mathcal{T}_{0}^{e},\dotsc,\mathcal{T}_{N+1}^{e},\mathcal{N}_{1}^{e},\dotsc,
\mathcal{N}_{N}^{e})}{\partial (b_{0}^{e},\dotsc,b_{N}^{e},a_{1}^{e},\dotsc,
a_{N+1}^{e})} \\
&=
\begin{vmatrix}
\frac{\partial \mathcal{T}_{0}^{e}}{\partial b_{0}^{e}} & \frac{\partial
\mathcal{T}_{0}^{e}}{\partial b_{1}^{e}} & \cdots & \frac{\partial \mathcal{
T}_{0}^{e}}{\partial b_{N}^{e}} & \frac{\partial \mathcal{T}_{0}^{e}}{\partial
a_{1}^{e}} & \frac{\partial \mathcal{T}_{0}^{e}}{\partial a_{2}^{e}} & \cdots
& \frac{\partial \mathcal{T}_{0}^{e}}{\partial a_{N+1}^{e}} \\
\frac{\partial \mathcal{T}_{1}^{e}}{\partial b_{0}^{e}} & \frac{\partial
\mathcal{T}_{1}^{e}}{\partial b_{1}^{e}} & \cdots & \frac{\partial \mathcal{
T}_{1}^{e}}{\partial b_{N}^{e}} & \frac{\partial \mathcal{T}_{1}^{e}}{\partial
a_{1}^{e}} & \frac{\partial \mathcal{T}_{1}^{e}}{\partial a_{2}^{e}} & \cdots
& \frac{\partial \mathcal{T}_{1}^{e}}{\partial a_{N+1}^{e}} \\
\vdots & \vdots & \ddots & \vdots & \vdots & \vdots & \ddots & \vdots \\
\frac{\partial \mathcal{T}_{N+1}^{e}}{\partial b_{0}^{e}} & \frac{\partial
\mathcal{T}_{N+1}^{e}}{\partial b_{1}^{e}} & \cdots & \frac{\partial \mathcal{
T}_{N+1}^{e}}{\partial b_{N}^{e}} & \frac{\partial \mathcal{T}_{N+1}^{e}}{
\partial a_{1}^{e}} & \frac{\partial \mathcal{T}_{N+1}^{e}}{\partial a_{2}^{
e}} & \cdots & \frac{\partial \mathcal{T}_{N+1}^{e}}{\partial a_{N+1}^{e}} \\
\frac{\partial \mathcal{N}_{1}^{e}}{\partial b_{0}^{e}} & \frac{\partial
\mathcal{N}_{1}^{e}}{\partial b_{1}^{e}} & \cdots & \frac{\partial \mathcal{
N}_{1}^{e}}{\partial b_{N}^{e}} & \frac{\partial \mathcal{N}_{1}^{e}}{\partial
a_{1}^{e}} & \frac{\partial \mathcal{N}_{1}^{e}}{\partial a_{2}^{e}} & \cdots
& \frac{\partial \mathcal{N}_{1}^{e}}{\partial a_{N+1}^{e}} \\
\frac{\partial \mathcal{N}_{2}^{e}}{\partial b_{0}^{e}} & \frac{\partial
\mathcal{N}_{2}^{e}}{\partial b_{1}^{e}} & \cdots & \frac{\partial \mathcal{
N}_{2}^{e}}{\partial b_{N}^{e}} & \frac{\partial \mathcal{N}_{2}^{e}}{\partial
a_{1}^{e}} & \frac{\partial \mathcal{N}_{2}^{e}}{\partial a_{2}^{e}} & \cdots
& \frac{\partial \mathcal{N}_{2}^{e}}{\partial a_{N+1}^{e}} \\
\vdots & \vdots & \ddots & \vdots & \vdots & \vdots & \ddots & \vdots \\
\frac{\partial \mathcal{N}_{N}^{e}}{\partial b_{0}^{e}} & \frac{\partial
\mathcal{N}_{N}^{e}}{\partial b_{1}^{e}} & \cdots & \frac{\partial \mathcal{
N}_{N}^{e}}{\partial b_{N}^{e}} & \frac{\partial \mathcal{N}_{N}^{e}}{\partial
a_{1}^{e}} & \frac{\partial \mathcal{N}_{N}^{e}}{\partial a_{2}^{e}} & \cdots
& \frac{\partial \mathcal{N}_{N}^{e}}{\partial a_{N+1}^{e}}
\end{vmatrix} \\
&= \dfrac{(-1)^{N}}{(4 \pi)^{N}} \! \left(\prod_{k=1}^{N+1} \dfrac{\partial
\mathcal{T}_{0}^{e}}{\partial b_{k-1}^{e}} \dfrac{\partial \mathcal{T}_{0}^{
e}}{\partial a_{k}^{e}} \right) \! \left(\prod_{j=1}^{N} \int_{a_{j}^{e}}^{
b_{j}^{e}}(R_{e}(s_{j}))^{1/2} \, \md s_{j} \right) \! \Delta_{d}^{e},
\end{align*}
where
\begin{equation*}
\Delta_{d}^{e} \! := \!
\begin{vmatrix}
1 & 1 & \cdots & 1 & 1 & 1 & \cdots & 1 \\
b_{0}^{e} & b_{1}^{e} & \cdots & b_{N}^{e} & a_{1}^{e} & a_{2}^{e} & \cdots &
a_{N+1}^{e} \\
\vdots & \vdots & \ddots & \vdots & \vdots & \vdots & \ddots & \vdots \\
(b_{0}^{e})^{N+1} & (b_{1}^{e})^{N+1} & \cdots & (b_{N}^{e})^{N+1} & (a_{1}^{
e})^{N+1} & (a_{2}^{e})^{N+1} & \cdots & (a_{N+1}^{e})^{N+1} \\
\frac{1}{s_{1}-b_{0}^{e}} & \frac{1}{s_{1}-b_{1}^{e}} & \cdots & \frac{1}{s_{1}
-b_{N}^{e}} & \frac{1}{s_{1}-a_{1}^{e}} & \frac{1}{s_{1}-a_{2}^{e}} & \cdots &
\frac{1}{s_{1}-a_{N+1}^{e}} \\
\frac{1}{s_{2}-b_{0}^{e}} & \frac{1}{s_{2}-b_{1}^{e}} & \cdots & \frac{1}{s_{2}
-b_{N}^{e}} & \frac{1}{s_{2}-a_{1}^{e}} & \frac{1}{s_{2}-a_{2}^{e}} & \cdots &
\frac{1}{s_{2}-a_{N+1}^{e}} \\
\vdots & \vdots & \ddots & \vdots & \vdots & \vdots & \ddots & \vdots \\
\frac{1}{s_{N}-b_{0}^{e}} & \frac{1}{s_{N}-b_{1}^{e}} & \cdots & \frac{1}{s_{N}
-b_{N}^{e}} & \frac{1}{s_{N}-a_{1}^{e}} & \frac{1}{s_{N}-a_{2}^{e}} & \cdots &
\frac{1}{s_{N}-a_{N+1}^{e}}
\end{vmatrix}.
\end{equation*}
The above determinant, that is, $\Delta_{d}^{e}$, has been calculated on 
pg.~780 of \cite{a92} (see, also, Section~5.3, Equations~(5.148) and~(5.149) 
of \cite{a62}), namely,
\begin{equation*}
\Delta_{d}^{e} \! = \! \dfrac{\left(\prod_{j=1}^{N+1} \prod_{k=1}^{N+1}(b_{k-
1}^{e} \! - \! a_{j}^{e}) \right) \! \left(\prod_{\substack{j,k=1\\j<k}}^{N+
1}(a_{k}^{e} \! - \! a_{j}^{e})(b_{k-1}^{e} \! - \! b_{j-1}^{e}) \right) \!
\left(\prod_{\substack{j,k=1\\j<k}}^{N}(s_{k} \! - \! s_{j}) \right)}{(-1)^{
N} \prod_{j=1}^{N} \prod_{k=1}^{N+1}(s_{j} \! - \! a_{k}^{e})(s_{j} \! - \!
b_{k-1}^{e})};
\end{equation*}
but, for $-\infty \! < \! b_{0}^{e} \! < \! a_{1}^{e} \! < \! s_{1} \! < \!
b_{1}^{e} \! < \! a_{2}^{e} \! < \! s_{2} \! < \! b_{2}^{e} \! < \! \dotsb \!
< \! b_{N-1}^{e} \! < \! a_{N}^{e} \! < \! s_{N} \! < \! b_{N}^{e} \! < \!
a_{N+1}^{e} \! < \! +\infty$, $\Delta_{d}^{e} \! \not= \! 0$ (which means that
it is of a fixed sign), and $\int_{a_{j}^{e}}^{b_{j}^{e}}(R_{e}(s_{j}))^{1/2}
\, \md s_{j} \! > \! 0$, $j \! = \! 1,\dotsc,N$, whence
\begin{equation*}
\left(\prod_{j=1}^{N} \int_{a_{j}^{e}}^{b_{j}^{e}}(R_{e}(s_{j}))^{1/2} \, \md
s_{j} \right) \! \Delta_{d}^{e} \! \not= \! 0.
\end{equation*}
It remains to show that $\partial \mathcal{T}_{0}^{e}/\partial b_{k-1}^{e}$
and $\partial \mathcal{T}_{0}^{e}/\partial a_{k}^{e}$, $k \! = \! 1,\dotsc,N
\! + \! 1$, too, are non-zero; for this purpose, one exploits the fact that
$\mathcal{T}_{0}^{e} \! = \! (\mi \pi)^{-1} \int_{J_{e}}(2s^{-1} \! + \!
\widetilde{V}^{\prime}(s))(R_{e}(s))^{-1/2}_{+} \, \md s$ is independent of
$z$. It follows {}from Equation~(3.7), the integral representation for $h_{
V}^{e}(z)$ given in the Lemma, and Equations~(F1) and~(F2) that
\begin{gather*}
\dfrac{(z \! - \! b_{k-1}^{e})}{\sqrt{\smash[b]{R_{e}(z)}}} \dfrac{\partial
\mathscr{F}^{e}(z)}{\partial b_{k-1}^{e}} \! = \! -\dfrac{1}{\mi \pi} \!
\left((z \! - \! b_{k-1}^{e}) \dfrac{\partial h_{V}^{e}(z)}{\partial b_{k-1}^{
e}} \! - \! \dfrac{1}{2}h_{V}^{e}(z) \right), \quad k \! = \! 1,\dotsc,N \! +
\! 1, \\
\dfrac{(z \! - \! a_{k}^{e})}{\sqrt{\smash[b]{R_{e}(z)}}} \dfrac{\partial
\mathscr{F}^{e}(z)}{\partial a_{k}^{e}} \! = \! -\dfrac{1}{\mi \pi} \! \left(
(z \! - \! a_{k}^{e}) \dfrac{\partial h_{V}^{e}(z)}{\partial a_{k}^{e}} \! -
\! \dfrac{1}{2}h_{V}^{e}(z) \right), \quad k \! = \! 1,\dotsc,N \! + \! 1:
\end{gather*}
using, now, the $z$-independence of $\mathcal{T}_{0}^{e}$, and the fact that,
for the case of regular $\widetilde{V} \colon \mathbb{R} \setminus \{0\} \!
\to \! \mathbb{R}$ satisfying conditions~(2.3)--(2.5), $h_{V}^{e}(b_{j-1}^{
e}),h_{V}^{e}(a_{j}^{e}) \! \not= \! 0$, $j \! = \! 1,\dotsc,N \! + \! 1$,
one shows that
\begin{gather*}
\left. \dfrac{(z \! - \! b_{k-1}^{e})}{\sqrt{\smash[b]{R_{e}(z)}}} \dfrac{
\partial \mathscr{F}^{e}(z)}{\partial b_{k-1}^{e}} \right\vert_{z=b_{k-1}^{
e}} \! = \! \dfrac{1}{2 \pi \mi}h_{V}^{e}(b_{k-1}^{e}) \! \not= \! 0, \quad
k \! = \! 1,\dotsc,N \! + \! 1, \\
\left. \dfrac{(z \! - \! a_{k}^{e})}{\sqrt{\smash[b]{R_{e}(z)}}} \dfrac{
\partial \mathscr{F}^{e}(z)}{\partial a_{k}^{e}} \right\vert_{z=a_{k}^{e}} \!
= \! \dfrac{1}{2 \pi \mi}h_{V}^{e}(a_{k}^{e}) \! \not= \! 0, \quad k \! = \!
1,\dotsc,N \! + \! 1;
\end{gather*}
thus, via Equations~(F1) and~(F2), one arrives at
\begin{equation*}
\dfrac{\partial \mathcal{T}_{0}^{e}}{\partial b_{k-1}^{e}} \! = \! -h_{V}^{e}
(b_{k-1}^{e}) \! \not= \! 0 \qquad \quad \text{and} \qquad \quad \dfrac{
\partial \mathcal{T}_{0}^{e}}{\partial a_{k}^{e}} \! = \! -h_{V}^{e}(a_{k}^{
e}) \! \not= \! 0, \quad k \! = \! 1,\dotsc,N \! + \! 1,
\end{equation*}
whence
\begin{equation*}
\prod_{k=1}^{N+1} \dfrac{\partial \mathcal{T}_{0}^{e}}{\partial b_{k-1}^{e}}
\dfrac{\partial \mathcal{T}_{0}^{e}}{\partial a_{k}^{e}} \! = \! \prod_{k=1}^{
N+1}h_{V}^{e}(b_{k-1}^{e})h_{V}^{e}(a_{k}^{e}) \! \not= \! 0.
\end{equation*}
Hence, $\operatorname{Jac}(\mathcal{T}_{0}^{e},\dotsc,\mathcal{T}_{N+1}^{e},
\mathcal{N}_{1}^{e},\dotsc,\mathcal{N}_{N}^{e}) \! \not= \! 0$.

It remains, still, to show that $\mathcal{T}_{j}^{e}$, $j \! = \! 0,\dotsc,
N \! + \! 1$, and $\mathcal{N}_{i}^{e}$, $i \! = \! 1,\dotsc,N$, are (real)
analytic functions of $\lbrace b_{j-1}^{e},a_{j}^{e} \rbrace_{j=1}^{N+1}$.
{}From the definition of $\mathcal{T}_{j}^{e}$, $j \! \in \! \mathbb{Z}_{0}^{
+}$, above, using the fact that they are independent of $z$, thus giving rise
to zero residue contributions, a straightforward residue calculus calculation
shows that, equivalently,
\begin{equation*}
\mathcal{T}_{j}^{e} \! = \! \dfrac{1}{2} \oint_{C_{\mathrm{R}}^{e}} \! \left(
\dfrac{2}{\mi \pi s} \! + \! \dfrac{\widetilde{V}^{\prime}(s)}{\mi \pi}
\right) \! \dfrac{s^{j}}{(R_{e}(s))^{1/2}} \, \md s, \quad j \! \in \!
\mathbb{Z}_{0}^{+},
\end{equation*}
where (the closed contour) $C_{\mathrm{R}}^{e}$ has been defined above: the
only factor depending on $\lbrace b_{k-1}^{e},a_{k}^{e} \rbrace_{k=1}^{N
+1}$ is $\sqrt{\smash[b]{R_{e}(z)}}$. As $\sqrt{\smash[b]{R_{e}(z)}}$ is
analytic $\forall \, \, z \! \in \! \mathbb{C} \setminus \cup_{j=1}^{N+1}
[b_{j-1}^{e},a_{j}^{e}]$, and since $C_{\mathrm{R}}^{e} \subset \mathbb{C}
\setminus \cup_{j=1}^{N+1}[b_{j-1}^{e},a_{j}^{e}]$, with $\operatorname{int}
(C_{\mathrm{R}}^{e}) \supset \overline{J_{e}} \cup \{z\}$, it follows that, in
particular, $\sqrt{\smash[b]{R_{e}(z)}} \! \upharpoonright_{C_{\mathrm{R}}^{
e}}$ is an analytic function of $\lbrace b_{j-1}^{e},a_{j}^{e} \rbrace_{j=
1}^{N+1}$, which implies, via the above (equivalent) contour integral
representation of $\mathcal{T}_{j}^{e}$, $j \! \in \! \mathbb{Z}_{0}^{+}$,
that $\mathcal{T}_{k}^{e}$, $k \! = \! 0,\dotsc,N \! + \! 1$, are (real)
analytic functions of $\lbrace b_{j-1}^{e},a_{j}^{e} \rbrace_{j=1}^{N+1}$.
Recalling that $(\mathcal{H} \psi_{V}^{e})(z) \! = \! \tfrac{1}{2 \pi}(\tfrac{
1}{z} \! + \! \tfrac{1}{2} \widetilde{V}^{\prime}(z) \! - \! (R_{e}(z))^{1/2}
h_{V}^{e}(z))$, it follows {}from the definition of $\mathcal{N}_{j}^{e}$,
$j \! = \! 1,\dotsc,N$, that
\begin{equation*}
\mathcal{N}_{j}^{e} \! = \! -\dfrac{1}{2 \pi} \int_{a_{j}^{e}}^{b_{j}^{e}}
(R_{e}(s))^{1/2}h_{V}^{e}(s) \, \md s, \quad j \! = \! 1,\dotsc,N:
\end{equation*}
making the linear change of variables $u_{j} \colon \mathbb{C} \! \to \!
\mathbb{C}$, $s \! \mapsto \! u_{j}(s) \! := \! (b_{j}^{e} \! - \! a_{j}^{
e})^{-1}(s \! - \! a_{j}^{e})$, $j \! = \! 1,\dotsc,N$, which take each of
the (compact) intervals $[a_{j}^{e},b_{j}^{e}]$, $j \! = \! 1,\dotsc,N$, onto
$[0,1]$, and setting
\begin{equation*}
\sqrt{\smash[b]{\widehat{R}_{e}(z)}}:= \! \left(\prod_{k_{1}=1}^{j}(z \! - \!
b_{k_{1}-1}^{e}) \prod_{k_{2}=1}^{j-1}(z \! - \! a_{k_{2}}^{e}) \prod_{k_{3}=j
+1}^{N+1}(a_{k_{3}}^{e} \! - \! z) \prod_{k_{4}=j+2}^{N+1}(b_{k_{4}-1}^{e} \!
- \! z) \right)^{1/2},
\end{equation*}
one arrives at
\begin{equation*}
\mathcal{N}_{j}^{e} \! = \! -\dfrac{1}{2 \pi} \! \left(b_{j}^{e} \! - \! a_{
j}^{e} \right)^{2} \! \int_{0}^{1} \! \left(u_{j}(1 \! - \! u_{j}) \right)^{
1/2} \left(\widehat{R}_{e}((b_{j}^{e} \! - \! a_{j}^{e})u_{j} \! + \! a_{j}^{
e}) \right)^{1/2}h_{V}^{e}((b_{j}^{e} \! - \! a_{j}^{e})u_{j} \! + \! a_{j}^{
e}) \, \md u_{j}, \quad j \! = \! 1,\dotsc,N.
\end{equation*}
Recalling that $h_{V}^{e}(z)$ is analytic on $\mathbb{R} \setminus \{0\}$, in
particular, $h_{V}^{e}(b_{j-1}^{e}),h_{V}^{e}(a_{j}^{e}) \! \not= \! 0$, $j \!
= \! 1,\dotsc,N \! + \! 1$, and that it is an analytic function of $\lbrace
b_{k-1}^{e}(z_{o}),a_{k}^{e}(z_{o}) \rbrace_{k=1}^{N+1}$ (since $-\infty \! <
\! b_{0}^{e} \! < \! a_{1}^{e} \! < \! b_{1}^{e} \! < \! a_{2}^{e} \! < \!
\dotsb \! < \! b_{N}^{e} \! < \! a_{N+1}^{e} \! < \! +\infty)$, and noting
{}from the definition of $\sqrt{\smash[b]{\widehat{R}_{e}(z)}}$ above that,
it, too, is an analytic function of $(b_{j}^{e} \! - \! a_{j}^{e})u_{j} \! +
\! a_{j}^{e}$, $(j,u_{j}) \! \in \! \{1,\dotsc,N\} \times [0,1]$, and thus an
analytic function of $\lbrace b_{j-1}^{e}(z_{o}),a_{j}^{e}(z_{o}) \rbrace_{j=
1}^{N+1}$, it follows that $\mathcal{N}_{j}^{e}$, $j \! = \! 1,\dotsc,N$, are
(real) analytic functions of $\lbrace b_{j-1}^{e}(z_{o}),a_{j}^{e}(z_{o})
\rbrace_{j=1}^{N+1}$.

Thus, as the Jacobian of the transformation $\lbrace b_{0}^{e}(z_{o}),\dotsc,
b_{N}^{e}(z_{o}),a_{1}^{e}(z_{o}),\dotsc,a_{N+1}^{e}(z_{o}) \rbrace \! \mapsto
\! \lbrace \mathcal{T}_{0}^{e},\dotsc,\linebreak[4]
\mathcal{T}_{N+1}^{e},\mathcal{N}_{1}^{e},\dotsc,\mathcal{N}_{N}^{e} \rbrace$
is non-zero whenever $\lbrace b_{j-1}^{e}(z_{o}),a_{j}^{e}(z_{o}) \rbrace_{j=
1}^{N+1}$, the end-points of the support of the `even' equilibrium measure,
are chosen so that, for regular $\widetilde{V} \colon \mathbb{R} \setminus \{
0\} \! \to \! \mathbb{R}$ satisfying conditions~(2.3)--(2.5), $\overline{J_{e}
}= \! \cup_{j=1}^{N+1}[b_{j-1}^{e},a_{j}^{e}]$, and $\mathcal{T}_{j}^{e}$, $j
\! = \! 0,\dotsc,N \! + \! 1$, and $\mathcal{N}_{k}^{e}$, $k \! = \! 1,\dotsc,
N$, are (real) analytic functions of $\lbrace b_{j-1}^{e}(z_{o}),a_{j}^{e}(z_{
o}) \rbrace_{j=1}^{N+1}$, it follows, via the Implicit Function Theorem, that
$b_{j-1}^{e}(z_{o}),a_{j}^{e}(z_{o})$, $j \! = \! 1,\dotsc,N \! + \! 1$, are
real analytic functions of $z_{o}$. \hfill $\qed$
\begin{eeeee}
It turns out that, for $\widetilde{V} \colon \mathbb{R} \setminus \{0\} \! 
\to \! \mathbb{R}$ (satisfying conditions~(2.3)--(2.5)) of the form
\begin{equation*}
\widetilde{V}(z) \! = \! \sum^{2m_{2}}_{k=-2m_{1}} \widetilde{\varrho}_{k}
z^{k},
\end{equation*}
with $\widetilde{\varrho}_{k} \! \in \! \mathbb{R}$, $k \! = \! -2m_{1},
\dotsc,2m_{2}$, $m_{1,2} \! \in \! \mathbb{N}$, and (since $\widetilde{V}
(\pm \infty),\widetilde{V}(0) \! > \! 0)$ $\widetilde{\varrho}_{-2m_{1}},
\widetilde{\varrho}_{2m_{2}} \! > \! 0$, the integral for $h_{V}^{e}
(z)$, that is, $h_{V}^{e}(z) \! = \! \tfrac{1}{2} \oint_{C_{\mathrm{R}}^{e}}
(R_{e}(s))^{-1/2}(\tfrac{\mi}{\pi s} \! + \! \tfrac{\mi \widetilde{V}^{\prime}
(s)}{2 \pi})(s \! - \! z)^{-1} \, \md s$, can be evaluated explicitly. Let
$C_{\mathrm{R}}^{e} \! = \! \widetilde{\Gamma}_{\infty}^{e} \cup \widetilde{
\Gamma}_{0}^{e}$, where $\widetilde{\Gamma}_{\infty}^{e} \! := \! \{\mathstrut
z^{\prime} \! = \! R \me^{\mi \vartheta}, \, R \! > \! 1/\varepsilon, \,
\vartheta \! \in \! [0,2 \pi]\}$ (oriented clockwise), and $\widetilde{
\Gamma}_{0}^{e} \! := \! \{\mathstrut z^{\prime} \! = \! r \me^{\mi
\vartheta}, \, 0 \! < \! r \! < \! \varepsilon, \, \vartheta \! \in \! [0,2
\pi]\}$ (oriented counter-clockwise), with $\varepsilon$ some arbitrarily
fixed, sufficiently small positive real number chosen such that: (i) $\partial
\{\mathstrut z^{\prime} \! \in \! \mathbb{C}; \, \vert z^{\prime} \vert \! =
\! \varepsilon\} \cap \partial \{\mathstrut z^{\prime} \! \in \! \mathbb{C};
\, \vert z^{\prime} \vert \! = \! 1/\varepsilon\} \! = \! \varnothing$; (ii)
$\{\mathstrut z^{\prime} \! \in \! \mathbb{C}; \, \vert z^{\prime} \vert \! <
\! \varepsilon\} \cap (J_{e} \cup \{z\}) \! = \! \varnothing$; (iii) $\{
\mathstrut z^{\prime} \! \in \! \mathbb{C}; \, \vert z^{\prime} \vert \! > \!
1/\varepsilon\} \cap (J_{e} \cup \{z\}) \! = \! \varnothing$; and (iv) $\{
\mathstrut z^{\prime} \! \in \! \mathbb{C}; \, \varepsilon \! < \! \vert
z^{\prime} \vert \! < \! 1/\varepsilon\} \! \supset \! J_{e} \cup \{z\}$. A
tedious, but otherwise straightforward, residue calculus calculation shows
that
\begin{align*}
h_{V}^{e}(z) \! =& \, \dfrac{1}{2}z^{2m_{2}-N-2} \, \, \sum_{j=0}^{2m_{2}-N-2}
\underset{\substack{0 \leqslant \vert k \vert + \vert l \vert \leqslant 2m_{2}
-j-N-2\\k_{i} \geqslant 0, \, \, l_{i} \geqslant 0, \, \, i \in \{0,\dotsc,N\}
}}{\sideset{}{'}{\sum}_{k_{0},\dotsc,k_{N}} \, \, \sideset{}{'}{\sum}_{l_{0},
\dotsc,l_{N}}}(2m_{2} \! - \! j) \widetilde{\varrho}_{2m_{2}-j} \! \left(
\prod_{p=0}^{N} \prod_{j_{p}=0}^{k_{p}-1} \! \left(\dfrac{1}{2} \! + \! j_{p}
\right) \! \right) \\
\times& \left(\prod_{q=0}^{N} \prod_{\widetilde{m}_{q}=0}^{l_{q}-1} \! \left(
\dfrac{1}{2} \! + \! \widetilde{m}_{q} \right) \! \right) \! \dfrac{\left(
\prod_{p^{\prime}=0}^{N}(b_{p^{\prime}}^{e})^{k_{p^{\prime}}} \right) \! \left(
\prod_{q^{\prime}=0}^{N}(a_{q^{\prime}+1}^{e})^{l_{q^{\prime}}} \right)}{\left(
\prod_{l^{\prime}=0}^{N}k_{l^{\prime}}! \right) \! \left(\prod_{\widetilde{m}^{
\prime}=0}^{N}l_{\widetilde{m}^{\prime}}! \right)} \, z^{-(j+\vert k \vert +
\vert l \vert)} \\
+& \, \dfrac{(-1)^{\mathcal{N}_{+}}(\prod_{k=1}^{N+1} \vert b_{k-1}^{e}a_{k}^{
e} \vert)^{-1/2}}{2z^{2m_{1}+1}} \, \, \sum_{j=-2m_{1}+1}^{0} \underset{
\substack{0 \leqslant \vert k \vert + \vert l \vert \leqslant 2m_{1}+j\\k_{i}
\geqslant 0, \, \, l_{i} \geqslant 0, \, \, i \in \{0,\dotsc,N\}}}{
\sideset{}{''}{\sum}_{k_{0},\dotsc,k_{N}} \, \, \sideset{}{''}{\sum}_{l_{0},
\dotsc,l_{N}}}(-2m_{1} \! - \! j) \widetilde{\varrho}_{-2m_{1}-j} \\
\times& \left(\prod_{p=0}^{N} \prod_{j_{p}=0}^{k_{p}-1} \! \left(\dfrac{1}{2}
\! + \! j_{p} \right) \! \right) \! \left(\prod_{q=0}^{N} \prod_{\widetilde{
m}_{q}=0}^{l_{q}-1} \! \left(\dfrac{1}{2} \! + \! \widetilde{m}_{q} \right) \!
\right) \! \dfrac{\left(\prod_{p^{\prime}=0}^{N}(b_{p^{\prime}}^{e})^{k_{p^{
\prime}}} \right)^{-1} \! \left(\prod_{q^{\prime}=0}^{N}(a_{q^{\prime}+1}^{e}
)^{l_{q^{\prime}}} \right)^{-1}}{\left(\prod_{l^{\prime}=0}^{N}k_{l^{\prime}}!
\right) \! \left(\prod_{\widetilde{m}^{\prime}=0}^{N}l_{\widetilde{m}^{\prime}
}! \right)} \\
\times& \, z^{\vert k \vert +\vert l \vert -j}+\dfrac{(-1)^{\mathcal{N}_{+}}
(\prod_{k=1}^{N+1} \vert b_{k-1}^{e}a_{k}^{e} \vert)^{-1/2}}{z},
\end{align*}
where $\mathcal{N}_{+} \! \in \! \lbrace 0,\dotsc,N \! + \! 1 \rbrace$ is the
number of bands to the right of $z \! = \! 0$, $\vert k \vert \! := \! k_{0}
\! + \! k_{1} \! + \! \dotsb \! + \! k_{N}$ $(\geqslant \! 0)$, $\vert l \vert
\! := \! l_{0} \! + \! l_{1} \! + \! \dotsb \! + \! l_{N}$ $(\geqslant \! 0)$,
and the primes (resp., double primes) on the summations mean that all possible
sums over $\{k_{l}\}_{l=0}^{N}$ and $\{l_{k}\}_{k=0}^{N}$ must be taken for
which $0 \! \leqslant \! k_{0} \! + \! \cdots \! + \! k_{N} \! + \! l_{0} \! +
\! \cdots \! + \! l_{N} \! \leqslant \! 2m_{2} \! - \! j \! - \! N \! - \! 2$,
$j \! = \! 0,\dotsc,2m_{2} \! - \! N \! - \! 2$, $k_{i} \! \geqslant \! 0$,
$l_{i} \! \geqslant \! 0$, $i \! = \! 0,\dotsc,N$ (resp., $0 \! \leqslant \!
k_{0} \! + \! \cdots \! + \! k_{N} \! + \! l_{0} \! + \! \cdots \! + \! l_{N}
\! \leqslant \! 2m_{1} \! + \! j$, $j \! = \! -2m_{1} \! + \! 1,\dotsc,0$,
$k_{i} \! \geqslant \! 0$, $l_{i} \! \geqslant \! 0$, $i \! = \! 0,\dotsc,
N)$. It is important to note that all of the above sums are finite sums: any
sums for which the upper limit is less than the lower limit are defined to be
zero, and any products in which the upper limit is less than the lower limit
are defined to be one; for example, $\sum_{j=0}^{-1}(\ast) \! := \! 0$ and
$\prod_{j=0}^{-1}(\ast) \! := \! 1$.

It is also interesting to note that one may derive explicit formulae for the
various moments of the `even' equilibrium measure, that is, $\int_{J_{e}}s^{
\pm m} \psi_{V}^{e}(s) \, \md s$, $m \! \in \! \mathbb{N}$, in terms of the
external field and the function $(R_{e}(z))^{1/2}$; without loss of
generality, and for demonstrative purposes only, consider, say, the following
moments: $\int_{J_{e}}s^{\pm j} \, \md \mu_{V}^{e}(s)$, $j \! = \! 1,2,3$ (the
calculations below straightforwardly generalise to $\int_{J_{e}}s^{\pm (k+3)}
\, \md \mu_{V}^{e}(s)$, $k \! \in \! \mathbb{N})$. Recall the following
formulae for $\mathscr{F}^{e}(z)$ given in the proof of Lemma~3.5:
\begin{gather*}
\mathscr{F}^{e}(z) \! = \! -\dfrac{1}{\pi \mi z} \! - \! \dfrac{2}{\pi \mi}
\int_{J_{e}} \dfrac{\md \mu_{V}^{e}(s)}{s \! - \! z}, \quad z \! \in \!
\mathbb{C} \setminus (J_{e} \cup \{0\}), \\
\mathscr{F}^{e}(z) \! = \! -\dfrac{1}{\pi \mi z} \! - \! (R_{e}(z))^{1/2}
\int_{J_{e}} \dfrac{(\frac{2 \mi}{\pi s} \! + \! \frac{\mi \widetilde{V}^{
\prime}(s)}{\pi})}{(R_{e}(s))^{1/2}_{+}(s \! - \! z)} \, \dfrac{\md s}{2 \pi
\mi}, \quad z \! \in \! \mathbb{C} \setminus (J_{e} \cup \{0\}).
\end{gather*}
One derives the following asymptotic expansions: (1) for $\mu_{V}^{e} \! \in
\! \mathcal{M}_{1}(\mathbb{R})$, in particular, $\int_{J_{e}}s^{-m} \, \md
\mu_{V}^{e}(s) \! < \! \infty$, $m \! \in \! \mathbb{N}$, $s \! \in \! J_{e}$
and $z \! \notin \! J_{e}$, with $\vert z/s \vert \! \ll \! 1$ (e.g., $\vert
z \vert \! \ll \! \min_{j=1,\dotsc,N+1}\{\vert b_{j-1}^{e} \! - \! a_{j}^{e}
\vert\})$, via the expansions $\tfrac{1}{z-s} \! = \! -\sum_{k=0}^{l} \tfrac{
z^{k}}{s^{k+1}} \! + \! \tfrac{z^{l+1}}{s^{l+1}(z-s)}$, $l \! \in \! \mathbb{
Z}_{0}^{+}$, and $\ln (1 \! - \! \ast) \! = \! -\sum_{k=1}^{\infty} \tfrac{
\ast^{k}}{k}$, $\vert \ast \vert \! \ll \! 1$,
\begin{equation*}
\mathscr{F}^{e}(z) \underset{z \to 0}{=} -\dfrac{1}{\pi \mi z} \! - \! \dfrac{
2}{\pi \mi} \int_{J_{e}}s^{-1} \, \md \mu_{V}^{e}(s) \! + \! z \! \left(-
\dfrac{2}{\pi \mi} \int_{J_{e}}s^{-2} \, \md \mu_{V}^{e}(s) \right) \! + \!
z^{2} \! \left(-\dfrac{2}{\pi \mi} \int_{J_{e}}s^{-3} \, \md \mu_{V}^{e}(s)
\right) \! + \! \mathcal{O}(z^{3}),
\end{equation*}
and
\begin{equation*}
\mathscr{F}^{e}(z) \underset{z \to 0}{=} -\dfrac{1}{\pi \mi z} \! + \!
\gamma_{V}^{e} \! \left(\check{Q}_{0}^{e} \! + \! z(\check{Q}_{1}^{e} \! - \!
\check{P}_{0}^{e} \check{Q}_{0}^{e}) \! + \! z^{2}(\check{Q}_{2}^{e} \! - \!
\check{P}_{0}^{e} \check{Q}_{1}^{e} \! + \! \check{P}_{1}^{e} \check{Q}_{0}^{
e}) \! + \! \mathcal{O}(z^{3}) \right),
\end{equation*}
where
\begin{gather*}
\gamma_{V}^{e} \! := \! (-1)^{\mathcal{N}_{+}} \left(\prod_{j=1}^{N+1}
\left\vert b_{j-1}^{e}a_{j}^{e} \right\vert \right)^{1/2}, \qquad \qquad
\check{P}_{0}^{e} \! := \! \dfrac{1}{2} \sum_{j=1}^{N+1} \! \left(\dfrac{1}{
b_{j-1}^{e}} \! + \! \dfrac{1}{a_{j}^{e}} \right), \\
\check{P}_{1}^{e} \! := \! \dfrac{1}{2}(\check{P}_{0}^{e})^{2} \! - \! \dfrac{
1}{4} \sum_{j=1}^{N+1} \! \left(\dfrac{1}{(b_{j-1}^{e})^{2}} \! + \! \dfrac{1}{
(a_{j}^{e})^{2}} \right), \qquad \check{Q}_{j}^{e} \! := \! -\int_{J_{e}}
\dfrac{(\frac{2 \mi}{\pi s} \! + \! \frac{\mi \widetilde{V}^{\prime}(s)}{\pi}
)}{(R_{e}(s))^{1/2}_{+}s^{j+1}} \, \dfrac{\md s}{2 \pi \mi}, \quad j \! = \!
0,1,2;
\end{gather*}
and (2) for $\mu_{V}^{e} \! \in \! \mathcal{M}_{1}(\mathbb{R})$, in
particular, $\int_{J_{e}} \md \mu_{V}^{e}(s) \! = \! 1$ and $\int_{J_{e}}s^{m}
\, \md \mu_{V}^{e}(s) \! < \! \infty$, $m \! \in \! \mathbb{N}$, $s \! \in \!
J_{e}$ and $z \! \notin \! J_{e}$, with $\vert s/z \vert \! \ll \! 1$ (e.g.,
$\vert z \vert \! \gg \! \max_{j=1,\dotsc,N+1}\{\vert b_{j-1}^{e} \! - \! a_{
j}^{e} \vert\})$, via the expansions $\tfrac{1}{s-z} \! = \! -\sum_{k=0}^{l}
\tfrac{s^{k}}{z^{k+1}} \! + \! \tfrac{s^{l+1}}{z^{l+1}(s-z)}$, $l \! \in \!
\mathbb{Z}_{0}^{+}$, and $\ln (1 \! - \! \ast) \! = \! -\sum_{k=1}^{\infty}
\tfrac{\ast^{k}}{k}$, $\vert \ast \vert \! \ll \! 1$,
\begin{equation*}
\mathscr{F}^{e}(z) \underset{z \to \infty}{=} \dfrac{1}{\pi \mi z} \! + \!
\dfrac{1}{z^{2}} \! \left(\dfrac{2}{\pi \mi} \int_{J_{e}}s \, \md \mu_{V}^{e}
(s) \right) \! + \! \dfrac{1}{z^{3}} \! \left(\dfrac{2}{\pi \mi} \int_{J_{e}}
s^{2} \, \md \mu_{V}^{e}(s) \right) \! + \! \dfrac{1}{z^{4}} \! \left(\dfrac{
2}{\pi \mi} \int_{J_{e}}s^{3} \, \md \mu_{V}^{e}(s) \right) \! + \! \mathcal{
O}(z^{-5}),
\end{equation*}
and
\begin{equation*}
\mathscr{F}^{e}(z) \underset{z \to \infty}{=} -\dfrac{1}{\pi \mi z} \! + \!
z^{N} \! \left(1 \! - \! \dfrac{\alpha_{V}^{e}}{z} \! + \! \dfrac{\widetilde{
P}_{0}^{e}}{z^{2}} \! + \! \dfrac{\widetilde{P}_{1}^{e}}{z^{3}} \! + \!
\dotsb \right) \! \int_{J_{e}} \dfrac{(\frac{2 \mi}{\pi s} \! + \! \frac{\mi
\widetilde{V}^{\prime}(s)}{\pi})}{(R_{e}(s))^{1/2}_{+}} \! \left(1 \! + \!
\cdots \! + \! \dfrac{s^{N}}{z^{N}} \! + \! \dfrac{s^{N+1}}{z^{N+1}} \! + \!
\dotsb \right) \! \dfrac{\md s}{2 \pi \mi},
\end{equation*}
where
\begin{gather*}
\alpha_{V}^{e} \! := \! \dfrac{1}{2} \sum_{j=1}^{N+1} \! \left(b_{j-1}^{e}
\! + \! a_{j}^{e} \right), \qquad \qquad \widetilde{P}_{0}^{e} \! := \! \dfrac{
1}{2}(\alpha_{V}^{e})^{2} \! - \! \dfrac{1}{4} \sum_{j=1}^{N+1} \! \left((b_{j
-1}^{e})^{2} \! + \! (a_{j}^{e})^{2} \right), \\
\widetilde{P}_{1}^{e} \! := \! -\dfrac{1}{3!} \sum_{j=1}^{N+1} \! \left((b_{j
-1}^{e})^{3} \! + \! (a_{j}^{e})^{3} \right) \! - \! \dfrac{(\alpha_{V}^{e})^{
3}}{3!} \! + \! \dfrac{\alpha_{V}^{e}}{4} \sum_{j=1}^{N+1} \! \left((b_{j-1}^{
e})^{2} \! + \! (a_{j}^{e})^{2} \right).
\end{gather*}
Recalling the following $(n$-dependent) $N \! + \! 2$ moment conditions stated
in Lemma~3.5,
\begin{equation*}
\int_{J_{e}} \dfrac{(\frac{2 \mi}{\pi s} \! + \! \frac{\mi \widetilde{V}^{
\prime}(s)}{\pi})s^{j}}{(R_{e}(s))^{1/2}_{+}} \, \md s \! = \! 0, \quad j \!
= \! 0,\dotsc,N, \qquad \text{and} \qquad \int_{J_{e}} \dfrac{(\frac{2 \mi}{
\pi s} \! + \! \frac{\mi \widetilde{V}^{\prime}(s)}{\pi})s^{N+1}}{(R_{e}(s)
)^{1/2}_{+}} \, \md s \! = \! 4,
\end{equation*}
and equating the respective pairs of asymptotic expansions above (as $z \! \to
\! 0$ and $z \! \to \! \infty)$ for $\mathscr{F}^{e}(z)$, one arrives at the
following expressions for the first three `positive' and `negative' moments of
the `even' equilibrium measure:
\begin{align*}
\int_{J_{e}}s \, \md \mu_{V}^{e}(s)=& \, \dfrac{1}{4} \int_{J_{e}} \dfrac{(
\frac{2 \mi}{\pi s} \! + \! \frac{\mi \widetilde{V}^{\prime}(s)}{\pi})s^{N+2}
}{(R_{e}(s))^{1/2}_{+}} \, \md s \! - \! \dfrac{1}{2} \sum_{j=1}^{N+1}(b_{j-
1}^{e} \! + \! a_{j}^{e}), \\
\int_{J_{e}}s^{2} \, \md \mu_{V}^{e}(s)=& \, \dfrac{1}{4} \int_{J_{e}} \dfrac{
(\frac{2 \mi}{\pi s} \! + \! \frac{\mi \widetilde{V}^{\prime}(s)}{\pi})s^{N+3}
}{(R_{e}(s))^{1/2}_{+}} \, \md s \! - \! \dfrac{1}{8} \! \left(\sum_{j=1}^{N
+1}(b_{j-1}^{e} \! + \! a_{j}^{e}) \right) \! \int_{J_{e}} \dfrac{(\frac{2
\mi}{\pi s} \! + \! \frac{\mi \widetilde{V}^{\prime}(s)}{\pi})s^{N+2}}{(R_{e}
(s))^{1/2}_{+}} \, \md s \\
+& \, \dfrac{1}{4} \! \left(\dfrac{1}{2} \! \left(\sum_{j=1}^{N+1}(b_{j-1}^{
e} \! + \! a_{j}^{e}) \right)^{2} \! - \! \sum_{j=1}^{N+1}((b_{j-1}^{e})^{2}
\! + \! (a_{j}^{e})^{2}) \right), \\
\int_{J_{e}}s^{3} \, \md \mu_{V}^{e}(s)=& \, \dfrac{1}{4} \int_{J_{e}} \dfrac{
(\frac{2 \mi}{\pi s} \! + \! \frac{\mi \widetilde{V}^{\prime}(s)}{\pi})s^{N+4}
}{(R_{e}(s))^{1/2}_{+}} \, \md s \! - \! \dfrac{1}{8} \! \left(\sum_{j=1}^{N
+1}(b_{j-1}^{e} \! + \! a_{j}^{e}) \right) \! \int_{J_{e}} \dfrac{(\frac{2
\mi}{\pi s} \! + \! \frac{\mi \widetilde{V}^{\prime}(s)}{\pi})s^{N+3}}{(R_{e}
(s))^{1/2}_{+}} \, \md s \\
+& \, \dfrac{1}{16} \! \left(\dfrac{1}{2} \! \left(\sum_{j=1}^{N+1}(b_{j-1}^{
e} \! + \! a_{j}^{e}) \right)^{2} \! - \! \sum_{j=1}^{N+1}((b_{j-1}^{e})^{2}
\! + \! (a_{j}^{e})^{2}) \right) \! \int_{J_{e}} \dfrac{(\frac{2 \mi}{\pi s}
\! + \! \frac{\mi \widetilde{V}^{\prime}(s)}{\pi})s^{N+2}}{(R_{e}(s))^{1/2}_{
+}} \, \md s \\
-& \, \dfrac{1}{8} \! \left(\dfrac{1}{3!} \! \left(\sum_{j=1}^{N+1}(b_{j-1}^{
e} \! + \! a_{j}^{e}) \right)^{3} \! + \! \dfrac{4}{3} \sum_{j=1}^{N+1}((b_{j
-1}^{e})^{3} \! + \! (a_{j}^{e})^{3}) \! - \! \sum_{j=1}^{N+1}(b_{j-1}^{e} \!
+ \! a_{j}^{e}) \right. \\
\times& \left. \sum_{k=1}^{N+1}((b_{k-1}^{e})^{2} \! + \! (a_{k}^{e})^{2})
\right), \\
\int_{J_{e}}s^{-1} \, \md \mu_{V}^{e}(s)=& \, \dfrac{1}{4}(-1)^{\mathcal{N}_{
+}} \! \left(\prod_{j=1}^{N+1} \vert b_{j-1}^{e}a_{j}^{e} \vert \right)^{1/2}
\! \int_{J_{e}} \dfrac{(\frac{2 \mi}{\pi s} \! + \! \frac{\mi \widetilde{V}^{
\prime}(s)}{\pi})}{(R_{e}(s))^{1/2}_{+}s} \, \md s, \\
\int_{J_{e}}s^{-2} \, \md \mu_{V}^{e}(s)=& \, \dfrac{1}{4}(-1)^{\mathcal{N}_{
+}} \! \left(\prod_{j=1}^{N+1} \vert b_{j-1}^{e}a_{j}^{e} \vert \right)^{1/2}
\! \left(\int_{J_{e}} \dfrac{(\frac{2 \mi}{\pi s} \! + \! \frac{\mi
\widetilde{V}^{\prime}(s)}{\pi})}{(R_{e}(s))^{1/2}_{+}s^{2}} \, \md s \! - \!
\dfrac{1}{2} \! \left(\sum_{j=1}^{N+1} \! \left(\dfrac{1}{b_{j-1}^{e}} \! + \!
\dfrac{1}{a_{j}^{e}} \right) \! \right) \right. \\
\times& \left. \int_{J_{e}} \dfrac{(\frac{2 \mi}{\pi s} \! + \! \frac{\mi
\widetilde{V}^{\prime}(s)}{\pi})}{(R_{e}(s))^{1/2}_{+}s} \, \md s \right), \\
\int_{J_{e}}s^{-3} \, \md \mu_{V}^{e}(s)=& \, \dfrac{1}{4}(-1)^{\mathcal{N}_{
+}} \! \left(\prod_{j=1}^{N+1} \vert b_{j-1}^{e}a_{j}^{e} \vert \right)^{1/2}
\! \left(\int_{J_{e}} \dfrac{(\frac{2 \mi}{\pi s} \! + \! \frac{\mi
\widetilde{V}^{\prime}(s)}{\pi})}{(R_{e}(s))^{1/2}_{+}s^{3}} \, \md s \! - \!
\dfrac{1}{2} \! \left(\sum_{j=1}^{N+1} \! \left(\dfrac{1}{b_{j-1}^{e}} \! + \!
\dfrac{1}{a_{j}^{e}} \right) \! \right) \right. \\
\times& \left. \int_{J_{e}} \dfrac{(\frac{2 \mi}{\pi s} \! + \! \frac{\mi
\widetilde{V}^{\prime}(s)}{\pi})}{(R_{e}(s))^{1/2}_{+}s^{2}} \, \md s \! + \!
\left(\dfrac{1}{8} \! \left(\sum_{j=1}^{N+1} \! \left(\dfrac{1}{b_{j-1}^{e}}
\! + \! \dfrac{1}{a_{j}^{e}} \right) \right)^{2} \! - \! \dfrac{1}{4} \sum_{
j=1}^{N+1} \! \left(\dfrac{1}{(b_{j-1}^{e})^{2}} \! + \! \dfrac{1}{(a_{j}^{
e})^{2}} \right) \right) \right. \\
\times& \left. \int_{J_{e}} \dfrac{(\frac{2 \mi}{\pi s} \! + \! \frac{\mi
\widetilde{V}^{\prime}(s)}{\pi})}{(R_{e}(s))^{1/2}_{+}s} \, \md s \right).
\end{align*}
It is important to note that all of the above integrals are real valued
(since, for $s \! \in \! \overline{J_{e}}$, $(R_{e}(s))^{1/2}_{+} \! = \! \mi
(\vert R_{e}(s) \vert)^{1/2} \! \in \! \mi \mathbb{R})$ and bounded (since,
for $j \! = \! 1,\dotsc,N \! + \! 1$, $(R_{e}(s))^{1/2} \! =_{s \downarrow
b_{j-1}^{e}} \! \mathcal{O}((s \! - \! b_{j-1}^{e})^{1/2})$ and $(R_{e}(s))^{
1/2} \! =_{s \uparrow a_{j}^{e}} \! \mathcal{O}((a_{j}^{e} \! - \! s)^{1/2})$,
that is, there are removable singularities at the end-points of the support of
the `even' equilibrium measure). \hfill $\blacksquare$
\end{eeeee}
\begin{ccccc}
Let the external field $\widetilde{V} \colon \mathbb{R} \setminus \{0\} \! 
\to \! \mathbb{R}$ satisfy conditions~{\rm (2.3)--(2.5)}. Let the `even' 
equilibrium measure, $\mu_{V}^{e}$, and its support, $\operatorname{supp}
(\mu_{V}^{e}) \! =: \! J_{e}$ $(\subset \overline{\mathbb{R}} \setminus 
\lbrace 0,\pm \infty \rbrace)$, be as described in Lemma~{\rm 3.5}, and 
let there exist $\ell_{e}$ $(\in \! \mathbb{R})$, the `even' variational 
constant, such that
\begin{equation}
\begin{aligned}
4 \int_{J_{e}} \ln (\vert x \! - \! s \vert) \psi_{V}^{e}(s) \, \md s \! - \!
2 \ln \vert x \vert \! - \! \widetilde{V}(x) \! - \! \ell_{e} =& \, 0, \quad
x \! \in \! \overline{J_{e}}, \\
4 \int_{J_{e}} \ln (\vert x \! - \! s \vert) \psi_{V}^{e}(s) \, \md s \! - \!
2 \ln \vert x \vert \! - \! \widetilde{V}(x) \! - \! \ell_{e} \leqslant& \, 0,
\quad x \! \in \! \mathbb{R} \setminus \overline{J_{e}},
\end{aligned}
\end{equation}
where, for $\widetilde{V}$ regular, the inequality in the second of
Equations~{\rm (3.9)} is strict. Then:
\begin{compactenum}
\item[{\rm (1)}] $g^{e}_{+}(z) \! + \! g^{e}_{-}(z) \! - \! \widetilde{V}(z)
\! - \! \ell_{e} \! + \! 2Q_{e} \! = \! 0$, $z \! \in \! \overline{J_{e}}$,
where $g^{e}_{\pm}(z) \! := \! \lim_{\varepsilon \downarrow 0}g^{e}(z \! \pm
\! \mi \varepsilon)$, and $Q_{e}$ is defined in Lemma~{\rm 3.4;}
\item[{\rm (2)}] $g^{e}_{+}(z) \! + \! g^{e}_{-}(z) \! - \! \widetilde{V}(z)
\! - \! \ell_{e} \! + \! 2Q_{e} \! \leqslant \! 0$, $z \! \in \! \mathbb{R}
\setminus \overline{J_{e}}$, where equality holds for at most a finite number
of points, and, for $\widetilde{V}$ regular, the inequality is strict;
\item[{\rm (3)}] $g^{e}_{+}(z) \! - \! g^{e}_{-}(z) \! \in \! \mi f_{g^{e}}^{
\mathbb{R}}(z)$, $z \! \in \! \mathbb{R}$, where $f_{g^{e}}^{\mathbb{R}}
\colon \mathbb{R} \! \to \! \mathbb{R}$ is some bounded function, and, in
particular, $g^{e}_{+}(z) \! - \! g^{e}_{-}(z) \! = \! \mi
\operatorname{const.}$, $z \! \in \! \mathbb{R} \setminus \overline{J_{e}}$,
where $\operatorname{const.} \! \in \! \mathbb{R};$
\item[{\rm (4)}] $\mi (g^{e}_{+}(z) \! - \! g^{e}_{-}(z))^{\prime} \! 
\geqslant \! 0$, $z \! \in \! J_{e}$, and where, for $\widetilde{V}$ regular,
equality holds for at most a finite number of points.
\end{compactenum}
\end{ccccc}

\emph{Proof.} Set (cf. Lemma~3.5) $J_{e} \! := \! \cup_{j=1}^{N+1}J_{j}^{e}$, 
where $J_{j}^{e} \! = \! (b_{j-1}^{e},a_{j}^{e}) \! =$ the $j$th `band', 
with $N \! \in \! \mathbb{N}$ and finite, $b_{0}^{e} \! := \! \min \{
\operatorname{supp}(\mu_{V}^{e})\} \! \notin \! \lbrace -\infty,0 \rbrace$, 
$a_{N+1}^{e} \! := \! \max \{\operatorname{supp}(\mu_{V}^{e})\} \! \notin \! 
\lbrace 0,+\infty \rbrace$, and $-\infty \! < \! b_{0}^{e} \! < \! a_{1}^{e} 
\! < \! b_{1}^{e} \! < \! a_{2}^{e} \! < \! \cdots \! < \! b_{N}^{e} \! < \! 
a_{N+1}^{e} \! < \! +\infty$, and $\lbrace b_{j-1}^{e},a_{j}^{e} \rbrace_{j=
1}^{N+1}$ satisfy the $n$-dependent and (locally) solvable system of $2(N \! 
+ \! 1)$ moment conditions given in Lemma~3.5. Consider the following cases: 
\textbf{(1)} $z \! \in \! \overline{J_{j}^{e}} \! := \! [b_{j-1}^{e},a_{j}^{
e}]$, $j \! = \! 1,\dotsc,N \! + \! 1$; \textbf{(2)} $z \! \in \! (a_{j}^{e},
b_{j}^{e}) \! =$ the $j$th `gap', $j \! = \! 1,\dotsc,N$; \textbf{(3)} $z 
\! \in \! (a_{N+1}^{e},+\infty)$; and \textbf{(4)} $z \! \in \! (-\infty,
b_{0}^{e})$.

$\mathbf{(1)}$ Recall the definition of $g^{e}(z)$ given in Lemma~3.4, namely,
$g^{e}(z) \! := \! \int_{J_{e}} \ln ((z \! - \! s)^{2}/zs) \psi_{V}^{e}(s) \,
\md s$, $z \! \in \! \mathbb{C} \setminus (-\infty,\max \{0,a_{N+1}^{e}\})$,
where the representation (cf. Lemma~3.5) $\md \mu_{V}^{e}(s) \! = \! \psi_{
V}^{e}(s) \, \md s$, $s \! \in \! J_{e}$, was substituted into the latter
formula. For $z \! \in \! \overline{J_{j}^{e}}$, $j \! = \! 1,\dotsc,N \! + \!
1$, one shows that
\begin{equation*}
g^{e}_{\pm}(z) \! = \! 2 \int_{J_{e}} \ln (\vert z \! - \! s \vert) \psi_{V}^{
e}(s) \, \md s \! \pm \! 2 \pi \mi \int_{z}^{a_{N+1}^{e}} \psi_{V}^{e}(s) \,
\md s \! - \! Q_{e} \! - \!
\begin{cases}
\ln \vert z \vert, &\text{$z \! > \! 0,$} \\
\ln \vert z \vert \! \pm \! \mi \pi, &\text{$z \! < \! 0$,}
\end{cases}
\end{equation*}
where $g^{e}_{\pm}(z) \! := \! \lim_{\varepsilon \downarrow 0}g^{e}(z \! \pm
\! \mi \varepsilon)$, and $Q_{e} \! := \! \int_{J_{e}} \ln (s) \psi_{V}^{e}(s)
\, \md s$, whence
\begin{equation*}
g^{e}_{+}(z) \! - \! g^{e}_{-}(z) \! = \! 4 \pi \mi \int_{z}^{a_{N+1}^{e}}
\psi_{V}^{e}(s) \, \md s \! + \!
\begin{cases}
0, &\text{$z \! > \! 0$,} \\
-2 \pi \mi, &\text{$z \! < \! 0$,}
\end{cases}
\end{equation*}
which shows that $g^{e}_{+}(z) \! - \! g^{e}_{-}(z) \! \in \! \mi \mathbb{R}$,
and $\Re (g^{e}_{+}(z) \! - \! g^{e}_{-}(z)) \! = \! 0$; moreover, using the
Fundamental Theorem of Calculus, one shows that $(g^{e}_{+}(z) \! - \! g^{e}_{
-}(z))^{\prime} \! = \! -4 \pi \mi \psi_{V}^{e}(z)$, whence $\mi (g^{e}_{+}(z)
\! - \! g^{e}_{-}(z))^{\prime} \! = \! 4 \pi \psi_{V}^{e}(z) \! \geqslant \!
0$, since $\psi_{V}^{e}(z) \! \geqslant \! 0 \, \, \forall \, \, z \! \in \!
\overline{J_{e}}$ $(\supset \overline{J_{j}^{e}}$, $j \! = \! 1,\dotsc,N \! +
\! 1)$. Furthermore, using the first of Equations~(3.9), one shows that
\begin{equation*}
g^{e}_{+}(z) \! + \! g^{e}_{-}(z) \! - \! \widetilde{V}(z) \! - \! l_{e} \! +
\! 2Q_{e} \! = \! 4 \int_{J_{e}} \ln (\vert z \! - \! s \vert) \psi_{V}^{e}(s)
\, \md s \! - \! 2 \ln \vert z \vert \! - \! \widetilde{V}(z) \! - \! \ell_{e}
\! = \! 0,
\end{equation*}
which gives the formula for the `even' variational constant $\ell_{e}$ $(\in 
\! \mathbb{R})$, which is the same \cite{a92,a97} (see, also, Section~7 of 
\cite{a56}) for each compact interval $\overline{J_{j}^{e}}$, $j \! = \! 1,
\dotsc,N \! + \! 1$; in particular,
\begin{equation*}
\ell_{e} \! = \! \dfrac{2}{\pi} \sum_{j=1}^{N+1} \int_{b_{j-1}^{e}}^{a_{j}^{
e}} \ln \! \left(\left\vert \tfrac{1}{2}(b_{N}^{e} \! + \! a_{N+1}^{e}) \! -
\! s \right\vert \right) \! \left(\vert R_{e}(s) \vert \right)^{1/2} \! h_{
V}^{e}(s) \, \md s \! - \! 2 \ln \! \left\vert \tfrac{1}{2}(b_{N}^{e} \! + \!
a_{N+1}^{e}) \right\vert \! - \! \widetilde{V} \! \left(\tfrac{1}{2}(b_{N}^{
e} \! + \! a_{N+1}^{e}) \right),
\end{equation*}
where $(\vert R_{e}(s) \vert)^{1/2}h_{V}^{e}(s) \! \geqslant \! 0$, $j \! = \!
1,\dotsc,N \! + \! 1$, and where there are no singularities in the integrand,
since, for (any) $r \! > \! 0$, $\lim_{\vert x \vert \to 0} \vert x \vert^{r}
\ln \vert x \vert \! = \! 0$.

$\mathbf{(2)}$ For $z \! \in \! (a_{j}^{e},b_{j}^{e})$, $j \! = \! 1,\dotsc,
N$, one shows that
\begin{equation*}
g^{e}_{\pm}(z) \! = \! 2 \int_{J_{e}} \ln (\vert z \! - \! s \vert) \psi_{V}^{
e}(s) \, \md s \! \pm \! 2 \pi \mi \sum_{k=j+1}^{N+1} \int_{b_{k-1}^{e}}^{a_{
k}^{e}} \psi_{V}^{e}(s) \, \md s \! - \! Q_{e} \! - \!
\begin{cases}
\ln \vert z \vert, &\text{$z \! > \! 0$,} \\
\ln \vert z \vert \! \pm \! \mi \pi, &\text{$z \! < \! 0$,}
\end{cases}
\end{equation*}
whence
\begin{equation*}
g^{e}_{+}(z) \! - \! g^{e}_{-}(z) \! = \! 4 \pi \mi \int_{b_{j}^{e}}^{a_{N+
1}^{e}} \psi_{V}^{e}(s) \, \md s \! + \!
\begin{cases}
0, &\text{$z \! > \! 0$,} \\
-2 \pi \mi, &\text{$z \! < \! 0$,}
\end{cases}
\end{equation*}
which shows that $g^{e}_{+}(z) \! - \! g^{e}_{-}(z) \! = \! \mi \operatorname{
const.}$, with $\operatorname{const.} \! \in \! \mathbb{R}$, and $\Re (g^{e}_{
+}(z) \! - \! g^{e}_{-}(z)) \! = \! 0$; moreover, $\mi (g^{e}_{+}(z) \! - \!
g^{e}_{-}(z))^{\prime} \! = \! 0$. One notes {}from the above formulae for
$g^{e}_{\pm}(z)$ that
\begin{equation*}
g^{e}_{+}(z) \! + \! g^{e}_{-}(z) \! - \! \widetilde{V}(z) \! - \! \ell_{e} \!
+ \! 2Q_{e} \! = \! 4 \int_{J_{e}} \ln (\vert z \! - \! s \vert) \psi_{V}^{e}
(s) \, \md s \! - \! 2 \ln \vert z \vert \! - \! \widetilde{V}(z) \! - \!
\ell_{e}.
\end{equation*}
Recalling that (cf. Lemma~3.5) $\mathcal{H} \colon \mathcal{L}^{2}_{\mathrm{
M}_{2}(\mathbb{C})} \! \to \! \mathcal{L}^{2}_{\mathrm{M}_{2}(\mathbb{C})}$,
$f \! \mapsto \! (\mathcal{H}f)(z) \! := \!
\pvi_{\raise-0.95ex\hbox{$\scriptstyle{} \mathbb{R}$}} \tfrac{f(s)}{z-s} \,
\tfrac{\md s}{\pi}$, where $\pvi_{}$ denotes the principle value integral, one
shows that, for $z \! \in \! (a_{j}^{e},b_{j}^{e})$, $j \! = \! 1,\dotsc,N$,
\begin{equation*}
4 \int_{J_{e}} \ln (\vert z \! - \! s \vert) \psi_{V}^{e}(s) \, \md s \! = \!
4 \pi \int_{a_{j}^{e}}^{z}(\mathcal{H} \psi_{V}^{e})(s) \, \md s \! + \! 4
\int_{J_{e}} \ln (\vert a_{j}^{e} \! - \! s \vert) \psi_{V}^{e}(s) \, \md s;
\end{equation*}
thus,
\begin{align*}
g^{e}_{+}(z) \! + \! g^{e}_{-}(z) \! - \! \widetilde{V}(z) \! - \! \ell_{e} \!
+ \! 2Q_{e}=& \, 4 \pi \int_{a_{j}^{e}}^{z}(\mathcal{H} \psi_{V}^{e})(s) \,
\md s \! + \! 4 \int_{J_{e}} \ln (\vert a_{j}^{e} \! - \! s \vert) \psi_{V}^{e}
(s) \, \md s \! - \! 2 \ln \vert z \vert \! - \! \widetilde{V}(z) \! - \!
\ell_{e} \\
=& \, 4 \int_{J_{e}} \ln (\vert a_{j}^{e} \! - \! s \vert) \psi_{V}^{e}(s) \,
\md s \! + \! 4 \pi \int_{a_{j}^{e}}^{z}(\mathcal{H} \psi_{V}^{e})(s) \, \md s
\! - \! 4 \pi \int_{a_{j}^{e}}^{z} \dfrac{\widetilde{V}^{\prime}(s)}{4 \pi} \,
\md s \\
-& \, 4 \pi \int_{a_{j}^{e}}^{z} \dfrac{1}{2 \pi s} \, \md s \! - \! 2 \ln
\vert a_{j}^{e} \vert \! - \! \widetilde{V}(a_{j}^{e}) \! - \! \ell_{e} \\
=& \, 4 \pi \int_{a_{j}^{e}}^{z} \! \left((\mathcal{H} \psi_{V}^{e})(s) \! -
\! \dfrac{\widetilde{V}^{\prime}(s)}{4 \pi} \! - \! \dfrac{1}{2 \pi s} \right)
\! \md s \\
+& \underbrace{\left(4 \int_{J_{e}} \ln (\vert a_{j}^{e} \! - \! s \vert)
\psi_{V}^{e}(s) \, \md s \! - \! 2 \ln \vert a_{j}^{e} \vert \! - \!
\widetilde{V}(a_{j}^{e}) \! - \! \ell_{e} \right)}_{= \, 0} \, \, \,
\Rightarrow
\end{align*}
\begin{equation*}
g^{e}_{+}(z) \! + \! g^{e}_{-}(z) \! - \! \widetilde{V}(z) \! - \! \ell_{e} \!
+ \! 2Q_{e} \! = \! 4 \pi \int_{a_{j}^{e}}^{z} \! \left((\mathcal{H} \psi_{V}^{
e})(s) \! - \! \dfrac{\widetilde{V}^{\prime}(s)}{4 \pi} \! - \! \dfrac{1}{2
\pi s} \right) \! \md s, \quad z \! \in \! (a_{j}^{e},b_{j}^{e}), \, \, \, \,
j \! = \! 1,\dotsc,N.
\end{equation*}
It was shown in the proof of Lemma~3.5 that $(\mathcal{H} \psi_{V}^{e})(s) \!
= \! \tfrac{\widetilde{V}^{\prime}(s)}{4 \pi} \! + \! \tfrac{1}{2 \pi s} \! -
\! \tfrac{1}{2 \pi}(R_{e}(s))^{1/2}h_{V}^{e}(s)$, $s \! \in \! (a_{j}^{e},b_{
j}^{e})$, $j \! = \! 1,\dotsc,N$, whence
\begin{equation*}
g^{e}_{+}(z) \! + \! g^{e}_{-}(z) \! - \! \widetilde{V}(z) \! - \! \ell_{e} \!
+ \! 2Q_{e} \! = \! -2 \int_{a_{j}^{e}}^{z}(R_{e}(s))^{1/2}h_{V}^{e}(s) \, \md
s \! < \! 0, \quad z \! \in \! \cup_{j=1}^{N}(a_{j}^{e},b_{j}^{e}):
\end{equation*}
since $h_{V}^{e}(z)$ is real analytic on $\mathbb{R} \setminus \{0\}$ and
$(R_{e}(s))^{1/2}h_{V}^{e}(s) \! > \! 0 \, \, \forall \, \, s \! \in \! \cup_{
j=1}^{N}(a_{j}^{e},b_{j}^{e})$, it follows that one has equality only at
points $z \! \in \! \cup_{j=1}^{N}(a_{j}^{e},b_{j}^{e})$ for which $h_{V}^{e}
(z) \! = \! 0$, which are finitely denumerable. (Note that, for $z \! \in \!
\cup_{j=1}^{N}(a_{j}^{e},b_{j}^{e})$, $(R_{e}(s))^{1/2}_{+} \! = \! (R_{e}
(s))^{1/2}_{-} \! = \! (R_{e}(s))^{1/2}$.)

$\mathbf{(3)}$ For $z \! \in \! (a_{N+1}^{e},+\infty)$, one shows that
\begin{equation*}
g^{e}_{\pm}(z) \! = \! 2 \int_{J_{e}} \ln (\vert z \! - \! s \vert) \psi_{V}^{
e}(s) \, \md s \! - \! Q_{e} \! - \!
\begin{cases}
\ln \vert z \vert, &\text{$z \! > \! 0$,} \\
\ln \vert z \vert \! \pm \! \mi \pi, &\text{$z \! < \! 0$,}
\end{cases}
\end{equation*}
whence
\begin{equation*}
g^{e}_{+}(z) \! - \! g^{e}_{-}(z) \! = \!
\begin{cases}
0, &\text{$z \! > \! 0$,} \\
-2 \pi \mi, &\text{$z \! < \! 0$,}
\end{cases}
\end{equation*}
which shows that $g^{e}_{+}(z) \! - \! g^{e}_{-}(z)$ is pure imaginary, and
$\mi (g^{e}_{+}(z) \! - \! g^{e}_{-}(z))^{\prime} \! = \! 0$. Also, one shows
that
\begin{equation*}
g^{e}_{+}(z) \! + \! g^{e}_{-}(z) \! - \! \widetilde{V}(z) \! - \! \ell_{e} \!
+ \! 2Q_{e} \! = \! 4 \int_{J_{e}} \ln (\vert z \! - \! s \vert) \psi_{V}^{e}
(s) \, \md s \! - \! 2 \ln \vert z \vert \! - \! \widetilde{V}(z) \! - \!
\ell_{e};
\end{equation*}
and, following the analysis of case~$\mathbf{(2)}$ above, one shows that, for
$z \! \in \! (a_{N+1}^{e},+\infty)$,
\begin{equation*}
4 \int_{J_{e}} \ln (\vert z \! - \! s \vert) \psi_{V}^{e}(s) \, \md s \! - \!
2 \ln \vert z \vert \! - \! \widetilde{V}(z) \! - \! \ell_{e} \! = \! 4 \pi
\int_{a_{N+1}^{e}}^{z} \! \left((\mathcal{H} \psi_{V}^{e})(s) \! - \! \dfrac{
\widetilde{V}^{\prime}(s)}{4 \pi} \! - \! \dfrac{1}{2 \pi s} \right) \md s,
\end{equation*}
thus, via the relation (cf. case~$\mathbf{(2)}$ above) $(\mathcal{H} \psi_{
V}^{e})(s) \! = \! \tfrac{\widetilde{V}^{\prime}(s)}{4 \pi} \! + \! \tfrac{1}{
2 \pi s} \! - \! \tfrac{1}{2 \pi}(R_{e}(s))^{1/2}h_{V}^{e}(s)$, $s \! \in \!
(a_{N+1}^{e},+\infty)$, one arrives at
\begin{equation*}
g^{e}_{+}(z) \! + \! g^{e}_{-}(z) \! - \! \widetilde{V}(z) \! - \! \ell_{e} \!
+ \! 2Q_{e} \! = \! -2 \int_{a_{N+1}^{e}}^{z}(R_{e}(s))^{1/2}h_{V}^{e}(s) \,
\md s \! < \! 0, \quad z \! \in \! (a_{N+1}^{e},+\infty).
\end{equation*}
If: (1) $z \! \to \! +\infty$ (e.g., $\vert z \vert \! \gg \! \max_{j=1,\dotsc,
N+1}\{\vert b_{j-1}^{e},a_{j}^{e} \vert\})$, $s \! \in \! J_{e}$, and $\vert
s/z \vert \! \ll \! 1$, {}from $\mu_{V}^{e} \! \in \! \mathcal{M}_{1}(\mathbb{
R})$, in particular, $\int_{J_{e}} \md \mu_{V}^{e}(s) \! = \! 1$ and $\int_{
J_{e}}s^{m} \, \md \mu_{V}^{e}(s) \! < \! \infty$, $m \! \in \! \mathbb{N}$,
the formula for $g^{e}_{+}(z) \! + \! g^{e}_{-}(z) \! - \! \widetilde{V}(z) \!
- \! \ell_{e} \! + \! 2Q_{e}$ above, and the expansions $\tfrac{1}{s-z} \! =
\! -\sum_{k=0}^{l} \tfrac{s^{k}}{z^{k+1}} \! + \! \tfrac{s^{l+1}}{z^{l+1}(s-z)
}$, $l \! \in \! \mathbb{Z}_{0}^{+}$, and $\ln (z \! - \! s) \! =_{\vert z
\vert \to \infty} \! \ln (z) \! - \! \sum_{k=1}^{\infty} \tfrac{1}{k}(\tfrac{
s}{z})^{k}$, one shows that
\begin{equation*}
g^{e}_{+}(z) \! + \! g^{e}_{-}(z) \! - \! \widetilde{V}(z) \! - \! \ell_{e}
\! + \! 2Q_{e} \underset{z \to +\infty}{=} \ln (z^{2} \! + \! 1) \! - \!
\widetilde{V}(z) \! + \! \mathcal{O}(1),
\end{equation*}
which, upon recalling that (cf. condition~(2.4)) $\lim_{\vert x \vert \to
\infty}(\widetilde{V}(x)/\ln (x^{2} \! + \! 1)) \! = \! +\infty$, shows that
$g^{e}_{+}(z) \! + \! g^{e}_{-}(z) \! - \! \widetilde{V}(z) \! - \! \ell_{e}
\! + \! 2Q_{e} \! < \! 0$; and (2) $\vert z \vert \! \to \! 0$ (e.g., $\vert z
\vert \! \ll \! \min_{j=1,\dotsc,N+1}\{\vert b_{j-1}^{e},a_{j}^{e} \vert\})$,
$s \! \in \! J_{e}$, and $\vert z/s \vert \! \ll \! 1$, {}from $\mu_{V}^{e} \!
\in \! \mathcal{M}_{1}(\mathbb{R})$, in particular, $\int_{J_{e}}s^{-m} \, \md
\mu_{V}^{e}(s) \! < \! \infty$, $m \! \in \! \mathbb{N}$, the above formula
for $g^{e}_{+}(z) \! + \! g^{e}_{-}(z) \! - \! \widetilde{V}(z) \! - \! \ell_{
e} \! + \! 2Q_{e}$, and the expansions $\tfrac{1}{z-s} \! = \! -\sum_{k=0}^{l}
\tfrac{z^{k}}{s^{k+1}} \! + \! \tfrac{z^{l+1}}{s^{l+1}(z-s)}$, $l \! \in \!
\mathbb{Z}_{0}^{+}$, and $\ln (s \! - \! z) \! =_{\vert z \vert \to 0} \! \ln
(s) \! - \! \sum_{k=1}^{\infty} \tfrac{1}{k}(\tfrac{z}{s})^{k}$, one shows that
\begin{equation*}
g^{e}_{+}(z) \! + \! g^{e}_{-}(z) \! - \! \widetilde{V}(z) \! - \! \ell_{e} \!
+ \! 2Q_{e} \underset{\vert z \vert \to 0}{=} \ln (z^{-2} \! + \! 1) \! - \!
\widetilde{V}(z) \! + \! \mathcal{O}(1),
\end{equation*}
whereupon, recalling that (cf. condition~(2.5)) $\lim_{\vert x \vert \to 0}
(\widetilde{V}(x)/\ln (x^{-2} \! + \! 1)) \! = \! +\infty$, it follows that
$g^{e}_{+}(z) \! + \! g^{e}_{-}(z) \! - \! \widetilde{V}(z) \! - \! \ell_{e}
\! + \! 2Q_{e} \! < \! 0$.

$\mathbf{(4)}$ For $z \! \in \! (-\infty,b_{0}^{e})$, one shows that
\begin{equation*}
g^{e}_{\pm}(z) \! = \! 2 \int_{J_{e}} \ln (\vert z \! - \! s \vert) \psi_{V}^{
e}(s) \, \md s \! \pm \! 2 \pi \mi \! - \! Q_{e} \! - \!
\begin{cases}
\ln \vert z \vert, &\text{$z \! > \! 0$,} \\
\ln \vert z \vert \! \pm \! \mi \pi, &\text{$z \! < \! 0$,}
\end{cases}
\end{equation*}
whence
\begin{equation*}
g^{e}_{+}(z) \! - \! g^{e}_{-}(z) \! = \! 4 \pi \mi \! + \!
\begin{cases}
0, &\text{$z \! > \! 0$,} \\
-2 \pi \mi, &\text{$z \! < \! 0$,}
\end{cases}
\end{equation*}
which shows that $g^{e}_{+}(z) \! - \! g^{e}_{-}(z)$ is pure imaginary, and
$\mi (g^{e}_{+}(z) \! - \! g^{e}_{-}(z))^{\prime} \! = \! 0$. Also,
\begin{equation*}
g^{e}_{+}(z) \! + \! g^{e}_{-}(z) \! - \! \widetilde{V}(z) \! - \! \ell_{e} \!
+ \! 2Q_{e} \! = \! 4 \int_{J_{e}} \ln (\vert z \! - \! s \vert) \psi_{V}^{e}
(s) \, \md s \! - \! 2 \ln \vert z \vert \! - \! \widetilde{V}(z) \! - \!
\ell_{e}:
\end{equation*}
proceeding as in the asymptotic analysis for case~$\mathbf{(3)}$ above, one
arrives at
\begin{equation*}
g^{e}_{+}(z) \! + \! g^{e}_{-}(z) \! - \! \widetilde{V}(z) \! - \! \ell_{e}
\! + \! 2Q_{e} \underset{z \to -\infty}{=} \ln (z^{2} \! + \! 1) \! - \!
\widetilde{V}(z) \! + \! \mathcal{O}(1),
\end{equation*}
and
\begin{equation*}
g^{e}_{+}(z) \! + \! g^{e}_{-}(z) \! - \! \widetilde{V}(z) \! - \! \ell_{e} \!
+ \! 2Q_{e} \underset{\vert z \vert \to 0}{=} \ln (z^{-2} \! + \! 1) \! - \!
\widetilde{V}(z) \! + \! \mathcal{O}(1),
\end{equation*}
whence, via conditions~(2.4) and~(2.5), $g^{e}_{+}(z) \! + \! g^{e}_{-}(z) \!
- \! \widetilde{V}(z) \! - \! \ell_{e} \! + \! 2Q_{e} \! < \! 0$, $z \! \in \!
(-\infty,b_{0}^{e})$. \hfill $\qed$
\section{The Model RHP and Parametrices}
In this section, the (normalised at infinity) auxiliary RHP for $\overset{e}{
\mathscr{M}} \colon \mathbb{C} \setminus \mathbb{R} \! \to \! 
\operatorname{SL}_{2}(\mathbb{C})$ formulated in Lemma~3.4 is augmented, 
by means of a sequence of contour deformations and transformations 
\emph{\`{a} la} Deift-Venakides-Zhou \cite{a1,a2,a3}, into simpler, `model' 
RHPs which, as $n \! \to \! \infty$, are solved explicitly (in closed form) 
in terms of Riemann theta functions (associated with the underlying genus-$N$ 
hyperelliptic Riemann surface) and Airy functions, and which give rise to the 
leading $(\mathcal{O}(1))$ terms of asymptotics for $\boldsymbol{\pi}_{2n}
(z)$, $\xi^{(2n)}_{n}$ and $\phi_{2n}(z)$ stated, respectively, in 
Theorems~2.3.1 and~2.3.2, and the asymptotic (as $n \! \to \! \infty)$ 
analysis of the parametrices, which are `approximate' solutions of the RHP 
for $\overset{e}{\mathscr{M}} \colon \mathbb{C} \setminus \mathbb{R} \! \to 
\! \operatorname{SL}_{2}(\mathbb{C})$ in neighbourhoods of the end-points of 
the support of the `even' equilibrium measure, and which give rise to the 
$\mathcal{O}(n^{-1})$ (and $\mathcal{O}(n^{-2}))$ corrections for 
$\boldsymbol{\pi}_{2n}(z)$, $\xi^{(2n)}_{n}$ and $\phi_{2n}(z)$ stated, 
respectively, in Theorems~2.3.1 and~2.3.2, is undertaken.
\begin{ccccc}
Let the external field $\widetilde{V} \colon \mathbb{R} \setminus \{0\} \! 
\to \! \mathbb{R}$ satisfy conditions~{\rm (2.3)--(2.5);} furthermore, let 
$\widetilde{V}$ be regular. Let the `even' equilibrium measure, $\mu_{V}^{e}$, 
and its support, $\operatorname{supp}(\mu_{V}^{e}) \! =: \! J_{e} \! = \! 
\cup_{j=1}^{N+1}J_{j}^{e} \! := \! \cup_{j=1}^{N+1}(b_{j-1}^{e},a_{j}^{e})$, 
be as described in Lemma~{\rm 3.5}, and, along with $\ell_{e}$ $(\in \! 
\mathbb{R})$, the `even' variational constant, satisfy the variational 
conditions stated in Lemma~{\rm 3.6}, Equations~{\rm (3.9);} moreover, let 
conditions~{\rm (1)}--{\rm (4)} stated in Lemma~{\rm 3.6} be valid. Then the 
{\rm RHP} for $\overset{e}{\mathscr{M}} \colon \mathbb{C} \setminus \mathbb{R} 
\! \to \! \mathrm{SL}_{2}(\mathbb{C})$ formulated in Lemma~{\rm 3.4} can be 
equivalently reformulated as follows: {\rm (1)} $\overset{e}{\mathscr{M}}(z)$ 
is holomorphic for $z \! \in \! \mathbb{C} \setminus \mathbb{R};$ {\rm (2)}
$\overset{e}{\mathscr{M}}_{\pm}(z) \! := \! \lim_{\underset{\pm \Im (z^{
\prime})>0}{z^{\prime} \to z}} \overset{e}{\mathscr{M}}(z^{\prime})$ satisfy 
the boundary condition
\begin{equation*}
\overset{e}{\mathscr{M}}_{+}(z) \! = \! \overset{e}{\mathscr{M}}_{-}(z)
\overset{e}{\upsilon}(z), \quad z \! \in \! \mathbb{R},
\end{equation*}
where
\begin{equation*}
\overset{e}{\upsilon}(z) \! = \!
\begin{cases}
\begin{pmatrix}
\me^{-4n \pi \mi \int_{z}^{a_{N+1}^{e}} \psi_{V}^{e}(s) \, \md s} & 1 \\
0 & \me^{4n \pi \mi \int_{z}^{a_{N+1}^{e}} \psi_{V}^{e}(s) \, \md s}
\end{pmatrix}, &\text{$z \! \in \! (b_{j-1}^{e},a_{j}^{e}), \quad j \! = \!
1,\dotsc,N \! + \! 1$,} \\
\begin{pmatrix}
\me^{-4n \pi \mi \int_{b_{k}^{e}}^{a_{N+1}^{e}} \psi_{V}^{e}(s) \, \md s} &
\me^{n(g^{e}_{+}(z)+g^{e}_{-}(z)-\widetilde{V}(z)-\ell_{e}+2Q_{e})} \\
0 & \me^{4n \pi \mi \int_{b_{k}^{e}}^{a_{N+1}^{e}} \psi_{V}^{e}(s) \, \md s}
\end{pmatrix}, &\text{$z \! \in \! (a_{k}^{e},b_{k}^{e}), \quad k \! = \! 1,
\dotsc,N$,} \\
\mathrm{I} \! + \! \me^{n(g^{e}_{+}(z)+g^{e}_{-}(z)-\widetilde{V}(z)-\ell_{e}+
2Q_{e})} \sigma_{+}, &\text{$z \! \in \! (-\infty,b_{0}^{e}) \! \cup \! (a_{N+
1}^{e},+\infty)$,}
\end{cases}
\end{equation*}
with $g^{e}(z)$ and $Q_{e}$ defined in Lemma~{\rm 3.4},
\begin{equation*}
\pm \Re \! \left(\mi \int_{z}^{a_{N+1}^{e}} \psi_{V}^{e}(s) \, \md s \right)
\! > \! 0, \quad z \! \in \! \mathbb{C}_{\pm} \cap (\cup_{j=1}^{N+1} \mathbb{
U}_{j}^{e}),
\end{equation*}
where $\mathbb{U}_{j}^{e} \! := \! \lbrace \mathstrut z \! \in \! \mathbb{C}^{
\ast}; \, \Re (z) \! \in \! (b_{j-1}^{e},a_{j}^{e}), \, \inf_{q \in J_{j}^{e}}
\vert z \! - \! q \vert \! < \! r_{j} \! \in \! (0,1) \rbrace$, $j \! = \! 1,
\dotsc,N \! + \! 1$, with $\mathbb{U}_{i}^{e} \cap \mathbb{U}_{j}^{e} \! = \!
\varnothing$, $i \! \not= \! j \! = \! 1,\dotsc,N \! + \! 1$, and $g^{e}_{+}
(z) \! + \! g^{e}_{-}(z) \! - \! \widetilde{V}(z) \! - \! \ell_{e} \! + \! 2
Q_{e} \! < \! 0$, $z \! \in \! (-\infty,b_{0}^{e}) \cup (a_{N+1}^{e},+\infty)
\cup (\cup_{j=1}^{N}(a_{j}^{e},b_{j}^{e}));$ {\rm (3)} $\overset{e}{\mathscr{
M}}(z) \! =_{\underset{z \in \mathbb{C} \setminus \mathbb{R}}{z \to \infty}}
\! \mathrm{I} \! + \! \mathcal{O}(z^{-1});$ and {\rm (4)} $\overset{e}{
\mathscr{M}}(z) \! =_{\underset{z \in \mathbb{C} \setminus \mathbb{R}}{z \to
0}} \! \mathcal{O}(1)$.
\end{ccccc}

\emph{Proof.} Item~(1) stated in the Lemma is simply a re-statement of
item~(1) of Lemma~3.4. Write $\mathbb{R} \! = \! (-\infty,b_{0}^{e}) \cup
(a_{N+1}^{e},+\infty) \cup (\cup_{j=1}^{N+1}J_{j}^{e}) \cup (\cup_{k=1}^{N}
(a_{k}^{e},b_{k}^{e})) \cup (\cup_{l=1}^{N+1}\{b_{l-1}^{e},a_{l}^{e}\})$,
where $J_{j}^{e} \! := \! (b_{j-1}^{e},a_{j}^{e})$, $j \! = \! 1,\dotsc,N \!
+ \! 1$. Recall {}from the proof of Lemma~3.6 that, for $\widetilde{V}$,
$\mu_{V}^{e}$, and $\ell_{e}$ described therein (and in the Lemma): (1)
\begin{equation*}
g^{e}_{+}(z) \! - \! g^{e}_{-}(z) \! = \!
\begin{cases}
4 \pi \mi \int_{z}^{a_{N+1}^{e}} \psi_{V}^{e}(s) \, \md s \! + \!
\begin{cases}
0, &\text{$z \! \in \! \mathbb{R}_{+} \cap \overline{J_{j}^{e}}, \quad j \! =
\! 1,\dotsc,N \! + \! 1$,} \\
-2 \pi \mi, &\text{$z \! \in \! \mathbb{R}_{-} \cap \overline{J_{j}^{e}},
\quad j \! = \! 1,\dotsc,N \! + \! 1$,}
\end{cases} \\
4 \pi \mi \int_{b_{j}^{e}}^{a_{N+1}^{e}} \psi_{V}^{e}(s) \, \md s \! + \!
\begin{cases}
0, &\text{$z \! \in \! \mathbb{R}_{+} \cap (a_{j}^{e},b_{j}^{e}), \quad j \!
= \! 1,\dotsc,N$,} \\
-2 \pi \mi, &\text{$z \! \in \! \mathbb{R}_{-} \cap (a_{j}^{e},b_{j}^{e}),
\quad j \! = \! 1,\dotsc,N$,}
\end{cases} \\
\begin{cases}
0, &\text{$z \! \in \! \mathbb{R}_{+} \cap (a_{N+1}^{e},+\infty)$,} \\
-2 \pi \mi, &\text{$z \! \in \! \mathbb{R}_{-} \cap (a_{N+1}^{e},+\infty)$,}
\end{cases} \\
4 \pi \mi \! + \!
\begin{cases}
0, &\text{$z \! \in \! \mathbb{R}_{+} \cap (-\infty,b_{0}^{e})$,} \\
-2 \pi \mi, &\text{$z \! \in \! \mathbb{R}_{-} \cap (-\infty,b_{0}^{e})$,}
\end{cases}
\end{cases}
\end{equation*}
where $\overline{J_{j}^{e}}$ $(:= \! J_{j}^{e} \cup \partial J_{j}^{e}) \! =
\! [b_{j-1}^{e},a_{j}^{e}]$, $j \! = \! 1,\dotsc,N \! + \! 1$; and (2)
\begin{equation*}
g^{e}_{+}(z) \! + \! g^{e}_{-}(z) \! - \! \widetilde{V}(z) \! - \! \ell_{e} \!
+ \! 2Q_{e} \! = \!
\begin{cases}
0, &\text{$z \! \in \! \cup_{j=1}^{N+1} \overline{J_{j}^{e}}$,} \\
-2 \int_{a_{j}^{e}}^{z}(R_{e}(s))^{1/2}h_{V}^{e}(s) \, \md s \! < \! 0,
&\text{$z \! \in \! (a_{j}^{e},b_{j}^{e}), \quad j \! = \! 1,\dotsc,N$,} \\
-2 \int_{a_{N+1}^{e}}^{z}(R_{e}(s))^{1/2}h_{V}^{e}(s) \, \md s \! < \! 0,
&\text{$z \! \in \! (a_{N+1}^{e},+\infty)$,} \\
2 \int_{z}^{b_{0}^{e}}(R_{e}(s))^{1/2}h_{V}^{e}(s) \, \md s \! < \! 0, &\text{
$z \! \in \! (-\infty,b_{0}^{e})$.}
\end{cases}
\end{equation*}
Recall, also, the formula for the `jump matrix' given in Lemma~3.4, namely,
\begin{equation*}
\begin{pmatrix}
\me^{-n(g^{e}_{+}(z)-g^{e}_{-}(z))} & \me^{n(g^{e}_{+}(z)+g^{e}_{-}(z)-
\widetilde{V}(z)-\ell_{e}+2Q_{e})} \\
0 & \me^{n(g^{e}_{+}(z)-g^{e}_{-}(z))}
\end{pmatrix}.
\end{equation*}
Partitioning $\mathbb{R}$ as given above, one obtains the formula for
$\overset{e}{\upsilon}(z)$ stated in the Lemma, thus item~(2); moreover,
items~(3) and~(4) are re-statements of the respective items of Lemma~3.4. It
remains, therefore, to show that $\Re (\mi \int_{z}^{a_{N+1}^{e}} \psi_{V}^{e}
(s) \, \md s)$ satisfies the inequalities stated in the Lemma. Recall {}from
the proof of Lemma~3.4 that $g^{e}(z)$ is uniformly Lipschitz continuous in
$\mathbb{C}_{\pm}$; moreover, via the Cauchy-Riemann conditions, item~(4) of
Lemma~3.6, that is, $\mi (g^{e}_{+}(z) \! - \! g^{e}_{-}(z))^{\prime} \!
\geqslant \! 0$, $z \! \in \! J_{e}$, implies that the quantity $g^{e}_{+}(z)
\! - \! g^{e}_{-}(z)$ has an analytic continuation, $\mathscr{G}^{e}(z)$, say,
to an open neighbourhood, $\mathbb{U}_{V}^{e}$, say, of $J_{e} \! = \! \cup_{
j=1}^{N+1}(b_{j-1}^{e},a_{j}^{e})$, where $\mathbb{U}_{V}^{e} \! := \! \cup_{
j=1}^{N+1} \mathbb{U}_{j}^{e}$, with $\mathbb{U}_{j}^{e} \! = \! \lbrace
\mathstrut z \! \in \! \mathbb{C}^{\ast}; \, \Re (z) \! \in \! (b_{j-1}^{e},
a_{j}^{e}), \, \inf_{q \in J_{j}^{e}} \vert z \! - \! q \vert \! < \! r_{j}
\! \in \! (0,1) \rbrace$, $j \! = \! 1,\dotsc,N \! + \! 1$, and $\mathbb{U}_{
i}^{e} \cap \mathbb{U}_{j}^{e} \! = \! \varnothing$, $i \! \not= \! j \! = \!
1,\dotsc,N \! + \! 1$, with the property that $\pm \Re (\mathscr{G}^{e}(z)) \!
> \! 0$ for $z \! \in \! \mathbb{C}_{\pm} \cap \mathbb{U}_{V}^{e}$. \hfill
$\qed$
\begin{eeeee}
Recalling that the external field $\widetilde{V} \colon \mathbb{R} \setminus
\{0\} \! \to \! \mathbb{R}$ is regular, that is, $h_{V}^{e}(z) \! \not\equiv
\! 0 \, \, \forall \, \, z \! \in \! \overline{J_{j}^{e}} \! := \! \cup_{j=1}^{
N+1}[b_{j-1}^{e},a_{j}^{e}]$, the second inequality of Equations~(3.9) is
strict, namely, $4 \int_{J_{e}} \ln (\vert x \! - \! s \vert) \psi_{V}^{e}(s)
\, \md s \! - \! 2 \ln \vert x \vert \! - \! \widetilde{V}(x) \! - \! \ell_{e}
\! + \! 2Q_{e} \! < \! 0$, $x \! \in \! \mathbb{R} \setminus \overline{J_{e}
}$, and ({}from the proof of Lemma~4.1) that $g^{e}_{+}(z) \! + \! g^{e}_{
-}(z) \! - \! \widetilde{V}(z) \! - \! \ell_{e} \! + \! 2Q_{e} \! < \! 0$, $z
\! \in \! (-\infty,b_{0}^{e}) \cup (a_{N+1}^{e},+\infty) \cup (\cup_{j=1}^{N}
(a_{j}^{e},b_{j}^{e}))$, it follows that
\begin{equation*}
\overset{e}{\upsilon}(z) \! \underset{n \to \infty}{=} \!
\begin{cases}
\me^{-(4n \pi \mi \int_{b_{j}^{e}}^{a_{N+1}^{e}} \psi_{V}^{e}(s) \, \md s)
\sigma_{3}} \! \left(\mathrm{I} \! + \! o(1) \sigma_{+} \right), &\text{$z \!
\in \! (a_{j}^{e},b_{j}^{e}), \quad j \! = \! 1,\dotsc,N$,} \\
\mathrm{I} \! + \! o(1) \sigma_{+}, &\text{$z \! \in \! (-\infty,b_{0}^{e})
\cup (a_{N+1}^{e},+\infty)$,}
\end{cases}
\end{equation*}
where $o(1)$ denotes terms that are exponentially small. \hfill $\blacksquare$
\end{eeeee}
\begin{ccccc}
Let the external field $\widetilde{V} \colon \mathbb{R} \setminus \{0\} \!
\to \! \mathbb{R}$ satisfy conditions~{\rm (2.3)--(2.5);} furthermore, let
$\widetilde{V}$ be regular. Let the `even' equilibrium measure, $\mu_{V}^{e}$,
and its support, $\operatorname{supp}(\mu_{V}^{e}) \! =: \! J_{e} \! = \!
\cup_{j=1}^{N+1}J_{j}^{e} \! := \! \cup_{j=1}^{N+1}(b_{j-1}^{e},a_{j}^{e})$,
be as described in Lemma~{\rm 3.5}, and, along with $\ell_{e}$ $(\in \!
\mathbb{R})$, the `even' variational constant, satisfy the variational
conditions stated in Lemma~{\rm 3.6}, Equations~{\rm (3.9);} moreover, let
conditions~{\rm (1)}--{\rm (4)} stated in Lemma~{\rm 3.6} be valid. Let
$\overset{e}{\mathscr{M}} \colon \mathbb{C} \setminus \mathbb{R} \! \to
\! \mathrm{SL}_{2}(\mathbb{C})$ solve the {\rm RHP} formulated in
Lemma~{\rm 4.1}, and let the deformed (and oriented) contour $\Sigma^{
\sharp}_{e} \! := \! \mathbb{R} \cup (\cup_{j=1}^{N+1}(J_{j}^{e,\smallfrown}
\cup J_{j}^{e,\smallsmile}))$ be as in Figure~{\rm 8} below; furthermore,
$\cup_{j=1}^{N+1}(\Omega_{j}^{e,\smallfrown} \cup \Omega_{j}^{e,\smallsmile}
\cup J_{j}^{e,\smallfrown} \cup J_{j}^{e,\smallsmile}) \subset \cup_{j=1}^{N
+1} \mathbb{U}_{j}^{e}$ (Figure~{\rm 8)}, where $\mathbb{U}_{j}^{e}$, $j \! =
\! 1,\dotsc,N \! + \! 1$, are defined in Lemma~{\rm 4.1}. Set
\begin{equation*}
\overset{e}{\mathscr{M}}^{\raise-1.0ex\hbox{$\scriptstyle \sharp$}}(z) \! :=
\!
\begin{cases}
\overset{e}{\mathscr{M}}(z), &\text{$z \! \in \! \mathbb{C} \setminus (\Sigma_{
e}^{\sharp} \cup (\cup_{j=1}^{N+1}(\Omega_{j}^{e,\smallfrown} \cup \Omega_{
j}^{e,\smallsmile})))$,} \\
\overset{e}{\mathscr{M}}(z) \! \left(\mathrm{I} \! - \! \me^{-4n \pi \mi \int_{
z}^{a_{N+1}^{e}} \psi_{V}^{e}(s) \, \md s} \, \sigma_{-} \right), &\text{$z \!
\in \! \mathbb{C}_{+} \cap (\cup_{j=1}^{N+1} \Omega_{j}^{e,\smallfrown})$,} \\
\overset{e}{\mathscr{M}}(z) \! \left(\mathrm{I} \! + \! \me^{4n \pi \mi \int_{
z}^{a_{N+1}^{e}} \psi_{V}^{e}(s) \, \md s} \, \sigma_{-} \right), &\text{$z \!
\in \! \mathbb{C}_{-} \cap (\cup_{j=1}^{N+1} \Omega_{j}^{e,\smallsmile})$.}
\end{cases}
\end{equation*}
Then $\overset{e}{\mathscr{M}}^{\raise-1.0ex\hbox{$\scriptstyle \sharp$}}
\colon \mathbb{C} \setminus \Sigma_{e}^{\sharp} \! \to \! \mathrm{SL}_{2}(
\mathbb{C})$ solves the following, equivalent {\rm RHP:} {\rm (1)} $\overset{
e}{\mathscr{M}}^{\raise-1.0ex\hbox{$\scriptstyle \sharp$}}(z)$ is holomorphic
for $z \! \in \! \mathbb{C} \setminus \Sigma_{e}^{\sharp};$ {\rm (2)}
$\overset{e}{\mathscr{M}}^{\raise-1.0ex\hbox{$\scriptstyle \sharp$}}_{\pm}
(z) \! := \! \lim_{\underset{z^{\prime} \, \in \, \pm \, \mathrm{side} \, 
\mathrm{of} \, \Sigma_{e}^{\sharp}}{z^{\prime} \to z}} 
\overset{e}{\mathscr{M}}^{\raise-1.0ex\hbox{$\scriptstyle \sharp$}}
(z^{\prime})$ satisfy the boundary condition
\begin{equation*}
\overset{e}{\mathscr{M}}^{\raise-1.0ex\hbox{$\scriptstyle \sharp$}}_{+}(z) \!
= \! \overset{e}{\mathscr{M}}^{\raise-1.0ex\hbox{$\scriptstyle \sharp$}}_{-}
(z) \overset{e}{\upsilon}^{\raise-1.0ex\hbox{$\scriptstyle \sharp$}}(z), \quad
z \! \in \! \textstyle \Sigma_{e}^{\sharp},
\end{equation*}
where
\begin{equation*}
\overset{e}{\upsilon}^{\raise-1.0ex\hbox{$\scriptstyle \sharp$}}(z) \! = \!
\begin{cases}
\mi \sigma_{2}, &\text{$z \! \in \! J_{j}^{e}, \quad j \! = \! 1,\dotsc,N \!
+ \! 1$,} \\
\mathrm{I} \! + \! \me^{-4n \pi \mi \int_{z}^{a_{N+1}^{e}} \psi_{V}^{e}(s) \,
\md s} \, \sigma_{-}, &\text{$z \! \in \! J_{j}^{e,\smallfrown}, \quad j \! =
\! 1,\dotsc,N \! + \! 1$,} \\
\mathrm{I} \! + \! \me^{4n \pi \mi \int_{z}^{a_{N+1}^{e}} \psi_{V}^{e}(s) \,
\md s} \, \sigma_{-}, &\text{$z \! \in \! J_{j}^{e,\smallsmile}, \quad j \! =
\! 1,\dotsc,N \! + \! 1$,} \\
\begin{pmatrix}
\me^{-4n \pi \mi \int_{b_{k}^{e}}^{a_{N+1}^{e}} \psi_{V}^{e}(s) \, \md s} &
\me^{n(g^{e}_{+}(z)+g^{e}_{-}(z)-\widetilde{V}(z)-\ell_{e}+2Q_{e})} \\
0 & \me^{4n \pi \mi \int_{b_{k}^{e}}^{a_{N+1}^{e}} \psi_{V}^{e}(s) \, \md s}
\end{pmatrix}, &\text{$z \! \in \! (a_{k}^{e},b_{k}^{e}), \quad k \! = \! 1,
\dotsc,N$,} \\
\mathrm{I} \! + \! \me^{n(g^{e}_{+}(z)+g^{e}_{-}(z)-\widetilde{V}(z)-\ell_{e}+
2Q_{e})} \sigma_{+}, &\text{$z \! \in \! (-\infty,b_{0}^{e}) \cup (a_{N+1}^{e},
+\infty)$,}
\end{cases}
\end{equation*}
with $\Re (\mi \int_{z}^{a_{N+1}^{e}} \psi_{V}^{e}(s) \, \md s) \! > \! 0$
(resp., $\Re (\mi \int_{z}^{a_{N+1}^{e}} \psi_{V}^{e}(s) \, \md s) \! < \!
0)$, $z \! \in \! \mathbb{C}_{+} \cap \Omega_{j}^{e,\smallfrown}$ (resp., $z
\! \in \! \mathbb{C}_{-} \cap \Omega_{j}^{e,\smallsmile})$, $j \! = \! 1,
\dotsc,N \! + \! 1;$ {\rm (3)}
\begin{equation*}
\overset{e}{\mathscr{M}}^{\raise-1.0ex\hbox{$\scriptstyle \sharp$}}(z)
\underset{\underset{z \in \mathbb{C} \setminus (\Sigma_{e}^{\sharp} \cup
(\cup_{j=1}^{N+1}(\Omega_{j}^{e,\smallfrown} \cup \Omega_{j}^{e,\smallsmile}
)))}{z \to \infty}}{=} \mathrm{I} \! + \! \mathcal{O}(z^{-1});
\end{equation*}
and {\rm (4)}
\begin{equation*}
\overset{e}{\mathscr{M}}^{\raise-1.0ex\hbox{$\scriptstyle \sharp$}}(z)
\underset{\underset{z \in \mathbb{C} \setminus (\Sigma_{e}^{\sharp} \cup
(\cup_{j=1}^{N+1}(\Omega_{j}^{e,\smallfrown} \cup \Omega_{j}^{e,\smallsmile})
))}{z \to 0}}{=} \mathcal{O}(1).
\end{equation*}
\end{ccccc}
\begin{figure}[tbh]
\begin{center}
\vspace{-0.40cm}
\begin{pspicture}(0,0)(14,5)
\psset{xunit=1cm,yunit=1cm}
\psarcn[linewidth=0.6pt,linestyle=solid,linecolor=magenta,arrowsize=1.5pt 5]%
{->}(2,1.5){1.8}{146}{90}
\psarcn[linewidth=0.6pt,linestyle=solid,linecolor=magenta](2,1.5){1.8}{90}{34}
\psarc[linewidth=0.6pt,linestyle=solid,linecolor=magenta,arrowsize=1.5pt 5]%
{->}(2,3.5){1.8}{214}{270}
\psarc[linewidth=0.6pt,linestyle=solid,linecolor=magenta](2,3.5){1.8}{270}{326}
\psline[linewidth=0.6pt,linestyle=solid,linecolor=black,arrowsize=1.5pt 4]%
{->}(0,2.5)(0.25,2.5)
\psline[linewidth=0.6pt,linestyle=solid,linecolor=black](0.25,2.5)(0.5,2.5)
\psline[linewidth=0.6pt,linestyle=solid,linecolor=magenta,arrowsize=1.5pt 5]%
{->}(0.5,2.5)(2,2.5)
\psline[linewidth=0.6pt,linestyle=solid,linecolor=magenta](2,2.5)(3.5,2.5)
\psline[linewidth=0.6pt,linestyle=solid,linecolor=black,arrowsize=1.5pt 4]%
{->}(3.5,2.5)(3.9,2.5)
\psline[linewidth=0.6pt,linestyle=solid,linecolor=black,arrowsize=1.5pt 4]%
{->}(5,2.5)(5.25,2.5)
\psline[linewidth=0.6pt,linestyle=solid,linecolor=black](5.25,2.5)(5.5,2.5)
\psline[linewidth=0.6pt,linestyle=solid,linecolor=magenta,arrowsize=1.5pt 5]%
{->}(5.5,2.5)(7,2.5)
\psline[linewidth=0.6pt,linestyle=solid,linecolor=magenta](7,2.5)(8.5,2.5)
\psline[linewidth=0.6pt,linestyle=solid,linecolor=black,arrowsize=1.5pt 4]%
{->}(8.5,2.5)(8.9,2.5)
\psarcn[linewidth=0.6pt,linestyle=solid,linecolor=magenta,arrowsize=1.5pt 5]%
{->}(7,1.5){1.8}{146}{90}
\psarcn[linewidth=0.6pt,linestyle=solid,linecolor=magenta](7,1.5){1.8}{90}{34}
\psarc[linewidth=0.6pt,linestyle=solid,linecolor=magenta,arrowsize=1.5pt 5]%
{->}(7,3.5){1.8}{214}{270}
\psarc[linewidth=0.6pt,linestyle=solid,linecolor=magenta](7,3.5){1.8}{270}{326}
\psline[linewidth=0.6pt,linestyle=solid,linecolor=black,arrowsize=1.5pt 4]%
{->}(10,2.5)(10.25,2.5)
\psline[linewidth=0.6pt,linestyle=solid,linecolor=black](10.25,2.5)(10.5,2.5)
\psline[linewidth=0.6pt,linestyle=solid,linecolor=magenta,arrowsize=1.5pt 5]%
{->}(10.5,2.5)(12,2.5)
\psline[linewidth=0.6pt,linestyle=solid,linecolor=magenta](12,2.5)(13.5,2.5)
\psline[linewidth=0.6pt,linestyle=solid,linecolor=black,arrowsize=1.5pt 4]%
{->}(13.5,2.5)(13.9,2.5)
\psarcn[linewidth=0.6pt,linestyle=solid,linecolor=magenta,arrowsize=1.5pt 5]%
{->}(12,1.5){1.8}{146}{90}
\psarcn[linewidth=0.6pt,linestyle=solid,linecolor=magenta](12,1.5){1.8}{90}{34}
\psarc[linewidth=0.6pt,linestyle=solid,linecolor=magenta,arrowsize=1.5pt 5]%
{->}(12,3.5){1.8}{214}{270}
\psarc[linewidth=0.6pt,linestyle=solid,linecolor=magenta](12,3.5){1.8}{270}%
{326}
\psline[linewidth=0.7pt,linestyle=dotted,linecolor=darkgray](3.95,2.5)%
(4.9,2.5)
\psline[linewidth=0.7pt,linestyle=dotted,linecolor=darkgray](8.95,2.5)%
(9.9,2.5)
\psdots[dotstyle=*,dotscale=1.5](0.5,2.5)
\psdots[dotstyle=*,dotscale=1.5](3.5,2.5)
\psdots[dotstyle=*,dotscale=1.5](5.5,2.5)
\psdots[dotstyle=*,dotscale=1.5](8.5,2.5)
\psdots[dotstyle=*,dotscale=1.5](10.5,2.5)
\psdots[dotstyle=*,dotscale=1.5](13.5,2.5)
\rput(0.5,2.1){\makebox(0,0){$\pmb{b_{0}^{e}}$}}
\rput(3.5,2.1){\makebox(0,0){$\pmb{a_{1}^{e}}$}}
\rput(5.5,2.1){\makebox(0,0){$\pmb{b_{j-1}^{e}}$}}
\rput(8.5,2.1){\makebox(0,0){$\pmb{a_{j}^{e}}$}}
\rput(10.5,2.1){\makebox(0,0){$\pmb{b_{N}^{e}}$}}
\rput(13.5,2.1){\makebox(0,0){$\pmb{a_{N+1}^{e}}$}}
\rput(2.7,2.5){\makebox(0,0){$\pmb{J_{1}^{e}}$}}
\rput(7.7,2.5){\makebox(0,0){$\pmb{J_{j}^{e}}$}}
\rput(12.7,2.5){\makebox(0,0){$\pmb{J_{N+1}^{e}}$}}
\rput(2.75,3.5){\makebox(0,0){$J_{1}^{e,\smallfrown}$}}
\rput(2.75,1.5){\makebox(0,0){$J_{1}^{e,\smallsmile}$}}
\rput(2,2.85){\makebox(0,0){$\Omega_{1}^{e,\smallfrown}$}}
\rput(2,2.15){\makebox(0,0){$\Omega_{1}^{e,\smallsmile}$}}
\rput(7.8,3.5){\makebox(0,0){$J_{j}^{e,\smallfrown}$}}
\rput(7.8,1.5){\makebox(0,0){$J_{j}^{e,\smallsmile}$}}
\rput(7,2.85){\makebox(0,0){$\Omega_{j}^{e,\smallfrown}$}}
\rput(7,2.15){\makebox(0,0){$\Omega_{j}^{e,\smallsmile}$}}
\rput(12.75,3.5){\makebox(0,0){$J_{N+1}^{e,\smallfrown}$}}
\rput(12.75,1.5){\makebox(0,0){$J_{N+1}^{e,\smallsmile}$}}
\rput(11.75,2.85){\makebox(0,0){$\Omega_{N+1}^{e,\smallfrown}$}}
\rput(11.75,2.15){\makebox(0,0){$\Omega_{N+1}^{e,\smallsmile}$}}
\end{pspicture}
\end{center}
\vspace{-1.00cm}
\caption{Oriented/deformed contour $\Sigma_{e}^{\sharp} \! := \! \mathbb{R}
\cup (\cup_{j=1}^{N+1}(J_{j}^{e,\smallfrown} \cup J_{j}^{e,\smallsmile}))$}
\end{figure}
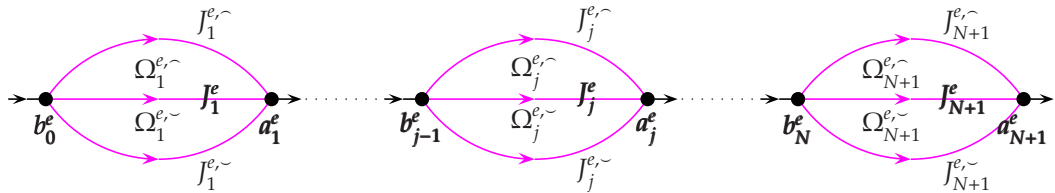

\emph{Proof.} Items~(1), (3), and~(4) in the formulation of the RHP for
$\overset{e}{\mathscr{M}}^{\raise-1.0ex\hbox{$\scriptstyle \sharp$}} \colon
\mathbb{C} \setminus \Sigma_{e}^{\sharp} \! \to \! \mathrm{SL}_{2}(\mathbb{
C})$ follow {}from the definition of
$\overset{e}{\mathscr{M}}^{\raise-1.0ex\hbox{$\scriptstyle \sharp$}}(z)$ (in
terms of $\overset{e}{\mathscr{M}}(z))$ given in the Lemma and the respective
items~(1), (3), and~(4) for the RHP for $\overset{e}{\mathscr{M}} \colon
\mathbb{C} \setminus \mathbb{R} \! \to \! \mathrm{SL}_{2}(\mathbb{C})$ stated
in Lemma~4.1; it remains, therefore, to verify item~(2), that is, the formula
for $\overset{e}{\upsilon}^{\raise-1.0ex\hbox{$\scriptstyle \sharp$}}(z)$.
Recall {}from item~(2) of Lemma~4.1 that, for $z \! \in \! (b_{j-1}^{e},a_{
j}^{e})$ $(\subset J_{e})$, $j \! = \! 1,\dotsc,N \! + \! 1$, $\overset{e}{
\mathscr{M}}_{+}(z) \! = \! \overset{e}{\mathscr{M}}_{-}(z) \overset{e}{
\upsilon}(z)$, where $\overset{e}{\upsilon}(z) \! = \!
\left(
\begin{smallmatrix}
\me^{-4n \pi \mi \int_{z}^{a_{N+1}^{e}} \psi_{V}^{e}(s) \, \md s} & 1 \\
0 & \me^{4n \pi \mi \int_{z}^{a_{N+1}^{e}} \psi_{V}^{e}(s) \, \md s}
\end{smallmatrix}
\right)$: noting the matrix factorisation
\begin{align*}
\begin{pmatrix}
\me^{-4n \pi \mi \int_{z}^{a_{N+1}^{e}} \psi_{V}^{e}(s) \, \md s} & 1 \\
0 & \me^{4n \pi \mi \int_{z}^{a_{N+1}^{e}} \psi_{V}^{e}(s) \, \md s}
\end{pmatrix} =& \,
\begin{pmatrix}
1 & 0 \\
\me^{4n \pi \mi \int_{z}^{a_{N+1}^{e}} \psi_{V}^{e}(s) \, \md s} & 1
\end{pmatrix} \!
\begin{pmatrix}
0 & 1 \\
-1 & 0
\end{pmatrix} \\
\times& \,
\begin{pmatrix}
1 & 0 \\
\me^{-4n \pi \mi \int_{z}^{a_{N+1}^{e}} \psi_{V}^{e}(s) \, \md s} & 1
\end{pmatrix},
\end{align*}
it follows that, for $z \! \in \! (b_{j-1}^{e},a_{j}^{e})$, $j \! = \! 1,
\dotsc,N \! + \! 1$,
\begin{equation*}
\overset{e}{\mathscr{M}}_{+}(z) \!
\begin{pmatrix}
1 & 0 \\
-\me^{-4n \pi \mi \int_{z}^{a_{N+1}^{e}} \psi_{V}^{e}(s) \, \md s} & 1
\end{pmatrix} \! = \overset{e}{\mathscr{M}}_{-}(z) \!
\begin{pmatrix}
1 & 0 \\
\me^{4n \pi \mi \int_{z}^{a_{N+1}^{e}} \psi_{V}^{e}(s) \, \md s} & 1
\end{pmatrix} \mi \sigma_{2}.
\end{equation*}
It was shown in Lemma~4.1 that $\pm \Re (\mi \int_{z}^{a_{N+1}^{e}} \psi_{V}^{
e}(s) \, \md s) \! > \! 0$ for $z \! \in \! \mathbb{C}_{\pm} \cap \mathbb{U}_{
j}^{e}$, where $\mathbb{U}_{j}^{e} \! := \! \lbrace \mathstrut z \! \in \!
\mathbb{C}^{\ast}; \, \Re (z) \! \in \! (b_{j-1}^{e},a_{j}^{e}), \, \inf_{q
\in J_{j}^{e}} \vert z \! - \! q \vert \! < \! r_{j} \! \in \! (0,1) \rbrace$,
$j \! = \! 1,\dotsc,N \! + \! 1$, with $\mathbb{U}_{i}^{e} \cap \mathbb{U}_{
j}^{e} \! = \! \varnothing$, $i \! \not= \! j \! = \! 1,\dotsc,N \! + \! 1$,
and $J_{j}^{e} \! := \! (b_{j-1}^{e},a_{j}^{e})$, $j \! = \! 1,\dotsc,N \! +
\! 1$. (One notes that the terms $\pm 4n \pi \mi \int_{z}^{a_{N+1}^{e}} \psi_{
V}^{e}(s) \, \md s$, which are pure imaginary for $z \! \in \! \mathbb{R}$,
and corresponding to which $\exp (\pm 4n \pi \mi \int_{z}^{a_{N+1}^{e}} \psi_{
V}^{e}(s) \, \md s)$ are undulatory, are continued analytically to $\mathbb{
C}_{\pm} \cap (\cup_{j=1}^{N+1} \mathbb{U}_{j}^{e})$, respectively,
corresponding to which $\exp (\pm 4n \pi \mi \int_{z}^{a_{N+1}^{e}} \psi_{V}^{
e}(s) \linebreak[4]
\md s)$ are exponentially decreasing as $n \! \to \! \infty)$. As
per the DZ non-linear steepest-descent method \cite{a1,a2} (see, also, the
extension \cite{a3}), one now `deforms' the original (and oriented) contour
$\mathbb{R}$ to the deformed, or extended, (and oriented) contour/skeleton
$\Sigma_{e}^{\sharp} \! := \! \mathbb{R} \cup (\cup_{j=1}^{N+1}(J_{j}^{e,
\smallfrown} \cup J_{j}^{e,\smallsmile}))$ (Figure~8) in such a way that the
upper (resp., lower) `lips' of the `lenses' $J_{j}^{e,\smallfrown}$ (resp.,
$J_{j}^{e,\smallsmile})$, $j \! = \! 1,\dotsc,N \! + \! 1$, which are the
boundaries of $\Omega_{j}^{e,\smallfrown}$ (resp., $\Omega_{j}^{e,\smallsmile}
)$, $j \! = \! 1,\dotsc,N \! + \! 1$, respectively, lie within the domain of
analytic continuation of $g^{e}_{+}(z) \! - \! g^{e}_{-}(z)$ (cf. the proof
of Lemma~4.1), that is, $\cup_{j=1}^{N+1}(\Omega_{j}^{e,\smallfrown} \cup
\Omega_{j}^{e,\smallsmile} \cup J_{j}^{e,\smallfrown} \cup J_{j}^{e,
\smallsmile}) \subset \cup_{j=1}^{N+1} \mathbb{U}_{j}^{e}$; in particular,
each (oriented) interval $J_{j}^{e} \! = \! (b_{j-1}^{e},a_{j}^{e})$, $j \! =
\! 1,\dotsc,N \! + \! 1$, in the original (and oriented) contour $\mathbb{R}$
is `split' (or branched) into three, and the new (and oriented) contour
$\Sigma_{e}^{\sharp}$ is the old contour $(\mathbb{R})$ together with the
(oriented) boundary of $N \! + \! 1$ lens-shaped regions, one region
surrounding each (bounded and oriented) interval $J_{j}^{e}$. Now,
recalling the definition of
$\overset{e}{\mathscr{M}}^{\raise-1.0ex\hbox{$\scriptstyle \sharp$}}(z)$
(in terms of $\overset{e}{\mathscr{M}}(z))$ stated in the Lemma, and the
expression for (the jump matrix) $\overset{e}{\upsilon}(z)$ given in
Lemma~4.1, one arrives at the formula for
$\overset{e}{\upsilon}^{\raise-1.0ex\hbox{$\scriptstyle \sharp$}}(z)$ given
in item~(2) of the Lemma. \hfill $\qed$
\begin{eeeee}
The jump condition stated in item~(2) of Lemma~4.2, that is,
$\overset{e}{\mathscr{M}}^{\raise-1.0ex\hbox{$\scriptstyle \sharp$}}_{+}(z) \!
= \! \overset{e}{\mathscr{M}}^{\raise-1.0ex\hbox{$\scriptstyle \sharp$}}_{-}
(z) \overset{e}{\upsilon}^{\raise-1.0ex\hbox{$\scriptstyle \sharp$}}(z)$,
$z \! \in \! \Sigma_{e}^{\sharp}$, with
$\overset{e}{\upsilon}^{\raise-1.0ex\hbox{$\scriptstyle \sharp$}}(z)$ given
therein, should, of course, be understood as follows: the $\mathrm{SL}_{2}
(\mathbb{C})$-valued functions
$\overset{e}{\mathscr{M}}^{\raise-1.0ex\hbox{$\scriptstyle \sharp$}} \! \!
\upharpoonright_{\mathbb{C}_{\pm} \setminus \Sigma_{e}^{\sharp}}$ have a
continuous extension to $\Sigma_{e}^{\sharp}$ with boundary values
$\overset{e}{\mathscr{M}}^{\raise-1.0ex\hbox{$\scriptstyle \sharp$}}_{\pm}
(z) \! := \! \lim_{\underset{z^{\prime} \, \in \, \pm \, \mathrm{side} \,
\mathrm{of} \, \Sigma_{e}^{\sharp}}{z^{\prime} \to z \in \Sigma_{e}^{\sharp}}}
\! \overset{e}{\mathscr{M}}^{\raise-1.0ex\hbox{$\scriptstyle \sharp$}}
(z^{\prime})$ satisfying the above jump relation
($\overset{e}{\mathscr{M}}^{\raise-1.0ex\hbox{$\scriptstyle \sharp$}}(z)$ is
continuous in each component of $\mathbb{C} \setminus \Sigma_{e}^{\sharp}$ up
to the boundary with boundary values
$\overset{e}{\mathscr{M}}^{\raise-1.0ex\hbox{$\scriptstyle \sharp$}}_{\pm}
(z)$ satisfying the above jump relation on $\Sigma_{e}^{\sharp})$. \hfill
$\blacksquare$
\end{eeeee}

Recalling {}from Lemma~4.1 that, for $z \! \in \! (-\infty,b_{0}^{e}) \cup
(a_{N+1}^{e},+\infty) \cup (\cup_{j=1}^{N}(a_{j}^{e},b_{j}^{e}))$, $g^{e}_{+}
(z) \! + \! g^{e}_{-}(z) \! - \! \widetilde{V}(z) \! - \! \ell_{e} \! + \! 2Q_{
e} \! < \! 0$, and, {}from Lemma~4.2, $\Re (\mi \int_{z}^{a_{N+1}^{e}} \psi_{
V}^{e}(s) \, \md s) \! > \! 0$ for $z \! \in \! J_{j}^{e,\smallfrown}$ (resp.,
$\Re (\mi \int_{z}^{a_{N+1}^{e}} \psi_{V}^{e}(s) \linebreak[4]
\md s) \! < \! 0$ for $z \! \in \! J_{j}^{e,\smallsmile})$, $j \! = \! 1,
\dotsc,N \! + \! 1$, one arrives at the following large-$n$ asymptotic
behaviour for the jump matrix
$\overset{e}{\upsilon}^{\raise-1.0ex\hbox{$\scriptstyle \sharp$}}(z)$:
\begin{equation*}
\overset{e}{\upsilon}^{\raise-1.0ex\hbox{$\scriptstyle \sharp$}}(z) \!
\underset{n \to \infty}{=} \!
\begin{cases}
\mi \sigma_{2}, &\text{$z \! \in \! J_{j}^{e}, \quad j \! = \! 1,\dotsc,N \! +
\! 1$,} \\
\mathrm{I} \! + \! \mathcal{O}(\me^{-nc \vert z \vert}) \sigma_{-}, &\text{$z
\! \in \! J_{j}^{e,\smallfrown} \cup J_{j}^{e,\smallsmile}, \quad j \! = \! 1,
\dotsc,N \! + \! 1$,} \\
\me^{-(4n \pi \mi \int_{b_{j}^{e}}^{a_{N+1}^{e}} \psi_{V}^{e}(s) \, \md s)
\sigma_{3}} \! \left(\mathrm{I} \! + \! \mathcal{O}(\me^{-nc \vert z-a_{j}^{
e} \vert}) \sigma_{+} \right), &\text{$z \! \in \! (a_{j}^{e},b_{j}^{e})
\setminus \widehat{\mathbb{U}}_{\delta_{0}^{e}}(0), \quad j \! = \! 1,\dotsc,
N$,} \\
\me^{-(4n \pi \mi \int_{b_{j}^{e}}^{a_{N+1}^{e}} \psi_{V}^{e}(s) \, \md s)
\sigma_{3}} \! \left(\mathrm{I} \! + \! \mathcal{O}(\me^{-nc \vert z \vert^{-
1}}) \sigma_{+} \right), &\text{$z \! \in \! (a_{j}^{e},b_{j}^{e}) \cap
\widehat{\mathbb{U}}_{\delta_{0}^{e}}(0), \quad j \! = \! 1,\dotsc,N$,} \\
\mathrm{I} \! + \! \mathcal{O}(\me^{-nc \vert z \vert}) \sigma_{+}, &\text{$z
\! \in \! ((-\infty,b_{0}^{e}) \cup (a_{N+1}^{e},+\infty)) \setminus \widehat{
\mathbb{U}}_{\delta_{0}^{e}}(0)$,} \\
\mathrm{I} \! + \! \mathcal{O}(\me^{-nc \vert z \vert^{-1}}) \sigma_{+},
&\text{$z \! \in \! ((-\infty,b_{0}^{e}) \cup (a_{N+1}^{e},+\infty)) \cap
\widehat{\mathbb{U}}_{\delta_{0}^{e}}(0)$,}
\end{cases}
\end{equation*}
where $c$ (some generic number) $> \! 0$, $\widehat{\mathbb{U}}_{\delta_{0}^{
e}}(0) \! := \! \lbrace \mathstrut z \! \in \! \mathbb{C}; \, \vert z \vert \!
< \! \delta_{0}^{e} \rbrace$, with $\delta_{0}^{e}$ some arbitrarily fixed,
sufficiently small positive real number, and where the respective convergences
are normal, that is, uniform in (respective) compact subsets (see Section~5
below).

Recall {}from Lemma~2.56 of \cite{a1} that, for an oriented skeleton in
$\mathbb{C}$ on which the jump matrix of an RHP is defined, one may always
choose to add or delete a portion of the skeleton on which the jump matrix
equals $\mathrm{I}$ without altering the RHP in the operator sense; hence,
neglecting those jumps on $\Sigma_{e}^{\sharp}$ tending exponentially quickly
(as $n \! \to \! \infty)$ to $\mathrm{I}$, and removing the corresponding
oriented skeletons {}from $\Sigma_{e}^{\sharp}$, it becomes more or less
transparent how to construct a parametrix, that is, an approximate solution,
of the (normalised at infinity) RHP for
$\overset{e}{\mathscr{M}}^{\raise-1.0ex\hbox{$\scriptstyle \sharp$}} \colon
\mathbb{C} \setminus \Sigma_{e}^{\sharp} \! \to \! \mathrm{SL}_{2}(\mathbb{C}
)$ stated in Lemma~4.2, namely, the large-$n$ solution of the RHP for
$\overset{e}{\mathscr{M}}^{\raise-1.0ex\hbox{$\scriptstyle \sharp$}}(z)$
formulated in Lemma~4.2 should be `close to' the solution of the
following (normalised at infinity) limiting, or model, RHP (for
$\overset{e}{m}^{\raise-1.0ex\hbox{$\scriptstyle \infty$}}(z))$.
\begin{ccccc}
Let the external field $\widetilde{V} \colon \mathbb{R} \setminus \{0\} \!
\to \! \mathbb{R}$ satisfy conditions~{\rm (2.3)--(2.5);} furthermore, let
$\widetilde{V}$ be regular. Let the `even' equilibrium measure, $\mu_{V}^{e}$,
and its support, $\operatorname{supp}(\mu_{V}^{e}) \! =: \! J_{e} \! = \!
\cup_{j=1}^{N+1}J_{j}^{e} \! := \! \cup_{j=1}^{N+1}(b_{j-1}^{e},a_{j}^{e})$,
be as described in Lemma~{\rm 3.5}, and, along with $\ell_{e}$ $(\in \!
\mathbb{R})$, the `even' variational constant, satisfy the variational
conditions stated in Lemma~{\rm 3.6}, Equations~{\rm (3.9);} moreover, let
conditions~{\rm (1)}--{\rm (4)} stated in Lemma~{\rm 3.6} be valid. Then
$\overset{e}{m}^{\raise-1.0ex\hbox{$\scriptstyle \infty$}} \colon \mathbb{C}
\setminus J_{e}^{\infty} \to \! \mathrm{SL}_{2}(\mathbb{C})$, where $J_{e}^{
\infty} \! := \! J_{e} \cup (\cup_{j=1}^{N}(a_{j}^{e},b_{j}^{e}))$, solves
the following (model) {\rm RHP:} {\rm (1)}
$\overset{e}{m}^{\raise-1.0ex\hbox{$\scriptstyle \infty$}}(z)$ is holomorphic
for $z \! \in \! \mathbb{C} \setminus J_{e}^{\infty};$ {\rm (2)}
$\overset{e}{m}^{\raise-1.0ex\hbox{$\scriptstyle \infty$}}_{\pm}(z) \! := \!
\lim_{\underset{z^{\prime} \, \in \, \pm \, \mathrm{side} \, \mathrm{of} \, 
J_{e}^{\infty}}{z^{\prime} \to z}} 
\overset{e}{m}^{\raise-1.0ex\hbox{$\scriptstyle \infty$}}(z^{\prime})$ 
satisfy the boundary condition
\begin{equation*}
\overset{e}{m}^{\raise-1.0ex\hbox{$\scriptstyle \infty$}}_{+}(z) \! = \!
\overset{e}{m}^{\raise-1.0ex\hbox{$\scriptstyle \infty$}}_{-}(z)
\overset{e}{\upsilon}^{\raise-1.0ex\hbox{$\scriptstyle \infty$}}(z), \quad z
\! \in \! J_{e}^{\infty},
\end{equation*}
where
\begin{equation*}
\overset{e}{\upsilon}^{\raise-1.0ex\hbox{$\scriptstyle \infty$}}(z) \! = \!
\begin{cases}
\mi \sigma_{2}, &\text{$z\! \in \! (b_{j-1}^{e},a_{j}^{e}), \quad j \! = \! 1,
\dotsc,N \! + \! 1$,} \\
\me^{-(4n \pi \mi \int_{b_{j}^{e}}^{a_{N+1}^{e}} \psi_{V}^{e}(s) \, \md s)
\sigma_{3}}, &\text{$z \! \in \! (a_{j}^{e},b_{j}^{e}), \quad j \! = \! 1,
\dotsc,N;$}
\end{cases}
\end{equation*}
{\rm (3)} $\overset{e}{m}^{\raise-1.0ex\hbox{$\scriptstyle \infty$}}(z) \! 
=_{\underset{z \in \mathbb{C} \setminus J_{e}^{\infty}}{z \to \infty}} \! 
\mathrm{I} \! + \! \mathcal{O}(z^{-1});$ and {\rm (4)} 
$\overset{e}{m}^{\raise-1.0ex\hbox{$\scriptstyle \infty$}}(z) \! =_{\underset{
z \in \mathbb{C} \setminus J_{e}^{\infty}}{z \to 0}} \! \mathcal{O}(1)$.
\end{ccccc}

The model RHP for $\overset{e}{m}^{\raise-1.0ex\hbox{$\scriptstyle \infty$}} 
\colon \mathbb{C} \setminus J_{e}^{\infty} \! \to \! \mathrm{SL}_{2}(\mathbb{
C})$ formulated in Lemma~4.3 is (explicitly) solvable in terms of Riemann 
theta functions (see, for example, Section~3 of \cite{a57}; see, also, 
Section~4.2 of \cite{a58}): the solution is succinctly presented below.
\begin{ccccc}
Let $\gamma^{e} \colon \mathbb{C} \setminus ((-\infty,b_{0}^{e}) \cup (a_{N+
1}^{e},+\infty) \cup (\cup_{j=1}^{N}(a_{j}^{e},b_{j}^{e}))) \! \to \! \mathbb{
C}$ be defined by
\begin{equation*}
\gamma^{e}(z) \! := \!
\begin{cases}
\left(\left(\dfrac{z \! - \! b_{0}^{e}}{z \! - \! a_{N+1}^{e}} \right)
\mathlarger{\prod_{k=1}^{N}} \! \left(\dfrac{z \! - \! b_{k}^{e}}{z \! - \!
a_{k}^{e}} \right) \right)^{1/4}, &\text{$z \! \in \! \mathbb{C}_{+}$,} \\
-\mi \left(\left(\dfrac{z \! - \! b_{0}^{e}}{z \! - \! a_{N+1}^{e}} \right)
\mathlarger{\prod_{k=1}^{N}} \! \left(\dfrac{z \! - \! b_{k}^{e}}{z \! - \!
a_{k}^{e}} \right) \right)^{1/4}, &\text{$z \! \in \! \mathbb{C}_{-}$.}
\end{cases}
\end{equation*}
Then, on the lower edge of each finite-length gap, that is, $(a_{j}^{e},b_{j}^{
e})^{-}$, $j \! = \! 1,\dotsc,N$, $\gamma^{e}(z) \! + \! (\gamma^{e}(z))^{-1}$
has exactly one root/zero, denoted $\left\lbrace z_{j}^{e,-} \! \in \! (a_{
j}^{e},b_{j}^{e})^{-} \! \subset \! \mathbb{C}_{-}, \, j \! = \! 1,\dotsc,N;
\, (\gamma^{e}(z) \! + \! (\gamma^{e}(z))^{-1}) \vert_{z=z_{j}^{e,-}} \! = \!
0 \right\rbrace$, and, on the upper edge of each finite-length gap, that is,
$(a_{j}^{e},b_{j}^{e})^{+}$, $j \! = \! 1,\dotsc,N$, $\gamma^{e}(z) \! - \!
(\gamma^{e}(z))^{-1}$ has exactly one root/zero, denoted $\left\lbrace z_{j}^{
e,+} \! \in \! (a_{j}^{e},b_{j}^{e})^{+} \! \subset \! \mathbb{C}_{+}, \, j \!
= \! 1,\dotsc,N; \, (\gamma^{e}(z) \! - \! (\gamma^{e}(z))^{-1}) \vert_{z=z_{
j}^{e,+}} \! = \! 0 \right\rbrace$ (in the plane, $z_{j}^{e,+} \! = \! z_{j}^{
e,-} \! := \! z_{j}^{e} \! \in \! (a_{j}^{e},b_{j}^{e})$, $j \! \in \! 1,
\dotsc,N)$. Furthermore, $\gamma^{e}(z)$ solves the following (scalar)
{\rm RHP:}
\begin{compactenum}
\item[{\rm (1)}] $\gamma^{e}(z)$ is holomorphic for $z \! \in \! \mathbb{C}
\setminus ((-\infty,b_{0}^{e}) \cup (a_{N+1}^{e},+\infty) \cup (\cup_{j=1}^{N}
(a_{j}^{e},b_{j}^{e})));$
\item[{\rm (2)}] $\gamma^{e}_{+}(z) \! = \! \gamma^{e}_{-}(z) \mi$, $z \! \in
\! (-\infty,b_{0}^{e}) \cup (a_{N+1}^{e},+\infty) \cup (\cup_{j=1}^{N}(a_{j}^{
e},b_{j}^{e}));$
\item[{\rm (3)}] $\gamma^{e}(z) \! =_{\underset{z \in \mathbb{C}_{\pm}}{z \to
\infty}} \! (-\mi)^{(1 \mp 1)/2}(1 \! + \! \mathcal{O}(z^{-1}));$ and
\item[{\rm (4)}] $\gamma^{e}(z) \! =_{\underset{z \in \mathbb{C}_{\pm}}{z \to
0}} \! \mathcal{O}(1)$.
\end{compactenum}
\end{ccccc}

\emph{Proof.} Define $\gamma^{e}(z)$ as in the Lemma. Then, one notes that 
$\gamma^{e}(z) \! \mp \! (\gamma^{e}(z))^{-1} \! = \! 0 \! \Leftrightarrow \! 
(\gamma^{e}(z))^{2} \! \mp \! 1 \! = \! 0 \! \Rightarrow \! (\gamma^{e}(z))^{
4} \! - \! 1 \! = \! 0 \! \Leftrightarrow \! \mathscr{Q}^{e}(z)$ $(\in \! 
\mathbb{R}[z])$ $:= \! (z \! - \! b_{0}^{e}) \prod_{k=1}^{N}(z \! - \! b_{k}^{
e}) \! - \! (z \! - \! a_{N+1}^{e}) \prod_{k=1}^{N}(z \! - \! a_{k}^{e}) \! = 
\! 0$, whence, via a straightforward calculation, one shows that $\mathscr{Q}^{
e}(a_{j}^{e}) \! = \! (-1)^{N-j+1} \widehat{\mathscr{Q}}^{e}_{a_{j}^{e}}$, $j 
\! = \! 1,\dotsc,N$, where $\widehat{\mathscr{Q}}^{e}_{a_{j}^{e}} \! := \! 
(b_{j}^{e} \! - \! a_{j}^{e})(a_{j}^{e} \! - \! b_{0}^{e}) \prod_{k=1}^{j-1}
(a_{j}^{e} \! - \! b_{k}^{e}) \prod_{l=j+1}^{N}(b_{l}^{e} \! - \! a_{j}^{e})$ 
$(> \! 0)$, and $\mathscr{Q}^{e}(b_{j}^{e}) \! = \! -(-1)^{N-j+1} \widehat{
\mathscr{Q}}^{e}_{b_{j}^{e}}$, $j \! = \! 1,\dotsc,N$, where $\widehat{
\mathscr{Q}}^{e}_{b_{j}^{e}} \! := \! (b_{j}^{e} \! - \! a_{j}^{e})(a_{N+1}^{
e} \! - \! b_{j}^{e}) \prod_{k=1}^{j-1}(b_{j}^{e} \! - \! a_{k}^{e}) \prod_{l
=j+1}^{N}(a_{l}^{e} \! - \! b_{j}^{e})$ $(> \! 0)$; thus, $\mathscr{Q}^{e}
(a_{j}^{e}) \mathscr{Q}^{e}(b_{j}^{e}) \! < \! 0$, $j \! = \! 1,\dotsc,N$, 
which shows that: (i) $\mathscr{Q}^{e}(z)$ has a root/zero, $z_{j}^{e}$, in 
each (open) interval $(a_{j}^{e},b_{j}^{e})$, $j \! = \! 1,\dotsc,N$; and (ii) 
since $\mathscr{Q}^{e}(z)$ is a unital polynomial with $\deg (\mathscr{Q}^{e}
(z)) \! = \! N$, $\{z_{j}^{e}\}_{j=1}^{N}$ are the only (simple) zeros/roots 
of $\mathscr{Q}^{e}(z)$. A straightforward analysis of the branch cuts shows 
that, for $z \! \in \! \cup_{j=1}^{N}(a_{j}^{e},b_{j}^{e})^{\pm}$, $\pm 
(\gamma^{e}(z))^{2} \! > \! 0$, whence $\lbrace z_{j}^{e,\pm} \rbrace_{j=1}^{
N} \! = \! \lbrace \mathstrut z^{\pm} \! \in \! (a_{j}^{e},b_{j}^{e})^{\pm} 
\subset \mathbb{C}_{\pm}, \, j \! = \! 1,\dotsc,N; \, (\gamma^{e}(z) \! \mp 
\! (\gamma^{e}(z))^{-1}) \vert_{z=z^{\pm}} \! = \! 0 \rbrace$. Setting 
$\widetilde{J}^{e} \! := \! (-\infty,b_{0}^{e}) \cup (a_{N+1}^{e},+\infty) 
\cup (\cup_{j=1}^{N}(a_{j}^{e},b_{j}^{e}))$, one shows that, upon performing 
a straig\-h\-t\-forward analysis of the branch cuts, $\gamma^{e}(z)$ solves 
the RHP $(\gamma^{e}(z),\mi,\widetilde{J}^{e})$ formulated in the Lemma. 
\hfill $\qed$
\begin{eeeee}
Recall {}from the proof of Lemma~4.4 that $\mathscr{Q}^{e}(z) \! := \! (z \!
- \! b_{0}^{e}) \prod_{k=1}^{N}(z \! - \! b_{k}^{e}) \! - \! (z \! - \! a_{N
+1}^{e}) \prod_{k=1}^{N}(z \linebreak[4]
- \! a_{k}^{e})$ $(\in \! \mathbb{R}[z])$, with $\deg
(\mathscr{Q}^{e}(z)) \! = \! N$; furthermore, $\mathscr{Q}^{e}(z_{j}^{e}) \! =
\! 0$, where $z_{j}^{e} \! \in \! (a_{j}^{e},b_{j}^{e})$, $j \! = \! 1,\dotsc,
N$. Writing $\mathscr{Q}^{e}(z) \! = \! \sum_{j=0}^{N} \mathfrak{q}_{j}^{e}z^{
j}$, where, in particular, $\mathfrak{q}_{0}^{e} \! = \! (-1)^{N}(\prod_{j=1}^{
N+1}a_{j}^{e} \! - \! \prod_{l=1}^{N+1}b_{l-1}^{e})$ and $\mathfrak{q}_{N}^{e}
\! = \! \sum_{j=1}^{N+1}(a_{j}^{e} \! - \! b_{j-1}^{e})$ $(\not= \! 0)$, one
uses a particular case of \emph{Gerschgorin's Circle Theorem} to arrive at the
following (upper) bound for the roots/zeros $z_{j}^{e}$, $j \! = \! 1,\dotsc,
N$: $\vert z_{j}^{e} \vert \! \leqslant \! \vert \mathfrak{q}_{N}^{e} \vert^{-
1} \sum_{l=0}^{N} \vert \mathfrak{q}_{l}^{e} \vert$, $j \! = \! 1,\dotsc,N$.
\hfill $\blacksquare$
\end{eeeee}

All of the notation/nomenclature used in Lemma~4.5 below has been defined at
the end of Subsection~2.1; the reader, therefore, is advised to peruse the
relevant notations(s), etc., before proceeding to Lemma~4.5. Let $\mathcal{
Y}_{e}$ denote the Riemann surface of $y^{2} \! = \! R_{e}(z) \! = \! \prod_{
k=1}^{N+1}(z \! - \! b_{k-1}^{e})(z \! - \! a_{k}^{e})$, where the
single-valued branch of the square root is chosen so that $z^{-(N+1)}(R_{e}
(z))^{1/2} \! \sim_{\underset{z \in \mathbb{C}_{\pm}}{z \to \infty}} \! \pm
1$. Let $\mathscr{P} \! := \! (y,z)$ denote a point on the Riemann surface
$\mathcal{Y}_{e}$ $(:= \! \lbrace \mathstrut (y,z); \, y^{2} \! = \! R_{e}(z)
\rbrace)$. The notation $\infty^{\pm}$ means: $\mathscr{P} \! \to \! \infty^{
\pm} \! \Leftrightarrow \! z \! \to \! \infty, y \! \sim \! \pm z^{N+1}$.
\begin{ccccc}
Let $\overset{e}{m}^{\raise-1.0ex\hbox{$\scriptstyle \infty$}} \colon \mathbb{
C} \setminus J_{e}^{\infty} \to \! \mathrm{SL}_{2}(\mathbb{C})$ solve the
{\rm RHP} formulated in Lemma~{\rm 4.3}. Then,
\begin{equation*}
\overset{e}{m}^{\raise-1.0ex\hbox{$\scriptstyle \infty$}}(z) \! = \!
\begin{cases}
\overset{e}{\mathfrak{M}}^{\raise-1.0ex\hbox{$\scriptstyle \infty$}}(z),
&\text{$z \! \in \! \mathbb{C}_{+}$,} \\
-\mi \, \overset{e}{\mathfrak{M}}^{\raise-1.0ex\hbox{$\scriptstyle \infty$}}
(z) \sigma_{2}, &\text{$z \in \! \mathbb{C}_{-}$,}
\end{cases}
\end{equation*}
where
\begin{equation*}
\overset{e}{\mathfrak{M}}^{\raise-1.0ex\hbox{$\scriptstyle \infty$}}(z) \! :=
\!
\begin{pmatrix}
\frac{\boldsymbol{\theta}^{e}(\boldsymbol{u}^{e}_{+}(\infty)+\boldsymbol{d}_{
e})}{\boldsymbol{\theta}^{e}(\boldsymbol{u}^{e}_{+}(\infty)-\frac{n}{2 \pi}
\boldsymbol{\Omega}^{e}+\boldsymbol{d}_{e})} & 0 \\
0 & \frac{\boldsymbol{\theta}^{e}(\boldsymbol{u}^{e}_{+}(\infty)+\boldsymbol{
d}_{e})}{\boldsymbol{\theta}^{e}(-\boldsymbol{u}^{e}_{+}(\infty)-\frac{n}{2
\pi} \boldsymbol{\Omega}^{e}-\boldsymbol{d}_{e})}
\end{pmatrix} \! \overset{e}{\boldsymbol{\Theta}}^{
\raise-1.0ex\hbox{$\scriptstyle \infty$}}(z),
\end{equation*}
and
\begin{equation*}
\overset{e}{\boldsymbol{\Theta}}^{\raise-1.0ex\hbox{$\scriptstyle \infty$}}(z)
\! = \!
\begin{pmatrix}
\frac{(\gamma^{e}(z)+(\gamma^{e}(z))^{-1})}{2} \frac{\boldsymbol{\theta}^{e}
(\boldsymbol{u}^{e}(z)-\frac{n}{2 \pi} \boldsymbol{\Omega}^{e}+\boldsymbol{d}_{
e})}{\boldsymbol{\theta}^{e}(\boldsymbol{u}^{e}(z)+\boldsymbol{d}_{e})} &
-\frac{(\gamma^{e}(z)-(\gamma^{e}(z))^{-1})}{2 \mi} \frac{\boldsymbol{\theta}^{
e}(-\boldsymbol{u}^{e}(z)-\frac{n}{2 \pi} \boldsymbol{\Omega}^{e}+\boldsymbol{
d}_{e})}{\boldsymbol{\theta}^{e}(-\boldsymbol{u}^{e}(z)+\boldsymbol{d}_{e})} \\
\frac{(\gamma^{e}(z)-(\gamma^{e}(z))^{-1})}{2 \mi} \frac{\boldsymbol{\theta}^{
e}(\boldsymbol{u}^{e}(z)-\frac{n}{2 \pi} \boldsymbol{\Omega}^{e}-\boldsymbol{
d}_{e})}{\boldsymbol{\theta}^{e}(\boldsymbol{u}^{e}(z)-\boldsymbol{d}_{e})} &
\frac{(\gamma^{e}(z)+(\gamma^{e}(z))^{-1})}{2} \frac{\boldsymbol{\theta}^{e}(-
\boldsymbol{u}^{e}(z)-\frac{n}{2 \pi} \boldsymbol{\Omega}^{e}-\boldsymbol{d}_{
e})}{\boldsymbol{\theta}^{e}(\boldsymbol{u}^{e}(z)+\boldsymbol{d}_{e})}
\end{pmatrix},
\end{equation*}
with $\gamma^{e}(z)$ characterised completely in Lemma~{\rm 4.4}, $\boldsymbol{
\Omega}^{e} \! := \! (\Omega_{1}^{e},\Omega_{2}^{e},\dotsc,\Omega_{N}^{e})^{
\mathrm{T}}$ $(\in \! \mathbb{R}^{N})$, where $\Omega_{j}^{e} \! = \! 4 \pi
\int_{b_{j}^{e}}^{a_{N+1}^{e}} \linebreak[4]
\psi_{V}^{e}(s) \, \md s$, $j \! = \! 1,\dotsc,N$, and ${}^{\mathrm{T}}$
denotes transposition, $\boldsymbol{d}_{e} \! \equiv \! -\sum_{j=1}^{N} \!
\int_{a_{j}^{e}}^{z_{j}^{e,-}} \! \boldsymbol{\omega}^{e}$ $(= \! \sum_{j=
1}^{N} \! \int_{a_{j}^{e}}^{z_{j}^{e,+}} \! \boldsymbol{\omega}^{e})$,
$\lbrace z_{j}^{e,\pm} \rbrace_{j=1}^{N}$ are characterised completely in
Lemma~{\rm 4.4}, $\boldsymbol{\omega}^{e}$ is the associated normalised basis
of holomorphic one-forms of $\mathcal{Y}_{e}$, $\boldsymbol{u}^{e}(z) \! := \!
\int_{a_{N+1}^{e}}^{z} \boldsymbol{\omega}^{e}$ $(\in \! \operatorname{Jac}
(\mathcal{Y}_{e}))$, and $\boldsymbol{u}^{e}_{+}(\infty) \! := \! \int_{a_{N+
1}^{e}}^{\infty^{+}} \boldsymbol{\omega}^{e}$ $(\infty^{+}$ is the point at
infinity in $\mathbb{C}_{+});$ furthermore, the solution is unique.
\end{ccccc}

\emph{Proof.} Let $\overset{e}{m}^{\raise-1.0ex\hbox{$\scriptstyle \infty$}}
\colon \mathbb{C} \setminus J_{e}^{\infty} \to \! \mathrm{SL}_{2}(\mathbb{C}
)$ solve the RHP formulated in Lemma~4.3, and define
$\overset{e}{m}^{\raise-1.0ex\hbox{$\scriptstyle \infty$}}(z)$, in terms of
$\overset{e}{\mathfrak{M}}^{\raise-1.0ex\hbox{$\scriptstyle \infty$}}(z)$,
as in the Lemma. A straightforward calculation shows that
$\overset{e}{\mathfrak{M}}^{\raise-1.0ex\hbox{$\scriptstyle \infty$}} \colon
\mathbb{C} \setminus \mathbb{R} \! \to \! \mathrm{SL}_{2}(\mathbb{C})$ solves
the following (normalised at infinity) `twisted' RHP: (i)
$\overset{e}{\mathfrak{M}}^{\raise-1.0ex\hbox{$\scriptstyle \infty$}}(z)$ is
holomorphic for $z \! \in \! \mathbb{C} \setminus \widetilde{J}^{e}$, where
$\widetilde{J}^{e} \! := \! (-\infty,b_{0}^{e}) \cup (a_{N+1}^{e},+\infty)
\cup (\cup_{j=1}^{N}(a_{j}^{e},b_{j}^{e}))$; (ii)
$\overset{e}{\mathfrak{M}}^{\raise-1.0ex\hbox{$\scriptstyle \infty$}}_{\pm}
(z) \! := \! \lim_{\underset{z^{\prime} \, \in \, \pm \, \mathrm{side} \, 
\mathrm{of} \, \widetilde{J}^{e}}{z^{\prime} \to z}} \! 
\overset{e}{\mathfrak{M}}^{\raise-1.0ex\hbox{$\scriptstyle \infty$}}
(z^{\prime})$ satisfy the boundary condition 
$\overset{e}{\mathfrak{M}}^{\raise-1.0ex\hbox{$\scriptstyle \infty$}}_{+}(z) 
\! = \! \overset{e}{\mathfrak{M}}^{\raise-1.0ex\hbox{$\scriptstyle \infty$}}_{
-}(z) \overset{e}{\mathscr{V}}^{\raise-1.0ex\hbox{$\scriptstyle \infty$}}(z)$,
$z \! \in \! \widetilde{J}^{e}$, where
\begin{equation}
\overset{e}{\mathscr{V}}^{\raise-1.0ex\hbox{$\scriptstyle \infty$}}(z) \! :=
\!
\begin{cases}
\mathrm{I}, &\text{$z \! \in \! J_{e}$,} \\
-\mi \sigma_{2}, &\text{$z \! \in \! (-\infty,b_{0}^{e}) \cup (a_{N+1}^{e},
+\infty)$,} \\
-\mi \sigma_{2} \me^{-\mi n \Omega_{j}^{e} \sigma_{3}}, &\text{$z \! \in \!
(a_{j}^{e},b_{j}^{e}), \quad j \! = \! 1,\dotsc,N$,}
\end{cases}
\end{equation}
with $\Omega_{j}^{e} \! = \! 4 \pi \int_{b_{j}^{e}}^{a_{N+1}^{e}} \psi_{V}^{
e}(s) \, \md s$, $j \! = \! 1,\dotsc,N$; (iii)
$\overset{e}{\mathfrak{M}}^{\raise-1.0ex\hbox{$\scriptstyle \infty$}}(z) \!
=_{\underset{z \in \mathbb{C}_{+}}{z \to \infty}} \! \mathrm{I} \! + \!
\mathcal{O}(z^{-1})$ and
$\overset{e}{\mathfrak{M}}^{\raise-1.0ex\hbox{$\scriptstyle \infty$}}(z) \!
=_{\underset{z \in \mathbb{C}_{-}}{z \to \infty}} \! \mi \sigma_{2} \! + \!
\mathcal{O}(z^{-1})$; and (iv)
$\overset{e}{\mathfrak{M}}^{\raise-1.0ex\hbox{$\scriptstyle \infty$}}(z) \!
=_{\underset{z \in \mathbb{C} \setminus \widetilde{J}^{e}}{z \to 0}} \!
\mathcal{O}(1)$. The solution of this latter (twisted) RHP for
$\overset{e}{\mathfrak{M}}^{\raise-1.0ex\hbox{$\scriptstyle \infty$}}(z)$ is
constructed out of the solution of two, simpler RHPs: $(\mathscr{N}^{e}(z),
-\mi \sigma_{2},\widetilde{J}^{e})$ and
$(\overset{e}{\mathfrak{m}}^{\raise-1.0ex\hbox{$\scriptstyle \infty$}}(z),
\overset{e}{\mathscr{U}}^{\raise-1.0ex\hbox{$\scriptstyle \infty$}}(z),
\widetilde{J}^{e})$, where
$\overset{e}{\mathscr{U}}^{\raise-1.0ex\hbox{$\scriptstyle \infty$}}(z)$
equals $\exp (\mi n \Omega_{j}^{e} \sigma_{3}) \sigma_{1}$ for $z \! \in \!
(a_{j}^{e},b_{j}^{e})$, $j \! = \! 1,\dotsc,N$, and equals $\mathrm{I}$
for $z \! \in \! (-\infty,b_{0}^{e}) \cup (a_{N+1}^{e},+\infty)$. The RHP
$(\mathscr{N}^{e}(z),-\mi \sigma_{2},\widetilde{J}^{e})$ is solved explicitly
by diagonalising the jump matrix, and thus reduced to two scalar RHPs
\cite{a2} (see, also, \cite{a57,a59,a90}): the solution is
\begin{equation*}
\mathscr{N}^{e}(z) \! = \!
\begin{pmatrix}
\frac{1}{2}(\gamma^{e}(z)+(\gamma^{e}(z))^{-1}) & -\frac{1}{2 \mi}(\gamma^{e}
(z)-(\gamma^{e}(z))^{-1}) \\
\frac{1}{2 \mi}(\gamma^{e}(z)-(\gamma^{e}(z))^{-1}) & \frac{1}{2}(\gamma^{e}
(z)+(\gamma^{e}(z))^{-1})
\end{pmatrix},
\end{equation*}
where $\gamma^{e} \colon \mathbb{C} \setminus \widetilde{J}^{e} \! \to \!
\mathbb{C}$ is characterised completely in Lemma~4.4; furthermore, $\mathscr{
N}^{e}(z)$ is piecewise holomorphic for $z \! \in \! \mathbb{C} \setminus
\widetilde{J}^{e}$, and $\mathscr{N}^{e}(z) \! =_{\underset{z \in \mathbb{
C}_{+}}{z \to \infty}} \! \mathrm{I} \! + \! \mathcal{O}(z^{-1})$ and
$\mathscr{N}^{e}(z) \! =_{\underset{z \in \mathbb{C}_{-}}{z \to \infty}} \!
\mi \sigma_{2} \! + \! \mathcal{O}(z^{-1})$\footnote{Note that, strictly
speaking, $\mathscr{N}^{e}(z)$, as given above, does not solve the RHP
$(\mathscr{N}^{e}(z),-\mi \sigma_{2},\widetilde{J}^{e})$ in the sense defined
heretofore, as $\mathscr{N}^{e} \! \! \upharpoonright_{\mathbb{C}_{\pm}}$ can
not be extended continuously to $\overline{\mathbb{C}}_{\pm}$; however,
$\mathscr{N}^{e}(\bm{\cdot} \! \pm \! \mi \varepsilon)$ converge in $\mathcal{
L}^{2}_{\mathrm{M}_{2}(\mathbb{C}),\mathrm{loc}}(\mathbb{R})$ as $\varepsilon
\! \downarrow \! 0$ to $\mathrm{SL}_{2}(\mathbb{C})$-valued functions
$\mathscr{N}^{e}(z)$ in $\mathcal{L}^{2}_{\mathrm{M}_{2}(\mathbb{C})}
(\widetilde{J}^{e})$ that satisfy $\mathscr{N}^{e}_{+}(z) \! = \! \mathscr{
N}^{e}_{-}(z)(-\mi \sigma_{2})$ $\mathrm{a.e.}$ on $\widetilde{J}^{e}$: one
then shows that $\mathscr{N}^{e}(z)$ is the unique solution of the RHP
$(\mathscr{N}^{e}(z),-\mi \sigma_{2},\widetilde{J}^{e})$, where the latter
boundary/jump condition is interpreted in the $\mathcal{L}^{2}_{\mathrm{M}_{2}
(\mathbb{C}),\mathrm{loc}}$ sense.}.

Consider, now, the functions $\boldsymbol{\theta}^{e}(\boldsymbol{u}^{e}(z) \!
\pm \! \boldsymbol{d}_{e})$, where $\boldsymbol{u}^{e}(z) \colon z \! \to \!
\operatorname{Jac}(\mathcal{Y}_{e})$, $z \! \mapsto \! \boldsymbol{u}^{e}(z)
\! := \! \int_{a_{N+1}^{e}}^{z} \boldsymbol{\omega}^{e}$, with $\boldsymbol{
\omega}^{e}$ the associated normalised basis of holomorphic one-forms of
$\mathcal{Y}_{e}$, $\boldsymbol{d}_{e} \! \equiv \! -\sum_{j=1}^{N} \int_{a_{
j}^{e}}^{z_{j}^{e,-}} \boldsymbol{\omega}^{e} \! = \! \sum_{j=1}^{N} \int_{a_{
j}^{e}}^{z_{j}^{e,+}} \boldsymbol{\omega}^{e}$, where $\equiv$ denotes
equivalence modulo the period lattice, and $\{z_{j}^{e,\pm}\}_{j=1}^{N}$ are
characterised completely in Lemma~4.4. {}From the general theory of theta
functions on Riemann surfaces (see, for example, \cite{a87,a88}), $\boldsymbol{
\theta}^{e}(\boldsymbol{u}^{e}(z) \! + \! \boldsymbol{d}_{e})$, for $z \! \in
\! \mathcal{Y}_{e} \! := \! \lbrace \mathstrut (y,z); \, y^{2} \! = \! \prod_{
k=1}^{N+1}(z \! - \! b_{k-1}^{e})(z \! - \! a_{k}^{e}) \rbrace$, is either
identically zero on $\mathcal{Y}_{e}$ or has precisely $N$ (simple) zeros (the
generic case). In this case, since the divisors $\prod_{j=1}^{N}z_{j}^{e,-}$
and $\prod_{j=1}^{N}z_{j}^{e,+}$ are non-special, one uses Lemma~3.27 of 
\cite{a57} (see, also, Lemma~4.2 of \cite{a58}) and the representation 
\cite{a88} $\boldsymbol{K}_{e} \! = \! \sum_{j=1}^{N} \int_{a_{j}^{e}}^{a_{N+
1}^{e}} \boldsymbol{\omega}^{e}$, for the `even' vector of Riemann constants,
with $2 \boldsymbol{K}_{e} \! = \! 0$ and $s \boldsymbol{K}_{e} \! \not= \!
0$, $0 \! < \! s \! < \! 2$, to arrive at
\begin{align*}
\boldsymbol{\theta}^{e}(\boldsymbol{u}^{e}(z) \! + \! \boldsymbol{d}_{e}) \!
=& \, \boldsymbol{\theta}^{e} \! \left(\boldsymbol{u}^{e}(z) \! - \! \sum_{j=
1}^{N} \int_{a_{j}^{e}}^{z_{j}^{e,-}} \boldsymbol{\omega}^{e} \right) \! =
\boldsymbol{\theta}^{e} \! \left(\int_{a_{N+1}^{e}}^{z} \boldsymbol{\omega}^{
e} \! - \! \boldsymbol{K}_{e} \! - \! \sum_{j=1}^{N} \int_{a_{N+1}^{e}}^{z_{
j}^{e,-}} \boldsymbol{\omega}^{e} \right) \! = \! 0 \\
\Leftrightarrow& \, z \! \in \! \left\{z_{1}^{e,-},z_{2}^{e,-},\dotsc,
z_{N}^{e,-} \right\}, \\
\boldsymbol{\theta}^{e}(\boldsymbol{u}^{e}(z) \! - \! \boldsymbol{d}_{e}) \!
=& \, \boldsymbol{\theta}^{e} \! \left(\boldsymbol{u}^{e}(z) \! - \! \sum_{j=
1}^{N} \int_{a_{j}^{e}}^{z_{j}^{e,+}} \boldsymbol{\omega}^{e} \right) \! =
\boldsymbol{\theta}^{e} \! \left(\int_{a_{N+1}^{e}}^{z} \boldsymbol{\omega}^{
e} \! - \! \boldsymbol{K}_{e} \! - \! \sum_{j=1}^{N} \int_{a_{N+1}^{e}}^{z_{
j}^{e,+}} \boldsymbol{\omega}^{e} \right) \! = \! 0 \\
\Leftrightarrow& \, z \! \in \! \left\{z_{1}^{e,+},z_{2}^{e,+},\dotsc,
z_{N}^{e,+} \right\}.
\end{align*}
Following Lemma~3.21 of \cite{a57}, set
\begin{equation*}
\overset{e}{\mathfrak{m}}^{\raise-1.0ex\hbox{$\scriptstyle \infty$}}(z) \! :=
\!
\begin{pmatrix}
\frac{\boldsymbol{\theta}^{e}(\boldsymbol{u}^{e}(z)-\frac{n}{2 \pi}
\boldsymbol{\Omega}^{e}+\boldsymbol{d}_{e})}{\boldsymbol{\theta}^{e}
(\boldsymbol{u}^{e}(z)+\boldsymbol{d}_{e})} & \frac{\boldsymbol{\theta}^{e}
(-\boldsymbol{u}^{e}(z)-\frac{n}{2 \pi} \boldsymbol{\Omega}^{e}+\boldsymbol{
d}_{e})}{\boldsymbol{\theta}^{e}(-\boldsymbol{u}^{e}(z)+\boldsymbol{d}_{e})} \\
\frac{\boldsymbol{\theta}^{e}(\boldsymbol{u}^{e}(z)-\frac{n}{2 \pi}
\boldsymbol{\Omega}^{e}-\boldsymbol{d}_{e})}{\boldsymbol{\theta}^{e}
(\boldsymbol{u}^{e}(z)-\boldsymbol{d}_{e})} & \frac{\boldsymbol{\theta}^{e}
(-\boldsymbol{u}^{e}(z)-\frac{n}{2 \pi} \boldsymbol{\Omega}^{e}-\boldsymbol{
d}_{e})}{\boldsymbol{\theta}^{e}(\boldsymbol{u}^{e}(z)+\boldsymbol{d}_{e})}
\end{pmatrix},
\end{equation*}
where $\boldsymbol{\Omega}^{e} \! := \! (\Omega_{1}^{e},\Omega_{2}^{e},\dotsc,
\Omega_{N}^{e})^{\mathrm{T}}$ $(\in \! \mathbb{R}^{N})$, with $\Omega_{j}^{
e}$, $j \! = \! 1,\dotsc,N$, given above, and ${}^{\mathrm{T}}$ denoting 
transposition. Using Lemma~3.18 of \cite{a57} (or, equivalently, 
Equations~(4.65) and~(4.66) of \cite{a58}), that is, for $z \! \in \! (a_{j}^{
e},b_{j}^{e})$, $j \! = \! 1,\dotsc,N$, $\boldsymbol{u}^{e}_{+}(z) \! + \!
\boldsymbol{u}^{e}_{-}(z) \! \equiv \! -\tau^{e}_{j}$ $(:= \! -\tau^{e}e_{j}
)$, $j \! = \! 1,\dotsc,N$, with $\tau^{e} \! := \! (\tau^{e})_{i,j=1,\dotsc,
N} \! := \! (\oint_{\boldsymbol{\beta}_{j}^{e}} \omega^{e}_{i})_{i,j=1,\dotsc,
N}$ (the associated matrix of Riemann periods), and, for $z \! \in \! (-\infty,
b_{0}^{e}) \cup (a_{N+1}^{e},+\infty)$, $\boldsymbol{u}^{e}_{+}(z) \! + \!
\boldsymbol{u}^{e}_{-}(z) \! \equiv \! 0$, where $\boldsymbol{u}^{e}_{\pm}(z)
\! := \! \int_{a_{N+1}^{e}}^{z^{\pm}} \boldsymbol{\omega}^{e}$, with $z^{\pm}
\! \in \! (a_{j}^{e},b_{j}^{e})^{\pm}$, $j \! = \! 1,\dotsc,N$, and the
evenness and (quasi-) periodicity properties of $\boldsymbol{\theta}^{e}(z)$,
one shows that, for $z \! \in \! (a_{j}^{e},b_{j}^{e})$, $j \! = \! 1,\dotsc,
N$,
\begin{gather*}
\dfrac{\boldsymbol{\theta}^{e}(\boldsymbol{u}^{e}_{+}(z) \! - \! \frac{n}{2
\pi} \boldsymbol{\Omega}^{e} \! + \! \boldsymbol{d}_{e})}{\boldsymbol{\theta}^{
e}(\boldsymbol{u}^{e}_{+}(z) \! + \! \boldsymbol{d}_{e})} \! = \! \me^{-\mi n
\Omega_{j}^{e}} \dfrac{\boldsymbol{\theta}^{e}(-\boldsymbol{u}^{e}_{-}(z) \! -
\! \frac{n}{2 \pi} \boldsymbol{\Omega}^{e} \! + \! \boldsymbol{d}_{e})}{
\boldsymbol{\theta}^{e}(-\boldsymbol{u}^{e}_{-}(z) \! + \! \boldsymbol{d}_{
e})}, \\
\dfrac{\boldsymbol{\theta}^{e}(\boldsymbol{u}^{e}_{+}(z) \! - \! \frac{n}{2
\pi} \boldsymbol{\Omega}^{e} \! - \! \boldsymbol{d}_{e})}{\boldsymbol{\theta}^{
e}(\boldsymbol{u}^{e}_{+}(z) \! - \! \boldsymbol{d}_{e})} \! = \! \me^{-\mi n
\Omega_{j}^{e}} \dfrac{\boldsymbol{\theta}^{e}(-\boldsymbol{u}^{e}_{-}(z) \! -
\! \frac{n}{2 \pi} \boldsymbol{\Omega}^{e} \! - \! \boldsymbol{d}_{e})}{
\boldsymbol{\theta}^{e}(\boldsymbol{u}^{e}_{-}(z) \! + \! \boldsymbol{d}_{
e})}, \\
\dfrac{\boldsymbol{\theta}^{e}(-\boldsymbol{u}^{e}_{+}(z) \! - \! \frac{n}{2
\pi} \boldsymbol{\Omega}^{e} \! + \! \boldsymbol{d}_{e})}{\boldsymbol{\theta}^{
e}(-\boldsymbol{u}^{e}_{+}(z) \! + \! \boldsymbol{d}_{e})} \! = \! \me^{\mi n
\Omega_{j}^{e}} \dfrac{\boldsymbol{\theta}^{e}(\boldsymbol{u}^{e}_{-}(z) \! -
\! \frac{n}{2 \pi} \boldsymbol{\Omega}^{e} \! + \! \boldsymbol{d}_{e})}{
\boldsymbol{\theta}^{e}(\boldsymbol{u}^{e}_{-}(z) \! + \! \boldsymbol{d}_{
e})}, \\
\dfrac{\boldsymbol{\theta}^{e}(-\boldsymbol{u}^{e}_{+}(z) \! - \! \frac{n}{2
\pi} \boldsymbol{\Omega}^{e} \! - \! \boldsymbol{d}_{e})}{\boldsymbol{\theta}^{
e}(\boldsymbol{u}^{e}_{+}(z) \! + \! \boldsymbol{d}_{e})} \! = \! \me^{\mi n
\Omega_{j}^{e}} \dfrac{\boldsymbol{\theta}^{e}(\boldsymbol{u}^{e}_{-}(z) \! -
\! \frac{n}{2 \pi} \boldsymbol{\Omega}^{e} \! - \! \boldsymbol{d}_{e})}{
\boldsymbol{\theta}^{e}(\boldsymbol{u}^{e}_{-}(z) \! - \! \boldsymbol{d}_{e})},
\end{gather*}
and, for $z \! \in \! (-\infty,b_{0}^{e}) \cup (a_{N+1}^{e},+\infty)$, one
obtains the same relations as above but with the changes $\exp (\mp \mi n
\Omega_{j}^{e}) \! \to \! 1$. Set, as in Proposition~3.31 of \cite{a57},
\begin{equation*}
\overset{e}{\mathcal{Q}}^{\raise-1.0ex\hbox{$\scriptstyle \infty$}}(z) \! :=
\!
\begin{pmatrix}
(\mathscr{N}^{e}(z))_{11}
(\overset{e}{\mathfrak{m}}^{\raise-1.0ex\hbox{$\scriptstyle \infty$}}(z))_{11}
& (\mathscr{N}^{e}(z))_{12}
(\overset{e}{\mathfrak{m}}^{\raise-1.0ex\hbox{$\scriptstyle \infty$}}(z))_{12}
\\(\mathscr{N}^{e}(z))_{21}
(\overset{e}{\mathfrak{m}}^{\raise-1.0ex\hbox{$\scriptstyle \infty$}}(z))_{21}
& (\mathscr{N}^{e}(z))_{22}
(\overset{e}{\mathfrak{m}}^{\raise-1.0ex\hbox{$\scriptstyle \infty$}}(z))_{22}
\end{pmatrix},
\end{equation*}
where $(\ast)_{ij}$, $i,j \! = \! 1,2$, denotes the $(i \, j)$-element of
$(\ast)$. Recalling that $\mathscr{N}^{e} \colon \mathbb{C} \setminus
\widetilde{J}^{e} \! \to \! \mathrm{SL}_{2}(\mathbb{C})$ solves the RHP
$(\mathscr{N}^{e}(z),-\mi \sigma_{2},\widetilde{J}^{e})$, using the above
theta-functional relations and the large-$z$ asymptotic expansion of
$\boldsymbol{u}^{e}(z)$ (see Section~5, the proof of Proposition~5.3), one
shows that
$\overset{e}{\mathcal{Q}}^{\raise-1.0ex\hbox{$\scriptstyle \infty$}}(z)$
solves the following RHP: (i)
$\overset{e}{\mathcal{Q}}^{\raise-1.0ex\hbox{$\scriptstyle \infty$}}(z)$ is
holomorphic for $z \! \in \! \mathbb{C} \setminus \widetilde{J}^{e}$; (ii) 
$\overset{e}{\mathcal{Q}}^{\raise-1.0ex\hbox{$\scriptstyle \infty$}}_{\pm}(z)
\! := \! \lim_{\underset{z^{\prime} \, \in \, \pm \, \mathrm{side} \, \mathrm{
of} \, \widetilde{J}^{e}}{z^{\prime} \to z}} 
\overset{e}{\mathcal{Q}}^{\raise-1.0ex\hbox{$\scriptstyle \infty$}}(z^{\prime}
)$ satisfy the boundary condition 
$\overset{e}{\mathcal{Q}}^{\raise-1.0ex\hbox{$\scriptstyle \infty$}}_{+}(z) \!
= \! \overset{e}{\mathcal{Q}}^{\raise-1.0ex\hbox{$\scriptstyle \infty$}}_{-}
(z) \overset{e}{\mathscr{V}}^{\raise-1.0ex\hbox{$\scriptstyle \infty$}}(z)$,
$z \! \in \! \widetilde{J}^{e}$, where
$\overset{e}{\mathscr{V}}^{\raise-1.0ex\hbox{$\scriptstyle \infty$}}(z)$ is
defined in Equation~(4.1); (iii)
\begin{align*}
\overset{e}{\mathcal{Q}}^{\raise-1.0ex\hbox{$\scriptstyle \infty$}}(z)
\underset{\underset{z \in \mathbb{C}_{+}}{z \to \infty}}{=}&
\begin{pmatrix}
\frac{\boldsymbol{\theta}^{e}(\boldsymbol{u}^{e}_{+}(\infty)-\frac{n}{2 \pi}
\boldsymbol{\Omega}^{e}+\boldsymbol{d}_{e})}{\boldsymbol{\theta}^{e}
(\boldsymbol{u}^{e}_{+}(\infty)+\boldsymbol{d}_{e})} & 0 \\
0 & \frac{\boldsymbol{\theta}^{e}(-\boldsymbol{u}^{e}_{+}(\infty)-\frac{n}{2
\pi} \boldsymbol{\Omega}^{e}-\boldsymbol{d}_{e})}{\boldsymbol{\theta}^{e}
(\boldsymbol{u}^{e}_{+}(\infty)+\boldsymbol{d}_{e})}
\end{pmatrix} \! + \! \mathcal{O}(z^{-1}), \\
\overset{e}{\mathcal{Q}}^{\raise-1.0ex\hbox{$\scriptstyle \infty$}}(z)
\underset{\underset{z \in \mathbb{C}_{-}}{z \to \infty}}{=}&
\begin{pmatrix}
0 & \frac{\boldsymbol{\theta}^{e}(-\boldsymbol{u}^{e}_{-}(\infty)-\frac{n}{2
\pi} \boldsymbol{\Omega}^{e}+\boldsymbol{d}_{e})}{\boldsymbol{\theta}^{e}
(-\boldsymbol{u}^{e}_{-}(\infty)+\boldsymbol{d}_{e})} \\
-\frac{\boldsymbol{\theta}^{e}(\boldsymbol{u}^{e}_{-}(\infty)-\frac{n}{2 \pi}
\boldsymbol{\Omega}^{e}-\boldsymbol{d}_{e})}{\boldsymbol{\theta}^{e}
(\boldsymbol{u}^{e}_{-}(\infty)-\boldsymbol{d}_{e})} & 0
\end{pmatrix} \! + \! \mathcal{O}(z^{-1}),
\end{align*}
where $\boldsymbol{u}^{e}_{\pm}(\infty) \! := \! \int_{a_{N+1}^{e}}^{\infty^{
\pm}} \boldsymbol{\omega}^{e}$ $(\infty^{\pm}$, respectively, are the points
at infinity in $\mathbb{C}_{\pm})$; and (iv)
$\overset{e}{\mathcal{Q}}^{\raise-1.0ex\hbox{$\scriptstyle \infty$}}(z) \!
=_{\underset{z \in \mathbb{C} \setminus \widetilde{J}^{e}}{z \to 0}} \!
\mathcal{O}(1)$. Now, using the fact that $\boldsymbol{u}^{e}_{-}(\infty) \! =
\! \int_{a_{N+1}^{e}}^{\infty^{-}} \boldsymbol{\omega}^{e} \! = \! -\int_{a_{N
+1}^{e}}^{\infty^{+}} \boldsymbol{\omega}^{e} \! = \! -\boldsymbol{u}^{e}_{+}
(\infty)$, upon multiplying
$\overset{e}{\mathcal{Q}}^{\raise-1.0ex\hbox{$\scriptstyle \infty$}}(z)$ on
the left by
\begin{equation*}
\operatorname{diag} \! \left(\dfrac{\boldsymbol{\theta}^{e}(\boldsymbol{u}^{
e}_{+}(\infty) \! + \! \boldsymbol{d}_{e})}{\boldsymbol{\theta}^{e}
(\boldsymbol{u}^{e}_{+}(\infty) \! - \! \frac{n}{2 \pi} \boldsymbol{\Omega}^{
e} \! + \! \boldsymbol{d}_{e})},\dfrac{\boldsymbol{\theta}^{e}(\boldsymbol{
u}^{e}_{+}(\infty) \! + \! \boldsymbol{d}_{e})}{\boldsymbol{\theta}^{e}(-
\boldsymbol{u}^{e}_{+}(\infty) \! - \! \frac{n}{2 \pi} \boldsymbol{\Omega}^{e}
\! - \! \boldsymbol{d}_{e})} \right) \! =: \!
\overset{e}{\mathfrak{c}}^{\raise-1.0ex\hbox{$\scriptstyle \infty$}},
\end{equation*}
that is, $\overset{e}{\mathcal{Q}}^{\raise-1.0ex\hbox{$\scriptstyle \infty$}}
(z) \! \to \! \overset{e}{\mathfrak{c}}^{\raise-1.0ex\hbox{$\scriptstyle
\infty$}} \overset{e}{\mathcal{Q}}^{\raise-1.0ex\hbox{$\scriptstyle \infty$}}
(z) \! =: \! \mathcal{M}^{\infty}_{e}(z)$, one easily shows that $\mathcal{
M}^{\infty}_{e} \colon \mathbb{C} \setminus \widetilde{J}^{e} \! \to \!
\mathrm{SL}_{2}(\mathbb{C})$ solves the RHP $(\mathcal{M}^{\infty}_{e}(z),
\overset{e}{\mathscr{V}}^{\raise-1.0ex\hbox{$\scriptstyle \infty$}}(z),
\widetilde{J}^{e})$. Using, finally, the formula
$\overset{e}{m}^{\raise-1.0ex\hbox{$\scriptstyle \infty$}}(z) \! = \!
\begin{cases}
\overset{e}{\mathfrak{M}}^{\raise-1.0ex\hbox{$\scriptstyle \infty$}}(z),
&\text{$z \! \in \! \mathbb{C}_{+}$,} \\
-\mi \, \overset{e}{\mathfrak{M}}^{\raise-1.0ex\hbox{$\scriptstyle \infty$}}
(z) \sigma_{2}, &\text{$z \! \in \! \mathbb{C}_{-}$,}
\end{cases}$ one shows that
$\overset{e}{m}^{\raise-1.0ex\hbox{$\scriptstyle \infty$}} \colon \mathbb{C}
\setminus J_{e}^{\infty} \! \to \! \operatorname{SL}_{2}(\mathbb{C})$ solves
the model RHP formulated in Lemma~4.3. One notes {}from the formula for
$\overset{e}{\mathfrak{M}}^{\raise-1.0ex\hbox{$\scriptstyle \infty$}}(z)$
stated in the Lemma that it is well defined for $\mathbb{C} \setminus \mathbb{
R}$; in particular, it is single valued and analytic (see below) for $z \! \in
\! \mathbb{C} \setminus \widetilde{J}^{e}$ (independently of the path in
$\mathbb{C} \setminus \widetilde{J}^{e}$ chosen to evaluate $\boldsymbol{u}^{
e}(z) \! = \! \int_{a_{N+1}^{e}}^{z} \boldsymbol{\omega}^{e})$. Furthermore
(cf. Lemma~4.4 and the analysis above), since $\lbrace \mathstrut z^{\prime}
\! \in \! \mathbb{C}; \, \boldsymbol{\theta}^{e}(\boldsymbol{u}^{e}(z^{\prime}
) \! \pm \! \boldsymbol{d}_{e}) \! = \! 0 \rbrace \! = \! \lbrace z_{j}^{e,
\mp} \rbrace_{j=1}^{N} \! = \! \lbrace \mathstrut z^{\prime} \! \in \!
\mathbb{C}; \, (\gamma^{e}(z) \! \pm \! (\gamma^{e}(z))^{-1}) \vert_{z=z^{
\prime}} \! = \! 0 \rbrace$, one notes that the (simple) poles of
$(\overset{e}{\mathfrak{m}}^{\raise-1.0ex\hbox{$\scriptstyle \infty$}}(z))_{1
1}$ and $(\overset{e}{\mathfrak{m}}^{\raise-1.0ex\hbox{$\scriptstyle \infty$}}
(z))_{22}$ (resp.,
$(\overset{e}{\mathfrak{m}}^{\raise-1.0ex\hbox{$\scriptstyle \infty$}}(z))_{1
2}$ and
$(\overset{e}{\mathfrak{m}}^{\raise-1.0ex\hbox{$\scriptstyle \infty$}}(z))_{
21})$, that is, $\lbrace \mathstrut z^{\prime} \! \in \! \mathbb{C}; \,
\boldsymbol{\theta}^{e}(\boldsymbol{u}^{e}(z^{\prime}) \! + \! \boldsymbol{
d}_{e}) \! = \! 0 \rbrace$ (resp., $\lbrace \mathstrut z^{\prime} \! \in \!
\mathbb{C}; \, \boldsymbol{\theta}^{e}(\boldsymbol{u}^{e}(z^{\prime}) \! - \!
\boldsymbol{d}_{e}) \! = \! 0 \rbrace)$, are exactly cancelled by the (simple)
zeros of $\gamma^{e}(z) \! + \! (\gamma^{e}(z))^{-1}$ (resp., $\gamma^{e}(z)
\! - \! (\gamma^{e}(z))^{-1})$; thus,
$\overset{e}{\mathfrak{M}}^{\raise-1.0ex\hbox{$\scriptstyle \infty$}}(z)$ has
only $\tfrac{1}{4}$-root singularities at the end-points of the support of the
`even' equilibrium measure, $\lbrace b_{j-1}^{e},a_{j}^{e} \rbrace_{j=1}^{N+
1}$. (This shows that
$\overset{e}{\mathfrak{M}}^{\raise-1.0ex\hbox{$\scriptstyle \infty$}}(z)$
obtains its boundary values,
$\overset{e}{\mathfrak{M}}^{\raise-1.0ex\hbox{$\scriptstyle \infty$}}_{\pm}(z)
\! := \! \lim_{\varepsilon \downarrow 0}
\overset{e}{\mathfrak{M}}^{\raise-1.0ex\hbox{$\scriptstyle \infty$}}(z \! \pm
\! \mi \varepsilon)$, in the $\mathcal{L}^{2}_{\mathrm{M}_{2}(\mathbb{C})}
(\mathbb{R})$ sense.) {}From the definition of
$\overset{e}{m}^{\raise-1.0ex\hbox{$\scriptstyle \infty$}}(z)$ in terms of
$\overset{e}{\mathfrak{M}}^{\raise-1.0ex\hbox{$\scriptstyle \infty$}}(z)$
given in the Lemma, the explicit formula for
$\overset{e}{\mathfrak{M}}^{\raise-1.0ex\hbox{$\scriptstyle \infty$}}(z)$, and
recalling that $\overset{e}{m}^{\raise-1.0ex\hbox{$\scriptstyle \infty$}}(z)$
solves the model RHP formulated in Lemma~4.3, one learns that, as $\det
(\overset{e}{\upsilon}^{\raise-1.0ex\hbox{$\scriptstyle \infty$}}(z)) \! = \!
1$, $\det (\overset{e}{m}^{\raise-1.0ex\hbox{$\scriptstyle \infty$}}_{+}(z))
\! = \! \det (\overset{e}{m}^{\raise-1.0ex\hbox{$\scriptstyle \infty$}}_{-}
(z))$, that is,
$\det (\overset{e}{m}^{\raise-1.0ex\hbox{$\scriptstyle \infty$}}(z))$ has no
`jumps', whence
$\det (\overset{e}{m}^{\raise-1.0ex\hbox{$\scriptstyle \infty$}}(z))$ has, at
worst, (isolated) $\tfrac{1}{2}$-root singularities at $\lbrace b_{j-1}^{e},
a_{j}^{e} \rbrace_{j=1}^{N+1}$, which are removable, which implies that
$\det (\overset{e}{m}^{\raise-1.0ex\hbox{$\scriptstyle \infty$}}(z))$ is
entire and bounded; hence, via a generalisation of Liouville's Theorem,
and the asymptotic relation $\det
(\overset{e}{m}^{\raise-1.0ex\hbox{$\scriptstyle \infty$}}(z)) \! =_{
\underset{z \in \mathbb{C} \setminus \mathbb{R}}{z \to \infty}} \! 1 \! + \!
\mathcal{O}(z^{-1})$, one arrives at
$\det (\overset{e}{m}^{\raise-1.0ex\hbox{$\scriptstyle \infty$}}(z)) \! = \!
1 \! \Rightarrow \! \overset{e}{m}^{\raise-1.0ex\hbox{$\scriptstyle \infty$}}
(z) \! \in \! \operatorname{SL}_{2}(\mathbb{C})$. Also, {}from the definition
of $\overset{e}{m}^{\raise-1.0ex\hbox{$\scriptstyle \infty$}}(z)$ in terms of
$\overset{e}{\mathfrak{M}}^{\raise-1.0ex\hbox{$\scriptstyle \infty$}}(z)$ and
the explicit formula for
$\overset{e}{\mathfrak{M}}^{\raise-1.0ex\hbox{$\scriptstyle \infty$}}(z)$, it
follows that both $\overset{e}{m}^{\raise-1.0ex\hbox{$\scriptstyle \infty$}}
(z)$ and $(\overset{e}{m}^{\raise-1.0ex\hbox{$\scriptstyle \infty$}}(z))^{-1}$
are uniformly bounded as functions of $n$ (as $n \! \to \! \infty)$ for $z$ in
compact subsets away {}from $\lbrace b_{j-1}^{e},a_{j}^{e} \rbrace_{j=1}^{N+
1}$.

Let $\mathscr{S}^{\infty}_{e} \colon \mathbb{C} \setminus \mathbb{R} \! \to
\! \operatorname{SL}_{2}(\mathbb{C})$ be another solution of the RHP
$(\overset{e}{\mathfrak{M}}^{\raise-1.0ex\hbox{$\scriptstyle \infty$}}(z),
\overset{e}{\mathscr{V}}^{\raise-1.0ex\hbox{$\scriptstyle \infty$}}(z),
\mathbb{R})$ formulated at the beginning of the proof, and set, as in the
\emph{Remark on Proposition}~3.43 of \cite{a57}, $\Delta^{e}(z) \! := \!
\overset{e}{\mathfrak{M}}^{\raise-1.0ex\hbox{$\scriptstyle \infty$}}(z) \! -
\! \mathscr{S}^{\infty}_{e}(z)$, whence $\Delta^{e} \! =_{\underset{z \in
\mathbb{C} \setminus \mathbb{R}}{z \to \infty}} \! \mathcal{O}(z^{-1})$; thus,
by Cauchy's Theorem, $\int_{\mathrm{C}_{\epsilon}^{e}} \Delta^{e}(s)(\Delta^{
e}(\overline{s}))^{\dagger} \, \md s \! = \! 0$, where ${}^{\dagger}$ denotes
the Hermitean adjoint, and $\mathrm{C}_{\epsilon}^{e}$ is the (closed and
simple) counter-clockwise-oriented contour $\mathrm{C}_{\epsilon}^{e} \!
:= \! \mathrm{C}_{\epsilon}^{e,\mathbb{R}} \cup \mathrm{C}_{\epsilon}^{e,
\smallfrown}$, where $\mathrm{C}_{\epsilon}^{e,\mathbb{R}} \! := \! \lbrace
\mathstrut x \! + \! \mi \epsilon; \, -\epsilon^{-1} \! \leqslant \! x \!
\leqslant \! \epsilon^{-1} \rbrace$ and $\mathrm{C}_{\epsilon}^{e,\smallfrown}
\! := \! \lbrace \mathstrut \epsilon^{-1} \me^{\mi \theta}; \, \theta \! \in
\! [\delta (\epsilon),\pi \! - \! \delta (\epsilon)], \, \delta (\epsilon) \!
:= \! \operatorname{tan}^{-1}(\epsilon^{2}) \rbrace$, with $\epsilon$ some
arbitrarily fixed, sufficiently small positive real number, and the principal
branch of $\operatorname{tan}^{-1}(\cdot)$ is taken. Writing $0 \! = \! \int_{
\mathrm{C}_{\epsilon}^{e}} \Delta^{e}(s)(\Delta^{e}(\overline{s}))^{\dagger}
\, \md s \! = \! (\int_{\mathrm{C}_{\epsilon}^{e,\mathbb{R}}} \! + \! \int_{
\mathrm{C}_{\epsilon}^{e,\smallfrown}}) \Delta^{e}(s)(\Delta^{e}(\overline{s}
))^{\dagger} \md s$, letting $\epsilon \downarrow 0$, in which case, since
$\Delta^{e}(z)(\Delta^{e}(\overline{z}))^{\dagger} \! =_{\underset{z \in
\mathbb{C} \setminus \mathbb{R}}{z \to \infty}} \! \mathcal{O}(z^{-2})$, an
application of Jordan's Lemma gives $\int_{\mathrm{C}_{\epsilon}^{e,
\smallfrown}} \Delta^{e}(s)(\Delta^{e}(\overline{s}))^{\dagger} \, \md s \!
=_{\epsilon \downarrow 0} \! 0$, one gets that
\begin{align*}
0 =& \, \int_{-\infty}^{+\infty} \Delta^{e}_{+}(s)(\Delta^{e}_{-}(s))^{
\dagger} \, \md s \! = \! \int_{-\infty}^{b_{0}^{e}} \Delta_{-}^{e}(s)(-\mi
\sigma_{2})(\Delta^{e}_{-}(s))^{\dagger} \, \md s \! + \! \int_{a_{N+1}^{e}}^{
+\infty} \Delta_{-}^{e}(s)(-\mi \sigma_{2})(\Delta^{e}_{-}(s))^{\dagger} \,
\md s \\
+& \, \sum_{j=1}^{N} \int_{a_{j}^{e}}^{b_{j}^{e}} \Delta_{-}^{e}(s)(-\mi
\sigma_{2} \me^{-\mi n \Omega_{j}^{e} \sigma_{3}})(\Delta^{e}_{-}(s))^{
\dagger} \, \md s \! + \! \int_{J_{e}} \Delta^{e}(s)(\Delta^{e}(s))^{\dagger}
\, \md s:
\end{align*}
adding the above to its Hermitean adjoint, that is,
\begin{align*}
0 =& \, \int_{-\infty}^{b_{0}^{e}} \Delta_{-}^{e}(s)(\mi \sigma_{2})(\Delta^{
e}_{-}(s))^{\dagger} \, \md s \! + \! \int_{a_{N+1}^{e}}^{+\infty} \Delta_{-}^{
e}(s)(\mi \sigma_{2})(\Delta^{e}_{-}(s))^{\dagger} \, \md s \\
+& \, \sum_{j=1}^{N} \int_{a_{j}^{e}}^{b_{j}^{e}} \Delta_{-}^{e}(s)(\mi \me^{
\mi n \Omega_{j}^{e} \sigma_{3}} \sigma_{2})(\Delta^{e}_{-}(s))^{\dagger} \,
\md s \! + \! \int_{J_{e}} \Delta^{e}(s)(\Delta^{e}(s))^{\dagger} \, \md s,
\end{align*}
one arrives at $2 \int_{J_{e}}\Delta^{e}(s)(\Delta^{e}(s))^{\dagger} \, \md
s \! = \! 0$; thus, $\Delta^{e}(z) \! = \! 0$, $z \! \in \! J_{e}$, which
implies that
$\overset{e}{\mathfrak{M}}^{\raise-1.0ex\hbox{$\scriptstyle \infty$}}(z) \! =
\! \mathscr{S}^{\infty}_{e}(z)$ for all $z$. \hfill $\qed$

In order to prove that there is a solution of the (full) RHP 
$(\overset{e}{\mathscr{M}}^{\raise-1.0ex\hbox{$\scriptstyle \sharp$}}(z),
\overset{e}{\upsilon}^{\raise-1.0ex\hbox{$\scriptstyle \sharp$}}(z),\Sigma_{
e}^{\sharp})$, formulated in Lemma~4.2, close to the parametrix, one needs 
to know that the parametrix is \emph{uniformly} bounded: more precisely, by 
(certain) general theorems (see, for example, \cite{a98}), one needs to know 
that $\overset{e}{\upsilon}^{\raise-1.0ex\hbox{$\scriptstyle \sharp$}}(z) \! 
\to \! \overset{e}{\upsilon}^{\raise-1.0ex\hbox{$\scriptstyle \infty$}}(z)$ 
as $n \! \to \! \infty$ uniformly for $z \! \in \! \Sigma_{e}^{\sharp}$ in 
the $\mathcal{L}^{2}_{\mathrm{M}_{2}(\mathbb{C})}(\Sigma_{e}^{\sharp}) \cap 
\mathcal{L}^{\infty}_{\mathrm{M}_{2}(\mathbb{C})}(\Sigma_{e}^{\sharp})$ sense, 
that is, uniformly,
\begin{equation*}
\lim_{n \to \infty} \norm{\overset{e}{\upsilon}^{
\raise-1.0ex\hbox{$\scriptstyle \sharp$}}(\cdot) \! - \!
\overset{e}{\upsilon}^{\raise-1.0ex\hbox{$\scriptstyle \infty$}}(\cdot)}_{
\mathcal{L}^{2}_{\mathrm{M}_{2}(\mathbb{C})}(\Sigma_{e}^{\sharp}) \cap
\mathcal{L}^{\infty}_{\mathrm{M}_{2}(\mathbb{C})}(\Sigma_{e}^{\sharp})} \! :=
\! \lim_{n \to \infty} \sum_{p \in \{2,\infty\}} \norm{\overset{e}{\upsilon}^{
\raise-1.0ex\hbox{$\scriptstyle \sharp$}}(\cdot) \! - \! \overset{e}{
\upsilon}^{\raise-1.0ex\hbox{$\scriptstyle \infty$}}(\cdot)}_{\mathcal{L}^{
p}_{\mathrm{M}_{2}(\mathbb{C})}(\Sigma_{e}^{\sharp})} \! = \! 0;
\end{equation*}
however, notwithstanding the fact that $\widetilde{V} \colon \mathbb{R}
\setminus \{0\} \! \to \! \mathbb{R}$ is regular $(h_{V}^{e}(b_{j-1}^{e}),
h_{V}^{e}(a_{j}^{e}) \! \not= \! 0$, $j \! = \! 1,\dotsc,N \! + \! 1)$, since
the strict inequalities $g^{e}_{+}(z) \! + \! g^{e}_{-}(z) \! - \! \widetilde{
V}(z) \! - \! \ell_{e} \! + \! 2Q_{e} \! < \! 0$, $z \! \in \! (-\infty,b_{
0}^{e}) \cup (a_{N+1}^{e},+\infty) \cup (\cup_{j=1}^{N}(a_{j}^{e},b_{j}^{e}
))$, and $\pm \Re (\mi \int_{z}^{a_{N+1}^{e}} \psi_{V}^{e}(s) \, \md s) \!
> \! 0$, $z \! \in \! \mathbb{C}_{\pm} \cap (\cup_{j=1}^{N+1} \mathbb{U}_{
j}^{e})$, fail at the end-points of the support of the `even' equilibrium
measure, this implies that
$\overset{e}{\upsilon}^{\raise-1.0ex\hbox{$\scriptstyle \sharp$}}(z) \! \to \!
\overset{e}{\upsilon}^{\raise-1.0ex\hbox{$\scriptstyle \infty$}}(z)$ as $n \!
\to \! \infty$ \emph{pointwise}, but not uniformly, for $z \! \in \! \Sigma_{
e}^{\sharp}$, whence, one can not conclude that
$\overset{e}{\mathscr{M}}^{\raise-1.0ex\hbox{$\scriptstyle \sharp$}}(z) \! \to
\! \overset{e}{m}^{\raise-1.0ex\hbox{$\scriptstyle \infty$}}(z)$ as $n \! \to
\! \infty$ uniformly for $z \! \in \! \Sigma_{e}^{\sharp}$. The resolution of
this lack of uniformity at the end-points of the support of the `even'
equilibrium measure constitutes, therefore, the essential analytical obstacle
remaining for the analysis of the RHP
$(\overset{e}{\mathscr{M}}^{\raise-1.0ex\hbox{$\scriptstyle \sharp$}}(z),
\overset{e}{\upsilon}^{\raise-1.0ex\hbox{$\scriptstyle \sharp$}}(z),\Sigma_{
e}^{\sharp})$, and a substantial part of the following analysis is devoted to
overcoming this problem.

The key necessary to remedy (and control) the above-mentioned analytical 
difficulty is to construct parametrices for the solution of the RHP 
$(\overset{e}{\mathscr{M}}^{\raise-1.0ex\hbox{$\scriptstyle \sharp$}}(z),
\overset{e}{\upsilon}^{\raise-1.0ex\hbox{$\scriptstyle \sharp$}}(z),\Sigma_{
e}^{\sharp})$ in `small' neighbourhoods (open discs) about $\{b_{j-1}^{e},
a_{j}^{e}\}_{j=1}^{N+1}$ (where the convergence of 
$\overset{e}{\upsilon}^{\raise-1.0ex\hbox{$\scriptstyle \sharp$}}(z)$ to 
$\overset{e}{\upsilon}^{\raise-1.0ex\hbox{$\scriptstyle \infty$}}(z)$ as 
$n \! \to \! \infty$ is not uniform) in such a way that, on the boundary 
of these open neighbourhoods, the parametrices `match' with the solution
of the model RHP, 
$\overset{e}{m}^{\raise-1.0ex\hbox{$\scriptstyle \infty$}}(z)$, up to $o(1)$
(in fact, $\mathcal{O}(n^{-1}))$ as $n \! \to \! \infty$; furthermore, in 
the generic framework considered in this work, namely, $\widetilde{V} \colon 
\mathbb{R} \setminus \{0\} \! \to \! \mathbb{R}$ is regular, in which case 
the (density of the) `even' equilibrium measure behaves as a square root at 
the end-points of $\operatorname{supp}(\mu_{V}^{e})$, that is, $\psi_{V}^{e}
(s) \! =_{s \downarrow b_{j-1}^{e}} \! \mathcal{O}((s \! - \! b_{j-1}^{e})^{
1/2})$ and $\psi_{V}^{e}(s) \! =_{s \uparrow a_{j}^{e}} \! \mathcal{O}
((a_{j}^{e} \! - \! s)^{1/2})$, $j \! = \! 1,\dotsc,N \! + \! 1$, it is well 
known \cite{a3,a59,a90,a99} that the parametrices can be expressed in terms 
of Airy functions. (The general method used to construct such parametrices 
is via a `Vanishing Lemma' \cite{a100}.) More precisely, one surrounds the 
end-points of the support of the `even' equilibrium measure, $\lbrace b_{j-
1}^{e},a_{j}^{e} \rbrace_{j=1}^{N+1}$, by `small', mutually disjoint open 
discs,
\begin{equation*}
\mathbb{D}_{\epsilon}(b_{j-1}^{e}) \! := \! \left\lbrace \mathstrut z \! \in
\! \mathbb{C}; \, \vert z \! - \! b_{j-1}^{e} \vert \! < \! \epsilon_{j}^{b}
\right\rbrace \qquad \text{and} \qquad \mathbb{D}_{\epsilon}(a_{j}^{e}) \! :=
\! \left\lbrace \mathstrut z \! \in \! \mathbb{C}; \, \vert z \! - \! a_{j}^{
e} \vert \! < \! \epsilon_{j}^{a} \right\rbrace, \quad j \! = \! 1,\dotsc,N \!
+ \! 1,
\end{equation*}
where $\epsilon_{j}^{b},\epsilon_{j}^{a}$ are arbitrarily fixed, sufficiently
small positive real numbers chosen so that $\mathbb{D}_{\epsilon}(b_{i-1}^{e})
\cap \mathbb{D}_{\epsilon}(a_{j}^{e}) \! = \! \varnothing$, $i,j \! = \! 1,
\dotsc,N \! + \! 1$, and defines $S_{p}^{e}(z)$, the parametrix for
$\overset{e}{\mathscr{M}}^{\raise-1.0ex\hbox{$\scriptstyle \sharp$}}(z)$,
by $\overset{e}{m}^{\raise-1.0ex\hbox{$\scriptstyle \infty$}}(z)$ for $z \!
\in \! \mathbb{C} \setminus (\cup_{j=1}^{N+1}(\mathbb{D}_{\epsilon}(b_{j-1}^{
e}) \cup \mathbb{D}_{\epsilon}(a_{j}^{e})))$, and by $m_{p}^{e}(z)$ for $z \!
\in \! \cup_{j=1}^{N+1}(\mathbb{D}_{\epsilon}(b_{j-1}^{e}) \cup \mathbb{D}_{
\epsilon}(a_{j}^{e}))$, and solves the local RHPs for $m_{p}^{e}(z)$ on
$\cup_{j=1}^{N+1}(\mathbb{D}_{\epsilon}(b_{j-1}^{e}) \cup \mathbb{D}_{\epsilon}
(a_{j}^{e}))$ in such a way (`optimal' in the nomenclature of \cite{a59}) 
that $m_{p}^{e}(z) \! \approx_{n \to \infty} \! 
\overset{e}{m}^{\raise-1.0ex\hbox{$\scriptstyle \infty$}}(z)$ (to $\mathcal{O}
(n^{-1}))$ for $z \! \in \! \cup_{j=1}^{N+1}(\partial \mathbb{D}_{\epsilon}
(b_{j-1}^{e}) \cup \partial \mathbb{D}_{\epsilon}(a_{j}^{e}))$, whence 
$\mathcal{R}^{e}(z) \! := \! 
\overset{e}{\mathscr{M}}^{\raise-1.0ex\hbox{$\scriptstyle \sharp$}}(z)(S_{p}^{
e}(z))^{-1} \colon \mathbb{C} \setminus \widetilde{\Sigma}_{e}^{\sharp} \! \to 
\! \operatorname{SL}_{2}(\mathbb{C})$, where $\widetilde{\Sigma}_{e}^{\sharp}
\! := \! \Sigma_{e}^{\sharp} \cup (\cup_{j=1}^{N+1}(\partial \mathbb{D}_{
\epsilon}(b_{j-1}^{e}) \cup \partial \mathbb{D}_{\epsilon}(a_{j}^{e})))$, 
solves the RHP $(\mathcal{R}^{e}(z),\upsilon^{e}_{\mathcal{R}}(z),\widetilde{
\Sigma}_{e}^{\sharp})$, with $\norm{\upsilon^{e}_{\mathcal{R}}(\cdot) \! - \!
\mathrm{I}}_{\cap_{p \in \{2,\infty\}} \mathcal{L}^{p}_{\mathrm{M}_{2}
(\mathbb{C})}(\widetilde{\Sigma}_{e}^{\sharp})} \! =_{n \to \infty} \!
\mathcal{O}(n^{-1})$ uniformly; in particular, the error term, which is
$\mathcal{O}(n^{-1})$ as $n \! \to \! \infty$, is uniform in $\cap_{p \in
\{1,2,\infty\}} \mathcal{L}^{p}_{\mathrm{M}_{2}(\mathbb{C})}(\widetilde{
\Sigma}_{e}^{\sharp})$. By general Riemann-Hilbert techniques (see, for
example, \cite{a98}), $\mathcal{R}^{e}(z)$ (and thus
$\overset{e}{\mathscr{M}}^{\raise-1.0ex\hbox{$\scriptstyle \sharp$}}(z)$ 
via the relation 
$\overset{e}{\mathscr{M}}^{\raise-1.0ex\hbox{$\scriptstyle \sharp$}}(z) \! = 
\! \mathcal{R}^{e}(z)S_{p}^{e}(z))$ can be computed to any order of $n^{-1}$ 
via a Neumann series expansion (of the corresponding resolvent kernel). In 
fact, at the very core of the above-mentioned discussion, and the analysis 
that follows, is the following Corollary (see, for example, \cite{a90}, 
Corollary~7.108):
\begin{fffff}[Deift {\rm \cite{a90}}]
For an oriented contour $\varSigma \subset \mathbb{C}$, let $m^{\infty} \colon
\mathbb{C} \setminus \varSigma \! \to \! \operatorname{SL}_{2}(\mathbb{C})$
and $m^{(n)} \colon \mathbb{C} \setminus \varSigma \! \to \!
\operatorname{SL}_{2}(\mathbb{C})$, $n \! \in \! \mathbb{N}$, respectively,
solve the following, equivalent {\rm RHPs}, $(m^{\infty}(z),\upsilon^{\infty}
(z),\varSigma)$ and $(m^{(n)}(z),\upsilon^{(n)}(z),\linebreak[4]
\varSigma)$, where
\begin{equation*}
\upsilon^{\infty} \colon \varSigma \! \to \! \operatorname{GL}_{2}(\mathbb{C}
), \, \, z \! \mapsto \! \left(\mathrm{I} \! - \! w^{\infty}_{-}(z) \right)^{
-1} \! \left(\mathrm{I} \! + \! w^{\infty}_{+}(z) \right)
\end{equation*}
and
\begin{equation*}
\upsilon^{(n)} \colon \varSigma \! \to \! \operatorname{GL}_{2}(\mathbb{C}),
\, \, z \! \mapsto \! \left(\mathrm{I} \! - \! w^{(n)}_{-}(z) \right)^{-1} \!
\left(\mathrm{I} \! + \! w^{(n)}_{+}(z) \right),
\end{equation*}
and suppose that $(\id \! - \! C^{\infty}_{w^{\infty}})^{-1}$ exists, where
\begin{equation*}
\mathcal{L}^{2}_{\mathrm{M}_{2}(\mathbb{C})}(\varSigma) \! \ni \! f \! \mapsto
\! C^{\infty}_{w^{\infty}}f \! := \! C^{\infty}_{+}(fw^{\infty}_{-}) \! + \!
C^{\infty}_{-}(fw^{\infty}_{+}),
\end{equation*}
with
\begin{equation*}
C^{\infty}_{\pm} \colon \mathcal{L}^{2}_{\mathrm{M}_{2}(\mathbb{C})}
(\varSigma) \! \to \! \mathcal{L}^{2}_{\mathrm{M}_{2}(\mathbb{C})}(\varSigma),
\, \, f \! \mapsto \! (C^{\infty}_{\pm}f)(z) := \lim_{\underset{z^{\prime} \,
\in \, \pm \, \mathrm{side} \, \mathrm{of} \, \varSigma}{z^{\prime} \to z}}
\int_{\varSigma} \dfrac{f(s)}{s \! - \! z^{\prime}} \, \dfrac{\md s}{2 \pi
\mi},
\end{equation*}
and $\norm{w^{(n)}_{l}(\cdot) \! - \! w^{\infty}_{l}(\cdot)}_{\cap_{p \in \{2,
\infty\}}\mathcal{L}^{p}_{\mathrm{M}_{2}(\mathbb{C})}(\varSigma)} \! \to \! 0$
as $n \! \to \! \infty$, $l \! = \! \pm 1$. Then, $\exists \, \, N^{\ast} \!
\in \! \mathbb{N}$ such that, $\forall \, \, n \! > \! N^{\ast}$, $m^{\infty}
(z)$ and $m^{(n)}(z)$ exist, and $\norm{m^{(n)}_{l}(\cdot) \! - \! m^{
\infty}_{l}(\cdot)}_{\mathcal{L}^{2}_{\mathrm{M}_{2}(\mathbb{C})}(\varSigma)}
\! \to \! 0$ as $n \! \to \! \infty$, $l \! = \! \pm 1$.
\end{fffff}

A detailed exposition, including further motivations, for the construction of 
parametrices of the above-mentioned type can be found in \cite{a3,a57,a58,%
a59,a61,a90}; rather than regurgitating, verbatim, much of the analysis that 
can be found in the latter references, the point of view taken here is that 
one follows the scheme presented therein to obtain the results stated below, 
that is, the parametrix for the RHP 
$(\overset{e}{\mathscr{M}}^{\raise-1.0ex\hbox{$\scriptstyle \sharp$}}(z),
\overset{e}{\upsilon}^{\raise-1.0ex\hbox{$\scriptstyle \sharp$}}(z),\Sigma_{
e}^{\sharp})$ formulated in Lemma~4.2. In the case of the right-most 
end-points of the support of the `even' equilibrium measure, $\lbrace 
a_{j}^{e} \rbrace_{j=1}^{N+1}$, a terse sketch of a proof is presented for 
the reader's convenience, and the remaining (left-most) end-points, namely, 
$b_{0}^{e},b_{1}^{e},\dotsc,b_{N}^{e}$, are analysed analogously.

The parametrix for the RHP
$(\overset{e}{\mathscr{M}}^{\raise-1.0ex\hbox{$\scriptstyle \sharp$}}
(z),\overset{e}{\upsilon}^{\raise-1.0ex\hbox{$\scriptstyle \sharp$}}(z),
\Sigma_{e}^{\sharp})$ is now presented. By a parametrix of the RHP
$(\overset{e}{\mathscr{M}}^{\raise-1.0ex\hbox{$\scriptstyle \sharp$}}(z),
\overset{e}{\upsilon}^{\raise-1.0ex\hbox{$\scriptstyle \sharp$}}(z),\Sigma_{
e}^{\sharp})$, in the neighbourhoods of the end-points of the support of the
`even' equilibrium measure, $\lbrace b_{j-1}^{e},a_{j}^{e} \rbrace_{j=1}^{N+
1}$, is meant the solution of the RHPs formulated in the following two Lemmas
(Lemmas~4.6 and~4.7). Define the `small', mutually disjoint (open) discs about
the end-points of the support of the `even' equilibrium measure as follows:
$\mathbb{U}_{\delta_{b_{j-1}}}^{e} \! := \! \lbrace \mathstrut z \! \in \!
\mathbb{C}; \, \vert z \! - \! b_{j-1}^{e} \vert \! < \! \delta_{b_{j-1}}^{
e} \! \in \! (0,1) \rbrace$ and $\mathbb{U}_{\delta_{a_{j}}}^{e} \! := \!
\lbrace \mathstrut z \! \in \! \mathbb{C}; \, \vert z \! - \! a_{j}^{e} \vert
\! < \! \delta_{a_{j}}^{e} \! \in \! (0,1) \rbrace$, $j \! = \! 1,\dotsc,N \!
+ \! 1$, where $\delta_{b_{j-1}}^{e}$ and $\delta_{a_{j}}^{e}$ are
sufficiently small, positive real numbers chosen (amongst other things: see
Lemmas~4.6 and~4.7 below) so that $\mathbb{U}_{\delta_{b_{i-1}}}^{e} \cap
\mathbb{U}_{\delta_{a_{j}}}^{e} \! = \! \varnothing$, $i,j \! = \! 1,\dotsc,
N \! + \! 1$ (the corresponding regions $\Omega_{b_{j-1}}^{e,l}$ and $\Omega_{
a_{j}}^{e,l}$, and arcs $\Sigma_{b_{j-1}}^{e,l}$ and $\Sigma_{a_{j}}^{e,l}$,
$j \! = \! 1,\dotsc,N \! + \! 1$, $l \! = \! 1,2,3,4$, respectively, are
defined more precisely below; see, also, Figures~5 and~6).
\begin{eeeee}
In order to simplify the results of Lemmas~4.6 and~4.7 (see below), it is
convenient to introduce the following notation: (i)
\begin{gather*}
\Psi^{e}_{1}(z) \! := \!
\begin{pmatrix}
\operatorname{Ai}(z) & \operatorname{Ai}(\omega^{2}z) \\
\operatorname{Ai}^{\prime}(z) & \omega^{2} \operatorname{Ai}^{\prime}
(\omega^{2}z)
\end{pmatrix} \! \me^{-\frac{\mi \pi}{6} \sigma_{3}}, \qquad \quad \Psi^{e}_{2}
(z) \! := \!
\begin{pmatrix}
\operatorname{Ai}(z) & \operatorname{Ai}(\omega^{2}z) \\
\operatorname{Ai}^{\prime}(z) & \omega^{2} \operatorname{Ai}^{\prime}
(\omega^{2}z)
\end{pmatrix} \! \me^{-\frac{\mi \pi}{6} \sigma_{3}}(\mathrm{I} \! - \!
\sigma_{-}), \\
\Psi^{e}_{3}(z) \! := \!
\begin{pmatrix}
\operatorname{Ai}(z) & -\omega^{2} \operatorname{Ai}(\omega z) \\
\operatorname{Ai}^{\prime}(z) & -\operatorname{Ai}^{\prime}(\omega z)
\end{pmatrix} \! \me^{-\frac{\mi \pi}{6} \sigma_{3}}(\mathrm{I} \! + \!
\sigma_{-}), \qquad \quad \Psi^{e}_{4}(z) \! := \!
\begin{pmatrix}
\operatorname{Ai}(z) & -\omega^{2} \operatorname{Ai}(\omega z) \\
\operatorname{Ai}^{\prime}(z) & -\operatorname{Ai}^{\prime}(\omega z)
\end{pmatrix} \! \me^{-\frac{\mi \pi}{6} \sigma_{3}},
\end{gather*}
where $\operatorname{Ai}(\pmb{\cdot})$ is the Airy function (cf.
Subsection~2.3), and $\omega \! = \! \exp (2 \pi \mi/3)$; and (ii)
\begin{equation*}
\mho^{e}_{j} \! := \!
\begin{cases}
\Omega^{e}_{j}, &\text{$j \! = \! 1,\dotsc,N$,} \\
0, &\text{$j \! = \! 0,N \! + \! 1$,}
\end{cases}
\end{equation*}
where $\Omega^{e}_{j} \! = \! 4 \pi \int_{b_{j}^{e}}^{a_{N+1}^{e}} \psi_{V}^{
e}(s) \, \md s$. \hfill $\blacksquare$
\end{eeeee}
\begin{ccccc}
Let $\overset{e}{\mathscr{M}}^{\raise-1.0ex\hbox{$\scriptstyle \sharp$}}
\colon \mathbb{C} \setminus \Sigma_{e}^{\sharp} \! \to \! \operatorname{SL}_{
2}(\mathbb{C})$ solve the {\rm RHP}
$(\overset{e}{\mathscr{M}}^{\raise-1.0ex\hbox{$\scriptstyle \sharp$}}(z),
\overset{e}{\upsilon}^{\raise-1.0ex\hbox{$\scriptstyle \sharp$}}(z),\Sigma_{
e}^{\sharp})$ formulated in Lemma~{\rm 4.2}, and set
\begin{equation*}
\mathbb{U}_{\delta_{b_{j-1}}}^{e} \! := \! \left\{ \mathstrut z \! \in \!
\mathbb{C}; \, \vert z \! - \! b_{j-1}^{e} \vert \! < \! \delta_{b_{j-1}}^{e}
\! \in \! (0,1) \right\}, \quad j \! = \! 1,\dotsc,N \! + \! 1.
\end{equation*}
Let
\begin{equation*}
\Phi_{b_{j-1}}^{e}(z) \! := \! \left(\dfrac{3n}{4} \xi_{b_{j-1}}^{e}(z)
\right)^{2/3}, \quad j \! = \! 1,\dotsc,N \! + \! 1,
\end{equation*}
with
\begin{equation*}
\xi_{b_{j-1}}^{e}(z) \! = \! -2 \int_{z}^{b_{j-1}^{e}}(R_{e}(s))^{1/2}
h_{V}^{e}(s) \, \md s,
\end{equation*}
where, for $z \! \in \! \mathbb{U}_{\delta_{b_{j-1}}}^{e} \setminus (-\infty,
b_{j-1}^{e})$, $\xi_{b_{j-1}}^{e}(z) \! = \! \mathfrak{b}(z \! - \! b_{j-
1}^{e})^{3/2}G_{b_{j-1}}^{e}(z)$, $j \! = \! 1,\dotsc,N \! + \! 1$, with
$\mathfrak{b} \! := \! \pm 1$ for $z \! \in \! \mathbb{C}_{\pm}$, and
$G_{b_{j-1}}^{e}(z)$ analytic, in particular,
\begin{equation*}
G_{b_{j-1}}^{e}(z) \underset{z \to b_{j-1}^{e}}{=} \dfrac{4}{3}f(b_{j-1}^{e})
\! + \! \dfrac{4}{5}f^{\prime}(b_{j-1}^{e})(z \! - \! b_{j-1}^{e}) \! + \!
\dfrac{2}{7}f^{\prime \prime}(b_{j-1}^{e})(z \! - \! b_{j-1}^{e})^{2} \! + \!
\mathcal{O} \! \left((z \! - \! b_{j-1}^{e})^{3} \right),
\end{equation*}
where
\begin{align*}
f(b_{0}^{e})=& \, \mi (-1)^{N}h_{V}^{e}(b_{0}^{e}) \eta_{b_{0}^{e}}, \\
f^{\prime}(b_{0}^{e})=& \, \mi (-1)^{N} \! \left(\dfrac{1}{2}h_{V}^{e}(b_{0}^{
e}) \eta_{b_{0}^{e}} \! \left(\sum_{l=1}^{N} \! \left(\dfrac{1}{b_{0}^{e} \! -
\! b_{l}^{e}} \! + \! \dfrac{1}{b_{0}^{e} \! - \! a_{l}^{e}} \right) \! + \!
\dfrac{1}{b_{0}^{e} \! - \! a_{N+1}^{e}} \right) \! + \! (h_{V}^{e}(b_{0}^{e}
))^{\prime} \eta_{b_{0}^{e}} \right), \\
f^{\prime \prime}(b_{0}^{e})=& \, \mi (-1)^{N} \! \left(\dfrac{h_{V}^{e}(b_{
0}^{e})(h_{V}^{e}(b_{0}^{e}))^{\prime \prime} \! - \! ((h_{V}^{e}(b_{0}^{e}))^{
\prime})^{2}}{h_{V}^{e}(b_{0}^{e})} \eta_{b_{0}^{e}} \! - \! \dfrac{1}{2}h_{
V}^{e}(b_{0}^{e}) \eta_{b_{0}^{e}} \right. \\
\times&\left. \, \left(\sum_{l=1}^{N} \! \left(\dfrac{1}{(b_{0}^{e} \! - \!
b_{l}^{e})^{2}} \! + \! \dfrac{1}{(b_{0}^{e} \! - \! a_{l}^{e})^{2}} \right)
\! + \! \dfrac{1}{(b_{0}^{e} \! - \! a_{N+1}^{e})^{2}} \right) \right. \\
+&\left. \, \left(\dfrac{1}{2} \! \left(\sum_{k=1}^{N} \! \left(\dfrac{1}{b_{
0}^{e} \! - \! b_{k}^{e}} \! + \! \dfrac{1}{b_{0}^{e} \! - \! a_{k}^{e}}
\right) \! + \! \dfrac{1}{b_{0}^{e} \! - \! a_{N+1}^{e}} \right) \! + \!
\dfrac{(h_{V}^{e}(b_{0}^{e}))^{\prime}}{h_{V}^{e}(b_{0}^{e})} \right) \right.
\\
\times&\left. \, \left(\dfrac{1}{2}h_{V}^{e}(b_{0}^{e}) \eta_{b_{0}^{e}} \!
\left(\sum_{l=1}^{N} \! \left(\dfrac{1}{b_{0}^{e} \! - \! b_{l}^{e}} \! + \!
\dfrac{1}{b_{0}^{e} \! - \! a_{l}^{e}} \right) \! + \! \dfrac{1}{b_{0}^{e} \!
- \! a_{N+1}^{e}} \right) \! + \! (h_{V}^{e}(b_{0}^{e}))^{\prime} \eta_{b_{0}^{
e}} \right) \right),
\end{align*}
with
\begin{equation*}
\eta_{b_{0}^{e}} \! := \! \left((a_{N+1}^{e} \! - \! b_{0}^{e}) \prod_{k=1}^{N}
(b_{k}^{e} \! - \! b_{0}^{e})(a_{k}^{e} \! - \! b_{0}^{e}) \right)^{1/2} \quad
(> \! 0),
\end{equation*}
and, for $j \! = \! 1,\dotsc,N$,
\begin{align*}
f(b_{j}^{e})=& \, \mi (-1)^{N-j}h_{V}^{e}(b_{j}^{e}) \eta_{b_{j}^{e}}, \\
f^{\prime}(b_{j}^{e})=& \, \mi (-1)^{N-j} \! \left(\dfrac{1}{2}h_{V}^{e}(b_{
j}^{e}) \eta_{b_{j}^{e}} \! \left(\sum_{\substack{k=1\\k \not= j}}^{N} \!
\left(\dfrac{1}{b_{j}^{e} \! - \! b_{k}^{e}} \! + \! \dfrac{1}{b_{j}^{e} \! -
\! a_{k}^{e}} \right) \! + \! \dfrac{1}{b_{j}^{e} \! - \! a_{j}^{e}} \! + \!
\dfrac{1}{b_{j}^{e} \! - \! a_{N+1}^{e}} \! + \! \dfrac{1}{b_{j}^{e} \! - \!
b_{0}^{e}} \right) \right. \\
+&\left. \, (h_{V}^{e}(b_{j}^{e}))^{\prime} \eta_{b_{j}^{e}} \right), \\
f^{\prime \prime}(b_{j}^{e})=& \, \mi (-1)^{N-j} \! \left(\dfrac{h_{V}^{e}(b_{
j}^{e})(h_{V}^{e}(b_{j}^{e}))^{\prime \prime} \! - \! ((h_{V}^{e}(b_{j}^{e}))^{
\prime})^{2}}{h_{V}^{e}(b_{j}^{e})} \eta_{b_{j}^{e}} \! - \! \dfrac{1}{2}h_{
V}^{e}(b_{j}^{e}) \eta_{b_{j}^{e}} \! \left(\sum_{\substack{k=1\\k \not= j}}^{
N} \! \left(\dfrac{1}{(b_{j}^{e} \! - \! b_{k}^{e})^{2}} \! + \! \dfrac{1}{(b_{
j}^{e} \! - \! a_{k}^{e})^{2}} \right) \right. \right. \\
+&\left. \left. \, \dfrac{1}{(b_{j}^{e} \! - \! a_{j}^{e})^{2}} \! + \! \dfrac{
1}{(b_{j}^{e} \! - \! a_{N+1}^{e})^{2}} \! + \! \dfrac{1}{(b_{j}^{e} \! - \!
b_{0}^{e})^{2}} \right) \! + \! \left(\dfrac{(h_{V}^{e}(b_{j}^{e}))^{\prime}}{
h_{V}^{e}(b_{j}^{e})} \! + \! \dfrac{1}{2} \! \left(\sum_{\substack{k=1\\k
\not= j}}^{N} \! \left(\dfrac{1}{b_{j}^{e} \! - \! b_{k}^{e}} \! + \! \dfrac{
1}{b_{j}^{e} \! - \! a_{k}^{e}} \right) \right. \right. \right. \\
+&\left. \left. \left. \, \dfrac{1}{b_{j}^{e} \! - \! a_{j}^{e}} \! + \!
\dfrac{1}{b_{j}^{e} \! - \! a_{N+1}^{e}} \! + \! \dfrac{1}{b_{j}^{e} \! - \!
b_{0}^{e}} \right) \! \right) \! \left(\dfrac{1}{2}h_{V}^{e}(b_{j}^{e}) \eta_{
b_{j}^{e}} \! \left(\sum_{\substack{k=1\\k \not= j}}^{N} \! \left(\dfrac{1}{b_{
j}^{e} \! - \! b_{k}^{e}} \! + \! \dfrac{1}{b_{j}^{e} \! - \! a_{k}^{e}}
\right) \right. \right. \right. \\
+&\left. \left. \left. \, \dfrac{1}{b_{j}^{e} \! - \! a_{j}^{e}} \! + \!
\dfrac{1}{b_{j}^{e} \! - \! a_{N+1}^{e}} \! + \! \dfrac{1}{b_{j}^{e} \! - \!
b_{0}^{e}} \right) \! + \!  (h_{V}^{e}(b_{j}^{e}))^{\prime} \eta_{b_{j}^{e}}
\right) \right),
\end{align*}
with
\begin{equation*}
\eta_{b_{j}^{e}} \! := \! \left((b_{j}^{e} \! - \! a_{j}^{e})(a_{N+1}^{e} \!
- \! b_{j}^{e})(b_{j}^{e} \! - \! b_{0}^{e}) \prod_{k=1}^{j-1}(b_{j}^{e} \! -
\! b_{k}^{e})(b_{j}^{e} \! - \! a_{k}^{e}) \prod_{l=j+1}^{N}(b_{l}^{e} \! - \!
b_{j}^{e})(a_{l}^{e} \! - \! b_{j}^{e}) \right)^{1/2} \quad (> \! 0),
\end{equation*}
and $((0,1) \! \ni)$ $\delta_{b_{j-1}}^{e}$, $j \! = \! 1,\dotsc,N \! + \! 1$,
are chosen sufficiently small so that $\Phi_{b_{j-1}}^{e}(z)$, which are
bi-holomorphic, conformal, and non-orientation preserving, map $\mathbb{U}_{
\delta_{b_{j-1}}}^{e}$ (and, thus, the oriented contours $\Sigma_{b_{j-1}}^{e}
\! := \! \cup_{l=1}^{4} \Sigma_{b_{j-1}}^{e,l}$, $j \! = \! 1,\dotsc,N \! + \!
1:$ Figure~{\rm 6)} injectively onto open $(n$-dependent) neighbourhoods
$\widehat{\mathbb{U}}_{\delta_{b_{j-1}}}^{e}$, $j \! = \! 1,\dotsc,N \! + \!
1$, of $0$ such that $\Phi_{b_{j-1}}^{e}(b_{j-1}^{e}) \! = \! 0$, $\Phi_{b_{j-
1}}^{e} \colon \mathbb{U}_{\delta_{b_{j-1}}}^{e} \! \to \! \widehat{\mathbb{
U}}_{\delta_{b_{j-1}}}^{e} \! := \! \Phi_{b_{j-1}}^{e}(\mathbb{U}_{\delta_{b_{
j-1}}}^{e})$, $\Phi_{b_{j-1}}^{e}(\mathbb{U}_{\delta_{b_{j-1}}}^{e} \cap
\Sigma_{b_{j-1}}^{e,l}) \! = \! \Phi_{b_{j-1}}^{e}(\mathbb{U}_{\delta_{b_{j-
1}}}^{e}) \cap \gamma_{b_{j-1}}^{e,l}$, and $\Phi_{b_{j-1}}^{e}(\mathbb{U}_{
\delta_{b_{j-1}}}^{e} \cap \Omega_{b_{j-1}}^{e,l}) \! = \! \Phi_{b_{j-1}}^{e}
(\mathbb{U}_{\delta_{b_{j-1}}}^{e}) \cap \widehat{\Omega}_{b_{j-1}}^{e,l}$,
$l \! = \! 1,2,3,4$, with $\widehat{\Omega}_{b_{j-1}}^{e,1} \! = \! \{
\mathstrut \zeta \! \in \! \mathbb{C}; \, \arg (\zeta) \! \in \! (0,2 \pi/3)
\}$, $\widehat{\Omega}_{b_{j-1}}^{e,2} \! = \! \{\mathstrut \zeta \! \in \!
\mathbb{C}; \, \arg (\zeta) \! \in \! (2 \pi/3,\pi)\}$, $\widehat{\Omega}_{
b_{j-1}}^{e,3} \! = \! \{\mathstrut \zeta \! \in \! \mathbb{C}; \, \arg
(\zeta) \! \in \! (-\pi,-2 \pi/3)\}$, and $\widehat{\Omega}_{b_{j-1}}^{e,4} \!
= \! \{\mathstrut \zeta \! \in \! \mathbb{C}; \, \arg (\zeta) \! \in \! (-2
\pi/3,0)\}$.

The parametrix for the {\rm RHP}
$(\overset{e}{\mathscr{M}}^{\raise-1.0ex\hbox{$\scriptstyle \sharp$}}(z),
\overset{e}{\upsilon}^{\raise-1.0ex\hbox{$\scriptstyle \sharp$}}(z),\Sigma_{
e}^{\sharp})$, for $z \! \in \! \mathbb{U}_{\delta_{b_{j-1}}}^{e}$, $j \! =
\! 1,\dotsc,N \! + \! 1$, is the solution of the following {\rm RHPs} for
$\mathcal{X}^{e} \colon \mathbb{U}_{\delta_{b_{j-1}}}^{e} \setminus \Sigma_{
b_{j-1}}^{e} \! \to \! \operatorname{SL}_{2}(\mathbb{C})$, $j \! = \! 1,\dotsc,
N \! + \! 1$, where $\Sigma_{b_{j-1}}^{e} \! := \! (\Phi_{b_{j-1}}^{e})^{-1}
(\gamma_{b_{j-1}}^{e})$, with $(\Phi_{b_{j-1}}^{e})^{-1}$ denoting the inverse
mapping, and $\gamma_{b_{j-1}}^{e} \! := \! \cup_{l=1}^{4} \gamma_{b_{j-1}}^{e,
l}:$ {\rm (i)} $\mathcal{X}^{e}(z)$ is holomorphic for $z \! \in \! \mathbb{
U}_{\delta_{b_{j-1}}}^{e} \setminus \Sigma_{b_{j-1}}^{e}$, $j \! = \! 1,
\dotsc,N \! + \! 1;$ {\rm (ii)} $\mathcal{X}^{e}_{\pm}(z) \! := \! \lim_{
\underset{z^{\prime} \, \in \, \pm \, \mathrm{side} \, \mathrm{of} \, \Sigma_{
b_{j-1}}^{e}}{z^{\prime} \to z}} \mathcal{X}^{e}(z^{\prime})$, $j \! = \! 1,
\dotsc,N \! + \! 1$, satisfy the boundary condition
\begin{equation*}
\mathcal{X}^{e}_{+}(z) \! = \! \mathcal{X}^{e}_{-}(z)
\overset{e}{\upsilon}^{\raise-1.0ex\hbox{$\scriptstyle \sharp$}}(z), \quad z
\! \in \! \mathbb{U}_{\delta_{b_{j-1}}}^{e} \cap \Sigma_{b_{j-1}}^{e}, \quad j
\! = \! 1,\dotsc,N \! + \! 1,
\end{equation*}
where $\overset{e}{\upsilon}^{\raise-1.0ex\hbox{$\scriptstyle \sharp$}}(z)$ is
given in Lemma~{\rm 4.2;} and {\rm (iii)} uniformly for $z \! \in \! \partial
\mathbb{U}_{\delta_{b_{j-1}}}^{e} \! := \! \left\lbrace \mathstrut z \! \in \!
\mathbb{C}; \, \vert z \! - \! b_{j-1}^{e} \vert \! = \! \delta_{b_{j-1}}^{e}
\right\rbrace$, $j \! = \! 1,\dotsc,N \! + \! 1$,
\begin{equation*}
\overset{e}{m}^{\raise-1.0ex\hbox{$\scriptstyle \infty$}}(z)(\mathcal{X}^{e}
(z))^{-1} \underset{\underset{z \in \partial \mathbb{U}_{\delta_{b_{j-1}}}^{
e}}{n \to \infty}}{=} \mathrm{I} \! + \! \mathcal{O}(n^{-1}), \quad j \! = \!
1,\dotsc,N \! + \! 1.
\end{equation*}
The solutions of the {\rm RHPs} $(\mathcal{X}^{e}(z),
\overset{e}{\upsilon}^{\raise-1.0ex\hbox{$\scriptstyle \sharp$}}(z),
\mathbb{U}_{\delta_{b_{j-1}}}^{e} \cap \Sigma_{b_{j-1}}^{e})$, $j \! = \! 1,
\dotsc,N \! + \! 1$, are:\\
{\rm \pmb{(1)}} for $z \! \in \! \Omega_{b_{j-1}}^{e,1} \! := \! \mathbb{U}_{
\delta_{b_{j-1}}}^{e} \cap (\Phi_{b_{j-1}}^{e})^{-1}(\widehat{\Omega}_{b_{j-
1}}^{e,1})$, $j \! = \! 1,\dotsc,N \! + \! 1$,
\begin{align*}
\mathcal{X}^{e}(z)=& \, \sqrt{\smash[b]{\pi}} \, \me^{-\frac{\mi \pi}{3}}
\overset{e}{m}^{\raise-1.0ex\hbox{$\scriptstyle \infty$}}(z) \sigma_{3} \me^{
\frac{\mi}{2}n \mho_{j-1}^{e} \operatorname{ad}(\sigma_{3})} \!
\begin{pmatrix}
\mi & -\mi \\
1 & 1
\end{pmatrix} \! (\Phi_{b_{j-1}}^{e}(z))^{\frac{1}{4} \sigma_{3}} \Psi^{e}_{1}
(\Phi_{b_{j-1}}^{e}(z)) \me^{\frac{1}{2}n \xi_{b_{j-1}}^{e}(z) \sigma_{3}}
\sigma_{3},
\end{align*}
where $\overset{e}{m}^{\raise-1.0ex\hbox{$\scriptstyle \infty$}}(z)$ is given
in Lemma~{\rm 4.5}, and $\Psi^{e}_{1}(z)$ and $\mho^{e}_{k}$ are defined in
Remark~{\rm 4.4;}\\
{\rm \pmb{(2)}} for $z \! \in \! \Omega_{b_{j-1}}^{e,2} \! := \! \mathbb{U}_{
\delta_{b_{j-1}}}^{e} \cap (\Phi_{b_{j-1}}^{e})^{-1}(\widehat{\Omega}_{b_{j-
1}}^{e,2})$, $j \! = \! 1,\dotsc,N \! + \! 1$,
\begin{align*}
\mathcal{X}^{e}(z)=& \, \sqrt{\smash[b]{\pi}} \, \me^{-\frac{\mi \pi}{3}}
\overset{e}{m}^{\raise-1.0ex\hbox{$\scriptstyle \infty$}}(z) \sigma_{3} \me^{
\frac{\mi}{2}n \mho_{j-1}^{e} \operatorname{ad}(\sigma_{3})} \!
\begin{pmatrix}
\mi & -\mi \\
1 & 1
\end{pmatrix} \! (\Phi_{b_{j-1}}^{e}(z))^{\frac{1}{4} \sigma_{3}} \Psi^{e}_{2}
(\Phi_{b_{j-1}}^{e}(z)) \me^{\frac{1}{2}n \xi_{b_{j-1}}^{e}(z) \sigma_{3}}
\sigma_{3},
\end{align*}
where $\Psi^{e}_{2}(z)$ is defined in Remark~{\rm 4.4;}\\
{\rm \pmb{(3)}} for $z \! \in \! \Omega_{b_{j-1}}^{e,3} \! := \! \mathbb{U}_{
\delta_{b_{j-1}}}^{e} \cap (\Phi_{b_{j-1}}^{e})^{-1}(\widehat{\Omega}_{b_{j-
1}}^{e,3})$, $j \! = \! 1,\dotsc,N \! + \! 1$,
\begin{align*}
\mathcal{X}^{e}(z)=& \, \sqrt{\smash[b]{\pi}} \, \me^{-\frac{\mi \pi}{3}}
\overset{e}{m}^{\raise-1.0ex\hbox{$\scriptstyle \infty$}}(z) \sigma_{3} \me^{-
\frac{\mi}{2}n \mho^{e}_{j-1} \operatorname{ad}(\sigma_{3})} \!
\begin{pmatrix}
\mi & -\mi \\
1 & 1
\end{pmatrix} \! (\Phi_{b_{j-1}}^{e}(z))^{\frac{1}{4} \sigma_{3}} \Psi^{e}_{3}
(\Phi_{b_{j-1}}^{e}(z)) \me^{\frac{1}{2}n \xi_{b_{j-1}}^{e}(z) \sigma_{3}}
\sigma_{3},
\end{align*}
where $\Psi^{e}_{3}(z)$ is defined in Remark~{\rm 4.4;}\\
{\rm \pmb{(4)}} for $z \! \in \! \Omega_{b_{j-1}}^{e,4} \! := \! \mathbb{U}_{
\delta_{b_{j-1}}}^{e} \cap (\Phi_{b_{j-1}}^{e})^{-1}(\widehat{\Omega}_{b_{j-
1}}^{e,4})$, $j \! = \! 1,\dotsc,N \! + \! 1$,
\begin{align*}
\mathcal{X}^{e}(z)=& \, \sqrt{\smash[b]{\pi}} \, \me^{-\frac{\mi \pi}{3}}
\overset{e}{m}^{\raise-1.0ex\hbox{$\scriptstyle \infty$}}(z) \sigma_{3} \me^{-
\frac{\mi}{2}n \mho_{j-1}^{e} \operatorname{ad}(\sigma_{3})} \!
\begin{pmatrix}
\mi & -\mi \\
1 & 1
\end{pmatrix} \! (\Phi_{b_{j-1}}^{e}(z))^{\frac{1}{4} \sigma_{3}} \Psi^{e}_{4}
(\Phi_{b_{j-1}}^{e}(z)) \me^{\frac{1}{2}n \xi_{b_{j-1}}^{e}(z) \sigma_{3}}
\sigma_{3},
\end{align*}
where $\Psi^{e}_{4}(z)$ is defined in Remark~{\rm 4.4}.
\end{ccccc}
\begin{ccccc}
Let $\overset{e}{\mathscr{M}}^{\raise-1.0ex\hbox{$\scriptstyle \sharp$}}
\colon \mathbb{C} \setminus \Sigma_{e}^{\sharp} \! \to \! \operatorname{SL}_{
2}(\mathbb{C})$ solve the {\rm RHP}
$(\overset{e}{\mathscr{M}}^{\raise-1.0ex\hbox{$\scriptstyle \sharp$}}(z),
\overset{e}{\upsilon}^{\raise-1.0ex\hbox{$\scriptstyle \sharp$}}(z),\Sigma_{
e}^{\sharp})$ formulated in Lemma~{\rm 4.2}, and set
\begin{equation*}
\mathbb{U}_{\delta_{a_{j}}}^{e} \! := \! \left\{ \mathstrut z \! \in \!
\mathbb{C}; \, \vert z \! - \! a_{j}^{e} \vert \! < \! \delta_{a_{j}}^{e} \!
\in \! (0,1) \right\}, \quad j \! = \! 1,\dotsc,N \! + \! 1.
\end{equation*}
Let
\begin{equation*}
\Phi_{a_{j}}^{e}(z) \! := \! \left(\dfrac{3n}{4} \xi_{a_{j}}^{e}(z)
\right)^{2/3}, \quad j \! = \! 1,\dotsc,N \! + \! 1,
\end{equation*}
with
\begin{equation*}
\xi_{a_{j}}^{e}(z) \! = \! 2 \int_{a_{j}^{e}}^{z}(R_{e}(s))^{1/2}h_{V}^{e}(s)
\, \md s,
\end{equation*}
where, for $z \! \in \! \mathbb{U}_{\delta_{a_{j}}}^{e} \setminus (-\infty,
a_{j}^{e})$, $\xi_{a_{j}}^{e}(z) \! = \! (z \! - \! a_{j}^{e})^{3/2}G_{a_{j}
}^{e}(z)$, $j \! = \! 1,\dotsc,N \! + \! 1$, with $G_{a_{j}}^{e}(z)$ analytic,
in particular,
\begin{equation*}
G_{a_{j}}^{e}(z) \! \underset{z \to a_{j}^{e}}{=} \! \dfrac{4}{3}f(a_{j}^{e})
\! + \! \dfrac{4}{5}f^{\prime}(a_{j}^{e})(z \! - \! a_{j}^{e}) \! + \! \dfrac{
2}{7}f^{\prime \prime}(a_{j}^{e})(z \! - \! a_{j}^{e})^{2} \! + \! \mathcal{O}
\! \left((z \! - \! a_{j}^{e})^{3} \right),
\end{equation*}
where
\begin{align*}
f(a_{N+1}^{e})=& \, h_{V}^{e}(a_{N+1}^{e}) \eta_{a_{N+1}^{e}}, \\
f^{\prime}(a_{N+1}^{e})=& \, \dfrac{1}{2}h_{V}^{e}(a_{N+1}^{e}) \eta_{a_{N+
1}^{e}} \! \left(\sum_{l=1}^{N} \! \left(\dfrac{1}{a_{N+1}^{e} \! - \! b_{l}^{
e}} \! + \! \dfrac{1}{a_{N+1}^{e} \! - \! a_{l}^{e}} \right) \! + \! \dfrac{
1}{a_{N+1}^{e} \! - \! b_{0}^{e}} \right) \\
+& \, (h_{V}^{e}(a_{N+1}^{e}))^{\prime} \eta_{a_{N+1}^{e}}, \\
f^{\prime \prime}(a_{N+1}^{e})=& \, \dfrac{h_{V}^{e}(a_{N+1}^{e})(h_{V}^{e}
(a_{N+1}^{e}))^{\prime \prime} \! - \! ((h_{V}^{e}(a_{N+1}^{e}))^{\prime})^{
2}}{h_{V}^{e}(a_{N+1}^{e})} \eta_{a_{N+1}^{e}} \! - \! \dfrac{1}{2}h_{V}^{e}
(a_{N+1}^{e}) \eta_{a_{N+1}^{e}} \\
\times& \, \left(\sum_{l=1}^{N} \! \left(\dfrac{1}{(a_{N+1}^{e} \! - \! b_{
l}^{e})^{2}} \! + \! \dfrac{1}{(a_{N+1}^{e} \! - \! a_{l}^{e})^{2}} \right) \!
+ \! \dfrac{1}{(a_{N+1}^{e} \! - \! b_{0}^{e})^{2}} \right) \\
+& \, \left(\dfrac{1}{2} \! \left(\sum_{k=1}^{N} \! \left(\dfrac{1}{a_{N+1}^{
e} \! - \! b_{k}^{e}} \! + \! \dfrac{1}{a_{N+1}^{e} \! - \! a_{k}^{e}} \right)
\! + \! \dfrac{1}{a_{N+1}^{e} \! - \! b_{0}^{e}} \right) \! + \! \dfrac{(h_{
V}^{e}(a_{N+1}^{e}))^{\prime}}{h_{V}^{e}(a_{N+1}^{e})} \right) \\
\times& \, \left(\dfrac{1}{2}h_{V}^{e}(a_{N+1}^{e}) \eta_{a_{N+1}^{e}} \!
\left(\sum_{l=1}^{N} \! \left(\dfrac{1}{a_{N+1}^{e} \! - \! a_{l}^{e}} \! + \!
\dfrac{1}{a_{N+1}^{e} \! - \! b_{l}^{e}} \right) \! + \! \dfrac{1}{a_{N+1}^{e}
\! - \! b_{0}^{e}} \right) \right. \\
+& \left. \, (h_{V}^{e}(a_{N+1}^{e}))^{\prime} \eta_{a_{N+1}^{e}} \right),
\end{align*}
with
\begin{equation*}
\eta_{a_{N+1}^{e}} \! := \! \left((a_{N+1}^{e} \! - \! b_{0}^{e}) \prod_{k=1
}^{N}(a_{N+1}^{e} \! - \! b_{k}^{e})(a_{N+1}^{e} \! - \! a_{k}^{e}) \right)^{
1/2} \quad (> \! 0),
\end{equation*}
and, for $j \! = \! 1,\dotsc,N$,
\begin{align*}
f(a_{j}^{e})=& \, (-1)^{N-j+1}h_{V}^{e}(a_{j}^{e}) \eta_{a_{j}^{e}}, \\
f^{\prime}(a_{j}^{e})=& \, (-1)^{N-j+1} \! \left(\dfrac{1}{2}h_{V}^{e}(a_{j}^{
e}) \eta_{a_{j}^{e}} \! \left(\sum_{\substack{k=1\\k \not= j}}^{N} \! \left(
\dfrac{1}{a_{j}^{e} \! - \! b_{k}^{e}} \! + \! \dfrac{1}{a_{j}^{e} \! - \!
a_{k}^{e}} \right) \! + \! \dfrac{1}{a_{j}^{e} \! - \! b_{j}^{e}} \! + \!
\dfrac{1}{a_{j}^{e} \! - \! a_{N+1}^{e}} \! + \! \dfrac{1}{a_{j}^{e} \! - \!
b_{0}^{e}} \right) \right. \\
+&\left. \, (h_{V}^{e}(a_{j}^{e}))^{\prime} \eta_{a_{j}^{e}} \right), \\
f^{\prime \prime}(a_{j}^{e})=& \, (-1)^{N-j+1} \! \left(\dfrac{h_{V}^{e}(a_{
j}^{e})(h_{V}^{e}(a_{j}^{e}))^{\prime \prime} \! - \! ((h_{V}^{e}(a_{j}^{e}))^{
\prime})^{2}}{h_{V}^{e}(a_{j}^{e})} \eta_{a_{j}^{e}} \! - \! \dfrac{1}{2}h_{
V}^{e}(a_{j}^{e}) \eta_{a_{j}^{e}} \! \left(\sum_{\substack{k=1\\k \not= j}}^{
N} \! \left(\dfrac{1}{(a_{j}^{e} \! - \! b_{k}^{e})^{2}} \! + \! \dfrac{1}{(a_{
j}^{e} \! - \! a_{k}^{e})^{2}} \right) \right. \right. \\
+&\left. \left. \, \dfrac{1}{(a_{j}^{e} \! - \! b_{j}^{e})^{2}} \! + \! \dfrac{
1}{(a_{j}^{e} \! - \! a_{N+1}^{e})^{2}} \! + \! \dfrac{1}{(a_{j}^{e} \! - \!
b_{0}^{e})^{2}} \right) \! + \! \left(\dfrac{(h_{V}^{e}(a_{j}^{e}))^{\prime}}{
h_{V}^{e}(a_{j}^{e})} \! + \! \dfrac{1}{2} \! \left(\sum_{\substack{k=1\\k
\not= j}}^{N} \! \left(\dfrac{1}{a_{j}^{e} \! - \! b_{k}^{e}} \! + \! \dfrac{
1}{a_{j}^{e} \! - \! a_{k}^{e}} \right) \right. \right. \right. \\
+&\left. \left. \left. \, \dfrac{1}{a_{j}^{e} \! - \! b_{j}^{e}} \! + \!
\dfrac{1}{a_{j}^{e} \! - \! a_{N+1}^{e}} \! + \! \dfrac{1}{a_{j}^{e} \! - \!
b_{0}^{e}} \right) \! \right) \! \left(\dfrac{1}{2}h_{V}^{e}(a_{j}^{e}) \eta_{
a_{j}^{e}} \! \left(\sum_{\substack{k=1\\k \not= j}}^{N} \! \left(\dfrac{1}{
a_{j}^{e} \! - \! b_{k}^{e}} \! + \! \dfrac{1}{a_{j}^{e} \! - \! a_{k}^{e}}
\right) \right. \right. \right. \\
+&\left. \left. \left. \, \dfrac{1}{a_{j}^{e} \! - \! b_{j}^{e}} \! + \!
\dfrac{1}{a_{j}^{e} \! - \! a_{N+1}^{e}} \! + \! \dfrac{1}{a_{j}^{e} \! - \!
b_{0}^{e}} \right) \! + \!  (h_{V}^{e}(a_{j}^{e}))^{\prime} \eta_{a_{j}^{e}}
\right) \right),
\end{align*}
with
\begin{equation*}
\eta_{a_{j}^{e}} \! := \! \left((b_{j}^{e} \! - \! a_{j}^{e})(a_{N+1}^{e} \!
- \! a_{j}^{e})(a_{j}^{e} \! - \! b_{0}^{e}) \prod_{k=1}^{j-1}(a_{j}^{e} \! -
\! b_{k}^{e})(a_{j}^{e} \! - \! a_{k}^{e}) \prod_{l=j+1}^{N}(b_{l}^{e} \! -
\! a_{j}^{e})(a_{l}^{e} \! - \! a_{j}^{e}) \right)^{1/2} \quad (> \! 0),
\end{equation*}
and $((0,1) \! \ni)$ $\delta_{a_{j}}^{e}$, $j \! = \! 1,\dotsc,N \! + \! 1$,
are chosen sufficiently small so that $\Phi_{a_{j}}^{e}(z)$, which are
bi-holomo\-r\-p\-h\-i\-c, conformal, and orientation preserving, map $\mathbb{
U}_{\delta_{a_{j}}}^{e}$ (and, thus, the oriented contours $\Sigma_{a_{j}}^{e}
\! := \! \cup_{l=1}^{4} \Sigma_{a_{j}}^{e,l}$, $j \! = \! 1,\dotsc,N \! + \!
1:$ Figure~{\rm 5)} injectively onto open $(n$-dependent) neighbourhoods
$\widehat{\mathbb{U}}_{\delta_{a_{j}}}^{e}$, $j \! = \! 1,\dotsc,N \! + \! 1$,
of $0$ such that $\Phi_{a_{j}}^{e}(a_{j}^{e}) \! = \! 0$, $\Phi_{a_{j}}^{e}
\colon \mathbb{U}_{\delta_{a_{j}}}^{e} \! \to \! \widehat{\mathbb{U}}_{\delta_{
a_{j}}}^{e} \! := \! \Phi_{a_{j}}^{e}(\mathbb{U}_{\delta_{a_{j}}}^{e})$,
$\Phi_{a_{j}}^{e}(\mathbb{U}_{\delta_{a_{j}}}^{e} \cap \Sigma_{a_{j}}^{e,l})
\! = \! \Phi_{a_{j}}^{e}(\mathbb{U}_{\delta_{a_{j}}}^{e}) \cap \gamma_{a_{j}
}^{e,l}$, and $\Phi_{a_{j}}^{e}(\mathbb{U}_{\delta_{a_{j}}}^{e} \cap \Omega_{
a_{j}}^{e,l}) \! = \! \Phi_{a_{j}}^{e}(\mathbb{U}_{\delta_{a_{j}}}^{e}) \cap
\widehat{\Omega}_{a_{j}}^{e,l}$, $l \! = \! 1,2,3,4$, with $\widehat{\Omega}_{
a_{j}}^{e,1} \! = \! \{\mathstrut \zeta \! \in \! \mathbb{C}; \, \arg (\zeta)
\! \in \! (0,2 \pi/3)\}$, $\widehat{\Omega}_{a_{j}}^{e,2} \! = \! \{\mathstrut
\zeta \! \in \! \mathbb{C}; \, \arg (\zeta) \! \in \! (2 \pi/3,\pi)\}$,
$\widehat{\Omega}_{a_{j}}^{e,3} \! = \! \{\mathstrut \zeta \! \in \! \mathbb{
C}; \, \arg (\zeta) \! \in \! (-\pi,-2 \pi/3)\}$, and $\widehat{\Omega}_{a_{
j}}^{e,4} \! = \! \{\mathstrut \zeta \! \in \! \mathbb{C}; \, \arg (\zeta) \!
\in \! (-2 \pi/3,0)\}$.

The parametrix for the {\rm RHP}
$(\overset{e}{\mathscr{M}}^{\raise-1.0ex\hbox{$\scriptstyle \sharp$}}(z),
\overset{e}{\upsilon}^{\raise-1.0ex\hbox{$\scriptstyle \sharp$}}(z),\Sigma_{
e}^{\sharp})$, for $z \! \in \! \mathbb{U}_{\delta_{a_{j}}}^{e}$, $j \! = \!
1,\dotsc,N \! + \! 1$, is the solution of the following {\rm RHPs} for
$\mathcal{X}^{e} \colon \mathbb{U}_{\delta_{a_{j}}}^{e} \setminus \Sigma_{a_{
j}}^{e} \! \to \! \operatorname{SL}_{2}(\mathbb{C})$, $j \! = \! 1,\dotsc,N \!
+ \! 1$, where $\Sigma_{a_{j}}^{e} \! := \! (\Phi_{a_{j}}^{e})^{-1}(\gamma_{a_{
j}}^{e})$, with $(\Phi_{a_{j}}^{e})^{-1}$ denoting the inverse mapping, and
$\gamma_{a_{j}}^{e} \! := \! \cup_{l=1}^{4} \gamma_{a_{j}}^{e,l}:$ {\rm (i)}
$\mathcal{X}^{e}(z)$ is holomorphic for $z \! \in \! \mathbb{U}_{\delta_{a_{
j}}}^{e} \setminus \Sigma_{a_{j}}^{e}$, $j \! = \! 1,\dotsc,N \! + \! 1;$
{\rm (ii)} $\mathcal{X}^{e}_{\pm}(z) \! := \! \lim_{\underset{z^{\prime} \,
\in \, \pm \, \mathrm{side} \, \mathrm{of} \, \Sigma_{a_{j}}^{e}}{z^{\prime}
\to z}} \mathcal{X}^{e}(z^{\prime})$, $j \! = \! 1,\dotsc,N \! + \! 1$,
satisfy the boundary condition
\begin{equation*}
\mathcal{X}^{e}_{+}(z) \! = \! \mathcal{X}^{e}_{-}(z)
\overset{e}{\upsilon}^{\raise-1.0ex\hbox{$\scriptstyle \sharp$}}(z), \quad z
\! \in \! \mathbb{U}_{\delta_{a_{j}}}^{e} \cap \Sigma_{a_{j}}^{e}, \quad j \!
= \! 1,\dotsc,N \! + \! 1,
\end{equation*}
where $\overset{e}{\upsilon}^{\raise-1.0ex\hbox{$\scriptstyle \sharp$}}(z)$ is
given in Lemma~{\rm 4.2;} and {\rm (iii)} uniformly for $z \! \in \! \partial
\mathbb{U}_{\delta_{a_{j}}}^{e} \! := \! \left\lbrace \mathstrut z \! \in \!
\mathbb{C}; \, \vert z \! - \! a_{j}^{e} \vert \! = \! \delta_{a_{j}}^{e}
\right\rbrace$, $j \! = \! 1,\dotsc,N \! + \! 1$,
\begin{equation*}
\overset{e}{m}^{\raise-1.0ex\hbox{$\scriptstyle \infty$}}(z)(\mathcal{X}^{e}
(z))^{-1} \underset{\underset{z \in \partial \mathbb{U}_{\delta_{a_{j}}}^{e}
}{n \to \infty}}{=} \mathrm{I} \! + \! \mathcal{O}(n^{-1}), \quad j \! = \!
1,\dotsc,N \! + \! 1.
\end{equation*}
The solutions of the {\rm RHPs} $(\mathcal{X}^{e}(z),
\overset{e}{\upsilon}^{\raise-1.0ex\hbox{$\scriptstyle \sharp$}}(z),
\mathbb{U}_{\delta_{a_{j}}}^{e} \cap \Sigma_{a_{j}}^{e})$, $j \! = \! 1,
\dotsc,N \! + \! 1$, are:\\
{\rm \pmb{(1)}} for $z \! \in \! \Omega_{a_{j}}^{e,1} \! := \! \mathbb{U}_{
\delta_{a_{j}}}^{e} \cap (\Phi_{a_{j}}^{e})^{-1}(\widehat{\Omega}_{a_{j}}^{e,
1})$, $j \! = \! 1,\dotsc,N \! + \! 1$,
\begin{align*}
\mathcal{X}^{e}(z)=& \, \sqrt{\smash[b]{\pi}} \, \me^{-\frac{\mi \pi}{3}}
\overset{e}{m}^{\raise-1.0ex\hbox{$\scriptstyle \infty$}}(z) \me^{\frac{\mi}{2}
n \mho_{j}^{e} \operatorname{ad}(\sigma_{3})} \!
\begin{pmatrix}
\mi & -\mi \\
1 & 1
\end{pmatrix} \! (\Phi_{a_{j}}^{e}(z))^{\frac{1}{4} \sigma_{3}} \Psi^{e}_{1}
(\Phi_{a_{j}}^{e}(z)) \me^{\frac{1}{2}n \xi_{a_{j}}^{e}(z) \sigma_{3}},
\end{align*}
where $\overset{e}{m}^{\raise-1.0ex\hbox{$\scriptstyle \infty$}}(z)$ is given
in Lemma~{\rm 4.5}, and $\Psi^{e}_{1}(z)$ and $\mho^{e}_{k}$ are defined
in Remark~{\rm 4.4;}\\
{\rm \pmb{(2)}} for $z \! \in \! \Omega_{a_{j}}^{e,2} \! := \! \mathbb{U}_{
\delta_{a_{j}}}^{e} \cap (\Phi_{a_{j}}^{e})^{-1}(\widehat{\Omega}_{a_{j}}^{e,
2})$, $j \! = \! 1,\dotsc,N \! + \! 1$,
\begin{align*}
\mathcal{X}^{e}(z)=& \, \sqrt{\smash[b]{\pi}} \, \me^{-\frac{\mi \pi}{3}}
\overset{e}{m}^{\raise-1.0ex\hbox{$\scriptstyle \infty$}}(z) \me^{\frac{\mi}{2}
n \mho_{j}^{e} \operatorname{ad}(\sigma_{3})} \!
\begin{pmatrix}
\mi & -\mi \\
1 & 1
\end{pmatrix} \! (\Phi_{a_{j}}^{e}(z))^{\frac{1}{4} \sigma_{3}} \Psi^{e}_{2}
(\Phi_{a_{j}}^{e}(z)) \me^{\frac{1}{2}n \xi_{a_{j}}^{e}(z) \sigma_{3}},
\end{align*}
where $\Psi^{e}_{2}(z)$ is defined in Remark~{\rm 4.4;}\\
{\rm \pmb{(3)}} for $z \! \in \! \Omega_{a_{j}}^{e,3} \! := \! \mathbb{U}_{
\delta_{a_{j}}}^{e} \cap (\Phi_{a_{j}}^{e})^{-1}(\widehat{\Omega}_{a_{j}}^{e,
3})$, $j \! = \! 1,\dotsc,N \! + \! 1$,
\begin{align*}
\mathcal{X}^{e}(z)=& \, \sqrt{\smash[b]{\pi}} \, \me^{-\frac{\mi \pi}{3}}
\overset{e}{m}^{\raise-1.0ex\hbox{$\scriptstyle \infty$}}(z) \me^{-\frac{\mi}{
2}n \mho^{e}_{j} \operatorname{ad}(\sigma_{3})} \!
\begin{pmatrix}
\mi & -\mi \\
1 & 1
\end{pmatrix} \! (\Phi_{a_{j}}^{e}(z))^{\frac{1}{4} \sigma_{3}} \Psi^{e}_{3}
(\Phi_{a_{j}}^{e}(z)) \me^{\frac{1}{2}n \xi_{a_{j}}^{e}(z) \sigma_{3}},
\end{align*}
where $\Psi^{e}_{3}(z)$ is defined in Remark~{\rm 4.4;}\\
{\rm \pmb{(4)}} for $z \! \in \! \Omega_{a_{j}}^{e,4} \! := \! \mathbb{U}_{
\delta_{a_{j}}}^{e} \cap (\Phi_{a_{j}}^{e})^{-1}(\widehat{\Omega}_{a_{j}}^{e,
4})$, $j \! = \! 1,\dotsc,N \! + \! 1$,
\begin{align*}
\mathcal{X}^{e}(z)=& \, \sqrt{\smash[b]{\pi}} \, \me^{-\frac{\mi \pi}{3}}
\overset{e}{m}^{\raise-1.0ex\hbox{$\scriptstyle \infty$}}(z) \me^{-\frac{\mi}{
2}n \mho_{j}^{e} \operatorname{ad}(\sigma_{3})} \!
\begin{pmatrix}
\mi & -\mi \\
1 & 1
\end{pmatrix} \! (\Phi_{a_{j}}^{e}(z))^{\frac{1}{4} \sigma_{3}} \Psi^{e}_{4}
(\Phi_{a_{j}}^{e}(z)) \me^{\frac{1}{2}n \xi_{a_{j}}^{e}(z) \sigma_{3}},
\end{align*}
where $\Psi^{e}_{4}(z)$ is defined in Remark~{\rm 4.4}.
\end{ccccc}
\begin{eeeee}
Perusing Lemmas~4.6 and~4.7, one notes that the normalisation condition at
infinity, which is needed in order to guarantee the existence of solutions to
the corresponding (parametrix) RHPs, is absent. The normalisation conditions
at infinity are replaced by the (uniform) matching conditions
$\overset{e}{m}^{\raise-1.0ex\hbox{$\scriptstyle \infty$}}(z)(\mathcal{X}^{e}
(z))^{-1} \! =_{\underset{z \in \partial \mathbb{U}_{\delta_{\ast_{j}}}^{e}}{n
\to \infty}} \! \mathrm{I} \! + \! \mathcal{O}(n^{-1})$, where $\ast_{j} \!
\in \! \{b_{j-1},a_{j}\}$, $j \! = \! 1,\dotsc,N \! + \! 1$, with $\partial
\mathbb{U}_{\delta_{\ast_{j}}}^{e}$ defined in Lemmas~4.6 and~4.7. \hfill
$\blacksquare$
\end{eeeee}

\emph{Sketch of proof of Lemma~{\rm 4.7}.} Let
$(\overset{e}{\mathscr{M}}^{\raise-1.0ex\hbox{$\scriptstyle \sharp$}}(z),
\overset{e}{\upsilon}^{\raise-1.0ex\hbox{$\scriptstyle \sharp$}}(z),\Sigma_{
e}^{\sharp})$ be the RHP formulated in Lemma~4.2, and recall the definitions
stated therein. For each $a_{j}^{e} \! \in \! \operatorname{supp}(\mu_{V}^{e}
)$, $j \! = \! 1,\dotsc,N \! + \! 1$, define $\mathbb{U}_{\delta_{a_{j}}}^{
e}$, $j \! = \! 1,\dotsc,N \! + \! 1$, as in the Lemma, that is, surround each
right-most end-point $a_{j}^{e}$ by open discs of radius $\delta_{a_{j}}^{e}
\! \in \! (0,1)$ centred at $a_{j}^{e}$. Recalling the formula for $\overset{
e}{\upsilon}^{\raise-1.0ex\hbox{$\scriptstyle \sharp$}}(z)$ given in
Lemma~4.2, one shows, via the proof of Lemma~4.1, that:
\begin{compactenum}
\item[(1)] $4 \pi \mi \int_{z}^{a_{N+1}^{e}} \psi_{V}^{e}(s) \, \md s \! = \!
4 \pi \mi (\int_{z}^{a_{j}^{e}} \! + \! \int_{a_{j}^{e}}^{b_{j}^{e}} \! + \!
\int_{b_{j}^{e}}^{a_{N+1}^{e}}) \psi_{V}^{e}(s) \, \md s$, whence, recalling
the expression for the density of the `even' equilibrium measure given in
Lemma~3.5, that is, $\md \mu_{V}^{e}(x) \! := \! \psi_{V}^{e}(x) \, \md x \! =
\! \tfrac{1}{2 \pi \mi}(R_{e}(x))^{1/2}_{+}h_{V}^{e}(x) \pmb{1}_{J_{e}}(x) \,
\md x$, one arrives at, upon considering the analytic continuation of $4 \pi
\mi \linebreak[4]
\cdot \int_{z}^{a_{N+1}^{e}} \psi_{V}^{e}(s) \, \md s$ to $\mathbb{C}
\setminus \mathbb{R}$ (cf. proof of Lemma~4.1), in particular, to the oriented
(open) skeletons $\mathbb{U}_{\delta_{a_{j}}}^{e} \cap (J_{j}^{e,\smallfrown}
\cup J_{j}^{e,\smallsmile})$, $j \! = \! 1,\dotsc,N \! + \! 1$, $4 \pi \mi
\int_{z}^{a_{N+1}^{e}} \psi_{V}^{e}(s) \, \md s \! = \! -\xi_{a_{j}}^{e}(z) \!
+ \! \mi \mho_{j}^{e}$, $j \! = \! 1,\dotsc,N \! + \! 1$, where $\xi_{a_{j}}^{
e}(z) \! = \! 2 \int_{a_{j}^{e}}^{z}(R_{e}(s))^{1/2}h_{V}^{e}(s) \, \md s$,
and $\mho_{j}^{e}$ are defined in Remark~4.4;
\item[(2)] $g^{e}_{+}(z) \! + \! g^{e}_{-}(z) \! - \! \widetilde{V}(z) \! - \!
\ell_{e} \! + \! 2Q_{e} \! = \! -2 \int_{a_{j}^{e}}^{z}(R_{e}(s))^{1/2}h_{V}^{
e}(s) \, \md s \! < \! 0$, $z \! \in \! (a_{N+1}^{e},+\infty) \cup (\cup_{j=
1}^{N}(a_{j}^{e},b_{j}^{e}))$.
\end{compactenum}
Via the latter formulae, which appear in the $(i \, j)$-elements, $i,j \!
= \! 1,2$, of the jump matrix
$\overset{e}{\upsilon}^{\raise-1.0ex\hbox{$\scriptstyle \sharp$}}(z)$,
denoting
$\overset{e}{\mathscr{M}}^{\raise-1.0ex\hbox{$\scriptstyle \sharp$}}(z)$ by
$\mathcal{X}^{e}(z)$ for $z \! \in \! \mathbb{U}_{\delta_{a_{j}}}^{e}$, $j \!
= \! 1,\dotsc,N \! + \! 1$, and defining
\begin{equation*}
\mathscr{P}_{a_{j}}^{e}(z) \! := \!
\begin{cases}
\mathcal{X}^{e}(z) \me^{-\frac{1}{2}n \xi_{a_{j}}^{e}(z) \sigma_{3}} \, \me^{
\frac{\mi}{2}n \mho_{j}^{e} \sigma_{3}}, &\text{$z \! \in \! \mathbb{C}_{+}
\cap \mathbb{U}_{\delta_{a_{j}}}^{e}, \quad j \! = \! 1,\dotsc,N \! + \! 1$,}
\\
\mathcal{X}^{e}(z) \me^{-\frac{1}{2}n \xi_{a_{j}}^{e}(z) \sigma_{3}} \, \me^{
-\frac{\mi}{2}n \mho_{j}^{e} \sigma_{3}}, &\text{$z \! \in \! \mathbb{C}_{-}
\cap \mathbb{U}_{\delta_{a_{j}}}^{e}, \quad j \! = \! 1,\dotsc,N \! + \! 1$,}
\end{cases}
\end{equation*}
one notes that $\mathscr{P}_{a_{j}}^{e} \colon \mathbb{U}_{\delta_{a_{j}}}^{e}
\setminus J_{a_{j}}^{e} \! \to \! \operatorname{GL}_{2}(\mathbb{C})$, where
$J_{a_{j}}^{e} \! := \! J_{j}^{e,\smallfrown} \cup J_{j}^{e,\smallsmile} \cup
(a_{j}^{e} \! - \! \delta_{a_{j}}^{e},a_{j}^{e} \! + \! \delta_{a_{j}}^{e})$,
$j \! = \! 1,\dotsc,N \! + \! 1$, solve the RHPs $(\mathscr{P}_{a_{j}}^{e}(z),
\upsilon^{e}_{\mathscr{P}_{a_{j}}}(z),J_{a_{j}}^{e})$, with constant jump
matrices $\upsilon_{\mathscr{P}_{a_{j}}}^{e}(z)$, $j \! = \! 1,\dotsc,N \! +
\! 1$, defined by
\begin{equation*}
\upsilon_{\mathscr{P}_{a_{j}}}^{e}(z) \! := \!
\begin{cases}
\mathrm{I} \! + \! \sigma_{-}, &\text{$z \! \in \! \mathbb{U}_{\delta_{a_{j}}
}^{e} \cap (J_{j}^{e,\smallfrown} \cup J_{j}^{e,\smallsmile}) \! = \! \Sigma^{
e,1}_{a_{j}} \cup \Sigma^{e,3}_{a_{j}}$,} \\
\mathrm{I} \! + \! \sigma_{+}, &\text{$z \! \in \! \mathbb{U}_{\delta_{a_{j}}
}^{e} \cap (a_{j}^{e},a_{j}^{e} \! + \! \delta_{a_{j}}^{e}) \! = \! \Sigma^{e,
4}_{a_{j}}$,} \\
\mi \sigma_{2}, &\text{$z \! \in \! \mathbb{U}_{\delta_{a_{j}}}^{e} \cap (a_{
j}^{e} \! - \! \delta_{a_{j}}^{e},a_{j}^{e}) \! = \! \Sigma^{e,2}_{a_{j}}$,}
\end{cases}
\end{equation*}
subject, still, to the asymptotic matching conditions
$\overset{e}{m}^{\raise-1.0ex\hbox{$\scriptstyle \infty$}}(z)(\mathcal{X}^{e}
(z))^{-1} \! =_{n \to \infty} \! \mathrm{I} \! + \! \mathcal{O}(n^{-1})$,
uniformly for $z \! \in \! \partial \mathbb{U}_{\delta_{a_{j}}}^{e}$, $j \! =
\! 1,\dotsc,N \! + \! 1$.

Set, as in the Lemma, $\Phi_{a_{j}}^{e}(z) \! := \! (\tfrac{3}{4}n \xi_{a_{j}
}^{e}(z))^{2/3}$, $j \! = \! 1,\dotsc,N \! + \! 1$, with $\xi_{a_{j}}^{e}(z)$
defined above: a careful analysis of the branch cuts shows that, for $z \!
\in \! \mathbb{U}_{\delta_{a_{j}}}^{e}$, $j \! = \! 1,\dotsc,N \! + \! 1$,
$\Phi_{a_{j}}^{e}(z)$ and $\xi_{a_{j}}^{e}(z)$ satisfy the properties stated
in the Lemma; in particular, for $\Phi_{a_{j}}^{e} \colon \mathbb{U}_{\delta_{
a_{j}}}^{e} \! \to \! \mathbb{C}$, $j \! = \! 1,\dotsc,N \! + \! 1$, $\Phi_{
a_{j}}^{e}(z) \! = \! (z \! - \! a_{j}^{e})^{3/2}G_{a_{j}}^{e}(z)$, with $G_{
a_{j}}^{e}(z)$ holomorphic for $z \! \in \! \mathbb{U}_{\delta_{a_{j}}}^{e}$
and characterised in the Lemma, $\Phi_{a_{j}}^{e}(a_{j}^{e}) \! = \! 0$,
$(\Phi_{a_{j}}^{e}(z))^{\prime} \! \not= \! 0$, $z \! \in \! \mathbb{U}_{
\delta_{a_{j}}}^{e}$, and where $(\Phi_{a_{j}}^{e}(a_{j}^{e}))^{\prime} \! =
\! (nf(a_{j}^{e}))^{2/3} \! > \! 0$, with $f(a_{j}^{e})$ given in the Lemma.
One now chooses $\delta_{a_{j}}^{e}$ $(\in \! (0,1))$, $j \! = \! 1,\dotsc,N
\! + \! 1$, and the oriented---open---skeletons (`near' $a_{j}^{e})$ $J_{a_{
j}}^{e}$, $j \! = \! 1,\dotsc,N \! + \! 1$, in such a way that their images
under the bi-holomorphic, conformal, and orientation-preserving mappings
$\Phi_{a_{j}}^{e}(z)$ are the union of the straight-line segments $\gamma_{
a_{j}}^{e,l}$, $l \! = \! 1,2,3,4$, $j \! = \! 1,\dotsc,N \! + \! 1$. Set
$\zeta \! := \! \Phi_{a_{j}}^{e}(z)$, $j \! = \! 1,\dotsc,N \! + \! 1$, and
consider $\mathcal{X}^{e}(\Phi_{a_{j}}^{e}(z)) \! := \! \Psi^{e}(\zeta)$.
Recalling the properties of $\Phi_{a_{j}}^{e}(z)$, a straightforward
calculation shows that $\Psi^{e} \colon \Phi_{a_{j}}^{e}(\mathbb{U}_{\delta_{
a_{j}}}^{e}) \setminus \cup_{l=1}^{4} \gamma_{a_{j}}^{e,l} \! \to \!
\operatorname{GL}_{2}(\mathbb{C})$, $j \! = \! 1,\dotsc,N \! + \! 1$, solves
the RHPs $(\Psi^{e}(\zeta),\upsilon_{\Psi^{e}}^{e}(\zeta),\cup_{l=1}^{4}
\gamma_{a_{j}}^{e,l})$, $j \! = \! 1,\dotsc,N \! + \! 1$, with constant jump
matrices $\upsilon_{\Psi^{e}}^{e}(\zeta)$, $j \! = \! 1,\dotsc,N \! + \! 1$,
defined by
\begin{equation*}
\upsilon_{\Psi^{e}}^{e}(\zeta) \! := \!
\begin{cases}
\mathrm{I} \! + \! \sigma_{-}, &\text{$\zeta \! \in \! \gamma_{a_{j}}^{e,1}
\cup \gamma_{a_{j}}^{e,3}$,} \\
\mathrm{I} \! + \! \sigma_{+}, &\text{$\zeta \! \in \! \gamma_{a_{j}}^{e,4}$,}
\\
\mi \sigma_{2}, &\text{$\zeta \! \in \! \gamma_{a_{j}}^{e,2}$.}
\end{cases}
\end{equation*}
The solution of the latter (yet-to-be normalised) RHPs is well known; in fact,
their solution is expressed in terms of the Airy function, and is given by
(see, for example, \cite{a3,a58,a59,a61,a90})
\begin{equation*}
\Psi^{e}(\zeta) \! = \!
\begin{cases}
\Psi^{e}_{1}(\zeta), &\text{$\zeta \! \in \! \widehat{\Omega}_{a_{j}}^{e,1},
\quad j \! = \! 1,\dotsc,N \! + \! 1$,} \\
\Psi^{e}_{2}(\zeta), &\text{$\zeta \! \in \! \widehat{\Omega}_{a_{j}}^{e,2},
\quad j \! = \! 1,\dotsc,N \! + \! 1$,} \\
\Psi^{e}_{3}(\zeta), &\text{$\zeta \! \in \! \widehat{\Omega}_{a_{j}}^{e,3},
\quad j \! = \! 1,\dotsc,N \! + \! 1$,} \\
\Psi^{e}_{4}(\zeta), &\text{$\zeta \! \in \! \widehat{\Omega}_{a_{j}}^{e,4},
\quad j \! = \! 1,\dotsc,N \! + \! 1$,}
\end{cases}
\end{equation*}
where $\Psi^{e}_{k}(z)$, $k \! = \! 1,2,3,4$, are defined in Remark~4.4. 
Recalling that $\Phi_{a_{j}}^{e}(z)$, $j \! = \! 1,\dotsc,N \! + \! 1$, are 
bi-holomorphic and orientation-preserving conformal mappings, with $\Phi_{a_{
j}}^{e}(a_{j}^{e}) \! = \! 0$ and $\Phi_{a_{j}}^{e}$ $(\colon \mathbb{U}_{
\delta_{a_{j}}}^{e} \! \to \! \Phi_{a_{j}}^{e}(\mathbb{U}_{\delta_{a_{j}}}^{
e}) \! =: \! \widehat{\mathbb{U}}_{\delta_{a_{j}}}^{e})$ $\colon \mathbb{U}_{
\delta_{a_{j}}}^{e} \cap J_{a_{j}}^{e} \! \to \! \Phi_{a_{j}}^{e}(\mathbb{U}_{
\delta_{a_{j}}}^{e} \cap J_{a_{j}}^{e}) \! = \! \widehat{\mathbb{U}}_{\delta_{
a_{j}}}^{e} \cap (\cup_{l=1}^{4} \gamma_{a_{j}}^{e,l})$, $j \! = \! 1,\dotsc,
N \! + \! 1$, one notes that, for any analytic maps $E_{a_{j}}^{e} \colon 
\mathbb{U}_{\delta_{a_{j}}}^{e} \! \to \! \operatorname{GL}_{2}(\mathbb{C})$, 
$j \! = \! 1,\dotsc,N \! + \! 1$, $\mathbb{U}_{\delta_{a_{j}}}^{e} \setminus 
J_{a_{j}}^{e} \! \ni \! \zeta \! \mapsto \! E_{a_{j}}^{e}(\zeta) \Psi^{e}
(\zeta)$ also solves the latter RHPs $(\Psi^{e}(\zeta),\upsilon^{e}_{\Psi^{
e}}(\zeta),\cup_{l=1}^{4} \gamma_{a_{j}}^{e,l})$, $j \! = \! 1,\dotsc,N \! 
+ \! 1$: one uses this `degree of freedom' of `multiplying on the left' by 
a non-degenerate, analytic, matrix-valued function in order to satisfy the 
remaining asymptotic (as $n \! \to \! \infty)$ matching condition for the 
parametrix, namely, 
$\overset{e}{m}^{\raise-1.0ex\hbox{$\scriptstyle \infty$}}(z)(\mathcal{X}^{e}
(z))^{-1} \! =_{\underset{z \in \partial \mathbb{U}_{\delta_{a_{j}}}^{e}}{n
\to \infty}} \! \mathrm{I} \! + \! \mathcal{O}(n^{-1})$, uniformly for $z \!
\in \! \partial \mathbb{U}_{\delta_{a_{j}}}^{e}$, $j \! = \! 1,\dotsc,N \! +
\! 1$.

Consider, say, and without loss of generality, the regions $\Omega_{a_{j}}^{e,
1} \! := \! (\Phi_{a_{j}}^{e})^{-1}(\widehat{\Omega}_{a_{j}}^{e,1})$, $j \! =
\! 1,\dotsc,N \! + \! 1$ (Figure~5). Re-tracing the above transformations, one
shows that, for $z \! \in \! \Omega_{a_{j}}^{e,1}$ $(\subset \! \mathbb{C}_{
+})$, $j \! = \! 1,\dotsc,N \! + \! 1$, $\mathcal{X}^{e}(z) \! = \! E_{a_{j}}^{
e}(z) \Psi^{e}((\tfrac{3}{4}n \xi_{a_{j}}^{e}(z))^{2/3}) \exp (\tfrac{n}{2}
(\xi_{a_{j}}^{e}(z) \! - \! \mi \mho_{j}^{e}) \sigma_{3})$, whence, using the
expression above for $\Psi^{e}(\zeta)$, $\zeta \! \in \! \mathbb{C}_{+} \cap
\widehat{\Omega}_{a_{j}}^{e,1}$, $j \! = \! 1,\dotsc,N \! + \! 1$, and the
asymptotic expansions for $\operatorname{Ai}(\pmb{\cdot})$ and $\operatorname{
Ai}^{\prime}(\pmb{\cdot})$ (as $n \! \to \! \infty)$ given in Equations~(2.6),
one arrives at
\begin{equation*}
\mathcal{X}^{e}(z) \underset{\underset{z \in \partial \Omega_{a_{j}}^{e,1}
\cap \partial \mathbb{U}_{\delta_{a_{j}}}^{e}}{n \to \infty}}{=} \dfrac{1}{
\sqrt{\smash[b]{2 \pi}}}E_{a_{j}}^{e}(z) \! \left(\! \left(\dfrac{3}{4}n
\xi_{a_{j}}^{e}(z) \right)^{2/3} \right)^{-\frac{1}{4} \sigma_{3}} \!
\begin{pmatrix}
\me^{-\frac{\mi \pi}{6}} & \me^{\frac{\mi \pi}{3}} \\
-\me^{-\frac{\mi \pi}{6}} & -\me^{\frac{4 \pi \mi}{3}}
\end{pmatrix} \! \me^{-\frac{\mi}{2}n \mho_{j}^{e} \sigma_{3}} \! \left(
\mathrm{I} \! + \! \mathcal{O}(n^{-1}) \right):
\end{equation*}
demanding that, for $z \! \in \! \partial \Omega_{a_{j}}^{e,1} \cap \partial
\mathbb{U}_{\delta_{a_{j}}}^{e}$, $j \! = \! 1,\dotsc,N \! + \! 1$,
$\overset{e}{m}^{\raise-1.0ex\hbox{$\scriptstyle \infty$}}(z)(\mathcal{X}^{e}
(z))^{-1} \! =_{n \to \infty} \! \mathrm{I} \! + \! \mathcal{O}(n^{-1})$, one
gets that
\begin{equation*}
E_{a_{j}}^{e}(z) \! = \! \dfrac{1}{\sqrt{\smash[b]{2 \mi}}}
\overset{e}{m}^{\raise-1.0ex\hbox{$\scriptstyle \infty$}}(z) \me^{\frac{\mi}{2}
n \mho_{j}^{e} \sigma_{3}} \!
\begin{pmatrix}
\mi & -\mi \\
1 & 1
\end{pmatrix} \! \left(\! \left(\dfrac{3}{4}n \xi_{a_{j}}^{e}(z) \right)^{2/3}
\right)^{\frac{1}{4} \sigma_{3}}, \quad j \! = \! 1,\dotsc,N \! + \! 1
\end{equation*}
(note that $\det (E_{a_{j}}^{e}(z)) \! = \! 1)$. One mimicks the above
paradigm for the remaining boundary skeletons $\partial \Omega_{a_{j}}^{e,l}
\cap \partial \mathbb{U}_{\delta_{a_{j}}}^{e}$, $l \! = \! 2,3,4$, $j \! = \!
1,\dotsc,N \! + \! 1$, and shows that the exact same formula for $E_{a_{j}}^{
e}(z)$ given above is obtained; thus, for $E_{a_{j}}^{e}(z)$, $j \! = \! 1,
\dotsc,N \! + \! 1$, as given above, one concludes that, uniformly for $z \!
\in \! \partial \mathbb{U}_{\delta_{a_{j}}}^{e}$, $j \! = \! 1,\dotsc,N \! +
\! 1$, $\overset{e}{m}^{\raise-1.0ex\hbox{$\scriptstyle \infty$}}(z)
(\mathcal{X}^{e}(z))^{-1} \! =_{\underset{z \in \partial \mathbb{U}_{\delta_{
a_{j}}}^{e}}{n \to \infty}} \! \mathrm{I} \! + \! \mathcal{O}(n^{-1})$. There
remains, however, the question of unimodularity, since
\begin{equation*}
\det \! \left(\mathcal{X}^{e}(z) \right) \! = \!
\begin{vmatrix}
\operatorname{Ai}(\Phi_{a_{j}}^{e}(z)) & \operatorname{Ai}(\omega^{2} \Phi_{
a_{j}}^{e}(z)) \\
\operatorname{Ai}^{\prime}(\Phi_{a_{j}}^{e}(z)) & \omega^{2} \operatorname{
Ai}^{\prime}(\omega^{2} \Phi_{a_{j}}^{e}(z))
\end{vmatrix} \qquad \text{or} \qquad
\begin{vmatrix}
\operatorname{Ai}(\Phi_{a_{j}}^{e}(z)) & -\omega^{2} \operatorname{Ai}(\omega
\Phi_{a_{j}}^{e}(z)) \\
\operatorname{Ai}^{\prime}(\Phi_{a_{j}}^{e}(z)) & -\operatorname{Ai}^{\prime}
(\omega \Phi_{a_{j}}^{e}(z))
\end{vmatrix}:
\end{equation*}
multiplying $\mathcal{X}^{e}(z)$ on the left by a constant, $\widetilde{c}$, 
say, using the Wronskian relations (see Chapter~10 of \cite{a93}) 
$\operatorname{W}(\operatorname{Ai}(\lambda),\operatorname{Ai}(\omega^{2}
\lambda)) \! = \! (2 \pi)^{-1} \exp (\mi \pi/6)$ and $\operatorname{W}
(\operatorname{Ai}(\lambda),\operatorname{Ai}(\omega \lambda)) \! = \! 
-(2 \pi)^{-1} \exp (-\mi \pi/6)$, and the linear dependence relation for 
Airy functions, $\operatorname{Ai}(\lambda) \! + \! \omega \operatorname{Ai}
(\omega \lambda) \! + \! \omega^{2} \operatorname{Ai}(\omega^{2} \lambda) \! 
= \! 0$, one shows that, upon imposing the condition $\det (\mathcal{X}^{e}
(z)) \! = \! 1$, $\widetilde{c} \! = \! (2 \pi)^{1/2} \exp (-\mi \pi/12)$. 
\hfill $\qed$

The above analyses lead to the following lemma.
\begin{ccccc}
Let $\overset{e}{\mathscr{M}}^{\raise-1.0ex\hbox{$\scriptstyle \sharp$}}
\colon \mathbb{C} \setminus \Sigma_{e}^{\sharp} \! \to \! \operatorname{SL}_{
2}(\mathbb{C})$ solve the {\rm RHP}
$(\overset{e}{\mathscr{M}}^{\raise-1.0ex\hbox{$\scriptstyle \sharp$}}(z),
\overset{e}{\upsilon}^{\raise-1.0ex\hbox{$\scriptstyle \sharp$}}(z),\Sigma_{
e}^{\sharp})$ formulated in Lemma~{\rm 4.2}. Define
\begin{equation*}
\mathscr{S}_{p}^{e}(z) \! := \!
\begin{cases}
\overset{e}{m}^{\raise-1.0ex\hbox{$\scriptstyle \infty$}}(z), &\text{$z \! \in
\! \mathbb{C} \setminus \cup_{j=1}^{N+1}(\mathbb{U}_{\delta_{b_{j-1}}}^{e}
\cup \mathbb{U}_{\delta_{a_{j}}}^{e})$,} \\
\mathcal{X}^{e}(z), &\text{$z \! \in \! \cup_{j=1}^{N+1}(\mathbb{U}_{\delta_{
b_{j-1}}}^{e} \cup \mathbb{U}_{\delta_{a_{j}}}^{e})$,}
\end{cases}
\end{equation*}
where $\overset{e}{m}^{\raise-1.0ex\hbox{$\scriptstyle \infty$}} \colon
\mathbb{C} \setminus J_{e}^{\infty} \! \to \! \operatorname{SL}_{2}(\mathbb{
C})$ is characterised completely in Lemma~{\rm 4.5}, and: {\rm (1)} for $z \!
\in \! \mathbb{U}_{\delta_{b_{j-1}}}^{e}$, $j \! = \! 1,\dotsc,N \! + \! 1$,
$\mathcal{X}^{e} \colon \mathbb{U}_{\delta_{b_{j-1}}}^{e} \setminus \Sigma_{
b_{j-1}}^{e} \! \to \! \operatorname{SL}_{2}(\mathbb{C})$ solve the {\rm RHPs}
$(\mathcal{X}^{e}(z),
\overset{e}{\upsilon}^{\raise-1.0ex\hbox{$\scriptstyle \sharp$}}(z),\Sigma_{
b_{j-1}}^{e})$, $j \! = \! 1,\dotsc,N \! + \! 1$, formulated in
Lemma~{\rm 4.6;} and {\rm (2)} for $z \! \in \! \mathbb{U}_{\delta_{a_{j}}}^{
e}$, $j \! = \! 1,\dotsc,N \! + \! 1$, $\mathcal{X}^{e} \colon \mathbb{U}_{
\delta_{a_{j}}}^{e} \setminus \Sigma_{a_{j}}^{e} \! \to \! \operatorname{SL}_{
2}(\mathbb{C})$ solve the {\rm RHPs} $(\mathcal{X}^{e}(z),
\overset{e}{\upsilon}^{\raise-1.0ex\hbox{$\scriptstyle \sharp$}}(z),\Sigma_{
a_{j}}^{e})$, $j \! = \! 1,\dotsc,N \! + \! 1$, formulated in Lemma~{\rm 4.7}.
Set
\begin{equation*}
\mathscr{R}^{e}(z) \! := \!
\overset{e}{\mathscr{M}}^{\raise-1.0ex\hbox{$\scriptstyle \sharp$}}(z) \!
\left(\mathscr{S}_{p}^{e}(z) \right)^{-1},
\end{equation*}
and define the augmented contour $\Sigma_{p}^{e} \! := \! \Sigma_{e}^{\sharp}
\cup (\cup_{j=1}^{N+1}(\partial \mathbb{U}_{\delta_{b_{j-1}}}^{e} \cup \mathbb{
U}_{\delta_{a_{j}}}^{e}))$, with the orientation given in Figure~{\rm 9}.
Then $\mathscr{R}^{e} \colon \mathbb{C} \setminus \Sigma_{p}^{e} \! \to \!
\operatorname{SL}_{2}(\mathbb{C})$ solves the following {\rm RHP:} {\rm (i)}
$\mathscr{R}^{e}(z)$ is holomorphic for $z \! \in \! \mathbb{C} \setminus
\Sigma_{p}^{e};$ {\rm (ii)} $\mathscr{R}^{e}_{\pm}(z) \! := \! \lim_{
\underset{z^{\prime} \! \in \, \pm \, \mathrm{side} \, \mathrm{of} \,
\Sigma_{p}^{e}}{z^{\prime} \to z}} \mathscr{R}^{e}(z^{\prime})$ satisfy the
boundary condition
\begin{equation*}
\mathscr{R}^{e}_{+}(z) \! = \! \mathscr{R}^{e}_{-}(z) \upsilon_{\mathscr{R}}^{
e}(z), \quad z \! \in \! \Sigma_{p}^{e},
\end{equation*}
where
\begin{equation*}
\upsilon^{e}_{\mathscr{R}}(z) \! := \!
\begin{cases}
\upsilon_{\mathscr{R}}^{e,1}(z), &\text{$z \! \in \! (-\infty,b_{0}^{e} \! -
\! \delta_{b_{0}}^{e}) \! \cup \! (a_{N+1}^{e} \! + \! \delta_{a_{N+1}}^{e},
+\infty) \! =: \! \Sigma_{p}^{e,1}$,} \\
\upsilon_{\mathscr{R}}^{e,2}(z), &\text{$z \! \in \! (a_{j}^{e} \! + \!
\delta_{a_{j}}^{e},b_{j}^{e} \! - \! \delta_{b_{j}}^{e}) \! =: \! \Sigma_{p,
j}^{e,2} \! \subset \! \cup_{l=1}^{N} \Sigma_{p,l}^{e,2} \! =: \! \Sigma_{p}^{
e,2}$,} \\
\upsilon_{\mathscr{R}}^{e,3}(z), &\text{$z \! \in \! \cup_{j=1}^{N+1}(J_{j}^{
e,\smallfrown} \setminus (\mathbb{C}_{+} \cap (\mathbb{U}_{\delta_{b_{j-1}}}^{
e} \cup \mathbb{U}_{\delta_{a_{j}}}^{e}))) \! =: \! \Sigma_{p}^{e,3}$,} \\
\upsilon_{\mathscr{R}}^{e,4}(z), &\text{$z \! \in \! \cup_{j=1}^{N+1}(J_{j}^{
e,\smallsmile} \setminus (\mathbb{C}_{-} \cap (\mathbb{U}_{\delta_{b_{j-1}}}^{
e} \cup \mathbb{U}_{\delta_{a_{j}}}^{e}))) \! =: \! \Sigma_{p}^{e,4}$,} \\
\upsilon_{\mathscr{R}}^{e,5}(z), &\text{$z \! \in \! \cup_{j=1}^{N+1}(\partial
\mathbb{U}_{\delta_{b_{j-1}}}^{e} \cup \mathbb{U}_{\delta_{a_{j}}}^{e}) \! =:
\! \Sigma_{p}^{e,5}$,} \\
\mathrm{I}, &\text{$z \! \in \! \Sigma_{p}^{e} \setminus \cup_{l=1}^{5}
\Sigma_{p}^{e,l}$,}
\end{cases}
\end{equation*}
with
\begin{align*}
\upsilon_{\mathscr{R}}^{e,1}(z) \! =& \, \mathrm{I} \! + \! \me^{n(g^{e}_{+}(z)
+g^{e}_{-}(z)-\widetilde{V}(z)-\ell_{e}+2Q_{e})} \, 
\overset{e}{m}^{\raise-1.0ex\hbox{$\scriptstyle \infty$}}(z) \sigma_{+}
(\overset{e}{m}^{\raise-1.0ex\hbox{$\scriptstyle \infty$}}(z))^{-1}, \\
\upsilon_{\mathscr{R}}^{e,2}(z) \! =& \, \mathrm{I} \! + \! \me^{-\mi n
\Omega_{j}^{e}+n(g^{e}_{+}(z)+g^{e}_{-}(z)-\widetilde{V}(z)-\ell_{e}+2Q_{e})} 
\, \overset{e}{m}^{\raise-1.0ex\hbox{$\scriptstyle \infty$}}_{-}(z) \sigma_{+}
(\overset{e}{m}^{\raise-1.0ex\hbox{$\scriptstyle \infty$}}_{-}(z))^{-1}, \\
\upsilon_{\mathscr{R}}^{e,3}(z) \! =& \, \mathrm{I} \! + \! \me^{-4 n \pi 
\mi \int_{z}^{a_{N+1}^{e}} \psi_{V}^{e}(s) \, \md s} \, 
\overset{e}{m}^{\raise-1.0ex\hbox{$\scriptstyle \infty$}}(z) \sigma_{-}
(\overset{e}{m}^{\raise-1.0ex\hbox{$\scriptstyle \infty$}}(z))^{-1}, \\
\upsilon_{\mathscr{R}}^{e,4}(z) \! =& \, \mathrm{I} \! + \! \me^{4 n \pi 
\mi \int_{z}^{a_{N+1}^{e}} \psi_{V}^{e}(s) \, \md s} \, 
\overset{e}{m}^{\raise-1.0ex\hbox{$\scriptstyle \infty$}}(z) \sigma_{-}
(\overset{e}{m}^{\raise-1.0ex\hbox{$\scriptstyle \infty$}}(z))^{-1}, \\
\upsilon_{\mathscr{R}}^{e,5}(z) \! =& \, \mathcal{X}^{e}(z)
(\overset{e}{m}^{\raise-1.0ex\hbox{$\scriptstyle \infty$}}(z))^{-1}:
\end{align*}
{\rm (iii)} $\mathscr{R}^{e}(z) \! =_{\underset{z \in \mathbb{C} \setminus
\Sigma_{p}^{e}}{z \to \infty}} \! \mathrm{I} \! + \! \mathcal{O}(z^{-1});$ and
{\rm (iv)} $\mathscr{R}^{e}(z) \! =_{\underset{z \in \mathbb{C} \setminus
\Sigma_{p}^{e}}{z \to 0}} \! \mathcal{O}(1)$.
\end{ccccc}
\begin{figure}[tbh]
\begin{center}
\begin{pspicture}(0,0)(15,5)
\psset{xunit=1cm,yunit=1cm}
\psline[linewidth=0.6pt,linestyle=solid,linecolor=black](0,2.5)(1,2.5)
\psline[linewidth=0.6pt,linestyle=solid,linecolor=magenta,arrowsize=1.5pt 5]%
{->}(1,2.5)(2.5,2.5)
\psline[linewidth=0.6pt,linestyle=solid,linecolor=magenta](2.4,2.5)(4,2.5)
\psline[linewidth=0.6pt,linestyle=solid,linecolor=black](4,2.5)(4.7,2.5)
\psarcn[linewidth=0.6pt,linestyle=solid,linecolor=magenta,arrowsize=1.5pt 5]%
{->}(2.5,1.5){1.8}{146}{90}
\psarcn[linewidth=0.6pt,linestyle=solid,linecolor=magenta](2.5,1.5){1.8}{90}%
{34}
\psarc[linewidth=0.6pt,linestyle=solid,linecolor=magenta,arrowsize=1.5pt 5]%
{->}(2.5,3.5){1.8}{214}{270}
\psarc[linewidth=0.6pt,linestyle=solid,linecolor=magenta](2.5,3.5){1.8}{270}%
{326}
\psarcn[linewidth=0.6pt,linestyle=solid,linecolor=cyan,arrowsize=1.5pt 5]%
{->}(4,2.5){0.6}{180}{45}
\psarcn[linewidth=0.6pt,linestyle=solid,linecolor=cyan](4,2.5){0.6}{45}{0}
\psarcn[linewidth=0.6pt,linestyle=solid,linecolor=cyan](4,2.5){0.6}{360}{180}
\psarcn[linewidth=0.6pt,linestyle=solid,linecolor=cyan,arrowsize=1.5pt 5]%
{->}(1,2.5){0.6}{180}{135}
\psarcn[linewidth=0.6pt,linestyle=solid,linecolor=cyan](1,2.5){0.6}{135}{0}
\psarcn[linewidth=0.6pt,linestyle=solid,linecolor=cyan](1,2.5){0.6}{360}{180}
\rput(2.5,2.9){\makebox(0,0){$\Omega^{e,\smallfrown}_{1}$}}
\rput(2.5,2.1){\makebox(0,0){$\Omega^{e,\smallsmile}_{1}$}}
\rput(3,3.6){\makebox(0,0){$J^{e,\smallfrown}_{1}$}}
\rput(3,1.4){\makebox(0,0){$J^{e,\smallsmile}_{1}$}}
\rput(1,3.4){\makebox(0,0){$\partial \mathbb{U}^{e}_{\delta_{b_{0}}}$}}
\rput(4.1,3.4){\makebox(0,0){$\partial \mathbb{U}^{e}_{\delta_{a_{1}}}$}}
\psline[linewidth=0.6pt,linestyle=solid,linecolor=black](5.3,2.5)(6,2.5)
\psline[linewidth=0.6pt,linestyle=solid,linecolor=magenta,arrowsize=1.5pt 5]%
{->}(6,2.5)(7.5,2.5)
\psline[linewidth=0.6pt,linestyle=solid,linecolor=magenta](7.4,2.5)(9,2.5)
\psline[linewidth=0.6pt,linestyle=solid,linecolor=black](9,2.5)(9.7,2.5)
\psarcn[linewidth=0.6pt,linestyle=solid,linecolor=magenta,arrowsize=1.5pt 5]%
{->}(7.5,1.5){1.8}{146}{90}
\psarcn[linewidth=0.6pt,linestyle=solid,linecolor=magenta](7.5,1.5){1.8}{90}%
{34}
\psarc[linewidth=0.6pt,linestyle=solid,linecolor=magenta,arrowsize=1.5pt 5]%
{->}(7.5,3.5){1.8}{214}{270}
\psarc[linewidth=0.6pt,linestyle=solid,linecolor=magenta](7.5,3.5){1.8}{270}%
{326}
\psarcn[linewidth=0.6pt,linestyle=solid,linecolor=cyan,arrowsize=1.5pt 5]%
{->}(6,2.5){0.6}{180}{135}
\psarcn[linewidth=0.6pt,linestyle=solid,linecolor=cyan](6,2.5){0.6}{135}{0}
\psarcn[linewidth=0.6pt,linestyle=solid,linecolor=cyan](6,2.5){0.6}{360}{180}
\psarcn[linewidth=0.6pt,linestyle=solid,linecolor=cyan,arrowsize=1.5pt 5]%
{->}(9,2.5){0.6}{180}{45}
\psarcn[linewidth=0.6pt,linestyle=solid,linecolor=cyan](9,2.5){0.6}{45}{0}
\psarcn[linewidth=0.6pt,linestyle=solid,linecolor=cyan](9,2.5){0.6}{360}{180}
\rput(7.5,2.9){\makebox(0,0){$\Omega^{e,\smallfrown}_{j}$}}
\rput(7.5,2.1){\makebox(0,0){$\Omega^{e,\smallsmile}_{j}$}}
\rput(8,3.6){\makebox(0,0){$J^{e,\smallfrown}_{j}$}}
\rput(8,1.4){\makebox(0,0){$J^{e,\smallsmile}_{j}$}}
\rput(6,3.4){\makebox(0,0){$\partial \mathbb{U}^{e}_{\delta_{b_{j-1}}}$}}
\rput(9.1,3.4){\makebox(0,0){$\partial \mathbb{U}^{e}_{\delta_{a_{j}}}$}}
\psline[linewidth=0.6pt,linestyle=solid,linecolor=black](10.3,2.5)(11,2.5)
\psline[linewidth=0.6pt,linestyle=solid,linecolor=magenta,arrowsize=1.5pt 5]%
{->}(11,2.5)(12.5,2.5)
\psline[linewidth=0.6pt,linestyle=solid,linecolor=magenta](12.4,2.5)(14,2.5)
\psline[linewidth=0.6pt,linestyle=solid,linecolor=black](14,2.5)(15,2.5)
\psarcn[linewidth=0.6pt,linestyle=solid,linecolor=magenta,arrowsize=1.5pt 5]%
{->}(12.5,1.5){1.8}{146}{90}
\psarcn[linewidth=0.6pt,linestyle=solid,linecolor=magenta](12.5,1.5){1.8}{90}%
{34}
\psarc[linewidth=0.6pt,linestyle=solid,linecolor=magenta,arrowsize=1.5pt 5]%
{->}(12.5,3.5){1.8}{214}{270}
\psarc[linewidth=0.6pt,linestyle=solid,linecolor=magenta](12.5,3.5){1.8}{270}%
{326}
\psarcn[linewidth=0.6pt,linestyle=solid,linecolor=cyan,arrowsize=1.5pt 5]%
{->}(11,2.5){0.6}{180}{135}
\psarcn[linewidth=0.6pt,linestyle=solid,linecolor=cyan](11,2.5){0.6}{135}{0}
\psarcn[linewidth=0.6pt,linestyle=solid,linecolor=cyan](11,2.5){0.6}{360}{180}
\psarcn[linewidth=0.6pt,linestyle=solid,linecolor=cyan,arrowsize=1.5pt 5]%
{->}(14,2.5){0.6}{180}{45}
\psarcn[linewidth=0.6pt,linestyle=solid,linecolor=cyan](14,2.5){0.6}{45}{0}
\psarcn[linewidth=0.6pt,linestyle=solid,linecolor=cyan](14,2.5){0.6}{360}{180}
\rput(12.4,2.95){\makebox(0,0){$\Omega^{e,\smallfrown}_{N+1}$}}
\rput(12.4,2.05){\makebox(0,0){$\Omega^{e,\smallsmile}_{N+1}$}}
\rput(13,3.6){\makebox(0,0){$J^{e,\smallfrown}_{N+1}$}}
\rput(13,1.4){\makebox(0,0){$J^{e,\smallsmile}_{N+1}$}}
\rput(11,3.4){\makebox(0,0){$\partial \mathbb{U}^{e}_{\delta_{b_{N}}}$}}
\rput(14.35,3.4){\makebox(0,0){$\partial \mathbb{U}^{e}_{\delta_{a_{N+1}}}$}}
\psline[linewidth=0.9pt,linestyle=dotted,linecolor=darkgray](4.8,2.5)(5.2,2.5)
\psline[linewidth=0.9pt,linestyle=dotted,linecolor=darkgray](9.8,2.5)(10.2,2.5)
\psdots[dotstyle=*,dotscale=1.5](1,2.5)
\psdots[dotstyle=*,dotscale=1.5](4,2.5)
\psdots[dotstyle=*,dotscale=1.5](6,2.5)
\psdots[dotstyle=*,dotscale=1.5](9,2.5)
\psdots[dotstyle=*,dotscale=1.5](11,2.5)
\psdots[dotstyle=*,dotscale=1.5](14,2.5)
\rput(1,2.2){\makebox(0,0){$\scriptstyle \pmb{b_{0}^{e}}$}}
\rput(4,2.2){\makebox(0,0){$\scriptstyle \pmb{a_{1}^{e}}$}}
\rput(6,2.2){\makebox(0,0){$\scriptstyle \pmb{b_{j-1}^{e}}$}}
\rput(9,2.2){\makebox(0,0){$\scriptstyle \pmb{a_{j}^{e}}$}}
\rput(11,2.2){\makebox(0,0){$\scriptstyle \pmb{b_{N}^{e}}$}}
\rput(14,2.2){\makebox(0,0){$\scriptstyle \pmb{a_{N+1}^{e}}$}}
\rput(3.15,2.5){\makebox(0,0){$\scriptstyle \pmb{J_{1}^{e}}$}}
\rput(8.15,2.5){\makebox(0,0){$\scriptstyle \pmb{J_{j}^{e}}$}}
\rput(13.05,2.5){\makebox(0,0){$\scriptstyle \pmb{J_{N+1}^{e}}$}}
\end{pspicture}
\end{center}
\vspace{-1.00cm}
\caption{The augmented contour $\Sigma_{p}^{e} \! := \! \Sigma_{e}^{\sharp}
\cup (\cup_{j=1}^{N+1}(\partial \mathbb{U}_{\delta_{b_{j-1}}}^{e} \cup
\partial \mathbb{U}_{\delta_{a_{j}}}^{e}))$}
\end{figure}
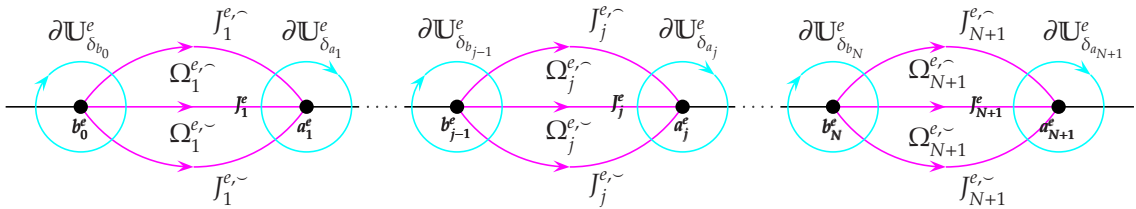

\emph{Proof.} Define the oriented, augmented skeleton $\Sigma_{p}^{e}$ as
in the Lemma: the RHP $(\mathscr{R}^{e}(z),\upsilon^{e}_{\mathscr{R}}(z),
\Sigma_{p}^{e})$ follows {}from the RHPs
$(\overset{e}{\mathscr{M}}^{\raise-1.0ex\hbox{$\scriptstyle \sharp$}}(z),
\overset{e}{\upsilon}^{\raise-1.0ex\hbox{$\scriptstyle \sharp$}}(z),\Sigma_{
e}^{\sharp})$ and $(\overset{e}{m}^{\raise-1.0ex\hbox{$\scriptstyle \infty$}}
(z),\overset{e}{\upsilon}^{\raise-1.0ex\hbox{$\scriptstyle \infty$}}(z),J_{e}^{
\infty})$ formulated in Lemmas~4.2 and~4.3, respectively, upon using the
definitions of $\mathscr{S}_{p}^{e}(z)$ and $\mathscr{R}^{e}(z)$ given in the
Lemma. \hfill $\qed$
\section{Asymptotic (as $n \! \to \! \infty)$ Solution of the RHP for 
$\stackrel{e}{\mathrm{Y}} \! (z)$}
In this section, via the Beals-Coifman (BC) construction \cite{a84}, the 
(normalised at infinity) RHP $(\mathscr{R}^{e}(z),\linebreak[4] 
\upsilon^{e}_{\mathscr{R}}(z),\Sigma_{p}^{e})$ formulated in Lemma~4.8 is 
solved asymptotically (as $n \! \to \! \infty)$; in particular, it is shown 
that, uniformly for $z \! \in \! \Sigma_{p}^{e}$,
\begin{equation*}
\norm{\upsilon^{e}_{\mathscr{R}}(\cdot) \! - \! \mathrm{I}}_{\cap_{p \in \{
1,2,\infty\}} \mathcal{L}^{p}_{\mathrm{M}_{2}(\mathbb{C})}(\Sigma_{p}^{e})}
\underset{n \to \infty}{=} \mathrm{I} \! + \! \mathcal{O} \! \left(f(n)n^{-1}
\right),
\end{equation*}
where $f(n) \! =_{n \to \infty} \! \mathcal{O}(1)$, and, subsequently, the
original \textbf{RHP1}, that is, $(\overset{e}{\mathrm{Y}}(z),\mathrm{I}
\! + \! \me^{-n \widetilde{V}(z)} \sigma_{+},\mathbb{R})$, is solved
asymptotically by re-tracing the finite sequence of RHP transformations
$\mathscr{R}^{e}(z)$ (Lemmas~5.3 and~4.8) $\to \!
\overset{e}{\mathscr{M}}^{\raise-1.0ex\hbox{$\scriptstyle \sharp$}}(z)$
(Lemma~4.2) $\to \! \overset{e}{\mathscr{M}}(z)$ (Lemma~3.4) $\to \!
\overset{e}{\mathrm{Y}}(z)$. The (unique) solution for $\overset{e}{\mathrm{
Y}}(z)$ then leads to the final asymptotic results for $\boldsymbol{\pi}_{2
n}(z)$ (in the entire complex plane), $\xi_{n}^{(2n)}$ and $\phi_{2n}(z)$ (in
the entire complex plane) stated, respectively, in Theorems~2.3.1 and~2.3.2.
\begin{bbbbb}
Let $\mathscr{R}^{e} \colon \mathbb{C} \setminus \Sigma_{p}^{e} \! \to \!
\operatorname{SL}_{2}(\mathbb{C})$ solve the {\rm RHP} $(\mathscr{R}^{e}(z),
\upsilon_{\mathscr{R}}^{e}(z),\Sigma_{p}^{e})$ formulated in Lemma {\rm 4.8}.
Then:
\begin{compactenum}
\item[{\rm (1)}] for $z \! \in \! (-\infty,b_{0}^{e} \! - \! \delta_{b_{0}}^{
e}) \cup (a_{N+1}^{e} \! + \! \delta_{a_{N+1}}^{e},+\infty) \! =: \! \Sigma_{
p}^{e,1}$,
\begin{equation*}
\upsilon_{\mathscr{R}}^{e}(z) \underset{n \to \infty}{=}
\begin{cases}
\mathrm{I} \! + \! \mathcal{O}(f_{\infty}(n) \me^{-nc_{\infty} \vert z
\vert}), &\text{$z \! \in \! \Sigma_{p}^{e,1} \setminus \mathbb{U}_{0}^{e}$,}
\\
\mathrm{I} \! + \! \mathcal{O}(f_{0}(n) \me^{-nc_{0} \vert z \vert^{-1}}),
&\text{$z \! \in \! \Sigma_{p}^{e} \cap \mathbb{U}_{0}^{e}$,}
\end{cases}
\end{equation*}
where $c_{0},c_{\infty} \! > \! 0$, $(f_{\infty}(n))_{ij} \! =_{n \to \infty}
\! \mathcal{O}(1)$, $(f_{0}(n))_{ij} \! =_{n \to \infty} \! \mathcal{O}(1)$,
$i,j \! = \! 1,2$, and $\mathbb{U}_{0}^{e} \! := \! \lbrace \mathstrut z \!
\in \! \mathbb{C}; \, \vert z \vert \! < \! \epsilon \rbrace$, with $\epsilon$
some arbitrarily fixed, sufficiently small positive real number;
\item[{\rm (2)}] for $z \! \in \! (a_{j}^{e} \! + \! \delta_{a_{j}}^{e},b_{
j}^{e} \! - \! \delta_{b_{j}}^{e}) \! =: \! \Sigma_{p,j}^{e,2} \subset \cup_{l
=1}^{N} \Sigma_{p,l}^{e,2} \! =: \! \Sigma_{p}^{e,2}$, $j \! = \! 1,\dotsc,N$,
\begin{equation*}
\upsilon_{\mathscr{R}}^{e}(z) \! \underset{n \to \infty}{=} \!
\begin{cases}
\mathrm{I} \! + \! \mathcal{O}(f_{j}(n) \me^{-nc_{j}(z \! - \! a_{j}^{e})}),
&\text{$z \! \in \! \Sigma^{e,2}_{p,j} \setminus \mathbb{U}^{e}_{0}$,} \\
\mathrm{I} \! + \! \mathcal{O}(\widetilde{f}_{j}(n) \me^{-n \widetilde{c}_{j}
\vert z \vert^{-1}}), &\text{$z \! \in \! \Sigma^{e,2}_{p,j} \cap \mathbb{
U}^{e}_{0}$,}
\end{cases}
\end{equation*}
where $c_{j},\widetilde{c}_{j} \! > \! 0$, $(f_{j}(n))_{kl} \! =_{n \to
\infty} \! \mathcal{O}(1)$, and $(\widetilde{f}_{j}(n))_{kl} \! =_{n \to
\infty} \! \mathcal{O}(1)$, $k,l \! = \! 1,2;$
\item[{\rm (3)}] for $z \! \in \! \cup_{j=1}^{N+1}(J_{j}^{e,\smallfrown}
\setminus (\mathbb{C}_{+} \cap (\mathbb{U}_{\delta_{b_{j-1}}}^{e} \cup
\mathbb{U}_{\delta_{a_{j}}}^{e}))) \! =: \! \Sigma_{p}^{e,3}$,
\begin{equation*}
\upsilon_{\mathscr{R}}^{e}(z) \underset{n \to \infty}{=} \mathrm{I} \! + \!
\mathcal{O}(\overset{\smallfrown}{f}(n) \me^{-n \overset{\smallfrown}{c} \vert
z \vert}),
\end{equation*}
where $\overset{\smallfrown}{c} \! > \! 0$ and $(\overset{\smallfrown}{f}(n))_{
ij} \! =_{n \to \infty} \! \mathcal{O}(1)$, $i,j \! = \! 1,2;$
\item[{\rm (4)}] for $z \! \in \! \cup_{j=1}^{N+1}(J_{j}^{e,\smallsmile}
\setminus (\mathbb{C}_{-} \cap (\mathbb{U}_{\delta_{b_{j-1}}}^{e} \cup
\mathbb{U}_{\delta_{a_{j}}}^{e}))) \! =: \! \Sigma_{p}^{e,4}$,
\begin{equation*}
\upsilon_{\mathscr{R}}^{e}(z) \underset{n \to \infty}{=} \mathrm{I} \! + \!
\mathcal{O}(\overset{\smallsmile}{f}(n) \me^{-n \overset{\smallsmile}{c} \vert
z \vert}),
\end{equation*}
where $\overset{\smallsmile}{c} \! > \! 0$ and $(\overset{\smallsmile}{f}(n))_{
ij} \! =_{n \to \infty} \! \mathcal{O}(1)$, $i,j \! = \! 1,2;$ and
\item[{\rm (5)}] for $z \! \in \! \cup_{j=1}^{N+1}(\partial \mathbb{U}_{
\delta_{b_{j-1}}}^{e} \cup \partial \mathbb{U}_{\delta_{a_{j}}}^{e}) \! =: \!
\Sigma_{p}^{e,5}$,
\begin{align*}
\upsilon_{\mathscr{R}}^{e}(z) \underset{\underset{z \in \mathbb{C}_{\pm} \cap
\partial \mathbb{U}_{\delta_{b_{j-1}}}^{e}}{n \to \infty}}{=}& \, \mathrm{I}
\! + \! \dfrac{1}{n \xi_{b_{j-1}}^{e}(z)} \,
\overset{e}{\mathfrak{M}}^{\raise-1.0ex\hbox{$\scriptstyle \infty$}}(z) \!
\begin{pmatrix}
\mp (s_{1}+t_{1}) & \mp \mi (s_{1}-t_{1}) \me^{\mi n \mho_{j-1}^{e}} \\
\mp \mi (s_{1}-t_{1}) \me^{-\mi n \mho_{j-1}^{e}} & \pm (s_{1}+t_{1})
\end{pmatrix} \!
(\overset{e}{\mathfrak{M}}^{\raise-1.0ex\hbox{$\scriptstyle \infty$}}(z))^{-1}
\\
+& \, \mathcal{O} \! \left(\dfrac{1}{(n \xi_{b_{j-1}}^{e}(z))^{2}} \,
\overset{e}{\mathfrak{M}}^{\raise-1.0ex\hbox{$\scriptstyle \infty$}}(z)f_{
b_{j-1}}^{e}(n)
(\overset{e}{\mathfrak{M}}^{\raise-1.0ex\hbox{$\scriptstyle \infty$}}(z))^{-1}
\right), \qquad j \! = \! 1,\dotsc,N \! + \! 1,
\end{align*}
where $\overset{e}{\mathfrak{M}}^{\raise-1.0ex\hbox{$\scriptstyle \infty$}}
(z)$ is characterised completely in Lemma {\rm 4.5}, $s_{1} \! = \! 5/72$,
$t_{1} \! = \! -7/72$, for $j \! = \! 1,\dotsc,N \! + \! 1$, $\xi_{b_{j-1}}^{
e}(z) \! = \! -2 \int_{z}^{b_{j-1}^{e}}(R_{e}(s))^{1/2}h_{V}^{e}(s) \, \md s
\! = \! (z \! - \! b_{j-1}^{e})^{3/2}G_{b_{j-1}}^{e}(z)$, with $G_{b_{j-1}}^{
e}(z)$ described completely in Lemma {\rm 4.6}, $\mho_{j-1}^{e}$ is defined in
Remark~{\rm 4.4}, and $(f_{b_{j-1}}^{e}(n))_{kl} \! =_{n \to \infty} \!
\mathcal{O}(1)$, $k,l \! = \! 1,2$, and
\begin{align*}
\upsilon_{\mathscr{R}}^{e}(z) \underset{\underset{z \in \mathbb{C}_{\pm} \cap
\partial \mathbb{U}_{\delta_{a_{j}}}^{e}}{n \to \infty}}{=}& \, \mathrm{I} \!
+ \! \dfrac{1}{n \xi_{a_{j}}^{e}(z)} \,
\overset{e}{\mathfrak{M}}^{\raise-1.0ex\hbox{$\scriptstyle \infty$}}(z) \!
\begin{pmatrix}
\mp (s_{1}+t_{1}) & \pm \mi (s_{1}-t_{1}) \me^{\mi n \mho_{j}^{e}} \\
\pm \mi (s_{1}-t_{1}) \me^{-\mi n \mho_{j}^{e}} & \pm (s_{1}+t_{1})
\end{pmatrix} \!
(\overset{e}{\mathfrak{M}}^{\raise-1.0ex\hbox{$\scriptstyle \infty$}}(z))^{-1}
\\
+& \, \mathcal{O} \! \left(\dfrac{1}{(n \xi_{a_{j}}^{e}(z))^{2}} \,
\overset{e}{\mathfrak{M}}^{\raise-1.0ex\hbox{$\scriptstyle \infty$}}(z)f_{a_{
j}}^{e}(n)
(\overset{e}{\mathfrak{M}}^{\raise-1.0ex\hbox{$\scriptstyle \infty$}}(z))^{-1}
\right), \qquad j \! = \! 1,\dotsc,N \! + \! 1,
\end{align*}
where, for $j \! = \! 1,\dotsc,N \! + \! 1$, $\xi_{a_{j}}^{e}(z) \! = \! 2
\int_{a_{j}^{e}}^{z}(R_{e}(s))^{1/2}h_{V}^{e}(s) \, \md s \! = \! (z \! -
\! a_{j}^{e})^{3/2}G_{a_{j}}^{e}(z)$, with $G_{a_{j}}^{e}(z)$ described
completely in Lemma {\rm 4.7}, and $(f_{a_{j}}^{e}(n))_{kl} \! =_{n \to
\infty} \! \mathcal{O}(1)$, $k,l \! = \! 1,2$.
\end{compactenum}
\end{bbbbb}

\emph{Proof.} Recall the definition of $\upsilon_{\mathscr{R}}^{e}(z)$ given
in Lemma~4.8. For $z \! \in \! \Sigma_{p}^{e,1} \! := \! (-\infty,b_{0}^{e} \!
- \! \delta_{b_{0}}^{e}) \cup (a_{N+1}^{e} \! + \! \delta_{a_{N+1}}^{e},
+\infty)$, recall {}from Lemma~4.8 that
\begin{equation*}
\upsilon_{\mathscr{R}}^{e}(z) \! := \! \upsilon_{\mathscr{R}}^{e,1}(z) \! = \!
\mathrm{I} \! + \! \exp \! \left(n(g^{e}_{+}(z) \! + \! g^{e}_{-}(z) \! - \!
\widetilde{V}(z) \! - \! \ell_{e} \! + \! 2Q_{e}) \right) \!
\overset{e}{m}^{\raise-1.0ex\hbox{$\scriptstyle \infty$}}(z) \sigma_{+}
(\overset{e}{m}^{\raise-1.0ex\hbox{$\scriptstyle \infty$}}(z))^{-1},
\end{equation*}
and, {}from the proof of Lemma~4.1, $g^{e}_{+}(z) \! + \! g^{e}_{-}(z) \! - \!
\widetilde{V}(z) \! - \! \ell_{e} \! + \! 2Q_{e}$ equals $-2 \int_{a_{N+1}^{e}
}^{z}(R_{e}(s))^{1/2}h_{V}^{e}(s) \, \md s$ $(< \! 0)$ for $z \! \in \! (a_{N+
1}^{e} \! + \! \delta_{a_{N+1}}^{e},+\infty)$ and equals $2 \int_{z}^{b_{0}^{
e}}(R_{e}(s))^{1/2}h_{V}^{e}(s) \, \md s$ $(< \! 0)$ for $z \! \in \!
(-\infty,b_{0}^{e} \! - \! \delta_{b_{0}}^{e})$; hence, recalling that
$\widetilde{V} \colon \mathbb{R} \setminus \{0\} \! \to \! \mathbb{R}$, which
is regular, satisfies conditions~(2.3)--(2.5), using the asymptotic expansions
(as $\vert z \vert \! \to \! \infty$ and $\vert z \vert \! \to \! 0)$ for $g^{
e}_{+}(z) \! + \! g^{e}_{-}(z) \! - \! \widetilde{V}(z) \! - \! \ell_{e} \! +
\! 2Q_{e}$ given in the proof of Lemma~3.6, that is, $g^{e}_{+}(z) \! + \! g^{
e}_{-}(z) \! - \! \widetilde{V}(z) \! - \! \ell_{e} \! + \! 2Q_{e} \! =_{\vert
z \vert \to \infty} \! \ln (z^{2} \! + \! 1) \! - \! \widetilde{V}(z) \! + \!
\mathcal{O}(1)$ and $g^{e}_{+}(z) \! + \! g^{e}_{-}(z) \! - \! \widetilde{V}
(z) \! - \! \ell_{e} \! + \! 2Q_{e} \! =_{\vert z \vert \to 0} \! \ln (z^{-2}
\! + \! 1) \! - \! \widetilde{V}(z) \! + \! \mathcal{O}(1)$, upon recalling
the expression for $\overset{e}{m}^{\raise-1.0ex\hbox{$\scriptstyle \infty$}}
(z)$ given in Lemma~4.5 and noting that the respective factors $\gamma^{e}(z)
\! \pm \! (\gamma^{e}(z))^{-1}$ and $\boldsymbol{\theta}^{e}(\pm \boldsymbol{
u}^{e}(z) \! - \! \tfrac{n}{2 \pi} \boldsymbol{\Omega}^{e} \! \pm \!
\boldsymbol{d}_{e})$ are uniformly bounded (with respect to $z)$ in compact
subsets outside the open intervals surrounding the end-points of the suppport
of the `even' equilibrium measure, defining $\mathbb{U}_{0}^{e}$ as in the
Proposition, one arrives at the asymptotic (as $n \! \to \! \infty)$ estimates
for $\upsilon_{\mathscr{R}}^{e}(z)$ on $\Sigma_{p}^{e,1} \setminus \mathbb{
U}_{0}^{e} \! \ni \! z$ and $\Sigma_{p}^{e,1} \cap \mathbb{U}_{0}^{e} \! \ni
\! z$ stated in item~(1) of the Proposition. (It should be noted that the
$n$-dependence of the $\operatorname{GL}_{2}(\mathbb{C})$-valued factors
$f_{\infty}(n)$ and $f_{0}(n)$ are inherited {}from the bounded $(\mathcal{O}
(1))$ $n$-dependence of the respective Riemann theta functions, whose
corresponding series converge absolutely and uniformly due to the fact that
the associated Riemann matrix of $\boldsymbol{\beta}^{e}$-periods, $\tau^{e}$,
is pure imaginary and $-\mi \tau^{e}$ is positive definite.)

For $z \! \in \! \Sigma_{p,j}^{e,2} \! := \! (a_{j}^{e} \! + \! \delta_{a_{j}
}^{e},b_{j}^{e} \! - \! \delta_{b_{j}}^{e})$, $j \! = \! 1,\dotsc,N$, recall
{}from Lemma~4.8 that
\begin{equation*}
\upsilon_{\mathscr{R}}^{e}(z) \! := \! \upsilon_{\mathscr{R}}^{e,2}(z) \! = \!
\mathrm{I} \! + \! \me^{-\mi n \Omega_{j}^{e}} \exp \! \left(n(g^{e}_{+}(z) \!
+ \! g^{e}_{-}(z) \! - \! \widetilde{V}(z) \! - \! \ell_{e} \! + \! 2Q_{e})
\right) \! \overset{e}{m}^{\raise-1.0ex\hbox{$\scriptstyle \infty$}}_{-}(z)
\sigma_{+}(\overset{e}{m}^{\raise-1.0ex\hbox{$\scriptstyle \infty$}}_{-}(z)
)^{-1},
\end{equation*}
and, {}from the proof of Lemma~4.1, $g^{e}_{+}(z) \! + \! g^{e}_{-}(z) \! - \!
\widetilde{V}(z) \! - \! \ell_{e} \! + \! 2Q_{e} \! = \! -2 \int_{a_{j}^{e}}^{
z}(R_{e}(s))^{1/2}h_{V}^{e}(s) \, \md s$ $(< \! 0)$. Recalling, also, that
$(R_{e}(z))^{1/2} \! := \! (\prod_{k=1}^{N+1}(z \! - \! b_{k-1}^{e})(z \! -
\! a_{k}^{e}))^{1/2}$ is continuous (and bounded) on the compact intervals
$[a_{j}^{e},b_{j}^{e}] \supset \Sigma_{p,j}^{e,2} \! \ni \! z$, $j \! = \! 1,
\dotsc,N$, vanishes at the end-points $\{a_{j}^{e}\}_{j=1}^{N}$ (resp., $\{
b_{j}^{e}\}_{j=1}^{N})$ like $(R_{e}(z))^{1/2} \! =_{z \downarrow a_{j}^{e}}
\! \mathcal{O}((z \! - \! a_{j}^{e})^{1/2})$ (resp., $(R_{e}(z))^{1/2} \! =_{
z \uparrow b_{j}^{e}} \! \mathcal{O}((b_{j}^{e} \! - \! z)^{1/2}))$, and is
differentiable on the open intervals $\Sigma_{p,j}^{e,2} \! \ni \! z$, and
$h_{V}^{e}(z) \! = \! \tfrac{1}{2} \oint_{C_{\mathrm{R}}^{e}}(\tfrac{\mi}{
\pi s} \! + \! \tfrac{\mi \widetilde{V}^{\prime}(s)}{2 \pi})(R_{e}(s))^{-1/2}
(s \! - \! z)^{-1} \, \md s$ is analytic, it follows that, for $z \! \in \!
\Sigma_{p,j}^{e,2}$,
\begin{equation*}
\inf_{z \in \Sigma_{p,j}^{e,2}}(R_{e}(z))^{1/2} \! =: \! \widehat{m}_{j} \!
\leqslant \! (R_{e}(z))^{1/2} \! \leqslant \! \widehat{M}_{j} \! := \! \sup_{
z \in \Sigma_{p,j}^{e,2}}(R_{e}(z))^{1/2}, \quad j \! = \! 1,\dotsc,N;
\end{equation*}
thus, recalling the expression for
$\overset{e}{m}^{\raise-1.0ex\hbox{$\scriptstyle \infty$}}(z)$ given in
Lemma~4.5 and noting that the respective factors $\gamma^{e}(z) \! \pm \!
(\gamma^{e}(z))^{-1}$ and $\boldsymbol{\theta}^{e}(\pm \boldsymbol{u}^{e}(z)
\! - \! \tfrac{n}{2 \pi} \boldsymbol{\Omega}^{e} \! \pm \! \boldsymbol{d}_{
e})$ are uniformly bounded (with respect to $z)$ in compact subsets outside
the open intervals surrounding the end-points of the suppport of the `even'
equilibrium measure, and defining $\mathbb{U}^{e}_{0}$ as in the Proposition,
after a straightforward integration argument, one arrives at the asymptotic
(as $n \! \to \! \infty)$ estimates for $\upsilon_{\mathscr{R}}^{e}(z)$ on
$\Sigma_{p,j}^{e,2} \setminus \mathbb{U}^{e}_{0} \! \ni \! z$ and $\Sigma^{e,
2}_{p,j} \cap \mathbb{U}^{e}_{0} \! \ni \! z$, $j \! = \! 1,\dotsc,N$, stated
in item~(2) of the Proposition (the $n$-dependence of the $\operatorname{GL}_{
2}(\mathbb{C})$-valued factors $f_{j}(n),\widetilde{f}_{j}(n)$, $j \! = \! 1,
\dotsc,N$, is inherited {}from the bounded $(\mathcal{O}(1))$ $n$-dependence
of the respective Riemann theta functions).

For $z \! \in \! \Sigma_{p}^{e,3} \! := \! \cup_{j=1}^{N+1}(J_{j}^{e,
\smallfrown} \setminus (\mathbb{C}_{+} \cap (\mathbb{U}_{\delta_{b_{j-1}}}^{e}
\cup \mathbb{U}_{\delta_{a_{j}}}^{e})))$, recall {}from the proof of
Lemma~4.8 that
\begin{equation*}
\upsilon_{\mathscr{R}}^{e}(z) \! := \! \upsilon_{\mathscr{R}}^{e,3}(z) \! = \!
\mathrm{I} \! + \! \exp \! \left(-4n \pi \mi \int_{z}^{a_{N+1}^{e}} \psi_{V}^{
e}(s) \, \md s \right) \!
\overset{e}{m}^{\raise-1.0ex\hbox{$\scriptstyle \infty$}}(z) \sigma_{-}
(\overset{e}{m}^{\raise-1.0ex\hbox{$\scriptstyle \infty$}}(z))^{-1},
\end{equation*}
and, {}from Lemma~4.1, $\Re (\mi \int_{z}^{a_{N+1}^{e}} \psi_{V}^{e}(s) \, \md
s) \! > \! 0$ for $z \! \in \! \mathbb{C}_{+} \cap (\cup_{j=1}^{N+1} \mathbb{
U}_{j}^{e}) \supset \Sigma_{p}^{e,3}$, where $\mathbb{U}_{j}^{e} \! := \!
\lbrace \mathstrut z \! \in \! \mathbb{C}^{\ast}; \, \Re (z) \! \in \! (b_{j
-1}^{e},a_{j}^{e}), \, \inf_{q \in (b_{j-1}^{e},a_{j}^{e})} \vert z \! - \! q
\vert \! < \! r_{j} \! \in \! (0,1) \rbrace$, $j \! = \! 1,\dotsc,N \! + \!
1$, with $\mathbb{U}_{i}^{e} \cap \mathbb{U}_{j}^{e} \! = \! \varnothing$, $i
\! \not= \! j \! = \! 1,\dotsc,N \! + \! 1$: using the expression for
$\overset{e}{m}^{\raise-1.0ex\hbox{$\scriptstyle \infty$}}(z)$
given in Lemma~4.5 and noting that the respective factors $\gamma^{e}(z) \!
\pm \! (\gamma^{e}(z))^{-1}$ and $\boldsymbol{\theta}^{e}(\pm \boldsymbol{u}^{
e}(z) \! - \! \tfrac{n}{2 \pi} \boldsymbol{\Omega}^{e} \! \pm \! \boldsymbol{
d}_{e})$ are uniformly bounded (with respect to $z)$ in compact subsets
outside the open intervals surrounding the end-points of the suppport of
the `even' equilibrium measure, an arc-length-parametrisation argument,
complemented by an application of the Maximum Length $(\mathrm{ML})$ Theorem,
leads one directly to the asymptotic (as $n \! \to \! \infty)$ estimate for
$\upsilon_{\mathscr{R}}^{e}(z)$ on $\Sigma_{p}^{e,3} \! \ni \! z$ stated in
item~(3) of the Proposition (the $n$-dependence of the $\operatorname{GL}_{2}
(\mathbb{C})$-valued factor $\overset{\smallfrown}{f}(n)$ is inherited {}from
the bounded $(\mathcal{O}(1))$ $n$-dependence of the respective Riemann theta
functions). The above argument applies, \emph{mutatis mutandis}, for the
asymptotic estimate of $\upsilon_{\mathscr{R}}^{e}(z)$ on $\Sigma_{p}^{e,4} \!
:= \! \cup_{j=1}^{N+1}(J_{j}^{e,\smallsmile} \setminus (\mathbb{C}_{-} \cap
(\mathbb{U}_{\delta_{b_{j-1}}}^{e} \cup \mathbb{U}_{\delta_{a_{j}}}^{e}))) \!
\ni \! z$ stated in item~(4) of the Proposition.

Since the estimates in item~(5) of the Proposition are similar, consider,
say, and without loss of generality, the asymptotic (as $n \! \to \! \infty)$
estimate for $\upsilon_{\mathscr{R}}^{e}(z)$ on $\partial \mathbb{U}_{\delta_{
a_{j}}}^{e} \! \ni \! z$, $j \! = \! 1,\dotsc,N \! + \! 1$: this argument
applies, \emph{mutatis mutandis}, for the large-$n$ asymptotics of $\upsilon^{
e}_{\mathscr{R}}(z)$ on $\cup_{j=1}^{N+1} \partial \mathbb{U}_{\delta_{b_{j-
1}}}^{e} \! \ni \! z$. For $z \! \in \! \partial \mathbb{U}_{\delta_{a_{j}}}^{
e}$, $j \! = \! 1,\dotsc,N \! + \! 1$, recall {}from the proof of Lemma~4.8
that $\upsilon^{e}_{\mathscr{R}}(z) \! := \! \upsilon_{\mathscr{R}}^{e,5}(z)
\! = \! \mathcal{X}^{e}(z)
(\overset{e}{m}^{\raise-1.0ex\hbox{$\scriptstyle \infty$}}(z))^{-1}$: using
the expression for the parametrix, $\mathcal{X}^{e}(z)$, given in Lemma~4.7,
and the large-argument asymptotics for the Airy function and its derivative
given in Equations~(2.6), one shows that, for $z \! \in \! \mathbb{C}_{+} \cap
\partial \mathbb{U}_{\delta_{a_{j}}}^{e}$, $j \! = \! 1,\dotsc,N \! + \! 1$,
\begin{align*}
\upsilon_{\mathscr{R}}^{e}(z) \underset{n \to \infty}{=}& \, \mathrm{I} \! +
\! \dfrac{\me^{-\frac{\mi \pi}{3}}}{n \xi_{a_{j}}^{e}(z)}
\overset{e}{m}^{\raise-1.0ex\hbox{$\scriptstyle \infty$}}(z) \!
\begin{pmatrix}
\mi \me^{\frac{\mi}{2}n \mho_{j}^{e}} & -\mi \me^{\frac{\mi}{2}n \mho_{j}^{e}}
\\
\me^{-\frac{\mi}{2}n \mho_{j}^{e}} & \me^{-\frac{\mi}{2}n \mho_{j}^{e}}
\end{pmatrix} \!
\begin{pmatrix}
-s_{1} \me^{-\frac{\mi \pi}{6}} \me^{-\frac{\mi}{2}n \mho_{j}^{e}} & s_{1}
\me^{\frac{\mi \pi}{3}} \me^{\frac{\mi}{2}n \mho_{j}^{e}} \\
t_{1} \me^{-\frac{\mi \pi}{6}} \me^{-\frac{\mi}{2}n \mho_{j}^{e}} & -t_{1}
\me^{\frac{4 \pi \mi}{3}} \me^{\frac{\mi}{2}n \mho_{j}^{e}}
\end{pmatrix} \\
\times& \, (\overset{e}{m}^{\raise-1.0ex\hbox{$\scriptstyle \infty$}}(z))^{-1}
\! + \! \mathcal{O} \! \left(\dfrac{1}{(n \xi_{a_{j}}^{e}(z))^{2}}
\overset{e}{m}^{\raise-1.0ex\hbox{$\scriptstyle \infty$}}(z) \!
\begin{pmatrix}
\ast & \ast \\
\ast & \ast
\end{pmatrix} \! (\overset{e}{m}^{\raise-1.0ex\hbox{$\scriptstyle \infty$}}
(z))^{-1} \right),
\end{align*}
where $\xi_{a_{j}}^{e}(z)$ and $\mho_{j}^{e}$, $j \! = \! 1,\dotsc,N \! + \!
1$, and $s_{1}$ and $t_{1}$ are defined in the Proposition,
$\overset{e}{m}^{\raise-1.0ex\hbox{$\scriptstyle \infty$}}(z)$ is given in
Lemma~4.5, and
$\left(
\begin{smallmatrix}
\ast & \ast \\
\ast & \ast
\end{smallmatrix}
\right) \! \in \! \operatorname{M}_{2}(\mathbb{C})$, and, for $z \! \in \!
\mathbb{C}_{-} \cap \partial \mathbb{U}_{\delta_{a_{j}}}^{e}$, $j \! = \! 1,
\dotsc,N \! + \! 1$,
\begin{align*}
\upsilon_{\mathscr{R}}^{e}(z) \underset{n \to \infty}{=}& \, \mathrm{I} \!
+ \! \dfrac{\me^{-\frac{\mi \pi}{3}}}{n \xi_{a_{j}}^{e}(z)}
\overset{e}{m}^{\raise-1.0ex\hbox{$\scriptstyle \infty$}}(z) \!
\begin{pmatrix}
\mi \me^{-\frac{\mi}{2}n \mho_{j}^{e}} & -\mi \me^{-\frac{\mi}{2}n \mho_{j}^{
e}} \\
\me^{\frac{\mi}{2}n \mho_{j}^{e}} & \me^{\frac{\mi}{2}n \mho_{j}^{e}}
\end{pmatrix} \!
\begin{pmatrix}
-s_{1} \me^{-\frac{\mi \pi}{6}} \me^{\frac{\mi}{2}n \mho_{j}^{e}} & s_{1}
\me^{\frac{\mi \pi}{3}} \me^{-\frac{\mi}{2}n \mho_{j}^{e}} \\
t_{1} \me^{-\frac{\mi \pi}{6}} \me^{\frac{\mi}{2}n \mho_{j}^{e}} & -t_{1}
\me^{\frac{4 \pi \mi}{3}} \me^{-\frac{\mi}{2}n \mho_{j}^{e}}
\end{pmatrix} \\
\times& \, (\overset{e}{m}^{\raise-1.0ex\hbox{$\scriptstyle \infty$}}(z))^{-1}
\! + \! \mathcal{O} \! \left(\dfrac{1}{(n \xi_{a_{j}}^{e}(z))^{2}}
\overset{e}{m}^{\raise-1.0ex\hbox{$\scriptstyle \infty$}}(z) \!
\begin{pmatrix}
\ast & \ast \\
\ast & \ast
\end{pmatrix} \! (\overset{e}{m}^{\raise-1.0ex\hbox{$\scriptstyle \infty$}}
(z))^{-1} \right).
\end{align*}
Upon recalling the formula for
$\overset{e}{m}^{\raise-1.0ex\hbox{$\scriptstyle \infty$}}(z)$ in terms of
$\overset{e}{\mathfrak{M}}^{\raise-1.0ex\hbox{$\scriptstyle \infty$}}(z)$
given in Lemma~4.5, and noting that the respective factors $\gamma^{e}(z) \!
\pm \! (\gamma^{e}(z))^{-1}$ and $\boldsymbol{\theta}^{e}(\pm \boldsymbol{u}^{
e}(z) \! - \! \tfrac{n}{2 \pi} \boldsymbol{\Omega}^{e} \! \pm \! \boldsymbol{
d}_{e})$ are uniformly bounded (with respect to $z)$ in compact subsets
outside the open intervals surrounding the end-points of the suppport of the
`even' equilibrium measure, after a straightforward matrix-multiplication
argument, one arrives at the asymptotic (as $n \! \to \! \infty)$ estimates
for $\upsilon^{e}_{\mathscr{R}}(z)$ on $\partial \mathbb{U}_{\delta_{a_{j}}
}^{e} \! \ni \! z$, $j \! = \! 1,\dotsc,N \! + \! 1$, stated in item~(5) of
the Proposition (the $n$-dependence of the $\operatorname{GL}_{2}
(\mathbb{C})$-valued factors $f_{a_{j}}^{e}(n)$, $j \! = \! 1,\dotsc,N \! +
\! 1$, is inherited {}from the bounded $(\mathcal{O}(1))$ $n$-dependence of
the respective Riemann theta functions). \hfill $\qed$
\begin{aaaaa}
For an oriented contour $D \! \subset \! \mathbb{C}$, let $\mathscr{N}_{q}
(D)$ denote the set of all bounded linear operators {}from $\mathcal{L}^{
q}_{\mathrm{M}_{2}(\mathbb{C})}(D)$ into $\mathcal{L}^{q}_{\mathrm{M}_{2}
(\mathbb{C})}(D)$, $q \! \in \! \lbrace 1,2,\infty \rbrace$.
\end{aaaaa}

Since the analysis that follows relies substantially on the BC \cite{a84} 
construction for the solution of a matrix (and suitably normalised) RHP on 
an oriented and unbounded contour, it is convenient to present, with some 
requisite preamble, a succinct and self-contained synopsis of it at this 
juncture. One agrees to call a contour $\Gamma^{\sharp}$ \emph{oriented} 
if:
\begin{compactenum}
\item[(1)] $\mathbb{C} \setminus \Gamma^{\sharp}$ has finitely many open
connected components;
\item[(2)] $\mathbb{C} \setminus \Gamma^{\sharp}$ is the disjoint union of
two, possibly disconnected, open regions, denoted by $\boldsymbol{\mho}^{+}$
and $\boldsymbol{\mho}^{-}$;
\item[(3)] $\Gamma^{\sharp}$ may be viewed as either the positively oriented
boundary for $\boldsymbol{\mho}^{+}$ or the negatively oriented boundary for
$\boldsymbol{\mho}^{-}$ ($\mathbb{C} \setminus \Gamma^{\sharp}$ is coloured
by two colours, $\pm)$.
\end{compactenum}
Let $\Gamma^{\sharp}$, as a closed set, be the union of finitely many 
oriented, simple, piecewise-smooth arcs. Denote the set of all 
self-intersections of $\Gamma^{\sharp}$ by $\widehat{\Gamma}^{\sharp}$ 
(with $\mathrm{card}(\widehat{\Gamma}^{\sharp}) \! < \! \infty$ assumed
throughout). Set $\widetilde{\Gamma}^{\sharp} \! := \! \Gamma^{\sharp} 
\setminus \widehat{\Gamma}^{\sharp}$. The BC \cite{a84} construction for the 
solution of a (matrix) RHP, in the absence of a discrete spectrum and spectral 
singularities \cite{a101} (see, also, \cite{a85,a86,a102,a103,a104}), on an 
oriented contour $\Gamma^{\sharp}$ consists of finding function $\mathcal{Y} 
\colon \mathbb{C} \setminus \Gamma^{\sharp} \! \to \! \mathrm{M}_{2}
(\mathbb{C})$ such that:
\begin{compactenum}
\item[(1)] $\mathcal{Y}(z)$ is holomorphic for $z \! \in \! \mathbb{C}
\setminus \Gamma^{\sharp}$, $\mathcal{Y} \! \! \upharpoonright_{\mathbb{C}
\setminus \Gamma^{\sharp}}$ has a continuous extension ({}from `above' and
`below') to $\widetilde{\Gamma}^{\sharp}$, and $\lim_{\genfrac{}{}{0pt}{2}{
z^{\prime} \to z}{z^{\prime} \, \in \, \pm \, \mathrm{side} \, \mathrm{of}
\, \widetilde{\Gamma}^{\sharp}}} \int_{\widetilde{\Gamma}^{\sharp}} \vert
\mathcal{Y}(z^{\prime}) \! - \! \mathcal{Y}_{\pm}(z) \vert^{2} \, \vert \md
z \vert \! = \! 0$;
\item[(2)] $\mathcal{Y}_{\pm}(z) \! := \! \lim_{\underset{z^{\prime} \, \in \,
\pm \, \mathrm{side} \, \mathrm{of} \, \widetilde{\Gamma}^{\sharp}}{z^{\prime}
\to z}} \mathcal{Y}(z^{\prime})$ satisfy $\mathcal{Y}_{+}(z) \! = \! \mathcal{
Y}_{-}(z) \upsilon (z)$, $z \! \in \! \widetilde{\Gamma}^{\sharp}$, for some
(smooth) `jump' matrix $\upsilon \colon \widetilde{\Gamma}^{\sharp} \! \to \!
\mathrm{GL}_{2}(\mathbb{C})$; and
\item[(3)] for arbitrarily fixed $\lambda_{o} \! \in \! \mathbb{C}$, and
uniformly with respect to $z$, $\mathcal{Y}(z) \! =_{\underset{z \in \mathbb{
C} \setminus \Gamma^{\sharp}}{z \to \lambda_{o}}} \! \mathrm{I} \! + \! o(1)$,
where $o(1) \! = \! \mathcal{O}(z \! - \! \lambda_{o})$ if $\lambda_{o}$ is
finite, and $o(1) \! = \! \mathcal{O}(z^{-1})$ if $\lambda_{o}$ is the point
at infinity).
\end{compactenum}
(Condition~(3) is referred to as the \emph{normalisation condition}, and is
necessary in order to prove uniqueness of the associated RHP: one says that
the RHP is `normalised at $\lambda_{o}$'.) Let $\upsilon (z) \! := \!
(\mathrm{I} \! - \! w_{-}(z))^{-1}(\mathrm{I} \! + \! w_{+}(z))$, $z \! \in
\! \widetilde{\Gamma}^{\sharp}$, be a (bounded algebraic) factorisation for
$\upsilon (z)$, where $w_{\pm}(z)$ are some upper/lower, or lower/upper,
triangular matrices (depending on the orientation of $\Gamma^{\sharp})$, and
$w_{\pm}(z) \! \in \! \cap_{p \in \{2,\infty\}} \mathcal{L}^{p}_{\mathrm{
M}_{2}(\mathbb{C})} \linebreak[4]
(\widetilde{\Gamma}^{\sharp})$ (if $\widetilde{\Gamma}^{\sharp}$ is unbounded,
one requires that $w_{\pm}(z) \! =_{\genfrac{}{}{0pt}{2}{z \to \infty}{z \in
\widetilde{\Gamma}^{\sharp}}} \! \pmb{0})$. Define $w(z) \! := \! w_{+}(z)
\! + \! w_{-}(z)$, and introduce the (normalised at $\lambda_{o})$ Cauchy
operators
\begin{equation*}
\mathcal{L}^{2}_{\mathrm{M}_{2}(\mathbb{C})}(\Gamma^{\sharp}) \! \ni \! f \!
\mapsto \! (C^{\lambda_{o}}_{\pm}f)(z) \! := \!
\lim_{\genfrac{}{}{0pt}{2}{z^{\prime} \to z}{z^{\prime} \, \in \, \pm \,
\mathrm{side} \, \mathrm{of} \, \Gamma^{\sharp}}} \int_{\Gamma^{\sharp}}
\dfrac{(z^{\prime} \! - \! \lambda_{o})f(\zeta)}{(\zeta \! - \! \lambda_{o})
(\zeta \! - \! z^{\prime})} \, \dfrac{\md \zeta}{2 \pi \mi},
\end{equation*}
where $\tfrac{(z-\lambda_{o})}{(\zeta -\lambda_{o})(\zeta -z)} \, \tfrac{\md
\zeta}{2 \pi \mi}$ is the Cauchy kernel normalised at $\lambda_{o}$ (which
reduces to the `standard' Cauchy kernel, that is, $\tfrac{1}{\zeta -z} \,
\tfrac{\md \zeta}{2 \pi \mi}$, in the limit $\lambda_{o} \! \to \! \infty)$,
with $C^{\lambda_{o}}_{\pm} \colon \mathcal{L}^{2}_{\mathrm{M}_{2}(\mathbb{C}
)}(\Gamma^{\sharp}) \! \to \! \mathcal{L}^{2}_{\mathrm{M}_{2}(\mathbb{C})}
(\Gamma^{\sharp})$ bounded in operator norm\footnote{$\vert \vert C^{\lambda_{
o}}_{\pm} \vert \vert_{\mathscr{N}_{2}(\Gamma^{\sharp})} \! < \! \infty$.},
and $\vert \vert (C^{\lambda_{o}}_{\pm}f)(\cdot) \vert \vert_{\mathcal{
L}^{2}_{\mathrm{M}_{2}(\mathbb{C})}(\Gamma^{\sharp})} \! \leqslant \!
\mathrm{const.} \vert \vert f(\cdot) \vert \vert_{\mathcal{L}^{2}_{\mathrm{
M}_{2}(\mathbb{C})}(\Gamma^{\sharp})}$. Introduce the BC operator $C^{
\lambda_{o}}_{w}$:
\begin{equation*}
\mathcal{L}^{2}_{\mathrm{M}_{2}(\mathbb{C})}(\Gamma^{\sharp}) \! \ni \! f \!
\mapsto \! C^{\lambda_{o}}_{w}f \! := \! C^{\lambda_{o}}_{+}(fw_{-}) \! + \!
C^{\lambda_{o}}_{-}(fw_{+}),
\end{equation*}
which, for $w_{\pm}\! \in \! \mathcal{L}^{\infty}_{\mathrm{M}_{2}(\mathbb{C})}
(\Gamma^{\sharp})$, is bounded {}from $\mathcal{L}^{2}_{\mathrm{M}_{2}(\mathbb{
C})}(\Gamma^{\sharp}) \! \to \! \mathcal{L}^{2}_{\mathrm{M}_{2}(\mathbb{C})}
(\Gamma^{\sharp})$, that is, $\vert \vert C^{\lambda_{o}}_{w} \vert \vert_{
\mathscr{N}_{2}(\Gamma^{\sharp})} \! < \! \infty$; furthermore, since 
$\mathbb{C} \setminus \Gamma^{\sharp}$ can be coloured by the two colours 
$\pm$, $C^{\lambda_{o}}_{\pm}$ are complementary projections \cite{a2,a85,%
a102,a103}, that is, $(C^{\lambda_{o}}_{+})^{2} \! = \! C^{\lambda_{o}}_{+}$, 
$(C^{\lambda_{o}}_{-})^{2} \! = \! -C^{\lambda_{o}}_{-}$, $C^{\lambda_{o}}_{
+}C^{\lambda_{o}}_{-} \! = \! C^{\lambda_{o}}_{-}C^{\lambda_{o}}_{+} \! = \! 
\underline{\pmb{0}}$ (the null operator), and $C^{\lambda_{o}}_{+} \! - \! 
C^{\lambda_{o}}_{-} \! = \! \id$ (the identity operator). (In the case that 
$C^{\lambda_{o}}_{+}$ and $-C^{\lambda_{o}}_{-}$ are complementary, the 
contour $\Gamma^{\sharp}$ can always be oriented in such a way that the $\pm$ 
regions lie on the $\pm$ sides of the contour, respectively.) The solution of 
the above (normalised at $\lambda_{o})$ RHP is given by the following integral 
representation.
\begin{ccccc}[Beals and Coifman {\rm \cite{a84}}]
Set
\begin{equation*}
\mu_{\lambda_{o}}(z) \! = \! \mathcal{Y}_{+}(z) \! \left(\mathrm{I} \! + \!
w_{+}(z) \right)^{-1} \! = \! \mathcal{Y}_{-}(z) \! \left(\mathrm{I} \! - \!
w_{-}(z) \right)^{-1}, \quad z \! \in \! \Gamma^{\sharp}.
\end{equation*}
If $\mu_{\lambda_{o}} \! \in \! \mathrm{I} \! + \! \mathcal{L}^{2}_{\mathrm{
M}_{2}(\mathbb{C})}(\Gamma^{\sharp})$ solves the linear singular integral
equation
\begin{equation*}
\left(\id \! - \! C^{\lambda_{o}}_{w} \right) \! \left(\mu_{\lambda_{o}}(z) \!
- \! \mathrm{I} \right) \! = \! C_{w}^{\lambda_{o}} \mathrm{I} \! = \! C^{
\lambda_{o}}_{+}(w_{-}(z)) \! + \! C^{\lambda_{o}}_{-}(w_{+}(z)), \quad z \!
\in \! \Gamma^{\sharp},
\end{equation*}
where $\id$ is the identity operator on $\mathcal{L}^{2}_{\mathrm{M}_{2}
(\mathbb{C})}(\Gamma^{\sharp})$, then the solution of the {\rm RHP}
$(\mathcal{Y}(z),\upsilon (z),\Gamma^{\sharp})$ is given by
\begin{equation*}
\mathcal{Y}(z) \! = \! \mathrm{I} \! + \! \int\nolimits_{\Gamma^{\sharp}}
\dfrac{(z \! - \! \lambda_{o}) \mu_{\lambda_{o}}(\zeta)w(\zeta)}{(\zeta \! -
\! \lambda_{o})(\zeta \! - \! z)} \, \dfrac{\md \zeta}{2 \pi \mi}, \quad z \!
\in \! \mathbb{C} \setminus \Gamma^{\sharp},
\end{equation*}
where $\mu_{\lambda_{o}}(z) \! := \! ((\id \! - \! C^{\lambda_{o}}_{w})^{-1}
\mathrm{I})(z)$\footnote{The linear singular integral equation for $\mu_{
\lambda_{o}}(\pmb{\cdot})$ stated in this Lemma~5.1 is well defined in
$\mathcal{L}^{2}_{\mathrm{M}_{2}(\mathbb{C})}(\Gamma^{\sharp})$ provided that
$w_{\pm}(\cdot) \! \in \! \mathcal{L}^{2}_{\mathrm{M}_{2}(\mathbb{C})}
(\Gamma^{\sharp}) \cap \mathcal{L}^{\infty}_{\mathrm{M}_{2}(\mathbb{C})}
(\Gamma^{\sharp})$; furthermore, it is assumed that the associated RHP
$(\mathcal{Y}(z),\upsilon (z),\Gamma^{\sharp})$ is solvable, that is,
$\mathrm{dim} \, \mathrm{ker}(\id \! - \! C^{\lambda_{o}}_{w}) \! = \!
\mathrm{dim} \left\lbrace \mathstrut \phi \! \in \! \mathcal{L}^{2}_{\mathrm{
M}_{2}(\mathbb{C})}(\Gamma^{\sharp}); \, (\id \! - \! C^{\lambda_{o}}_{w})
\phi \! = \! \underline{\pmb{0}} \right\rbrace \! = \! \mathrm{dim} \,
\varnothing \! = \! 0$ $(\Rightarrow  (\id \! - \! C^{\lambda_{o}}_{w})^{-1}
\! \! \upharpoonright_{\mathcal{L}^{2}_{\mathrm{M}_{2}(\mathbb{C})}(\Gamma^{
\sharp})}$ exists).}.
\end{ccccc}

Recall that $\mathscr{R}^{e} \colon \mathbb{C} \setminus \Sigma_{p}^{e} \! 
\to \! \operatorname{SL}_{2}(\mathbb{C})$, which solves the RHP $(\mathscr{
R}^{e}(z),\upsilon_{\mathscr{R}}^{e}(z),\Sigma_{p}^{e})$ formulated in Lemma 
4.8, is normalised at infinity, that is, $\mathscr{R}^{e}(\infty) \! = \! 
\mathrm{I}$. Removing {}from the specification of the RHP $(\mathscr{R}^{e}
(z),\upsilon_{\mathscr{R}}^{e} \linebreak[4]
(z),\Sigma_{p}^{e})$ the oriented skeletons on which the jump matrix, 
$\upsilon_{\mathscr{R}}^{e}(z)$, is equal to $\mathrm{I}$, in particular (cf. 
Lemma~4.8), the oriented skeleton $\Sigma_{p}^{e} \setminus \cup_{l=1}^{5} 
\Sigma_{p}^{e,l}$, and setting $\Sigma_{p}^{e} \setminus (\Sigma_{p}^{e} 
\setminus \cup_{l=1}^{5} \Sigma_{p}^{e,l}) \! =: \! \widetilde{\Sigma}_{
p}^{e}$ (see Figure~10),
\begin{figure}[tbh]
\begin{center}
\begin{pspicture}(0,0)(15,5)
\psset{xunit=1cm,yunit=1cm}
\psdots[dotstyle=*,dotscale=1.5](1.3,2.5)
\psdots[dotstyle=*,dotscale=1.5](3.7,2.5)
\psdots[dotstyle=*,dotscale=1.5](6.3,2.5)
\psdots[dotstyle=*,dotscale=1.5](8.7,2.5)
\psdots[dotstyle=*,dotscale=1.5](11.3,2.5)
\psdots[dotstyle=*,dotscale=1.5](13.7,2.5)
\psline[linewidth=0.6pt,linestyle=solid,linecolor=black](0.4,2.5)(0.7,2.5)
\psline[linewidth=0.6pt,linestyle=solid,linecolor=black,arrowsize=1.5pt 3]%
{->}(0,2.5)(0.4,2.5)
\psline[linewidth=0.6pt,linestyle=solid,linecolor=black,arrowsize=1.5pt 3]%
{->}(4.3,2.5)(4.55,2.5)
\psline[linewidth=0.6pt,linestyle=solid,linecolor=black](4.55,2.5)(4.75,2.5)
\psarcn[linewidth=0.6pt,linestyle=solid,linecolor=magenta,arrowsize=1.5pt 5]%
{->}(2.5,1.5){1.8}{121}{89}
\psarcn[linewidth=0.6pt,linestyle=solid,linecolor=magenta](2.5,1.5){1.8}{90}%
{59}
\psarc[linewidth=0.6pt,linestyle=solid,linecolor=magenta,arrowsize=1.5pt 5]%
{->}(2.5,3.5){1.8}{239}{271}
\psarc[linewidth=0.6pt,linestyle=solid,linecolor=magenta](2.5,3.5){1.8}{270}%
{301}
\psarcn[linewidth=0.6pt,linestyle=solid,linecolor=cyan,arrowsize=1.5pt 5]%
{->}(1.3,2.5){0.6}{180}{135}
\psarcn[linewidth=0.6pt,linestyle=solid,linecolor=cyan](1.3,2.5){0.6}{135}{0}
\psarcn[linewidth=0.6pt,linestyle=solid,linecolor=cyan](1.3,2.5){0.6}{360}{180}
\psarcn[linewidth=0.6pt,linestyle=solid,linecolor=cyan,arrowsize=1.5pt 5]%
{->}(3.7,2.5){0.6}{180}{45}
\psarcn[linewidth=0.6pt,linestyle=solid,linecolor=cyan](3.7,2.5){0.6}{45}{0}
\psarcn[linewidth=0.6pt,linestyle=solid,linecolor=cyan](3.7,2.5){0.6}{360}{180}
\psline[linewidth=0.6pt,linestyle=solid,linecolor=black](5.5,2.5)(5.7,2.5)
\psline[linewidth=0.6pt,linestyle=solid,linecolor=black,arrowsize=1.5pt 3]%
{->}(5.25,2.5)(5.5,2.5)
\psline[linewidth=0.6pt,linestyle=solid,linecolor=black,arrowsize=1.5pt 3]%
{->}(9.3,2.5)(9.55,2.5)
\psline[linewidth=0.6pt,linestyle=solid,linecolor=black](9.55,2.5)(9.75,2.5)
\psarcn[linewidth=0.6pt,linestyle=solid,linecolor=magenta,arrowsize=1.5pt 5]%
{->}(7.5,1.5){1.8}{121}{89}
\psarcn[linewidth=0.6pt,linestyle=solid,linecolor=magenta](7.5,1.5){1.8}{90}%
{59}
\psarc[linewidth=0.6pt,linestyle=solid,linecolor=magenta,arrowsize=1.5pt 5]%
{->}(7.5,3.5){1.8}{239}{271}
\psarc[linewidth=0.6pt,linestyle=solid,linecolor=magenta](7.5,3.5){1.8}{270}%
{301}
\psarcn[linewidth=0.6pt,linestyle=solid,linecolor=cyan,arrowsize=1.5pt 5]%
{->}(6.3,2.5){0.6}{180}{135}
\psarcn[linewidth=0.6pt,linestyle=solid,linecolor=cyan](6.3,2.5){0.6}{135}{0}
\psarcn[linewidth=0.6pt,linestyle=solid,linecolor=cyan](6.3,2.5){0.6}{360}{180}
\psarcn[linewidth=0.6pt,linestyle=solid,linecolor=cyan,arrowsize=1.5pt 5]%
{->}(8.7,2.5){0.6}{180}{45}
\psarcn[linewidth=0.6pt,linestyle=solid,linecolor=cyan](8.7,2.5){0.6}{45}{0}
\psarcn[linewidth=0.6pt,linestyle=solid,linecolor=cyan](8.7,2.5){0.6}{360}{180}
\psline[linewidth=0.6pt,linestyle=solid,linecolor=black](10.5,2.5)(10.7,2.5)
\psline[linewidth=0.6pt,linestyle=solid,linecolor=black,arrowsize=1.5pt 3]%
{->}(10.25,2.5)(10.5,2.5)
\psline[linewidth=0.6pt,linestyle=solid,linecolor=black,arrowsize=1.5pt 3]%
{->}(14.3,2.5)(14.65,2.5)
\psline[linewidth=0.6pt,linestyle=solid,linecolor=black](14.65,2.5)(14.95,2.5)
\psarcn[linewidth=0.6pt,linestyle=solid,linecolor=magenta,arrowsize=1.5pt 5]%
{->}(12.5,1.5){1.8}{121}{89}
\psarcn[linewidth=0.6pt,linestyle=solid,linecolor=magenta](12.5,1.5){1.8}{90}%
{59}
\psarc[linewidth=0.6pt,linestyle=solid,linecolor=magenta,arrowsize=1.5pt 5]%
{->}(12.5,3.5){1.8}{239}{271}
\psarc[linewidth=0.6pt,linestyle=solid,linecolor=magenta](12.5,3.5){1.8}{270}%
{301}
\psarcn[linewidth=0.6pt,linestyle=solid,linecolor=cyan,arrowsize=1.5pt 5]%
{->}(11.3,2.5){0.6}{180}{135}
\psarcn[linewidth=0.6pt,linestyle=solid,linecolor=cyan](11.3,2.5){0.6}{135}{0}
\psarcn[linewidth=0.6pt,linestyle=solid,linecolor=cyan](11.3,2.5){0.6}{360}%
{180}
\psarcn[linewidth=0.6pt,linestyle=solid,linecolor=cyan,arrowsize=1.5pt 5]%
{->}(13.7,2.5){0.6}{180}{45}
\psarcn[linewidth=0.6pt,linestyle=solid,linecolor=cyan](13.7,2.5){0.6}{45}{0}
\psarcn[linewidth=0.6pt,linestyle=solid,linecolor=cyan](13.7,2.5){0.6}{360}%
{180}
\psline[linewidth=0.9pt,linestyle=dotted,linecolor=darkgray](4.8,2.5)(5.2,2.5)
\psline[linewidth=0.9pt,linestyle=dotted,linecolor=darkgray](9.8,2.5)(10.2,2.5)
\rput(1.3,2.2){\makebox(0,0){$\scriptstyle \pmb{b_{0}^{e}}$}}
\rput(3.7,2.2){\makebox(0,0){$\scriptstyle \pmb{a_{1}^{e}}$}}
\rput(6.3,2.2){\makebox(0,0){$\scriptstyle \pmb{b_{j-1}^{e}}$}}
\rput(8.7,2.2){\makebox(0,0){$\scriptstyle \pmb{a_{j}^{e}}$}}
\rput(11.3,2.2){\makebox(0,0){$\scriptstyle \pmb{b_{N}^{e}}$}}
\rput(13.7,2.2){\makebox(0,0){$\scriptstyle \pmb{a_{N+1}^{e}}$}}
\end{pspicture}
\end{center}
\vspace{-1.05cm}
\caption{Oriented skeleton $\widetilde{\Sigma}_{p}^{e} \! := \! \Sigma_{p}^{e}
\setminus (\Sigma_{p}^{e} \setminus \cup_{l=1}^{5} \Sigma_{p}^{e,l})$}
\end{figure}
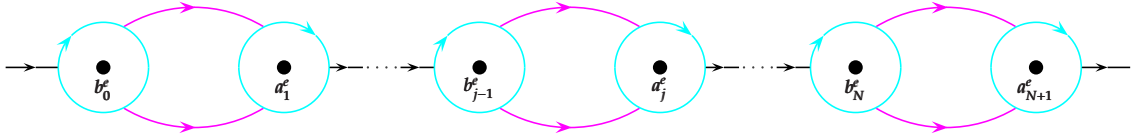
one arrives at the equivalent RHP $(\mathscr{R}^{e}(z),\upsilon_{\mathscr{
R}}^{e}(z),\widetilde{\Sigma}_{p}^{e})$ for $\mathscr{R}^{e} \colon \mathbb{C}
\setminus \widetilde{\Sigma}_{p}^{e} \! \to \! \operatorname{SL}_{2}(\mathbb{
C})$ (the normalisation at infinity, of course, remains unchanged). Via the BC 
\cite{a84} construction discussed above, write, for $\upsilon_{\mathscr{R}}^{
e} \colon \widetilde{\Sigma}_{p}^{e} \! \to \! \operatorname{SL}_{2}(\mathbb{
C})$, the (bounded algebraic) factorisation
\begin{equation*}
\upsilon_{\mathscr{R}}^{e}(z) \! := \! \left(\mathrm{I} \! - \! w_{-}^{
\Sigma_{\mathscr{R}}^{e}}(z) \right)^{-1} \! \left(\mathrm{I} \! + \! w_{+}^{
\Sigma_{\mathscr{R}}^{e}}(z) \right), \quad z \! \in \! \widetilde{\Sigma}_{
p}^{e}:
\end{equation*}
taking the (so-called) trivial factorisation \cite{a86} (see pp.~293 and~294,
\emph{Proof of Theorem}~3.14 and \emph{Proposition}~1.9; see, also, 
\cite{a103,a104}) $w^{\Sigma_{\mathscr{R}}^{e}}_{-}(z) \! \equiv \! \pmb{0}$, 
whence $\upsilon_{\mathscr{R}}^{e}(z) \! = \! \mathrm{I} \! + \! w^{\Sigma_{
\mathscr{R}}^{e}}_{+}(z)$, $z \! \in \! \widetilde{\Sigma}_{p}^{e}$, it 
follows {}from Lemma 5.1 that, upon normalising the Cauchy (integral) 
operator(s) at infinity (take the limit $\lambda_{o} \! \to \! \infty$ in 
Lemma 5.1), the $(\operatorname{SL}_{2}(\mathbb{C})$-valued) integral 
representation for the---unique---solution of the equivalent RHP $(\mathscr{
R}^{e}(z),\upsilon_{\mathscr{R}}^{e}(z),\widetilde{\Sigma}_{p}^{e})$ is
\begin{equation}
\mathscr{R}^{e}(z) \! = \! \mathrm{I} \! + \! \int_{\widetilde{\Sigma}_{p}^{
e}} \dfrac{\mu^{\Sigma_{\mathscr{R}}^{e}}(s) w^{\Sigma_{\mathscr{R}}^{e}}_{+}
(s)}{s \! - \! z} \, \dfrac{\md s}{2 \pi \mi}, \quad z \! \in \! \mathbb{C}
\setminus \widetilde{\Sigma}_{p}^{e},
\end{equation}
where $\mu^{\Sigma_{\mathscr{R}}^{e}}(\bm{\cdot}) \! \in \! \mathrm{I} \! + \!
\mathcal{L}^{2}_{\mathrm{M}_{2}(\mathbb{C})}(\widetilde{\Sigma}_{p}^{e})$
solves the (linear) singular integral equation
\begin{equation*}
(\id \! - \! C^{\infty}_{w^{\Sigma_{\mathscr{R}}^{e}}}) \mu^{\Sigma_{\mathscr{
R}}^{e}}(z) \! = \! \mathrm{I}, \quad z \! \in \! \widetilde{\Sigma}_{p}^{e},
\end{equation*}
with
\begin{equation*}
\mathcal{L}^{2}_{\mathrm{M}_{2}(\mathbb{C})}(\widetilde{\Sigma}_{p}^{e}) \!
\ni \! f \! \mapsto \! C^{\infty}_{w^{\Sigma_{\mathscr{R}}^{e}}}f \! := \!
C_{-}^{\infty}(fw^{\Sigma_{\mathscr{R}}^{e}}_{+}),
\end{equation*}
and
\begin{equation*}
\mathcal{L}^{2}_{\mathrm{M}_{2}(\mathbb{C})}(\widetilde{\Sigma}_{p}^{e}) \!
\ni \! f \! \mapsto \! (C_{\pm}^{\infty}f)(z) \! := \lim_{\underset{z^{\prime}
\, \in \, \pm \, \mathrm{side} \, \mathrm{of} \, \widetilde{\Sigma}_{p}^{e}}{
z^{\prime} \to z}} \int_{\widetilde{\Sigma}_{p}^{e}} \dfrac{f(s)}{s \! - \!
z^{\prime}} \, \dfrac{\md s}{2 \pi \mi};
\end{equation*}
furthermore, $\norm{C_{\pm}^{\infty}}_{\mathscr{N}_{2}(\widetilde{\Sigma}_{
p}^{e})} \! < \! \infty$.
\begin{bbbbb}
Let $\mathscr{R}^{e} \colon \mathbb{C} \setminus \widetilde{\Sigma}_{p}^{e}
\! \to \! \operatorname{SL}_{2}(\mathbb{C})$ solve the following, equivalent
{\rm RHP:} {\rm (i)} $\mathscr{R}^{e}(z)$ is holomorphic for $z \! \in \!
\mathbb{C} \setminus \widetilde{\Sigma}_{p}^{e};$ {\rm (ii)} $\mathscr{R}^{
e}_{\pm}(z) \! := \! \lim_{\underset{z^{\prime} \, \in \, \pm \, \mathrm{side}
\, \mathrm{of} \, \widetilde{\Sigma}_{p}^{e}}{z^{\prime} \to z}} \! \mathscr{
R}^{e}(z^{\prime})$ satisfy the boundary condition
\begin{equation*}
\mathscr{R}_{+}^{e}(z) \! = \! \mathscr{R}_{-}^{e}(z) \upsilon_{\mathscr{R}}^{
e}(z), \quad z \! \in \! \widetilde{\Sigma}_{p}^{e},
\end{equation*}
where $\upsilon_{\mathscr{R}}^{e}(z)$, for $z \! \in \! \widetilde{\Sigma}_{
p}^{e}$, is defined in Lemma~{\rm 4.8} and satisfies the asymptotic (as $n \!
\to \! \infty)$ estimates given in Proposition~{\rm 5.1;} {\rm (iii)}
$\mathscr{R}^{e}(z) \! =_{\underset{z \in \mathbb{C} \setminus \widetilde{
\Sigma}_{p}^{e}}{z \to \infty}} \! \mathrm{I} \! + \! \mathcal{O}(z^{-1});$
and {\rm (iv)} $\mathscr{R}^{e}(z) \! =_{\underset{z \in \mathbb{C} \setminus
\widetilde{\Sigma}_{p}^{e}}{z \to 0}} \! \mathcal{O}(1)$. Then:
\begin{compactenum}
\item[{\rm (1)}] for $z \! \in \! (-\infty,b_{0}^{e} \! - \! \delta_{b_{0}}^{
e}) \cup (a_{N+1}^{e} \! + \! \delta_{a_{N+1}}^{e},+\infty) \! =: \! \Sigma_{
p}^{e,1}$,
\begin{equation*}
\norm{w_{+}^{\Sigma_{\mathscr{R}}^{e}}(\cdot)}_{\mathcal{L}^{q}_{\mathrm{M}_{2}
(\mathbb{C})}(\Sigma_{p}^{e,1})} \! \underset{n \to \infty}{=} \! \mathcal{O}
\! \left(\dfrac{f(n) \me^{-nc}}{n^{1/q}}\right), \quad q \! = \! 1,2, \qquad
\norm{w_{+}^{\Sigma_{\mathscr{R}}^{e}}(\cdot)}_{\mathcal{L}^{\infty}_{\mathrm{
M}_{2}(\mathbb{C})}(\Sigma_{p}^{e,1})} \! \underset{n \to \infty}{=} \!
\mathcal{O} \! \left(f(n) \me^{-nc}\right),
\end{equation*}
where $c \! > \! 0$ and $f(n) \! =_{n \to \infty} \! \mathcal{O}(1);$
\item[{\rm (2)}] for $z \! \in \! (a_{j}^{e} \! + \! \delta_{a_{j}}^{e},b_{
j}^{e} \! - \! \delta_{b_{j}}^{e}) \! =: \! \Sigma_{p,j}^{e,2} \! \subset \!
\cup_{l=1}^{N} \Sigma_{p,l}^{e,2} \! =: \Sigma_{p}^{e,2}$, $j \! = \! 1,
\dotsc,N$,
\begin{equation*}
\norm{w_{+}^{\Sigma_{\mathscr{R}}^{e}}(\cdot)}_{\mathcal{L}^{q}_{\mathrm{M}_{2}
(\mathbb{C})}(\Sigma_{p,j}^{e,2})} \! \underset{n \to \infty}{=} \! \mathcal{
O} \! \left(\dfrac{f_{j}(n) \me^{-nc_{j}}}{n^{1/q}}\right), \quad q \! = \!
1,2, \quad \norm{w_{+}^{\Sigma_{\mathscr{R}}^{e}}(\cdot)}_{\mathcal{L}^{
\infty}_{\mathrm{M}_{2}(\mathbb{C})}(\Sigma_{p,j}^{e,2})} \! \underset{n \to
\infty}{=} \! \mathcal{O} \! \left(f_{j}(n) \me^{-nc_{j}}\right),
\end{equation*}
where $c_{j} \! > \! 0$ and $f_{j}(n) \! =_{n \to \infty} \! \mathcal{O}(1);$
\item[{\rm (3)}] for $z \! \in \! \cup_{j=1}^{N+1}(J_{j}^{e,\smallfrown}
\setminus (\mathbb{C}_{+} \cap (\mathbb{U}_{\delta_{b_{j-1}}}^{e} \cup
\mathbb{U}_{\delta_{a_{j}}}^{e}))) \! =: \! \Sigma_{p}^{e,3}$,
\begin{equation*}
\norm{w_{+}^{\Sigma_{\mathscr{R}}^{e}}(\cdot)}_{\mathcal{L}^{q}_{\mathrm{M}_{2}
(\mathbb{C})}(\Sigma_{p}^{e,3})} \! \underset{n \to \infty}{=} \! \mathcal{O}
\! \left(\dfrac{f(n) \me^{-nc}}{n^{1/q}}\right), \quad q \! = \! 1,2, \qquad
\norm{w_{+}^{\Sigma_{\mathscr{R}}^{e}}(\cdot)}_{\mathcal{L}^{\infty}_{\mathrm{
M}_{2}(\mathbb{C})}(\Sigma_{p}^{e,3})} \! \underset{n \to \infty}{=} \!
\mathcal{O} \! \left(f(n) \me^{-nc}\right),
\end{equation*}
where $c \! > \! 0$ and $f(n) \! =_{n \to \infty} \! \mathcal{O}(1);$
\item[{\rm (4)}] for $z \! \in \! \cup_{j=1}^{N+1}(J_{j}^{e,\smallsmile}
\setminus (\mathbb{C}_{-} \cap (\mathbb{U}_{\delta_{b_{j-1}}}^{e} \cup
\mathbb{U}_{\delta_{a_{j}}}^{e}))) \! =: \! \Sigma_{p}^{e,4}$,
\begin{equation*}
\norm{w_{+}^{\Sigma_{\mathscr{R}}^{e}}(\cdot)}_{\mathcal{L}^{q}_{\mathrm{M}_{2}
(\mathbb{C})}(\Sigma_{p}^{e,4})} \! \underset{n \to \infty}{=} \! \mathcal{O}
\! \left(\dfrac{f(n) \me^{-nc}}{n^{1/q}}\right), \quad q \! = \! 1,2, \qquad
\norm{w_{+}^{\Sigma_{\mathscr{R}}^{e}}(\cdot)}_{\mathcal{L}^{\infty}_{\mathrm{
M}_{2}(\mathbb{C})}(\Sigma_{p}^{e,4})} \! \underset{n \to \infty}{=} \!
\mathcal{O} \! \left(f(n) \me^{-nc}\right),
\end{equation*}
where $c \! > \! 0$ and $f(n) \! =_{n \to \infty} \! \mathcal{O}(1);$ and
\item[{\rm (5)}] for $z \! \in \! \cup_{j=1}^{N+1}(\partial \mathbb{U}_{
\delta_{b_{j-1}}}^{e} \cup \partial \mathbb{U}_{\delta_{a_{j}}}^{e}) \! =: \!
\Sigma_{p}^{e,5}$,
\begin{equation*}
\norm{w_{+}^{\Sigma_{\mathscr{R}}^{e}}(\cdot)}_{\mathcal{L}^{q}_{\mathrm{M}_{2}
(\mathbb{C})}(\Sigma_{p}^{e,5})} \! \underset{n \to \infty}{=} \! \mathcal{O}
\! \left(n^{-1}f(n)\right), \quad q \! \in \! \lbrace 1,2,\infty \rbrace,
\end{equation*}
where $f(n) \! =_{n \to \infty} \! \mathcal{O}(1)$.
\end{compactenum}
Furthermore,
\begin{equation*}
\norm{C^{\infty}_{w^{\Sigma_{\mathscr{R}}^{e}}}}_{\mathscr{N}_{r}(\widetilde{
\Sigma}_{p}^{e})} \! \underset{n \to \infty}{=} \! \mathcal{O} \! \left(n^{-1}
f(n)\right), \quad r \! \in \! \lbrace 2,\infty \rbrace,
\end{equation*}
where $f(n) \! =_{n \to \infty} \! \mathcal{O}(1);$ in particular, $(\id \! -
\! C^{\infty}_{w^{\Sigma_{\mathscr{R}}^{e}}})^{-1} \! \! \upharpoonright_{
\mathcal{L}^{2}_{\mathrm{M}_{2}(\mathbb{C})}(\widetilde{\Sigma}_{p}^{e})}$
exists, that is,
\begin{equation*}
\norm{(\id \! - \! C^{\infty}_{w^{\Sigma_{\mathscr{R}}^{e}}})^{-1}}_{\mathscr{
N}_{2}(\widetilde{\Sigma}_{p}^{e})} \! \underset{n \to \infty}{=} \! \mathcal{
O}(1),
\end{equation*}
and it can be expanded in a Neumann series.
\end{bbbbb}

\emph{Proof.} Without loss of generality, assume that $0 \! \in \! \Sigma^{e,
1}_{p}$ (cf. Proposition~5.1). Recall that $w_{+}^{\Sigma_{\mathscr{R}}^{e}}
(z) \! = \! \upsilon_{\mathscr{R}}^{e}(z) \! - \! \mathrm{I}$, $z \! \in \!
\widetilde{\Sigma}_{p}^{e}$. For $z \! \in \! \Sigma_{p}^{e,1}$, using the
asymptotic (as $n \! \to \! \infty)$ estimate for $\upsilon_{\mathscr{R}}^{e}
(z)$ given in item~(1) of Proposition~5.1, one gets that
\begin{align*}
\norm{w_{+}^{\Sigma_{\mathscr{R}}^{e}}(\cdot)}_{\mathcal{L}^{\infty}_{\mathrm{
M}_{2}(\mathbb{C})}(\Sigma_{p}^{e,1})} &:= \max_{i,j=1,2} \, \, \sup_{z \in
\Sigma_{p}^{e,1}} \vert (w_{+}^{\Sigma_{\mathscr{R}}^{e}}(z))_{ij} \vert \!
\underset{n \to \infty}{=} \! \mathcal{O} \! \left(f(n) \me^{-nc}\right), \\
\norm{w_{+}^{\Sigma_{\mathscr{R}}^{e}}(\cdot)}_{\mathcal{L}^{1}_{\mathrm{M}_{
2}(\mathbb{C})}(\Sigma_{p}^{e,1})} &:= \int_{\Sigma_{p}^{e,1}} \vert w_{+}^{
\Sigma_{\mathscr{R}}^{e}}(z) \vert \, \vert \md z \vert \! = \! \int_{(\Sigma_{
p}^{e,1} \setminus \mathbb{U}_{0}^{e}) \cup \mathbb{U}_{0}^{e}} \vert w_{+}^{
\Sigma_{\mathscr{R}}^{e}}(z) \vert \, \vert \md z \vert \\
&= \, \left(\int_{\Sigma_{p}^{e,1} \setminus \mathbb{U}_{0}^{e}} \! + \!
\int_{\mathbb{U}_{0}^{e}} \right) \! \! \left(\sum_{i,j=1}^{2}(w_{+}^{\Sigma_{
\mathscr{R}}^{e}}(z))_{ij} \overline{(w_{+}^{\Sigma_{\mathscr{R}}^{e}}(z))_{
ij}}\right)^{1/2} \, \vert \md z \vert \\
& \, \underset{n \to \infty}{=} \! \mathcal{O}(n^{-1}f(n) \me^{-nc}) \! + \!
\mathcal{O}(n^{-1}f(n) \me^{-nc}) \! \underset{n \to \infty}{=} \! \mathcal{O}
(n^{-1}f(n) \me^{-nc})
\end{align*}
$(\vert \md z \vert$ denotes arc length), and
\begin{align*}
\norm{w_{+}^{\Sigma_{\mathscr{R}}^{e}}(\cdot)}_{\mathcal{L}^{2}_{\mathrm{M}_{
2}(\mathbb{C})}(\Sigma_{p}^{e,1})} &:= \left(\int_{\Sigma_{p}^{e,1}} \vert
w_{+}^{\Sigma_{\mathscr{R}}^{e}}(z) \vert^{2} \, \vert \md z \vert \right)^{
1/2} \! = \! \left(\int_{(\Sigma_{p}^{e,1} \setminus \mathbb{U}_{0}^{e}) \cup
\mathbb{U}_{0}^{e}} \vert w_{+}^{\Sigma_{\mathscr{R}}^{e}}(z) \vert^{2} \,
\vert \md z \vert \right)^{1/2} \\
&= \, \left( \! \left(\int_{\Sigma_{p}^{e,1} \setminus \mathbb{U}_{0}^{e}} \!
+ \! \int_{\mathbb{U}_{0}^{e}} \right) \! \! \left(\sum_{i,j=1}^{2}(w_{+}^{
\Sigma_{\mathscr{R}}^{e}}(z))_{ij} \overline{(w_{+}^{\Sigma_{\mathscr{R}}^{e}}
(z))_{ij}} \right) \vert \md z \vert \right)^{1/2} \\
& \, \underset{n \to \infty}{=} \! \left(\mathcal{O}(n^{-1}f(n) \me^{-nc}) \!
+ \! \mathcal{O}(n^{-1}f(n) \me^{-nc}) \right)^{1/2} \! \underset{n \to
\infty}{=} \! \mathcal{O}(n^{-1/2}f(n) \me^{-nc}),
\end{align*}
where $c \! > \! 0$ and $f(n) \! =_{n \to \infty} \! \mathcal{O}(1)$.

For $z \! \in \! (a_{j}^{e} \! + \! \delta_{a_{j}}^{e},b_{j}^{e} \! - \!
\delta_{b_{j}}^{e}) \! =: \! \Sigma_{p,j}^{e,2} \! \subset \! \cup_{l=1}^{N}
\Sigma_{p,l}^{e,2} \! =: \! \Sigma_{p}^{e,2}$ $(\subset \! \widetilde{\Sigma}_{
p}^{e})$, $j \! = \! 1,\dotsc,N$, using the asymptotic (as $n \! \to \!
\infty)$ estimate for $\upsilon_{\mathscr{R}}^{e}(z)$ given in item~(2) of
Proposition~5.1, one gets that
\begin{align*}
\norm{w_{+}^{\Sigma_{\mathscr{R}}^{e}}(\cdot)}_{\mathcal{L}^{\infty}_{\mathrm{
M}_{2}(\mathbb{C})}(\Sigma_{p,j}^{e,2})} &:= \max_{k,m=1,2} \, \, \sup_{z \in
\Sigma_{p,j}^{e,2}} \vert (w_{+}^{\Sigma_{\mathscr{R}}^{e}}(z))_{km} \vert \!
\underset{n \to \infty}{=} \! \mathcal{O} \! \left(f_{j}(n) \me^{-nc_{j}}
\right), \\
\norm{w_{+}^{\Sigma_{\mathscr{R}}^{e}}(\cdot)}_{\mathcal{L}^{1}_{\mathrm{M}_{
2}(\mathbb{C})}(\Sigma_{p,j}^{e,2})} &:= \int_{\Sigma_{p,j}^{e,2}} \vert w_{
+}^{\Sigma_{\mathscr{R}}^{e}}(z) \vert \, \vert \md z \vert \! = \! \int_{
\Sigma_{p,j}^{e,2}} \! \left(\sum_{i,j=1}^{2}(w_{+}^{\Sigma_{\mathscr{R}}^{e}}
(z))_{ij} \overline{(w_{+}^{\Sigma_{\mathscr{R}}^{e}}(z))_{ij}} \right)^{1/2}
\, \vert \md z \vert \\
& \, \underset{n \to \infty}{=} \! \mathcal{O}(n^{-1}f_{j}(n) \me^{-nc_{j}}),
\end{align*}
and
\begin{align*}
\norm{w_{+}^{\Sigma_{\mathscr{R}}^{e}}(\cdot)}_{\mathcal{L}^{2}_{\mathrm{M}_{
2}(\mathbb{C})}(\Sigma_{p,j}^{e,2})} &:= \left(\int_{\Sigma_{p,j}^{e,2}} \vert
w_{+}^{\Sigma_{\mathscr{R}}^{e}}(z) \vert^{2} \, \vert \md z \vert \right)^{
1/2} \! = \! \left(\int_{\Sigma_{p,j}^{e,2}} \sum_{k,l=1}^{2}(w_{+}^{\Sigma_{
\mathscr{R}}^{e}}(z))_{kl} \overline{(w_{+}^{\Sigma_{\mathscr{R}}^{e}}(z))_{k
l}} \, \vert \md z \vert \right)^{1/2} \\
& \, \underset{n \to \infty}{=} \! \mathcal{O}(n^{-1/2}f_{j}(n) \me^{-nc_{j}}),
\end{align*}
where $c_{j} \! > \! 0$ and $f_{j}(n) \! =_{n \to \infty} \! \mathcal{O}(1)$,
$j \! = \! 1,\dotsc,N$.

For $z \! \in \! \cup_{j=1}^{N+1}(J_{j}^{e,\smallfrown} \setminus (\mathbb{C}_{
+} \cap (\mathbb{U}_{\delta_{b_{j-1}}}^{e} \cup \mathbb{U}_{\delta_{a_{j}}}^{
e}))) \! =: \! \Sigma_{p}^{e,3}$ $(\subset \! \widetilde{\Sigma}_{p}^{e})$,
using the asymptotic (as $n \! \to \! \infty)$ estimate for $\upsilon_{
\mathscr{R}}^{e}(z)$ given in item~(3) of Proposition~5.1, one gets that
\begin{align*}
\norm{w_{+}^{\Sigma_{\mathscr{R}}^{e}}(\cdot)}_{\mathcal{L}^{\infty}_{\mathrm{
M}_{2}(\mathbb{C})}(\Sigma_{p}^{e,3})} &:= \max_{i,j=1,2} \, \, \sup_{z \in
\Sigma_{p}^{e,3}} \vert (w_{+}^{\Sigma_{\mathscr{R}}^{e}}(z))_{ij} \vert \!
\underset{n \to \infty}{=} \! \mathcal{O} \! \left(f(n) \me^{-nc} \right), \\
\norm{w_{+}^{\Sigma_{\mathscr{R}}^{e}}(\cdot)}_{\mathcal{L}^{1}_{\mathrm{M}_{
2}(\mathbb{C})}(\Sigma_{p}^{e,3})} &:= \int_{\Sigma_{p}^{e,3}} \vert w_{+}^{
\Sigma_{\mathscr{R}}^{e}}(z) \vert \, \vert \md z \vert \! = \! \int_{\Sigma_{
p}^{e,3}} \! \left(\sum_{i,j=1}^{2}(w_{+}^{\Sigma_{\mathscr{R}}^{e}}(z))_{ij}
\overline{(w_{+}^{\Sigma_{\mathscr{R}}^{e}}(z))_{ij}} \right)^{1/2} \, \vert
\md z \vert \\
& \, \underset{n \to \infty}{=} \! \mathcal{O}(n^{-1}f(n) \me^{-nc}),
\end{align*}
and
\begin{align*}
\norm{w_{+}^{\Sigma_{\mathscr{R}}^{e}}(\cdot)}_{\mathcal{L}^{2}_{\mathrm{M}_{
2}(\mathbb{C})}(\Sigma_{p}^{e,3})} &:= \left(\int_{\Sigma_{p}^{e,3}} \vert w_{
+}^{\Sigma_{\mathscr{R}}^{e}}(z) \vert^{2} \, \vert \md z \vert \right)^{1/2}
\! = \! \left(\int_{\Sigma_{p}^{e,3}} \sum_{i,j=1}^{2}(w_{+}^{\Sigma_{
\mathscr{R}}^{e}}(z))_{ij} \overline{(w_{+}^{\Sigma_{\mathscr{R}}^{e}}(z))_{
ij}} \, \vert \md z \vert \right)^{1/2} \\
& \, \underset{n \to \infty}{=} \! \mathcal{O}(n^{-1/2}f(n) \me^{-nc}),
\end{align*}
where $c \! > \! 0$ and $f(n) \! =_{n \to \infty} \! \mathcal{O}(1)$: the
above analysis applies, \emph{mutatis mutandis}, for the analogous estimates
on $\Sigma_{p}^{e,4} \! := \! \cup_{j=1}^{N+1}(J_{j}^{e,\smallsmile} \setminus
(\mathbb{C}_{-} \cap (\mathbb{U}_{\delta_{b_{j-1}}}^{e} \cup \mathbb{U}_{
\delta_{a_{j}}}^{e}))) \! \ni \! z$.

For $z \! \in \! \cup_{j=1}^{N+1}(\partial \mathbb{U}_{\delta_{b_{j-1}}}^{e}
\cup \partial \mathbb{U}_{\delta_{a_{j}}}^{e}) \! =: \! \Sigma_{p}^{e,5}$
$(\subset \! \widetilde{\Sigma}_{p}^{e})$, using the $(2(N \! + \! 1))$
asymptotic (as $n \! \to \! \infty)$ estimates for $\upsilon_{\mathscr{R}}^{
e}(z)$ given in item~(5) of Proposition~5.1, one gets that
\begin{align*}
\norm{w_{+}^{\Sigma_{\mathscr{R}}^{e}}(\cdot)}_{\mathcal{L}^{\infty}_{\mathrm{
M}_{2}(\mathbb{C})}(\Sigma_{p}^{e,5})} &:= \max_{i,j=1,2} \, \, \sup_{z \in
\Sigma_{p}^{e,5}} \vert (w_{+}^{\Sigma_{\mathscr{R}}^{e}}(z))_{ij} \vert \!
\underset{n \to \infty}{=} \! \mathcal{O} \! \left(n^{-1}f(n) \right), \\
\norm{w_{+}^{\Sigma_{\mathscr{R}}^{e}}(\cdot)}_{\mathcal{L}^{1}_{\mathrm{M}_{
2}(\mathbb{C})}(\Sigma_{p}^{e,5})} &:= \int_{\Sigma_{p}^{e,5}} \vert w_{+}^{
\Sigma_{\mathscr{R}}^{e}}(z) \vert \, \vert \md z \vert \! = \! \int_{\cup_{j=
1}^{N+1}(\partial \mathbb{U}_{\delta_{b_{j-1}}}^{e} \cup \partial \mathbb{U}_{
\delta_{a_{j}}}^{e})} \vert w_{+}^{\Sigma_{\mathscr{R}}^{e}}(z) \vert \, \vert
\md z \vert \\
&= \, \sum_{k=1}^{N+1} \! \left(\int_{\partial \mathbb{U}_{\delta_{b_{k-1}}}^{
e}} \! + \! \int_{\partial \mathbb{U}_{\delta_{a_{k}}}^{e}} \right) \! \!
\left(\sum_{i,j=1}^{2}(w_{+}^{\Sigma_{\mathscr{R}}^{e}}(z))_{ij} \overline{
(w_{+}^{\Sigma_{\mathscr{R}}^{e}}(z))_{ij}}\right)^{1/2} \, \vert \md z \vert,
\end{align*}
whence, (cf. Lemma~4.5) using the fact that the respective factors $\gamma^{e}
(z) \! \pm \! (\gamma^{e}(z))^{-1}$ and $\boldsymbol{\theta}^{e}(\pm
\boldsymbol{u}^{e}(z) \! - \! \tfrac{n}{2 \pi} \boldsymbol{\Omega}^{e} \! \pm
\! \boldsymbol{d}_{e})$ are uniformly bounded (with respect to $z)$ in compact
intervals outside open intervals surrounding the end-points of the support of
the `even' equilibrium measure, one arrives at
\begin{align*}
\norm{w_{+}^{\Sigma_{\mathscr{R}}^{e}}(\cdot)}_{\mathcal{L}^{1}_{\mathrm{M}_{
2}(\mathbb{C})}(\Sigma_{p}^{e,5})} &\underset{n \to \infty}{=} \dfrac{1}{n}
\sum_{k=1}^{N+1} \! \left(\int_{\partial \mathbb{U}_{\delta_{b_{k-1}}}^{e}}
\dfrac{\vert \star_{b_{k-1}}^{e}(z;n) \vert}{(z \! - \! b_{k-1}^{e})^{3/2}} \,
\vert \md z \vert \! + \! \int_{\partial \mathbb{U}_{\delta_{a_{k}}}^{e}}
\dfrac{\vert \star_{a_{k}}^{e}(z;n) \vert}{(z \! - \! a_{k}^{e})^{3/2}} \,
\vert \md z \vert \right) \\
&\underset{n \to \infty}{=} \, \mathcal{O}(n^{-1}f(n)),
\end{align*}
and, similarly,
\begin{align*}
\norm{w_{+}^{\Sigma_{\mathscr{R}}^{e}}(\cdot)}_{\mathcal{L}^{2}_{\mathrm{M}_{
2}(\mathbb{C})}(\Sigma_{p}^{e,5})} &:= \left(\int_{\Sigma_{p}^{e,5}} \vert w_{
+}^{\Sigma_{\mathscr{R}}^{e}}(z) \vert^{2} \, \vert \md z \vert \right)^{1/2}
\! = \! \left(\int_{\cup_{j=1}^{N+1}(\partial \mathbb{U}_{\delta_{b_{j-1}}}^{
e} \cup \partial \mathbb{U}_{\delta_{a_{j}}}^{e})} \vert w_{+}^{\Sigma_{
\mathscr{R}}^{e}}(z) \vert^{2} \, \vert \md z \vert \right)^{1/2} \\
&= \, \left(\sum_{k=1}^{N+1} \! \left(\int_{\partial \mathbb{U}_{\delta_{b_{k
-1}}}^{e}} \! + \! \int_{\partial \mathbb{U}_{\delta_{a_{k}}}^{e}} \right) \!
\sum_{i,j=1}^{2}(w_{+}^{\Sigma_{\mathscr{R}}^{e}}(z))_{ij} \overline{(w_{+}^{
\Sigma_{\mathscr{R}}^{e}}(z))_{ij}} \, \vert \md z \vert \right)^{1/2} \\
&\underset{n \to \infty}{=} \, \left(\dfrac{1}{n^{2}} \sum_{k=1}^{N+1} \!
\left(\int_{\partial \mathbb{U}_{\delta_{b_{k-1}}}^{e}} \dfrac{\vert \star_{
b_{k-1}}^{e}(z;n) \vert}{(z \! - \! b_{k-1}^{e})^{3}} \, \vert \md z \vert \!
+ \! \int_{\partial \mathbb{U}_{\delta_{a_{k}}}^{e}} \dfrac{\vert \star_{a_{
k}}^{e}(z;n) \vert}{(z \! - \! a_{k}^{e})^{3}} \, \vert \md z \vert \right)
\right)^{1/2} \\
&\underset{n \to \infty}{=} \, \mathcal{O}(n^{-1}f(n)),
\end{align*}
where $f(n) \! =_{n \to \infty} \! \mathcal{O}(1)$.

Recall that $C^{\infty}_{w^{\Sigma_{\mathscr{R}}^{e}}}f \! := \! C^{\infty}_{-}
(fw^{\Sigma_{\mathscr{R}}^{e}}_{+})$, where $(C^{\infty}_{-}g)(z) \! := \!
\lim_{\underset{z^{\prime} \in -\widetilde{\Sigma}_{p}^{e}}{z^{\prime} \to z}}
\! \int_{\widetilde{\Sigma}_{p}^{e}} \tfrac{g(s)}{s-z^{\prime}} \, \tfrac{\md
s}{2 \pi \mi}$, with $-\widetilde{\Sigma}_{p}^{e}$ shorthand for `the $-$ side
of $\widetilde{\Sigma}_{p}^{e}$'. For the $\norm{C^{\infty}_{w^{\Sigma_{
\mathscr{R}}^{e}}}}_{\mathscr{N}_{\infty}(\widetilde{\Sigma}_{p}^{e})}$ norm,
one proceeds as follows:
\begin{align*}
\norm{C^{\infty}_{w^{\Sigma_{\mathscr{R}}^{e}}}g}_{\mathcal{L}^{\infty}_{
\mathrm{M}_{2}(\mathbb{C})}(\widetilde{\Sigma}_{p}^{e})} &:= \max_{j,l=1,2}
\, \, \sup_{z \in \widetilde{\Sigma}_{p}^{e}} \vert (C^{\infty}_{w^{\Sigma_{
\mathscr{R}}^{e}}}g)_{jl}(z) \vert \! = \! \max_{j,l=1,2} \, \, \sup_{z \in
\widetilde{\Sigma}_{p}^{e}} \! \left\vert \lim_{\underset{z^{\prime} \in
-\widetilde{\Sigma}_{p}^{e}}{z^{\prime} \to z}} \int_{\widetilde{\Sigma}_{p}^{
e}} \dfrac{(g(s)w^{\Sigma_{\mathscr{R}}^{e}}_{+}(s))_{jl}}{s \! - \! z^{\prime}
} \, \dfrac{\md s}{2 \pi \mi} \right\vert \\
&\leqslant \, \norm{g(\cdot)}_{\mathcal{L}^{\infty}_{\mathrm{M}_{2}(\mathbb{C}
)}(\widetilde{\Sigma}_{p}^{e})} \max_{j,l=1,2} \, \, \sup_{z \in \widetilde{
\Sigma}_{p}^{e}} \! \left\vert \lim_{\underset{z^{\prime} \in -\widetilde{
\Sigma}_{p}^{e}}{z^{\prime} \to z}} \int_{\widetilde{\Sigma}_{p}^{e}} \dfrac{
(w^{\Sigma_{\mathscr{R}}^{e}}_{+}(s))_{jl}}{s \! - \! z^{\prime}} \, \dfrac{
\md s}{2 \pi \mi} \right\vert \\
&\leqslant \, \norm{g(\cdot)}_{\mathcal{L}^{\infty}_{\mathrm{M}_{2}(\mathbb{C}
)}(\widetilde{\Sigma}_{p}^{e})} \max_{j,l=1,2} \, \, \sup_{z \in \widetilde{
\Sigma}_{p}^{e}} \! \left\vert \lim_{\underset{z^{\prime} \in -\widetilde{
\Sigma}_{p}^{e}}{z^{\prime} \to z}} \! \left(\int_{\Sigma_{p}^{e,1} \setminus
\mathbb{U}_{0}^{e}} \! + \! \int_{\mathbb{U}_{0}^{e}} \! + \! \sum_{k=1}^{N+1}
\int_{\Sigma_{p,k}^{e,2}} \! + \! \int_{\Sigma_{p}^{e,3}} \right. \right. \\
&+ \left. \left. \, \int_{\Sigma_{p}^{e,4}} \! + \! \sum_{k=1}^{N+1} \! \left(
\int_{\partial \mathbb{U}_{\delta_{b_{k-1}}}^{e}} \! + \! \int_{\partial
\mathbb{U}_{a_{k}}^{e}} \right) \right) \! \dfrac{(w^{\Sigma_{\mathscr{R}}^{
e}}_{+}(s))_{jl}}{s \! - \! z^{\prime}} \, \dfrac{\md s}{2 \pi \mi}
\right\vert \\
&\underset{n \to \infty}{\leqslant} \norm{g(\cdot)}_{\mathcal{L}^{\infty}_{
\mathrm{M}_{2}(\mathbb{C})}(\widetilde{\Sigma}_{p}^{e})} \max_{j,l=1,2} \,
\, \sup_{z \in \widetilde{\Sigma}_{p}^{e}} \! \left\vert \lim_{\underset{z^{
\prime} \in -\widetilde{\Sigma}_{p}^{e}}{z^{\prime} \to z}} \! \left(\int_{
\Sigma_{p}^{e,1} \setminus \mathbb{U}_{0}^{e}} \dfrac{(\mathcal{O}(\me^{-n
c_{\infty} \vert s \vert}f_{\infty}(n)))_{jl}}{s \! - \! z^{\prime}} \,
\dfrac{\md s}{2 \pi \mi} \right. \right. \\
&+ \left. \left. \, \int_{\mathbb{U}_{0}^{e}} \dfrac{(\mathcal{O}(\me^{-nc_{0}
\vert s \vert^{-1}}f_{0}(n)))_{jl}}{s \! - \! z^{\prime}} \, \dfrac{\md s}{2
\pi \mi} \! + \! \sum_{k=1}^{N} \int_{\Sigma_{p,k}^{e,2}} \dfrac{(\mathcal{O}
(\me^{-nc_{k}(s-a_{k}^{e})}f_{k}(n)))_{jl}}{s \! - \! z^{\prime}} \, \dfrac{
\md s}{2 \pi \mi} \right. \right. \\
&+ \left. \left. \, \int_{\Sigma_{p}^{e,3}} \dfrac{(\mathcal{O}(\me^{-n
\overset{\smallfrown}{c} \vert s \vert} \overset{\smallfrown}{f}(n)))_{jl}}{s
\! - \! z^{\prime}} \, \dfrac{\md s}{2 \pi \mi} \! + \! \int_{\Sigma_{p}^{e,4}
} \dfrac{(\mathcal{O}(\me^{-n \overset{\smallsmile}{c} \vert s \vert} \overset{
\smallsmile}{f}(n)))_{jl}}{s \! - \! z^{\prime}} \, \dfrac{\md s}{2 \pi \mi}
\right. \right. \\
&+ \left. \left. \, \sum_{k=1}^{N+1} \! \left(\int_{\partial \mathbb{U}_{
\delta_{b_{k-1}}}^{e}} \! \left(\! \mathcal{O} \! \left(\dfrac{
\overset{e}{\mathfrak{M}}^{\raise-1.0ex\hbox{$\scriptstyle \infty$}}(s)
\left(
\begin{smallmatrix}
\ast & \ast \\
\ast & \ast
\end{smallmatrix}
\right)
(\overset{e}{\mathfrak{M}}^{\raise-1.0ex\hbox{$\scriptstyle \infty$}}(s))^{-1}
}{n(s \! - \! z^{\prime})(s \! - \! b_{k-1}^{e})^{3/2}G_{b_{k-1}}^{e}(s)}
\right) \right)_{jl} \dfrac{\md s}{2 \pi \mi} \right. \right. \right. \\
&+ \left. \left. \left. \, \int_{\partial \mathbb{U}_{\delta_{a_{k}}}^{e}} \!
\left(\! \mathcal{O} \! \left(\dfrac{
\overset{e}{\mathfrak{M}}^{\raise-1.0ex\hbox{$\scriptstyle \infty$}}(s)
\left(
\begin{smallmatrix}
\ast & \ast \\
\ast & \ast
\end{smallmatrix}
\right)
(\overset{e}{\mathfrak{M}}^{\raise-1.0ex\hbox{$\scriptstyle \infty$}}(s))^{-1}
}{n(s \! - \! z^{\prime})(s \! - \! a_{k}^{e})^{3/2}G_{a_{k}}^{e}(s)} \right)
\right)_{jl} \dfrac{\md s}{2 \pi \mi} \right) \right) \right\vert,
\end{align*}
whence, (cf. Lemma~4.5) using the fact that the respective factors $\gamma^{e}
(z) \! \pm \! (\gamma^{e}(z))^{-1}$ and $\boldsymbol{\theta}^{e}(\pm
\boldsymbol{u}^{e}(z) \! - \! \tfrac{n}{2 \pi} \boldsymbol{\Omega}^{e} \! \pm
\! \boldsymbol{d}_{e})$ are uniformly bounded (with respect to $z)$ in compact
intervals outside open intervals surrounding the end-points of the support of
the `even' equilibrium measure, one arrives at, after a straightforward
integration argument and an application of the Maximum Length (ML) Theorem,
\begin{align*}
\norm{C^{\infty}_{w^{\Sigma_{\mathscr{R}}^{e}}}g}_{\mathcal{L}^{\infty}_{
\mathrm{M}_{2}(\mathbb{C})}(\widetilde{\Sigma}_{p}^{e})} &\underset{n \to
\infty}{\leqslant} \norm{g(\cdot)}_{\mathcal{L}^{\infty}_{\mathrm{M}_{2}
(\mathbb{C})}(\widetilde{\Sigma}_{p}^{e})} \! \left(\! \mathcal{O} \! \left(
\dfrac{f(n) \me^{-nc}}{n \operatorname{dist}(z,\widetilde{\Sigma}_{p}^{e})}
\right) \! + \! \mathcal{O} \! \left(\dfrac{f(n) \me^{-nc}}{\operatorname{dist}
(z,\widetilde{\Sigma}_{p}^{e})} \right) \right. \\
&+ \left. \mathcal{O} \! \left(\dfrac{f(n)}{n \operatorname{dist}(z,\widetilde{
\Sigma}_{p}^{e})} \right) \right) \! \underset{n \to \infty}{=} \! \norm{g
(\cdot)}_{\mathcal{L}^{\infty}_{\mathrm{M}_{2}(\mathbb{C})}(\widetilde{
\Sigma}_{p}^{e})} \mathcal{O}(n^{-1}f(n)),
\end{align*}
where $\operatorname{dist}(z,\widetilde{\Sigma}_{p}^{e}) \! := \! \inf
\left\lbrace \mathstrut \vert z \! - \! r \vert; \, r \! \in \! \widetilde{
\Sigma}_{p}^{e}, \, z \! \in \! \mathbb{C} \setminus \widetilde{\Sigma}_{p}^{
e} \right\rbrace$ $(> \! 0)$, and $f(n) \! =_{n \to \infty} \! \mathcal{O}
(1)$, whence one obtains the asymptotic (as $n \to \infty)$ estimate for
$\norm{C^{\infty}_{w^{\Sigma_{\mathscr{R}}^{e}}}}_{\mathscr{N}_{\infty}
(\widetilde{\Sigma}_{p}^{e})}$ stated in the Proposition. Similarly, for
$\norm{C^{\infty}_{w^{\Sigma_{\mathscr{R}}^{e}}}}_{\mathscr{N}_{2}(\widetilde{
\Sigma}_{p}^{e})}$:
\begin{align*}
\norm{C^{\infty}_{w^{\Sigma_{\mathscr{R}}^{e}}}g}_{\mathcal{L}^{2}_{\mathrm{
M}_{2}(\mathbb{C})}(\widetilde{\Sigma}_{p}^{e})} &:= \left(\int_{\widetilde{
\Sigma}_{p}^{e}} \vert (C^{\infty}_{w^{\Sigma_{\mathscr{R}}^{e}}}g)(z) \vert^{
2} \, \vert \md z \vert \right)^{1/2} \! = \! \left(\int_{\widetilde{\Sigma}_{
p}^{e}} \sum_{j,l=1}^{2}(C^{\infty}_{w^{\Sigma_{\mathscr{R}}^{e}}}g)_{jl}(z)
\overline{(C^{\infty}_{w^{\Sigma_{\mathscr{R}}^{e}}}g)_{jl}(z)} \, \vert \md z
\vert \right)^{1/2} \\
&= \, \left(\int_{\widetilde{\Sigma}_{p}^{e}} \sum_{j,l=1}^{2} \! \left\vert
\lim_{\underset{z^{\prime} \in -\widetilde{\Sigma}_{p}^{e}}{z^{\prime} \to z}}
\int_{\widetilde{\Sigma}_{p}^{e}} \dfrac{(g(s)w^{\Sigma_{\mathscr{R}}^{e}}_{+}
(s))_{jl}}{s \! - \! z^{\prime}} \, \dfrac{\md s}{2 \pi \mi} \right\vert^{2}
\vert \md z \vert \right)^{1/2} \\
&\leqslant \, \norm{g(\cdot)}_{\mathcal{L}^{2}_{\mathrm{M}_{2}(\mathbb{C})}
(\widetilde{\Sigma}_{p}^{e})} \! \left(\int_{\widetilde{\Sigma}_{p}^{e}}
\sum_{j,l=1}^{2} \! \left\vert \lim_{\underset{z^{\prime} \in -\widetilde{
\Sigma}_{p}^{e}}{z^{\prime} \to z}} \int_{\widetilde{\Sigma}_{p}^{e}} \dfrac{
(w^{\Sigma_{\mathscr{R}}^{e}}_{+}(s))_{jl}}{s \! - \! z^{\prime}} \, \dfrac{
\md s}{2 \pi \mi} \right\vert^{2} \vert \md z \vert \right)^{1/2} \\
&\leqslant \, \norm{g(\cdot)}_{\mathcal{L}^{2}_{\mathrm{M}_{2}(\mathbb{C})}
(\widetilde{\Sigma}_{p}^{e})} \! \left(\int_{\widetilde{\Sigma}_{p}^{e}} \sum_{
j,l=1}^{2} \! \left\vert \lim_{\underset{z^{\prime} \in -\widetilde{\Sigma}_{
p}^{e}}{z^{\prime} \to z}} \! \left(\int_{\Sigma_{p}^{e,1} \setminus \mathbb{
U}_{0}^{e}} \! + \! \int_{\mathbb{U}_{0}^{e}} \! + \! \sum_{k=1}^{N+1} \int_{
\Sigma_{p,k}^{e,2}} \! + \! \int_{\Sigma_{p}^{e,3}} \right. \right. \right. \\
&+ \left. \left. \left. \, \int_{\Sigma_{p}^{e,4}} \! + \! \sum_{k=1}^{N+1} \!
\left(\int_{\partial \mathbb{U}_{\delta_{b_{k-1}}}^{e}} \! + \! \int_{\partial
\mathbb{U}_{a_{k}}^{e}} \right) \right) \! \dfrac{(w^{\Sigma_{\mathscr{R}}^{e}
}_{+}(s))_{jl}}{s \! - \! z^{\prime}} \, \dfrac{\md s}{2 \pi \mi} \right\vert^{
2} \vert \md z \vert \right)^{1/2} \\
&\underset{n \to \infty}{\leqslant} \norm{g(\cdot)}_{\mathcal{L}^{2}_{\mathrm{
M}_{2}(\mathbb{C})}(\widetilde{\Sigma}_{p}^{e})} \! \left(\int_{\widetilde{
\Sigma}_{p}^{e}} \sum_{j,l=1}^{2} \! \left\vert \lim_{\underset{z^{\prime} \in
-\widetilde{\Sigma}_{p}^{e}}{z^{\prime} \to z}} \! \left(\int_{\Sigma_{p}^{e,
1} \setminus \mathbb{U}_{0}^{e}} \dfrac{(\mathcal{O}(\me^{-nc_{\infty} \vert s
\vert}f_{\infty}(n)))_{jl}}{s \! - \! z^{\prime}} \, \dfrac{\md s}{2 \pi \mi}
\right. \right. \right. \\
&+ \left. \left. \left. \, \int_{\mathbb{U}_{0}^{e}} \dfrac{(\mathcal{O}(\me^{
-nc_{0} \vert s \vert^{-1}}f_{0}(n)))_{jl}}{s \! - \! z^{\prime}} \, \dfrac{
\md s}{2 \pi \mi} \! + \! \sum_{k=1}^{N} \int_{\Sigma_{p,k}^{e,2}} \dfrac{(
\mathcal{O}(\me^{-nc_{k}(s-a_{k}^{e})}f_{k}(n)))_{jl}}{s \! - \! z^{\prime}}
\, \dfrac{\md s}{2 \pi \mi} \right. \right. \right. \\
&+ \left. \left. \left. \, \int_{\Sigma_{p}^{e,3}} \dfrac{(\mathcal{O}(\me^{-n
\overset{\smallfrown}{c} \vert s \vert} \overset{\smallfrown}{f}(n)))_{jl}}{s
\! - \! z^{\prime}} \, \dfrac{\md s}{2 \pi \mi} \! + \! \int_{\Sigma_{p}^{e,4}
} \dfrac{(\mathcal{O}(\me^{-n \overset{\smallsmile}{c} \vert s \vert} \overset{
\smallsmile}{f}(n)))_{jl}}{s \! - \! z^{\prime}} \, \dfrac{\md s}{2 \pi \mi}
\right. \right. \right. \\
&+ \left. \left. \left. \, \sum_{k=1}^{N+1} \! \left(\int_{\partial \mathbb{
U}_{\delta_{b_{k-1}}}^{e}} \! \left(\! \mathcal{O} \! \left(\dfrac{
\overset{e}{\mathfrak{M}}^{\raise-1.0ex\hbox{$\scriptstyle \infty$}}(s)
\left(
\begin{smallmatrix}
\ast & \ast \\
\ast & \ast
\end{smallmatrix}
\right)
(\overset{e}{\mathfrak{M}}^{\raise-1.0ex\hbox{$\scriptstyle \infty$}}(s))^{-1}
}{n(s \! - \! z^{\prime})(s \! - \! b_{k-1}^{e})^{3/2}G_{b_{k-1}}^{e}(s)}
\right) \right)_{jl} \dfrac{\md s}{2 \pi \mi} \right. \right. \right. \right.
\\
&+ \left. \left. \left. \left. \, \int_{\partial \mathbb{U}_{\delta_{a_{k}}}^{
e}} \! \left(\! \mathcal{O} \! \left(\dfrac{
\overset{e}{\mathfrak{M}}^{\raise-1.0ex\hbox{$\scriptstyle \infty$}}(s)
\left(
\begin{smallmatrix}
\ast & \ast \\
\ast & \ast
\end{smallmatrix}
\right)
(\overset{e}{\mathfrak{M}}^{\raise-1.0ex\hbox{$\scriptstyle \infty$}}(s))^{-1}
}{n(s \! - \! z^{\prime})(s \! - \! a_{k}^{e})^{3/2}G_{a_{k}}^{e}(s)} \right)
\right)_{jl} \dfrac{\md s}{2 \pi \mi} \right) \right) \right\vert^{2} \vert
\md z \vert \right)^{1/2},
\end{align*}
whence, (cf. Lemma~4.5) using the fact that the respective factors $\gamma^{e}
(z) \! \pm \! (\gamma^{e}(z))^{-1}$ and $\boldsymbol{\theta}^{e}(\pm
\boldsymbol{u}^{e}(z) \! - \! \tfrac{n}{2 \pi} \boldsymbol{\Omega}^{e} \! \pm
\! \boldsymbol{d}_{e})$ are uniformly bounded (with respect to $z)$ in compact
intervals outside open intervals surrounding the end-points of the support
of the `even' equilibrium measure, one arrives at, after a straightforward
integration argument and an application of the ML Theorem,
\begin{align*}
\norm{C^{\infty}_{w^{\Sigma_{\mathscr{R}}^{e}}}g}_{\mathcal{L}^{2}_{\mathrm{
M}_{2}(\mathbb{C})}(\widetilde{\Sigma}_{p}^{e})} &\underset{n \to \infty}{
\leqslant} \norm{g(\cdot)}_{\mathcal{L}^{2}_{\mathrm{M}_{2}(\mathbb{C})}
(\widetilde{\Sigma}_{p}^{e})} \! \left(\! \mathcal{O} \! \left(\dfrac{f(n)
\me^{-nc}}{n \operatorname{dist}(z,\widetilde{\Sigma}_{p}^{e})} \right) \! +
\! \mathcal{O} \! \left(\dfrac{f(n) \me^{-nc}}{\operatorname{dist}(z,
\widetilde{\Sigma}_{p}^{e})} \right) \right. \\
&+ \left. \mathcal{O} \! \left(\dfrac{f(n)}{n \operatorname{dist}(z,\widetilde{
\Sigma}_{p}^{e})} \right) \right) \! \underset{n \to \infty}{=} \! \norm{g
(\cdot)}_{\mathcal{L}^{2}_{\mathrm{M}_{2}(\mathbb{C})}(\widetilde{\Sigma}_{
p}^{e})} \mathcal{O}(n^{-1}f(n)),
\end{align*}
where $f(n) \! =_{n \to \infty} \! \mathcal{O}(1)$, whence one obtains the
asymptotic (as $n \to \infty)$ estimate for $\norm{C^{\infty}_{w^{\Sigma_{
\mathscr{R}}^{e}}}}_{\mathscr{N}_{2}(\widetilde{\Sigma}_{p}^{e})}$ stated in
the Proposition. The above analysis establishes the fact that, as $n \! \to \!
\infty$, $C^{\infty}_{w^{\Sigma_{\mathscr{R}}^{e}}} \! \in \! \mathscr{N}_{2}
(\widetilde{\Sigma}_{p}^{e})$, with operator norm $\norm{C^{\infty}_{w^{
\Sigma_{\mathscr{R}}^{e}}}}_{\mathscr{N}_{2}(\widetilde{\Sigma}_{p}^{e})} \!
=_{n \to \infty} \! \mathcal{O}(n^{-1}f(n))$, where $f(n) \! =_{n \to \infty}
\! \mathcal{O}(1)$; due to a well-known result for bounded linear operators
in Hilbert space \cite{a105}, it follows, thus, that $(\id \! - \! C^{\infty}_{
w^{\Sigma_{\mathscr{R}}^{e}}})^{-1} \! \! \upharpoonright_{\mathcal{L}^{2}_{
\mathrm{M}_{2}(\mathbb{C})}(\widetilde{\Sigma}_{p}^{e})}$ exists, and $(\id
\! - \! C^{\infty}_{w^{\Sigma_{\mathscr{R}}^{e}}}) \! \! \upharpoonright_{
\mathcal{L}^{2}_{\mathrm{M}_{2}(\mathbb{C})}(\widetilde{\Sigma}_{p}^{e})}$ can
be inverted by a Neumann series (as $n \! \to \! \infty)$, with $\norm{(\id
\! - \! C^{\infty}_{w^{\Sigma_{\mathscr{R}}^{e}}})^{-1}}_{\mathscr{N}_{2}
(\widetilde{\Sigma}_{p}^{e})} \! \leqslant_{n \to \infty} \! (1 \! - \!
\norm{C^{\infty}_{w^{\Sigma_{\mathscr{R}}^{e}}}}_{\mathscr{N}_{2}(\widetilde{
\Sigma}_{p}^{e})})^{-1} \! =_{n \to \infty} \! \mathcal{O}(1)$. \hfill $\qed$
\begin{ccccc}
Set $\Sigma_{\circlearrowright}^{e} \! := \! \Sigma_{p}^{e,5}$ $(= \! \cup_{j=
1}^{N+1}(\partial \mathbb{U}_{\delta_{b_{j-1}}}^{e} \cup \partial \mathbb{U}_{
\delta_{a_{j}}}^{e}))$ and $\Sigma_{\scriptscriptstyle \blacksquare}^{e} \! :=
\! \widetilde{\Sigma}_{p}^{e} \setminus \Sigma_{\circlearrowright}^{e}$, and
let $\mathscr{R}^{e} \colon \mathbb{C} \setminus \widetilde{\Sigma}_{p}^{e}
\! \to \! \operatorname{SL}_{2}(\mathbb{C})$ solve the (equivalent) {\rm RHP}
$(\mathscr{R}^{e}(z),\upsilon_{\mathscr{R}}^{e}(z),\widetilde{\Sigma}_{p}^{e}
)$ formulated in Proposition~{\rm 5.2} with integral representation given by
Equation~{\rm (5.1)}. Let the asymptotic (as $n \! \to \! \infty)$ estimates
and bounds given in Propositions~{\rm 5.1} and~{\rm 5.2} be valid. Then,
uniformly for compact subsets of $\mathbb{C} \setminus \widetilde{\Sigma}_{
p}^{e} \! \ni \! z$,
\begin{equation*}
\mathscr{R}^{e}(z) \underset{\underset{z \in \mathbb{C} \setminus \widetilde{
\Sigma}^{e}_{p}}{n \to \infty}}{=} \mathrm{I} \! + \! \int_{\Sigma_{
\circlearrowright}^{e}} \dfrac{w^{\Sigma_{\circlearrowright}^{e}}_{+}(s)}{s
\! - \! z} \, \dfrac{\md s}{2 \pi \mi} \! + \! \mathcal{O} \! \left(\dfrac{
f(n)}{n^{2} \operatorname{dist}(z,\widetilde{\Sigma}_{p}^{e})} \right),
\end{equation*}
where $w_{+}^{\Sigma_{\circlearrowright}^{e}}(z) \! := \! w_{+}^{\Sigma_{
\mathscr{R}}^{e}}(z) \! \! \upharpoonright_{\Sigma_{\circlearrowright}^{e}}$,
and $(f(n))_{ij} \! =_{n \to \infty} \! \mathcal{O}(1)$, $i,j \! = \! 1,2$.
\end{ccccc}

\emph{Proof.} Define $\Sigma_{\circlearrowright}^{e}$ and $\Sigma^{e}_{
\scriptscriptstyle \blacksquare}$ as in the Lemma, and write $\widetilde{
\Sigma}_{p}^{e} \! = \! (\widetilde{\Sigma}_{p}^{e} \setminus \Sigma^{e}_{
\circlearrowright}) \cup \Sigma_{\circlearrowright}^{e} \! := \! \Sigma^{e}_{
\scriptscriptstyle \blacksquare} \cup \Sigma_{\circlearrowright}^{e}$ (with
$\Sigma^{e}_{\scriptscriptstyle \blacksquare} \cap \Sigma_{\circlearrowright}^{
e} \! = \! \varnothing)$. Recall, {}from Equation~(5.1), the integral
representation for $\mathscr{R}^{e} \colon \mathbb{C} \setminus \widetilde{
\Sigma}_{p}^{e} \! \to \! \operatorname{SL}_{2}(\mathbb{C})$:
\begin{equation*}
\mathscr{R}^{e}(z) \! = \! \mathrm{I} \! + \! \int_{\widetilde{\Sigma}_{p}^{e}}
\dfrac{\mu^{\Sigma_{\mathscr{R}}^{e}}(s)w_{+}^{\Sigma_{\mathscr{R}}^{e}}(s)}{
s \! - \! z} \, \dfrac{\md s}{2 \pi \mi}, \quad z \! \in \! \mathbb{C}
\setminus \widetilde{\Sigma}_{p}^{e}.
\end{equation*}
Using the linearity property of the Cauchy integral operator $C^{\infty}_{w^{
\Sigma_{\mathscr{R}}^{e}}}$, one shows that $C^{\infty}_{w^{\Sigma_{\mathscr{
R}}^{e}}} \! = \! C^{\infty}_{w^{\Sigma_{\circlearrowright}^{e}}} \! + \! C^{
\infty}_{w^{\Sigma_{\scriptscriptstyle \blacksquare}^{e}}}$. Via a repeated
application of the second resolvent identity\footnote{For general operators
$\mathscr{A}$ and $\mathscr{B}$, if $(\id \! - \! \mathscr{A})^{-1}$ and $(\id
\! - \! \mathscr{B})^{-1}$ exist, then $(\id \! - \! \mathscr{B})^{-1} \! - \!
(\id \! - \! \mathscr{A})^{-1} \! = \! (\id \! - \! \mathscr{B})^{-1}(\mathscr{
B} \! - \! \mathscr{A})(\id \! - \! \mathscr{A})^{-1}$ \cite{a105}.}:
\begin{align*}
\mu^{\Sigma_{\mathscr{R}}^{e}}(z) &= \mathrm{I} \! + \! ((\id \! - \! C^{
\infty}_{w^{\Sigma_{\mathscr{R}}^{e}}})^{-1}C^{\infty}_{w^{\Sigma_{\mathscr{
R}}^{e}}} \mathrm{I})(z) \! = \! \mathrm{I} \! + \! ((\id \! - \! C^{\infty}_{
w^{\Sigma_{\circlearrowright}^{e}}} \! - \! C^{\infty}_{w^{\Sigma_{
\scriptscriptstyle \blacksquare}^{e}}})^{-1}(C^{\infty}_{w^{\Sigma_{
\circlearrowright}^{e}}} \! + \! C^{\infty}_{w^{\Sigma_{\scriptscriptstyle
\blacksquare}^{e}}}) \mathrm{I})(z) \\
&=\mathrm{I} \! + \! ((\id \! - \! C^{\infty}_{w^{\Sigma_{\circlearrowright}^{
e}}} \! - \! C^{\infty}_{w^{\Sigma_{\scriptscriptstyle \blacksquare}^{e}}})^{-
1}C^{\infty}_{w^{\Sigma_{\circlearrowright}^{e}}} \mathrm{I})(z) \! + \! ((\id
\! - \! C^{\infty}_{w^{\Sigma_{\circlearrowright}^{e}}} \! - \! C^{\infty}_{w^{
\Sigma_{\scriptscriptstyle \blacksquare}^{e}}})^{-1}C^{\infty}_{w^{\Sigma_{
\scriptscriptstyle \blacksquare}^{e}}} \mathrm{I})(z) \\
&=\mathrm{I} \! + \! (((\id \! - \! C^{\infty}_{w^{\Sigma_{\circlearrowright}^{
e}}})(\id \! - \! (\id \! - \! C^{\infty}_{w^{\Sigma_{\circlearrowright}^{e}}}
)^{-1}C^{\infty}_{w^{\Sigma_{\scriptscriptstyle \blacksquare}^{e}}}))^{-1}C^{
\infty}_{w^{\Sigma_{\circlearrowright}^{e}}} \mathrm{I})(z) \\
&+(((\id \! - \! C^{\infty}_{w^{\Sigma_{\scriptscriptstyle \blacksquare}^{e}}})
(\id \! - \! (\id \! - \! C^{\infty}_{w^{\Sigma_{\scriptscriptstyle
\blacksquare}^{e}}})^{-1}C^{\infty}_{w^{\Sigma_{\circlearrowright}^{e}}}))^{-1}
C^{\infty}_{w^{\Sigma_{\scriptscriptstyle \blacksquare}^{e}}} \mathrm{I})(z) \\
&=\mathrm{I} \! + \! ((\id \! - \! (\id \! - \! C^{\infty}_{w^{\Sigma_{
\circlearrowright}^{e}}})^{-1}C^{\infty}_{w^{\Sigma_{\scriptscriptstyle
\blacksquare}^{e}}})^{-1}(\id \! + \! (\id \! - \! C^{\infty}_{w^{\Sigma_{
\circlearrowright}^{e}}})^{-1}C^{\infty}_{w^{\Sigma_{\circlearrowright}^{e}}})
C^{\infty}_{w^{\Sigma_{\circlearrowright}^{e}}} \mathrm{I})(z) \\
&+((\id \! - \! (\id \! - \! C^{\infty}_{w^{\Sigma_{\scriptscriptstyle
\blacksquare}^{e}}})^{-1}C^{\infty}_{w^{\Sigma_{\circlearrowright}^{e}}})^{-1}
(\id \! + \! (\id \! - \! C^{\infty}_{w^{\Sigma_{\scriptscriptstyle
\blacksquare}^{e}}})^{-1}C^{\infty}_{w^{\Sigma_{\scriptscriptstyle
\blacksquare}^{e}}})C^{\infty}_{w^{\Sigma_{\scriptscriptstyle \blacksquare}^{
e}}} \mathrm{I})(z) \\
&=\mathrm{I} \! + \! ((\id \! - \! (\id \! - \! C^{\infty}_{w^{\Sigma_{
\circlearrowright}^{e}}})^{-1}C^{\infty}_{w^{\Sigma_{\scriptscriptstyle
\blacksquare}^{e}}})^{-1}((\id \! - \! C^{\infty}_{w^{\Sigma_{
\circlearrowright}^{e}}})^{-1}C^{\infty}_{w^{\Sigma_{\circlearrowright}^{e}}})
(C^{\infty}_{w^{\Sigma_{\circlearrowright}^{e}}} \mathrm{I}))(z) \\
&+((\id \! - \! (\id \! - \! C^{\infty}_{w^{\Sigma_{\scriptscriptstyle
\blacksquare}^{e}}})^{-1}C^{\infty}_{w^{\Sigma_{\circlearrowright}^{e}}})^{-1}
((\id \! - \! C^{\infty}_{w^{\Sigma_{\scriptscriptstyle \blacksquare}^{e}}})^{
-1}C^{\infty}_{w^{\Sigma_{\scriptscriptstyle \blacksquare}^{e}}})(C^{\infty}_{
w^{\Sigma_{\scriptscriptstyle \blacksquare}^{e}}} \mathrm{I}))(z) \\
&+((\id \! - \! (\id \! - \! C^{\infty}_{w^{\Sigma_{\circlearrowright}^{e}}})^{
-1}C^{\infty}_{w^{\Sigma_{\scriptscriptstyle \blacksquare}^{e}}})^{-1}(C^{
\infty}_{w^{\Sigma_{\circlearrowright}^{e}}} \mathrm{I}))(z) \! + \! ((\id \!
- \! (\id \! - \! C^{\infty}_{w^{\Sigma_{\scriptscriptstyle \blacksquare}^{e}}}
)^{-1}C^{\infty}_{w^{\Sigma_{\circlearrowright}^{e}}})^{-1}(C^{\infty}_{w^{
\Sigma_{\scriptscriptstyle \blacksquare}^{e}}} \mathrm{I}))(z) \\
&=\mathrm{I} \! + \! ((\id \! + \! (\id \! - \! (\id \! - \! C^{\infty}_{w^{
\Sigma_{\circlearrowright}^{e}}})^{-1}C^{\infty}_{w^{\Sigma_{
\scriptscriptstyle \blacksquare}^{e}}})^{-1}(\id \! - \! C^{\infty}_{w^{
\Sigma_{\circlearrowright}^{e}}})^{-1}C^{\infty}_{w^{\Sigma_{
\scriptscriptstyle \blacksquare}^{e}}})(C^{\infty}_{w^{\Sigma_{
\circlearrowright}^{e}}} \mathrm{I}))(z) \\
&+((\id \! + \! (\id \! - \! (\id \! - \! C^{\infty}_{w^{\Sigma_{
\scriptscriptstyle \blacksquare}^{e}}})^{-1}C^{\infty}_{w^{\Sigma_{
\circlearrowright}^{e}}})^{-1}(\id \! - \! C^{\infty}_{w^{
\Sigma_{\scriptscriptstyle \blacksquare}^{e}}})^{-1}C^{\infty}_{w^{\Sigma_{
\circlearrowright}^{e}}})(C^{\infty}_{w^{\Sigma_{\scriptscriptstyle
\blacksquare}^{e}}} \mathrm{I}))(z) \\
&+((\id \! + \! (\id \! - \! (\id \! - \! C^{\infty}_{w^{\Sigma_{
\circlearrowright}^{e}}})^{-1}C^{\infty}_{w^{\Sigma_{\scriptscriptstyle
\blacksquare}^{e}}})^{-1}(\id \! - \! C^{\infty}_{w^{\Sigma_{
\circlearrowright}^{e}}})^{-1}C^{\infty}_{w^{\Sigma_{\scriptscriptstyle
\blacksquare}^{e}}})(\id \! - \! C^{\infty}_{w^{\Sigma_{\circlearrowright}^{
e}}})^{-1}C^{\infty}_{w^{\Sigma_{\circlearrowright}^{e}}}(C^{\infty}_{w^{
\Sigma_{\circlearrowright}^{e}}} \mathrm{I}))(z) \\
&+((\id \! + \! (\id \! - \! (\id \! - \! C^{\infty}_{w^{\Sigma_{
\scriptscriptstyle \blacksquare}^{e}}})^{-1}C^{\infty}_{w^{\Sigma_{
\circlearrowright}^{e}}})^{-1}(\id \! - \! C^{\infty}_{w^{\Sigma_{
\scriptscriptstyle \blacksquare}^{e}}})^{-1}C^{\infty}_{w^{\Sigma_{
\circlearrowright}^{e}}})(\id \! - \! C^{\infty}_{w^{\Sigma_{
\scriptscriptstyle \blacksquare}^{e}}})^{-1}C^{\infty}_{w^{\Sigma_{
\scriptscriptstyle \blacksquare}^{e}}}(C^{\infty}_{w^{\Sigma_{
\scriptscriptstyle \blacksquare}^{e}}} \mathrm{I}))(z) \\
&=\mathrm{I} \! + \! (C^{\infty}_{w^{\Sigma_{\circlearrowright}^{e}}} \mathrm{
I})(z) \! + \! (C^{\infty}_{w^{\Sigma_{\scriptscriptstyle \blacksquare}^{e}}}
\mathrm{I})(z) \! + \! ((\id \! - \! C^{\infty}_{w^{\Sigma_{
\circlearrowright}^{e}}})^{-1}C^{\infty}_{w^{\Sigma_{\circlearrowright}^{e}}}
(C^{\infty}_{w^{\Sigma_{\circlearrowright}^{e}}} \mathrm{I}))(z) \\
&+((\id \! - \! C^{\infty}_{w^{\Sigma_{\scriptscriptstyle \blacksquare}^{e}}}
)^{-1}C^{\infty}_{w^{\Sigma_{\scriptscriptstyle \blacksquare}^{e}}}(C^{
\infty}_{w^{\Sigma_{\scriptscriptstyle \blacksquare}^{e}}} \mathrm{I}))(z) \!
+ \! ((\id \! - \! (\id \! - \! C^{\infty}_{w^{\Sigma_{\circlearrowright}^{e}}
})^{-1}C^{\infty}_{w^{\Sigma_{\scriptscriptstyle \blacksquare}^{e}}})^{-1}(\id
\! - \! C^{\infty}_{w^{\Sigma_{\circlearrowright}^{e}}})^{-1}C^{\infty}_{w^{
\Sigma_{\scriptscriptstyle \blacksquare}^{e}}} \\
&\times (C^{\infty}_{w^{\Sigma_{\circlearrowright}^{e}}} \mathrm{I}))(z) \!
+ \! ((\id \! - \! (\id \! - \! C^{\infty}_{w^{\Sigma_{\scriptscriptstyle
\blacksquare}^{e}}})^{-1}C^{\infty}_{w^{\Sigma_{\circlearrowright}^{e}}})^{-1}
(\id \! - \! C^{\infty}_{w^{\Sigma_{\scriptscriptstyle \blacksquare}^{e}}})^{-
1}C^{\infty}_{w^{\Sigma_{\circlearrowright}^{e}}}(C^{\infty}_{w^{\Sigma_{
\scriptscriptstyle \blacksquare}^{e}}} \mathrm{I}))(z) \\
&+((\id \! - \! (\id \! - \! C^{\infty}_{w^{\Sigma_{\circlearrowright}^{e}}}
)^{-1}C^{\infty}_{w^{\Sigma_{\scriptscriptstyle \blacksquare}^{e}}})^{-1}(\id
\! - \! C^{\infty}_{w^{\Sigma_{\circlearrowright}^{e}}})^{-1}C^{\infty}_{w^{
\Sigma_{\scriptscriptstyle \blacksquare}^{e}}}(\id \! - \! C^{\infty}_{w^{
\Sigma_{\circlearrowright}^{e}}})^{-1}C^{\infty}_{w^{\Sigma_{
\circlearrowright}^{e}}}(C^{\infty}_{w^{\Sigma_{\circlearrowright}^{e}}}
\mathrm{I}))(z) \\
&+((\id \! - \! (\id \! - \! C^{\infty}_{w^{\Sigma_{\scriptscriptstyle
\blacksquare}^{e}}})^{-1}C^{\infty}_{w^{\Sigma_{\circlearrowright}^{e}}})^{-1}
(\id \! - \! C^{\infty}_{w^{\Sigma_{\scriptscriptstyle \blacksquare}^{e}}})^{-
1}C^{\infty}_{w^{\Sigma_{\circlearrowright}^{e}}}(\id \! - \! C^{\infty}_{w^{
\Sigma_{\scriptscriptstyle \blacksquare}^{e}}})^{-1}C^{\infty}_{w^{\Sigma_{
\scriptscriptstyle \blacksquare}^{e}}}(C^{\infty}_{w^{\Sigma_{
\scriptscriptstyle \blacksquare}^{e}}} \mathrm{I}))(z) \\
&=\mathrm{I} \! + \! (C^{\infty}_{w^{\Sigma_{\circlearrowright}^{e}}} \mathrm{
I})(z) \! + \! (C^{\infty}_{w^{\Sigma_{\scriptscriptstyle \blacksquare}^{e}}}
\mathrm{I})(z) \! + \! ((\id \! - \! C^{\infty}_{w^{\Sigma_{
\circlearrowright}^{e}}})^{-1}C^{\infty}_{w^{\Sigma_{\circlearrowright}^{e}}}
(C^{\infty}_{w^{\Sigma_{\circlearrowright}^{e}}} \mathrm{I}))(z) \! + \! ((\id
\! - \! C^{\infty}_{w^{\Sigma_{\scriptscriptstyle \blacksquare}^{e}}})^{-1}C^{
\infty}_{w^{\Sigma_{\scriptscriptstyle \blacksquare}^{e}}} \\
&\times \! (C^{\infty}_{w^{\Sigma_{\scriptscriptstyle \blacksquare}^{e}}}
\mathrm{I}))(z) \! + \! ((\id \! - \! (\id \! - \! C^{\infty}_{w^{\Sigma_{
\circlearrowright}^{e}}})^{-1}(\id \! - \! C^{\infty}_{w^{\Sigma_{
\scriptscriptstyle \blacksquare}^{e}}})^{-1}C^{\infty}_{w^{\Sigma_{
\scriptscriptstyle \blacksquare}^{e}}}C^{\infty}_{w^{\Sigma_{
\circlearrowright}^{e}}})^{-1}(\id \! - \! C^{\infty}_{w^{\Sigma_{
\circlearrowright}^{e}}})^{-1}(\id \! - \! C^{\infty}_{w^{\Sigma_{
\scriptscriptstyle \blacksquare}^{e}}})^{-1} \\
&\times C^{\infty}_{w^{\Sigma_{\circlearrowright}^{e}}}(C^{\infty}_{w^{\Sigma_{
\scriptscriptstyle \blacksquare}^{e}}} \mathrm{I}))(z) \! + \! ((\id \! - \!
(\id \! - \! C^{\infty}_{w^{\Sigma_{\scriptscriptstyle \blacksquare}^{e}}})^{
-1}(\id \! - \! C^{\infty}_{w^{\Sigma_{\circlearrowright}^{e}}})^{-1}C^{
\infty}_{w^{\Sigma_{\circlearrowright}^{e}}}C^{\infty}_{w^{\Sigma_{
\scriptscriptstyle \blacksquare}^{e}}})^{-1}(\id \! - \! C^{\infty}_{w^{
\Sigma_{\scriptscriptstyle \blacksquare}^{e}}})^{-1} \\
&\times (\id \! - \! C^{\infty}_{w^{\Sigma_{\circlearrowright}^{e}}})^{-1}C^{
\infty}_{w^{\Sigma_{\scriptscriptstyle \blacksquare}^{e}}}(C^{\infty}_{w^{
\Sigma_{\circlearrowright}^{e}}} \mathrm{I}))(z) \! + \! ((\id \! - \! (\id \!
- \! C^{\infty}_{w^{\Sigma_{\circlearrowright}^{e}}})^{-1}(\id \! - \! C^{
\infty}_{w^{\Sigma_{\scriptscriptstyle \blacksquare}^{e}}})^{-1}C^{\infty}_{
w^{\Sigma_{\scriptscriptstyle \blacksquare}^{e}}}C^{\infty}_{w^{\Sigma_{
\circlearrowright}^{e}}})^{-1} \\
&\times (\id \! - \! C^{\infty}_{w^{\Sigma_{\circlearrowright}^{e}}})^{-1}(\id
\! - \! C^{\infty}_{w^{\Sigma_{\scriptscriptstyle \blacksquare}^{e}}})^{-1}C^{
\infty}_{w^{\Sigma_{\circlearrowright}^{e}}}(\id \! - \! C^{\infty}_{w^{
\Sigma_{\scriptscriptstyle \blacksquare}^{e}}})^{-1}C^{\infty}_{w^{\Sigma_{
\scriptscriptstyle \blacksquare}^{e}}}(C^{\infty}_{w^{\Sigma_{
\scriptscriptstyle \blacksquare}^{e}}} \mathrm{I}))(z) \! + \! ((\id \! - \!
(\id \! - \! C^{\infty}_{w^{\Sigma_{\scriptscriptstyle \blacksquare}^{e}}})^{
-1} \\
&\times (\id \! - \! C^{\infty}_{w^{\Sigma_{\circlearrowright}^{e}}})^{-1}
C^{\infty}_{w^{\Sigma_{\circlearrowright}^{e}}}C^{\infty}_{w^{\Sigma_{
\scriptscriptstyle \blacksquare}^{e}}})^{-1}(\id \! - \! C^{\infty}_{w^{
\Sigma_{\scriptscriptstyle \blacksquare}^{e}}})^{-1}(\id \! - \! C^{\infty}_{
w^{\Sigma_{\circlearrowright}^{e}}})^{-1}C^{\infty}_{w^{\Sigma_{
\scriptscriptstyle \blacksquare}^{e}}}(\id \! - \! C^{\infty}_{w^{\Sigma_{
\circlearrowright}^{e}}})^{-1}C^{\infty}_{w^{\Sigma_{\circlearrowright}^{e}}}
\\
&\times (C^{\infty}_{w^{\Sigma_{\circlearrowright}^{e}}} \mathrm{I}))(z);
\end{align*}
hence, recalling the integral representation for $\mathscr{R}^{e}(z)$ given
above, one arrives at, for $\mathbb{C} \setminus \widetilde{\Sigma}_{p}^{e}
\! \ni \! z$,
\begin{equation}
\mathscr{R}^{e}(z) \! - \! \mathrm{I} \! - \! \int_{\Sigma_{\circlearrowright}
^{e}} \dfrac{w^{\Sigma_{\circlearrowright}^{e}}_{+}(s)}{s \! - \! z} \, \dfrac{
\md s}{2 \pi \mi} \! = \! \int_{\Sigma_{\scriptscriptstyle \blacksquare}^{e}}
\dfrac{w_{+}^{\Sigma_{\scriptscriptstyle \blacksquare}^{e}}(s)}{s \! - \! z}
\, \dfrac{\md s}{2 \pi \mi} \! + \! \sum_{k=1}^{8}I_{k}^{e},
\end{equation}
where $w_{+}^{\Sigma_{\circlearrowright}^{e}}(z) \! := \! w_{+}^{\Sigma_{
\mathscr{R}}^{e}}(z) \! \! \upharpoonright_{\Sigma_{\circlearrowright}^{e}}$,
$w_{+}^{\Sigma_{\scriptscriptstyle \blacksquare}^{e}}(z) \! := \! w_{+}^{
\Sigma_{\mathscr{R}}^{e}}(z) \! \! \upharpoonright_{\Sigma_{\scriptscriptstyle
\blacksquare}^{e}}$,
\begin{align*}
I_{1}^{e} :=& \int_{\widetilde{\Sigma}_{p}^{e}} \dfrac{(C^{\infty}_{w^{\Sigma_{
\scriptscriptstyle \blacksquare}^{e}}} \mathrm{I})(s)w_{+}^{\Sigma_{\mathscr{
R}}^{e}}(s)}{s \! - \! z} \, \dfrac{\md s}{2 \pi \mi}, \qquad \qquad \quad
I_{2}^{e} \! := \! \int_{\widetilde{\Sigma}_{p}^{e}} \dfrac{(C^{\infty}_{w^{
\Sigma_{\circlearrowright}^{e}}} \mathrm{I})(s)w_{+}^{\Sigma_{\mathscr{R}}^{e}}
(s)}{s \! - \! z} \, \dfrac{\md s}{2 \pi \mi}, \\
I_{3}^{e} :=& \int_{\widetilde{\Sigma}_{p}^{e}} \dfrac{((\id \! - \! C^{
\infty}_{w^{\Sigma_{\scriptscriptstyle \blacksquare}^{e}}})^{-1}C^{\infty}_{
w^{\Sigma_{\scriptscriptstyle \blacksquare}^{e}}}(C^{\infty}_{w^{\Sigma_{
\scriptscriptstyle \blacksquare}^{e}}} \mathrm{I}))(s)w_{+}^{\Sigma_{\mathscr{
R}}^{e}}(s)}{s \! - \! z} \, \dfrac{\md s}{2 \pi \mi}, \\
I_{4}^{e} :=& \int_{\widetilde{\Sigma}_{p}^{e}} \dfrac{((\id \! - \! C^{
\infty}_{w^{\Sigma_{\circlearrowright}^{e}}})^{-1}C^{\infty}_{w^{\Sigma_{
\circlearrowright}^{e}}}(C^{\infty}_{w^{\Sigma_{\circlearrowright}^{e}}}
\mathrm{I}))(s)w_{+}^{\Sigma_{\mathscr{R}}^{e}}(s)}{s \! - \! z} \, \dfrac{
\md s}{2 \pi \mi}, \\
I_{5}^{e} :=& \int_{\widetilde{\Sigma}_{p}^{e}}((\id \! - \! (\id \! - \! C^{
\infty}_{w^{\Sigma_{\circlearrowright}^{e}}})^{-1}(\id \! - \! C^{\infty}_{w^{
\Sigma_{\scriptscriptstyle \blacksquare}^{e}}})^{-1}C^{\infty}_{w^{\Sigma_{
\scriptscriptstyle \blacksquare}^{e}}}C^{\infty}_{w^{\Sigma_{
\circlearrowright}^{e}}})^{-1}(\id \! - \! C^{\infty}_{w^{\Sigma_{
\circlearrowright}^{e}}})^{-1}(\id \! - \! C^{\infty}_{w^{\Sigma_{
\scriptscriptstyle \blacksquare}^{e}}})^{-1} \\
\times& \, \dfrac{C^{\infty}_{w^{\Sigma_{\circlearrowright}^{e}}}(C^{\infty}_{
w^{\Sigma_{\scriptscriptstyle \blacksquare}^{e}}} \mathrm{I}))(s)w_{+}^{
\Sigma_{\mathscr{R}}^{e}}(s)}{s \! - \! z} \, \dfrac{\md s}{2 \pi \mi}, \\
I_{6}^{e} :=& \int_{\widetilde{\Sigma}_{p}^{e}}((\id \! - \! (\id \! - \! C^{
\infty}_{w^{\Sigma_{\scriptscriptstyle \blacksquare}^{e}}})^{-1}(\id \! - \!
C^{\infty}_{w^{\Sigma_{\circlearrowright}^{e}}})^{-1}C^{\infty}_{w^{\Sigma_{
\circlearrowright}^{e}}}C^{\infty}_{w^{\Sigma_{\scriptscriptstyle
\blacksquare}^{e}}})^{-1}(\id \! - \! C^{\infty}_{w^{\Sigma_{
\scriptscriptstyle \blacksquare}^{e}}})^{-1}(\id \! - \! C^{\infty}_{w^{
\Sigma_{\circlearrowright}^{e}}})^{-1} \\
\times& \, \dfrac{C^{\infty}_{w^{\Sigma_{\scriptscriptstyle \blacksquare}^{e}}}
(C^{\infty}_{w^{\Sigma_{\circlearrowright}^{e}}} \mathrm{I}))(s)w_{+}^{\Sigma_{
\mathscr{R}}^{e}}(s)}{s \! - \! z} \, \dfrac{\md s}{2 \pi \mi}, \\
I_{7}^{e} :=& \int_{\widetilde{\Sigma}_{p}^{e}}((\id \! - \! (\id \! - \! C^{
\infty}_{w^{\Sigma_{\circlearrowright}^{e}}})^{-1}(\id \! - \! C^{\infty}_{w^{
\Sigma_{\scriptscriptstyle \blacksquare}^{e}}})^{-1}C^{\infty}_{w^{\Sigma_{
\scriptscriptstyle \blacksquare}^{e}}}C^{\infty}_{w^{\Sigma_{
\circlearrowright}^{e}}})^{-1}(\id \! - \! C^{\infty}_{w^{\Sigma_{
\circlearrowright}^{e}}})^{-1}(\id \! - \! C^{\infty}_{w^{\Sigma_{
\scriptscriptstyle \blacksquare}^{e}}})^{-1} \\
\times& \, \dfrac{C^{\infty}_{w^{\Sigma_{\circlearrowright}^{e}}}(\id \! - \!
C^{\infty}_{w^{\Sigma_{\scriptscriptstyle \blacksquare}^{e}}})^{-1}C^{\infty}_{
w^{\Sigma_{\scriptscriptstyle \blacksquare}^{e}}}(C^{\infty}_{w^{\Sigma_{
\scriptscriptstyle \blacksquare}^{e}}} \mathrm{I}))(s)w_{+}^{\Sigma_{\mathscr{
R}}^{e}}(s)}{s \! - \! z} \, \dfrac{\md s}{2 \pi \mi}, \\
I_{8}^{e} :=& \int_{\widetilde{\Sigma}_{p}^{e}}((\id \! - \! (\id \! - \! C^{
\infty}_{w^{\Sigma_{\scriptscriptstyle \blacksquare}^{e}}})^{-1}(\id \! - \!
C^{\infty}_{w^{\Sigma_{\circlearrowright}^{e}}})^{-1}C^{\infty}_{w^{\Sigma_{
\circlearrowright}^{e}}}C^{\infty}_{w^{\Sigma_{\scriptscriptstyle
\blacksquare}^{e}}})^{-1}(\id \! - \! C^{\infty}_{w^{\Sigma_{
\scriptscriptstyle \blacksquare}^{e}}})^{-1}(\id \! - \! C^{\infty}_{w^{
\Sigma_{\circlearrowright}^{e}}})^{-1} \\
\times& \, \dfrac{C^{\infty}_{w^{\Sigma_{\scriptscriptstyle \blacksquare}^{e}}}
(\id \! - \! C^{\infty}_{w^{\Sigma_{\circlearrowright}^{e}}})^{-1}C^{\infty}_{
w^{\Sigma_{\circlearrowright}^{e}}}(C^{\infty}_{w^{\Sigma_{\circlearrowright}^{
e}}} \mathrm{I}))(s)w_{+}^{\Sigma_{\mathscr{R}}^{e}}(s)}{s \! - \! z} \,
\dfrac{\md s}{2 \pi \mi}.
\end{align*}
One now proceeds to estimate, as $n \! \to \! \infty$, the respective terms
on the right-hand side of Equation~(5.2) using the estimates and bounds given
in Propositions~5.1 and 5.2.
\begin{equation*}
\left\vert \int_{\Sigma_{\scriptscriptstyle \blacksquare}^{e}} \dfrac{w_{+}^{
\Sigma_{\scriptscriptstyle \blacksquare}^{e}}(s)}{s \! - \! z} \, \dfrac{\md
s}{2 \pi \mi} \right\vert \! \leqslant \! \int_{\Sigma_{\scriptscriptstyle
\blacksquare}^{e}} \dfrac{\vert w_{+}^{\Sigma_{\scriptscriptstyle
\blacksquare}^{e}}(s) \vert}{\vert s \! - \! z \vert} \, \dfrac{\vert \md
s \vert}{2 \pi} \! \leqslant \! \dfrac{\norm{w_{+}^{\Sigma_{\scriptscriptstyle
\blacksquare}^{e}}(\cdot)}_{\mathcal{L}^{1}_{\mathrm{M}_{2}(\mathbb{C})}
(\Sigma_{\scriptscriptstyle \blacksquare}^{e})}}{2 \pi \operatorname{dist}
(z,\Sigma_{\scriptscriptstyle \blacksquare}^{e})} \! \underset{n \to \infty}{
\leqslant} \! \mathcal{O} \! \left(\dfrac{f(n) \me^{-nc}}{n \operatorname{dist}
(z,\widetilde{\Sigma}_{p}^{e})} \right),
\end{equation*}
where, here and below, $(f(n) \! > \! 0$ and) $f(n) \! =_{n \to \infty} \!
\mathcal{O}(1)$ and $c \! > \! 0$. One estimates $I_{1}^{e}$ as follows:
\begin{align*}
\vert I_{1}^{e} \vert &\leqslant \int_{\widetilde{\Sigma}_{p}^{e}} \dfrac{
\vert (C^{\infty}_{w^{\Sigma_{\scriptscriptstyle \blacksquare}^{e}}} \mathrm{
I})(s) \vert \vert w_{+}^{\Sigma_{\mathscr{R}}^{e}}(s) \vert}{\vert s \! - \!
z \vert} \, \dfrac{\vert \md s \vert}{2 \pi} \! \leqslant \! \dfrac{\norm{(C^{
\infty}_{w^{\Sigma_{\scriptscriptstyle \blacksquare}^{e}}} \mathrm{I})(\cdot)
}_{\mathcal{L}^{2}_{\mathrm{M}_{2}(\mathbb{C})}(\widetilde{\Sigma}_{p}^{e})}
\norm{w_{+}^{\Sigma_{\mathscr{R}}^{e}}(\cdot)}_{\mathcal{L}^{2}_{\mathrm{M}_{2}
(\mathbb{C})}(\widetilde{\Sigma}_{p}^{e})}}{2 \pi \operatorname{dist}(z,
\widetilde{\Sigma}_{p}^{e})} \\
&\leqslant \dfrac{\operatorname{const.} \, \norm{w_{+}^{\Sigma_{
\scriptscriptstyle \blacksquare}^{e}}(\cdot)}_{\mathcal{L}^{2}_{\mathrm{M}_{2}
(\mathbb{C})}(\Sigma^{e}_{\scriptscriptstyle \blacksquare})}(\norm{w_{+}^{
\Sigma_{\circlearrowright}^{e}}(\cdot)}_{\mathcal{L}^{2}_{\mathrm{M}_{2}
(\mathbb{C})}(\Sigma^{e}_{\circlearrowright})} \! + \! \norm{w_{+}^{\Sigma_{
\scriptscriptstyle \blacksquare}^{e}}(\cdot)}_{\mathcal{L}^{2}_{\mathrm{M}_{2}
(\mathbb{C})}(\Sigma^{e}_{\scriptscriptstyle \blacksquare})})}{2 \pi
\operatorname{dist}(z,\widetilde{\Sigma}_{p}^{e})} \\
&\underset{n \to \infty}{\leqslant} \mathcal{O} \! \left(\dfrac{f(n) \me^{-nc}
}{\sqrt{\smash[b]{n}} \, \operatorname{dist}(z,\widetilde{\Sigma}_{p}^{e})}
\right) \! \left(\! \mathcal{O} \! \left(\dfrac{f(n)}{n} \right) \! + \!
\mathcal{O} \! \left(\dfrac{f(n) \me^{-nc}}{\sqrt{\smash[b]{n}}} \right)
\right) \! \underset{n \to \infty}{\leqslant} \! \mathcal{O} \! \left(\dfrac{
f(n) \me^{-nc}}{n \operatorname{dist}(z,\widetilde{\Sigma}_{p}^{e})} \right),
\end{align*}
where, here and below, $\operatorname{const.}$ denotes some positive,
$\mathcal{O}(1)$ constant; in going {}from the second line to the third line
in the above asymptotic (as $n \! \to \! \infty)$ estimation for $I_{1}^{e}$,
one uses the fact that, for $a,b \! > \! 0$, $\sqrt{\smash[b]{a^{2} \! + \!
b^{2}}} \! \leqslant \! \sqrt{\smash[b]{a^{2}}} \! + \! \sqrt{\smash[b]{b^{
2}}}$ (a fact used repeatedly below). One estimates $I_{2}^{e}$ as follows:
\begin{align*}
\vert I_{2}^{e} \vert &\leqslant \int_{\widetilde{\Sigma}_{p}^{e}} \dfrac{
\vert (C^{\infty}_{w^{\Sigma_{\circlearrowright}^{e}}} \mathrm{I})(s) \vert
\vert w_{+}^{\Sigma_{\mathscr{R}}^{e}}(s) \vert}{\vert s \! - \! z \vert} \,
\dfrac{\vert \md s \vert}{2 \pi} \! \leqslant \! \dfrac{\norm{(C^{\infty}_{w^{
\Sigma_{\circlearrowright}^{e}}} \mathrm{I})(\cdot)}_{\mathcal{L}^{2}_{\mathrm{
M}_{2}(\mathbb{C})}(\widetilde{\Sigma}_{p}^{e})} \norm{w_{+}^{\Sigma_{\mathscr{
R}}^{e}}(\cdot)}_{\mathcal{L}^{2}_{\mathrm{M}_{2}(\mathbb{C})}(\widetilde{
\Sigma}_{p}^{e})}}{2 \pi \operatorname{dist}(z,\widetilde{\Sigma}_{p}^{e})} \\
&\leqslant \dfrac{\operatorname{const.} \, \norm{w_{+}^{\Sigma_{
\circlearrowright}^{e}}(\cdot)}_{\mathcal{L}^{2}_{\mathrm{M}_{2}(\mathbb{C})}
(\Sigma^{e}_{\circlearrowright})}(\norm{w_{+}^{\Sigma_{\circlearrowright}^{e}}
(\cdot)}_{\mathcal{L}^{2}_{\mathrm{M}_{2}(\mathbb{C})}(\Sigma^{e}_{
\circlearrowright})} \! + \! \norm{w_{+}^{\Sigma_{\scriptscriptstyle
\blacksquare}^{e}}(\cdot)}_{\mathcal{L}^{2}_{\mathrm{M}_{2}(\mathbb{C})}
(\Sigma^{e}_{\scriptscriptstyle \blacksquare})})}{2 \pi \operatorname{dist}(z,
\widetilde{\Sigma}_{p}^{e})} \\
&\underset{n \to \infty}{\leqslant} \mathcal{O} \! \left(\dfrac{f(n)}{n
\operatorname{dist}(z,\widetilde{\Sigma}_{p}^{e})} \right) \! \left(\!
\mathcal{O} \! \left(\dfrac{f(n)}{n} \right) \! + \! \mathcal{O} \! \left(
\dfrac{f(n) \me^{-nc}}{\sqrt{\smash[b]{n}}} \right) \right) \! \underset{n
\to \infty}{\leqslant} \! \mathcal{O} \! \left(\dfrac{f(n)}{n^{2}
\operatorname{dist}(z,\widetilde{\Sigma}_{p}^{e})} \right).
\end{align*}
One estimates $I_{3}^{e}$ as follows:
\begin{align*}
\vert I_{3}^{e} \vert &\leqslant \int_{\widetilde{\Sigma}_{p}^{e}} \dfrac{
\vert ((\id \! - \! C^{\infty}_{w^{\Sigma_{\scriptscriptstyle \blacksquare}^{
e}}})^{-1}C^{\infty}_{w^{\Sigma_{\scriptscriptstyle \blacksquare}^{e}}}(C^{
\infty}_{w^{\Sigma_{\scriptscriptstyle \blacksquare}^{e}}} \mathrm{I}))(s)
\vert \vert w_{+}^{\Sigma_{\mathscr{R}}^{e}}(s) \vert}{\vert s \! - \! z
\vert} \, \dfrac{\vert \md s \vert}{2 \pi} \\
&\leqslant \dfrac{\norm{((\id \! - \! C^{\infty}_{w^{\Sigma_{
\scriptscriptstyle \blacksquare}^{e}}})^{-1}C^{\infty}_{w^{\Sigma_{
\scriptscriptstyle \blacksquare}^{e}}}(C^{\infty}_{w^{\Sigma_{
\scriptscriptstyle \blacksquare}^{e}}} \mathrm{I}))(\cdot)}_{\mathcal{L}^{2}_{
\mathrm{M}_{2}(\mathbb{C})}(\widetilde{\Sigma}_{p}^{e})} \norm{w_{+}^{\Sigma_{
\mathscr{R}}^{e}}(\cdot)}_{\mathcal{L}^{2}_{\mathrm{M}_{2}(\mathbb{C})}
(\widetilde{\Sigma}_{p}^{e})}}{2 \pi \operatorname{dist}(z,\widetilde{\Sigma}_{
p}^{e})} \\
&\leqslant \dfrac{\norm{(\id \! - \! C^{\infty}_{w^{\Sigma_{\scriptscriptstyle
\blacksquare}^{e}}})^{-1}}_{\mathscr{N}_{2}(\widetilde{\Sigma}_{p}^{e})}
\norm{C^{\infty}_{w^{\Sigma_{\scriptscriptstyle \blacksquare}^{e}}}}_{\mathscr{
N}_{2}(\widetilde{\Sigma}_{p}^{e})} \norm{(C^{\infty}_{w^{\Sigma_{
\scriptscriptstyle \blacksquare}^{e}}} \mathrm{I})(\cdot)}_{\mathcal{L}^{2}_{
\mathrm{M}_{2}(\mathbb{C})}(\widetilde{\Sigma}_{p}^{e})} \norm{w_{+}^{\Sigma_{
\mathscr{R}}^{e}}(\cdot)}_{\mathcal{L}^{2}_{\mathrm{M}_{2}(\mathbb{C})}
(\widetilde{\Sigma}_{p}^{e})}}{2 \pi \operatorname{dist}(z,\widetilde{\Sigma}_{
p}^{e})} \\
&\leqslant \dfrac{\operatorname{const.} \, \norm{(\id \! - \! C^{\infty}_{w^{
\Sigma_{\scriptscriptstyle \blacksquare}^{e}}})^{-1}}_{\mathscr{N}_{2}
(\widetilde{\Sigma}_{p}^{e})} \norm{C^{\infty}_{w^{\Sigma_{\scriptscriptstyle
\blacksquare}^{e}}}}_{\mathscr{N}_{2}(\widetilde{\Sigma}_{p}^{e})} \norm{w_{
+}^{\Sigma_{\scriptscriptstyle \blacksquare}^{e}}(\cdot)}_{\mathcal{L}^{2}_{
\mathrm{M}_{2}(\mathbb{C})}(\Sigma_{\scriptscriptstyle \blacksquare}^{e})}}{2
\pi \operatorname{dist}(z,\widetilde{\Sigma}_{p}^{e})} \\
&\times \left(\norm{w_{+}^{\Sigma_{\circlearrowright}^{e}}(\cdot)}_{\mathcal{
L}^{2}_{\mathrm{M}_{2}(\mathbb{C})}(\Sigma^{e}_{\circlearrowright})} \! + \!
\norm{w_{+}^{\Sigma_{\scriptscriptstyle \blacksquare}^{e}}(\cdot)}_{\mathcal{
L}^{2}_{\mathrm{M}_{2}(\mathbb{C})}(\Sigma^{e}_{\scriptscriptstyle
\blacksquare})} \right);
\end{align*}
using the fact that (cf. Proposition~5.2) $\norm{(\id \! - \! C^{\infty}_{w^{
\Sigma_{\scriptscriptstyle \blacksquare}^{e}}})^{-1}}_{\mathscr{N}_{2}
(\widetilde{\Sigma}_{p}^{e})} \! =_{n \to \infty} \! \mathcal{O}(1)$ (via
a Neuman series inversion argument, since $\norm{C^{\infty}_{w^{\Sigma_{
\scriptscriptstyle \blacksquare}^{e}}}}_{\mathscr{N}_{2}(\widetilde{\Sigma}_{
p}^{e})} \! =_{n \to \infty} \! \mathcal{O}(n^{-1}f(n) \me^{-nc}))$, one gets
that
\begin{equation*}
\vert I_{3}^{e} \vert \! \underset{n \to \infty}{\leqslant} \! \mathcal{O} \!
\left(\dfrac{f(n) \me^{-nc}}{n \operatorname{dist}(z,\widetilde{\Sigma}_{p}^{
e})} \right) \! \mathcal{O} \! \left(\dfrac{f(n) \me^{-nc}}{\sqrt{\smash[b]{n}}
} \right) \! \left(\! \mathcal{O} \! \left(\dfrac{f(n)}{n} \right) \! + \!
\mathcal{O} \! \left(\dfrac{f(n) \me^{-nc}}{\sqrt{\smash[b]{n}}} \right)
\right) \! \underset{n \to \infty}{\leqslant} \! \mathcal{O} \! \left(\dfrac{
f(n) \me^{-nc}}{n^{2} \operatorname{dist}(z,\widetilde{\Sigma}_{p}^{e})}
\right).
\end{equation*}
One estimates $I_{4}^{e}$ as follows:
\begin{align*}
\vert I_{4}^{e} \vert &\leqslant \int_{\widetilde{\Sigma}_{p}^{e}} \dfrac{
\vert ((\id \! - \! C^{\infty}_{w^{\Sigma_{\circlearrowright}^{e}}})^{-1}
C^{\infty}_{w^{\Sigma_{\circlearrowright}^{e}}}(C^{\infty}_{w^{\Sigma_{
\circlearrowright}^{e}}} \mathrm{I}))(s) \vert \vert w_{+}^{\Sigma_{\mathscr{
R}}^{e}}(s) \vert}{\vert s \! - \! z \vert} \, \dfrac{\vert \md s \vert}{2
\pi} \\
&\leqslant \dfrac{\norm{((\id \! - \! C^{\infty}_{w^{\Sigma_{
\circlearrowright}^{e}}})^{-1}C^{\infty}_{w^{\Sigma_{\circlearrowright}^{e}}}
(C^{\infty}_{w^{\Sigma_{\circlearrowright}^{e}}} \mathrm{I}))(\cdot)}_{
\mathcal{L}^{2}_{\mathrm{M}_{2}(\mathbb{C})}(\widetilde{\Sigma}_{p}^{e})}
\norm{w_{+}^{\Sigma_{\mathscr{R}}^{e}}(\cdot)}_{\mathcal{L}^{2}_{\mathrm{M}_{2}
(\mathbb{C})}(\widetilde{\Sigma}_{p}^{e})}}{2 \pi \operatorname{dist}(z,
\widetilde{\Sigma}_{p}^{e})} \\
&\leqslant \dfrac{\norm{(\id \! - \! C^{\infty}_{w^{\Sigma_{
\circlearrowright}^{e}}})^{-1}}_{\mathscr{N}_{2}(\widetilde{\Sigma}_{p}^{e})}
\norm{C^{\infty}_{w^{\Sigma_{\circlearrowright}^{e}}}}_{\mathscr{N}_{2}
(\widetilde{\Sigma}_{p}^{e})} \norm{(C^{\infty}_{w^{\Sigma_{
\circlearrowright}^{e}}} \mathrm{I})(\cdot)}_{\mathcal{L}^{2}_{\mathrm{M}_{2}
(\mathbb{C})}(\widetilde{\Sigma}_{p}^{e})} \norm{w_{+}^{\Sigma_{\mathscr{R}}^{
e}}(\cdot)}_{\mathcal{L}^{2}_{\mathrm{M}_{2}(\mathbb{C})}(\widetilde{\Sigma}_{
p}^{e})}}{2 \pi \operatorname{dist}(z,\widetilde{\Sigma}_{p}^{e})} \\
&\leqslant \dfrac{\operatorname{const.} \, \norm{(\id \! - \! C^{\infty}_{w^{
\Sigma_{\circlearrowright}^{e}}})^{-1}}_{\mathscr{N}_{2}(\widetilde{\Sigma}_{
p}^{e})} \norm{C^{\infty}_{w^{\Sigma_{\circlearrowright}^{e}}}}_{\mathscr{N}_{
2}(\widetilde{\Sigma}_{p}^{e})} \norm{w_{+}^{\Sigma_{\circlearrowright}^{e}}
(\cdot)}_{\mathcal{L}^{2}_{\mathrm{M}_{2}(\mathbb{C})}(\Sigma_{
\circlearrowright}^{e})}}{2 \pi \operatorname{dist}(z,\widetilde{\Sigma}_{p}^{
e})} \\
&\times \left(\norm{w_{+}^{\Sigma_{\circlearrowright}^{e}}(\cdot)}_{\mathcal{
L}^{2}_{\mathrm{M}_{2}(\mathbb{C})}(\Sigma^{e}_{\circlearrowright})} \! + \!
\norm{w_{+}^{\Sigma_{\scriptscriptstyle \blacksquare}^{e}}(\cdot)}_{\mathcal{
L}^{2}_{\mathrm{M}_{2}(\mathbb{C})}(\Sigma^{e}_{\scriptscriptstyle
\blacksquare})} \right);
\end{align*}
using the fact that (cf. Proposition~5.2) $\norm{(\id \! - \! C^{\infty}_{w^{
\Sigma_{\circlearrowright}^{e}}})^{-1}}_{\mathscr{N}_{2}(\widetilde{\Sigma}_{
p}^{e})} \! =_{n \to \infty} \! \mathcal{O}(1)$ (via a Neuman series inversion
argument, since $\norm{C^{\infty}_{w^{\Sigma_{\circlearrowright}^{e}}}}_{
\mathscr{N}_{2}(\widetilde{\Sigma}_{p}^{e})} \! =_{n \to \infty} \! \mathcal{O}
(n^{-1}f(n)))$, one gets that
\begin{equation*}
\vert I_{4}^{e} \vert \! \underset{n \to \infty}{\leqslant} \! \mathcal{O} \!
\left(\dfrac{f(n)}{n \operatorname{dist}(z,\widetilde{\Sigma}_{p}^{e})}
\right) \! \mathcal{O} \! \left(\dfrac{f(n)}{n} \right) \! \left(\! \mathcal{
O} \! \left(\dfrac{f(n)}{n} \right) \! + \! \mathcal{O} \! \left(\dfrac{f(n)
\me^{-nc}}{\sqrt{\smash[b]{n}}} \right) \right) \! \underset{n \to \infty}{
\leqslant} \! \mathcal{O} \! \left(\dfrac{f(n)}{n^{3} \operatorname{dist}(z,
\widetilde{\Sigma}_{p}^{e})} \right).
\end{equation*}
One estimates $I_{5}^{e}$ as follows:
\begin{align*}
\vert I_{5}^{e} \vert &\leqslant \int_{\widetilde{\Sigma}_{p}^{e}} \vert ((\id
\! - \! (\id \! - \! C^{\infty}_{w^{\Sigma_{\circlearrowright}^{e}}})^{-1}(\id
\! - \! C^{\infty}_{w^{\Sigma_{\scriptscriptstyle \blacksquare}^{e}}})^{-1}
C^{\infty}_{w^{\Sigma_{\scriptscriptstyle \blacksquare}^{e}}}C^{\infty}_{w^{
\Sigma_{\circlearrowright}^{e}}})^{-1}(\id \! - \! C^{\infty}_{w^{\Sigma_{
\circlearrowright}^{e}}})^{-1}(\id \! - \! C^{\infty}_{w^{\Sigma_{
\scriptscriptstyle \blacksquare}^{e}}})^{-1} \\
&\times \, \dfrac{C^{\infty}_{w^{\Sigma_{\circlearrowright}^{e}}}(C^{\infty}_{
w^{\Sigma_{\scriptscriptstyle \blacksquare}^{e}}} \mathrm{I}))(s) \vert \vert
w_{+}^{\Sigma_{\mathscr{R}}^{e}}(s) \vert}{\vert s \! - \! z \vert} \, \dfrac{
\vert \md s \vert}{2 \pi} \\
&\leqslant \, \vert \vert ((\id \! - \! (\id \! - \! C^{\infty}_{w^{\Sigma_{
\circlearrowright}^{e}}})^{-1}(\id \! - \! C^{\infty}_{w^{\Sigma_{
\scriptscriptstyle \blacksquare}^{e}}})^{-1}C^{\infty}_{w^{\Sigma_{
\scriptscriptstyle \blacksquare}^{e}}}C^{\infty}_{w^{\Sigma_{
\circlearrowright}^{e}}})^{-1}(\id \! - \! C^{\infty}_{w^{\Sigma_{
\circlearrowright}^{e}}})^{-1}(\id \! - \! C^{\infty}_{w^{\Sigma_{
\scriptscriptstyle \blacksquare}^{e}}})^{-1} \\
&\times \, \dfrac{C^{\infty}_{w^{\Sigma_{\circlearrowright}^{e}}}(C^{\infty}_{
w^{\Sigma_{\scriptscriptstyle \blacksquare}^{e}}} \mathrm{I}))(\cdot) \vert
\vert_{\mathcal{L}^{2}_{\mathrm{M}_{2}(\mathbb{C})}(\widetilde{\Sigma}_{p}^{
e})} \norm{w_{+}^{\Sigma_{\mathscr{R}}^{e}}(\cdot)}_{\mathcal{L}^{2}_{\mathrm{
M}_{2}(\mathbb{C})}(\widetilde{\Sigma}_{p}^{e})}}{2 \pi \operatorname{dist}
(z,\widetilde{\Sigma}_{p}^{e})} \\
&\leqslant \, \norm{(\id \! - \! (\id \! - \! C^{\infty}_{w^{\Sigma_{
\circlearrowright}^{e}}})^{-1}(\id \! - \! C^{\infty}_{w^{\Sigma_{
\scriptscriptstyle \blacksquare}^{e}}})^{-1}C^{\infty}_{w^{\Sigma_{
\scriptscriptstyle \blacksquare}^{e}}}C^{\infty}_{w^{\Sigma_{
\circlearrowright}^{e}}})^{-1}}_{\mathscr{N}_{2}(\widetilde{\Sigma}_{p}^{e})}
\norm{(\id \! - \! C^{\infty}_{w^{\Sigma_{\circlearrowright}^{e}}})^{-1}}_{
\mathscr{N}_{2}(\widetilde{\Sigma}_{p}^{e})} \\
&\times \, \dfrac{\norm{(\id \! - \! C^{\infty}_{w^{\Sigma_{\scriptscriptstyle
\blacksquare}^{e}}})^{-1}}_{\mathscr{N}_{2}(\widetilde{\Sigma}_{p}^{e})}
\norm{C^{\infty}_{w^{\Sigma_{\circlearrowright}^{e}}}}_{\mathscr{N}_{2}
(\widetilde{\Sigma}_{p}^{e})} \norm{(C^{\infty}_{w^{\Sigma_{\scriptscriptstyle
\blacksquare}^{e}}} \mathrm{I})(\cdot)}_{\mathcal{L}^{2}_{\mathrm{M}_{2}
(\mathbb{C})}(\widetilde{\Sigma}_{p}^{e})} \norm{w_{+}^{\Sigma_{\mathscr{R}}^{
e}}(\cdot)}_{\mathcal{L}^{2}_{\mathrm{M}_{2}(\mathbb{C})}(\widetilde{\Sigma}_{
p}^{e})}}{2 \pi \operatorname{dist}(z,\widetilde{\Sigma}_{p}^{e})} \\
&\leqslant \, \norm{(\id \! - \! (\id \! - \! C^{\infty}_{w^{\Sigma_{
\circlearrowright}^{e}}})^{-1}(\id \! - \! C^{\infty}_{w^{\Sigma_{
\scriptscriptstyle \blacksquare}^{e}}})^{-1}C^{\infty}_{w^{\Sigma_{
\scriptscriptstyle \blacksquare}^{e}}}C^{\infty}_{w^{\Sigma_{
\circlearrowright}^{e}}})^{-1}}_{\mathscr{N}_{2}(\widetilde{\Sigma}_{p}^{e})}
\norm{(\id \! - \! C^{\infty}_{w^{\Sigma_{\circlearrowright}^{e}}})^{-1}}_{
\mathscr{N}_{2}(\widetilde{\Sigma}_{p}^{e})} \\
&\times \, \operatorname{const.} \, \norm{(\id \! - \! C^{\infty}_{w^{\Sigma_{
\scriptscriptstyle \blacksquare}^{e}}})^{-1}}_{\mathscr{N}_{2}(\widetilde{
\Sigma}_{p}^{e})} \norm{C^{\infty}_{w^{\Sigma_{\circlearrowright}^{e}}}}_{
\mathscr{N}_{2}(\widetilde{\Sigma}_{p}^{e})} \norm{w_{+}^{\Sigma_{
\scriptscriptstyle \blacksquare}^{e}}(\cdot)}_{\mathcal{L}^{2}_{\mathrm{M}_{2}
(\mathbb{C})}(\Sigma_{\scriptscriptstyle \blacksquare}^{e})} \\
&\times \dfrac{(\norm{w_{+}^{\Sigma_{\circlearrowright}^{e}}(\cdot)}_{\mathcal{
L}^{2}_{\mathrm{M}_{2}(\mathbb{C})}(\Sigma^{e}_{\circlearrowright})} \! + \!
\norm{w_{+}^{\Sigma_{\scriptscriptstyle \blacksquare}^{e}}(\cdot)}_{\mathcal{
L}^{2}_{\mathrm{M}_{2}(\mathbb{C})}(\Sigma^{e}_{\scriptscriptstyle
\blacksquare})})}{2 \pi \operatorname{dist}(z,\widetilde{\Sigma}_{p}^{e})};
\end{align*}
using the fact that (cf. Proposition~5.2) $\norm{(\id \! - \! (\id \! - \! C^{
\infty}_{w^{\Sigma_{\circlearrowright}^{e}}})^{-1}(\id \! - \! C^{\infty}_{w^{
\Sigma_{\scriptscriptstyle \blacksquare}^{e}}})^{-1}C^{\infty}_{w^{\Sigma_{
\scriptscriptstyle \blacksquare}^{e}}}C^{\infty}_{w^{\Sigma_{
\circlearrowright}^{e}}})^{-1}}_{\mathscr{N}_{2}(\widetilde{\Sigma}_{p}^{e})}
\! =_{n \to \infty} \! \mathcal{O}(1)$ (via a Neuman series inversion 
argument, since $\norm{C^{\infty}_{w^{\Sigma_{\scriptscriptstyle 
\blacksquare}^{e}}}}_{\mathscr{N}_{2}(\widetilde{\Sigma}_{p}^{e})} \! =_{n 
\to \infty} \! \mathcal{O}(n^{-1}f(n) \me^{-nc})$ and $\norm{C^{\infty}_{w^{
\Sigma_{\circlearrowright}^{e}}}}_{\mathscr{N}_{2}(\widetilde{\Sigma}_{p}^{
e})} \linebreak[4]
=_{n \to \infty} \! \mathcal{O}(n^{-1}f(n)))$, one gets that
\begin{equation*}
\vert I_{5}^{e} \vert \! \underset{n \to \infty}{\leqslant} \! \mathcal{O}
\! \left(\dfrac{f(n)}{n \operatorname{dist}(z,\widetilde{\Sigma}_{p}^{e})}
\right) \! \mathcal{O} \! \left(\dfrac{f(n) \me^{-nc}}{\sqrt{\smash[b]{n}}}
\right) \! \left(\! \mathcal{O} \! \left(\dfrac{f(n)}{n} \right) \! + \!
\mathcal{O} \! \left(\dfrac{f(n) \me^{-nc}}{\sqrt{\smash[b]{n}}} \right)
\right) \! \underset{n \to \infty}{\leqslant} \! \mathcal{O} \! \left(
\dfrac{f(n) \me^{-nc}}{n^{2} \operatorname{dist}(z,\widetilde{\Sigma}_{p}^{
e})} \right).
\end{equation*}
One estimates $I_{6}^{e}$ as follows:
\begin{align*}
\vert I_{6}^{e} \vert &\leqslant \int_{\widetilde{\Sigma}_{p}^{e}} \vert ((\id
\! - \! (\id \! - \! C^{\infty}_{w^{\Sigma_{\scriptscriptstyle \blacksquare}^{
e}}})^{-1}(\id \! - \! C^{\infty}_{w^{\Sigma_{\circlearrowright}^{e}}})^{-1}
C^{\infty}_{w^{\Sigma_{\circlearrowright}^{e}}}C^{\infty}_{w^{\Sigma_{
\scriptscriptstyle \blacksquare}^{e}}})^{-1}(\id \! - \! C^{\infty}_{w^{
\Sigma_{\scriptscriptstyle \blacksquare}^{e}}})^{-1}(\id \! - \! C^{\infty}_{
w^{\Sigma_{\circlearrowright}^{e}}})^{-1} \\
&\times \, \dfrac{C^{\infty}_{w^{\Sigma_{\scriptscriptstyle \blacksquare}^{e}}}
(C^{\infty}_{w^{\Sigma_{\circlearrowright}^{e}}} \mathrm{I}))(s) \vert \vert
w_{+}^{\Sigma_{\mathscr{R}}^{e}}(s) \vert}{\vert s \! - \! z \vert} \, \dfrac{
\vert \md s \vert}{2 \pi} \\
&\leqslant \, \vert \vert ((\id \! - \! (\id \! - \! C^{\infty}_{w^{\Sigma_{
\scriptscriptstyle \blacksquare}^{e}}})^{-1}(\id \! - \! C^{\infty}_{w^{
\Sigma_{\circlearrowright}^{e}}})^{-1}C^{\infty}_{w^{\Sigma_{
\circlearrowright}^{e}}}C^{\infty}_{w^{\Sigma_{\scriptscriptstyle
\blacksquare}^{e}}})^{-1}(\id \! - \! C^{\infty}_{w^{\Sigma_{
\scriptscriptstyle \blacksquare}^{e}}})^{-1}(\id \! - \! C^{\infty}_{w^{
\Sigma_{\circlearrowright}^{e}}})^{-1} \\
&\times \, \dfrac{C^{\infty}_{w^{\Sigma_{\scriptscriptstyle \blacksquare}^{e}}}
(C^{\infty}_{w^{\Sigma_{\circlearrowright}^{e}}} \mathrm{I}))(\cdot) \vert
\vert_{\mathcal{L}^{2}_{\mathrm{M}_{2}(\mathbb{C})}(\widetilde{\Sigma}_{p}^{
e})} \norm{w_{+}^{\Sigma_{\mathscr{R}}^{e}}(\cdot)}_{\mathcal{L}^{2}_{\mathrm{
M}_{2}(\mathbb{C})}(\widetilde{\Sigma}_{p}^{e})}}{2 \pi \operatorname{dist}
(z,\widetilde{\Sigma}_{p}^{e})} \\
&\leqslant \, \norm{(\id \! - \! (\id \! - \! C^{\infty}_{w^{\Sigma_{
\scriptscriptstyle \blacksquare}^{e}}})^{-1}(\id \! - \! C^{\infty}_{w^{
\Sigma_{\circlearrowright}^{e}}})^{-1}C^{\infty}_{w^{\Sigma_{
\circlearrowright}^{e}}}C^{\infty}_{w^{\Sigma_{\scriptscriptstyle
\blacksquare}^{e}}})^{-1}}_{\mathscr{N}_{2}(\widetilde{\Sigma}_{p}^{e})} \norm{
(\id \! - \! C^{\infty}_{w^{\Sigma_{\scriptscriptstyle \blacksquare}^{e}}})^{-
1}}_{\mathscr{N}_{2}(\widetilde{\Sigma}_{p}^{e})} \\
&\times \, \dfrac{\norm{(\id \! - \! C^{\infty}_{w^{\Sigma_{
\circlearrowright}^{e}}})^{-1}}_{\mathscr{N}_{2}(\widetilde{\Sigma}_{p}^{e})}
\norm{C^{\infty}_{w^{\Sigma_{\scriptscriptstyle \blacksquare}^{e}}}}_{\mathscr{
N}_{2}(\widetilde{\Sigma}_{p}^{e})} \norm{(C^{\infty}_{w^{\Sigma_{
\circlearrowright}^{e}}} \mathrm{I})(\cdot)}_{\mathcal{L}^{2}_{\mathrm{M}_{2}
(\mathbb{C})}(\widetilde{\Sigma}_{p}^{e})} \norm{w_{+}^{\Sigma_{\mathscr{R}}^{
e}}(\cdot)}_{\mathcal{L}^{2}_{\mathrm{M}_{2}(\mathbb{C})}(\widetilde{\Sigma}_{
p}^{e})}}{2 \pi \operatorname{dist}(z,\widetilde{\Sigma}_{p}^{e})} \\
&\leqslant \, \norm{(\id \! - \! (\id \! - \! C^{\infty}_{w^{\Sigma_{
\scriptscriptstyle \blacksquare}^{e}}})^{-1}(\id \! - \! C^{\infty}_{w^{
\Sigma_{\circlearrowright}^{e}}})^{-1}C^{\infty}_{w^{\Sigma_{
\circlearrowright}^{e}}}C^{\infty}_{w^{\Sigma_{\scriptscriptstyle
\blacksquare}^{e}}})^{-1}}_{\mathscr{N}_{2}(\widetilde{\Sigma}_{p}^{e})} \norm{
(\id \! - \! C^{\infty}_{w^{\Sigma_{\scriptscriptstyle \blacksquare}^{e}}})^{
-1}}_{\mathscr{N}_{2}(\widetilde{\Sigma}_{p}^{e})} \\
&\times \, \operatorname{const.} \, \norm{(\id \! - \! C^{\infty}_{w^{\Sigma_{
\circlearrowright}^{e}}})^{-1}}_{\mathscr{N}_{2}(\widetilde{\Sigma}_{p}^{e})}
\norm{C^{\infty}_{w^{\Sigma_{\scriptscriptstyle \blacksquare}^{e}}}}_{\mathscr{
N}_{2}(\widetilde{\Sigma}_{p}^{e})} \norm{w_{+}^{\Sigma_{\circlearrowright}^{
e}}(\cdot)}_{\mathcal{L}^{2}_{\mathrm{M}_{2}(\mathbb{C})}(\Sigma_{
\circlearrowright}^{e})} \\
&\times \, \dfrac{(\norm{w_{+}^{\Sigma_{\circlearrowright}^{e}}(\cdot)}_{
\mathcal{L}^{2}_{\mathrm{M}_{2}(\mathbb{C})}(\Sigma^{e}_{\circlearrowright})}
\! + \! \norm{w_{+}^{\Sigma_{\scriptscriptstyle \blacksquare}^{e}}(\cdot)}_{
\mathcal{L}^{2}_{\mathrm{M}_{2}(\mathbb{C})}(\Sigma^{e}_{\scriptscriptstyle
\blacksquare})})}{2 \pi \operatorname{dist}(z,\widetilde{\Sigma}_{p}^{e})};
\end{align*}
using the fact that (cf. Proposition~5.2) $\norm{(\id \! - \! (\id \! - \! C^{
\infty}_{w^{\Sigma_{\scriptscriptstyle \blacksquare}^{e}}})^{-1}(\id \! - \!
C^{\infty}_{w^{\Sigma_{\circlearrowright}^{e}}})^{-1}C^{\infty}_{w^{\Sigma_{
\circlearrowright}^{e}}}C^{\infty}_{w^{\Sigma_{\scriptscriptstyle 
\blacksquare}^{e}}})^{-1}}_{\mathscr{N}_{2}(\widetilde{\Sigma}_{p}^{e})} 
\! =_{n \to \infty} \! \mathcal{O}(1)$ (via a Neuman series inversion 
argument, since $\norm{C^{\infty}_{w^{\Sigma_{\scriptscriptstyle 
\blacksquare}^{e}}}}_{\mathscr{N}_{2}(\widetilde{\Sigma}_{p}^{e})} \! =_{n 
\to \infty} \! \mathcal{O}(n^{-1}f(n) \me^{-nc})$ and $\norm{C^{\infty}_{w^{
\Sigma_{\circlearrowright}^{e}}}}_{\mathscr{N}_{2}(\widetilde{\Sigma}_{p}^{
e})} \linebreak[4]
=_{n \to \infty} \! \mathcal{O}(n^{-1}f(n)))$, one gets that
\begin{equation*}
\vert I_{6}^{e} \vert \! \underset{n \to \infty}{\leqslant} \! \mathcal{O} \!
\left(\dfrac{f(n) \me^{-nc}}{n \operatorname{dist}(z,\widetilde{\Sigma}_{p}^{
e})} \right) \! \mathcal{O} \! \left(\dfrac{f(n)}{n} \right) \! \left(\!
\mathcal{O} \! \left(\dfrac{f(n)}{n} \right) \! + \! \mathcal{O} \! \left(
\dfrac{f(n) \me^{-nc}}{\sqrt{\smash[b]{n}}} \right) \right) \! \underset{n
\to \infty}{\leqslant} \! \mathcal{O} \! \left(\dfrac{f(n) \me^{-nc}}{n^{3}
\operatorname{dist}(z,\widetilde{\Sigma}_{p}^{e})} \right).
\end{equation*}
One estimates $I_{7}^{e}$, succinctly, as follows:
\begin{align*}
\vert I_{7}^{e} \vert &\leqslant \int_{\widetilde{\Sigma}_{p}^{e}} \vert ((\id
\! - \! (\id \! - \! C^{\infty}_{w^{\Sigma_{\circlearrowright}^{e}}})^{-1}(\id
\! - \! C^{\infty}_{w^{\Sigma_{\scriptscriptstyle \blacksquare}^{e}}})^{-1}
C^{\infty}_{w^{\Sigma_{\scriptscriptstyle \blacksquare}^{e}}}C^{\infty}_{w^{
\Sigma_{\circlearrowright}^{e}}})^{-1}(\id \! - \! C^{\infty}_{w^{\Sigma_{
\circlearrowright}^{e}}})^{-1}(\id \! - \! C^{\infty}_{w^{\Sigma_{
\scriptscriptstyle \blacksquare}^{e}}})^{-1} \\
&\times \, \dfrac{C^{\infty}_{w^{\Sigma_{\circlearrowright}^{e}}}(\id \! - \!
C^{\infty}_{w^{\Sigma_{\scriptscriptstyle \blacksquare}^{e}}})^{-1}C^{\infty}_{
w^{\Sigma_{\scriptscriptstyle \blacksquare}^{e}}}(C^{\infty}_{w^{\Sigma_{
\scriptscriptstyle \blacksquare}^{e}}} \mathrm{I}))(s) \vert \vert w_{+}^{
\Sigma_{\mathscr{R}}^{e}}(s) \vert}{\vert s \! - \! z \vert} \, \dfrac{\vert
\md s \vert}{2 \pi} \\
&\leqslant \, \dfrac{\norm{(\id \! - \! (\id \! - \! C^{\infty}_{w^{\Sigma_{
\circlearrowright}^{e}}})^{-1}(\id \! - \! C^{\infty}_{w^{\Sigma_{
\scriptscriptstyle \blacksquare}^{e}}})^{-1}C^{\infty}_{w^{\Sigma_{
\scriptscriptstyle \blacksquare}^{e}}}C^{\infty}_{w^{\Sigma_{
\circlearrowright}^{e}}})^{-1}}_{\mathscr{N}_{2}(\widetilde{\Sigma}_{p}^{e})}
\norm{(\id \! - \! C^{\infty}_{w^{\Sigma_{\circlearrowright}^{e}}})^{-1}}_{
\mathscr{N}_{2}(\widetilde{\Sigma}_{p}^{e})}}{2 \pi \operatorname{dist}(z,
\widetilde{\Sigma}_{p}^{e})} \\
&\times \, \norm{(\id \! - \! C^{\infty}_{w^{\Sigma_{\scriptscriptstyle
\blacksquare}^{e}}})^{-1}}_{\mathscr{N}_{2}(\widetilde{\Sigma}_{p}^{e})}
\norm{C^{\infty}_{w^{\Sigma_{\circlearrowright}^{e}}}}_{\mathscr{N}_{2}
(\widetilde{\Sigma}_{p}^{e})} \norm{(\id \! - \! C^{\infty}_{w^{\Sigma_{
\scriptscriptstyle \blacksquare}^{e}}})^{-1}}_{\mathscr{N}_{2}(\widetilde{
\Sigma}_{p}^{e})} \norm{C^{\infty}_{w^{\Sigma_{\scriptscriptstyle
\blacksquare}^{e}}}}_{\mathscr{N}_{2}(\widetilde{\Sigma}_{p}^{e})} \\
&\times \, \operatorname{const.} \, \norm{w_{+}^{\Sigma_{\scriptscriptstyle
\blacksquare}^{e}}(\cdot)}_{\mathcal{L}^{2}_{\mathrm{M}_{2}(\mathbb{C})}
(\Sigma^{e}_{\scriptscriptstyle \blacksquare})} \! \left(\norm{w_{+}^{\Sigma_{
\scriptscriptstyle \blacksquare}^{e}}(\cdot)}_{\mathcal{L}^{2}_{\mathrm{M}_{2}
(\mathbb{C})}(\Sigma^{e}_{\scriptscriptstyle \blacksquare})} \! + \! \norm{w_{
+}^{\Sigma_{\circlearrowright}^{e}}(\cdot)}_{\mathcal{L}^{2}_{\mathrm{M}_{2}
(\mathbb{C})}(\Sigma^{e}_{\circlearrowright})} \right);
\end{align*}
using the fact that (established above) $\norm{(\id \! - \! (\id \! - \! C^{
\infty}_{w^{\Sigma_{\circlearrowright}^{e}}})^{-1}(\id \! - \! C^{\infty}_{w^{
\Sigma_{\scriptscriptstyle \blacksquare}^{e}}})^{-1}C^{\infty}_{w^{\Sigma_{
\scriptscriptstyle \blacksquare}^{e}}}C^{\infty}_{w^{\Sigma_{
\circlearrowright}^{e}}})^{-1}}_{\mathscr{N}_{2}(\widetilde{\Sigma}_{p}^{e})}
\! =_{n \to \infty} \! \mathcal{O}(1)$, one gets that
\begin{align*}
\vert I_{7}^{e} \vert &\underset{n \to \infty}{\leqslant} \mathcal{O} \! \left(
\dfrac{f(n)}{n \operatorname{dist}(z,\widetilde{\Sigma}_{p}^{e})} \right) \!
\mathcal{O} \! \left(\dfrac{f(n) \me^{-nc}}{n} \right) \! \mathcal{O} \! \left(
\dfrac{f(n) \me^{-nc}}{\sqrt{\smash[b]{n}}} \right) \! \left(\! \mathcal{O} \!
\left(\dfrac{f(n)}{n} \right) \! + \! \mathcal{O} \! \left(\dfrac{f(n) \me^{-
nc}}{\sqrt{\smash[b]{n}}} \right) \right) \\
&\underset{n \to \infty}{\leqslant} \mathcal{O} \! \left(\dfrac{f(n) \me^{-
nc}}{n^{3} \operatorname{dist}(z,\widetilde{\Sigma}_{p}^{e})} \right).
\end{align*}
One estimates $I_{8}^{e}$, succinctly, as follows:
\begin{align*}
\vert I_{8}^{e} \vert &\leqslant \int_{\widetilde{\Sigma}_{p}^{e}} \vert ((\id
\! - \! (\id \! - \! C^{\infty}_{w^{\Sigma_{\scriptscriptstyle \blacksquare}^{
e}}})^{-1}(\id \! - \! C^{\infty}_{w^{\Sigma_{\circlearrowright}^{e}}})^{-1}
C^{\infty}_{w^{\Sigma_{\circlearrowright}^{e}}}C^{\infty}_{w^{\Sigma_{
\scriptscriptstyle \blacksquare}^{e}}})^{-1}(\id \! - \! C^{\infty}_{w^{
\Sigma_{\scriptscriptstyle \blacksquare}^{e}}})^{-1}(\id \! - \! C^{\infty}_{
w^{\Sigma_{\circlearrowright}^{e}}})^{-1} \\
&\times \, \dfrac{C^{\infty}_{w^{\Sigma_{\scriptscriptstyle \blacksquare}^{e}}}
(\id \! - \! C^{\infty}_{w^{\Sigma_{\circlearrowright}^{e}}})^{-1}C^{\infty}_{
w^{\Sigma_{\circlearrowright}^{e}}}(C^{\infty}_{w^{\Sigma_{\circlearrowright}^{
e}}} \mathrm{I}))(s) \vert \vert w_{+}^{\Sigma_{\mathscr{R}}^{e}}(s) \vert}{
\vert s \! - \! z \vert} \, \dfrac{\vert \md s \vert}{2 \pi} \\
&\leqslant \, \dfrac{\norm{(\id \! - \! (\id \! - \! C^{\infty}_{w^{\Sigma_{
\scriptscriptstyle \blacksquare}^{e}}})^{-1}(\id \! - \! C^{\infty}_{w^{
\Sigma_{\circlearrowright}^{e}}})^{-1}C^{\infty}_{w^{\Sigma_{
\circlearrowright}^{e}}}C^{\infty}_{w^{\Sigma_{\scriptscriptstyle
\blacksquare}^{e}}})^{-1}}_{\mathscr{N}_{2}(\widetilde{\Sigma}_{p}^{e})}
\norm{(\id \! - \! C^{\infty}_{w^{\Sigma_{\scriptscriptstyle \blacksquare}^{
e}}})^{-1}}_{\mathscr{N}_{2}(\widetilde{\Sigma}_{p}^{e})}}{2 \pi
\operatorname{dist}(z,\widetilde{\Sigma}_{p}^{e})} \\
&\times \, \norm{(\id \! - \! C^{\infty}_{w^{\Sigma_{\circlearrowright}^{e}}}
)^{-1}}_{\mathscr{N}_{2}(\widetilde{\Sigma}_{p}^{e})} \norm{C^{\infty}_{w^{
\Sigma_{\scriptscriptstyle \blacksquare}^{e}}}}_{\mathscr{N}_{2}(\widetilde{
\Sigma}_{p}^{e})} \norm{(\id \! - \! C^{\infty}_{w^{\Sigma_{
\circlearrowright}^{e}}})^{-1}}_{\mathscr{N}_{2}(\widetilde{\Sigma}_{p}^{e})}
\norm{C^{\infty}_{w^{\Sigma_{\circlearrowright}^{e}}}}_{\mathscr{N}_{2}
(\widetilde{\Sigma}_{p}^{e})} \\
&\times \, \operatorname{const.} \, \norm{w_{+}^{\Sigma_{\circlearrowright}^{
e}}(\cdot)}_{\mathcal{L}^{2}_{\mathrm{M}_{2}(\mathbb{C})}(\Sigma^{e}_{
\circlearrowright})} \! \left(\norm{w_{+}^{\Sigma_{\scriptscriptstyle
\blacksquare}^{e}}(\cdot)}_{\mathcal{L}^{2}_{\mathrm{M}_{2}(\mathbb{C})}
(\Sigma^{e}_{\scriptscriptstyle \blacksquare})} \! + \! \norm{w_{+}^{\Sigma_{
\circlearrowright}^{e}}(\cdot)}_{\mathcal{L}^{2}_{\mathrm{M}_{2}(\mathbb{C})}
(\Sigma^{e}_{\circlearrowright})} \right);
\end{align*}
using the fact that (established above) $\norm{(\id \! - \! (\id \! - \! C^{
\infty}_{w^{\Sigma_{\scriptscriptstyle \blacksquare}^{e}}})^{-1}(\id \! - \!
C^{\infty}_{w^{\Sigma_{\circlearrowright}^{e}}})^{-1}C^{\infty}_{w^{\Sigma_{
\circlearrowright}^{e}}}C^{\infty}_{w^{\Sigma_{\scriptscriptstyle
\blacksquare}^{e}}})^{-1}}_{\mathscr{N}_{2}(\widetilde{\Sigma}_{p}^{e})} \!
=_{n \to \infty} \! \mathcal{O}(1)$, one gets that
\begin{align*}
\vert I_{8}^{e} \vert &\underset{n \to \infty}{\leqslant} \mathcal{O} \! \left(
\dfrac{f(n) \me^{-nc}}{n \operatorname{dist}(z,\widetilde{\Sigma}_{p}^{e})}
\right) \! \mathcal{O} \! \left(\dfrac{f(n)}{n} \right) \! \mathcal{O} \!
\left(\dfrac{f(n)}{n} \right) \! \left(\! \mathcal{O} \! \left(\dfrac{f(n)}{n}
\right) \! + \! \mathcal{O} \! \left(\dfrac{f(n) \me^{-nc}}{\sqrt{\smash[b]{
n}}} \right) \right) \\
&\underset{n \to \infty}{\leqslant} \mathcal{O} \! \left(\dfrac{f(n) \me^{-
nc}}{n^{4} \operatorname{dist}(z,\widetilde{\Sigma}_{p}^{e})} \right).
\end{align*}
Gathering the above-derived bounds, one arrives at the result stated in the 
Lemma. \hfill $\qed$
\begin{ccccc}
Let $\mathscr{R}^{e} \colon \mathbb{C} \setminus \widetilde{\Sigma}_{p}^{e} 
\! \to \! \operatorname{SL}_{2}(\mathbb{C})$ be the solution of the {\rm RHP} 
$(\mathscr{R}^{e}(z),\upsilon_{\mathscr{R}}^{e}(z),\widetilde{\Sigma}_{p}^{
e})$ formulated in Proposition~{\rm 5.2} with the $n \! \to \! \infty$ 
integral representation given in Lemma~{\rm 5.2}. Then, uniformly for compact 
subsets of $\mathbb{C} \setminus \widetilde{\Sigma}_{p}^{e} \! \ni \! z$,
\begin{equation*}
\mathscr{R}^{e}(z) \underset{\underset{z \in \mathbb{C} \setminus \widetilde{
\Sigma}_{p}^{e}}{n \to \infty}}{=} \mathrm{I} \! + \! \dfrac{1}{n} \! \left(
\mathscr{R}^{e}_{\infty}(z) \! - \! \widetilde{\mathscr{R}}^{e}_{\infty}(z)
\right) \! + \! \mathcal{O} \! \left(\dfrac{f(z;n)}{n^{2}} \right),
\end{equation*}
where $\mathscr{R}^{e}_{\infty}(z)$ is defined in Theorem~{\rm 2.3.1}, 
Equations~{\rm (2.23)--(2.57)}, $\widetilde{\mathscr{R}}^{e}_{\infty}(z)$ 
is defined in Theorem~{\rm 2.3.1}, Equations~{\rm (2.14)--(2.20)} 
and~{\rm (2.70)--(2.74)}, and $f(z;n)$, where the $n$-dependence arises due 
to the $n$-dependence of the associated Riemann theta functions, is a bounded 
(with respect to $z$ and $n)$, $\operatorname{GL}_{2}(\mathbb{C})$-valued 
function which is analytic (with respect to $z)$ for $z \! \in \! \mathbb{C} 
\setminus \widetilde{\Sigma}_{p}^{e}$, and $(f(\bm{\cdot};n))_{kl} \! =_{n 
\to \infty} \! \mathcal{O}(1)$, $k,l \! = \! 1,2$.
\end{ccccc}
\begin{eeeee}
Note {}from the formulation of Lemma~5.3 above that (cf. Theorem~2.3.1, 
Equations (2.24)--(2.27)), for $j \! = \! 1,\dotsc,N \! + \! 1$, 
$\operatorname{tr}(\mathscr{A}^{e}(a_{j}^{e})) \! = \! \operatorname{tr}
(\mathscr{A}^{e}(b_{j-1}^{e})) \! = \! \operatorname{tr}(\mathscr{B}^{e}
(a_{j}^{e})) \! = \! \operatorname{tr}(\mathscr{B}^{e}(b_{j-1}^{e})) \! = \! 
0$. \hfill $\blacksquare$
\end{eeeee}

\emph{Proof.} Recall the integral representation for $\mathscr{R}^{e} \colon
\mathbb{C} \setminus \widetilde{\Sigma}_{p}^{e} \! \to \! \operatorname{SL}_{
2}(\mathbb{C})$ given in Lemma~5.2:
\begin{equation*}
\mathscr{R}^{e}(z) \underset{n \to \infty}{=} \mathrm{I} \! + \! \int_{
\Sigma^{e}_{\circlearrowright}} \dfrac{w_{+}^{\Sigma^{e}_{\circlearrowright}}
(s)}{s \! - \! z} \, \dfrac{\md s}{2 \pi \mi} \! + \! \mathcal{O} \! \left(\!
\dfrac{f(n)}{n^{2} \operatorname{dist}(z,\widetilde{\Sigma}_{p}^{e})} \right),
\quad z \! \in \! \mathbb{C} \setminus \widetilde{\Sigma}_{p}^{e},
\end{equation*}
where $\Sigma^{e}_{\circlearrowright} \! := \! \cup_{j=1}^{N+1}(\partial
\mathbb{U}^{e}_{\delta_{b_{j-1}}} \cup \partial \mathbb{U}^{e}_{\delta_{a_{
j}}})$, and $(f(n))_{kl} \! =_{n \to \infty} \! \mathcal{O}(1)$, $k,l \! =
\! 1,2$. Recalling that the radii of the open discs $\mathbb{U}^{e}_{\delta_{
b_{j-1}}},\mathbb{U}^{e}_{\delta_{a_{j}}}$, $j \! = \! 1,\dotsc,N \! + \! 1$,
are chosen, amongst other factors (cf. Lemmas~4.6 and~4.7), so that $\mathbb{
U}^{e}_{\delta_{b_{j-1}}} \cap \mathbb{U}^{e}_{\delta_{a_{k}}} \! = \!
\varnothing$, $j,k \! = \! 1,\dotsc,N \! + \! 1$, it follows {}from the
above integral representation that
\begin{equation*}
\mathscr{R}^{e}(z) \underset{n \to \infty}{=} \mathrm{I} \! - \! \sum_{j=1}^{
N+1} \! \left(\! \oint_{\partial \mathbb{U}^{e}_{\delta_{b_{j-1}}}}+\oint_{
\partial \mathbb{U}^{e}_{\delta_{a_{j}}}} \right) \! \dfrac{w_{+}^{\Sigma^{
e}_{\circlearrowright}}(s)}{s \! - \! z} \, \dfrac{\md s}{2 \pi \mi} \! + \!
\mathcal{O} \! \left(\! \dfrac{f(n)}{n^{2} \operatorname{dist}(z,\widetilde{
\Sigma}_{p}^{e})} \right), \quad z \! \in \! \mathbb{C} \setminus \widetilde{
\Sigma}_{p}^{e},
\end{equation*}
where $\oint_{\partial \mathbb{U}^{e}_{\delta_{b_{j-1}}}},\oint_{\partial
\mathbb{U}^{e}_{\delta_{a_{j}}}}$, $j \! = \! 1,\dotsc,N \! + \! 1$, are
counter-clockwise-oriented, closed (contour) integrals (Figure~10) about the
end-points of the support of the `even' equilibrium measure, $\lbrace b_{j-
1}^{e},a_{j}^{e} \rbrace_{j=1}^{N+1}$. The evaluation of these $2(N \! + \!
1)$ contour integrals requires the application of the Cauchy and Residue
Theorems; and, since the evaluation of the respective integrals entails
analogous calculations, consider, say, and without loss of generality, the
evaluation of the integrals about the right-most end-points $a_{j}^{e}$, $j
\! = \! 1,\dotsc,N$, namely:
\begin{equation*}
\oint_{\partial \mathbb{U}^{e}_{\delta_{a_{j}}}} \dfrac{w_{+}^{\Sigma^{e}_{
\circlearrowright}}(s)}{s \! - \! z} \, \dfrac{\md s}{2 \pi \mi}, \quad j \!
= \! 1,\dotsc,N.
\end{equation*}
Recalling {}from Lemma~4.7 that $\xi_{a_{j}}^{e}(z) \! = \! (z \! - \! a_{j}^{
e})^{3/2}G_{a_{j}}^{e}(z)$, $z \! \in \! \mathbb{U}^{e}_{\delta_{a_{j}}}
\setminus (-\infty,a_{j}^{e})$, $j \! = \! 1,\dotsc,N$, it follows {}from
item~(5) of Proposition~5.1 that, since $w_{+}^{\Sigma^{e}_{\circlearrowright}}
(z) \! = \! \upsilon_{\mathscr{R}}^{e}(z) \! - \! \mathrm{I}$,
\begin{align*}
w_{+}^{\Sigma^{e}_{\circlearrowright}}(z) \underset{\underset{z \in \mathbb{
C}_{\pm} \cap \partial \mathbb{U}_{\delta_{a_{j}}}^{e}}{n \to \infty}}{=}& \,
\dfrac{1}{n(z \! - \! a_{j}^{e})^{3/2}G_{a_{j}}^{e}(z)}
\overset{e}{\mathfrak{M}}^{\raise-1.0ex\hbox{$\scriptstyle \infty$}}(z) \!
\begin{pmatrix}
\mp (s_{1}+t_{1}) & \pm \mi (s_{1}-t_{1}) \me^{\mi n \Omega_{j}^{e}} \\
\pm \mi (s_{1}-t_{1}) \me^{-\mi n \Omega_{j}^{e}} & \pm (s_{1}+t_{1})
\end{pmatrix} \!
(\overset{e}{\mathfrak{M}}^{\raise-1.0ex\hbox{$\scriptstyle \infty$}}(z))^{-1}
\\
+& \, \mathcal{O} \! \left(\dfrac{1}{n^{2}(z \! - \! a_{j}^{e})^{3}(G_{a_{
j}}^{e}(z))^{2}}
\overset{e}{\mathfrak{M}}^{\raise-1.0ex\hbox{$\scriptstyle \infty$}}(z)f_{a_{
j}}^{e}(n)(\overset{e}{\mathfrak{M}}^{\raise-1.0ex\hbox{$\scriptstyle \infty$}}
(z))^{-1} \right), \quad j \! = \! 1,\dotsc,N,
\end{align*}
where $\overset{e}{\mathfrak{M}}^{\raise-1.0ex\hbox{$\scriptstyle \infty$}}
(z)$ and $\Omega_{j}^{e}$ are defined in Lemma~4.5, and $(f_{a_{j}}^{e}
(n))_{kl} \! =_{n \to \infty} \! \mathcal{O}(1)$, $k,l \! = \! 1,2$.
A matrix-multipl\-i\-c\-a\-t\-i\-o\-n argument shows that
$\overset{e}{\mathfrak{M}}^{\raise-1.0ex\hbox{$\scriptstyle \infty$}}(z) \!
\left(
\begin{smallmatrix}
\mp (s_{1}+t_{1}) & \pm \mi (s_{1}-t_{1}) \me^{\mi n \Omega_{j}^{e}} \\
\pm \mi (s_{1}-t_{1}) \me^{-\mi n \Omega_{j}^{e}} & \pm (s_{1}+t_{1})
\end{smallmatrix}
\right) \!
(\overset{e}{\mathfrak{M}}^{\raise-1.0ex\hbox{$\scriptstyle \infty$}}(z))^{-
1}$ is given by
\begin{equation*}
\begin{pmatrix}
\boxed{\begin{matrix}
\mp \frac{1}{4}(s_{1} \! + \! t_{1}) \! \left(\frac{(\gamma^{e}(z))^{2}+1}{
\gamma^{e}(z)} \right)^{2} \! \mathfrak{m}^{e}_{11}(z) \mathfrak{m}^{e}_{22}
(z) \\
\mp \frac{1}{4}(s_{1} \! + \! t_{1}) \! \left(\frac{(\gamma^{e}(z))^{2}-1}{
\gamma^{e}(z)} \right)^{2} \! \mathfrak{m}^{e}_{12}(z) \mathfrak{m}^{e}_{21}
(z) \\
\mp \frac{1}{4}(s_{1} \! - \! t_{1}) \! \left(\frac{(\gamma^{e}(z))^{4}-1}{(
\gamma^{e}(z))^{2}} \right) \! \mathfrak{m}^{e}_{11}(z) \mathfrak{m}^{e}_{21}
(z) \me^{\mi n \Omega_{j}^{e}} \\
\mp \frac{1}{4}(s_{1} \! - \! t_{1}) \! \left(\frac{(\gamma^{e}(z))^{4}-1}{(
\gamma^{e}(z))^{2}} \right) \! \mathfrak{m}^{e}_{12}(z) \mathfrak{m}^{e}_{22}
(z) \me^{-\mi n \Omega_{j}^{e}}
\end{matrix}} &
\boxed{\begin{matrix}
\pm \frac{\mi}{2}(s_{1} \! + \! t_{1}) \! \left(\frac{(\gamma^{e}(z))^{4}-1}{(
\gamma^{e}(z))^{2}} \right) \! \mathfrak{m}^{e}_{11}(z) \mathfrak{m}^{e}_{12}
(z) \\
\pm \frac{\mi}{4}(s_{1} \! - \! t_{1}) \! \left(\frac{(\gamma^{e}(z))^{2}+1}{
\gamma^{e}(z)} \right)^{2} \! (\mathfrak{m}^{e}_{11}(z))^{2} \me^{\mi n
\Omega_{j}^{e}} \\
\pm \frac{\mi}{4}(s_{1} \! - \! t_{1}) \! \left(\frac{(\gamma^{e}(z))^{2}-1}{
\gamma^{e}(z)} \right)^{2} \! (\mathfrak{m}^{e}_{12}(z))^{2} \me^{-\mi n
\Omega_{j}^{e}}
\end{matrix}} \\
\boxed{\begin{matrix}
\pm \frac{\mi}{2}(s_{1} \! + \! t_{1}) \! \left(\frac{(\gamma^{e}(z))^{4}-1}{(
\gamma^{e}(z))^{2}} \right) \! \mathfrak{m}^{e}_{21}(z) \mathfrak{m}^{e}_{22}
(z) \\
\pm \frac{\mi}{4}(s_{1} \! - \! t_{1}) \! \left(\frac{(\gamma^{e}(z))^{2}-1}{
\gamma^{e}(z)} \right)^{2} \! (\mathfrak{m}^{e}_{21}(z))^{2} \me^{\mi n
\Omega_{j}^{e}} \\
\pm \frac{\mi}{4}(s_{1} \! - \! t_{1}) \! \left(\frac{(\gamma^{e}(z))^{2}+1}{
\gamma^{e}(z)} \right)^{2} \! (\mathfrak{m}^{e}_{22}(z))^{2} \me^{-\mi n
\Omega_{j}^{e}}
\end{matrix}} &
\boxed{\begin{matrix}
\pm \frac{1}{4}(s_{1} \! + \! t_{1}) \! \left(\frac{(\gamma^{e}(z))^{2}+1}{
\gamma^{e}(z)} \right)^{2} \! \mathfrak{m}^{e}_{11}(z) \mathfrak{m}^{e}_{22}
(z) \\
\pm \frac{1}{4}(s_{1} \! + \! t_{1}) \! \left(\frac{(\gamma^{e}(z))^{2}-1}{
\gamma^{e}(z)} \right)^{2} \! \mathfrak{m}^{e}_{12}(z) \mathfrak{m}^{e}_{21}
(z) \\
\pm \frac{1}{4}(s_{1} \! - \! t_{1}) \! \left(\frac{(\gamma^{e}(z))^{4}-1}{(
\gamma^{e}(z))^{2}} \right) \! \mathfrak{m}^{e}_{11}(z) \mathfrak{m}^{e}_{21}
(z) \me^{\mi n \Omega_{j}^{e}} \\
\pm \frac{1}{4}(s_{1} \! - \! t_{1}) \! \left(\frac{(\gamma^{e}(z))^{4}-1}{(
\gamma^{e}(z))^{2}} \right) \! \mathfrak{m}^{e}_{12}(z) \mathfrak{m}^{e}_{22}
(z) \me^{-\mi n \Omega_{j}^{e}}
\end{matrix}}
\end{pmatrix},
\end{equation*}
where $s_{1}$ and $t_{1}$ are given in Theorem~2.3.1, Equations~(2.28), 
$\gamma^{e}(z)$ is defined in Lemma~4.4, and $\mathfrak{m}^{e}_{kl}(z)$, 
$k,l \! = \! 1,2$, are defined in Theorem~2.3.1, Equations~(2.17)--(2.20). 
Recall that, for $j \! = \! 1,\dotsc,N$, $\omega_{j}^{e} \! = \! \sum_{k=1}^{
N}c_{jk}^{e}(R_{e}(z))^{-1/2}z^{N-k} \, \md z$, where $c_{jk}^{e}$, $j,k 
\! = \! 1,\dotsc,N$, are obtained {}from Equations~(E1) and~(E2), and (the 
multi-valued function) $(R_{e}(z))^{1/2}$ is defined in Theorem~2.3.1, 
Equation~(2.8). One shows that
\begin{equation*}
\omega_{m}^{e} \underset{\underset{j=1,\dotsc,N}{z \to a_{j}^{e}}}{=} \dfrac{
(\mathfrak{f}_{e}(a_{j}^{e}))^{-1}}{\sqrt{\smash[b]{z \! - \! a_{j}^{e}}}}
\! \left(\mathfrak{p}_{m}^{\natural}(a_{j}^{e}) \! + \! \mathfrak{q}_{m}^{
\natural}(a_{j}^{e})(z \! - \! a_{j}^{e}) \! + \! \mathfrak{r}_{m}^{\natural}
(a_{j}^{e})(z \! - \! a_{j}^{e})^{2} \! + \! \mathcal{O}((z \! - \! a_{j}^{
e})^{3}) \right) \md z, \quad m \! = \! 1,\dotsc,N,
\end{equation*}
where
\begin{gather*}
\mathfrak{f}_{e}(\xi) \! = \! (-1)^{N-j+1} \! \left((a_{N+1}^{e} \! - \! \xi)
(\xi \! - \! b_{0}^{e})(b_{j}^{e} \! - \! \xi) \prod_{k=1}^{j-1}(\xi \! - \! 
b_{k}^{e})(\xi \! - \! a_{k}^{e}) \prod_{l=j+1}^{N}(b_{l}^{e} \! - \! \xi)
(a_{l}^{e} \! - \! \xi) \right)^{1/2}, \\
\mathfrak{p}_{m}^{\natural}(\xi) \! = \! \sum_{k=1}^{N}c_{mk}^{e} \xi^{N-k}, 
\qquad \qquad \quad \mathfrak{q}_{m}^{\natural}(\xi) \! = \! \sum_{k=1}^{N}
c_{mk}^{e} \xi^{N-k-1} \! \left(N \! - \! k \! - \! \dfrac{\xi \mathfrak{
f}_{e}^{\prime}(\xi)}{\mathfrak{f}_{e}(\xi)} \right), \\
\mathfrak{r}_{m}^{\natural}(\xi) \! = \! \sum_{k=1}^{N}c_{mk}^{e} \xi^{N-k-2} 
\! \left(\dfrac{(N \! - \! k)(N \! - \! k \! - \! 1)}{2} \! - \! \dfrac{(N 
\! - \! k) \xi \mathfrak{f}_{e}^{\prime}(\xi)}{\mathfrak{f}_{e}(\xi)} \! + \! 
\xi^{2} \! \left(\left(\dfrac{\mathfrak{f}_{e}^{\prime}(\xi)}{\mathfrak{f}_{
e}(\xi)} \right)^{2} \! - \! \dfrac{\mathfrak{f}_{e}^{\prime \prime}(\xi)}{2 
\mathfrak{f}_{e}(\xi)} \right) \right),
\end{gather*}
with $(-1)^{-N+j-1} \mathfrak{f}_{e}(a_{j}^{e}) \! > \! 0$,
\begin{align*}
\mathfrak{f}_{e}^{\prime}(\xi) =& \, \dfrac{1}{2} \mathfrak{f}_{e}(\xi) \! 
\left(\sum_{\substack{k=1\\k \not= j}}^{N} \left(\dfrac{1}{\xi\! - \! b_{k}^{
e}} \! + \! \dfrac{1}{\xi \! - \! a_{k}^{e}} \right) \! + \! \dfrac{1}{\xi \! 
- \! b_{j}^{e}} \! + \! \dfrac{1}{\xi \! - \! a_{N+1}^{e}} \! + \! \dfrac{
1}{\xi \! - \! b_{0}^{e}} \right), \\
\mathfrak{f}_{e}^{\prime \prime}(\xi) =& -\dfrac{1}{2} \mathfrak{f}_{e}(\xi) 
\! \left(\sum_{\substack{k=1\\k \not= j}}^{N} \left(\dfrac{1}{(\xi \! - \! 
b_{k}^{e})^{2}} \! + \! \dfrac{1}{(\xi \! - \! a_{k}^{e})^{2}} \right) \! + 
\! \dfrac{1}{(\xi \! - \! b_{j}^{e})^{2}} \! + \! \dfrac{1}{(\xi \! - \! 
a_{N+1}^{e})^{2}} \! + \! \dfrac{1}{(\xi \! - \! b_{0}^{e})^{2}} \right) \\
+& \, \dfrac{1}{4} \mathfrak{f}_{e}(\xi) \! \left(\sum_{\substack{k=1\\k\not= 
j}}^{N} \left(\dfrac{1}{\xi \! - \! b_{k}^{e}} \! + \! \dfrac{1}{\xi \! - \! 
a_{k}^{e}} \right) \! + \! \dfrac{1}{\xi \! - \! b_{j}^{e}} \! + \! \dfrac{
1}{\xi \! - \! a_{N+1}^{e}} \! + \! \dfrac{1}{\xi \! - \! b_{0}^{e}} 
\right)^{2}.
\end{align*}
Recall (cf. Lemma~4.5), also, that $\bm{u}^{e} \! \equiv \! \int_{a_{N+1}^{e}
}^{z} \bm{\omega}^{e}$ $(\in \! \operatorname{Jac}(\mathcal{Y}_{e}))$, where 
$\equiv$ denotes congruence modulo the period lattice, with $\bm{\omega}^{
e} \! := \! (\omega_{1}^{e},\omega_{2}^{e},\dotsc,\omega_{N}^{e})$; hence, 
via the above expansion (as $z \! \to \! a_{j}^{e}$, $j \! = \! 1,\dotsc,N)$ 
for $\omega_{m}^{e}$, $m \! = \! 1,\dotsc,N$, one arrives at
\begin{equation*}
\int_{a_{j}^{e}}^{z} \omega_{m}^{e} \underset{\underset{j=1,\dotsc,N}{z 
\to a_{j}^{e}}}{\equiv} \dfrac{2 \mathfrak{p}_{m}^{\natural}(a_{j}^{e})}{
\mathfrak{f}_{e}(a_{j}^{e})}(z \! - \! a_{j}^{e})^{1/2} \! + \! \dfrac{2 
\mathfrak{q}_{m}^{\natural}(a_{j}^{e})}{3 \mathfrak{f}_{e}(a_{j}^{e})}(z \! 
- \! a_{j}^{e})^{3/2} \! + \! \dfrac{2 \mathfrak{r}_{m}^{\natural}(a_{j}^{
e})}{5 \mathfrak{f}_{e}(a_{j}^{e})}(z \! - \! a_{j}^{e})^{5/2} \! + \! 
\mathcal{O}((z \! - \! a_{j}^{e})^{7/2}).
\end{equation*}
{}From the definition of $\mathfrak{m}^{e}_{kl}(z)$, $k,l \! = \! 1,2$, given 
in Theorem~2.3.1, Equations~(2.17)--(2.20), the definition of the `even' 
Riemann theta function given by Equation~(2.1), and recalling that $\mathfrak{
m}^{e}_{kl}(z)$, $k,l \! = \! 1,2$, satisfy the jump relation (cf. Lemma~4.5) 
$\mathfrak{m}^{e}_{+}(z) \! = \! \mathfrak{m}^{e}_{-}(z)(\exp (-\mi n \Omega_{
j}^{e}) \sigma_{-} \! + \! \exp (\mi n \Omega_{j}^{e}) \sigma_{+})$, via the 
above asymptotic expansion (as $z \! \to \! a_{j}^{e}$, $j \! = \! 1,\dotsc,
N)$ for $\int_{a_{j}^{e}}^{z} \omega_{m}^{e}$, $m \! = \! 1,\dotsc,N$, one 
arrives at
\begin{align*}
\mathfrak{m}^{e}_{11}(z) \underset{\underset{j=1,\dotsc,N}{z \to a_{j}^{e}}}{
=}& \, \varkappa_{1}^{e}(a_{j}^{e}) \! \left(1 \! + \! \mi \aleph^{1}_{1}(a_{
j}^{e})(z \! - \! a_{j}^{e})^{1/2} \! + \! \daleth^{1}_{1}(a_{j}^{e})(z \! -
\! a_{j}^{e}) \! + \! \mi \beth^{1}_{1}(a_{j}^{e})(z \! - \! a_{j}^{e})^{3/2}
\! + \! \gimel^{1}_{1}(a_{j}^{e})(z \! - \! a_{j}^{e})^{2} \right. \\
+&\left. \, \mathcal{O}((z \! - \! a_{j}^{e})^{5/2}) \right), \\
\mathfrak{m}^{e}_{12}(z) \underset{\underset{j=1,\dotsc,N}{z \to a_{j}^{e}}}{
=}& \, \varkappa_{1}^{e}(a_{j}^{e}) \! \left(1 \! - \! \mi \aleph^{-1}_{1}
(a_{j}^{e})(z \! - \! a_{j}^{e})^{1/2} \! + \! \daleth^{-1}_{1}(a_{j}^{e})(z
\! - \! a_{j}^{e}) \! - \! \mi \beth^{-1}_{1}(a_{j}^{e})(z \! - \! a_{j}^{e}
)^{3/2} \! + \! \gimel^{-1}_{1}(a_{j}^{e})(z \! - \! a_{j}^{e})^{2} \right. \\
+&\left. \, \mathcal{O}((z \! - \! a_{j}^{e})^{5/2}) \right) \! \exp \! 
\left(\mi n \Omega_{j}^{e} \right), \\
\mathfrak{m}^{e}_{21}(z) \underset{\underset{j=1,\dotsc,N}{z \to a_{j}^{e}}}{
=}& \, \varkappa_{2}^{e}(a_{j}^{e}) \! \left(1 \! + \! \mi \aleph^{1}_{-1}(a_{
j}^{e})(z \! - \! a_{j}^{e})^{1/2} \! + \! \daleth^{1}_{-1}(a_{j}^{e})(z \! -
\! a_{j}^{e}) \! + \! \mi \beth^{1}_{-1}(a_{j}^{e})(z \! - \! a_{j}^{e})^{3/2}
\! + \! \gimel^{1}_{-1}(a_{j}^{e})(z \! - \! a_{j}^{e})^{2} \right. \\
+&\left. \, \mathcal{O}((z \! - \! a_{j}^{e})^{5/2}) \right), \\
\mathfrak{m}^{e}_{22}(z) \underset{\underset{j=1,\dotsc,N}{z \to a_{j}^{e}}}{
=}& \, \varkappa_{2}^{e}(a_{j}^{e}) \! \left(1 \! - \! \mi \aleph^{-1}_{-1}
(a_{j}^{e})(z \! - \! a_{j}^{e})^{1/2} \! + \! \daleth^{-1}_{-1}(a_{j}^{e})
(z \! - \! a_{j}^{e}) \! - \! \mi \beth^{-1}_{-1}(a_{j}^{e})(z \! - \! a_{
j}^{e})^{3/2} \! + \! \gimel^{-1}_{-1}(a_{j}^{e})(z \! - \! a_{j}^{e})^{2}
\right. \\
+&\left. \, \mathcal{O}((z \! - \! a_{j}^{e})^{5/2}) \right) \! \exp \! 
\left(\mi n \Omega_{j}^{e} \right),
\end{align*}
where, for $\varepsilon_{1},\varepsilon_{2} \! = \! \pm 1$,
\begin{align*}
\varkappa_{1}^{e}(\xi)=& \, \dfrac{\bm{\theta}^{e}(\bm{u}^{e}_{+}(\infty) \! 
+ \! \bm{d}_{e}) \bm{\theta}^{e}(\bm{u}^{e}_{+}(\xi) \! - \! \frac{n}{2 \pi} 
\bm{\Omega}^{e} \! + \! \bm{d}_{e})}{\bm{\theta}^{e}(\bm{u}^{e}_{+}(\infty) 
\! - \! \frac{n}{2 \pi} \bm{\Omega}^{e} \! + \! \bm{d}_{e}) \bm{\theta}^{e}
(\bm{u}^{e}_{+}(\xi) \! + \! \bm{d}_{e})}, \\
\varkappa_{2}^{e}(\xi)=& \, \dfrac{\bm{\theta}^{e}(-\bm{u}^{e}_{+}(\infty) \! 
- \! \bm{d}_{e}) \bm{\theta}^{e}(\bm{u}^{e}_{+}(\xi) \! - \! \frac{n}{2 \pi} 
\bm{\Omega}^{e} \! - \! \bm{d}_{e})}{\bm{\theta}^{e}(-\bm{u}^{e}_{+}(\infty) 
\! - \! \frac{n}{2 \pi} \bm{\Omega}^{e} \! - \! \bm{d}_{e}) \bm{\theta}^{e}
(\bm{u}^{e}_{+}(\xi) \! - \! \bm{d}_{e})}, \\
\aleph^{\varepsilon_{1}}_{\varepsilon_{2}}(\xi)=& \, -\dfrac{\mathfrak{u}^{e}
(\varepsilon_{1},\varepsilon_{2},\bm{0};\xi)}{\bm{\theta}^{e}(\varepsilon_{1} 
\bm{u}^{e}_{+}(\xi) \! + \! \varepsilon_{2} \bm{d}_{e})} \! + \! \dfrac{
\mathfrak{u}^{e}(\varepsilon_{1},\varepsilon_{2},\bm{\Omega}^{e};\xi)}{\bm{
\theta}^{e}(\varepsilon_{1} \bm{u}^{e}_{+}(\xi) \! - \! \frac{n}{2 \pi} 
\bm{\Omega}^{e} \! + \! \varepsilon_{2} \bm{d}_{e})}, \\
\daleth^{\varepsilon_{1}}_{\varepsilon_{2}}(\xi)=& \, -\dfrac{\mathfrak{v}^{
e}(\varepsilon_{1},\varepsilon_{2},\bm{0};\xi)}{\bm{\theta}^{e}(\varepsilon_{
1} \bm{u}^{e}_{+}(\xi) \! + \! \varepsilon_{2} \bm{d}_{e})} \! + \! \dfrac{
\mathfrak{v}^{e}(\varepsilon_{1},\varepsilon_{2},\bm{\Omega}^{e};\xi)}{\bm{
\theta}^{e}(\varepsilon_{1} \bm{u}^{e}_{+}(\xi) \! - \! \frac{n}{2 \pi} 
\bm{\Omega}^{e} \! + \! \varepsilon_{2} \bm{d}_{e})} \! - \! \left(\dfrac{
\mathfrak{u}^{e}(\varepsilon_{1},\varepsilon_{2},\bm{0};\xi)}{\bm{\theta}^{e}
(\varepsilon_{1} \bm{u}^{e}_{+}(\xi) \! + \! \varepsilon_{2} \bm{d}_{e})} 
\right)^{2} \\
+& \, \dfrac{\mathfrak{u}^{e}(\varepsilon_{1},\varepsilon_{2},\bm{0};\xi) 
\mathfrak{u}^{e}(\varepsilon_{1},\varepsilon_{2},\bm{\Omega}^{e};\xi)}{\bm{
\theta}^{e}(\varepsilon_{1} \bm{u}^{e}_{+}(\xi) \! + \! \varepsilon_{2} 
\bm{d}_{e}) \bm{\theta}^{e}(\varepsilon_{1} \bm{u}^{e}_{+}(\xi) \! - \! 
\frac{n}{2 \pi} \bm{\Omega}^{e} \! + \! \varepsilon_{2} \bm{d}_{e})}, \\
\beth^{\varepsilon_{1}}_{\varepsilon_{2}}(\xi)=& \, -\dfrac{\mathfrak{w}^{
e}(\varepsilon_{1},\varepsilon_{2},\bm{0};\xi)}{\bm{\theta}^{e}(\varepsilon_{
1} \bm{u}^{e}_{+}(\xi) \! + \! \varepsilon_{2} \bm{d}_{e})} \! + \! \dfrac{
\mathfrak{w}^{e}(\varepsilon_{1},\varepsilon_{2},\bm{\Omega}^{e};\xi)}{\bm{
\theta}^{e}(\varepsilon_{1} \bm{u}^{e}_{+}(\xi) \! - \! \frac{n}{2 \pi} \bm{
\Omega}^{e} \! + \! \varepsilon_{2} \bm{d}_{e})} \! + \! \dfrac{2 \mathfrak{
u}^{e}(\varepsilon_{1},\varepsilon_{2},\bm{0};\xi) \mathfrak{v}^{e}
(\varepsilon_{1},\varepsilon_{2},\bm{0};\xi)}{(\bm{\theta}^{e}(\varepsilon_{
1} \bm{u}^{e}_{+}(\xi) \! + \! \varepsilon_{2} \bm{d}_{e}))^{2}} \\
-& \, \dfrac{\mathfrak{v}^{e}(\varepsilon_{1},\varepsilon_{2},\bm{0};\xi) 
\mathfrak{u}^{e}(\varepsilon_{1},\varepsilon_{2},\bm{\Omega}^{e};\xi)}{\bm{
\theta}^{e}(\varepsilon_{1} \bm{u}^{e}_{+}(\xi) \! + \! \varepsilon_{2} \bm{
d}_{e}) \bm{\theta}^{e}(\varepsilon_{1} \bm{u}^{e}_{+}(\xi) \! - \! \frac{
n}{2 \pi} \bm{\Omega}^{e} \! + \! \varepsilon_{2} \bm{d}_{e})} \! + \! \left(
\dfrac{\mathfrak{u}^{e}(\varepsilon_{1},\varepsilon_{2},\bm{0};\xi)}{\bm{
\theta}^{e}(\varepsilon_{1} \bm{u}^{e}_{+}(\xi) \! + \! \varepsilon_{2} \bm{
d}_{e})} \right)^{3} \\
-& \, \dfrac{\mathfrak{u}^{e}(\varepsilon_{1},\varepsilon_{2},\bm{0};\xi)
\mathfrak{v}^{e}(\varepsilon_{1},\varepsilon_{2},\bm{\Omega}^{e};\xi)}{\bm{
\theta}^{e}(\varepsilon_{1} \bm{u}^{e}_{+}(\xi) \! + \! \varepsilon_{2} \bm{
d}_{e}) \bm{\theta}^{e}(\varepsilon_{1} \bm{u}^{e}_{+}(\xi) \! - \! \frac{n}{
2 \pi} \bm{\Omega}^{e} \! + \! \varepsilon_{2} \bm{d}_{e})} \! - \! \dfrac{
(\mathfrak{u}^{e}(\varepsilon_{1},\varepsilon_{2},\bm{0};\xi))^{2}}{(\bm{
\theta}^{e}(\varepsilon_{1} \bm{u}^{e}_{+}(\xi) \! + \! \varepsilon_{2} \bm{
d}_{e}))^{2}} \\
\times& \, \dfrac{\mathfrak{u}^{e}(\varepsilon_{1},\varepsilon_{2},\bm{
\Omega}^{e};\xi)}{\bm{\theta}^{e}(\varepsilon_{1} \bm{u}^{e}_{+}(\xi) \! - 
\! \frac{n}{2 \pi} \bm{\Omega}^{e} \! + \! \varepsilon_{2} \bm{d}_{e})}, \\
\gimel^{\varepsilon_{1}}_{\varepsilon_{2}}(\xi)=& \, -\dfrac{\mathfrak{z}^{
e}(\varepsilon_{1},\varepsilon_{2},\bm{0};\xi)}{\bm{\theta}^{e}(\varepsilon_{
1} \bm{u}^{e}_{+}(\xi) \! + \! \varepsilon_{2} \bm{d}_{e})} \! + \! \dfrac{
\mathfrak{z}^{e}(\varepsilon_{1},\varepsilon_{2},\bm{\Omega}^{e};\xi)}{\bm{
\theta}^{e}(\varepsilon_{1} \bm{u}^{e}_{+}(\xi) \! - \! \frac{n}{2 \pi} 
\bm{\Omega}^{e} \! + \! \varepsilon_{2} \bm{d}_{e})} \! + \! \left(\dfrac{
\mathfrak{v}^{e}(\varepsilon_{1},\varepsilon_{2},\bm{0};\xi)}{\bm{\theta}^{e}
(\varepsilon_{1} \bm{u}^{e}_{+}(\xi) \! + \! \varepsilon_{2} \bm{d}_{e})} 
\right)^{2} \\
-& \, \dfrac{\mathfrak{v}^{e}(\varepsilon_{1},\varepsilon_{2},\bm{0};\xi) 
\mathfrak{v}^{e}(\varepsilon_{1},\varepsilon_{2},\bm{\Omega}^{e};\xi)}{\bm{
\theta}^{e}(\varepsilon_{1} \bm{u}^{e}_{+}(\xi) \! + \! \varepsilon_{2} 
\bm{d}_{e}) \bm{\theta}^{e}(\varepsilon_{1} \bm{u}^{e}_{+}(\xi) \! - \! 
\frac{n}{2 \pi} \bm{\Omega}^{e} \! + \! \varepsilon_{2} \bm{d}_{e})} \! - 
\! \dfrac{2 \mathfrak{u}^{e}(\varepsilon_{1},\varepsilon_{2},\bm{0};\xi) 
\mathfrak{w}^{e}(\varepsilon_{1},\varepsilon_{2},\bm{0};\xi)}{(\bm{\theta}^{e}
(\varepsilon_{1} \bm{u}^{e}_{+}(\xi) \! + \! \varepsilon_{2} \bm{d}_{e}))^{
2}} \\
+& \, \dfrac{\mathfrak{w}^{e}(\varepsilon_{1},\varepsilon_{2},\bm{0};\xi) 
\mathfrak{u}^{e}(\varepsilon_{1},\varepsilon_{2},\bm{\Omega}^{e};\xi)}{\bm{
\theta}^{e}(\varepsilon_{1} \bm{u}^{e}_{+}(\xi) \! + \! \varepsilon_{2} \bm{
d}_{e}) \bm{\theta}^{e}(\varepsilon_{1} \bm{u}^{e}_{+}(\xi) \! - \! \frac{
n}{2 \pi} \bm{\Omega}^{e} \! + \! \varepsilon_{2} \bm{d}_{e})} \! + \! \dfrac{
3(\mathfrak{u}^{e}(\varepsilon_{1},\varepsilon_{2},\bm{0};\xi))^{2} \mathfrak{
v}^{e}(\varepsilon_{1},\varepsilon_{2},\bm{0};\xi)}{(\bm{\theta}^{e}
(\varepsilon_{1} \bm{u}^{e}_{+}(\xi) \! + \! \varepsilon_{2} \bm{d}_{e}))^{
3}} \\
+& \, \dfrac{\mathfrak{u}^{e}(\varepsilon_{1},\varepsilon_{2},\bm{0};\xi) 
\mathfrak{w}^{e}(\varepsilon_{1},\varepsilon_{2},\bm{\Omega}^{e};\xi)}{\bm{
\theta}^{e}(\varepsilon_{1} \bm{u}^{e}_{+}(\xi) \! + \! \varepsilon_{2} \bm{
d}_{e}) \bm{\theta}^{e}(\varepsilon_{1} \bm{u}^{e}_{+}(\xi) \! - \! \frac{n}{
2 \pi} \bm{\Omega}^{e} \! + \! \varepsilon_{2} \bm{d}_{e})} \! + \! \left(
\dfrac{\mathfrak{u}^{e}(\varepsilon_{1},\varepsilon_{2},\bm{0};\xi)}{\bm{
\theta}^{e}(\varepsilon_{1} \bm{u}^{e}_{+}(\xi) \! + \! \varepsilon_{2} 
\bm{d}_{e})} \right)^{4} \\
-& \, \dfrac{2 \mathfrak{u}^{e}(\varepsilon_{1},\varepsilon_{2},\bm{0};\xi)
\mathfrak{v}^{e}(\varepsilon_{1},\varepsilon_{2},\bm{0};\xi) \mathfrak{u}^{
e}(\varepsilon_{1},\varepsilon_{2},\bm{\Omega}^{e};\xi)}{(\bm{\theta}^{e}
(\varepsilon_{1} \bm{u}^{e}_{+}(\xi) \! + \! \varepsilon_{2} \bm{d}_{e}))^{2}
\bm{\theta}^{e}(\varepsilon_{1} \bm{u}^{e}_{+}(\xi) \! - \! \frac{n}{2 \pi}
\bm{\Omega}^{e} \! + \! \varepsilon_{2} \bm{d}_{e})} \! - \! \dfrac{
(\mathfrak{u}^{e}(\varepsilon_{1},\varepsilon_{2},\bm{0};\xi))^{2}}{(\bm{
\theta}^{e}(\varepsilon_{1} \bm{u}^{e}_{+}(\xi) \! + \! \varepsilon_{2} 
\bm{d}_{e}))^{2}} \\
\times& \, \dfrac{\mathfrak{v}^{e}(\varepsilon_{1},\varepsilon_{2},\bm{
\Omega}^{e};\xi)}{\bm{\theta}^{e}(\varepsilon_{1} \bm{u}^{e}_{+}(\xi) \! - \! 
\frac{n}{2\pi} \bm{\Omega}^{e} \! + \! \varepsilon_{2} \bm{d}_{e})} \! - \! 
\dfrac{(\mathfrak{u}^{e}(\varepsilon_{1},\varepsilon_{2},\bm{0};\xi))^{3} 
\mathfrak{u}^{e}(\varepsilon_{1},\varepsilon_{2},\bm{\Omega}^{e};\xi)}{(\bm{
\theta}^{e}(\varepsilon_{1} \bm{u}^{e}_{+}(\xi) \! + \! \varepsilon_{2} 
\bm{d}_{e}))^{3} \bm{\theta}^{e}(\varepsilon_{1} \bm{u}^{e}_{+}(\xi) \! - 
\! \frac{n}{2 \pi} \bm{\Omega}^{e} \! + \! \varepsilon_{2} \bm{d}_{e})},
\end{align*}
with $\bm{0} \! := \! (0,0,\dotsc,0)^{\operatorname{T}}$ $(\in \! \mathbb{
R}^{N})$,
\begin{gather*}
\mathfrak{u}^{e}(\varepsilon_{1},\varepsilon_{2},\bm{\Omega}^{e};\xi) \! 
:= \! 2 \pi \overset{e}{\Lambda}^{\raise-1.0ex\hbox{$\scriptstyle 1$}}_{0}
(\varepsilon_{1},\varepsilon_{2},\bm{\Omega}^{e};\xi), \qquad \qquad 
\mathfrak{v}^{e}(\varepsilon_{1},\varepsilon_{2},\bm{\Omega}^{e};\xi) \! := 
\! -2 \pi^{2} \overset{e}{\Lambda}^{\raise-1.0ex\hbox{$\scriptstyle 2$}}_{0}
(\varepsilon_{1},\varepsilon_{2},\bm{\Omega}^{e};\xi), \\
\mathfrak{w}^{e}(\varepsilon_{1},\varepsilon_{2},\bm{\Omega}^{e};\xi) \! := 
\! 2 \pi \! \left(\overset{e}{\Lambda}^{\raise-1.0ex\hbox{$\scriptstyle 0$}}_{
1}(\varepsilon_{1},\varepsilon_{2},\bm{\Omega}^{e};\xi) \! - \! \dfrac{2 \pi^{
2}}{3} \overset{e}{\Lambda}^{\raise-1.0ex\hbox{$\scriptstyle 3$}}_{0}
(\varepsilon_{1},\varepsilon_{2},\bm{\Omega}^{e};\xi) \right), \\
\mathfrak{z}^{e}(\varepsilon_{1},\varepsilon_{2},\bm{\Omega}^{e};\xi) \! := \! 
-(2 \pi)^{2} \! \left(\overset{e}{\Lambda}^{\raise-1.0ex\hbox{$\scriptstyle 
1$}}_{1}(\varepsilon_{1},\varepsilon_{2},\bm{\Omega}^{e};\xi) \! - \! \dfrac{
\pi^{2}}{6} \overset{e}{\Lambda}^{\raise-1.0ex\hbox{$\scriptstyle 4$}}_{0}
(\varepsilon_{1},\varepsilon_{2},\bm{\Omega}^{e};\xi) \right), \\
\overset{e}{\Lambda}^{\raise-1.0ex\hbox{$\scriptstyle j_{1}$}}_{j_{2}}
(\varepsilon_{1},\varepsilon_{2},\bm{\Omega}^{e};\xi) \! = \! \sum_{m \in
\mathbb{Z}^{N}}(\mathfrak{r}_{e}(\xi))^{j_{1}}(\mathfrak{s}_{e}(\xi))^{j_{2}}
\, \me^{2 \pi \mi (m,\varepsilon_{1} \bm{u}^{e}_{+}(\xi)-\frac{n}{2 \pi}
\bm{\Omega}^{e}+ \varepsilon_{2} \bm{d}_{e})+ \pi \mi (m,\bm{\tau}^{e}m)}, \\
\mathfrak{r}_{e}(\xi) \! := \! \dfrac{2(m,\vec{\moo}^{e}_{1}(\xi))}{\mathfrak{
f}_{e}(\xi)}, \qquad \qquad \mathfrak{s}_{e}(\xi) \! := \! \dfrac{2(m,\vec{
\moo}^{e}_{2}(\xi))}{3 \mathfrak{f}_{e}(\xi)}, \\
\vec{\moo}^{e}_{1}(\xi) \! = \! \left(\mathfrak{p}_{1}^{\natural}(\xi),
\mathfrak{p}_{2}^{\natural}(\xi),\dotsc,\mathfrak{p}_{N}^{\natural}(\xi) 
\right), \qquad \qquad \vec{\moo}^{e}_{2}(\xi) \! = \! \left(\mathfrak{q}_{
1}^{\natural}(\xi),\mathfrak{q}_{2}^{\natural}(\xi),\dotsc,\mathfrak{q}_{N}^{
\natural}(\xi) \right).
\end{gather*}

Recall the definition of $\gamma^{e}(z)$ given in Lemma~4.4: a careful 
analysis of the branch cuts shows that, for $j \! = \! 1,\dotsc,N$,
\begin{align*}
(\gamma^{e}(z))^{2} &\underset{z \in \mathbb{C}_{\pm}}{=} \pm \, \dfrac{\left(
\! (z \! - \! b_{j}^{e}) \prod_{\substack{k=1\\k \not= j}}^{N} \! \left(\!
\frac{z-b_{k}^{e}}{z-a_{k}^{e}} \right) \! \left(\! \frac{z-b_{0}^{e}}{z-a_{N
+1}^{e}} \right) \right)^{1/2}}{\sqrt{\smash[b]{z \! - \! a_{j}^{e}}}} \\
&\underset{\mathbb{C}_{\pm} \ni z \to a_{j}^{e}}{=} \pm \, \dfrac{\left(
Q_{0}^{e}(a_{j}^{e}) \! + \! Q_{1}^{e}(a_{j}^{e})(z \! - \! a_{j}^{e}) \!
+ \! \frac{1}{2}Q_{2}^{e}(a_{j}^{e})(z \! - \! a_{j}^{e})^{2} \! + \!
\mathcal{O}((z \! - \! a_{j}^{e})^{3}) \right)}{\sqrt{\smash[b]{z \! - \!
a_{j}^{e}}}},
\end{align*}
where $Q_{0}^{e}(a_{j}^{e}),Q_{1}^{e}(a_{j}^{e})$, $j \! = \! 1,\dotsc,N$,
are given in Theorem~2.3.1, Equations~(2.35) and~(2.36), and
\begin{align*}
Q_{2}^{e}(a_{j}^{e})=& -\dfrac{1}{2}Q_{0}^{e}(a_{j}^{e}) \! \left(\sum_{
\substack{k=1\\k \not= j}}^{N} \! \left(\dfrac{1}{(a_{j}^{e} \! - \! b_{k}^{e}
)^{2}} \! - \! \dfrac{1}{(a_{j}^{e} \! - \! a_{k}^{e})^{2}} \right) \! + \!
\dfrac{1}{(a_{j}^{e} \! - \! b_{0}^{e})^{2}} \! - \! \dfrac{1}{(a_{j}^{e} \!
- \! a_{N+1}^{e})^{2}} \! + \! \dfrac{1}{(a_{j}^{e} \! - \! b_{j}^{e})^{2}}
\right) \\
+& \, \dfrac{1}{4}Q_{0}^{e}(a_{j}^{e}) \! \left(\sum_{\substack{k=1\\k \not=
j}}^{N} \! \left(\dfrac{1}{a_{j}^{e} \! - \! b_{k}^{e}} \! - \! \dfrac{1}{a_{
j}^{e} \! - \! a_{k}^{e}} \right) \! + \! \dfrac{1}{a_{j}^{e} \! - \! b_{0}^{
e}} \! - \! \dfrac{1}{a_{j}^{e} \! - \! a_{N+1}^{e}} \! + \! \dfrac{1}{a_{j}^{
e} \! - \! b_{j}^{e}} \right)^{2}, \quad j \! = \! 1,\dotsc,N.
\end{align*}
Recall the above formula for 
$\overset{e}{\mathfrak{M}}^{\raise-1.0ex\hbox{$\scriptstyle \infty$}}(z) \!
\left(
\begin{smallmatrix}
\mp (s_{1}+t_{1}) & \pm \mi (s_{1}-t_{1}) \me^{\mi n \Omega_{j}^{e}} \\
\pm \mi (s_{1}-t_{1}) \me^{-\mi n \Omega_{j}^{e}} & \pm (s_{1}+t_{1})
\end{smallmatrix}
\right) \!
(\overset{e}{\mathfrak{M}}^{\raise-1.0ex\hbox{$\scriptstyle \infty$}}(z))^{-
1}$. Substituting the a\-b\-o\-v\-e expansions (as $z \! \to \! a_{j}^{e}$, 
$j \! = \! 1,\dotsc,N)$ into this formula, equating coefficients of like 
powers of $(z \! - \! a_{j}^{e})^{-p/2}(nG_{a_{j}}^{e}(z))^{-1}$, $p \! \in \! 
\{4,3,2,1,0\}$, and considering, say, the $(1 \, 1)$-element of the resulting 
(asymptotic) expansions, one arrives at, up to terms that are $\mathcal{O} 
\! \left((n^{2}(z \! - \! a_{j}^{e})^{3}(G_{a_{j}}^{e}(z))^{2})^{-1} 
\overset{e}{\mathfrak{M}}^{\raise-1.0ex\hbox{$\scriptstyle \infty$}}(z)
f_{a_{j}}^{e}(n)(\overset{e}{\mathfrak{M}}^{\raise-1.0ex\hbox{$\scriptstyle 
\infty$}}(z))^{-1} \right)$ (modulo a minus sign, this result is equally 
applicable to the $(2 \, 2)$-element, since $\operatorname{tr}(w_{+}^{
\Sigma^{e}_{\circlearrowright}}(z)) \! = \! 0)$, upon setting, for economy 
of notation, $Q_{q}^{e}(a_{j}^{e}) \! =: \! Q_{q}$, $q \! = \! 0,1,2$, 
$\varkappa_{1}^{e}(a_{j}^{e}) \! =: \! \varkappa_{1}^{e}$, $\varkappa_{2}^{
e}(a_{j}^{e}) \! =: \! \varkappa_{2}^{e}$, $\aleph^{\varepsilon_{1}}_{
\varepsilon_{2}}(a_{j}^{e}) \! =: \! \aleph^{\varepsilon_{1}}_{\varepsilon_{
2}}$, $\daleth^{\varepsilon_{1}}_{\varepsilon_{2}}(a_{j}^{e}) \! =: \! 
\daleth^{\varepsilon_{1}}_{\varepsilon_{2}}$, $\beth^{\varepsilon_{1}}_{
\varepsilon_{2}}(a_{j}^{e}) \! =: \! \beth^{\varepsilon_{1}}_{\varepsilon_{
2}}$, and $\gimel^{\varepsilon_{1}}_{\varepsilon_{2}}(a_{j}^{e}) \! =: \! 
\gimel^{\varepsilon_{1}}_{\varepsilon_{2}}$:
\begin{align*}
\mathcal{O} \! \left(\dfrac{(z \! - \! a_{j}^{e})^{-2} \, \me^{\mi n \Omega_{
j}^{e}}}{nG_{a_{j}}^{e}(z)} \right):& -\dfrac{(s_{1} \! + \! t_{1}) 
\varkappa_{1}^{e} \varkappa_{2}^{e}Q_{0}}{4} \! - \! \dfrac{(s_{1} \! + \! 
t_{1}) \varkappa_{1}^{e} \varkappa_{2}^{e}Q_{0}}{ 4} \! - \! \dfrac{(s_{1} 
\! - \! t_{1}) \varkappa_{1}^{e} \varkappa_{2}^{e}Q_{0}}{4} \! - \! \dfrac{
(s_{1} \! - \! t_{1}) \varkappa_{1}^{e} \varkappa_{2}^{e}Q_{0}}{4}; \\
\mathcal{O} \! \left(\dfrac{(z \! - \! a_{j}^{e})^{-3/2} \, \me^{\mi n 
\Omega_{j}^{e}}}{nG_{a_{j}}^{e}(z)} \right):& -\dfrac{\mi (s_{1} \! + \! 
t_{1}) \varkappa_{1}^{e} \varkappa_{2}^{e}Q_{0}(\aleph_{1}^{1} \! - \! 
\aleph_{-1}^{-1})}{4} \! - \! \dfrac{\mi (s_{1} \! + \! t_{1}) \varkappa_{
1}^{e} \varkappa_{2}^{e}Q_{0}(\aleph^{1}_{-1} \! - \! \aleph^{-1}_{1})}{4} \\
&-\dfrac{\mi (s_{1} \! - \! t_{1}) \varkappa_{1}^{e} \varkappa_{2}^{e}Q_{0}
(\aleph^{1}_{-1} \! + \! \aleph^{1}_{1})}{4} \! + \! \dfrac{\mi (s_{1} \! - 
\! t_{1}) \varkappa_{1}^{e} \varkappa_{2}^{e}Q_{0}(\aleph^{-1}_{-1} \! + \! 
\aleph^{-1}_{1})}{4} \\
&-\dfrac{(s_{1} \! + \! t_{1}) \varkappa_{1}^{e} \varkappa_{2}^{e}}{2} \! + 
\! \dfrac{(s_{1} \! + \! t_{1}) \varkappa_{1}^{e} \varkappa_{2}^{e}}{2}; \\
\mathcal{O} \! \left(\dfrac{(z \! - \! a_{j}^{e})^{-1} \, \me^{\mi n 
\Omega_{j}^{e}}}{nG_{a_{j}}^{e}(z)} \right):& -\dfrac{(s_{1} \! + \! t_{1}) 
\varkappa_{1}^{e} \varkappa_{2}^{e}}{4} \! \left(Q_{1} \! + \! Q_{0} \! 
\left(\daleth^{-1}_{-1} \! + \! \daleth^{1}_{1} \! + \! \aleph^{1}_{1} 
\aleph^{-1}_{-1} \right) \right) \! - \! \dfrac{(s_{1} \! + \! t_{1}) 
\varkappa_{1}^{e} \varkappa_{2}^{e}}{4Q_{0}} \\
&-\dfrac{(s_{1} \! + \! t_{1}) \varkappa_{1}^{e} \varkappa_{2}^{e}}{4} \! 
\left(Q_{1} \! + \! Q_{0} \! \left(\daleth^{1}_{-1} \! + \! \daleth^{-1}_{1} 
\! + \! \aleph^{-1}_{1} \aleph^{1}_{-1} \right) \right) \! - \! \dfrac{(s_{1} 
\! + \! t_{1}) \varkappa_{1}^{e} \varkappa_{2}^{e}}{4Q_{0}} \\
&-\dfrac{(s_{1} \! - \! t_{1}) \varkappa_{1}^{e} \varkappa_{2}^{e}}{4} \! 
\left(Q_{1} \! + \! Q_{0} \! \left(\daleth^{1}_{-1} \! + \! \daleth^{1}_{1} 
\! - \! \aleph^{1}_{1} \aleph^{1}_{-1} \right) \right) \! + \! \dfrac{(s_{1} 
\! - \! t_{1}) \varkappa_{1}^{e} \varkappa_{2}^{e}}{4Q_{0}} \\
&-\dfrac{(s_{1} \! - \! t_{1}) \varkappa_{1}^{e} \varkappa_{2}^{e}}{4} \! 
\left(Q_{1} \! + \! Q_{0} \! \left(\daleth^{-1}_{-1} \! + \! \daleth^{-1}_{
1} \! - \! \aleph^{-1}_{1} \aleph^{-1}_{-1} \right) \right) \! + \! \dfrac{
(s_{1} \! - \! t_{1}) \varkappa_{1}^{e} \varkappa_{2}^{e}}{4Q_{0}} \\
&-\dfrac{\mi (s_{1} \! + \! t_{1}) \varkappa_{1}^{e} \varkappa_{2}^{e}}{2} 
\! \left(\aleph^{1}_{1} \! - \! \aleph^{-1}_{-1} \right) \! + \! \dfrac{\mi 
(s_{1} \! + \! t_{1}) \varkappa_{1}^{e} \varkappa_{2}^{e}}{2} \! \left(
\aleph^{1}_{-1} \! - \! \aleph^{-1}_{1} \right); \\
\mathcal{O} \! \left(\dfrac{(z \! - \! a_{j}^{e})^{-1/2} \, \me^{\mi n 
\Omega_{j}^{e}}}{nG_{a_{j}}^{e}(z)} \right):& -\dfrac{\mi (s_{1} \! + \! 
t_{1}) \varkappa_{1}^{e} \varkappa_{2}^{e}}{4} \! \left(Q_{1} \! \left(
\aleph^{1}_{1} \! - \! \aleph^{-1}_{-1} \right) \! + \! Q_{0} \! \left(
\beth^{1}_{1} \! - \! \beth^{-1}_{-1} \! + \! \aleph^{1}_{1} \daleth^{-1}_{
-1} \! - \! \aleph^{-1}_{-1} \daleth^{1}_{1} \right) \right) \\
&-\dfrac{\mi (s_{1} \! + \! t_{1}) \varkappa_{1}^{e} \varkappa_{2}^{e}}{4} 
\! \left(Q_{1} \! \left(\aleph^{1}_{-1} \! - \! \aleph^{-1}_{1} \right) \! + 
\! Q_{0} \! \left(\beth^{1}_{-1} \! - \! \beth^{-1}_{1} \! + \! \aleph^{1}_{
-1} \daleth^{-1}_{1} \! - \! \aleph^{-1}_{1} \daleth^{1}_{-1} \right) \right) 
\\
&-\dfrac{\mi (s_{1} \! - \! t_{1}) \varkappa_{1}^{e} \varkappa_{2}^{e}}{4} \! 
\left(Q_{1} \! \left(\aleph^{1}_{-1} \! + \! \aleph^{1}_{1} \right) \! + \! 
Q_{0} \! \left(\beth^{1}_{-1} \! + \! \beth^{1}_{1} \! + \! \aleph^{1}_{1} 
\daleth^{1}_{-1} \! + \! \aleph^{1}_{-1} \daleth^{1}_{1} \right) \right) \\
&+\dfrac{\mi (s_{1} \! - \! t_{1}) \varkappa_{1}^{e} \varkappa_{2}^{e}}{4} \! 
\left(Q_{1} \! \left(\aleph^{-1}_{-1} \! + \! \aleph^{-1}_{1} \right) \! + \! 
Q_{0} \! \left(\beth^{-1}_{-1} \! + \! \beth^{-1}_{1} \! + \! \aleph^{-1}_{1} 
\daleth^{-1}_{-1} \! + \! \aleph^{-1}_{-1} \daleth^{-1}_{1} \right) \right) \\
&-\dfrac{\mi (s_{1} \! + \! t_{1}) \varkappa_{1}^{e} \varkappa_{2}^{e}}{4Q_{
0}} \! \left(\aleph^{1}_{1} \! - \! \aleph^{-1}_{-1} \right) \! - \! \dfrac{
(s_{1} \! + \! t_{1}) \varkappa_{1}^{e} \varkappa_{2}^{e}}{2} \! \left(
\daleth^{-1}_{-1} \! + \! \daleth^{1}_{1} \! + \! \aleph^{1}_{1} \aleph^{-
1}_{-1} \right) \\
&-\dfrac{\mi (s_{1} \! + \! t_{1}) \varkappa_{1}^{e} \varkappa_{2}^{e}}{4Q_{
0}} \! \left(\aleph^{1}_{-1} \! - \! \aleph^{-1}_{1} \right) \! + \! \dfrac{
(s_{1} \! + \! t_{1}) \varkappa_{1}^{e} \varkappa_{2}^{e}}{2} \! \left(
\daleth^{1}_{-1} \! + \! \daleth^{-1}_{1} \! + \! \aleph^{-1}_{1} \aleph^{
1}_{-1} \right) \\
&+\dfrac{\mi (s_{1} \! - \! t_{1}) \varkappa_{1}^{e} \varkappa_{2}^{e}}{4Q_{
0}} \! \left(\aleph^{1}_{-1} \! + \! \aleph^{1}_{1} \right) \! - \! \dfrac{
\mi (s_{1} \! - \! t_{1}) \varkappa_{1}^{e} \varkappa_{2}^{e}}{4Q_{0}} \! 
\left(\aleph^{-1}_{1} \! + \! \aleph^{-1}_{-1} \right); \\
\mathcal{O} \! \left(\dfrac{\me^{\mi n \Omega_{j}^{e}}}{nG_{a_{j}}^{e}(z)} 
\right):& -\dfrac{(s_{1} \! + \! t_{1}) \varkappa_{1}^{e} \varkappa_{2}^{e}
}{4} \! \left(Q_{0} \! \left(\gimel^{-1}_{-1} \! + \! \gimel^{1}_{1} \! + \! 
\daleth^{1}_{1} \daleth^{-1}_{-1} \! + \! \aleph^{-1}_{-1} \beth^{1}_{1} \! 
+ \! \aleph^{1}_{1} \beth^{-1}_{-1} \right) \! + \! Q_{1} \! \left(\daleth^{
-1}_{-1} \! + \! \daleth^{1}_{1} \! + \! \aleph^{1}_{1} \aleph^{-1}_{-1} 
\right) \right. \\
&\left. + \, \dfrac{1}{2}Q_{2} \right) \! - \! \dfrac{(s_{1} \! + \! t_{1}) 
\varkappa_{1}^{e} \varkappa_{2}^{ e}}{4} \! \left(Q_{0} \! \left(\gimel^{1}_{
-1} \! + \! \gimel^{-1}_{1} \! + \! \daleth^{-1}_{1} \daleth^{1}_{-1} \! + \! 
\aleph^{-1}_{1} \beth^{1}_{-1} \! + \! \aleph^{1}_{-1} \beth^{-1}_{1} \right) 
\! + \! \dfrac{1}{2}Q_{2} \right. \\
&\left. + \, Q_{1} \! \left( \daleth^{1}_{-1} \! + \! \daleth^{-1}_{1} \! + 
\! \aleph^{-1}_{1} \aleph^{1}_{-1} \right) \right) \! - \! \dfrac{(s_{1} \! 
- \! t_{1}) \varkappa_{1}^{e} \varkappa_{2}^{e}}{4} \! \left(Q_{0} \! \left(
\gimel^{1}_{-1} \! + \! \gimel^{1}_{1} \! + \! \daleth^{1}_{1} \daleth^{1}_{
-1} \! - \! \aleph^{1}_{1} \beth^{1}_{-1} \! - \! \aleph^{1}_{-1} \beth^{1}_{
1} \right) \right. \\
&\left. + \, \dfrac{1}{2}Q_{2} \! + \! Q_{1} \! \left(\daleth^{1}_{-1} \! + \! 
\daleth^{1}_{1} \! - \! \aleph^{1}_{1} \aleph^{1}_{-1} \right) \right) \! - 
\! \dfrac{(s_{1} \! - \! t_{1}) \varkappa_{1}^{e} \varkappa_{2}^{ e}}{4} \! 
\left(Q_{0} \! \left(\gimel^{-1}_{-1} \! + \! \gimel^{-1}_{1} \! + \! 
\daleth^{-1}_{1} \daleth^{-1}_{-1} \! - \! \aleph^{-1}_{1} \beth^{-1}_{-1} 
\right. \right. \\
&\left. \left. - \, \aleph^{-1}_{-1} \beth^{-1}_{1} \right) \! + \! \dfrac{1}{
2}Q_{2} \! + \! Q_{1} \! \left( \daleth^{-1}_{-1} \! + \! \daleth^{-1}_{1} \! 
- \! \aleph^{-1}_{1} \aleph^{-1}_{-1} \right) \right) \! - \! \dfrac{(s_{1} 
\! + \! t_{1}) \varkappa_{1}^{e} \varkappa_{2}^{e}}{4Q_{0}} \! \left(
\daleth^{-1}_{-1} \! + \! \daleth^{1}_{1} \! + \! \aleph^{1}_{1} \aleph^{-
1}_{-1} \right. \\
&\left. - \, Q_{1}(Q_{0})^{-1} \right) \! - \! \dfrac{(s_{1} \! + \! t_{1}) 
\varkappa_{1}^{e} \varkappa_{2}^{ e}}{4Q_{0}} \! \left(\daleth^{1}_{-1} \! 
+ \! \daleth^{-1}_{1} \! + \! \aleph^{-1}_{1} \aleph^{1}_{-1} \! - \! Q_{1}
(Q_{0})^{-1} \right) \! - \! \dfrac{\mi (s_{1} \! + \! t_{1}) \varkappa_{
1}^{e} \varkappa_{2}^{e}}{2} \\
&\times \left(\beth^{1}_{1} \! - \! \beth^{-1}_{-1} \! + \! \aleph^{1}_{1} 
\daleth^{-1}_{-1} \! - \! \aleph^{-1}_{ -1} \daleth^{1}_{1} \right) \! + \! 
\dfrac{\mi (s_{1} \! + \! t_{1}) \varkappa_{1}^{e} \varkappa_{2}^{e}}{2} \! 
\left(\beth^{1}_{-1} \! - \! \beth^{-1}_{1} \! + \! \aleph^{1}_{-1} 
\daleth^{-1}_{1} \! - \! \aleph^{-1}_{1} \daleth^{1}_{-1} \right) \\
&+ \dfrac{(s_{1} \! - \! t_{1}) \varkappa_{1}^{e} \varkappa_{2}^{e}}{4Q_{0}} 
\! \left( \daleth^{1}_{-1} \! + \! \daleth^{1}_{1} \! - \! \aleph^{1}_{1} 
\aleph^{1}_{-1} \! - \! Q_{1}(Q_{0})^{-1} \right) \! + \! \dfrac{(s_{1} \! - 
\! t_{1}) \varkappa_{1}^{e} \varkappa_{2}^{ e}}{4Q_{0}} \! \left(\daleth^{-
1}_{-1} \! + \! \daleth^{-1}_{1} \right. \\
&\left. - \, \aleph^{-1}_{1} \aleph^{-1}_{-1} \! - \!Q_{1}(Q_{0})^{-1} 
\right).
\end{align*}
Repeating the above analysis, \emph{mutatis mutandis}, for the $(1 \, 2)$- and 
$(2 \, 1)$-elements, substituting $\aleph^{-1}_{1} \! = \! \aleph^{1}_{1}$, 
$\aleph^{-1}_{-1} \! = \! \aleph^{1}_{-1}$, $\daleth^{-1}_{1} \! = \! 
\daleth^{1}_{1}$, $\daleth^{-1}_{-1} \! = \! \daleth^{1}_{-1}$, $\beth^{-1}_{
1} \! = \! \beth^{1}_{1}$, $\beth^{-1}_{-1} \! = \! \beth^{1}_{-1}$, $\gimel^{
-1}_{1} \! = \! \gimel^{1}_{1}$, and $\gimel^{-1}_{-1} \! = \! \gimel^{1}_{-
1}$ into the above (and resulting) `coefficient equations', and simplifying, 
one shows that: (i) the coefficients of the terms that are $\mathcal{O}((z \! 
- \! a_{j}^{e})^{-p/2}(nG_{a_{j}}^{e}(z))^{-1} \exp (\mi n \Omega_{j}^{e}))$, 
$p \! = \! 1,3$, are equal to zero; and (ii) recalling {}from Lemma~4.7 that, 
for $z \! \in \! \mathbb{U}^{e}_{\delta_{a_{j}}} \setminus (-\infty,a_{j}^{
e})$, $j \! = \! 1, \dotsc,N$, $G_{a_{j}}^{e}(z) \! =_{z \to a_{j}^{e}} \! 
\widehat{\alpha}_{0} \! + \! \widehat{\alpha}_{1}(z \! - \! a_{j}^{e}) \! + 
\! \widehat{\alpha}_{2}(z \! - \! a_{j}^{e})^{2} \! + \! \mathcal{O}((z \! - 
\! a_{j}^{e})^{3})$, where $\widehat{\alpha}_{0} \! = \! \widehat{\alpha}_{
0}^{e}(a_{j}^{e}) \! := \! \tfrac{4}{3}f(a_{j}^{e})$, $\widehat{\alpha}_{1} 
\! = \! \widehat{\alpha}_{1}^{e}(a_{j}^{e}) \! := \! \tfrac{4}{5}f^{\prime}
(a_{j}^{e})$, and $\widehat{ \alpha}_{2} \! = \! \widehat{\alpha}_{2}^{e}
(a_{j}^{e}) \! := \! \tfrac{2}{7}f^{\prime \prime}(a_{j}^{e})$, with $f(a_{
j}^{e})$, $f^{\prime}(a_{j}^{e})$, and $f^{\prime \prime}(a_{j}^{e})$ given 
in Lemma~4.7, substituting the expansion for $G^{e}_{a_{j}}(z)$ (as $z \! 
\to \! a_{j}^{e}$, $j \! = \! 1,\dotsc,N)$ into the remaining non-zero 
coefficient equations, collecting coefficients of like powers of $(z \! - 
\! a_{j}^{e})^{-p}$, $p \! = \! 0,1,2$, and continuing, analytically, the 
resulting (rational) expressions to $\partial \mathbb{U}^{e}_{\delta_{a_{
j}}}$, $j \! = \! 1,\dotsc,N$, one arrives at, after a lengthy calculation 
and reinserting explicit $a_{j}^{e}$, $j \! = \! 1,\dotsc,N$, dependencies,
\begin{align}
w_{+}^{\Sigma^{e}_{\circlearrowright}}(z) \underset{\underset{z \in \partial
\mathbb{U}^{e}_{\delta_{a_{j}}}}{n \to \infty}}{=}& \dfrac{1}{n} \! \left(
\dfrac{\mathscr{A}^{e}(a_{j}^{e})}{\widehat{\alpha}_{0}^{e}(a_{j}^{e})(z \! 
- \! a_{j}^{e})^{2}} \! + \! \dfrac{(\mathscr{B}^{e}(a_{j}^{e}) \widehat{
\alpha}_{0}^{e}(a_{j}^{e}) \! - \! \mathscr{A}^{e}(a_{j}^{e}) \widehat{
\alpha}_{1}^{e}(a_{j}^{e}))}{(\widehat{\alpha}_{0}^{e}(a_{j}^{e}))^{2}(z \! 
- \! a_{j}^{e})} \right. \nonumber \\
&\left. + \, \dfrac{\left(\mathscr{A}^{e}(a_{j}^{e}) \widehat{\alpha}_{0}^{e}
(a_{j}^{e}) \! \left(\! \left(\frac{\widehat{\alpha}_{1}^{e}(a_{j}^{e})}{
\widehat{\alpha}_{0}^{e}(a_{j}^{e})} \right)^{2} \! - \! \frac{\widehat{
\alpha}_{2}^{e}(a_{j}^{e})}{\widehat{\alpha}_{0}^{e}(a_{j}^{e})} \right) \! 
- \! \mathscr{B}^{e}(a_{j}^{e}) \widehat{\alpha}_{1}^{e}(a_{j}^{e}) \! + \!
\mathscr{C}^{e}(a_{j}^{e}) \widehat{\alpha}_{0}^{e}(a_{j}^{e}) \right)}{(
\widehat{\alpha}_{0}^{e}(a_{j}^{e}))^{2}} \right) \nonumber \\
&+ \, \mathcal{O} \! \left(\dfrac{1}{n} \sum_{k=1}^{\infty}f_{k}^{e}(n)(z 
\! - \! a_{j}^{e})^{k} \right) \! + \! \mathcal{O} \! \left(\dfrac{
\overset{e}{\mathfrak{M}}^{\raise-1.0ex\hbox{$\scriptstyle \infty$}}(z)f^{e}_{
a_{j}}(n)(\overset{e}{\mathfrak{M}}^{\raise-1.0ex\hbox{$\scriptstyle \infty$}}
(z))^{-1}}{n^{2}(z \! - \! a_{j}^{e})^{3}(G^{e}_{a_{j}}(z))^{2}} \right),
\end{align}
where $\mathscr{A}^{e}(a_{j}^{e}),\mathscr{B}^{e}(a_{j}^{e})$, $j \! = \! 1,
\dotsc,N$, are defined in Theorem~2.3.1, Equations~(2.25), (2.27), (2.28), 
(2.35)--(2.45), (2.49), (2.56), and~(2.57),
\begin{equation*}
\dfrac{\mathscr{C}^{e}(a_{j}^{e})}{\me^{\mi n \Omega_{j}^{e}}} \! := \!
\begin{pmatrix}
\boxed{\begin{matrix} \varkappa_{1}^{e}(a_{j}^{e}) \varkappa_{2}^{e}(a_{j}^{
e}) \! \left(-s_{1} \! \left\{Q_{0}^{e}(a_{j}^{e}) \! \left[\gimel^{1}_{-1}
(a_{j}^{e}) \right. \right. \right. \\
\left. \left. \left. + \, \gimel^{1}_{1}(a_{j}^{e}) \! + \! \daleth^{1}_{1}
(a_{j}^{e}) \daleth^{1}_{-1}(a_{j}^{e}) \right] \! + \! Q_{1}^{e}(a_{j}^{e}) 
\right. \right. \\
\left. \left. \times \left[\daleth^{1}_{1}(a_{j}^{e}) \! + \! \daleth^{1}_{-
1}(a_{j}^{e}) \right] \! + \! \frac{1}{2}Q_{2}^{e}(a_{j}^{e}) \right. \right. 
\\
\left. \left. + \, (Q_{0}^{e}(a_{j}^{e}))^{-1} \aleph^{1}_{1}(a_{j}^{e}) 
\aleph^{1}_{-1}(a_{j}^{e}) \right\} \right. \\
\left. - \, t_{1} \! \left\{Q_{0}^{e}(a_{j}^{e}) \! \left[\aleph^{1}_{-1}
(a_{j}^{e}) \beth^{1}_{1}(a_{j}^{e}) \right. \right. \right. \\
\left. \left. \left. + \, \aleph^{1}_{1}(a_{j}^{e}) \beth^{1}_{-1}(a_{j}^{e}) 
\right] \! + \! (Q_{0}^{e}(a_{j}^{e}))^{-1} \right. \right. \\
\left. \left. \times \left[\daleth^{1}_{-1}(a_{j}^{e}) \! + \! \daleth^{1}_{
1}(a_{j}^{e}) \! - \! Q_{1}^{e}(a_{j}^{e})(Q_{0}^{e}(a_{j}^{e}))^{-1} \right] 
\right. \right. \\
\left. \left. + \, Q_{1}^{e}(a_{j}^{e}) \aleph^{1}_{1}(a_{j}^{e}) \aleph^{
1}_{-1}(a_{j}^{e}) \right\} \right. \\
\left. + \, \mi (s_{1} \! + \! t_{1}) \! \left\{\beth^{1}_{-1}(a_{j}^{e}) \! 
+ \! \aleph^{1}_{-1}(a_{j}^{e}) \daleth^{1}_{1}(a_{j}^{e}) \right. \right. \\
\left. \left. - \, \beth^{1}_{1}(a_{j}^{e}) \! - \! \aleph^{1}_{1}(a_{j}^{e}) 
\daleth^{1}_{-1}(a_{j}^{e}) \right\} \right)
\end{matrix}} & 
\boxed{\begin{matrix}
(\varkappa_{1}^{e}(a_{j}^{e}))^{2} \! \left(\mi s_{1} \! \left\{Q_{0}^{e}
(a_{j}^{e}) \! \left[2 \gimel^{1}_{1}(a_{j}^{e}) \right. \right. \right. \\
\left. \left. \left. + \, (\daleth^{1}_{1}(a_{j}^{e}))^{2} \right] \! + \! 
2Q_{1}^{e}(a_{j}^{e}) \daleth^{1}_{1}(a_{j}^{e}) \right. \right. \\
\left. \left. - \, (Q_{0}^{e}(a_{j}^{e}))^{-1}(\aleph^{1}_{1}(a_{j}^{e}))^{2} 
\! + \! \frac{1}{2}Q_{2}^{e}(a_{j}^{e}) \right\} \right. \\
\left. + \, \mi t_{1} \! \left\{2Q_{0}^{e}(a_{j}^{e}) \aleph^{1}_{1}(a_{j}^{
e}) \beth^{1}_{1}(a_{j}^{e}) \right. \right. \\
\left. \left. + \, Q_{1}^{e}(a_{j}^{e})(\aleph^{1}_{1}(a_{j}^{e}))^{2} \! + 
\! (Q_{0}^{e}(a_{j}^{e}))^{-1} \right. \right. \\
\left. \left. \times \left[Q_{1}^{e}(a_{j}^{e})(Q_{0}^{e}(a_{j}^{e}))^{-1} \! 
- \! 2 \daleth^{1}_{1}(a_{j}^{e}) \right] \right\} \right. \\
\left. - \, 2(s_{1} \! - \! t_{1}) \! \left\{\beth^{1}_{1}(a_{j}^{e}) \! + \! 
\aleph^{1}_{1}(a_{j}^{e}) \daleth^{1}_{1}(a_{j}^{e}) \right\} \right)
\end{matrix}} \\
\boxed{\begin{matrix} (\varkappa_{2}^{e}(a_{j}^{e}))^{2} \! \left(\mi s_{1} 
\! \left\{Q_{0}^{e}(a_{j}^{e}) \! \left[2 \gimel^{1}_{-1}(a_{j}^{e}) \right. 
\right. \right. \\
\left. \left. \left. + \, (\daleth^{1}_{-1}(a_{j}^{e}))^{2} \right] \! + \! 
2Q_{1}^{e}(a_{j}^{e}) \daleth^{1}_{-1}(a_{j}^{e}) \right. \right. \\
\left. \left. - \, (Q_{0}^{e}(a_{j}^{e}))^{-1}(\aleph^{1}_{-1}(a_{j}^{e}))^{
2} \! + \! \frac{1}{2}Q_{2}^{e}(a_{j}^{e}) \right\} \right. \\
\left. + \, \mi t_{1} \! \left\{2Q_{0}^{e}(a_{j}^{e}) \aleph^{1}_{-1}(a_{j}^{
e}) \beth^{1}_{-1}(a_{j}^{e}) \right. \right. \\
\left. \left. + \, Q_{1}^{e}(a_{j}^{e})(\aleph^{1}_{-1}(a_{j}^{e}))^{2} \! + 
\! (Q_{0}^{e}(a_{j}^{e}))^{-1} \right. \right. \\
\left. \left. \times \left[Q_{1}^{e}(a_{j}^{e})(Q_{0}^{e}(a_{j}^{e}))^{-1} \! 
- \! 2 \daleth^{1}_{-1}(a_{j}^{e}) \right] \right\} \right. \\
\left. + \, 2(s_{1} \! - \! t_{1}) \! \left\{\beth^{1}_{-1}(a_{j}^{e}) \! + 
\! \aleph^{1}_{-1}(a_{j}^{e}) \daleth^{1}_{-1}(a_{j}^{e}) \right\} \right)
\end{matrix}} & 
\boxed{\begin{matrix} \varkappa_{1}^{e}(a_{j}^{e}) \varkappa_{2}^{e}(a_{j}^{
e}) \! \left(s_{1} \! \left\{Q_{0}^{e}(a_{j}^{e}) \! \left[\gimel^{1}_{-1}
(a_{j}^{e}) \right. \right. \right. \\
\left. \left. \left. + \, \gimel^{1}_{1}(a_{j}^{e}) \! + \! \daleth^{1}_{1}
(a_{j}^{e}) \daleth^{1}_{-1}(a_{j}^{e}) \right] \! + \! Q_{1}^{e}(a_{j}^{e}) 
\right. \right. \\
\left. \left. \times \left[\daleth^{1}_{1}(a_{j}^{e}) \! + \! \daleth^{1}_{-
1}(a_{j}^{e}) \right] \! + \! \frac{1}{2}Q_{2}^{e}(a_{j}^{e}) \right. \right. 
\\
\left. \left. + \, (Q_{0}^{e}(a_{j}^{e}))^{-1}\aleph^{1}_{1}(a_{j}^{e}) 
\aleph^{1}_{-1}(a_{j}^{e}) \right\} \right. \\
\left. + \, t_{1} \! \left\{Q_{0}^{e}(a_{j}^{e}) \! \left[\aleph^{1}_{-1}
(a_{j}^{e}) \beth^{1}_{1}(a_{j}^{e}) \right. \right. \right. \\
\left. \left. \left. + \, \aleph^{1}_{1}(a_{j}^{e}) \beth^{1}_{-1}(a_{j}^{e}) 
\right] \! + \! (Q_{0}^{e}(a_{j}^{e}))^{-1} \right. \right. \\
\left. \left. \times \left[\daleth^{1}_{-1}(a_{j}^{e}) \! + \! \daleth^{1}_{
1}(a_{j}^{e}) \! - \! Q_{1}^{e}(a_{j}^{e})(Q_{0}^{e}(a_{j}^{e}))^{-1} \right] 
\right. \right. \\
\left. \left. + \, Q_{1}^{e}(a_{j}^{e}) \aleph^{1}_{1}(a_{j}^{e}) \aleph^{
1}_{-1}(a_{j}^{e}) \right\} \right. \\
\left. + \, \mi (s_{1} \! + \! t_{1}) \! \left\{\beth^{1}_{1}(a_{j}^{e}) \! 
+ \! \aleph^{1}_{1}(a_{j}^{e}) \daleth^{1}_{-1}(a_{j}^{e}) \right. \right. \\
\left. \left. - \, \beth^{1}_{-1}(a_{j}^{e}) \! - \! \aleph^{1}_{-1}(a_{j}^{
e}) \daleth^{1}_{1}(a_{j}^{e}) \right\} \right)
\end{matrix}}
\end{pmatrix},
\end{equation*}
(with $\operatorname{tr}(\mathscr{C}^{e}(a_{j}^{e})) \! = \! 0)$, and $(f^{
e}_{k}(n))_{ij} \! =_{n \to \infty} \! \mathcal{O}(1)$, $k \! \in \! \mathbb{
N}$, $i,j \! = \! 1,2$. (The expression for $\mathscr{C}^{e}(a_{j}^{e})$ is 
necessary for obtaining asymptotics \textbf{at} the end-points $\lbrace a_{
j}^{e} \rbrace_{j=1}^{N}$, as well as for Remark~5.2 below.) Returning to the 
counter-clockwise-oriented integrals $\oint_{\partial \mathbb{U}^{e}_{\delta_{
a_{j}}}} \tfrac{w_{+}^{\Sigma^{e}_{\circlearrowright}}(s)}{s-z} \, \tfrac{\md 
s}{2 \pi \mi}$, $z \! \in \! \mathbb{C} \setminus \widetilde{\Sigma}_{p}^{e}$, 
it follows, via the Residue and Cauchy Theorems, that, for $j \! = \! 1,\dotsc,
N$,
\begin{equation*}
\oint_{\partial \mathbb{U}^{e}_{\delta_{a_{j}}}} \dfrac{w_{+}^{\Sigma^{e}_{
\circlearrowright}}(s)}{s \! - \! z} \, \dfrac{\md s}{2 \pi \mi} \underset{
n \to \infty}{=}
\begin{cases}
-\dfrac{\widehat{\mathscr{A}}^{e}(a_{j}^{e})}{n(z \! - \! a_{j}^{e})^{2}} \!
- \! \dfrac{\widehat{\mathscr{B}}^{e}(a_{j}^{e})}{n(z \! - \! a_{j}^{e})} \!
+ \! \mathcal{O} \! \left(\dfrac{\widehat{f}^{e}(z;n)}{n^{2}} \right),
&\text{$z \! \in \! \mathbb{C} \setminus (\mathbb{U}^{e}_{\delta_{a_{j}}}
\cup \partial \mathbb{U}^{e}_{\delta_{a_{j}}})$,} \\
-\dfrac{\widehat{\mathscr{A}}^{e}(a_{j}^{e})}{n(z \! - \! a_{j}^{e})^{2}} \!
- \! \dfrac{\widehat{\mathscr{B}}^{e}(a_{j}^{e})}{n(z \! - \! a_{j}^{e})} \!
+ \! \dfrac{\mathscr{R}^{\infty}_{a_{j}^{e}}(z)}{n} \! + \! \mathcal{O} \!
\left(\dfrac{\widehat{f}^{e}(z;n)}{n^{2}} \right), &\text{$z \! \in \!
\mathbb{U}^{e}_{\delta_{a_{j}}}$,}
\end{cases}
\end{equation*}
where $\widehat{\mathscr{A}}^{e}(a_{j}^{e}) \! := \! (\widehat{\alpha}_{0}^{e}
(a_{j}^{e}))^{-1} \mathscr{A}^{e}(a_{j}^{e})$, $\widehat{\mathscr{B}}^{e}(a_{
j}^{e}) \! := \! (\widehat{\alpha}_{0}^{e}(a_{j}^{e}))^{-2}(\mathscr{B}^{e}
(a_{j}^{e}) \widehat{\alpha}_{0}^{e}(a_{j}^{e}) \! - \! \mathscr{A}^{e}(a_{
j}^{e}) \widehat{\alpha}_{1}^{e}(a_{j}^{e}))$, $\mathscr{R}^{\infty}_{a_{j}^{
e}} (z)$ is given in Theorem~2.3.1, Equations~(2.73) and~(2.74), and 
$\widehat{f}^{e}(z;n)$, where the $n$-dependence arises due to the 
$n$-dependence of the associated Riemann theta functions, denotes some 
bounded (with respect to both $z$ and $n)$, analytic (for $\mathbb{C} 
\setminus \widetilde{\Sigma}^{e}_{p} \! \ni \! z)$, $\operatorname{GL}_{2}
(\mathbb{C})$-valued function for which $(\widehat{f}^{e}(z;n))_{kl} \! 
=_{\underset{z \in \mathbb{C} \setminus \widetilde{\Sigma}^{e}_{p}}{n \to 
\infty}} \! \mathcal{O}(1)$, $k,l \! = \! 1,2$. Repeating the above analysis 
for the remaining end-points of the support of the `even' equilibrium measure, 
that is, $\lbrace b_{0}^{e},\dotsc,b_{N}^{e},a_{N+1}^{e} \rbrace$, one arrives 
at the result stated in the Lemma. \hfill $\qed$
\begin{eeeee}
A brisk perusing of the asymptotic (as $n \! \to \! \infty)$ result for 
$\mathscr{R}^{e}(z)$ stated in Lemma~5.3 seems to imply that, at first glance, 
there are second-order poles at $\lbrace b_{j-1}^{e},a_{j}^{e} \rbrace_{j=1}^{
N+1}$; however, this is not the case. As the proof of Lemma~5.3 demonstrates 
(cf. the analysis leading up to Equations~(5.3)), Laurent series expansions 
about $\lbrace b_{j-1}^{e},a_{j}^{e} \rbrace_{j=1}^{N+1}$ show that, as $n \! 
\to \! \infty$, all expansions are, indeed, analytic; in particular: (i) for 
$z \! \in \! \mathbb{U}^{e}_{\delta_{a_{j}}}$, $j \! = \! 1,\dotsc,N \!+ \! 
1$ (all contour integrals are counter-clockwise oriented),
\begin{align*}
\oint_{\partial \mathbb{U}^{e}_{\delta_{a_{j}}}} \dfrac{w_{+}^{\Sigma^{e}_{
\circlearrowright}}(s)}{s \! - \! z} \, \dfrac{\md s}{2 \pi \mi} \underset{
\underset{z \in \mathbb{U}^{e}_{\delta_{a_{j}}}}{n \to \infty}}{=}& \, \dfrac{
\left(\mathscr{A}^{e}(a_{j}^{e}) \widehat{\alpha}_{0}^{e}(a_{j}^{e}) \! \left(
\! \left(\frac{\widehat{\alpha}_{1}^{e}(a_{j}^{e})}{\widehat{\alpha}_{0}^{e}
(a_{j}^{e})} \right)^{2} \! - \! \frac{\widehat{\alpha}_{2}^{e}(a_{j}^{e})}{
\widehat{\alpha}_{0}^{e}(a_{j}^{e})} \right) \! - \! \mathscr{B}^{e}(a_{j}^{
e}) \widehat{\alpha}_{1}^{e}(a_{j}^{e}) \! + \! \mathscr{C}^{e}(a_{j}^{e})
\widehat{\alpha}_{0}^{e}(a_{j}^{e}) \right)}{n(\widehat{\alpha}_{0}^{e}(a_{
j}^{e}))^{2}} \\
&+ \, \dfrac{1}{n} \sum_{k=1}^{\infty}f_{k}^{a_{j}^{e}}(n)(z \! - \! a_{j}^{
e})^{k} \! + \! \mathcal{O} \! \left(\dfrac{\widehat{f}^{e}(z;n)}{n^{2}}
\right),
\end{align*}
where, for $j \! = \! 1,\dotsc,N$, $\mathscr{A}^{e}(a_{j}^{e})$, $\mathscr{
B}^{e}(a_{j}^{e})$, $\mathscr{C}^{e}(a_{j}^{e})$, $\widehat{\alpha}_{0}^{e}
(a_{j}^{e})$, $\widehat{\alpha}_{1}^{e}(a_{j}^{e})$, and $\widehat{\alpha}_{
2}^{e}(a_{j}^{e})$ are given in (the proof of) Lemma~5.3, $\mathscr{A}^{e}
(a_{N+1}^{e})$, $\mathscr{B}^{e}(a_{N+1}^{e})$ are given in Theorem~2.3.1, 
Equations~(2.25), (2.27), (2.28), (2.31), (2.32), (2.37)--(2.45), (2.47), 
(2.52), and~(2.53), $\widehat{\alpha}_{0}^{e}(a_{N+1}^{e}) \! := \! \tfrac{
4}{3}f(a_{N+1}^{e})$, $\widehat{\alpha}_{1}^{e}(a_{N+1}^{e}) \! := \! 
\tfrac{4}{5} f^{\prime}(a_{N+1}^{e})$, and $\widehat{\alpha}_{2}^{e}(a_{N+
1}^{e}) \! := \! \tfrac{2}{7}f^{\prime \prime}(a_{N+1}^{e})$, with $f(a_{N+
1}^{e})$, $f^{\prime}(a_{N+1}^{e})$, and $f^{\prime \prime}(a_{N+1}^{e})$ 
given in Lemma~4.7, $\mathscr{C}^{e}(a_{N+1}^{e})$ is given by the same 
expression as $\mathscr{ C}^{e}(a_{j}^{e})$ above subject to the modifications 
$\Omega_{j}^{e} \! \to \! 0$, $a_{j}^{e} \! \to \! a_{N+1}^{e}$, $Q_{0}^{e}
(a_{j}^{e}) \! \to \! Q_{0}^{e}(a_{N+1}^{e})$, $Q_{1}^{e}(a_{j}^{e}) \! \to 
\! Q_{1}^{e}(a_{N+1}^{e})$, with $Q_{0}^{e}(a_{N+1}^{e})$, $Q_{1}^{e}(a_{N+
1}^{e})$ given in Theorem~2.3.1, Equations~(2.31) and~(2.32), and $Q_{2}^{e}
(a_{j}^{e}) \! \to \! Q_{2}^{e}(a_{N+1}^{e})$, where
\begin{align*}
Q_{2}^{e}(a_{N+1}^{e})&= -\dfrac{1}{2}Q_{0}^{e}(a_{N+1}^{e}) \! \left(\sum_{k
=1}^{N} \! \left(\dfrac{1}{(a_{N+1}^{e} \! - \! b_{k}^{e})^{2}} \! - \! \dfrac{
1}{(a_{N+1}^{e} \! - \! a_{k}^{e})^{2}} \right) \! + \! \dfrac{1}{(a_{N+1}^{e}
\! - \! b_{0}^{e})^{2}} \right) \\
&+ \, \dfrac{1}{4}Q_{0}^{e}(a_{N+1}^{e}) \! \left(\sum_{k=1}^{N} \! \left(
\dfrac{1}{a_{N+1}^{e} \! - \! b_{k}^{e}} \! - \! \dfrac{1}{a_{N+1}^{e} \! - \!
a_{k}^{e}} \right) \! + \! \dfrac{1}{a_{N+1}^{e} \! - \! b_{0}^{e}} \right)^{
2},
\end{align*}
$\widehat{f}^{e}(z;n)$ is characterised completely at the end of the proof of
Lemma~5.3, and $(f_{k}^{a_{j}^{e}}(n))_{l_{1}l_{2}} \! =_{n \to \infty} \!
\mathcal{O}(1)$, $k \! \in \! \mathbb{N}$, $l_{1},l_{2} \! = \! 1,2$; and (ii)
for $z \! \in \! \mathbb{U}^{e}_{\delta_{b_{j}}}$, $j \! = \! 0,\dotsc,N$,
\begin{align*}
\oint_{\partial \mathbb{U}^{e}_{\delta_{b_{j}}}} \dfrac{w_{+}^{\Sigma^{e}_{
\circlearrowright}}(s)}{s \! - \! z} \, \dfrac{\md s}{2 \pi \mi} \underset{
\underset{z \in \mathbb{U}^{e}_{\delta_{b_{j}}}}{n \to \infty}}{=}& \, \dfrac{
\left(\mathscr{A}^{e}(b_{j}^{e}) \widehat{\alpha}_{0}^{e}(b_{j}^{e}) \! \left(
\! \left(\frac{\widehat{\alpha}_{1}^{e}(b_{j}^{e})}{\widehat{\alpha}_{0}^{e}
(b_{j}^{e})} \right)^{2} \! - \! \frac{\widehat{\alpha}_{2}^{e}(b_{j}^{e})}{
\widehat{\alpha}_{0}^{e}(b_{j}^{e})} \right) \! - \! \mathscr{B}^{e}(b_{j}^{
e}) \widehat{\alpha}_{1}^{e}(b_{j}^{e}) \! + \! \mathscr{C}^{e}(b_{j}^{e})
\widehat{\alpha}_{0}^{e}(b_{j}^{e}) \right)}{n(\widehat{\alpha}_{0}^{e}(b_{
j}^{e}))^{2}} \\
&+ \, \dfrac{1}{n} \sum_{k=1}^{\infty}f_{k}^{b_{j}^{e}}(n)(z \! - \! b_{j}^{
e})^{k} \! + \! \mathcal{O} \! \left(\dfrac{\widehat{f}^{e}(z;n)}{n^{2}}
\right),
\end{align*}
where, for $j \! = \! 1,\dotsc,N \! + \! 1$, $\mathscr{A}^{e}(b_{j-1}^{e})$, 
$\mathscr{B}^{e}(b_{j-1}^{e})$ are given in Theorem~2.3.1, Equations~(2.24), 
(2.26), (2.28), (2.29), (2.30), (2.33), (2.34), (2.37)--(2.45), (2.46), 
(2.48), (2.50), (2.51), (2.54), and~(2.55), $\widehat{\alpha}_{0}^{e}(b_{j-
1}^{e}) \! := \! \tfrac{4}{3}f(b_{j-1}^{e})$, $\widehat{\alpha}_{1}^{e}(b_{j
-1}^{e}) \! := \! \tfrac{4}{5}f^{\prime}(b_{j-1}^{e})$, and $\widehat{
\alpha}_{2}^{e}(b_{j-1}^{e}) \! := \! \tfrac{2}{7}f^{\prime \prime}(b_{j-1}^{
e})$, with $f(b_{j-1}^{e})$, $f^{\prime}(b_{j-1}^{e})$, and $f^{\prime \prime}
(b_{j-1}^{e})$ given in Lemma~4.6, for $j \! = \! 1,\dotsc,N$,
\begin{equation*}
\dfrac{\mathscr{C}^{e}(b_{j}^{e})}{\me^{\mi n \Omega_{j}^{e}}} \! := \! 
\begin{pmatrix} \boxed{\begin{matrix} \varkappa_{1}^{e}(b_{j}^{e}) 
\varkappa_{2}^{e}(b_{j}^{e}) \! \left(s_{1} \! \left\{-Q_{0}^{e}(b_{j}^{e}) 
\aleph^{1}_{1}(b_{j}^{e}) \aleph^{1}_{-1}(b_{j}^{e}) \right. \right. \\
\left. \left. - \, (Q_{0}^{e}(b_{j}^{e}))^{-3}(Q_{1}^{e}(b_{j}^{e}))^{2} \! 
+ \! \frac{1}{2}Q_{2}^{e}(b_{j}^{e})(Q_{0}^{e}(b_{j}^{e}))^{-2} \right. 
\right. \\
\left. \left. + \, Q_{1}^{e}(b_{j}^{e})(Q_{0}^{e}(b_{j}^{e}))^{-2} \! \left[
\daleth^{1}_{1}(b_{j}^{e}) \! + \! \daleth^{1}_{-1}(b_{j}^{e}) \right] 
\right. \right. \\
\left. \left. - \, (Q_{0}^{e}(b_{j}^{e}))^{-1} \! \left[\gimel^{1}_{1}
(b_{j}^{e}) \! + \! \gimel^{1}_{-1}(b_{j}^{e}) \! + \! \daleth^{1}_{1}
(b_{j}^{e}) \daleth^{1}_{-1}(b_{j}^{e}) \right] \right\} \right. \\
\left. + \, t_{1} \! \left\{-Q_{1}^{e}(b_{j}^{e}) \! - \! Q_{0}^{e}(b_{j}^{e}) 
\! \left[\daleth^{1}_{-1}(b_{j}^{e}) \! + \! \daleth^{1}_{1}(b_{j}^{e}) 
\right] \right. \right. \\
\left. \left. + \, Q_{1}^{e}(b_{j}^{e})(Q_{0}^{e}(b_{j}^{e}))^{-2} \aleph^{
1}_{1}(b_{j}^{e}) \aleph^{1}_{-1}(b_{j}^{e}) \! - \! (Q_{0}^{e}(b_{j}^{e})
)^{-1} \right. \right. \\
\left. \left. \times \left[\aleph^{1}_{1}(b_{j}^{e}) \beth^{1}_{-1}(b_{j}^{e}) 
\! + \! \aleph^{1}_{-1}(b_{j}^{e}) \beth^{1}_{1}(b_{j}^{e}) \right] \right\} 
\right. \\
\left. + \, \mi (s_{1} \! + \! t_{1}) \! \left[\beth^{1}_{-1}(b_{j}^{e}) \! 
- \! \aleph^{1}_{1}(b_{j}^{e}) \daleth^{1}_{-1}(b_{j}^{e}) \right. \right. \\
\left. \left. + \, \aleph^{1}_{-1}(b_{j}^{e}) \daleth^{1}_{1}(b_{j}^{e}) \! 
- \! \beth^{1}_{1}(b_{j}^{e}) \right] \right)
\end{matrix}} & 
\boxed{\begin{matrix} (\varkappa_{1}^{e}(b_{j}^{e}))^{2} \! \left(\mi s_{1} 
\! \left\{Q_{0}^{e}(b_{j}^{e})(\aleph^{1}_{1}(b_{j}^{e}))^{2} \! - \! (Q_{
0}^{e}(b_{j}^{e}))^{-3} \right. \right. \\
\left. \left. \times \, (Q_{1}^{e}(b_{j}^{e}))^{2} \! + \! \frac{1}{2}Q_{2}^{
e}(b_{j}^{e})(Q_{0}^{e}(b_{j}^{e}))^{-2} \! + \! 2Q_{1}^{e}(b_{j}^{e}) 
\right. \right. \\
\left. \left. \times \, (Q_{0}^{e}(b_{j}^{e}))^{-2} \daleth^{1}_{1}(b_{j}^{e}) 
\! - \! (Q_{0}^{e}(b_{j}^{e}))^{-1} \! \left[2 \gimel^{1}_{1}(b_{j}^{e}) 
\right. \right. \right. \\
\left. \left. \left. + \, (\daleth^{1}_{1}(b_{j}^{e}))^{2} \right] \right\} 
\! + \! \mi t_{1} \! \left\{2Q_{0}^{e}(b_{j}^{e}) \daleth^{1}_{1}(b_{j}^{e}) 
\! + \! Q_{1}^{e}(b_{j}^{e}) \right. \right. \\
\left. \left. + \, Q_{1}^{e}(b_{j}^{e})(Q_{0}^{e}(b_{j}^{e}))^{-2}(\aleph^{
1}_{1}(b_{j}^{e}))^{2} \! - \! 2(Q_{0}^{e}(b_{j}^{e}))^{-1} \right. \right. 
\\
\left. \left. \times \, \aleph^{1}_{1}(b_{j}^{e}) \beth^{1}_{1}(b_{j}^{e}) 
\right\} \! + \! 2(s_{1} \! - \! t_{1}) \right. \\
\left. \times \left\{\beth^{1}_{1}(b_{j}^{e}) \! + \! \aleph^{1}_{1}(b_{j}^{
e}) \daleth^{1}_{1}(b_{j}^{e}) \right\} \right)
\end{matrix}} \\
\boxed{\begin{matrix} (\varkappa_{2}^{e}(b_{j}^{e}))^{2} \! \left(\mi s_{1} 
\! \left\{Q_{0}^{e}(b_{j}^{e})(\aleph^{1}_{-1}(b_{j}^{e}))^{2} \! - \! (Q_{
0}^{e}(b_{j}^{e}))^{-3} \right. \right. \\
\left. \left. \times \, (Q_{1}^{e}(b_{j}^{e}))^{2} \! + \! \frac{1}{2}Q_{2}^{
e}(b_{j}^{e})(Q_{0}^{e}(b_{j}^{e}))^{-2} \! + \! 2Q_{1}^{e}(b_{j}^{e}) 
\right. \right. \\
\left. \left. \times \, (Q_{0}^{e}(b_{j}^{e}))^{-2} \daleth^{1}_{-1}(b_{j}^{
e}) \! - \! (Q_{0}^{e}(b_{j}^{e}))^{-1} \! \left[2 \gimel^{1}_{-1}(b_{j}^{e}) 
\right. \right. \right. \\
\left. \left. \left. + \, (\daleth^{1}_{-1}(b_{j}^{e}))^{2} \right] \right\} 
\! + \! \mi t_{1} \! \left\{2Q_{0}^{e}(b_{j}^{e}) \daleth^{1}_{-1}(b_{j}^{e}) 
\! + \! Q_{1}^{e}(b_{j}^{e}) \right. \right. \\
\left. \left. + \, Q_{1}^{e}(b_{j}^{e})(Q_{0}^{e}(b_{j}^{e}))^{-2}(\aleph^{
1}_{-1}(b_{j}^{e}))^{2} \! - \! 2(Q_{0}^{e}(b_{j}^{e}))^{-1} \right. \right. 
\\
\left. \left. \times \, \aleph^{1}_{-1}(b_{j}^{e}) \beth^{1}_{-1}(b_{j}^{e}) 
\right\} \! - \! 2(s_{1} \! - \! t_{1}) \right. \\
\left. \times \left\{\beth^{1}_{-1}(b_{j}^{e}) \! + \! \aleph^{1}_{-1}
(b_{j}^{e}) \daleth^{1}_{-1}(b_{j}^{e}) \right\} \right)
\end{matrix}} & \boxed{\begin{matrix} \varkappa_{1}^{e}(b_{j}^{e}) 
\varkappa_{2}^{e}(b_{j}^{e}) \! \left(s_{1} \! \left\{Q_{0}^{e}(b_{j}^{e}) 
\aleph^{1}_{1}(b_{j}^{e}) \aleph^{1}_{-1}(b_{j}^{e}) \right. \right. \\
\left. \left. + \, (Q_{0}^{e}(b_{j}^{e}))^{-3}(Q_{1}^{e}(b_{j}^{e}))^{2} \! 
- \! \frac{1}{2}Q_{2}^{e}(b_{j}^{e})(Q_{0}^{e}(b_{j}^{e}))^{-2} \right. 
\right. \\
\left. \left. - \, Q_{1}^{e}(b_{j}^{e})(Q_{0}^{e}(b_{j}^{e}))^{-2} \! \left[
\daleth^{1}_{1}(b_{j}^{e}) \! + \! \daleth^{1}_{-1}(b_{j}^{e}) \right] 
\right. \right. \\
\left. \left. + \, (Q_{0}^{e}(b_{j}^{e}))^{-1} \! \left[\gimel^{1}_{1}(b_{
j}^{e}) \! + \! \gimel^{1}_{-1}(b_{j}^{e}) \! + \! \daleth^{1}_{1}(b_{j}^{e}) 
\daleth^{1}_{-1}(b_{j}^{e}) \right] \right\} \right. \\
\left. + \, t_{1} \! \left\{Q_{1}^{e}(b_{j}^{e}) \! + \! Q_{0}^{e}(b_{j}^{e}) 
\! \left[\daleth^{1}_{-1}(b_{j}^{e}) \! + \! \daleth^{1}_{1}(b_{j}^{e}) 
\right] \right. \right. \\
\left. \left. - \, Q_{1}^{e}(b_{j}^{e})(Q_{0}^{e}(b_{j}^{e}))^{-2} \aleph^{
1}_{1}(b_{j}^{e}) \aleph^{1}_{-1}(b_{j}^{e}) \! + \! (Q_{0}^{e}(b_{j}^{e})
)^{-1} \right. \right. \\
\left. \left. \times \left[\aleph^{1}_{1}(b_{j}^{e}) \beth^{1}_{-1}(b_{j}^{
e}) \! + \! \aleph^{1}_{-1}(b_{j}^{e}) \beth^{1}_{1}(b_{j}^{e}) \right] 
\right\} \right. \\
\left. + \, \mi (s_{1} \! + \! t_{1}) \! \left[\beth^{1}_{1}(b_{j}^{e}) \! - 
\! \aleph^{1}_{-1}(b_{j}^{e}) \daleth^{1}_{1}(b_{j}^{e}) \right. \right. \\
\left. \left. + \, \aleph^{1}_{1}(b_{j}^{e}) \daleth^{1}_{-1}(b_{j}^{e}) \! 
- \! \beth^{1}_{-1}(b_{j}^{e}) \right] \right)
\end{matrix}}
\end{pmatrix}
\end{equation*}
(with $\operatorname{tr}(\mathscr{C}^{e}(b_{j}^{e})) \! = \! 0)$, where 
$Q_{0}^{e}(b_{j}^{e})$, $Q_{1}^{e}(b_{j}^{e})$ are given in Theorem~2.3.1, 
Equations~(2.33) and~(2.34),
\begin{align*}
Q_{2}^{e}(b_{j}^{e})&= -\dfrac{1}{2}Q_{0}^{e}(b_{j}^{e}) \! \left(\sum_{
\substack{k=1\\k \not= j}}^{N} \! \left(\dfrac{1}{(b_{j}^{e} \! - \! b_{k}^{
e})^{2}} \! - \! \dfrac{1}{(b_{j}^{e} \! - \! a_{k}^{e})^{2}} \right) \! + \!
\dfrac{1}{(b_{j}^{e} \! - \! b_{0}^{e})^{2}} \! - \! \dfrac{1}{(b_{j}^{e} \!
- \! a_{N+1}^{e})^{2}} \! - \! \dfrac{1}{(b_{j}^{e} \! - \! a_{j}^{e})^{2}}
\right) \\
&+ \, \dfrac{1}{4}Q_{0}^{e}(b_{j}^{e}) \! \left(\sum_{\substack{k=1\\k \not=
j}}^{N} \! \left(\dfrac{1}{b_{j}^{e} \! - \! b_{k}^{e}} \! - \! \dfrac{1}{b_{
j}^{e} \! - \! a_{k}^{e}} \right) \! + \! \dfrac{1}{b_{j}^{e} \! - \! b_{0}^{
e}} \! - \! \dfrac{1}{b_{j}^{e} \! - \! a_{N+1}^{e}} \! - \! \dfrac{1}{b_{j}^{
e} \! - \! a_{j}^{e}} \right)^{2},
\end{align*}
$\mathscr{C}^{e}(b_{0}^{e})$ is given by the same expression as
$\mathscr{C}^{ e}(b_{j}^{e})$ above subject to the modifications
$\Omega_{j}^{e} \! \to \! 0$, $b_{ j}^{e} \! \to \! b_{0}^{e}$,
$Q_{0}^{e}(b_{j}^{e}) \! \to \! Q_{0}^{e}(b_{0}^{ e})$,
$Q_{1}^{e}(b_{j}^{e}) \! \to \! Q_{1}^{e}(b_{0}^{e})$, with
$Q_{0}^{e} (b_{0}^{e})$, $Q_{1}^{e}(b_{0}^{e})$ given in
Theorem~2.3.1, Equations~(2.29) and~(2.30), and $Q_{2}^{e}(b_{j}^{e}) \! \to 
\! Q_{2}^{e}(b_{0}^{e})$, where
\begin{align*}
Q_{2}^{e}(b_{0}^{e})&= -\dfrac{1}{2}Q_{0}^{e}(b_{0}^{e}) \! \left(\sum_{k=1}^{
N} \! \left(\dfrac{1}{(b_{0}^{e} \! - \! b_{k}^{e})^{2}} \! - \! \dfrac{1}{(
b_{0}^{e} \! - \! a_{k}^{e})^{2}} \right) \! - \! \dfrac{1}{(b_{0}^{e} \! - \!
a_{N+1}^{e})^{2}} \right) \\
&+ \, \dfrac{1}{4}Q_{0}^{e}(b_{0}^{e}) \! \left(\sum_{k=1}^{N} \! \left(
\dfrac{1}{b_{0}^{e} \! - \! b_{k}^{e}} \! - \! \dfrac{1}{b_{0}^{e} \! - \!
a_{k}^{e}} \right) \! - \! \dfrac{1}{b_{0}^{e} \! - \! a_{N+1}^{e}} \right)^{
2},
\end{align*}
and $(f_{k}^{b_{j-1}^{e}}(n))_{l_{1}l_{2}} \! =_{n \to \infty} \! \mathcal{O}
(1)$, $j \! = \! 1,\dotsc,N \! + \! 1$, $k \! \in \! \mathbb{N}$, $l_{1},l_{2}
\! = \! 1,2$. \hfill $\blacksquare$
\end{eeeee}

Re-tracing the finite sequence of RHP transformations (all of which are 
invertible) and definitions, namely, $\mathscr{R}^{e}(z)$ (Lemmas~5.3 and 
4.8) and $\mathscr{S}^{e}_{p}(z)$ (Lemma~4.8) $\to$ $\mathcal{X}^{e}(z)$
(Lemmas~4.6 and~4.7) $\to$
$\overset{e}{m}^{\raise-1.0ex\hbox{$\scriptstyle \infty$}}(z)$ (Lemma 4.5) 
$\to$ $\overset{e}{\mathscr{M}}^{\raise-1.0ex\hbox{$\scriptstyle \sharp$}}(z)$
(Lemma~4.2) $\to$ $\overset{e}{\mathscr{M}}(z)$ (Lemma~4.1) $\to$ $\overset{e}{
\mathrm{Y}}(z)$ (Lemma~3.4), the asymptotic (as $n \! \to \! \infty)$ solution
of the original \textbf{RHP1}, that is, $(\overset{e}{\mathrm{Y}}(z),\mathrm{
I} \! + \! \exp (-n \widetilde{V}(z)) \sigma_{+},\mathbb{R})$, in the various 
bounded and unbounded regions (Figure~7), is given by:
\begin{compactenum}
\item[(1)] for $z \! \in \! \Upsilon^{e}_{1} \cup \Upsilon^{e}_{2}$,
\begin{equation*}
\overset{e}{\mathrm{Y}}(z) \! = \! \me^{\frac{n \ell_{e}}{2} \operatorname{ad}
(\sigma_{3})} \mathscr{R}^{e}(z)
\overset{e}{m}^{\raise-1.0ex\hbox{$\scriptstyle \infty$}}(z) \me^{n
(g^{e}(z)+\int_{J_{e}} \ln (s) \psi_{V}^{e}(s) \, \md s) \sigma_{3}},
\end{equation*}
where
$g^{e}(z)$, $\psi_{V}^{e}(z)$, $\ell_{e}$,
$\overset{e}{m}^{\raise-1.0ex\hbox{$\scriptstyle \infty$}}(z)$, and $\mathscr{
R}^{e}(z)$ are given in Lemmas~3.4, 3.5, 3.6, 4.5, and~5.3, respectively;
\item[(2)] for $z \! \in \! \Upsilon^{e}_{3}$,
\begin{equation*}
\overset{e}{\mathrm{Y}}(z) \! = \! \me^{\frac{n \ell_{e}}{2} \operatorname{ad}
(\sigma_{3})} \mathscr{R}^{e}(z)
\overset{e}{m}^{\raise-1.0ex\hbox{$\scriptstyle \infty$}}(z) \!
\left(\mathrm{I} \! + \! \me^{-4n \pi \mi \int_{z}^{a_{N+1}^{e}} \psi_{V}^{e}
(s) \, \md s} \, \sigma_{-} \right) \! \me^{n(g^{e}(z)+\int_{J_{e}} \ln (s)
\psi_{V}^{e}(s) \, \md s) \sigma_{3}};
\end{equation*}
\item[(3)] for $z \! \in \! \Upsilon^{e}_{4}$,
\begin{equation*}
\overset{e}{\mathrm{Y}}(z) \! = \! \me^{\frac{n \ell_{e}}{2} \operatorname{ad}
(\sigma_{3})} \mathscr{R}^{e}(z)
\overset{e}{m}^{\raise-1.0ex\hbox{$\scriptstyle \infty$}}(z) \!
\left(\mathrm{I} \! - \! \me^{4n \pi \mi \int_{z}^{a_{N+1}^{e}} \psi_{V}^{e}
(s) \, \md s} \, \sigma_{-} \right) \! \me^{n(g^{e}(z)+\int_{J_{e}} \ln (s)
\psi_{V}^{e}(s) \, \md s) \sigma_{3}};
\end{equation*}
\item[(4)] for $z \! \in \! (\Omega^{e,1}_{b_{j-1}} \cup \Omega^{e,4}_{b_{j-
1}}) \cup (\Omega^{e,1}_{a_{j}} \cup \Omega^{e,4}_{a_{j}})$, $j \! = \! 1,
\dotsc,N \! + \! 1$,
\begin{equation*}
\overset{e}{\mathrm{Y}}(z) \! = \! \me^{\frac{n \ell_{e}}{2} \operatorname{ad}
(\sigma_{3})} \mathscr{R}^{e}(z)
\mathcal{X}^{e}(z) \me^{n(g^{e}(z)+\int_{J_{e}} \ln (s) \psi_{V}^{e}(s) \, \md
s) \sigma_{3}},
\end{equation*}
where, for $z \! \in \! \mathbb{U}^{e}_{\delta_{b_{j-1}}}$ $(\supset \Omega^{
e,1}_{b_{j-1}} \cup \Omega^{e,4}_{b_{j-1}})$, $\mathcal{X}^{e}(z)$ is given by
Lemma~4.6, and, for $z \! \in \! \mathbb{U}^{e}_{\delta_{a_{j}}}$ $(\supset
\Omega^{e,1}_{a_{j}} \cup \Omega^{e,4}_{a_{j}})$, $\mathcal{X}^{e}(z)$ is
given by Lemma~4.7;
\item[(5)] for $z \! \in \! \Omega^{e,2}_{b_{j-1}} \cup \Omega^{e,2}_{a_{j}}$,
$j \! = \! 1,\dotsc,N \! + \! 1$,
\begin{equation*}
\overset{e}{\mathrm{Y}}(z) \! = \! \me^{\frac{n \ell_{e}}{2} \operatorname{ad}
(\sigma_{3})} \mathscr{R}^{e}(z)
\mathcal{X}^{e}(z) \! \left(\mathrm{I} \! + \! \me^{-4n \pi \mi
\int_{z}^{a_{N+1}^{e}} \psi_{V}^{e}(s) \, \md s} \, \sigma_{-} \right) \!
\me^{n(g^{e}(z)+\int_{J_{e}} \ln (s) \psi_{V}^{e}(s) \, \md s) \sigma_{3}};
\end{equation*}
\item[(6)] for $z \! \in \! \Omega^{e,3}_{b_{j-1}} \cup \Omega^{e,3}_{a_{j}}$,
$j \! = \! 1,\dotsc,N \! + \! 1$,
\begin{equation*}
\overset{e}{\mathrm{Y}}(z) \! = \! \me^{\frac{n \ell_{e}}{2} \operatorname{ad}
(\sigma_{3})} \mathscr{R}^{e}(z)
\mathcal{X}^{e}(z) \! \left(\mathrm{I} \! - \! \me^{4n \pi \mi \int_{
z}^{a_{N+1}^{e}} \psi_{V}^{e}(s) \, \md s} \, \sigma_{-} \right) \! \me^{n
(g^{e}(z)+\int_{J_{e}} \ln (s) \psi_{V}^{e}(s) \, \md s) \sigma_{3}}.
\end{equation*}
\end{compactenum}
Multiplying the respective matrices in items~(1)--(6) above and collecting
$(1 \, 1)$- and $(1 \, 2)$-elements, one arrives at, finally, the asymptotic
(as $n \! \to \! \infty)$ results for $\boldsymbol{\pi}_{2n}(z)$ and $\int_{
\mathbb{R}} \tfrac{\boldsymbol{\pi}_{2n}(s) \exp (-n \widetilde{V}(s))}{s-z}
\, \tfrac{\md s}{2 \pi \mi}$ (in the entire complex plane) stated in
Theorem~2.3.1.

In order to obtain asymptotics (as $n \! \to \! \infty)$ for $\xi^{(2n)}_{n}$
$(= \! \norm{\boldsymbol{\pi}_{2n}(\cdot)}^{-1}_{\mathscr{L}} \! = \! (H^{(-2
n)}_{2n}/H^{(-2n)}_{2n+1})^{1/2})$ and $\phi_{2n}(z)$ $(= \! \xi^{(2n)}_{n}
\boldsymbol{\pi}_{2n}(z))$ stated in Theorem~2.3.2, large-$z$ asymptotics
for $\overset{e}{\operatorname{Y}}(z)$ are necessary.
\begin{bbbbb}
Let $\mathscr{R}^{e} \colon \mathbb{C} \setminus \widetilde{\Sigma}_{p}^{e}
\! \to \! \operatorname{SL}_{2}(\mathbb{C})$ be the solution of the {\rm RHP}
$(\mathscr{R}^{e}(z),\upsilon_{\mathscr{R}}^{e}(z),\widetilde{\Sigma}_{p}^{e}
)$ formulated in Proposition~{\rm 5.2} with $n \! \to \! \infty$ asymptotics
given in Lemma~{\rm 5.3}. Then,
\begin{equation*}
\mathscr{R}^{e}(z) \underset{z \to \infty}{=} \mathrm{I} \! + \! \dfrac{1}{
z} \mathscr{R}^{e,\infty}_{1}(n) \! + \! \dfrac{1}{z^{2}} \mathscr{R}^{e,
\infty}_{2}(n) \! + \! \mathcal{O} \! \left(\dfrac{1}{z^{3}} \right),
\end{equation*}
where, for $k \! = \! 1,2$,
\begin{equation*}
\mathscr{R}^{e,\infty}_{k}(n) \! := \! -\int_{\Sigma^{e}_{\circlearrowright}}
s^{k-1}w^{\Sigma^{e}_{\circlearrowright}}_{+}(s) \, \dfrac{\md s}{2 \pi \mi}
\! = \! \sum_{j=1}^{N+1} \, \sum_{q \in \{b_{j-1}^{e},a_{j}^{e}\}}
\operatorname{Res} \! \left(z^{k-1}w^{\Sigma^{e}_{\circlearrowright}}_{+}(z);
q \right),
\end{equation*}
with, in particular,
\begin{align*}
\mathscr{R}^{e,\infty}_{1}(n) \underset{n \to \infty}{=}& \, \dfrac{1}{n}
\sum_{j=1}^{N+1} \! \left(\dfrac{(\mathscr{B}^{e}(a_{j}^{e}) \widehat{\alpha}_{
0}^{e}(a_{j}^{e}) \! - \! \mathscr{A}^{e}(a_{j}^{e}) \widehat{\alpha}_{1}^{e}
(a_{j}^{e}))}{(\widehat{\alpha}_{0}^{e}(a_{j}^{e}))^{2}} \! + \! \dfrac{(
\mathscr{B}^{e}(b_{j-1}^{e}) \widehat{\alpha}_{0}^{e}(b_{j-1}^{e}) \! - \!
\mathscr{A}^{e}(b_{j-1}^{e}) \widehat{\alpha}_{1}^{e}(b_{j-1}^{e}))}{(
\widehat{\alpha}_{0}^{e}(b_{j-1}^{e}))^{2}} \right) \\
+& \, \mathcal{O} \! \left(\dfrac{1}{n^{2}} \right), \\
\mathscr{R}^{e,\infty}_{2}(n) \underset{n \to \infty}{=}& \, \dfrac{1}{n}
\sum_{j=1}^{N+1} \! \left(\dfrac{(\widehat{\alpha}_{0}^{e}(b_{j-1}^{e})
\mathscr{A}^{e}(b_{j-1}^{e}) \! + \! b_{j-1}^{e}(\mathscr{B}^{e}(b_{j-1}^{e})
\widehat{\alpha}_{0}^{e}(b_{j-1}^{e}) \! - \! \mathscr{A}^{e}(b_{j-1}^{e})
\widehat{\alpha}_{1}^{e}(b_{j-1}^{e})))}{(\widehat{\alpha}_{0}^{e}(b_{j-1}^{
e}))^{2}} \right. \\
+&\left. \, \dfrac{(\widehat{\alpha}_{0}^{e}(a_{j}^{e}) \mathscr{A}^{e}(a_{
j}^{e}) \! + \! a_{j}^{e}(\mathscr{B}^{e}(a_{j}^{e}) \widehat{\alpha}_{0}^{e}
(a_{j}^{e}) \! - \! \mathscr{A}^{e}(a_{j}^{e}) \widehat{\alpha}_{1}^{e}(a_{
j}^{e})))}{(\widehat{\alpha}_{0}^{e}(a_{j}^{e}))^{2}} \right) \! + \!
\mathcal{O} \! \left(\dfrac{1}{n^{2}} \right),
\end{align*}
and all the parameters defined in Lemma~{\rm 5.3}.

Let $\overset{e}{m}^{\raise-1.0ex\hbox{$\scriptstyle \infty$}} \colon \mathbb{
C} \setminus J_{e}^{\infty} \! \to \! \operatorname{SL}_{2}(\mathbb{C})$ solve
the {\rm RHP} $(\overset{e}{m}^{\raise-1.0ex\hbox{$\scriptstyle \infty$}}(z),
J_{e}^{\infty},
\overset{e}{\upsilon}^{\raise-1.0ex\hbox{$\scriptstyle \infty$}}(z))$
formulated in Lemma~{\rm 4.3} with (unique) solution given by
Lemma~{\rm 4.5}. For $\varepsilon_{1},\varepsilon_{2} \! = \! \pm 1$, set
\begin{gather*}
\theta^{e}_{\infty}(\varepsilon_{1},\varepsilon_{2},\boldsymbol{\Omega}^{e})
\! := \! \boldsymbol{\theta}^{e}(\varepsilon_{1} \boldsymbol{u}^{e}_{+}
(\infty) \! - \! \tfrac{n}{2 \pi} \boldsymbol{\Omega}^{e} \! + \!
\varepsilon_{2} \boldsymbol{d}_{e}), \\
\alpha^{e}_{\infty}(\varepsilon_{1},\varepsilon_{2},\boldsymbol{\Omega}^{e})
\! := \! 2 \pi \mi \varepsilon_{1} \sum_{m \in \mathbb{Z}^{N}} (m,\widehat{
\boldsymbol{\alpha}}^{e}_{\infty}) \me^{2 \pi \mi (m,\varepsilon_{1}
\boldsymbol{u}^{e}_{+}(\infty)-\frac{n}{2 \pi} \boldsymbol{\Omega}^{e}+
\varepsilon_{2} \boldsymbol{d}_{e})+ \pi \mi (m,\tau^{e}m)},
\end{gather*}
where $\widehat{\boldsymbol{\alpha}}^{e}_{\infty} \! = \! (\widehat{\alpha}_{
\infty,1}^{e},\widehat{\alpha}^{e}_{\infty,2},\dotsc,\widehat{\alpha}^{e}_{
\infty,N})$, with $\widehat{\alpha}^{e}_{\infty,j} \! := \! c_{j1}^{e}$, $j
\! = \! 1,\dotsc,N$, and
\begin{gather*}
\beta^{e}_{\infty}(\varepsilon_{1},\varepsilon_{2},\boldsymbol{\Omega}^{e}) \!
:= \! 2 \pi \sum_{m \in \mathbb{Z}^{N}} \! \left(\pi (m,\widehat{\boldsymbol{
\alpha}}^{e}_{\infty})^{2} \! + \! \mi \varepsilon_{1}(m,\widehat{\boldsymbol{
\beta}}^{e}_{\infty}) \right) \! \me^{2 \pi \mi (m,\varepsilon_{1} \boldsymbol{
u}^{e}_{+}(\infty)-\frac{n}{2 \pi} \boldsymbol{\Omega}^{e}+\varepsilon_{2}
\boldsymbol{d}_{e})+ \pi \mi (m,\tau^{e}m)},
\end{gather*}
where $\widehat{\boldsymbol{\beta}}^{e}_{\infty} \! = \! (\widehat{\beta}^{e}_{
\infty,1},\widehat{\beta}^{e}_{\infty,2},\dotsc,\widehat{\beta}^{e}_{\infty,
N})$, with $\widehat{\beta}^{e}_{\infty,j} \! := \! \tfrac{1}{2}(c^{e}_{j2} \!
+ \! \tfrac{1}{2}c_{j1}^{e} \sum_{k=1}^{N+1}(b_{k-1}^{e} \! + \! a_{k}^{e}))$,
$j \! = \! 1,\dotsc,N$, where $c_{j1}^{e},c_{j2}^{e}$, $j \! = \! 1,\dotsc,
N$, are obtained {}from Equations~{\rm (E1)} and~{\rm (E2)}. Then,
\begin{gather*}
\overset{e}{m}^{\raise-1.0ex\hbox{$\scriptstyle \infty$}}(z) \underset{z \to
\infty}{=} \mathrm{I} \! + \! \dfrac{1}{z}
\overset{e}{m}_{1}^{\raise-1.0ex\hbox{$\scriptstyle \infty$}} \! + \! \dfrac{
1}{z^{2}} \overset{e}{m}_{2}^{\raise-1.0ex\hbox{$\scriptstyle \infty$}} \! +
\! \mathcal{O} \! \left(\dfrac{1}{z^{3}} \right),
\end{gather*}
where
\begin{align*}
(\overset{e}{m}_{1}^{\raise-1.0ex\hbox{$\scriptstyle \infty$}})_{11} &=\dfrac{
\boldsymbol{\theta}^{e}(\boldsymbol{u}^{e}_{+}(\infty) \! + \! \boldsymbol{d}_{
e})}{\boldsymbol{\theta}^{e}(\boldsymbol{u}^{e}_{+}(\infty) \! - \! \frac{n}{2
\pi} \boldsymbol{\Omega}^{e} \! + \! \boldsymbol{d}_{e})} \! \left(\dfrac{
\theta^{e}_{\infty}(1,1,\boldsymbol{\Omega}^{e}) \alpha^{e}_{\infty}(1,1,\vec{
\pmb{0}}) \! - \! \alpha^{e}_{\infty}(1,1,\boldsymbol{\Omega}^{e}) \theta^{
e}_{\infty}(1,1,\vec{\pmb{0}})}{(\theta^{e}_{\infty}(1,1,\vec{\pmb{0}}))^{2}}
\right), \\
(\overset{e}{m}_{1}^{\raise-1.0ex\hbox{$\scriptstyle \infty$}})_{12} &=\dfrac{
1}{4 \mi} \! \left(\sum_{k=1}^{N+1}(b_{k-1}^{e} \! - \! a_{k}^{e}) \right) \!
\dfrac{\boldsymbol{\theta}^{e}(\boldsymbol{u}^{e}_{+}(\infty) \! + \!
\boldsymbol{d}_{e}) \theta^{e}_{\infty}(-1,1,\boldsymbol{\Omega}^{e})}{
\boldsymbol{\theta}^{e}(\boldsymbol{u}^{e}_{+}(\infty) \! - \! \frac{n}{2 \pi}
\boldsymbol{\Omega}^{e} \! + \! \boldsymbol{d}_{e}) \theta^{e}_{\infty}(-1,1,
\vec{\pmb{0}})}, \\
(\overset{e}{m}_{1}^{\raise-1.0ex\hbox{$\scriptstyle \infty$}})_{21}&=-\dfrac{
1}{4 \mi} \! \left(\sum_{k=1}^{N+1}(b_{k-1}^{e} \! - \! a_{k}^{e}) \right)
\! \dfrac{\boldsymbol{\theta}^{e}(\boldsymbol{u}^{e}_{+}(\infty) \! + \!
\boldsymbol{d}_{e}) \theta^{e}_{\infty}(1,-1,\boldsymbol{\Omega}^{e})}{
\boldsymbol{\theta}^{e}(-\boldsymbol{u}^{e}_{+}(\infty) \! - \! \frac{n}{2
\pi} \boldsymbol{\Omega}^{e} \! - \! \boldsymbol{d}_{e}) \theta^{e}_{\infty}
(1,-1,\vec{\pmb{0}})}, \\
(\overset{e}{m}_{1}^{\raise-1.0ex\hbox{$\scriptstyle \infty$}})_{22} &=\left(
\dfrac{\theta^{e}_{\infty}(-1,-1,\boldsymbol{\Omega}^{e}) \alpha^{e}_{\infty}
(-1,-1,\vec{\pmb{0}}) \! - \! \alpha^{e}_{\infty}(-1,-1,\boldsymbol{\Omega}^{
e}) \theta^{e}_{\infty}(-1,-1,\vec{\pmb{0}})}{(\theta^{e}_{\infty}(-1,-1,\vec{
\pmb{0}}))^{2}} \right) \\
&\times \dfrac{\boldsymbol{\theta}^{e}(\boldsymbol{u}^{e}_{+}(\infty) \! + \!
\boldsymbol{d}_{e})}{\boldsymbol{\theta}^{e}(-\boldsymbol{u}^{e}_{+}(\infty)
\! - \! \frac{n}{2 \pi} \boldsymbol{\Omega}^{e} \! - \! \boldsymbol{d}_{e})},
\\
(\overset{e}{m}_{2}^{\raise-1.0ex\hbox{$\scriptstyle \infty$}})_{11} &=\left(
\theta^{e}_{\infty}(1,1,\boldsymbol{\Omega}^{e}) \! \left(\beta^{e}_{\infty}
(1,1,\vec{\pmb{0}}) \theta^{e}_{\infty}(1,1,\vec{\pmb{0}}) \! + \! (\alpha^{
e}_{\infty}(1,1,\vec{\pmb{0}}))^{2} \right) \! - \! \alpha^{e}_{\infty}(1,1,
\boldsymbol{\Omega}^{e}) \alpha^{e}_{\infty}(1,1,\vec{\pmb{0}}) \theta^{e}_{
\infty}(1,1,\vec{\pmb{0}}) \right. \\
&\left. - \, \beta^{e}_{\infty}(1,1,\boldsymbol{\Omega}^{e})(\theta^{e}_{
\infty}(1,1,\vec{\pmb{0}}))^{2} \right) \! \dfrac{(\theta^{e}_{\infty}(1,1,
\vec{\pmb{0}}))^{-3} \boldsymbol{\theta}^{e}(\boldsymbol{u}^{e}_{+}(\infty)
\! + \! \boldsymbol{d}_{e})}{\boldsymbol{\theta}^{e}(\boldsymbol{u}^{e}_{+}
(\infty) \! - \! \frac{n}{2 \pi} \boldsymbol{\Omega}^{e} \! + \! \boldsymbol{
d}_{e})} \! + \! \dfrac{1}{32} \! \left(\sum_{k=1}^{N+1}(b_{k-1}^{e} \! - \!
a_{k}^{e}) \right)^{2}, \\
(\overset{e}{m}_{2}^{\raise-1.0ex\hbox{$\scriptstyle \infty$}})_{12} &=
\dfrac{\boldsymbol{\theta}^{e}(\boldsymbol{u}^{e}_{+}(\infty) \! + \!
\boldsymbol{d}_{e})}{\boldsymbol{\theta}^{e}(\boldsymbol{u}^{e}_{+}(\infty) \!
- \! \frac{n}{2 \pi} \boldsymbol{\Omega}^{e} \! + \! \boldsymbol{d}_{e})} \!
\left(\! \left(\dfrac{\theta^{e}_{\infty}(-1,1,\boldsymbol{\Omega}^{e})
\alpha^{e}_{\infty}(-1,1,\vec{\pmb{0}}) \! - \! \alpha^{e}_{\infty}(-1,1,
\boldsymbol{\Omega}^{e}) \theta^{e}_{\infty}(-1,1,\vec{\pmb{0}})}{(\theta^{e}_{
\infty}(-1,1,\vec{\pmb{0}}))^{2}} \right) \right. \\
&\left. \times \, \dfrac{1}{4 \mi} \! \left(\sum_{k=1}^{N+1}(b_{k-1}^{e} \! -
\! a_{k}^{e}) \right) \! + \! \dfrac{1}{8 \mi} \! \left(\sum_{k=1}^{N+1}((b_{k
-1}^{e})^{2} \! - \! (a_{k}^{e})^{2}) \right) \! \dfrac{\theta^{e}_{\infty}(-1,
1,\boldsymbol{\Omega}^{e})}{\theta^{e}_{\infty}(-1,1,\vec{\pmb{0}})} \right),
\\
(\overset{e}{m}_{2}^{\raise-1.0ex\hbox{$\scriptstyle \infty$}})_{21} &=-
\dfrac{\boldsymbol{\theta}^{e}(\boldsymbol{u}^{e}_{+}(\infty) \! + \!
\boldsymbol{d}_{e})}{\boldsymbol{\theta}^{e}(-\boldsymbol{u}^{e}_{+}(\infty)
\! - \! \frac{n}{2 \pi} \boldsymbol{\Omega}^{e} \! - \! \boldsymbol{d}_{e})}
\! \left(\! \left(\dfrac{\theta^{e}_{\infty}(1,-1,\boldsymbol{\Omega}^{e})
\alpha^{e}_{\infty}(1,-1,\vec{\pmb{0}}) \! - \! \alpha^{e}_{\infty}(1,-1,
\boldsymbol{\Omega}^{e}) \theta^{e}_{\infty}(1,-1,\vec{\pmb{0}})}{(\theta^{e}_{
\infty}(1,-1,\vec{\pmb{0}}))^{2}} \right) \right. \\
&\left. \times \, \dfrac{1}{4 \mi} \! \left(\sum_{k=1}^{N+1}(b_{k-1}^{e} \! -
\! a_{k}^{e}) \right) \! + \! \dfrac{1}{8 \mi} \! \left(\sum_{k=1}^{N+1}((b_{k
-1}^{e})^{2} \! - \! (a_{k}^{e})^{2}) \right) \! \dfrac{\theta^{e}_{\infty}(1,
-1,\boldsymbol{\Omega}^{e})}{\theta^{e}_{\infty}(1,-1,\vec{\pmb{0}})} \right),
\\
(\overset{e}{m}_{2}^{\raise-1.0ex\hbox{$\scriptstyle \infty$}})_{22} &=\left(
\theta^{e}_{\infty}(-1,-1,\boldsymbol{\Omega}^{e}) \! \left(\beta^{e}_{\infty}
(-1,-1,\vec{\pmb{0}}) \theta^{e}_{\infty}(-1,-1,\vec{\pmb{0}}) \! + \!
(\alpha^{e}_{\infty}(-1,-1,\vec{\pmb{0}}))^{2} \right) \! - \! \alpha^{e}_{
\infty}(-1,-1,\boldsymbol{\Omega}^{e}) \right. \\
&\left. \times \, \alpha^{e}_{\infty}(-1,-1,\vec{\pmb{0}}) \theta^{e}_{\infty}
(-1,-1,\vec{\pmb{0}}) \! - \! \beta^{e}_{\infty}(-1,-1,\boldsymbol{\Omega}^{e})
(\theta^{e}_{\infty}(-1,-1,\vec{\pmb{0}}))^{2} \right) \! (\theta^{e}_{\infty}
(-1,-1,\vec{\pmb{0}}))^{-3} \\
&\times \dfrac{\boldsymbol{\theta}^{e}(\boldsymbol{u}^{e}_{+}(\infty) \! + \!
\boldsymbol{d}_{e})}{\boldsymbol{\theta}^{e}(-\boldsymbol{u}^{e}_{+}(\infty)
\! - \! \frac{n}{2 \pi} \boldsymbol{\Omega}^{e} \! - \! \boldsymbol{d}_{e})}
\! + \! \dfrac{1}{32} \! \left(\sum_{k=1}^{N+1}(b_{k-1}^{e} \! - \! a_{k}^{e})
\right)^{2},
\end{align*}
with $(\star)_{ij}$, $i,j \! = \! 1,2$, denoting the $(i \, j)$-element of
$\star$, and $\vec{\pmb{0}} \! := \! (0,0,\dotsc,0)^{\mathrm{T}}$ $(\in \!
\mathbb{R}^{N})$.

Let $\overset{e}{\mathrm{Y}} \colon \mathbb{C} \setminus \mathbb{R} \! \to \!
\operatorname{SL}_{2}(\mathbb{C})$ be the solution of {\rm \pmb{RHP1}}. Then,
\begin{equation*}
\overset{e}{\mathrm{Y}}(z)z^{-n \sigma_{3}} \underset{z \to \infty}{=}
\mathrm{I} \! + \! \dfrac{1}{z} \mathrm{Y}^{e,\infty}_{1} \! + \! \dfrac{1}{
z^{2}} \mathrm{Y}^{e,\infty}_{2} \! + \! \mathcal{O} \! \left(\dfrac{1}{z^{
3}} \right),
\end{equation*}
where
\begin{align*}
(\mathrm{Y}^{e,\infty}_{1})_{11} &= -2n \int_{J_{e}}s \psi_{V}^{e}(s) \, \md s
\! + \! (\overset{e}{m}_{1}^{\raise-1.0ex\hbox{$\scriptstyle \infty$}})_{11}
\! + \! (\mathscr{R}^{e,\infty}_{1}(n))_{11}, \\
(\mathrm{Y}^{e,\infty}_{1})_{12} &= \me^{n \ell_{e}} \! \left(
(\overset{e}{m}_{1}^{\raise-1.0ex\hbox{$\scriptstyle \infty$}})_{12} \! + \!
(\mathscr{R}^{e,\infty}_{1}(n))_{12} \right), \\
(\mathrm{Y}^{e,\infty}_{1})_{21} &= \me^{-n \ell_{e}} \! \left(
(\overset{e}{m}_{1}^{\raise-1.0ex\hbox{$\scriptstyle \infty$}})_{21} \! + \!
(\mathscr{R}^{e,\infty}_{1}(n))_{21} \right), \\
(\mathrm{Y}^{e,\infty}_{1})_{22} &= 2n \int_{J_{e}}s \psi_{V}^{e}(s) \, \md s
\! + \! (\overset{e}{m}_{1}^{\raise-1.0ex\hbox{$\scriptstyle \infty$}})_{22}
\! + \! (\mathscr{R}^{e,\infty}_{1}(n))_{22}, \\
(\mathrm{Y}^{e,\infty}_{2})_{11} &= 2n^{2} \! \left(\int_{J_{e}}s \psi_{V}^{e}
(s) \, \md s \right)^{2} \! - \! n \int_{J_{e}}s^{2} \psi_{V}^{e}(s) \, \md s
\! - \! 2n \! \left(
(\overset{e}{m}_{1}^{\raise-1.0ex\hbox{$\scriptstyle \infty$}})_{11} \! + \!
(\mathscr{R}^{e,\infty}_{1}(n))_{11} \right) \! \int_{J_{e}}s \psi_{V}^{e}(s)
\, \md s \\
&+(\overset{e}{m}_{2}^{\raise-1.0ex\hbox{$\scriptstyle \infty$}})_{11} \! + \!
(\mathscr{R}^{e,\infty}_{2}(n))_{11} \! + \! (\mathscr{R}^{e,\infty}_{1}(n))_{
11}(\overset{e}{m}_{1}^{\raise-1.0ex\hbox{$\scriptstyle \infty$}})_{11} \! +
\! (\mathscr{R}^{e,\infty}_{1}(n))_{12}
(\overset{e}{m}_{1}^{\raise-1.0ex\hbox{$\scriptstyle \infty$}})_{21}, \\
(\mathrm{Y}^{e,\infty}_{2})_{12} &= \me^{n \ell_{e}} \! \left(2n \! \left(
(\overset{e}{m}_{1}^{\raise-1.0ex\hbox{$\scriptstyle \infty$}})_{12} \! + \!
(\mathscr{R}^{e,\infty}_{1}(n))_{12} \right) \! \int_{J_{e}}s \psi_{V}^{e}(s)
\, \md s \! + \!
(\overset{e}{m}_{2}^{\raise-1.0ex\hbox{$\scriptstyle \infty$}})_{12} \! + \!
(\mathscr{R}^{e,\infty}_{2}(n))_{12} \right. \\
&\left. + \, (\mathscr{R}^{e,\infty}_{1}(n))_{11}
(\overset{e}{m}_{1}^{\raise-1.0ex\hbox{$\scriptstyle \infty$}})_{12} \! + \!
(\mathscr{R}^{e,\infty}_{1}(n))_{12}
(\overset{e}{m}_{1}^{\raise-1.0ex\hbox{$\scriptstyle \infty$}})_{22} \right),
\\
(\mathrm{Y}^{e,\infty}_{2})_{21} &= \me^{-n \ell_{e}} \! \left(-2n \! \left(
(\overset{e}{m}_{1}^{\raise-1.0ex\hbox{$\scriptstyle \infty$}})_{21} \! + \!
(\mathscr{R}^{e,\infty}_{1}(n))_{21} \right) \! \int_{J_{e}}s \psi_{V}^{e}(s)
\, \md s \! + \!
(\overset{e}{m}_{2}^{\raise-1.0ex\hbox{$\scriptstyle \infty$}})_{21} \! + \!
(\mathscr{R}^{e,\infty}_{2}(n))_{21} \right. \\
&\left. + \, (\mathscr{R}^{e,\infty}_{1}(n))_{21}
(\overset{e}{m}_{1}^{\raise-1.0ex\hbox{$\scriptstyle \infty$}})_{11} \! + \!
(\mathscr{R}^{e,\infty}_{1}(n))_{22}
(\overset{e}{m}_{1}^{\raise-1.0ex\hbox{$\scriptstyle \infty$}})_{21} \right),
\\
(\mathrm{Y}^{e,\infty}_{2})_{22} &= 2n^{2} \! \left(\int_{J_{e}}s \psi_{V}^{e}
(s) \, \md s \right)^{2} \! + \! n \int_{J_{e}}s^{2} \psi_{V}^{e}(s) \, \md s
\! + \! 2n \! \left(
(\overset{e}{m}_{1}^{\raise-1.0ex\hbox{$\scriptstyle \infty$}})_{22} \! + \!
(\mathscr{R}^{e,\infty}_{1}(n))_{22} \right) \! \int_{J_{e}}s \psi_{V}^{e}(s)
\, \md s \\
&+(\overset{e}{m}_{2}^{\raise-1.0ex\hbox{$\scriptstyle \infty$}})_{22} \! + \!
(\mathscr{R}^{e,\infty}_{2}(n))_{22} \! + \! (\mathscr{R}^{e,\infty}_{1}(n))_{
21}(\overset{e}{m}_{1}^{\raise-1.0ex\hbox{$\scriptstyle \infty$}})_{12} \! +
\! (\mathscr{R}^{e,\infty}_{1}(n))_{22}
(\overset{e}{m}_{1}^{\raise-1.0ex\hbox{$\scriptstyle \infty$}})_{22}.
\end{align*}
\end{bbbbb}

\emph{Proof.} Let $\mathscr{R}^{e} \colon \mathbb{C} \setminus \widetilde{
\Sigma}^{e}_{p} \! \to \! \operatorname{SL}_{2}(\mathbb{C})$ be the solution
of the RHP $(\mathscr{R}^{e}(z),\upsilon_{\mathscr{R}}^{e}(z),\widetilde{
\Sigma}^{e}_{p})$ formulated in Proposition~5.2 with $n \! \to \! \infty$
asymptotics given in Lemma~5.3. For $\vert z \vert \! \gg \! \max_{j=1,
\dotsc,N+1} \lbrace \vert b_{j-1}^{e} \! - \! a_{j}^{e} \vert \rbrace$, via
the expansion $\tfrac{1}{s-z} \! = \! -\sum_{k=0}^{l} \tfrac{s^{k}}{z^{k+1}}
\! + \! \tfrac{s^{l+1}}{z^{l+1}(s-z)}$, $l \! \in \! \mathbb{Z}_{0}^{+}$,
where $s \! \in \! \lbrace b_{j-1}^{e},a_{j}^{e} \rbrace$, $j \! = \! 1,
\dotsc,N \! + \! 1$, one obtains the asymptotics for $\mathscr{R}^{e}(z)$
stated in the Proposition.

Let $\overset{e}{m}^{\raise-1.0ex\hbox{$\scriptstyle \infty$}} \colon \mathbb{
C} \setminus J_{e}^{\infty} \! \to \! \operatorname{SL}_{2}(\mathbb{C})$ solve
the RHP $(\overset{e}{m}^{\raise-1.0ex\hbox{$\scriptstyle \infty$}}(z),J_{e}^{
\infty},\overset{e}{\upsilon}^{\raise-1.0ex\hbox{$\scriptstyle \infty$}}(z))$
formulated in Lemma~4.3 with (unique) solution given by Lemma~4.5.
In order to obtain large-$z$ asymptotics of
$\overset{e}{m}^{\raise-1.0ex\hbox{$\scriptstyle \infty$}}(z)$, one needs
large-$z$ asymptotics of $(\gamma^{e}(z))^{\pm 1}$ and $\tfrac{\boldsymbol{
\theta}^{e}(\varepsilon_{1} \boldsymbol{u}^{e}(z)-\frac{n}{2 \pi} \boldsymbol{
\Omega}^{e}+\varepsilon_{2} \boldsymbol{d}_{e})}{\boldsymbol{\theta}^{e}(
\varepsilon_{1} \boldsymbol{u}^{e}(z)+\varepsilon_{2} \boldsymbol{d}_{e})}$,
$\varepsilon_{1},\varepsilon_{2} \! = \! \pm 1$. Consider, say, and without
loss of generality, $z \! \to \! \infty$ asymptotics for $z \! \in \! \mathbb{
C}_{+}$, that is, $z \! \to \! \infty^{+}$, where, by definition, $\sqrt{
\smash[b]{\star (z)}} \! := \! +\sqrt{\smash[b]{\star (z)}}$: equivalently,
one may consider $z \! \to \! \infty$ asymptotics for $z \! \in \! \mathbb{
C}_{-}$, that is, $z \! \to \! \infty^{-}$; however, recalling that $\sqrt{
\smash[b]{\star (z)}} \! \upharpoonright_{\mathbb{C}_{+}} \! = \! -\sqrt{
\smash[b]{\star (z)}} \! \upharpoonright_{\mathbb{C}_{-}}$, one obtains (in
either case, and via the sheet-interchange index) the same $z \! \to \!
\infty$ asymptotics (for
$\overset{e}{m}^{\raise-1.0ex\hbox{$\scriptstyle \infty$}}(z))$. Recall the
expression for $\gamma^{e}(z)$ given in Lemma~4.4: for $\vert z \vert \! \gg
\! \max_{j=1,\dotsc,N+1} \lbrace \vert b_{j-1}^{e} \! - \! a_{j}^{e} \vert
\rbrace $, via the expansions $\tfrac{1}{s-z} \! = \! -\sum_{k=0}^{l} \tfrac{
s^{k}}{z^{k+1}} \! + \! \tfrac{s^{l+1}}{z^{l+1}(s-z)}$, $l \! \in \! \mathbb{
Z}_{0}^{+}$, and $\ln (z \! - \! s) \! =_{\vert z \vert \to \infty} \! \ln
(z) \! - \! \sum_{k=1}^{\infty} \tfrac{1}{k}(\tfrac{s}{z})^{k}$, where $s \!
\in \! \lbrace b_{j-1}^{e},a_{j}^{e} \rbrace$, $j \! = \! 1,\dotsc,N \! + \!
1$, one shows that
\begin{align*}
(\gamma^{e}(z))^{\pm 1} \underset{z \to \infty^{+}}{=}& \, 1 \! + \! \dfrac{
1}{z} \! \left(\pm \dfrac{1}{4} \sum_{k=1}^{N+1}(a_{k}^{e} \! - \! b_{k-1}^{
e}) \right) \! + \! \dfrac{1}{z^{2}} \! \left(\pm \dfrac{1}{8} \sum_{k=1}^{N
+1} \! \left((a_{k}^{e})^{2} \! - \! (b_{k-1}^{e})^{2} \right) \right. \\
+&\left. \dfrac{1}{32} \! \left(\sum_{k=1}^{N+1}(a_{k}^{e} \! - \! b_{k-1}^{
e}) \right)^{2} \, \right) \! + \! \mathcal{O} \! \left(\dfrac{1}{z^{3}}
\right),
\end{align*}
whence
\begin{equation*}
\dfrac{1}{2}(\gamma^{e}(z) \! + \! (\gamma^{e}(z))^{-1}) \underset{z \to
\infty^{+}}{=} 1 \! + \! \dfrac{1}{z^{2}} \! \left(\dfrac{1}{32} \! \left(
\sum_{k=1}^{N+1}(a_{k}^{e} \! - \! b_{k-1}^{e}) \right)^{2} \, \right) \! + \!
\mathcal{O} \! \left(\dfrac{1}{z^{3}} \right),
\end{equation*}
and
\begin{equation*}
\dfrac{1}{2 \mi}(\gamma^{e}(z) \! - \! (\gamma^{e}(z))^{-1}) \underset{z \to
\infty^{+}}{=} \dfrac{1}{z} \! \left(\dfrac{1}{4 \mi} \sum_{k=1}^{N+1}(a_{
k}^{e} \! - \! b_{k-1}^{e}) \right) \! + \! \dfrac{1}{z^{2}} \! \left(
\dfrac{1}{8 \mi} \sum_{k=1}^{N+1} \! \left((a_{k}^{e})^{2} \! - \! (b_{k-
1}^{e})^{2} \right) \right) \! + \! \mathcal{O} \! \left(\dfrac{1}{z^{3}}
\right).
\end{equation*}
Recall {}from Lemma~4.5 that $\boldsymbol{u}^{e}(z) \! := \! \int_{a_{N+1}^{
e}}^{z} \boldsymbol{\omega}^{e}$ $(\in \operatorname{Jac}(\mathcal{Y}_{e})$, 
with $\mathcal{Y}_{e} \! := \! \lbrace \mathstrut (y,z); \, y^{2} \! = \! 
R_{e}(z) \rbrace)$, where $\boldsymbol{\omega}^{e}$, the associated normalised 
basis of holomorphic one-forms of $\mathcal{Y}_{e}$, is given by $\boldsymbol{
\omega}^{e} \! = \! (\omega_{1}^{e},\omega_{2}^{e},\dotsc,\omega_{N}^{e})$, 
with $\omega_{j}^{e} \! := \! \sum_{k=1}^{N}c_{jk}^{e}(\prod_{i=1}^{N+1}
(z \! - \! b_{i-1}^{e})(z \! - \! a_{i}^{e}))^{-1/2}z^{N-k} \, \md z$, $j \! 
= \! 1,\dotsc,N$, where $c_{jk}^{e}$, $j,k \! = \! 1,\dotsc,N$, are obtained 
{}from Equations~(E1) and~(E2). Writing
\begin{equation*}
\boldsymbol{u}^{e}(z) \! = \! \left(\int_{a_{N+1}^{e}}^{\infty^{+}} \! + \!
\int_{\infty^{+}}^{z} \right) \! \boldsymbol{\omega}^{e} \! = \! \boldsymbol{
u}^{e}_{+}(\infty) \! + \! \int_{\infty^{+}}^{z} \boldsymbol{\omega}^{e},
\end{equation*}
where $\boldsymbol{u}^{e}_{+}(\infty) \! := \! \int_{a_{N+1}^{e}}^{\infty^{+}
} \boldsymbol{\omega}^{e}$ (cf. Lemma~4.5), for $\vert z \vert \! \gg \!
\max_{j=1,\dotsc,N+1} \lbrace \vert b_{j-1}^{e} \! - \! a_{j}^{e} \vert
\rbrace$, via the expansions $\tfrac{1}{s-z} \! = \! -\sum_{k=0}^{l} \tfrac{
s^{k}}{z^{k+1}} \! + \! \tfrac{s^{l+1}}{z^{l+1}(s-z)}$, $l \! \in \! \mathbb{
Z}_{0}^{+}$, and $\ln (z \! - \! s) \! =_{\vert z \vert \to \infty} \! \ln
(z) \! - \! \sum_{k=1}^{\infty} \tfrac{1}{k}(\tfrac{s}{z})^{k}$, where $s \!
\in \! \lbrace b_{k-1}^{e},a_{k}^{e} \rbrace$, $k \! = \! 1,\dotsc,N \! + \!
1$, one shows that, for $j \! = \! 1,\dotsc,N$,
\begin{equation*}
\omega_{j}^{e} \underset{z \to \infty^{+}}{=} \dfrac{c_{j1}^{e}}{z^{2}} \,
\md z \! + \! \dfrac{(c_{j2}^{e} \! + \! \frac{1}{2}c_{j1}^{e} \sum_{i=1}^{
N+1}(a_{i}^{e} \! + \! b_{i-1}^{e}))}{z^{3}} \, \md z \! + \! \mathcal{O} \!
\left(\dfrac{\md z}{z^{4}} \right),
\end{equation*}
whence
\begin{align*}
\int_{\infty^{+}}^{z} \omega_{j}^{e} \underset{z \to \infty^{+}}{=}& -\dfrac{
c_{j1}^{e}}{z} \! - \! \dfrac{\frac{1}{2}(c_{j2}^{e} \! + \! \frac{1}{2}c_{j
1}^{e} \sum_{i=1}^{N+1}(a_{i}^{e} \! + \! b_{i-1}^{e}))}{z^{2}} \! + \!
\mathcal{O} \! \left(\dfrac{1}{z^{3}} \right) \\
=:& -\dfrac{\widehat{\alpha}^{e}_{\infty,j}}{z} \! - \! \dfrac{\widehat{
\beta}^{e}_{\infty,j}}{z^{2}} \! + \! \mathcal{O} \! \left(\dfrac{1}{z^{3}}
\right).
\end{align*}
Defining $\theta^{e}_{\infty}(\varepsilon_{1},\varepsilon_{2},\boldsymbol{
\Omega}^{e})$, $\alpha^{e}_{\infty}(\varepsilon_{1},\varepsilon_{2},
\boldsymbol{\Omega}^{e})$, and $\beta^{e}_{\infty}(\varepsilon_{1},
\varepsilon_{2},\boldsymbol{\Omega}^{e})$, $\varepsilon_{1},\varepsilon_{2}
\! = \! \pm 1$, as in the Proposition, recalling that $\boldsymbol{\omega}^{
e} \! = \! (\omega_{1}^{e},\omega_{2}^{e},\dotsc,\omega_{N}^{e})$, and that
the associated $N \! \times \! N$ Riemann matrix of
$\boldsymbol{\beta}^{e}$-periods, that is, $\tau^{e} \! = \! (\tau^{e})_{i,j
=1,\dotsc,N} \! := \! (\oint_{\boldsymbol{\beta}^{e}_{j}} \omega^{e}_{i})_{i,
j=1,\dotsc,N}$, is non-degenerate, symmetric, and $-\mi \tau^{e}$ is positive
definite, via the above asymptotic (as $z \! \to \! \infty^{+})$ expansion for
$\int_{\infty^{+}}^{z} \omega^{e}_{j}$, $j \! = \! 1,\dotsc,N$, one shows that
\begin{equation*}
\dfrac{\boldsymbol{\theta}^{e}(\varepsilon_{1} \boldsymbol{u}^{e}(z) \! - \!
\frac{n}{2 \pi} \boldsymbol{\Omega}^{e} \! + \! \varepsilon_{2} \boldsymbol{
d}_{e})}{\boldsymbol{\theta}^{e}(\varepsilon_{1} \boldsymbol{u}^{e}(z) \!
+ \! \varepsilon_{2} \boldsymbol{d}_{e})} \underset{z \to \infty^{+}}{=}
\varTheta_{0}^{e} \! + \! \dfrac{1}{z} \varTheta_{1}^{e} \! + \! \dfrac{1}{z^{
2}} \varTheta_{2}^{e} \! + \! \mathcal{O} \! \left(\dfrac{1}{z^{3}} \right),
\end{equation*}
where
\begin{align*}
\varTheta_{0}^{e} :=& \, \dfrac{\theta^{e}_{\infty}(\varepsilon_{1},
\varepsilon_{2},\boldsymbol{\Omega}^{e})}{\theta^{e}_{\infty}(\varepsilon_{1},
\varepsilon_{2},\vec{\pmb{0}})}, \\
\varTheta_{1}^{e} :=& \, \dfrac{\theta^{e}_{\infty}(\varepsilon_{1},
\varepsilon_{2},\boldsymbol{\Omega}^{e}) \alpha^{e}_{\infty}(\varepsilon_{1},
\varepsilon_{2},\vec{\pmb{0}}) \! - \! \alpha^{e}_{\infty}(\varepsilon_{1},
\varepsilon_{2},\boldsymbol{\Omega}^{e}) \theta^{e}_{\infty}(\varepsilon_{1},
\varepsilon_{2},\vec{\pmb{0}})}{(\theta^{e}_{\infty}(\varepsilon_{1},
\varepsilon_{2},\vec{\pmb{0}}))^{2}}, \\
\varTheta^{e}_{2} :=& \left(\theta^{e}_{\infty}(\varepsilon_{1},\varepsilon_{
2},\boldsymbol{\Omega}^{e}) \! \left(\beta^{e}_{\infty}(\varepsilon_{1},
\varepsilon_{2},\vec{\pmb{0}}) \theta^{e}_{\infty}(\varepsilon_{1},
\varepsilon_{2},\vec{\pmb{0}}) \! + \! (\alpha^{e}_{\infty}(\varepsilon_{1},
\varepsilon_{2},\vec{\pmb{0}}))^{2} \right) \! - \! \alpha^{e}_{\infty}
(\varepsilon_{1},\varepsilon_{2},\boldsymbol{\Omega}^{e}) \right. \\
\times&\left. \alpha^{e}_{\infty}(\varepsilon_{1},\varepsilon_{2},\vec{\pmb{0}}
) \theta^{e}_{\infty}(\varepsilon_{1},\varepsilon_{2},\vec{\pmb{0}}) \! - \!
\beta^{e}_{\infty}(\varepsilon_{1},\varepsilon_{2},\boldsymbol{\Omega}^{e})
(\theta^{e}_{\infty}(\varepsilon_{1},\varepsilon_{2},\vec{\pmb{0}}))^{2}
\right) \! (\theta^{e}_{\infty}(\varepsilon_{1},\varepsilon_{2},\vec{\pmb{0}}
))^{-3},
\end{align*}
with $\vec{\pmb{0}} \! := \! (0,0,\dotsc,0)^{\mathrm{T}}$ $(\in \! \mathbb{
R}^{N})$. Via the above asymptotic (as $z \! \to \! \infty^{+})$ expansions
for $\tfrac{1}{2}(\gamma^{e}(z) \! + \! (\gamma^{e}(z))^{-1})$, $\tfrac{1}{2
\mi}(\gamma^{e}(z) \! - \! (\gamma^{e}(z))^{-1})$, and $\tfrac{\boldsymbol{
\theta}^{e}(\varepsilon_{1} \boldsymbol{u}^{e}(z)-\frac{n}{2 \pi} \boldsymbol{
\Omega}^{e}+\varepsilon_{2} \boldsymbol{d}_{e})}{\boldsymbol{\theta}^{e}
(\varepsilon_{1} \boldsymbol{u}^{e}(z)+\varepsilon_{2} \boldsymbol{d}_{e})}$,
one arrives at, upon recalling the expression for
$\overset{e}{m}^{\raise-1.0ex\hbox{$\scriptstyle \infty$}}(z)$ given in
Lemma~4.5, the asymptotic expansion for
$\overset{e}{m}^{\raise-1.0ex\hbox{$\scriptstyle \infty$}}(z)$ stated in the
Proposition.

Let $\overset{e}{\operatorname{Y}} \colon \mathbb{C} \setminus \mathbb{R} \! 
\to \! \operatorname{SL}_{2}(\mathbb{C})$ be the (unique) solution of \textbf{
RHP1}, that is, $(\overset{e}{\operatorname{Y}}(z),\mathrm{I} \! + \! \exp 
(-n \widetilde{V}(z)) \sigma_{+},\mathbb{R})$. Recall, also, that, for $z \! 
\in \! \Upsilon^{e}_{1} \cup \Upsilon^{e}_{2}$ (Figure~7),
\begin{equation*}
\overset{e}{\operatorname{Y}}(z) \! = \! \me^{\tfrac{n \ell_{e}}{2}
\operatorname{ad}(\sigma_{3})} \mathscr{R}^{e}(z)
\overset{e}{m}^{\raise-1.0ex\hbox{$\scriptstyle \infty$}}
(z) \me^{n(g^{e}(z)+\int_{J_{e}} \ln (s) \psi_{V}^{e}(s) \, \md s) \sigma_{
3}}:\end{equation*}
consider, say, and without loss of generality, large-$z$ asymptotics for
$\overset{e}{\operatorname{Y}}(z)$ for $z \! \in \! \Upsilon_{1}^{e}$.
Recalling the definition of $g^{e}(z)$ given in Lemma~3.4, that is, $g^{e}(z)
\! := \! \int_{J_{e}} \ln ((z \! - \! s)^{2}(zs)^{-1}) \psi_{V}^{e}(s) \, \md
s$, $z \! \in \! \mathbb{C} \setminus (-\infty,\max \lbrace 0,\linebreak[4]
a_{N+1}^{e} \rbrace)$, for $\vert z \vert \! \gg \! \max_{j=1,\dotsc,N+1}
\lbrace \vert b_{j-1}^{e} \! - \! a_{j}^{e} \vert \rbrace$, in particular,
$\vert s/z \vert \! \ll \! 1$ with $s \! \in \! J_{e}$, and noting that
$\int_{J_{e}} \psi_{V}^{e}(s) \, \md s \! = \! 1$ and $\int_{J_{e}}s^{m}
\psi_{V}^{e}(s) \, \md s \! < \! \infty$, $m \! \in \! \mathbb{N}$, via the
expansions $\tfrac{1}{s-z} \! = \! -\sum_{k=0}^{l} \tfrac{s^{k}}{z^{k+1}} \!
+ \! \tfrac{s^{l+1}}{z^{l+1}(s-z)}$, $l \! \in \! \mathbb{Z}_{0}^{+}$, and
$\ln (z \! - \! s) \! =_{\vert z \vert \to \infty} \! \ln (z) \! - \! \sum_{k
=1}^{\infty} \tfrac{1}{k}(\tfrac{s}{z})^{k}$, one shows that
\begin{equation*}
g^{e}(z) \underset{z \to \infty}{=} \ln (z) \! - \! \int_{J_{e}} \ln (s)
\psi_{V}^{e}(s) \, \md s \! + \! \dfrac{1}{z} \! \left(-2 \int_{J_{e}}s \psi_{
V}^{e}(s) \, \md s \right) \! + \! \dfrac{1}{z^{2}} \! \left(-\int_{J_{e}}
s^{2} \psi_{V}^{e}(s) \, \md s \right) \! + \! \mathcal{O} \! \left(\dfrac{
1}{z^{3}}\right)
\end{equation*}
(explicit expressions for $\int_{J_{e}}s^{k} \psi_{V}^{e}(s) \, \md s$, $k \!
= \! 1,2$, are given in Remark~3.2): using the asymptotic (as $z \! \to \!
\infty)$ expansions for $g^{e}(z)$, $\mathscr{R}^{e}(z)$, and
$\overset{e}{m}^{\raise-1.0ex\hbox{$\scriptstyle \infty$}}(z)$ derived above,
upon recalling the formula for $\overset{e}{\operatorname{Y}}(z)$, one arrives
at, after a matrix-multiplication argument, the asymptotic expansion for
$\overset{e}{\operatorname{Y}}(z)z^{-n \sigma_{3}}$ stated in the Proposition.
\hfill $\qed$
\begin{bbbbb}
Let $\overset{e}{\operatorname{Y}} \colon \mathbb{C} \setminus \mathbb{R} \!
\to \! \operatorname{SL}_{2}(\mathbb{C})$ be the solution of {\rm \pmb{RHP1}}
with $z$ $(\in \! \mathbb{C} \setminus \mathbb{R})$ $\to \! \infty$
asymptotics given in Proposition~{\rm 5.3}. Then,
\begin{equation*}
\xi^{(2n)}_{n} \! = \! \dfrac{1}{\norm{\boldsymbol{\pi}_{2n}(\pmb{\cdot})}_{
\mathscr{L}}} \! = \sqrt{\dfrac{H^{(-2n)}_{2n}}{H^{(-2n)}_{2n+1}}}= \! \left(
-\dfrac{1}{2 \pi \mi (\operatorname{Y}^{e,\infty}_{1})_{12}} \right)^{1/2}
\quad (> \! 0),
\end{equation*}
where $(\operatorname{Y}^{e,\infty}_{1})_{12} \! = \! \me^{n \ell_{e}} \!
\left((\overset{e}{m}_{1}^{\raise-1.0ex\hbox{$\scriptstyle \infty$}})_{12} \!
+ \! (\mathscr{R}^{e,\infty}_{1}(n))_{12} \right)$, with
$(\overset{e}{m}_{1}^{\raise-1.0ex\hbox{$\scriptstyle \infty$}})_{12}$ and
$(\mathscr{R}^{e,\infty}_{1}(n))_{12}$ given in Proposition~{\rm 5.3}.
\end{bbbbb}

\emph{Proof.} Recall {}from Lemma~2.2.1 that $\boldsymbol{\pi}_{2n}(z) \! :=
\! (\overset{e}{\operatorname{Y}}(z))_{11}$ and $(\overset{e}{\operatorname{
Y}}(z))_{12} \! = \! \int_{\mathbb{R}} \tfrac{\boldsymbol{\pi}_{2n}(s) \exp
(-n \widetilde{V}(s))}{s-z} \tfrac{\md s}{2 \pi \mi}$. Using (for $\vert s/z
\vert \! \ll \! 1)$ the expansion $\tfrac{1}{s-z} \! = \! -\sum_{k=0}^{l}
\tfrac{s^{k}}{z^{k+1}} \! + \! \tfrac{s^{l+1}}{z^{l+1}(s-z)}$, $l \! \in \!
\mathbb{Z}_{0}^{+}$, and recalling that $\langle \boldsymbol{\pi}_{2n},z^{j}
\rangle_{\mathscr{L}} \! = \! 0$, $j \! = \! -n,\dotsc,n \! - \! 1$, and
$\phi_{2n}(z) \! = \! \xi^{(2n)}_{n} \boldsymbol{\pi}_{2n}(z)$, one proceeds
as follows:
\begin{align*}
\left(\overset{e}{\operatorname{Y}}(z) \right)_{12} \underset{\underset{z \in
\mathbb{C} \setminus \mathbb{R}}{z \to \infty}}{=}& \, -\dfrac{1}{z} \int_{
\mathbb{R}} \boldsymbol{\pi}_{2n}(s) \! \left(1 \! + \! \dfrac{s}{z} \! + \!
\cdots \! + \! \dfrac{s^{n-1}}{z^{n-1}} \! + \! \dfrac{s^{n}}{z^{n}} \! + \!
\cdots \right) \! \me^{-n \widetilde{V}(s)} \, \dfrac{\md s}{2 \pi \mi} \\
\underset{\underset{z \in \mathbb{C} \setminus \mathbb{R}}{z \to \infty}}{=}&
\, -\dfrac{1}{z} \int_{\mathbb{R}} \boldsymbol{\pi}_{2n}(s) \! \left(\dfrac{
s^{n}}{z^{n}} \right) \! \me^{-n \widetilde{V}(s)} \, \dfrac{\md s}{2 \pi \mi}
\! + \! \mathcal{O} \! \left(\dfrac{1}{z^{n+2}} \right) \\
\underset{\underset{z \in \mathbb{C} \setminus \mathbb{R}}{z \to \infty}}{=}&
\, -\dfrac{z^{-(n+1)}}{\xi^{(2n)}_{n}} \int_{\mathbb{R}} \underbrace{\xi_{n}^{
(2n)} \boldsymbol{\pi}_{2n}(s)}_{= \, \phi_{2n}(s)} \dfrac{\me^{-n \widetilde{
V}(s)}}{\xi_{n}^{(2n)}} \! \underbrace{\left(\xi_{n}^{(2n)}s^{n} \! + \!
\cdots \! + \! \dfrac{\xi_{-n}^{(2n)}}{s^{n}} \right)}_{= \, \phi_{2n}(s)}
\dfrac{\md s}{2 \pi \mi} \\
+& \, \mathcal{O} \! \left(z^{-(n+2)} \right) \\
\underset{\underset{z \in \mathbb{C} \setminus \mathbb{R}}{z \to \infty}}{=}&
\, -\dfrac{z^{-(n+1)}}{2 \pi \mi (\xi_{n}^{(2n)})^{2}} \underbrace{\int_{
\mathbb{R}} \phi_{2n}(s) \phi_{2n}(s) \me^{-n \widetilde{V}(s)} \, \md s}_{=
\, 1}+ \, \mathcal{O} \! \left(z^{-(n+2)} \right) \Rightarrow \\
\left(\overset{e}{\operatorname{Y}}(z)z^{-n \sigma_{3}} \right)_{12} \underset{
\underset{z \in \mathbb{C} \setminus \mathbb{R}}{z \to \infty}}{=}& \, \dfrac{
1}{z} \! \left(-\dfrac{1}{2 \pi \mi (\xi_{n}^{(2n)})^{2}} \right) \! + \!
\mathcal{O} \! \left(\dfrac{1}{z^{2}} \right);
\end{align*}
but, noting {}from Proposition~5.3 that
\begin{equation*}
\left(\overset{e}{\operatorname{Y}}(z)z^{-n \sigma_{3}} \right)_{12} \underset{
\underset{z \in \mathbb{C} \setminus \mathbb{R}}{z \to \infty}}{=} \dfrac{1}{
z} \! \left(\operatorname{Y}^{e,\infty}_{1} \right)_{12} \! + \!\mathcal{O} \!
\left(\dfrac{1}{z^{2}} \right),
\end{equation*}
upon equating the above two asymptotic expansions for $(\overset{e}{
\operatorname{Y}}(z)z^{-n \sigma_{3}})_{12}$, one arrives at the result
stated in the Proposition. \hfill $\qed$

Using the results of Propositions~5.3 and~5.4, one obtains the $n \! \to \!
\infty$ asymptotics for $\xi^{(2n)}_{n}$ and $\phi_{2n}(z)$ (in the entire
complex plane) stated in Theorem~2.3.2.

Small-$z$ asymptotics for $\overset{e}{\mathrm{Y}}(z)$ are given in the 
Appendix (see Lemma~A.1): these latter asymptotics are necessary for the 
results of \cite{a52}.

\vspace*{1.50cm}
\textbf{\Large Acknowledgements}

K.~T.-R.~McLaughlin was supported, in part, by National Science Foundation
Grant Nos.~DMS--9970328 and~DMS--0200749. X.~Zhou was supported, in part, by
National Science Foundation Grant No.~DMS--0300844. The authors are grateful
to P.~Gonz\'{a}lez-Vera for comments related to the two-point Pad\'{e}
approximants error formula~$(\mathrm{TPA1})$.
\clearpage
\section*{Appendix: Small-$z$ Asymptotics for $\overset{e}{\mathrm{Y}}(z)$}
\setcounter{section}{1}
\setcounter{z0}{1}
Even though the results of Lemma~A.1 below, namely, small-$z$ asymptotics (as 
$(\mathbb{C} \setminus \mathbb{R} \! \ni)$ $z \! \to \! 0)$ of $\overset{e}{
\operatorname{Y}}(z)$, are not necessary in order to prove Theorems~2.3.1 
and~2.3.2, they are essential for the results of \cite{a52}, related to 
asymptotics of the coefficients of the system of three- and five-term 
recurrence relations and the corresponding Laurent-Jacobi matrices (cf. 
Section~1). For the sake of completeness, therefore, and in order to eschew 
any duplication of the analysis of this paper, $(\mathbb{C} \setminus 
\mathbb{R} \! \ni)$ $z \! \to \! 0$ asymptotics for $\overset{e}{
\operatorname{Y}}(z)$ are presented here.
\begin{ay}
Let $\mathscr{R}^{e} \colon \mathbb{C} \setminus \widetilde{\Sigma}_{p}^{e} 
\! \to \! \operatorname{SL}_{2}(\mathbb{C})$ be the solution of the {\rm RHP} 
$(\mathscr{R}^{e}(z),\upsilon_{\mathscr{R}}^{e}(z),\widetilde{\Sigma}_{p}^{
e})$ formulated in Proposition~{\rm 5.2} with $n \! \to \! \infty$ asymptotics 
given in Lemma~{\rm 5.3}. Then,
\begin{equation*}
\mathscr{R}^{e}(z) \underset{z \to 0}{=} \mathrm{I} \! + \! \mathscr{R}^{e,
0}_{0}(n) \! + \! \mathscr{R}^{e,0}_{1}(n)z \! + \! \mathscr{R}^{e,0}_{2}
(n)z^{2} \! + \! \mathcal{O}(z^{3}),
\end{equation*}
where, for $k \! = \! 1,2,3$,
\begin{equation*}
\mathscr{R}^{e,0}_{k-1}(n) \! := \! \int_{\Sigma_{\circlearrowright}^{e}}s^{-k}
w^{\Sigma^{e}_{\circlearrowright}}_{+}(s) \, \dfrac{\md s}{2 \pi \mi} \! = \!
-\sum_{j=1}^{N+1} \, \sum_{q \in \{b_{j-1}^{e},a_{j}^{e}\}} \operatorname{Res}
\! \left(z^{-k}w^{\Sigma^{e}_{\circlearrowright}}_{+}(z);q \right),
\end{equation*}
with, in particular,
\begin{align*}
\mathscr{R}^{e,0}_{k-1}(n) \underset{n \to \infty}{=}& \, \dfrac{1}{n} \sum_{
j=1}^{N+1} \! \left(\dfrac{(\mathscr{A}^{e}(b_{j-1}^{e})(\widehat{\alpha}_{1}^{
e}(b_{j-1}^{e}) \! + \! k(b_{j-1}^{e})^{-1} \widehat{\alpha}_{0}^{e}(b_{j-1}^{
e})) \! - \! \mathscr{B}^{e}(b_{j-1}^{e}) \widehat{\alpha}_{0}^{e}(b_{j-1}^{
e}))}{(b_{j-1}^{e})^{k}(\widehat{\alpha}_{0}^{e}(b_{j-1}^{e}))^{2}} \right. \\
+&\left. \, \dfrac{(\mathscr{A}^{e}(a_{j}^{e})(\widehat{\alpha}_{1}^{e}(a_{j}^{
e}) \! + \! k(a_{j}^{e})^{-1} \widehat{\alpha}_{0}^{e}(a_{j}^{e})) \! - \!
\mathscr{B}^{e}(a_{j}^{e}) \widehat{\alpha}_{0}^{e}(a_{j}^{e}))}{(a_{j}^{e})^{
k}(\widehat{\alpha}_{0}^{e}(a_{j}^{e}))^{2}} \right) \! + \! \mathcal{O} \!
\left(\dfrac{1}{n^{2}} \right),
\end{align*}
and all the parameters defined in Lemma~{\rm 5.3}.

Let $\overset{e}{m}^{\raise-1.0ex\hbox{$\scriptstyle \infty$}} \colon \mathbb{
C} \setminus J_{e}^{\infty} \! \to \! \operatorname{SL}_{2}(\mathbb{C})$ solve
the {\rm RHP} $(\overset{e}{m}^{\raise-1.0ex\hbox{$\scriptstyle \infty$}}(z),
J_{e}^{\infty},
\overset{e}{\upsilon}^{\raise-1.0ex\hbox{$\scriptstyle \infty$}}(z))$
formulated in Lemma~{\rm 4.3} with (unique) solution given by
Lemma~{\rm 4.5}. For $\varepsilon_{1},\varepsilon_{2} \! = \! \pm 1$, set
\begin{gather*}
\theta^{e}_{0}(\varepsilon_{1},\varepsilon_{2},\boldsymbol{\Omega}^{e}) \! :=
\! \boldsymbol{\theta}^{e}(\varepsilon_{1} \boldsymbol{u}^{e}_{+}(0) \! - \!
\tfrac{n}{2 \pi} \boldsymbol{\Omega}^{e} \! + \! \varepsilon_{2} \boldsymbol{
d}_{e}),
\end{gather*}
where $\boldsymbol{u}^{e}_{+}(0) \! = \! \int_{a_{N+1}^{e}}^{0^{+}}
\boldsymbol{\omega}^{e}$ $(0^{+} \! \in \! \mathbb{C}_{+})$,
\begin{gather*}
\widetilde{\alpha}^{e}_{0}(\varepsilon_{1},\varepsilon_{2},\boldsymbol{
\Omega}^{e}) \! := \! 2 \pi \mi \varepsilon_{1} \sum_{m \in \mathbb{Z}^{N}}
(m,\widehat{\boldsymbol{\alpha}}^{e}_{0}) \me^{2 \pi \mi (m,\varepsilon_{1}
\boldsymbol{u}^{e}_{+}(0)-\frac{n}{2 \pi} \boldsymbol{\Omega}^{e}+\varepsilon_{
2} \boldsymbol{d}_{e})+ \pi \mi (m,\tau^{e}m)},
\end{gather*}
where $\widehat{\boldsymbol{\alpha}}^{e}_{0} \! = \! (\widehat{\alpha}_{0,
1}^{e},\widehat{\alpha}^{e}_{0,2},\dotsc,\widehat{\alpha}^{e}_{0,N})$, with
$\widehat{\alpha}^{e}_{0,j} \! := \! (-1)^{\mathcal{N}_{+}}(\prod_{i=1}^{N+1}
\vert b_{i-1}^{e}a_{i}^{e} \vert)^{-1/2}c_{jN}^{e}$, $j \! = \! 1,\dotsc,N$,
where $\mathcal{N}_{+} \! \in \! \lbrace 0,\dotsc,N \! + \! 1 \rbrace$ is
the number of bands to the right of $z \! = \! 0$,
\begin{gather*}
\beta^{e}_{0}(\varepsilon_{1},\varepsilon_{2},\boldsymbol{\Omega}^{e}) \! :=
\! 2 \pi \sum_{m \in \mathbb{Z}^{N}} \! \left(\mi \varepsilon_{1}(m,\widehat{
\boldsymbol{\beta}}^{e}_{0}) \! - \!  \pi (m,\widehat{\boldsymbol{\alpha}}^{
e}_{0})^{2} \right) \! \me^{2 \pi \mi (m,\varepsilon_{1} \boldsymbol{u}^{e}_{+}
(0)-\frac{n}{2 \pi} \boldsymbol{\Omega}^{e}+\varepsilon_{2} \boldsymbol{d}_{e})
+ \pi \mi (m,\tau^{e}m)},
\end{gather*}
where $\widehat{\boldsymbol{\beta}}^{e}_{0} \! = \! (\widehat{\beta}^{e}_{0,1},
\widehat{\beta}^{e}_{0,2},\dotsc,\widehat{\beta}^{e}_{0,N})$, with $\widehat{
\beta}^{e}_{0,j} \! := \! \tfrac{1}{2}(-1)^{\mathcal{N}_{+}}(\prod_{i=1}^{N+
1} \vert b_{i-1}^{e}a_{i}^{e} \vert)^{-1/2}(c^{e}_{jN-1} \! + \! \tfrac{1}{2}
c_{jN}^{e} \sum_{k=1}^{N+1}((a_{k}^{e})^{-1} \! + \! (b_{k-1}^{e})^{-1}))$,
$j \! = \! 1,\dotsc,N$, where $c_{jN}^{e},c_{jN-1}^{e}$, $j \! = \! 1,\dotsc,
N$, are obtained {}from Equations~{\rm (E1)} and~{\rm (E2)}. Set $\gamma^{e}_{
0} \! := \! \gamma^{e}(0) \! = \! (\prod_{k=1}^{N+1}b_{k-1}^{e}(a_{k}^{e})^{-
1})^{1/4}$ $(> \! 0)$. Then,
\begin{gather*}
\overset{e}{m}^{\raise-1.0ex\hbox{$\scriptstyle \infty$}}(z) \underset{z \to
0}{=} \overset{e}{m}_{0}^{\raise-1.0ex\hbox{$\scriptstyle 0$}} \! + \! z
\overset{e}{m}_{1}^{\raise-1.0ex\hbox{$\scriptstyle 0$}} \! + \! z^{2}
\overset{e}{m}_{2}^{\raise-1.0ex\hbox{$\scriptstyle 0$}} \! + \!
\mathcal{O}(z^{3}),
\end{gather*}
where
\begin{align*}
(\overset{e}{m}_{0}^{\raise-1.0ex\hbox{$\scriptstyle 0$}})_{11} &= \dfrac{
\boldsymbol{\theta}^{e}(\boldsymbol{u}^{e}_{+}(\infty) \! + \! \boldsymbol{d}_{
e})}{\boldsymbol{\theta}^{e}(\boldsymbol{u}^{e}_{+}(\infty) \! - \! \frac{n}{2
\pi} \boldsymbol{\Omega}^{e} \! + \! \boldsymbol{d}_{e})} \! \left(\dfrac{
\gamma^{e}_{0} \! + \! (\gamma^{e}_{0})^{-1}}{2} \right) \! \dfrac{\theta^{
e}_{0}(1,1,\boldsymbol{\Omega}^{e})}{\theta^{e}_{0}(1,1,\vec{\pmb{0}})}, \\
(\overset{e}{m}_{0}^{\raise-1.0ex\hbox{$\scriptstyle 0$}})_{12} &= -\dfrac{
\boldsymbol{\theta}^{e}(\boldsymbol{u}^{e}_{+}(\infty) \! + \! \boldsymbol{d}_{
e})}{\boldsymbol{\theta}^{e}(\boldsymbol{u}^{e}_{+}(\infty) \! - \! \frac{n}{2
\pi} \boldsymbol{\Omega}^{e} \! + \! \boldsymbol{d}_{e})} \! \left(\dfrac{
\gamma^{e}_{0} \! - \! (\gamma^{e}_{0})^{-1}}{2 \mi} \right) \! \dfrac{
\theta^{e}_{0}(-1,1,\boldsymbol{\Omega}^{e})}{\theta^{e}_{0}(-1,1,\vec{\pmb{0}
})}, \\
(\overset{e}{m}_{0}^{\raise-1.0ex\hbox{$\scriptstyle 0$}})_{21} &= \dfrac{
\boldsymbol{\theta}^{e}(\boldsymbol{u}^{e}_{+}(\infty) \! + \! \boldsymbol{d}_{
e})}{\boldsymbol{\theta}^{e}(-\boldsymbol{u}^{e}_{+}(\infty) \! - \! \frac{n}{
2 \pi} \boldsymbol{\Omega}^{e} \! - \! \boldsymbol{d}_{e})} \! \left(\dfrac{
\gamma^{e}_{0} \! - \! (\gamma^{e}_{0})^{-1}}{2 \mi} \right) \! \dfrac{\theta^{
e}_{0}(1,-1,\boldsymbol{\Omega}^{e})}{\theta^{e}_{0}(1,-1,\vec{\pmb{0}})}, \\
(\overset{e}{m}_{0}^{\raise-1.0ex\hbox{$\scriptstyle 0$}})_{22} &= \dfrac{
\boldsymbol{\theta}^{e}(\boldsymbol{u}^{e}_{+}(\infty) \! + \! \boldsymbol{d}_{
e})}{\boldsymbol{\theta}^{e}(-\boldsymbol{u}^{e}_{+}(\infty) \! - \! \frac{n}{
2 \pi} \boldsymbol{\Omega}^{e} \! - \! \boldsymbol{d}_{e})} \! \left(\dfrac{
\gamma^{e}_{0} \! + \! (\gamma^{e}_{0})^{-1}}{2} \right) \! \dfrac{\theta^{e}_{
0}(-1,-1,\boldsymbol{\Omega}^{e})}{\theta^{e}_{0}(-1,-1,\vec{\pmb{0}})}, \\
(\overset{e}{m}_{1}^{\raise-1.0ex\hbox{$\scriptstyle 0$}})_{11} &= \dfrac{
\boldsymbol{\theta}^{e}(\boldsymbol{u}^{e}_{+}(\infty) \! + \! \boldsymbol{d}_{
e})}{\boldsymbol{\theta}^{e}(\boldsymbol{u}^{e}_{+}(\infty) \! - \! \frac{n}{2
\pi} \boldsymbol{\Omega}^{e} \! + \! \boldsymbol{d}_{e})} \! \left(\! \left(\!
\dfrac{\widetilde{\alpha}_{0}^{e}(1,1,\boldsymbol{\Omega}^{e}) \theta^{e}_{0}
(1,1,\vec{\pmb{0}}) \! - \! \widetilde{\alpha}_{0}^{e}(1,1,\vec{\pmb{0}})
\theta^{e}_{0}(1,1,\boldsymbol{\Omega}^{e})}{(\theta^{e}_{0}(1,1,\vec{\pmb{0}}
))^{2}} \right) \right. \\
&\left. \times \left(\dfrac{\gamma^{e}_{0} \! + \! (\gamma^{e}_{0})^{-1}}{2}
\right) \! + \! \left(\dfrac{\gamma^{e}_{0} \! - \! (\gamma^{e}_{0})^{-1}}{8}
\right) \! \left(\sum_{k=1}^{N+1} \! \left(\dfrac{1}{a_{k}^{e}} \! - \! \dfrac{
1}{b_{k-1}^{e}} \right) \right) \! \dfrac{\theta^{e}_{0}(1,1,\boldsymbol{
\Omega}^{e})}{\theta^{e}_{0}(1,1,\vec{\pmb{0}})} \right), \\
(\overset{e}{m}_{1}^{\raise-1.0ex\hbox{$\scriptstyle 0$}})_{12} &= -\dfrac{
\boldsymbol{\theta}^{e}(\boldsymbol{u}^{e}_{+}(\infty) \! + \! \boldsymbol{d}_{
e})}{\boldsymbol{\theta}^{e}(\boldsymbol{u}^{e}_{+}(\infty) \! - \! \frac{n}{2
\pi} \boldsymbol{\Omega}^{e} \! + \! \boldsymbol{d}_{e})} \! \left(\! \left(\!
\dfrac{\widetilde{\alpha}_{0}^{e}(-1,1,\boldsymbol{\Omega}^{e}) \theta^{e}_{0}
(-1,1,\vec{\pmb{0}}) \! - \! \widetilde{\alpha}_{0}^{e}(-1,1,\vec{\pmb{0}})
\theta^{e}_{0}(-1,1,\boldsymbol{\Omega}^{e})}{(\theta^{e}_{0}(-1,1,\vec{\pmb{0}
}))^{2}} \right) \right. \\
&\left. \times \left(\dfrac{\gamma^{e}_{0} \! - \! (\gamma^{e}_{0})^{-1}}{2
\mi} \right) \! + \! \left(\dfrac{\gamma^{e}_{0} \! + \! (\gamma^{e}_{0})^{-1}
}{8 \mi} \right) \! \left(\sum_{k=1}^{N+1} \! \left(\dfrac{1}{a_{k}^{e}} \! -
\! \dfrac{1}{b_{k-1}^{e}} \right) \right) \! \dfrac{\theta^{e}_{0}(-1,1,
\boldsymbol{\Omega}^{e})}{\theta^{e}_{0}(-1,1,\vec{\pmb{0}})} \right), \\
(\overset{e}{m}_{1}^{\raise-1.0ex\hbox{$\scriptstyle 0$}})_{21} &= \dfrac{
\boldsymbol{\theta}^{e}(\boldsymbol{u}^{e}_{+}(\infty) \! + \! \boldsymbol{d}_{
e})}{\boldsymbol{\theta}^{e}(-\boldsymbol{u}^{e}_{+}(\infty) \! - \! \frac{n}{
2 \pi} \boldsymbol{\Omega}^{e} \! - \! \boldsymbol{d}_{e})} \! \left(\! \left(
\! \dfrac{\widetilde{\alpha}_{0}^{e}(1,-1,\boldsymbol{\Omega}^{e}) \theta^{e}_{
0}(1,-1,\vec{\pmb{0}}) \! - \! \widetilde{\alpha}_{0}^{e}(1,-1,\vec{\pmb{0}})
\theta^{e}_{0}(1,-1,\boldsymbol{\Omega}^{e})}{(\theta^{e}_{0}(1,-1,\vec{\pmb{
0}}))^{2}} \right) \right. \\
&\left. \times \left(\dfrac{\gamma^{e}_{0} \! - \! (\gamma^{e}_{0})^{-1}}{2
\mi} \right) \! + \! \left(\dfrac{\gamma^{e}_{0} \! + \! (\gamma^{e}_{0})^{-1}
}{8 \mi} \right) \! \left(\sum_{k=1}^{N+1} \! \left(\dfrac{1}{a_{k}^{e}} \!
- \! \dfrac{1}{b_{k-1}^{e}} \right) \right) \! \dfrac{\theta^{e}_{0}(1,-1,
\boldsymbol{\Omega}^{e})}{\theta^{e}_{0}(1,-1,\vec{\pmb{0}})} \right), \\
(\overset{e}{m}_{1}^{\raise-1.0ex\hbox{$\scriptstyle 0$}})_{22} &= \dfrac{
\boldsymbol{\theta}^{e}(\boldsymbol{u}^{e}_{+}(\infty) \! + \! \boldsymbol{d}_{
e})}{\boldsymbol{\theta}^{e}(-\boldsymbol{u}^{e}_{+}(\infty) \! - \! \frac{n}{
2 \pi} \boldsymbol{\Omega}^{e} \! - \! \boldsymbol{d}_{e})} \! \left(\! \left(
\! \dfrac{\widetilde{\alpha}_{0}^{e}(-1,-1,\boldsymbol{\Omega}^{e}) \theta^{
e}_{0}(-1,-1,\vec{\pmb{0}}) \! - \! \widetilde{\alpha}_{0}^{e}(-1,-1,\vec{\pmb{
0}}) \theta^{e}_{0}(-1,-1,\boldsymbol{\Omega}^{e})}{(\theta^{e}_{0}(-1,-1,\vec{
\pmb{0}}))^{2}} \right) \right. \\
&\left. \times \left(\dfrac{\gamma^{e}_{0} \! + \! (\gamma^{e}_{0})^{-1}}{2}
\right) \! + \! \left(\dfrac{\gamma^{e}_{0} \! - \! (\gamma^{e}_{0})^{-1}}{8}
\right) \! \left(\sum_{k=1}^{N+1} \! \left(\dfrac{1}{a_{k}^{e}} \! - \! \dfrac{
1}{b_{k-1}^{e}} \right) \right) \! \dfrac{\theta^{e}_{0}(-1,-1,\boldsymbol{
\Omega}^{e})}{\theta^{e}_{0}(-1,-1,\vec{\pmb{0}})} \right), \\
(\overset{e}{m}_{2}^{\raise-1.0ex\hbox{$\scriptstyle 0$}})_{11} &= \dfrac{
\boldsymbol{\theta}^{e}(\boldsymbol{u}^{e}_{+}(\infty) \! + \! \boldsymbol{d}_{
e})}{\boldsymbol{\theta}^{e}(\boldsymbol{u}^{e}_{+}(\infty) \! - \! \frac{n}{2
\pi} \boldsymbol{\Omega}^{e} \! + \! \boldsymbol{d}_{e})} \! \left(\! \left(
\theta^{e}_{0}(1,1,\boldsymbol{\Omega}^{e}) \! \left((\widetilde{\alpha}_{0}^{
e}(1,1,\vec{\pmb{0}}))^{2} \! - \! \beta^{e}_{0}(1,1,\vec{\pmb{0}}) \theta^{
e}_{0}(1,1,\vec{\pmb{0}}) \right) \right. \right. \\
&\left. \left. - \, \widetilde{\alpha}_{0}^{e}(1,1,\boldsymbol{\Omega}^{e})
\widetilde{\alpha}_{0}^{e}(1,1,\vec{\pmb{0}}) \theta^{e}_{0}(1,1,\vec{\pmb{0}}
) \! + \! \beta^{e}_{0}(1,1,\boldsymbol{\Omega}^{e})(\theta^{e}_{0}(1,1,\vec{
\pmb{0}}))^{2} \right) \! \tfrac{1}{(\theta^{e}_{0}(1,1,\vec{\pmb{0}}))^{3}}
\right. \\
&\left. \times \left(\dfrac{\gamma^{e}_{0} \! + \! (\gamma^{e}_{0})^{-1}}{2}
\right) \! + \! \left(\dfrac{\widetilde{\alpha}_{0}^{e}(1,1,\boldsymbol{
\Omega}^{e}) \theta^{e}_{0}(1,1,\vec{\pmb{0}}) \! - \! \widetilde{\alpha}_{0}^{
e}(1,1,\vec{\pmb{0}}) \theta^{e}_{0}(1,1,\boldsymbol{\Omega}^{e})}{(\theta^{
e}_{0}(1,1,\vec{\pmb{0}}))^{2}} \right) \right. \\
&\left. \times \left(\sum_{k=1}^{N+1} \! \left(\dfrac{1}{a_{k}^{e}} \! - \!
\dfrac{1}{b_{k-1}^{e}} \right) \right) \! \left(\dfrac{\gamma^{e}_{0} \! - \!
(\gamma^{e}_{0})^{-1}}{8} \right) \! + \! \left(\! \left(\dfrac{\gamma^{e}_{0}
\! - \! (\gamma^{e}_{0})^{-1}}{16} \right) \! \sum_{k=1}^{N+1} \! \left(
\dfrac{1}{(a_{k}^{e})^{2}} \! - \! \dfrac{1}{(b_{k-1}^{e})^{2}} \right)
\right. \right. \\
&\left. \left. +\left(\dfrac{\gamma^{e}_{0} \! + \! (\gamma^{e}_{0})^{-1}}{64}
\right) \! \left(\sum_{k=1}^{N+1} \! \left(\dfrac{1}{a_{k}^{e}} \! - \!
\dfrac{1}{b_{k-1}^{e}} \right) \right)^{2} \right) \! \dfrac{\theta^{e}_{0}
(1,1,\boldsymbol{\Omega}^{e})}{\theta^{e}_{0}(1,1,\vec{\pmb{0}})} \right), \\
(\overset{e}{m}_{2}^{\raise-1.0ex\hbox{$\scriptstyle 0$}})_{12} &= -\dfrac{
\boldsymbol{\theta}^{e}(\boldsymbol{u}^{e}_{+}(\infty) \! + \! \boldsymbol{d}_{
e})}{\boldsymbol{\theta}^{e}(\boldsymbol{u}^{e}_{+}(\infty) \! - \! \frac{n}{2
\pi} \boldsymbol{\Omega}^{e} \! + \! \boldsymbol{d}_{e})} \! \left(\! \left(
\theta^{e}_{0}(-1,1,\boldsymbol{\Omega}^{e}) \! \left((\widetilde{\alpha}_{0}^{
e}(-1,1,\vec{\pmb{0}}))^{2} \! - \! \beta^{e}_{0}(-1,1,\vec{\pmb{0}}) \theta^{
e}_{0}(-1,1,\vec{\pmb{0}}) \right) \right. \right. \\
&\left. \left. - \, \widetilde{\alpha}_{0}^{e}(-1,1,\boldsymbol{\Omega}^{e})
\widetilde{\alpha}_{0}^{e}(-1,1,\vec{\pmb{0}}) \theta^{e}_{0}(-1,1,\vec{\pmb{
0}}) \! + \! \beta^{e}_{0}(-1,1,\boldsymbol{\Omega}^{e})(\theta^{e}_{0}(-1,1,
\vec{\pmb{0}}))^{2} \right) \! \tfrac{1}{(\theta^{e}_{0}(-1,1,\vec{\pmb{0}})
)^{3}} \right. \\
&\left. \times \left(\dfrac{\gamma^{e}_{0} \! - \! (\gamma^{e}_{0})^{-1}}{2
\mi} \right) \! + \! \left(\dfrac{\widetilde{\alpha}_{0}^{e}(-1,1,\boldsymbol{
\Omega}^{e}) \theta^{e}_{0}(-1,1,\vec{\pmb{0}}) \! - \! \widetilde{\alpha}_{
0}^{e}(-1,1,\vec{\pmb{0}}) \theta^{e}_{0}(-1,1,\boldsymbol{\Omega}^{e})}{
(\theta^{e}_{0}(-1,1,\vec{\pmb{0}}))^{2}} \right) \right. \\
&\left. \times \left(\sum_{k=1}^{N+1} \! \left(\dfrac{1}{a_{k}^{e}} \! - \!
\dfrac{1}{b_{k-1}^{e}} \right) \right) \! \left(\dfrac{\gamma^{e}_{0} \! + \!
(\gamma^{e}_{0})^{-1}}{8 \mi} \right) \! + \! \left(\! \left(\dfrac{\gamma^{
e}_{0} \! + \! (\gamma^{e}_{0})^{-1}}{16 \mi} \right) \! \sum_{k=1}^{N+1} \!
\left(\dfrac{1}{(a_{k}^{e})^{2}} \! - \! \dfrac{1}{(b_{k-1}^{e})^{2}} \right)
\right. \right. \\
&\left. \left. +\left(\dfrac{\gamma^{e}_{0} \! - \! (\gamma^{e}_{0})^{-1}}{64
\mi} \right) \! \left(\sum_{k=1}^{N+1} \! \left(\dfrac{1}{a_{k}^{e}} \! - \!
\dfrac{1}{b_{k-1}^{e}} \right) \right)^{2} \right) \! \dfrac{\theta^{e}_{0}
(-1,1,\boldsymbol{\Omega}^{e})}{\theta^{e}_{0}(-1,1,\vec{\pmb{0}})} \right), \\
(\overset{e}{m}_{2}^{\raise-1.0ex\hbox{$\scriptstyle 0$}})_{21} &= \dfrac{
\boldsymbol{\theta}^{e}(\boldsymbol{u}^{e}_{+}(\infty) \! + \! \boldsymbol{d}_{
e})}{\boldsymbol{\theta}^{e}(-\boldsymbol{u}^{e}_{+}(\infty) \! - \! \frac{n}{
2 \pi} \boldsymbol{\Omega}^{e} \! - \! \boldsymbol{d}_{e})} \! \left(\! \left(
\theta^{e}_{0}(1,-1,\boldsymbol{\Omega}^{e}) \! \left((\widetilde{\alpha}_{0}^{
e}(1,-1,\vec{\pmb{0}}))^{2} \! - \! \beta^{e}_{0}(1,-1,\vec{\pmb{0}}) \theta^{
e}_{0}(1,-1,\vec{\pmb{0}}) \right) \right. \right. \\
&\left. \left. - \, \widetilde{\alpha}_{0}^{e}(1,-1,\boldsymbol{\Omega}^{e})
\widetilde{\alpha}_{0}^{e}(1,-1,\vec{\pmb{0}}) \theta^{e}_{0}(1,-1,\vec{\pmb{
0}}) \! + \! \beta^{e}_{0}(1,-1,\boldsymbol{\Omega}^{e})(\theta^{e}_{0}(1,-1,
\vec{\pmb{0}}))^{2} \right) \! \tfrac{1}{(\theta^{e}_{0}(1,-1,\vec{\pmb{0}})
)^{3}} \right. \\
&\left. \times \left(\dfrac{\gamma^{e}_{0} \! - \! (\gamma^{e}_{0})^{-1}}{2
\mi} \right) \! + \! \left(\dfrac{\widetilde{\alpha}_{0}^{e}(1,-1,\boldsymbol{
\Omega}^{e}) \theta^{e}_{0}(1,-1,\vec{\pmb{0}}) \! - \! \widetilde{\alpha}_{
0}^{e}(1,-1,\vec{\pmb{0}}) \theta^{e}_{0}(1,-1,\boldsymbol{\Omega}^{e})}{
(\theta^{e}_{0}(1,-1,\vec{\pmb{0}}))^{2}} \right) \right. \\
&\left. \times \left(\sum_{k=1}^{N+1} \! \left(\dfrac{1}{a_{k}^{e}} \! - \!
\dfrac{1}{b_{k-1}^{e}} \right) \right) \! \left(\dfrac{\gamma^{e}_{0} \! + \!
(\gamma^{e}_{0})^{-1}}{8 \mi} \right) \! + \! \left(\! \left(\dfrac{\gamma^{
e}_{0} \! + \! (\gamma^{e}_{0})^{-1}}{16 \mi} \right) \! \sum_{k=1}^{N+1} \!
\left(\dfrac{1}{(a_{k}^{e})^{2}} \! - \! \dfrac{1}{(b_{k-1}^{e})^{2}} \right)
\right. \right. \\
&\left. \left. +\left(\dfrac{\gamma^{e}_{0} \! - \! (\gamma^{e}_{0})^{-1}}{64
\mi} \right) \! \left(\sum_{k=1}^{N+1} \! \left(\dfrac{1}{a_{k}^{e}} \! - \!
\dfrac{1}{b_{k-1}^{e}} \right) \right)^{2} \right) \! \dfrac{\theta^{e}_{0}
(1,-1,\boldsymbol{\Omega}^{e})}{\theta^{e}_{0}(1,-1,\vec{\pmb{0}})} \right), \\
(\overset{e}{m}_{2}^{\raise-1.0ex\hbox{$\scriptstyle 0$}})_{22} &= \dfrac{
\boldsymbol{\theta}^{e}(\boldsymbol{u}^{e}_{+}(\infty) \! + \! \boldsymbol{d}_{
e})}{\boldsymbol{\theta}^{e}(-\boldsymbol{u}^{e}_{+}(\infty) \! - \! \frac{n}{
2 \pi} \boldsymbol{\Omega}^{e} \! - \! \boldsymbol{d}_{e})} \! \left(\! \left(
\theta^{e}_{0}(-1,-1,\boldsymbol{\Omega}^{e}) \! \left((\widetilde{\alpha}_{
0}^{e}(-1,-1,\vec{\pmb{0}}))^{2} \! - \! \beta^{e}_{0}(-1,-1,\vec{\pmb{0}})
\theta^{e}_{0}(-1,-1,\vec{\pmb{0}}) \right) \right. \right. \\
&\left. \left. - \, \widetilde{\alpha}_{0}^{e}(-1,-1,\boldsymbol{\Omega}^{e})
\widetilde{\alpha}_{0}^{e}(-1,-1,\vec{\pmb{0}}) \theta^{e}_{0}(-1,-1,\vec{\pmb{
0}}) \! + \! \beta^{e}_{0}(-1,-1,\boldsymbol{\Omega}^{e})(\theta^{e}_{0}(-1,-1,
\vec{\pmb{0}}))^{2} \right) \! \tfrac{1}{(\theta^{e}_{0}(-1,-1,\vec{\pmb{0}})
)^{3}} \right. \\
&\left. \times \left(\dfrac{\gamma^{e}_{0} \! + \! (\gamma^{e}_{0})^{-1}}{2}
\right) \! + \! \left(\dfrac{\widetilde{\alpha}_{0}^{e}(-1,-1,\boldsymbol{
\Omega}^{e}) \theta^{e}_{0}(-1,-1,\vec{\pmb{0}}) \! - \! \widetilde{\alpha}_{
0}^{e}(-1,-1,\vec{\pmb{0}}) \theta^{e}_{0}(-1,-1,\boldsymbol{\Omega}^{e})}{
(\theta^{e}_{0}(-1,-1,\vec{\pmb{0}}))^{2}} \right) \right. \\
&\left. \times \left(\sum_{k=1}^{N+1} \! \left(\dfrac{1}{a_{k}^{e}} \! - \!
\dfrac{1}{b_{k-1}^{e}} \right) \right) \! \left(\dfrac{\gamma^{e}_{0} \! - \!
(\gamma^{e}_{0})^{-1}}{8} \right) \! + \! \left(\! \left(\dfrac{\gamma^{e}_{0}
\! - \! (\gamma^{e}_{0})^{-1}}{16} \right) \! \sum_{k=1}^{N+1} \! \left(
\dfrac{1}{(a_{k}^{e})^{2}} \! - \! \dfrac{1}{(b_{k-1}^{e})^{2}} \right)
\right. \right. \\
&\left. \left. +\left(\dfrac{\gamma^{e}_{0} \! + \! (\gamma^{e}_{0})^{-1}}{64}
\right) \! \left(\sum_{k=1}^{N+1} \! \left(\dfrac{1}{a_{k}^{e}} \! - \!
\dfrac{1}{b_{k-1}^{e}} \right) \right)^{2} \right) \! \dfrac{\theta^{e}_{0}
(-1,-1,\boldsymbol{\Omega}^{e})}{\theta^{e}_{0}(-1,-1,\vec{\pmb{0}})} \right),
\end{align*}
with $(\star)_{ij}$, $i,j \! = \! 1,2$, denoting the $(i \, j)$-element of
$\star$, and $\vec{\pmb{0}} \! := \! (0,0,\dotsc,0)^{\mathrm{T}}$ $(\in \!
\mathbb{R}^{N})$. Set
\begin{align*}
(\widehat{Q}^{e}_{0})_{11} :=& \,
(\overset{e}{m}_{0}^{\raise-1.0ex\hbox{$\scriptstyle 0$}})_{11} \! \left(1 \!
+ \! (\mathscr{R}^{e,0}_{0}(n))_{11} \right) \! + \! (\mathscr{R}^{e,0}_{0}
(n))_{12}(\overset{e}{m}_{0}^{\raise-1.0ex\hbox{$\scriptstyle 0$}})_{21}, \\
(\widehat{Q}^{e}_{0})_{12} :=& \,
(\overset{e}{m}_{0}^{\raise-1.0ex\hbox{$\scriptstyle 0$}})_{12} \! \left(1 \!
+ \! (\mathscr{R}^{e,0}_{0}(n))_{11} \right) \! + \! (\mathscr{R}^{e,0}_{0}
(n))_{12}(\overset{e}{m}_{0}^{\raise-1.0ex\hbox{$\scriptstyle 0$}})_{22}, \\
(\widehat{Q}^{e}_{0})_{21} :=& \,
(\overset{e}{m}_{0}^{\raise-1.0ex\hbox{$\scriptstyle 0$}})_{21} \! \left(1 \!
+ \! (\mathscr{R}^{e,0}_{0}(n))_{22} \right) \! + \! (\mathscr{R}^{e,0}_{0}
(n))_{21}(\overset{e}{m}_{0}^{\raise-1.0ex\hbox{$\scriptstyle 0$}})_{11}, \\
(\widehat{Q}^{e}_{0})_{22} :=& \,
(\overset{e}{m}_{0}^{\raise-1.0ex\hbox{$\scriptstyle 0$}})_{22} \! \left(1 \!
+ \! (\mathscr{R}^{e,0}_{0}(n))_{22} \right) \! + \! (\mathscr{R}^{e,0}_{0}
(n))_{21}(\overset{e}{m}_{0}^{\raise-1.0ex\hbox{$\scriptstyle 0$}})_{12}, \\
(\widehat{Q}^{e}_{1})_{11} :=& \,
(\overset{e}{m}_{1}^{\raise-1.0ex\hbox{$\scriptstyle 0$}})_{11} \! \left(1 \!
+ \! (\mathscr{R}^{e,0}_{0}(n))_{11} \right) \! + \! (\mathscr{R}^{e,0}_{0}
(n))_{12}(\overset{e}{m}_{1}^{\raise-1.0ex\hbox{$\scriptstyle 0$}})_{21} \! +
\! (\mathscr{R}^{e,0}_{1}(n))_{11}
(\overset{e}{m}_{0}^{\raise-1.0ex\hbox{$\scriptstyle 0$}})_{11} \\
+& \, (\mathscr{R}^{e,0}_{1}(n))_{12}
(\overset{e}{m}_{0}^{\raise-1.0ex\hbox{$\scriptstyle 0$}})_{21}, \\
(\widehat{Q}^{e}_{1})_{12} :=& \,
(\overset{e}{m}_{1}^{\raise-1.0ex\hbox{$\scriptstyle 0$}})_{12} \! \left(1 \!
+ \! (\mathscr{R}^{e,0}_{0}(n))_{11} \right) \! + \! (\mathscr{R}^{e,0}_{0}
(n))_{12}(\overset{e}{m}_{1}^{\raise-1.0ex\hbox{$\scriptstyle 0$}})_{22} \! +
\! (\mathscr{R}^{e,0}_{1}(n))_{11}
(\overset{e}{m}_{0}^{\raise-1.0ex\hbox{$\scriptstyle 0$}})_{12} \\
+& \, (\mathscr{R}^{e,0}_{1}(n))_{12}
(\overset{e}{m}_{0}^{\raise-1.0ex\hbox{$\scriptstyle 0$}})_{22}, \\
(\widehat{Q}^{e}_{1})_{21} :=& \,
(\overset{e}{m}_{1}^{\raise-1.0ex\hbox{$\scriptstyle 0$}})_{21} \! \left(1 \!
+ \! (\mathscr{R}^{e,0}_{0}(n))_{22} \right) \! + \! (\mathscr{R}^{e,0}_{0}
(n))_{21}(\overset{e}{m}_{1}^{\raise-1.0ex\hbox{$\scriptstyle 0$}})_{11} \! +
\! (\mathscr{R}^{e,0}_{1}(n))_{21}
(\overset{e}{m}_{0}^{\raise-1.0ex\hbox{$\scriptstyle 0$}})_{11} \\
+& \, (\mathscr{R}^{e,0}_{1}(n))_{22}
(\overset{e}{m}_{0}^{\raise-1.0ex\hbox{$\scriptstyle 0$}})_{21}, \\
(\widehat{Q}^{e}_{1})_{22} :=& \,
(\overset{e}{m}_{1}^{\raise-1.0ex\hbox{$\scriptstyle 0$}})_{22} \! \left(1 \!
+ \! (\mathscr{R}^{e,0}_{0}(n))_{22} \right) \! + \! (\mathscr{R}^{e,0}_{0}
(n))_{21}(\overset{e}{m}_{1}^{\raise-1.0ex\hbox{$\scriptstyle 0$}})_{12} \! +
\! (\mathscr{R}^{e,0}_{1}(n))_{21}
(\overset{e}{m}_{0}^{\raise-1.0ex\hbox{$\scriptstyle 0$}})_{12} \\
+& \, (\mathscr{R}^{e,0}_{1}(n))_{22}
(\overset{e}{m}_{0}^{\raise-1.0ex\hbox{$\scriptstyle 0$}})_{22}, \\
(\widehat{Q}^{e}_{2})_{11} :=& \,
(\overset{e}{m}_{2}^{\raise-1.0ex\hbox{$\scriptstyle 0$}})_{11} \! \left(1 \!
+ \! (\mathscr{R}^{e,0}_{0}(n))_{11} \right) \! + \! (\mathscr{R}^{e,0}_{0}
(n))_{12}(\overset{e}{m}_{2}^{\raise-1.0ex\hbox{$\scriptstyle 0$}})_{21} \! +
\! (\mathscr{R}^{e,0}_{1}(n))_{11}
(\overset{e}{m}_{1}^{\raise-1.0ex\hbox{$\scriptstyle 0$}})_{11} \\
+& \, (\mathscr{R}^{e,0}_{1}(n))_{12}
(\overset{e}{m}_{1}^{\raise-1.0ex\hbox{$\scriptstyle 0$}})_{21} \! + \!
(\mathscr{R}^{e,0}_{2}(n))_{11}
(\overset{e}{m}_{0}^{\raise-1.0ex\hbox{$\scriptstyle 0$}})_{11} \! + \!
(\mathscr{R}^{e,0}_{2}(n))_{12}
(\overset{e}{m}_{0}^{\raise-1.0ex\hbox{$\scriptstyle 0$}})_{21}, \\
(\widehat{Q}^{e}_{2})_{12} :=& \,
(\overset{e}{m}_{2}^{\raise-1.0ex\hbox{$\scriptstyle 0$}})_{12} \! \left(1 \!
+ \! (\mathscr{R}^{e,0}_{0}(n))_{11} \right) \! + \! (\mathscr{R}^{e,0}_{0}
(n))_{12}(\overset{e}{m}_{2}^{\raise-1.0ex\hbox{$\scriptstyle 0$}})_{22} \! +
\! (\mathscr{R}^{e,0}_{1}(n))_{11}
(\overset{e}{m}_{1}^{\raise-1.0ex\hbox{$\scriptstyle 0$}})_{12} \\
+& \, (\mathscr{R}^{e,0}_{1}(n))_{12}
(\overset{e}{m}_{1}^{\raise-1.0ex\hbox{$\scriptstyle 0$}})_{22} \! + \!
(\mathscr{R}^{e,0}_{2}(n))_{11}
(\overset{e}{m}_{0}^{\raise-1.0ex\hbox{$\scriptstyle 0$}})_{12} \! + \!
(\mathscr{R}^{e,0}_{2}(n))_{12}
(\overset{e}{m}_{0}^{\raise-1.0ex\hbox{$\scriptstyle 0$}})_{22}, \\
(\widehat{Q}^{e}_{2})_{21} :=& \,
(\overset{e}{m}_{2}^{\raise-1.0ex\hbox{$\scriptstyle 0$}})_{21} \! \left(1 \!
+ \! (\mathscr{R}^{e,0}_{0}(n))_{22} \right) \! + \! (\mathscr{R}^{e,0}_{0}
(n))_{21}(\overset{e}{m}_{2}^{\raise-1.0ex\hbox{$\scriptstyle 0$}})_{11} \! +
\! (\mathscr{R}^{e,0}_{1}(n))_{21}
(\overset{e}{m}_{1}^{\raise-1.0ex\hbox{$\scriptstyle 0$}})_{11} \\
+& \, (\mathscr{R}^{e,0}_{1}(n))_{22}
(\overset{e}{m}_{1}^{\raise-1.0ex\hbox{$\scriptstyle 0$}})_{21} \! + \!
(\mathscr{R}^{e,0}_{2}(n))_{21}
(\overset{e}{m}_{0}^{\raise-1.0ex\hbox{$\scriptstyle 0$}})_{11} \! + \!
(\mathscr{R}^{e,0}_{2}(n))_{22}
(\overset{e}{m}_{0}^{\raise-1.0ex\hbox{$\scriptstyle 0$}})_{21}, \\
(\widehat{Q}^{e}_{2})_{22} :=& \,
(\overset{e}{m}_{2}^{\raise-1.0ex\hbox{$\scriptstyle 0$}})_{22} \! \left(1 \!
+ \! (\mathscr{R}^{e,0}_{0}(n))_{22} \right) \! + \! (\mathscr{R}^{e,0}_{0}
(n))_{21}(\overset{e}{m}_{2}^{\raise-1.0ex\hbox{$\scriptstyle 0$}})_{12} \! +
\! (\mathscr{R}^{e,0}_{1}(n))_{21}
(\overset{e}{m}_{1}^{\raise-1.0ex\hbox{$\scriptstyle 0$}})_{12} \\
+& \, (\mathscr{R}^{e,0}_{1}(n))_{22}
(\overset{e}{m}_{1}^{\raise-1.0ex\hbox{$\scriptstyle 0$}})_{22} \! + \!
(\mathscr{R}^{e,0}_{2}(n))_{21}
(\overset{e}{m}_{0}^{\raise-1.0ex\hbox{$\scriptstyle 0$}})_{12} \! + \!
(\mathscr{R}^{e,0}_{2}(n))_{22}
(\overset{e}{m}_{0}^{\raise-1.0ex\hbox{$\scriptstyle 0$}})_{22}.
\end{align*}
Let $\overset{e}{\mathrm{Y}} \colon \mathbb{C} \setminus \mathbb{R} \! \to \!
\operatorname{SL}_{2}(\mathbb{C})$ be the solution of {\rm \pmb{RHP1}}. Then,
\begin{equation*}
\overset{e}{\mathrm{Y}}(z)z^{n \sigma_{3}} \underset{z \to 0}{=} \mathrm{Y}^{
e,0}_{0} \! + \! z \mathrm{Y}^{e,0}_{1} \! + \! z^{2} \mathrm{Y}^{e,0}_{2} \!
+ \! \mathcal{O}(z^{3}),
\end{equation*}
where
\begin{align*}
(\mathrm{Y}^{e,0}_{0})_{11} =& \, (\widehat{Q}^{e}_{0})_{11} \me^{2n(\int_{J_{
e}} \ln (\lvert s \rvert) \psi_{V}^{e}(s) \, \md s + \mi \pi \int_{J_{e} \cap
\mathbb{R}_{+}} \psi_{V}^{e}(s) \, \md s)}, \\
(\mathrm{Y}^{e,0}_{0})_{12} =& \, (\widehat{Q}^{e}_{0})_{12} \me^{n(\ell_{e}
-2 \int_{J_{e}} \ln (\lvert s \rvert) \psi_{V}^{e}(s) \, \md s - 2 \pi \mi
\int_{J_{e} \cap \mathbb{R}_{+}} \psi_{V}^{e}(s) \, \md s)}, \\
(\mathrm{Y}^{e,0}_{0})_{21} =& \, (\widehat{Q}^{e}_{0})_{21} \me^{-n(\ell_{e}
-2 \int_{J_{e}} \ln (\lvert s \rvert) \psi_{V}^{e}(s) \, \md s - 2 \pi \mi
\int_{J_{e} \cap \mathbb{R}_{+}} \psi_{V}^{e}(s) \, \md s)}, \\
(\mathrm{Y}^{e,0}_{0})_{22} =& \, (\widehat{Q}^{e}_{0})_{22} \me^{-2n(\int_{J_{
e}} \ln (\lvert s \rvert) \psi_{V}^{e}(s) \, \md s + \mi \pi \int_{J_{e} \cap
\mathbb{R}_{+}} \psi_{V}^{e}(s) \, \md s)}, \\
(\mathrm{Y}^{e,0}_{1})_{11} =& \, \left(\! (\widehat{Q}^{e}_{1})_{11} \! - \!
2n(\widehat{Q}^{e}_{0})_{11} \int_{J_{e}}s^{-1} \psi_{V}^{e}(s) \, \md s
\right) \! \me^{2n(\int_{J_{e}} \ln (\lvert s \rvert) \psi_{V}^{e}(s) \, \md s
+ \mi \pi \int_{J_{e} \cap \mathbb{R}_{+}} \psi_{V}^{e}(s) \, \md s)}, \\
(\mathrm{Y}^{e,0}_{1})_{12} =& \, \left(\! (\widehat{Q}^{e}_{1})_{12} \! + \!
2n(\widehat{Q}^{e}_{0})_{12} \int_{J_{e}}s^{-1} \psi_{V}^{e}(s) \, \md s
\right) \! \me^{n(\ell_{e}-2 \int_{J_{e}} \ln (\lvert s \rvert) \psi_{V}^{e}
(s) \, \md s - 2 \pi \mi \int_{J_{e} \cap \mathbb{R}_{+}} \psi_{V}^{e}(s) \,
\md s)}, \\
(\mathrm{Y}^{e,0}_{1})_{21} =& \, \left(\! (\widehat{Q}^{e}_{1})_{21} \! - \!
2n(\widehat{Q}^{e}_{0})_{21} \int_{J_{e}}s^{-1} \psi_{V}^{e}(s) \, \md s
\right) \! \me^{-n(\ell_{e}-2 \int_{J_{e}} \ln (\lvert s \rvert) \psi_{V}^{e}
(s) \, \md s - 2 \pi \mi \int_{J_{e} \cap \mathbb{R}_{+}} \psi_{V}^{e}(s) \,
\md s)}, \\
(\mathrm{Y}^{e,0}_{1})_{22} =& \, \left(\! (\widehat{Q}^{e}_{1})_{22} \! + \!
2n(\widehat{Q}^{e}_{0})_{22} \int_{J_{e}}s^{-1} \psi_{V}^{e}(s) \, \md s
\right) \! \me^{-2n(\int_{J_{e}} \ln (\lvert s \rvert) \psi_{V}^{e}(s) \, \md
s + \mi \pi \int_{J_{e} \cap \mathbb{R}_{+}} \psi_{V}^{e}(s) \, \md s)}, \\
(\mathrm{Y}^{e,0}_{2})_{11} =& \left(\! (\widehat{Q}^{e}_{2})_{11} \! - \! 2n
(\widehat{Q}^{e}_{1})_{11} \int_{J_{e}}s^{-1} \psi_{V}^{e}(s) \, \md s \! + \!
(\widehat{Q}^{e}_{0})_{11} \! \left(2n^{2} \! \left(\int_{J_{e}}s^{-1} \psi_{
V}^{e}(s) \, \md s \right)^{2} \right. \right. \\
-& \left. \left. n \int_{J_{e}}s^{-2} \psi_{V}^{e}(s) \, \md s \right) \right)
\! \me^{2n(\int_{J_{e}} \ln (\lvert s \rvert) \psi_{V}^{e}(s) \, \md s + \mi
\pi \int_{J_{e} \cap \mathbb{R}_{+}} \psi_{V}^{e}(s) \, \md s)}, \\
(\mathrm{Y}^{e,0}_{2})_{12} =& \left(\! (\widehat{Q}^{e}_{2})_{12} \! + \! 2n
(\widehat{Q}^{e}_{1})_{12} \int_{J_{e}}s^{-1} \psi_{V}^{e}(s) \, \md s \! + \!
(\widehat{Q}^{e}_{0})_{12} \! \left(2n^{2} \! \left(\int_{J_{e}}s^{-1} \psi_{
V}^{e}(s) \, \md s \right)^{2} \right. \right. \\
+& \left. \left. n \int_{J_{e}}s^{-2} \psi_{V}^{e}(s) \, \md s \right)
\right) \! \me^{n(\ell_{e}-2 \int_{J_{e}} \ln (\lvert s \rvert) \psi_{V}^{e}
(s) \, \md s - 2 \pi \mi \int_{J_{e} \cap \mathbb{R}_{+}} \psi_{V}^{e}(s) \,
\md s)}, \\
(\mathrm{Y}^{e,0}_{2})_{21} =& \left(\! (\widehat{Q}^{e}_{2})_{21} \! - \! 2n
(\widehat{Q}^{e}_{1})_{21} \int_{J_{e}}s^{-1} \psi_{V}^{e}(s) \, \md s \! + \!
(\widehat{Q}^{e}_{0})_{21} \! \left(2n^{2} \! \left(\int_{J_{e}}s^{-1} \psi_{
V}^{e}(s) \, \md s \right)^{2} \right. \right. \\
-& \left. \left. n \int_{J_{e}}s^{-2} \psi_{V}^{e}(s) \, \md s \right)
\right) \! \me^{-n(\ell_{e}-2 \int_{J_{e}} \ln (\lvert s \rvert) \psi_{V}^{e}
(s) \, \md s - 2 \pi \mi \int_{J_{e} \cap \mathbb{R}_{+}} \psi_{V}^{e}(s) \,
\md s)}, \\
(\mathrm{Y}^{e,0}_{2})_{22} =& \left(\! (\widehat{Q}^{e}_{2})_{22} \! + \! 2n
(\widehat{Q}^{e}_{1})_{22} \int_{J_{e}}s^{-1} \psi_{V}^{e}(s) \, \md s \! + \!
(\widehat{Q}^{e}_{0})_{22} \! \left(2n^{2} \! \left(\int_{J_{e}}s^{-1} \psi_{
V}^{e}(s) \, \md s \right)^{2} \right. \right. \\
+& \left. \left. n \int_{J_{e}}s^{-2} \psi_{V}^{e}(s) \, \md s \right)
\right) \! \me^{-2n(\int_{J_{e}} \ln (\lvert s \rvert) \psi_{V}^{e}(s) \, \md
s + \mi \pi \int_{J_{e} \cap \mathbb{R}_{+}} \psi_{V}^{e}(s) \, \md s)}.
\end{align*}
\end{ay}

\emph{Proof.} Let $\mathscr{R}^{e} \colon \mathbb{C} \setminus \widetilde{
\Sigma}^{e}_{p} \! \to \! \operatorname{SL}_{2}(\mathbb{C})$ be the solution
of the RHP $(\mathscr{R}^{e}(z),\upsilon_{\mathscr{R}}^{e}(z),\widetilde{
\Sigma}^{e}_{p})$ formulated in Proposition~5.2 with $n \! \to \! \infty$
asymptotics given in Lemma~5.3. For $\vert z \vert \! \ll \! \min_{j=1,\dotsc,
N+1} \lbrace \vert b_{j-1}^{e} \! - \! a_{j}^{e} \vert \rbrace$, via the
expansion $\tfrac{1}{z-s} \! = \! -\sum_{k=0}^{l} \tfrac{z^{k}}{s^{k+1}} \! +
\! \tfrac{z^{l+1}}{s^{l+1}(z-s)}$, $l \! \in \! \mathbb{Z}_{0}^{+}$, where $s
\! \in \! \lbrace b_{j-1}^{e},a_{j}^{e} \rbrace$, $j \! = \! 1,\dotsc,N \! +
\! 1$, one obtains the asymptotics for $\mathscr{R}^{e}(z)$ stated in the
Proposition.

Let $\overset{e}{m}^{\raise-1.0ex\hbox{$\scriptstyle \infty$}} \colon \mathbb{
C} \setminus J_{e}^{\infty} \! \to \! \operatorname{SL}_{2}(\mathbb{C})$ solve
the RHP $(\overset{e}{m}^{\raise-1.0ex\hbox{$\scriptstyle \infty$}}(z),J_{e}^{
\infty},\overset{e}{\upsilon}^{\raise-1.0ex\hbox{$\scriptstyle \infty$}}(z))$
formulated in Lemma~4.3 with (unique) solution given by Lemma~4.5. In order
to obtain small-$z$ asymptotics of
$\overset{e}{m}^{\raise-1.0ex\hbox{$\scriptstyle \infty$}}(z)$, one needs
small-$z$ asymptotics of $(\gamma^{e}(z))^{\pm 1}$ and $\tfrac{\boldsymbol{
\theta}^{e}(\varepsilon_{1} \boldsymbol{u}^{e}(z)-\frac{n}{2 \pi} \boldsymbol{
\Omega}^{e}+\varepsilon_{2} \boldsymbol{d}_{e})}{\boldsymbol{\theta}^{e}
(\varepsilon_{1} \boldsymbol{u}^{e}(z)+\varepsilon_{2} \boldsymbol{d}_{e})}$,
$\varepsilon_{1},\varepsilon_{2} \! = \! \pm 1$. Consider, say, and without
loss of generality, $z \! \to \! 0$ asymptotics for $z \! \in \! \mathbb{C}_{
+}$ (designated $z \! \to \! 0^{+})$, where, by definition, $\sqrt{\smash[b]{
\star (z)}} \! := \! + \sqrt{\smash[b]{\star (z)}}$: equivalently, one may
consider $z \! \to \! 0$ asymptotics for $z \! \in \! \mathbb{C}_{-}$
(designated $z \! \to \! 0^{-})$; however, recalling that $\sqrt{\smash[b]{
\star (z)}} \! \upharpoonright_{\mathbb{C}_{+}} \! = \! -\sqrt{\smash[b]{\star
(z)}} \! \upharpoonright_{\mathbb{C}_{-}}$, one obtains (in either case, and
via the sheet-interchange index) the same $z \! \to \! 0$ asymptotics (for
$\overset{e}{m}^{\raise-1.0ex\hbox{$\scriptstyle \infty$}}(z))$. Recall the
expression for $\gamma^{e}(z)$ given in Lemma~4.4: for $\vert z \vert \! \ll
\! \min_{j=1,\dotsc,N+1} \lbrace \vert b_{j-1}^{e} \! - \! a_{j}^{e} \vert
\rbrace$, via the expansions $\tfrac{1}{z-s} \! = \! -\sum_{k=0}^{l} \tfrac{
z^{k}}{s^{k+1}} \! + \! \tfrac{z^{l+1}}{s^{l+1}(z-s)}$, $l \! \in \! \mathbb{
Z}_{0}^{+}$, and $\ln (s \! - \! z) \! =_{\vert z \vert \to 0} \! \ln (s) \!
- \! \sum_{k=1}^{\infty} \tfrac{1}{k}(\tfrac{z}{s})^{k}$, where $s \! \in \!
\lbrace b_{j-1}^{e},a_{j}^{e} \rbrace$, $j \! = \! 1,\dotsc,N \! + \! 1$, one
shows that, upon defining $\gamma^{e}(0)$ as in the Proposition,
\begin{align*}
(\gamma^{e}(z))^{\pm 1} \underset{z \to 0^{+}}{=}& \, (\gamma^{e}_{0})^{\pm
1} \! \left(1 \! + \! z \! \left(\pm \dfrac{1}{4} \sum_{k=1}^{N+1} \! \left(
\dfrac{1}{a_{k}^{e}} \! - \! \dfrac{1}{b_{k-1}^{e}} \right) \right) \! + \!
z^{2} \! \left(\pm \dfrac{1}{8} \sum_{k=1}^{N+1} \! \left(\dfrac{1}{(a_{k}^{
e})^{2}} \! - \! \dfrac{1}{(b_{k-1}^{e})^{2}} \right) \right. \right. \\
+&\left. \left. \dfrac{1}{32} \! \left(\sum_{k=1}^{N+1} \! \left(\dfrac{1}{a_{
k}^{e}} \! - \! \dfrac{1}{b_{k-1}^{e}} \right) \right)^{2} \, \right) \! + \!
\mathcal{O}(z^{3}) \right),
\end{align*}
whence
\begin{align*}
\dfrac{1}{2}(\gamma^{e}(z) \! + \! (\gamma^{e}(z))^{-1}) \underset{z \to 0^{
+}}{=}& \, \dfrac{(\gamma^{e}_{0} \! + \! (\gamma^{e}_{0})^{-1})}{2} \! + \!
z \! \left(\! \left(\dfrac{\gamma^{e}_{0} \! - \! (\gamma^{e}_{0})^{-1}}{8}
\right) \sum_{k=1}^{N+1} \! \left(\dfrac{1}{a_{k}^{e}} \! - \! \dfrac{1}{b_{k
-1}^{e}} \right) \right) \\
+& \, z^{2} \! \left(\! \left(\dfrac{\gamma^{e}_{0} \! - \! (\gamma^{e}_{0}
)^{-1}}{16} \right) \sum_{k=1}^{N+1} \! \left(\dfrac{1}{(a_{k}^{e})^{2}} \! -
\! \dfrac{1}{(b_{k-1}^{e})^{2}} \right) \! + \! \left(\dfrac{\gamma^{e}_{0} \!
+ \! (\gamma^{e}_{0})^{-1}}{64} \right) \right. \\
\times&\left. \left(\sum_{k=1}^{N+1} \! \left(\dfrac{1}{a_{k}^{e}} \! - \!
\dfrac{1}{b_{k-1}^{e}} \right) \right)^{2} \, \right) \! + \! \mathcal{O}
(z^{3}),
\end{align*}
and
\begin{align*}
\dfrac{1}{2 \mi}(\gamma^{e}(z) \! - \! (\gamma^{e}(z))^{-1}) \underset{z \to
0^{+}}{=}& \, \dfrac{(\gamma^{e}_{0} \! - \! (\gamma^{e}_{0})^{-1})}{2 \mi} \!
+ \! z \! \left(\! \left(\dfrac{\gamma^{e}_{0} \! + \! (\gamma^{e}_{0})^{-1}}{
8 \mi} \right) \sum_{k=1}^{N+1} \! \left(\dfrac{1}{a_{k}^{e}} \! - \! \dfrac{
1}{b_{k-1}^{e}} \right) \right) \\
+& \, z^{2} \! \left(\! \left(\dfrac{\gamma^{e}_{0} \! + \! (\gamma^{e}_{0})^{
-1}}{16 \mi} \right) \sum_{k=1}^{N+1} \! \left(\dfrac{1}{(a_{k}^{e})^{2}} \! -
\! \dfrac{1}{(b_{k-1}^{e})^{2}} \right) \! + \! \left(\dfrac{\gamma^{e}_{0} \!
- \! (\gamma^{e}_{0})^{-1}}{64 \mi} \right) \right. \\
\times&\left. \left(\sum_{k=1}^{N+1} \! \left(\dfrac{1}{a_{k}^{e}} \! - \!
\dfrac{1}{b_{k-1}^{e}} \right) \right)^{2} \, \right) \! + \! \mathcal{O}
(z^{3}).
\end{align*}
Recall {}from Lemma~4.5 that $\boldsymbol{u}^{e}(z) \! := \! \int_{a_{N+1}^{
e}}^{z} \boldsymbol{\omega}^{e}$ $(\in \operatorname{Jac}(\mathcal{Y}_{e})$,
with $\mathcal{Y}_{e} \! := \! \lbrace \mathstrut (y,z); \, y^{2} \! = \!
R_{e}(z) \rbrace)$, where $\boldsymbol{\omega}^{e}$, the associated normalised
basis of holomorphic one-forms of $\mathcal{Y}_{e}$, is given by $\boldsymbol{
\omega}^{e} \! = \! (\omega_{1}^{e},\omega_{2}^{e},\dotsc,\omega_{N}^{e})$,
with $\omega_{j}^{e} \! := \! \sum_{k=1}^{N}c_{jk}^{e}(\prod_{i=1}^{N+1}
(z \! - \! b_{i-1}^{e})(z \! - \! a_{i}^{e}))^{-1/2}z^{N-k} \, \md z$, $j \! 
= \! 1,\dotsc,N$, where $c_{jk}^{e}$, $j,k \! = \! 1,\dotsc,N$, are obtained 
{}from Equations~(E1) and~(E2). Writing
\begin{equation*}
\boldsymbol{u}^{e}(z) \! = \! \left(\int_{a_{N+1}^{e}}^{0^{+}} \! + \! \int_{
0^{+}}^{z} \right) \! \boldsymbol{\omega}^{e} \! = \! \boldsymbol{u}^{e}_{+}
(0) \! + \! \int_{0^{+}}^{z} \boldsymbol{\omega}^{e},
\end{equation*}
where $\boldsymbol{u}^{e}_{+}(0) \! := \! \int_{a_{N+1}^{e}}^{0^{+}} 
\boldsymbol{\omega}^{e}$, for $\vert z \vert \! \ll \! \min_{j=1,\dotsc,N+1} 
\lbrace \vert b_{j-1}^{e} \! - \! a_{j}^{e} \vert \rbrace$, via the expansions 
$\tfrac{1}{z-s} \! = \! -\sum_{k=0}^{l} \tfrac{z^{k}}{s^{k+1}} \! + \! \tfrac{
z^{l+1}}{s^{l+1}(z-s)}$, $l \! \in \! \mathbb{Z}_{0}^{+}$, and $\ln (s \! - \! 
z) \! =_{\vert z \vert \to 0} \! \ln (s) \! - \! \sum_{k=1}^{\infty} \tfrac{
1}{k}(\tfrac{z}{s})^{k}$, where $s \! \in \! \lbrace b_{k-1}^{e},a_{k}^{e} 
\rbrace$, $k \! = \! 1,\dotsc,N \! + \! 1$, one shows that, for $j \! = \! 1,
\dotsc,N$,
\begin{equation*}
\omega_{j}^{e} \underset{z \to 0^{+}}{=} (-1)^{\mathcal{N}_{+}} \! \left(
\prod_{i=1}^{N+1} \vert b_{i-1}^{e}a_{i}^{e} \vert \right)^{-1/2} \! \left(
c_{jN}^{e} \, \md z \! + \! \left(c_{jN-1}^{e} \! + \! \dfrac{c_{jN}^{e}}{2}
\sum_{k=1}^{N+1} \! \left(\dfrac{1}{b_{k-1}^{e}} \! + \! \dfrac{1}{a_{k}^{e}}
\right) \right) \! z \, \md z \! + \! \mathcal{O}(z^{2} \, \md z) \right),
\end{equation*}
where $\mathcal{N}_{+} \! \in \! \lbrace 0,\dotsc,N \! + \! 1 \rbrace$ is 
the number of bands to the right of $z \! = \! 0$, whence
\begin{align*}
\int_{0^{+}}^{z} \omega_{j}^{e} \underset{z \to 0^{+}}{=}& (-1)^{\mathcal{N}_{
+}} \! \left(\prod_{i=1}^{N+1} \vert b_{i-1}^{e}a_{i}^{e} \vert \right)^{-1/2}
\! \left(c_{jN}^{e}z \! + \! \dfrac{1}{2} \! \left(c_{jN-1}^{e} \! + \!
\dfrac{c_{jN}^{e}}{2} \sum_{k=1}^{N+1} \! \left(\dfrac{1}{b_{k-1}^{e}} \! +
\! \dfrac{1}{a_{k}^{e}} \right) \right) \! z^{2} \! + \! \mathcal{O}(z^{3})
\right) \\
=:& \, \widehat{\alpha}^{e}_{0,j}z \! + \! \widehat{\beta}^{e}_{0,j}z^{2} \! +
\! \mathcal{O}(z^{3}).
\end{align*}
Defining $\theta^{e}_{0}(\varepsilon_{1},\varepsilon_{2},\boldsymbol{\Omega}^{
e})$, $\widetilde{\alpha}^{e}_{0}(\varepsilon_{1},\varepsilon_{2},\boldsymbol{
\Omega}^{e})$, and $\beta^{e}_{0}(\varepsilon_{1},\varepsilon_{2},\boldsymbol{
\Omega}^{e})$, $\varepsilon_{1},\varepsilon_{2} \! = \! \pm 1$, as in the
Proposition, recalling that $\boldsymbol{\omega}^{e} \! = \! (\omega_{1}^{e},
\omega_{2}^{e},\dotsc,\omega_{N}^{e})$, and that the associated $N \! \times
\! N$ Riemann matrix of $\boldsymbol{\beta}^{e}$-periods, that is, $\tau^{e}
\! = \! (\tau^{e})_{i,j=1,\dotsc,N} \! := \! (\oint_{\boldsymbol{\beta}^{e}_{
j}} \omega^{e}_{i})_{i,j=1,\dotsc,N}$, is non-degenerate, symmetric, and $-
\mi \tau^{e}$ is positive definite, via the above asymptotic (as $z \! \to \!
0^{+})$ expansion for $\int_{0^{+}}^{z} \omega^{e}_{j}$, $j \! = \! 1,\dotsc,
N$, one shows that
\begin{equation*}
\dfrac{\boldsymbol{\theta}^{e}(\varepsilon_{1} \boldsymbol{u}^{e}(z) \! - \!
\frac{n}{2 \pi} \boldsymbol{\Omega}^{e} \! + \! \varepsilon_{2} \boldsymbol{
d}_{e})}{\boldsymbol{\theta}^{e}(\varepsilon_{1} \boldsymbol{u}^{e}(z) \! + \!
\varepsilon_{2} \boldsymbol{d}_{e})} \underset{z \to 0^{+}}{=} \digamma_{0}^{
e} \! + \! z \digamma_{1}^{e} \! + \! z^{2} \digamma_{2}^{e} \! + \! \mathcal{
O}(z^{3}),
\end{equation*}
where
\begin{align*}
\digamma_{0}^{e} :=& \, \dfrac{\theta^{e}_{0}(\varepsilon_{1},\varepsilon_{2},
\boldsymbol{\Omega}^{e})}{\theta^{e}_{0}(\varepsilon_{1},\varepsilon_{2},\vec{
\pmb{0}})}, \\
\digamma_{1}^{e} :=& \, \dfrac{\widetilde{\alpha}^{e}_{0}(\varepsilon_{1},
\varepsilon_{2},\boldsymbol{\Omega}^{e}) \theta^{e}_{0}(\varepsilon_{1},
\varepsilon_{2},\vec{\pmb{0}}) \! - \! \theta^{e}_{0}(\varepsilon_{1},
\varepsilon_{2},\boldsymbol{\Omega}^{e}) \widetilde{\alpha}^{e}_{0}
(\varepsilon_{1},\varepsilon_{2},\vec{\pmb{0}})}{(\theta^{e}_{0}(\varepsilon_{
1},\varepsilon_{2},\vec{\pmb{0}}))^{2}}, \\
\digamma^{e}_{2} :=& \left(\theta^{e}_{0}(\varepsilon_{1},\varepsilon_{2},
\boldsymbol{\Omega}^{e}) \! \left((\widetilde{\alpha}^{e}_{0}(\varepsilon_{1},
\varepsilon_{2},\vec{\pmb{0}}))^{2} \! - \! \beta^{e}_{0}(\varepsilon_{1},
\varepsilon_{2},\vec{\pmb{0}}) \theta^{e}_{0}(\varepsilon_{1},\varepsilon_{2},
\vec{\pmb{0}}) \right) \! - \! \widetilde{\alpha}^{e}_{0}(\varepsilon_{1},
\varepsilon_{2},\boldsymbol{\Omega}^{e}) \right. \\
\times&\left. \widetilde{\alpha}^{e}_{0}(\varepsilon_{1},\varepsilon_{2},\vec{
\pmb{0}}) \theta^{e}_{0}(\varepsilon_{1},\varepsilon_{2},\vec{\pmb{0}}) \! +
\! \beta^{e}_{0}(\varepsilon_{1},\varepsilon_{2},\boldsymbol{\Omega}^{e})
(\theta^{e}_{0}(\varepsilon_{1},\varepsilon_{2},\vec{\pmb{0}}))^{2} \right) \!
(\theta^{e}_{0}(\varepsilon_{1},\varepsilon_{2},\vec{\pmb{0}}))^{-3},
\end{align*}
with $\vec{\pmb{0}} \! := \! (0,0,\dotsc,0)^{\mathrm{T}}$ $(\in \! \mathbb{
R}^{N})$. Via the above asymptotic (as $z \! \to \! 0^{+})$ expansions for
$\tfrac{1}{2}(\gamma^{e}(z) \! + \! (\gamma^{e}(z))^{-1})$, $\tfrac{1}{2 \mi}
(\gamma^{e}(z) \! - \! (\gamma^{e}(z))^{-1})$, and $\tfrac{\boldsymbol{
\theta}^{e}(\varepsilon_{1} \boldsymbol{u}^{e}(z)-\frac{n}{2 \pi} \boldsymbol{
\Omega}^{e}+\varepsilon_{2} \boldsymbol{d}_{e})}{\boldsymbol{\theta}^{e}
(\varepsilon_{1} \boldsymbol{u}^{e}(z)+\varepsilon_{2} \boldsymbol{d}_{e})}$,
one arrives at, upon recalling the expression for
$\overset{e}{m}^{\raise-1.0ex\hbox{$\scriptstyle \infty$}}(z)$ given in
Lemma~4.5, the asymptotic expansion for
$\overset{e}{m}^{\raise-1.0ex\hbox{$\scriptstyle \infty$}}(z)$ stated in the
Proposition.

Let $\overset{e}{\operatorname{Y}} \colon \mathbb{C} \setminus \mathbb{R}
\! \to \! \operatorname{SL}_{2}(\mathbb{C})$ be the (unique) solution of
\textbf{RHP1}, that is, $(\overset{e}{\operatorname{Y}}(z),\mathrm{I} \! + \!
\exp (-n \widetilde{V}(z)) \sigma_{+},\mathbb{R})$. Recall, also, that, for
$z \! \in \! \Upsilon^{e}_{1} \cup \Upsilon^{e}_{2}$ (Figure~7),
\begin{equation*}
\overset{e}{\operatorname{Y}}(z) \! = \! \me^{\tfrac{n \ell_{e}}{2}
\operatorname{ad}(\sigma_{3})} \mathscr{R}^{e}(z)
\overset{e}{m}^{\raise-1.0ex\hbox{$\scriptstyle \infty$}}
(z) \me^{n(g^{e}(z)+\int_{J_{e}} \ln (s) \psi_{V}^{e}(s) \, \md s) \sigma_{3}}:
\end{equation*}
consider, say, and without loss of generality, small-$z$ asymptotics for
$\overset{e}{\operatorname{Y}}(z)$ for $z \! \in \! \Upsilon_{1}^{e}$.
Recalling {}from Lemma~3.4 that $g^{e}(z) \! := \! \int_{J_{e}} \ln ((z \! -
\! s)^{2}(zs)^{-1}) \psi_{V}^{e}(s) \, \md s$, $z \! \in \! \mathbb{C}
\setminus (-\infty,\max \{0,a_{N+1}^{e}\})$, for $\vert z \vert \! \ll \!
\min_{j=1,\dotsc,N+1} \lbrace \vert b_{j-1}^{e} \! - \! a_{j}^{e} \vert
\rbrace$, in particular, $\vert z/s \vert \! \ll \! 1$ with $s \! \in \!
J_{e}$, and noting that $\int_{J_{e}} \psi_{V}^{e}(s) \, \md s \! = \! 1$ and
$\int_{J_{e}}s^{-m} \psi_{V}^{e}(s) \, \md s \! < \! \infty$, $m \! \in \!
\mathbb{N}$, via the expansions $\tfrac{1}{z-s} \! = \! -\sum_{k=0}^{l}
\tfrac{z^{k}}{s^{k+1}} \! + \! \tfrac{z^{l+1}}{s^{l+1}(z-s)}$, $l \! \in \!
\mathbb{Z}_{0}^{+}$, and $\ln (s \! - \! z) \! =_{\vert z \vert \to 0} \!
\ln (s) \! - \! \sum_{k=1}^{\infty} \tfrac{1}{k}(\tfrac{z}{s})^{k}$, one
shows that
\begin{align*}
g^{e}(z) \underset{\mathbb{C}_{\pm} \ni z \to 0}{=}& \, -\ln (z) \! - \! 
Q_{e} \! + \! 2 \int_{J_{e}} \ln (\lvert s \rvert) \psi_{V}^{e}(s) \, \md s 
\! \pm \! 2 \pi \mi \int_{J_{e} \cap \mathbb{R}_{+}} \psi_{V}^{e}(s) \, \md 
s \\
+& \, z \! \left(-2 \int_{J_{e}}s^{-1} \psi_{V}^{e}(s) \, \md s \right) \! 
+ \! z^{2} \! \left(-\int_{J_{e}}s^{-2} \psi_{V}^{e}(s) \, \md s \right) \! 
+ \! \mathcal{O}(z^{3}),
\end{align*}
where $\int_{J_{e} \cap \mathbb{R}_{+}} \psi_{V}^{e}(s) \, \md s$ is given in
the proof of Lemma~3.4. (Explicit expressions for $\int_{J_{e}}s^{-k} \psi_{
V}^{e}(s) \, \md s$, $k \! = \! 1,2$, are given in Remark~3.2.) Using the
asymptotic (as $z \! \to \! 0)$ expansions for $g^{e}(z)$, $\mathscr{R}^{e}
(z)$, and $\overset{e}{m}^{\raise-1.0ex\hbox{$\scriptstyle \infty$}}(z)$
derived above, upon recalling the formula for $\overset{e}{\operatorname{
Y}}(z)$, one arrives at, after a matrix-multiplication argument, the
asymptotic expansion for $\overset{e}{\operatorname{Y}}(z)z^{n \sigma_{3}}$
stated in the Proposition. \hfill $\qed$

\end{document}